\documentclass[12pt]{article}
\usepackage[utf8]{inputenc}
\usepackage{amscd}
\usepackage{amsthm}
\usepackage{extarrows}
\usepackage{latexsym}
\usepackage{mathrsfs}
\usepackage{pst-all,multido,ifthen}
\usepackage{amsmath}
\usepackage{amssymb}
\usepackage[all]{xy}
\usepackage[english]{babel}
\usepackage{mathtools}
\usepackage{upgreek}
\usepackage{scalerel}
\usepackage{stackengine,wasysym}
\stackMath
\usepackage{tikz}
\newcommand*\circled[1]{\tikz[baseline=(char.base)]{
 \node[shape=circle,draw,inner sep=0.1pt] (char) {#1};}}
\newcommand\reallywidehat[1]{%
\savestack{\tmpbox}{\stretchto{%
 \scaleto{%
 \scalerel*[\widthof{\ensuremath{#1}}]{\kern-.6pt\bigwedge\kern-.6pt}%
 {\rule[-\textheight/2]{1ex}{\textheight}}
 }{\textheight}%
}{0.5ex}}%
\stackon[1pt]{#1}{\tmpbox}%
}
\newcommand\reallywidetilde[1]{\ThisStyle{%
 \setbox0=\hbox{$\SavedStyle#1$}%
 \stackengine{-.1\LMpt}{\SavedStyle#1}{%
 \stretchto{\scaleto{\SavedStyle\mkern.2mu\AC}{.5150\wd0}}{.6\ht0}%
 }{O}{c}{F}{T}{S}%
}}

\newcommand{\xdownarrow}[1]{%
 {\left\downarrow\vbox to #1{}\right.\kern-\nulldelimiterspace}
}

\newtheorem{thm}{Theorem}
\newtheorem{prop}{Proposition}
\newtheorem{cor}{Corollary}
\newtheorem{lemma}{Lemma}
\newtheorem{Fact}{Fact}
\newtheorem{claim}{Claim}

\theoremstyle{definition}
\newtheorem{df}{Definition}

\theoremstyle{remark}
\newtheorem{rem}{Remark}
\newtheorem{ex}{Example}
\newtheorem{prob}{Problem}
\def\endproof{$\hfill \square$}
\def\alt{\textup{alt}}
\def\alleven{\textup{all-even}}
\def\asydim{\textup{asy.dim}}
\def\asytrdeg{\textup{asy.tr.deg}}
\def\big-for{\textup{big-for}}
\def\can{\textup{can}}
\def\ch{\textup{char}}

\def\CL{\textup{CL}}
\def\CM{\textup{CM}}
\def\con{\textup{conj}}
\def\cong{\textup{cong}}
\def\C{\textup{C}}
\def\crys{\textup{crys}}
\def\cyc{\textup{Cyc}}
\def\deg{\textup{deg}}

\def\Det{\textup{Det}}
\def\Diag{\textup{Diag}}
\def\Diag{\textup{Diag}}

\def\dR{\textup{dR}}
\def\End{\textup{End}}
\def\even{\textup{even}}
\def\for{\textup{for}}
\def\Frac{\textup{Frac}}
\def\G{\textup{G}}
\def\GL{\textup{GL}}
\def\Gr{\textup{Gr}}
\def\GSp{\textup{GSp}}
\def\GT{\textup{GT}}
\def\Hom{\textup{Hom}}
\def\Index{\textup{Index}}
\def\Im{\textup{Im}}
\def\Ker{\textup{Ker}}
\def\L{\textup{L}}
\def\Lie{\textup{Lie}}
\def\M{\textup{M}}
\def\Maps{\textup{Maps}}

\def\n{\textup{n}}
\def\Mat{\textup{Mat}}
\def\O{\textup{O}}
\def\Ob{\textup{Ob}}
\def\odd{\textup{odd}}
\def\ord{\textup{ord}}
\def\Par{\textup{Par}}
\def\Perm{\textup{Perm}}
\def\PGL{\textup{PGL}}
\def\phidim{\textup{$\phi$-tr.deg}}
\def\phirank{\textup{$\phi$-rank}}
\def\Q{\textup{Q}}
\def\QCL{\textup{QCL}}
\def\Proj{\textup{Proj}}
\def\Prol{\widehat{\textup{\bf FilRing}_{\delta}}}
\def\rank{\textup{rank}}
\def\red{\textup{red}}
\def\s{\textup{s}}
\def\sep{\textup{sep}}
\def\S{\textup{S}}
\def\SL{\textup{SL}}
\def\small-can{\textup{small-can}}
\def\smooth{\textup{sm}}
\def\SO{\textup{SO}}
\def\Span{\textup{Span}}
\def\Spec{\textup{Spec}}
\def\Spf{\textup{Spf}}
\def\Spf{\textup{Spf}}
\def\Stab{\textup{Stab}}
\def\Sym{\textup{Sym}}
\def\t{\textup{t}}
\def\T{\textup{T}}
\def\torus{\textup{torus}}
\def\Trace{\textup{Tr}}
\def\trdeg{\textup{tr.deg}}
\def\U{\textup{U}}

\def\uni{\textup{uni}}
\def\univ{\textup{univ}}
\def\V{\textup{V}}
\def\W{\textup{W}}
\def\X{\textup{X}}
\def\Y{\textup{Y}}
\def\Z{\textup{Z}}
\def\f{\textup{\bf f}}
\def\tot{\textup{tot}}

\begin{document}

\title{The $\delta$-invariant theory of Hecke correspondences on $\mathcal A_g$}
\author{Alexandru Buium and Adrian Vasiu}
\maketitle
\noindent
{\bf ABSTRACT.} 
Let $p$ be a prime, let $N\geq 3$ be an integer prime to $p$, let $R$ be the ring of $p$-typical Witt vectors with coefficients in an algebraic closure of $\mathbb F_p$, and consider the correspondence $\mathcal A'_{g,1,N,R}\rightrightarrows \mathcal A_{g,1,N,R}$ obtained by taking the union of all prime to $p$ Hecke correspondences on Mumford's moduli scheme of principally polarized abelian schemes of relative dimension $g$ endowed with symplectic similitude level-$N$ structure over $R$-schemes. It is well-known that the coequalizer $\mathcal A_{g,1,N,R}/\mathcal A'_{g,1,N,R}$ of the above correspondence exists and is trivial in the category of schemes, i.e., is $\Spec(R)$. We construct and study in detail such a coequalizer (categorical quotient) in a more refined geometry (category) referred to as {\it $\delta$-geometry}. This geometry is in essence obtained from the usual algebraic geometry by equipping all $R$-algebras with {\it $p$-derivations}. In particular, we prove that our substitute of $\mathcal A_{g,1,N,R}/\mathcal A'_{g,1,N,R}$ in $\delta$-geometry has the same `dimension' as $\mathcal A_{g,1,N,R}$, thus solving a main open problem in the work of Barc\u{a}u--Buium. We also give applications to the study of various Zariski dense loci in $\mathcal A_{g,1,N,R}$ such as of isogeny classes and of points with complex multiplication. To prove our results we develop a Serre--Tate expansion theory for {\it Siegel $\delta$-modular forms} of arbitrary genus which we then combine with old and new results from the geometric invariant theory of multiple quadratic forms and of multiple endomorphisms.

\bigskip\noindent
{\bf KEY WORDS:} abelian scheme, group scheme, cohomology, $\delta$-modular form, field, flow, Fourier expansion, Hilbert series, invariant, Hecke covariant form, module, quasi-linearity, quotient, polarization, polynomial, reductive group, representation, ring, scheme, semi-simplicial object, Serre--Tate expansion, Siegel modular variety, stable point, and transcendence basis. 

\bigskip\noindent
{\bf MSC 2020:} 11E04, 11F06, 11F46, 11F55, 11F60, 11F99, 11G10, 11G15, 11G18, 14G35, 14K02, 14K99, 14L24, 14L30, 14L35, 14R10, 14R20, 15A15, 15A18, 15A24, 15A69, 15A72, 20C20, 20G05

\newpage

\tableofcontents

\newpage\section{Introduction}\label{S1}

\subsection{Motivation}\label{S11} 
This paper is motivated by the following problem. Let $p$ be a prime, let $N\geq 3$ be an integer prime to $p$, and let $R$ be a $\mathbb Z_{(p)}$-algebra. Suppose one wants to consider the coequalizer\footnote{We recall the notion of {\it coequalizer}, alternatively {\it categorical quotient}, 
 of a pair of morphisms $\pi_1,\pi_2:X'\rightrightarrows X$ in a category: it is a morphism $\pi:X\rightarrow Y$ of the category such that $\pi\circ \pi_1=\pi\circ \pi_2$ and it is universal with this property.} (categorical quotient),
 $\mathcal A_{g,1,N,R}/\mathcal A'_{g,1,N,R}$, in the category of schemes, of the correspondence $\mathcal A'_{g,1,N,R}\rightrightarrows \mathcal A_{g,1,N,R}$ obtained by taking the union $\mathcal A'_{g,1,N,R}$ of all prime to $p$ Hecke correspondences on Mumford's moduli scheme $\mathcal A_{g,1,N,R}$ of principally polarized abelian schemes of relative dimension $g$ endowed with symplectic similitude level-$N$ structure over $R$-schemes.\footnote{It is an open closed subscheme of the $R$-scheme constructed in \cite{MFK}, Ch. 7, Thm. 7.9 in the context of level-$N$ structures. It is Serre's lemma of \cite{Mum}, Ch. IV, Sect. 21, Thm. 5 which allows to apply \cite{MFK}, Ch. 7, Thm. 7.9 with $N\geq 3$ and not just with $N>6^g\sqrt{g!}$.}  The categorical quotient $\mathcal A_{g,1,N,R}/\mathcal A'_{g,1,N,R}$ in
algebraic geometry (exists and) is $\Spec(R)$ but this has no applications.

In our paper we  construct a categorical quotient as above in a `richer' geometry 
(category) that is  referred to as {\it $\delta$-geometry}. The ground ring $R$ of $\delta$-geometry is the ring of $p$-typical Witt vectors with coefficients in an algebraic closure $k$ of $\mathbb F_p$. Then, in $\delta$-geometry, our substitute for the quotient $\mathcal A_{g,1,N,R}/\mathcal A'_{g,1,N,R}$ turns out to be highly non-trivial and its structure is studied in detail. We  introduce the $\delta$-geometry in Subsection \ref{S21}; our version of it is related to but slightly simpler and more flexible than that in \cite{Bu05}, Sect. 2.1. Then the substitute of $\mathcal A_{g,1,N,R}/\mathcal A'_{g,1,N,R}$ appears in the $\delta$-geometry as a certain `$\Proj_{\delta}$' object attached to the graded $R$-algebra $\mathbb I_g$ of {\it Hecke covariant Siegel $\delta$-modular forms of genus $g$}. 

Siegel $\delta$-modular forms were introduced in \cite{Bu00}, Sect. 1.3 for genus $g=1$ and in \cite{BB} for all genera $g\geq 1$. The present paper can be viewed as an independent continuation of \cite{Bu05, BB}: only few facts from \cite{Bu05, BB} are needed here and we  briefly review those facts in the body of the present paper. Note that our theory here is worked out over the $\mathbb Z_p$-algebra $R$ instead of, as in \cite{BB}, over $\mathbb Z_p$; moreover, in \cite{Bu00}, \cite{BB}, \cite{Bu05} the prime $p$ was assumed to be $\geq 5$ but the methods of the present paper work for $p\in\{2,3\}$ as well. We also note that, as in \cite{BB}, 
our definition of Siegel $\delta$-modular forms does not involve 
the classical symplectic similitude level-$N$ structures. A more refined concept of Siegel $\delta$-modular form could be introduced that uses these level structures but 
we expect that the theory is independent of any level (for $g=1$ one can check that this is so based on \cite{Bu05}, Sects. 8.5 and 8.6). So we  often work with the algebraic  stacks $\mathcal A_{g,R}$ of principally polarized abelian schemes of relative dimension $g$ over $R$-schemes and with their respective unions of correspondences $\mathcal A'_{g,R}$. 

It is worth mentioning that the $\delta$-geometry referred to above and based on \cite{Bu95a, Bu05} is an arithmetic analogue of the differential algebraic geometry of Ritt and Kolchin \cite{Ko}
and of the difference algebraic geometry of Ritt and Cohn \cite{Co, Le}: the derivations in \cite{Ko} and the difference operators in \cite{Co, Le} are replaced in the $\delta$-geometry by {\it Fermat quotient operators} ({\it $p$-derivations}), see \cite{Jo, Bu95a}. Moreover, the theory of the present paper can be viewed as an arithmetic analogue of the differential algebraic theory developed in \cite{Bu95b}. In particular, for instance, the maps of sets in Equation (\ref{EQ046}) of Subsubsection \ref{S237}
defined by tuples of Hecke covariant Siegel $\delta$-modular forms  can be viewed as arithmetic analogues of the maps defined by `Schwarzian type' differential equation  considered in \cite{Bu95b}, Thm. 0.2. Our paper is logically independent from \cite{Bu95b} but some of our motivation and ideas have their origin in it. 

We note that a variant of $\delta$-geometry was used in \cite{Borg} as a candidate for the geometry over the field with one element. Also $\delta$-rings recently played a role in the foundation of prismatic cohomology (see \cite{BS}). It would be interesting to interrelate these topics and the topic of the present paper.

Here is an informal description of the philosophy behind our paper (and \cite{BB, Bu05}). Classical modular forms and $p$-adic modular forms (see \cite{Ka73a}) can be envisioned as sections of certain line bundles on the moduli stack $\mathcal A_{g,R}$; there exist some technical caveats in this description (related to non-representability issues) which are not accounted for here. On the other hand, to each scheme $X$ one attaches its $p$-jet spaces $J^r(X)$ indexed by an integer $r\geq 0$ (see \cite{Bu95a}, Sect. 1.2, where $J^r(X)$ are denoted by $X^r$). Then Siegel $\delta$-modular forms of orders $\leq r$ can be informally described as sections of certain line bundles on the `$p$-jet stacks $J^r(\mathcal A_{g,R})$'; similar caveats as above apply in this description. Weights of Siegel $\delta$-modular forms are no longer assumed to be integers but rather elements in the ring $\mathbb Z[\phi]$, where $\phi$ is an indeterminate thought of as a Frobenius lift. The $\delta$-geometric approach to the study of Hecke correspondences is based on replacing classical or $p$-adic Siegel modular forms by  Siegel $\delta$-modular forms.

 \subsection{Strategy}\label{S12}
 Among Siegel $\delta$-modular forms we  single out certain forms that have a covariance property with respect to the Hecke correspondences on $\mathcal A_{g,R}$, to be referred to as the {\it Hecke covariance}. It is a property stronger than that of being a {\it Hecke eigenform}; in fact, no classical or $p$-adic modular forms can be Hecke covariant, i.e., our theory has no classical counterpart. Besides the $\mathbb Z[\phi]$-graded $R$-algebra $\mathbb I_g$ of Hecke covariant Siegel $\delta$-modular forms of genus $g$, we will also study a larger $R$-algebra $\mathbb I_{g,\ord}$ which is the variant of $\mathbb I_g$ in which $\mathcal A_{g,R}$ is replaced by its open substack $\mathcal A_{g,R,\ord}$ whose fiber over $\Spec(k)$ is the ordinary locus of $\mathcal A_{g,k}$. 
The rings $\mathbb I_g$ and $\mathbb I_{g,\ord}$ are integral domains.
One can consider  the `$\Proj_{\delta}$' objects 
$\Proj_{\delta}(\mathbb I_g,\f)$ and $\Proj_{\delta}(\mathbb I_{g,\ord},\f)$ 
over $\Spec(R)$ attached to these $\mathbb Z[\phi]$-graded $R$-algebras and to `linear systems' $\f=(f_0,\ldots,f_n)$ with $f_0,\ldots, f_n$ of the same weight in 
either $\mathbb I_g$ or $\mathbb I_{g,\ord}$.
These `$\Proj_{\delta}$' objects
are semisimplicial objects in the dual category of rings equipped with $p$-derivations. The category of such semisimplicial objects
is  referred to as the category of {\it semi-simplicial affine $\delta$-spaces} or shortly of {\it ssa-$\delta$-spaces} and geometry in this category is what we refer to, in this paper, as {\it $\delta$-geometry}. Then 
$\Proj_{\delta}(\mathbb I_g,\f)$ and $\Proj_{\delta}(\mathbb I_{g,\ord},\f)$ are the categorical quotients, in the category of ssa-$\delta$-spaces, of suitable correspondences that depend on $\f$ and
that play, in $\delta$-geometry, the roles of the correspondences
$\mathcal A'_{g,R}\rightrightarrows \mathcal A'_{g,R}$ and their restrictions $\mathcal A'_{g,R,\ord}\rightrightarrows \mathcal A_{g,R,\ord}$ (respectively).
 This formalism is somewhat parallel to the formalism of (classical) geometric invariant theory as in \cite{MFK}; instead of group actions we consider Hecke correspondences and instead of usual varieties or schemes we consider {\it ssa-$\delta$-spaces}.

The main tool in \cite{BB} was a theory of Fourier expansions of Siegel $\delta$-modular forms (at infinity). For $g=1$ one has a Fourier expansion theory (at infinity) and also a Serre--Tate expansion theory (at ordinary $k$-valued points), see \cite{Bu00}, Sect. 7, \cite{Ba}, Sect. 5, \cite{Bu05}, Ch. 8, Sects. 8.1 and 8.3. One of the somewhat surprising outcomes was that for Hecke covariant forms the two types of expansions coincide up to a natural identification. 

The main tool  developed  in this paper is a Serre--Tate expansion theory for Siegel $\delta$-modular forms of genus $g\geq 1$. The
 main new idea is  
 to construct a series of `comparison maps' between 
 the $R$-algebras $\mathbb I_g$ and $\mathbb I_{g,\ord}$ and certain $R$-algebras that appear in the geometric invariant theory of multiple quadratic forms and multiple endomorphisms over $R$; we  then transfer information from the latter $R$-algebras to the former $R$-algebras. We need some new results on the algebras of geometric invariant theory which are proved in Section 4. We emphasize that invariant theory shows up in our paper in two distinct ways:

\medskip\noindent
{\bf 1.} First, the $R$-algebras $\mathbb I_g$ and $\mathbb I_{g,\ord}$ are (in the context of Hecke correspondences on $\mathcal A_{g,R}$) two `$\delta$-geometric analogues' of the classical algebras of invariants of linear algebraic group actions.

\smallskip\noindent
{\bf 2.} Second, we use concrete instances of classical geometric invariant theory (in the context of the congruent $\pmb{\SL}_g$-actions on vector spaces of multiple quadratic forms (i.e., of symmetric matrices) and of conjugate $\pmb{\SL}_g$-actions on vector spaces of multiple endomorphisms over either $k$ or an algebraic closure of the field of fractions $K:=\Frac(R)$ of $R$, see Section \ref{S4}) to study the $R$-algebras $\mathbb I_g$ and $\mathbb I_{g,\ord}$.

\bigskip

\subsection{Results} 

We now outline our main results that are described in Subsection \ref{S24}. 

First, we prove a Serre--Tate expansion principle similar to the classical Fourier expansion principle (see Theorem \ref{T1}). Generalizing the known result for $g=1$, we prove that the Serre--Tate expansions of the {\it basic} Hecke covariant forms $f^r_{g,\crys}$ of genus $g$, size $g$, order $r$ and weight $(\phi^r,1)$ (introduced in \cite{BB}, Subsect. 4.1 and reviewed in Subsection \ref{S34} using the language recalled in Subsection \ref{S23} and) coincide, up to a natural identification, again, with their Fourier expansions (see Theorem \ref{T2}). For special generating properties of $f^1_{g,\crys}$ see 
Theorem \ref{T3} (a).
Furthermore, generalizing a construction introduced in \cite{Ba} in case $g=1$, we  define (see Subsection \ref{S34}) for all $g\geq 1$ a remarkable  ordinary Siegel $\delta$-modular form $f^{\partial}_{g,\crys}$ of genus $g$, size $g$, and order $1$ (also referred to as {\it basic}) which is invertible and whose Serre--Tate expansion equals the identity $g\times g$ matrix (see Theorem \ref{T4} (a)). 

Second, using the Serre--Tate expansion theory and the basic forms  
we  show that the $R$-algebra $\mathbb I_{g,\ord}$ coincides, up to torsion, with a ring of Laurent polynomials in a countable number of indeterminates over the $R$-algebra of invariants of multiple quadratic forms (see Theorem \ref{T5}). The analogue of the latter $R$-algebra over an algebraically closed field is a classical object and is revisited and further studied in detail in Subsections \ref{S41} to \ref{S43}. As a consequence of this link with invariant theory, if we denote by $\mathbb I^r_{g,\ord}$ the $R$-subalgebra of elements of $\mathbb I_{g,\ord}$ of order $\leq r$, we  obtain an exact formula for the transcendence degree of the fraction field $\Frac(\mathbb I^r_{g,\ord})$ over $K$ (see Remark \ref{R14}). This transcendence degree divided by $r$ tends to a limit, referred to as the {\it asymptotic transcendence degree} of $\mathbb I_{g,\ord}$ and denoted by
 $\asytrdeg(\mathbb I_{g,\ord})$, which we prove to be equal to $\dim(\mathcal A_{g,K})+1=\dim(\mathcal A_{g,1,N,K})+1=\frac{g(g+1)}{2}+1$.
We  also prove that  $\Frac(\mathbb I_{g,\ord})$ is `$\phi$-generated' by finitely many elements in $\mathbb I^4_{g,\ord}$.
Finally, we  prove that the `$\phi$-transcendence degree' of  
$\Frac(\mathbb I_{g,\ord})$ over $K$, i.e., the maximum number of $\phi$-algebraically independent elements of $\Frac(\mathbb I_{g,\ord})$, equals $\frac{g(g+1)}{2}+1$ (see Theorem \ref{T6} for the last three results). 

Third, defining $\mathbb I_g^r$ similarly (so $\mathbb I_g^r=\mathbb I_g\cap \mathbb I_{g,\ord}^r$) and using the invariant theory of multiple endomorphisms acted upon via conjugation, we  show that the sequence obtained by dividing the transcendence degree of $\Frac(\mathbb I^r_g)$ over $K$ by $r$ tends to a limit $\asytrdeg(\mathbb I_g)$ which we prove to be equal, again, to $\frac{g(g+1)}{2}+1$ (see Theorem \ref{T8}). 
Note that, by a general result in difference algebra (see \cite{Le}, Thm. 4.4.1, p. 292), the field
 $\Frac(\mathbb I_g)$ is also $\phi$-finitely generated. We prove that `$\phi$-transcendence degree' of $\Frac(\mathbb I_g)$ over $K$ equals $\frac{g(g+1)}{2}+1$ (see Theorem \ref{T8}). 
The latter result was first obtained for $g\in\{1,2\}$ in \cite{BB}, Thms. 1.13 and 1.14, and for $g\geq 3$ solves the `main open problem' formulated in the paragraph after \cite{BB}, Thm. 1.13. We also prove that the transcendence degree of $\Frac(\mathbb I_{g,\ord})$ over $\Frac(\mathbb I_g)$ is $\leq \frac{g(g+1)}{2}+2$ (see Remark \ref{R15}). 
 We can view the latter as saying that the $R$-subalgebra $\mathbb I_g$ of $\mathbb I_{g,\ord}$ is a `transcendentally close' to the algebra $\mathbb I_{g,\ord}$.
As an application of our structure theory for $\mathbb I_g$ we  show that various Zariski dense loci in $\mathcal A_{g,1,N,R}(R)$ are contained in zero loci of (many) Hecke covariant Siegel $\delta$-modular forms that are proper subsets of $\mathcal A_{g,1,N,R}(R)$. This is the case with the locus (set): 

\medskip
{\bf (i)} of  abelian schemes with complex multiplication (see Remark \ref{R12});

\smallskip
{\bf (ii)} of a  polarized isogeny prime to $p$ class (see Remark \ref{R11}); 

\smallskip
{\bf (iii)} of  principally polarized abelian schemes with non-commutative `prime to $p$' ring of endomorphisms (see Theorem \ref{T9});

\smallskip
{\bf (iv)} of  principally polarized abelian schemes whose polarized isogeny prime to $p$ class is decomposable (see Theorem \ref{T10}).

\smallskip

\medskip\noindent We also prove a  `$p$-adic approximation' result saying, informally speaking, that each element of $\mathbb I_{g,\ord}$ is a $p$-adic limit of fractions of elements in $\mathbb I_g$ (see Theorem \ref{T11}). This means that the $R$-subalgebra $\mathbb I_g$ of $\mathbb I_{g,\ord}$ is also `big' in the sense of being a `close $p$-adic approximation' of $\mathbb I_{g,\ord}$. 

Fourth, we  prove a {\it $\delta$-flow} structure theorem for the arithmetic differential equation defined by the basic form $f^1_{g,\crys}$ (see Theorem \ref{T12}). The concept of $\delta$-flow appearing here is an arithmetic analogue of the classical concept of flow on a manifold and played a role, for instance, in \cite{BM}, Def. 5.3, p. 117 and \cite{BP}, Lem. 4.45. Using our $\delta$-flow structure theorem plus the arithmetic analogue (see \cite{Bu97}, Cor. 1.7) of Manin's Theorem of the Kernel (see \cite{Man}, Sect. 5, Thm. 2), we  prove (see Theorem \ref{T12.12}) an extension to arbitrary genus $g\geq 1$ of the `Reciprocity Theorem for $\CL$ (canonical lifts) points' in \cite{BP}, Thm. 3.5; we do not obtain, however, the conclusion, in the latter theorem, that our `reciprocity map' is non-constant. Theorem \ref{T12.12} is a purely diophantine result about 
 $\mathbb Z$-linear dependence 
in the group of $K$-rational points of an abelian variety over $K$ of points arising from $\CM$ (complex multiplication) points on a Siegel modular variety.

As an addition, we would like to mention that Propositions \ref{P8} to \ref{P12} of Subsection \ref{S310.66} compute cohomology groups of our ssa-$\delta$-spaces  $\Proj_{\delta}(\mathbb I_g,\f)$ and $\Proj_{\delta}(\mathbb I_{g,\ord},\f)$   in some simple cases.

\subsection{Open problems}  Our results suggest a number of problems and directions of investigations, including the following ones:

\smallskip

{\bf (a)}  Find generators and relations for the rings $\mathbb I_{g,\ord}$ and $\mathbb I_g$.

\smallskip

{\bf (b)} Compute the exact transcendence degree of $\Frac(\mathbb I_{g,\ord})$ over $\Frac(\mathbb I_g)$.

\smallskip

{\bf (c)} Compute the exact ranks of the $R$-modules of forms in  $\mathbb I_{g,\ord}$ and $\mathbb I_g$, of given order and weight. 

\smallskip

{\bf (d)} Compute the ideals of all  forms in the rings
$\mathbb I_{g,\ord}$ and $\mathbb I_g$
that vanish on the various loci
(i) to (iv) of $\mathcal A_{g,1,N,R}(R)$. 

\smallskip

{\bf (e)} Compute the cohomology groups of all $\Proj_{\delta}(\mathbb I_g,\f)$ and  $\Proj_{\delta}(\mathbb I_{g,\ord},\f)$.

\smallskip

{\bf (f)} Develop an analogue of our theory for more general Shimura varieties.
 (The case of Shimura curves was developed in \cite{Bu03}, \cite{Bu04}, \cite{Bu05}.)
 
 \smallskip
 
 {\bf (g)} Recalling that in the case of $\mathcal A_1$ and of Shimura curves (see \cite{Ba}, \cite{Bu04}, \cite{Bu05}) an important role was played by certain  operators, called {\it $\delta$-Serre operators}, that are analogues  of the Serre differential operators acting on classical or $p$-adic modular forms, develop a theory  of  $\delta$-Serre operators for stacks of higher dimensional Shimura varieties, in particular for $\mathcal A_g$.
 
\smallskip

For our partial results and more details and comments on questions (a) to (e) we refer to 
Theorems \ref{T3},  \ref{T4},  \ref{T7} and  \ref{T8} and Remarks \ref{R7.4}, \ref{R9}, \ref{R11}, \ref{R13}, \ref{R14}, \ref{R15}, \ref{R16} and \ref{R17} of Subsection \ref{S24} .

\bigskip

\subsection{Plan of the paper} Subsection \ref{S21} reviews some notation, terminology, and concepts including $p$-jet spaces introduced in \cite{Bu95a}, p. 315. Subsection \ref{S22} introduces some dimension concepts that are useful in the paper. The basic concepts on Siegel $\delta$-forms \cite{BB} are recalled and expanded in Subsection \ref{S23}. The main results are stated in Subsection \ref{S24} and are proved in Section \ref{S3}. For a detailed description of the contents of Section 3 and, in particular, for our detailed strategy of proofs, we refer to the beginning of Section 3. Section \ref{S4} is devoted to revisiting the classical  invariant theory of multiple quadratic forms and multiple endomorphisms and proving new results that are being often used and cited in Sections \ref{S1} to \ref{S3}: it is organized as an appendix, i.e., it is completely a self-contained section (including the notation) which can be read out at any time. For  more details on Section 4 we refer to its beginning. 

\newpage\section{Main concepts and results}\label{S2}

\subsection{Concepts of $\delta$-geometry}\label{S21} 

We begin by reviewing some ring theoretic terminology and notation.
Unless otherwise stated, all rings will always be assumed commutative with $1$.

By a {\it filtered} ring we mean a ring $B$ equipped with an increasing exhausting filtration $(B^r)_{r\geq 0}$
of subrings; hence $B^r\subset B^{r+1}$ for all integers $r\geq 0$ and 
$B=\bigcup _{r=0}^{\infty} B^r$. 
Filtered rings are the objects of a category whose morphisms are ring homomorphisms $u:B\rightarrow \tilde{B}$ that preserve the filtrations, i.e., we have $u(B^r)\subset \tilde{B}^r$ for all integers $r\geq 0$.

By a {\it graded} ring we mean a ring $B$ equipped with a direct sum decomposition $B=\oplus_{w\in W}B(w)$ where $W$ is an abelian semigroup and $B(w)$ are subgroups such that we have $B(w)B(v)\subset B(w+v)$ for all $w,v\in W$. The elements of $B(w)$ will be referred to as having {\it degree} $w$ or {\it weight} $w$. We also say that $B$ is a $W$-graded ring. The $W$-graded rings are the objects of a category whose morphisms are ring homomorphisms $u:B\rightarrow \tilde{B}$ that preserve the gradings, i.e., we have $u(B(w))\subset \tilde{B}(w)$ for all $w\in W$.

\begin{df}\label{df1}
A {\it filtered graded ring} is a 
ring $B$ equipped with a filtered ring structure, given 
by a filtration $(B^r)_{r\geq 0}$, and also equipped with a graded ring structure, given by $B=\oplus_{w\in W} B(w)$ with $W$ as an abelian semigroup, the two structures being compatible in the sense that for all integers $r\geq 0$ we have a direct sum decomposition $B^r=\oplus_{w\in W} B^r\cap [B(w)]$. We also call $B$ a {\it filtered $W$-graded ring}
and we set $B^r(w):= B^r\cap [B(w)]$.
\end{df}

Filtered $W$-graded rings are the objects of a category whose morphisms are ring homomorphisms that preserve the filtrations and the gradings.

Assume one is given a sequence of $W$-graded rings $B^r=\oplus_{w\in W}B^r(w)$, $r\geq 0$,
such that $B^r$ is a subring of $B^{r+1}$ and $B^r(w)\subset B^{r+1}(w)$ for all integers  $r\geq 0$ and $w\in W$; let us refer to this data as a {\it filtered $W$-graded system} or just as a {\it filtered graded system}. To each filtered graded system as above we can attach a filtered graded ring as follows: define $B:=\bigcup_{r\geq 0} B^r$ and $B(w):=\bigcup_{r\geq 0}B^r(w)$. Then
it is easy to check that 
$B=\oplus_{w\in W}B(w)$ and
$B^r(w)=B^r\cap B(w)$ for all integers  $r\geq 0$ and $w\in W$, so
 $B$ with the above grading and with filtration $(B^r)_{r\geq 0}$ is a filtered graded ring. All filtered graded rings can be obtained in this way.

Throughout the paper $p$ will denote a field characteristic: so it is either $0$ or a prime natural number. Until Section \ref{S4} we will assume that $p>0$. 

If $S$ is a ring and $L$ is an $S$-module, let $\widehat{L}$ be its $p$-adic completion. We recall that $L$ is called $p$-adically complete
if the functorial $S$-linear map $L\rightarrow \widehat{L}$ is an isomorphism. We write $\overline{S}:=S/(p)=S/pS$ and we denote by $S^{\times}$ the group of units of $S$. If $S$ is local, we denote by $\mathfrak m_S$ its maximal ideal; we have $S^{\times}=S\setminus\mathfrak m_S$. 

If $X$ is a scheme, we denote by $\mathcal O_X$ its structure ringed sheaf and we use the shorter notation $\overline{X}:=X\times_{\Spec(\mathbb Z)}\Spec(\mathbb F_p)$
and $\mathcal O(X):=H^0(X,\mathcal O_X)$.
If $X=\Spec(S)$ is affine, then $\overline{X}=\Spec(\overline{S})$ and for a  scheme $Z$ over $X$ we let $Z(S)$ be the set of all $X$-valued points of $Z$. 

Similarly, if $Y$ is a $p$-adic formal scheme (i.e., a formal scheme with ideal of definition generated by $p$), we denote by $\mathcal O_Y$ its structure ringed sheaf and we use the shorter notation $\overline{Y}:=Y\times_{\Spf(\mathbb Z_p)}\Spf(\mathbb F_p)$
and $\mathcal O(Y):=H^0(Y,\mathcal O_Y)$.
If $Y=\Spf(S)$ is affine, then $\overline{Y}=\Spf(\overline{S})=\Spec(\overline{S})$ and for a $p$-adic formal scheme $Z$ over $Y$ we let $Z(S)$ be the set of all $Y$-valued points of $Z$. 

For a noetherian scheme $X$ we denote by $\widehat{X}$ the $p$-adic formal scheme obtained by $p$-adically  completing $X$. If $X$ is affine, then $\overline{X}=\overline{\widehat{X}}$. 
 
Let $\pmb{\Mat}_g$ be the ring $\mathbb Z$-scheme of $g\times g$ matrices. Let $1_g\in\pmb{\Mat}_g(S)$ be the identity matrix. Let $\mathbb E_g$ be the matrix with all entries equal to $1$. For each matrix
$M\in\pmb{\Mat}_g(S)$ let
 $M^{\t}$ and $M^*$ denote the transpose and the adjugate of $M$ (respectively); so $MM^*=\det(M)\cdot 1_g$.
 
\begin{df}\label{df2}
Let $\rho:S_1\rightarrow S_2$ be a ring homomorphism between two $\mathbb Z_{(p)}$-algebras and let $\bar\rho:S_1/pS_1\rightarrow S_2/pS_2$ be its reduction modulo $p$. By a {\it Frobenius lift relative to $\rho$} we mean a homomorphism $\phi: S_1\rightarrow S_2$ whose reduction modulo $p$ maps $x\in S_1/pS_1$ to $\bar\rho(x)^p\in S_2/pS_2$.
By a {\it Frobenius lift} on a $\mathbb Z_{(p)}$-algebra $S$ we understand a Frobenius lift $\phi:S \rightarrow S$ relative to the identity $1_S$. \end{df}

Throughout the paper we let $k$, $R$, and $K$ be as in the Introduction; so:

\medskip

$\bullet$ $k$ is an algebraic closure of $\mathbb F_p$; 

\smallskip

$\bullet$ $R=W(k)$ is the ring of $p$-typical Witt vectors with coefficients in $k$; 

\smallskip

$\bullet$ $K=R[\frac{1}{p}]=\Frac(R)$. 

\medskip

Recall that there exists a unique Frobenius lift $\phi$ on $R$.

We next review the concepts of $p$-derivations \cite{Jo, Bu95a} and of $p$-jet spaces
\cite{Bu95a}. Consider the polynomial
$$C_p(z_1,z_2):=\frac{z_1^p+z_2^p-(z_1+z_2)^p}{p} \in\mathbb Z[z_1,z_2].$$

\begin{df}\label{df3}
If $\iota:S\rightarrow B$ is a ring homomorphism (in practice it is injective) between two $\mathbb Z_{(p)}$-algebras, then by a {\it $p$-derivation on $B$ relative to $\iota$} (or simply a {\it $p$-derivation relative to $\iota$}) we mean a function $\delta:S \rightarrow B$ such that the following two axioms hold:

\medskip

{\bf (DAX1)} For all $x,y\in S$ we have $\delta(x+y) = \delta x + \delta y +\iota(C_p(x,y))$.

\smallskip

{\bf (DAX2)} For all $x,y\in S$ we have $\delta(xy) = \iota(x)^p \delta y + \iota(y)^p \delta x +p \delta x \delta y$.
\end{df}

If $\delta$ is a $p$-derivation relative to $\iota$, then the homomorphism $\phi:S \rightarrow B$ defined by the rule 
$$\phi(x):=\iota(x)^p+p \delta(x)$$ 
(with $x\in S$) is a Frobenius lift relative to $\iota$, to be called the {\it Frobenius lift attached to $\delta$}. Conversely, given a Frobenius lift $\phi:S\rightarrow B$ relative to $\iota:S\rightarrow B$, if $p$ is a non-zero divisor in $B$, then $\delta:S\rightarrow B$, defined by the rule $\delta(x):=\frac{\phi(x)-\iota(x)^p}{p}$, is a $p$-derivation relative to $\iota$. 

By a {\it $p$-derivation} on a $\mathbb Z_{(p)}$-algebra $S$ we understand a $p$-derivation $\phi:S \rightarrow S$ relative to the identity $1_S$. Hence the $R$-algebra $R$ has a unique $p$-derivation $\delta:R\rightarrow R$.

\begin{df}\label{df4}

{\bf (a)}
By a $\delta$-{\it ring} we mean a $\mathbb Z_{(p)}$-algebra $S$ equipped with a $p$-derivation $\delta:S\rightarrow S$. 

{\bf (b)} By a {\it filtered $\delta$-ring} we mean a $\delta$-ring $S$ 
which is also a filtered ring with filtration $(S^r)_{r\geq 0}$ 
such that we have $\delta(S^r)\subset S^{r+1}$ for all integers $r\geq 0$.\end{df}

Note that a filtered $\delta$-ring in the above sense 
is not a filtered object in the category of $\delta$-rings; rather it is an object of the category of $\delta$-rings and also an object of the category of filtered rings and the two structures are subject to an appropriate compatibility condition.

A morphism of $\delta$-rings is a ring homomorphism that commutes with the corresponding $p$-derivations. 
The class of $\delta$-rings and their morphisms form a category $\textup{\bf Ring}_{\delta}$ which has fiber coproducts and equalizers. If $\textup{\bf Ring}$ is the category of rings over $\mathbb Z_{(p)}$, then
the forgetful functor $\textup{\bf Ring}_{\delta}\rightarrow\textup{\bf Ring}$, which forgets the $p$-derivations, commutes with coproducts and equalizers.

A morphism of filtered $\delta$-rings $u:S\rightarrow \tilde{S}$ is a morphism of $\delta$-rings 
 such that $u(S^r)\subset \tilde{S}^r$ for all integers  $r\geq 0$.
The class of filtered $\delta$-rings and their morphisms form a category $\textup{\bf FilRing}_{\delta}$ that has fiber coproducts and equalizers.

The dual category $\textup{\bf Ring}_{\delta}^{\textup{op}}$ will be denoted by $\textup{\bf Aff}_{\delta}$ and its objects will be referred to as {\it affine $\delta$-spaces}.
The object in  $\textup{\bf Aff}_{\delta}$ that corresponds to a $\delta$-ring $B\in \Ob(\textup{\bf Ring}_{\delta})$ will be denoted by $\Spec_{\delta}(B)\in \Ob(\textup{\bf Aff}_{\delta})$. We often drop mentioning the $p$-derivations $\delta$ when referring to objects or morphisms of either $\textup{\bf Ring}_{\delta}$ or $\textup{\bf Aff}_{\delta}$.

\begin{rem}\label{R0}
Note that in the definition above the word `affine' is not an adjective as we did not introduce a  concept of {\it $\delta$-space}. One can introduce such a concept as follows, 
cf. \cite{Borg}, p. 7.
We equip the category $\textup{\bf Aff}:=\textup{\bf Ring}^{\textup{op}}$ with  the \'{e}tale topology and 
we denote by  $\textup{\bf Sp}$  the category of sheaves of sets   on  $\textup{\bf Aff}$.
The infinite-length $p$-typical Witt vector functor 
$W^*:\textup{\bf Sp}\rightarrow \textup{\bf Sp}$  possesses a natural monad structure.
 Then 
one  defines a {\it $\delta$-space} to be 
an object of $\textup{\bf Sp}$   equipped  with a $W^*$-action.
 (For the details of this definition we refer to \cite{Borg}; note however that in loc. cit.  one works over $\mathbb Z$ instead of over $\mathbb Z_{(p)}$, one uses  the big Witt vectors instead of the $p$-typical Witt vectors, and one uses the terminology of `$\Lambda$-spaces' instead of `$\delta$-spaces'.)   With morphisms of $\delta$-spaces defined in the obvious way we obtain the {\it category $\textup{\bf Sp}_{\delta}$ of $\delta$-spaces}.
Then  $\textup{\bf Aff}_{\delta}$ can be  identified with the full subcategory 
  of $\textup{\bf Sp}_{\delta}$ whose objects are the objects of
  $\textup{\bf Sp}_{\delta}$ with the property that their
 image in $\textup{\bf Sp}$ is represented by an object of $\textup{\bf Aff}$. 
The category $\textup{\bf Sp}_{\delta}$ can be chosen as  a natural  framework for $\delta$-geometry as in \cite{Borg}. 
However, in the present paper, we found it is more helpful and expedient to present
  $\delta$-geometry  
 in a semi-simplicial (rather than sheaf theoretic)  framework, see Definition \ref{df8} below. This framework is a  variant of the one in \cite{Bu05}, Ch. 2. 
 In particular, the category $\textup{\bf Sp}_{\delta}$ will play no role in the present paper.
\end{rem}

For $d\in\mathbb N$, a basic example of filtered $\delta$-ring is given as follows. Let
 $$S^r=R[y,y',\ldots,y^{(r)}]:=R[y_i^{(l)}|1\leq i\leq d, 0\leq l\leq r],\ \ r\geq 0,$$
 where 
 $$\begin{array}{l}
 y=y^{(0)}=(y_1,\ldots,y_d)=(y_1^{(0)},\ldots,y_d^{(0)}),\\
 y'=y^{(1)}=(y_1',\ldots,y_d')=(y_1^{(1)},
 \ldots,
 y_d^{(1)}),\\
\ldots,\;  y^{(r)}=(y_1^{(r)},\ldots,y_d^{(r)}),\;
 \ldots
 \end{array}$$ 
 are $d$-tuples of indeterminates. We have inclusions
 $S^r\rightarrow S^{r+1}$ 
 and, for each $r$, we let $\delta:S^r\rightarrow S^{r+1}$
 be the unique $p$-derivation relative to the inclusion $S^r\subset S^{r+1}$ such that we have $\delta y=y',\ldots,\delta y^{(r)}=y^{(r+1)}$ (i.e., for integers $1\leq i\leq d$ and $0\leq l\leq r$ we have $\delta y_i^{(l)})=y_i^{(l+1)}$). We define the $R$-algebra
 of {\it $\delta$-polynomials},
\begin{equation}\label{EQ001}
 R[y,y',\ldots,y^{(r)},\ldots]:=\bigcup_{r=0}^{\infty} S^r;
\end{equation} this algebra has an induced structure of a filtered $\delta$-ring with filtration $(S^r)_{r\geq 0}$.

\begin{df}\label{df5}
Let $\Prol$ be the category whose objects are the filtered $\delta$-rings $S$
 with
each $S^r$ a $p$-adically complete, 
flat $R$-algebra and with the inclusions $S^r\subset S^{r+1}$ assumed to be $R$-algebra homomorphisms; the  morphisms in the category are the morphisms of filtered $\delta$-rings  that induce $R$-algebra homomorphisms on the rings in the filtrations. 
\end{df}

For $S\in \textup{Ob}(\Prol)$, $S$ is always flat over $R$ but in general it is not $p$-adically complete (not even $p$-adically separated); more precisely, $S$ is $p$-adically complete if and only if there exists a pair $(r_0,n_0)\in\mathbb Z_{\geq 0}^2$ such that $p^{n_0}S\subset S^{r_0}$ as one can easily check based on Baire's category theorem.\footnote{The following example due to Ofer Gabber shows that in general we cannot assume take $n_0=0$ for $r_0>>0$. If for an integer $r\geq 0$, $S^r$ is the inverse image in the ring $R\{x,y\}$ of restricted formal series in two indeterminates $x$ and $y$ of $k[x,xy,\ldots,xy^r]$, then for $S=\cup_{r\geq 0} S^r$ we can take $(r_0,n_0)=(0,1)$ but $n_0$ cannot be $0$ no matter what $r_0$ is.} In particular, if $S/pS$ is a finitely generated $k$-algebra or if $S^{r+1}/S^r$ is torsion free for each integer $r\geq 0$, then $S$ is $p$-adically complete if and only if there exists $r_0\in\mathbb Z_{\geq 0}$ such that $S=S^{r_0}$.
 The category $\Prol$ has an initial object $R$ with filtration given by $R^r=R$ for all integers  $r\geq 0$. 

If $X=\Spec(B)$ is an affine smooth scheme over $\Spec(R)$,
then we recall from \cite{Bu95a}, Sect. 1, 
 that there exists a filtered $\delta$-ring $J(B) \in\Ob(\Prol)$, with filtration $(J^r(B))_{r\geq 0}$ starting with $J^0(B)=\widehat{B}$, that satisfies the following universal property: for each $S\in\Ob(\Prol)$ and for every
$R$-algebra homomorphism $B \rightarrow S^0$ there exists a unique morphism of filtered $\delta$-rings $J(B) \rightarrow
S$ over $R$
that extends $B \rightarrow S^0$. 
The affine $p$-adic formal schemes $J^r(X):=\Spf(J^r(B))$ are $p$-adic completions of affine smooth schemes over $\Spec(R)$.
This follows from the fact that if $y=(y_1,\ldots,y_d) \subset B^d=\mathcal O(X)^d$ define \'{e}tale coordinates on $X$ (i.e, the $R$-algebra homomorphism  from the ring of polynomials $R[x_1,\ldots,x_n]$ into $B$ that maps each indeterminate $x_i$ to $y_i$ is \'etale), then we have a natural $B$-algebra isomorphism
\begin{equation}\label{EQ002}
\reallywidehat{B[y_1',\ldots,y_d',\ldots, y_1^{(r)},\ldots,y_d^{(r)}]}\simeq J^r(B)=\mathcal O(J^r(X))
\end{equation}
defined by $y_i^{(l)}\to\delta^l y_i$ for all integers $1\leq i\leq d$ and $1\leq l\leq r$ (see \cite{Bu95a}, Prop. 1.4). As explained in \cite{Bu95a}, Sect. 1, the construction can be globalized to attach formal schemes $J^r(X)$ to each separated (not necessarily affine) smooth scheme $X$ over $\Spec(R)$:
if $X=\cup_{i\in I} U_i$ is an affine open cover, then $J^r(X)=\cup J^r(U_i)$ is an open cover with $J^r(U_i)$ glued naturally so that for all $i,j\in I$ we have an identity $J^r(U_i)\cap J_r(U_j)=J^r(U_i\cap U_j)$.

\begin{df}\label{df6}
The $p$-adic formal scheme $J^r(X)$ is called the $r$-th $p$-{\it jet space} of the separated smooth scheme $X$ over $\Spec(R)$.\end{df}

For a separated 
smooth scheme $X$ over $\Spec(R)$ we may consider the filtered $\delta$-ring
$\mathcal O^{\infty}(X):=\bigcup_{r=0}^{\infty}\mathcal O(J^r(X))$ with filtration given by the rings in the union.

By the universal property mentioned above applied with $S=R$, for every $R$-valued point $P\in X(R)$, $P:\Spec(R)\rightarrow X$, there exists a unique sequence $(J^r(P))_{r\geq 0}$ of $R$-valued points $J^r(P)\in J^r(X)(R)$, $J^r(P):\Spf(R)\rightarrow J^r(X)$, 
with $J^r(P)$ lifting $J^{r-1}(P)$ for $r\geq 1$,
such 
that $J^0(P)$ is the $p$-adic completion of $P$ and such that for every affine open subscheme $U\subset X$ with $P\in U(R)$ the induced sequence of $R$-algebra homomorphisms $\mathcal O(J^r(U))\rightarrow R$ defines a morphism of  $\delta$-rings $\mathcal O^{\infty}(U)\rightarrow R$. The $R$-valued points $J^r(P)$ are covariantly functorial in $X$. 
 For every element $\varphi\in \mathcal O(J^r(X))$ 
 we consider the induced morphism (still denoted by) $\varphi:J^r(X)\rightarrow \widehat{\mathbb A^1_R}$
 and we continue to denote by $\varphi:J^r(X)(R)\rightarrow \widehat{\mathbb A^1_R}(R)=R$
 the induced map between the sets of $\Spf(R)$-valued points. Then 
 we have a well-defined map
$$
\varphi_R:X(R)\rightarrow R,\ \ \varphi_R(P):=\varphi(J^r(P)).
 $$
 The rule $\varphi\mapsto \varphi_R$ defines an $R$-algebra homomorphism
 \begin{equation}\label{EQ003}
 \mathcal O(J^r(X))\rightarrow\Maps(X(R), R)\ \ \ 
 \end{equation}
 which is injective (see \cite{Bu05}, Prop. 3.19) and thus viewed as an inclusion that allows us to write $\varphi$ in place of $\varphi_R$. Here $\Maps(X(R), R)$ denotes the $R$-algebra of set theoretic maps from $X(R)$ to $R$.

 \begin{df}\label{df7}
 A map of sets $X(R)\rightarrow R$ is called a {\it $\delta$-function} of order $\leq r$ on $X(R)$ if it belongs to the image of the $R$-algebra homomorphism (\ref{EQ003}).
 \end{df}

We say that a $\delta$-function of order $\leq r$ has order $r$ if either $r=0$ or $r\geq 1$ and it does not have order $\leq r-1$. The $\delta$-functions of order $0$ are exactly the functions induced by formal functions on $\widehat{X}=J^0(X)$. For $r\geq 1$
 the $\delta$-functions of order $r$ are not induced by formal functions on $\widehat{X}$; rather, they are given locally in the Zariski topology by restricted power series in the affine coordinates and their $p$-derivatives up to oder $r$.

One can ask for a geometry based on $\delta$-rings that could be referred to as {\it $\delta$-geometry}. For our practical purposes we would like a version of $\delta$-geometry that requires a minimum amount of foundational things. Such a formalism was introduced and used in \cite{Bu05}. We shall employ here a slightly different (although related) formalism which is, arguably, more flexible and at the same time simpler than the one in \cite{Bu05}. Here are the relevant concepts. 

First we recall (see \cite{W}, Sect. 8.1) that a semi-simplicial object in a category is
a collection of objects $X_{\bullet}=(X_d)_{d\geq 0}$  and morphisms, called {\it face maps},
$f_{d,0},\ldots,f_{d,d+1}:X_{d+1}\rightarrow X_d$ that satisfy the {\it semi-simplicial face relations}
$$f_{d,i}\circ f_{d+1,s}=f_{d,s-1}\circ f_{d+1,i} \;\;\; \forall\;\; \textup{integers satisfying}\;\; 0\leq i<s\leq d+2\geq 2.$$
Consider now the following example. Given a separated scheme $X$ and an affine open cover $X=U_0\cup \cdots \cup U_n$ one can define a semi-simplicial object $X_{\bullet}=(X_d)_{d\geq 0}$ in the category of affine schemes
with 
$$X_d:=\coprod_{0\leq i_0<\cdots <i_d\leq n} U_{i_0}\cap \cdots \cap U_{i_d}$$
and with the obvious face maps. Then $X_{\bullet}$ can be viewed as a `substitute' for $X$. We would like to take this simple example as a blueprint for $\delta$-geometry. 

\begin{df}\label{df8}
 A {\it semi-simplicial affine $\delta$-space} (or simply an {\it ssa-$\delta$-space})  is a semi-simplicial object in the category ${\bf Aff}_{\delta}={\bf Ring}_{\delta}^{\textup{op}}$.\end{df}
 
 In concrete terms, an ssa-$\delta$-space is a sequence $\Spec_{\delta}(B_{\bullet})=(\Spec_{\delta}(B_d))_{d\geq 0}$ where $B_d$ are $\delta$-rings equipped with $\delta$-ring homomorphisms $u_{d,0},\ldots,u_{d,d+1}:B_d\rightarrow B_{d+1}$ that satisfy the {\it coface} relations
$$u_{d+1,s}\circ u_{d,i}=u_{d+1,i}\circ u_{d,s-1}\;\;\; \forall\;\; \textup{integers that satisfy}\;\; 0\leq i<s\leq d+2.$$
Morphisms of ssa-$\delta$-spaces are defined as morphisms of semi-simplicial objects in ${\bf Aff}_{\delta}$. In this way we get the {\it category of semi-simplicial affine $\delta$-spaces} denoted by 
$\textup{\bf ssAff}_{\delta}$. 
This category has fiber products and coequalizers: the coequalizer 
of a pair of morphisms $\Spec_{\delta}(B'_{\bullet})
\rightrightarrows \Spec_{\delta}(B_{\bullet})$ in $\textup{\bf ssAff}_{\delta}$
is the ssa-$\delta$-space $\Spec_{\delta}(C_{\bullet})$
where for each $d\geq 0$, $C_d$ is the equalizer in $\textup{\bf Ring}_{\delta}$ of the corresponding pair of morphisms $B_d \rightrightarrows B'_d$.

Each $\delta$-ring $S$ defines an ssa-$\delta$-space 
$\Spec_{\delta}(S_{\bullet})=(\Spec_{\delta}(S_d))_{d\geq 0}$ where $S_d=S$ for all integers  $d\geq 0$ and where all face maps are the identity. The functor $\textup{\bf Aff}_{\delta}\rightarrow \textup{\bf ssAff}_{\delta}$ defined by $\Spec_{\delta}(S)\mapsto \Spec_{\delta}(S_{\bullet})$ identifies $\textup{\bf Aff}_{\delta}$ with a full subcategory of $\textup{\bf ssAff}_{\delta}$.

We will usually denote ssa-$\delta$-spaces by symbols such as $X_{\delta},Y_{\delta}$, etc.

An ssa-$\delta$-space $X_{\delta}$ is said to be {\it over $S$}, if it is equipped with a morphism $X_{\delta}\rightarrow \Spec_{\delta}(S_{\bullet})$. We will be especially interested in the category of ssa-$\delta$-spaces over $R$. Then $\Spec_{\delta}(R_{\bullet})$ is a final object
in the category of ssa-$\delta$-spaces over $R$.

The $d$-th cohomology group of an ssa-$\delta$-space $X_{\delta}=\Spec_{\delta}(B_{\bullet})$ is defined as usual:
$$H^d(X_{\delta}):=\frac{\Ker(\partial_{d})}{\textup{Im}(\partial_{d-1})},\ \ \ \partial_d:=\sum_{i=0}^{d+1} (-1)^i u_{d,i}:B_d\rightarrow B_{d+1},\ \ \ \ \partial_{-1}:=0.$$
The group $H^0(X_{\delta})$ has a structure of $\delta$-ring which we call the {\it ring of global functions} on $X_{\delta}$. If $X_{\delta}$ is an ssa-$\delta$-space over $R$, then its cohomology groups are $R$-modules equipped with semilinear endomorphisms induced by the Frobenius lift $\phi:R\to R$   and denoted also by $\phi$. (We recall that, if $M$ and $N$ are $R$-modules, an additive map $f:M\rightarrow N$ is called {\it semilinear} if it satisfies $f(\lambda x)=\phi(\lambda)f(x)$ for all $(\lambda,x)\in R\times M$.)
If $X_{\delta}\rightarrow Y_{\delta}$ is a morphism of ssa-$\delta$-spaces over $R$, then for each integer $d\geq 0$ we have a naturally induced $R$-linear map $H^d(Y_{\delta})\rightarrow H^d(X_{\delta})$ compatible with the semilinear endomorphisms $\phi$. The map $H^0(Y_{\delta})\rightarrow H^0(X_{\delta})$
is an $R$-algebra homomorphism.

If $A$ is a $\delta$-ring we denote by $D:=D(A)\subset A$ the multiplicative  set of all $a\in A$ such that the image of $a$ in $A/pA$ is a non-zero divisor; then $\phi(D)\subset D$ and the ring $D^{-1}A$ has a natural structure of $\delta$-ring. The latter ring is a natural substitute, in our framework, for the total ring of fractions. Note that if $B$ is the zero ring then $D(B)=B$.

\begin{df}\label{df9}
{\bf (a)} An ssa-$\delta$-space $X_{\delta}$ is {\it connected} if the ring $H^0(X_{\delta})$ has no non-trivial idempotents. 

\smallskip

{\bf (b)} An  ssa-$\delta$-space $X_{\delta}$ is {\it acyclic} if $H^d(X_{\delta})=0$ for all integers $d\geq 1$.

\smallskip

{\bf (c)} A morphism of ssa-$\delta$-spaces over $R$, $X_{\delta}\rightarrow Y_{\delta}$, is a {\it quasi-isomorphism} if the homomorphisms $H^d(Y_{\delta})\rightarrow H^d(X_{\delta})$ are isomorphisms for all integers  $d\geq 0$.

{\bf (d)} An ssa-$\delta$-space $\Spec_{\delta}(B_{\bullet})$ over $R$ is {\it discrete} if for every integer $d\geq 0$ the $\delta$-ring $B_d$ is a power of $R$, i.e., $B_d=R^{\Sigma_d}$, for some set $\Sigma_d$, with the convention $R^{\emptyset}:=0$. (Then $(\Sigma_d)_{d\geq 0}$
has a natural structure of semi-simplicial set and, conversely, each semi-simplicial set gives rise to a discrete ssa-$\delta$-space as above.)

\smallskip
{\bf (e)} A {\it $\delta$-point} is an ssa-$\delta$-space over $R$ that is connected, discrete, and acyclic.

\smallskip
{\bf (f)} Given an ssa-$\delta$-space $X_{\delta}$ over $R$, a {\it $\delta$-point of} $X_{\delta}$  is a  morphism  over $R$ 
of the form $P_{\delta}\rightarrow X_{\delta}$ where $P_{\delta}$ is a $\delta$-point.

\smallskip

{\bf (g)} For an ssa-$\delta$-space $X_{\delta}=\Spec_{\delta}(B_{\bullet})$ let $D_d:=D(B_d)\subset B_d$ for $d\geq 0$ and consider the following condition:
 for all integers  $d\geq 0$ and all $i\in\{0,\ldots,d+1\}$ we have $u_{d,i}(D_d)\subset D_{d+1}$. If this condition is satisfied we obtain a naturally associated ssa-$\delta$-space $D^{-1}X_{\delta}:=\Spec_{\delta}(D^{-1}_{\bullet}B_{\bullet})$. The $\delta$-ring $H^0(D^{-1}X_{\delta})$ is called the $\delta$-{\it ring of rational functions} on $X_{\delta}$.
\end{df}

For each $\delta$-ring $S$, the ssa-$\delta$-space $\Spec_{\delta}(S_{\bullet})$ is acyclic and moreover we have $H^0(\Spec_{\delta}(S_{\bullet}))=S$. In particular, 
the ssa-$\delta$-space $\Spec_{\delta}(R_{\bullet})$ is a $\delta$-point. Other $\delta$-points will play a role later as well, see Equation (\ref{EQ038.5}). For each $\delta$-point $P_{\delta}$, the canonical morphism $P_{\delta}\rightarrow
\Spec_{\delta}(R_{\bullet})$ is a quasi-isomorphism. 
 We denote by $X_{\delta}(\delta\textup{-pts})$ the set of all $\delta$-points of $X_{\delta}$. 
 
 For each morphism $X_{\delta}\rightarrow Y_{\delta}$ of ssa-$\delta$-spaces over $R$ we have an induced map of sets $X_{\delta}(\delta\textup{-pts})\rightarrow Y_{\delta}(\delta\textup{-pts})$.
 
 One can enhance ssa-$\delta$-spaces by considering filtrations as follows:
 
 \begin{df}\label{df10}
 A {\it filtered ssa-$\delta$-space} is a semi-simplicial object in ${\bf FilRing}_{\delta}^{\textup{op}}$.
 \end{df}
 
Filtered ssa-$\delta$-spaces form a category, $\textup{\bf FilssAff}_{\delta}$, in the obvious way.
This category has fiber products and coequalizers and there exists a forgetful functor 
$\textup{\bf FilssAff}_{\delta}\rightarrow \textup{\bf ssAff}_{\delta}$ which forgets the filtrations; so we shall view
 filtered ssa-$\delta$-spaces as ssa-$\delta$-spaces with `additional structure' given by filtrations. The forgetful functor above commutes with fiber products and coequalizers. 
 The ssa-$\delta$-space $\Spec(R_{\bullet})$ has a natural structure of a filtered ssa-$\delta$-space
 with the constant filtration equal to $R$ in each degree; so we will speak about filtered ssa-$\delta$-spaces {\it over} $R$.
 If $X_{\delta}$ is a filtered ssa-$\delta$-space, then all cohomology groups
 $H^d(X_{\delta})$ inherit an exhaustive increasing filtration by subgroups 
 $H^d(X_{\delta})^r$
 which, if $X_{\delta}$ is over $R$, are $R$-submodules equipped with semilinear homomorphisms $H^d(X_{\delta})^r\rightarrow H^d(X_{\delta})^{r+1}$. If in addition we have $X_{\delta}=\Spec_{\delta}(B_{\bullet})$ and the condition in Definition
 \ref{df9} (g) is satisfied, then $D^{-1}X_{\delta}$ has a structure of a filtered ssa-$\delta$-space given by the filtration 
 $((D^{-1}_dB_d)^r)_{r\geq 0}=((D^r_d)^{-1}B^r_d)_{r\geq 0}$ of each $D^{-1}_dB_d$, where $D^r_d:=B^r_d\cap D_d$.
 
 \medskip

 A `first' example of (filtered) ssa-$\delta$-space over $R$ is attached to a quasi-compact separated smooth scheme $X$ over $\Spec(R)$ as follows. Given $n\in\mathbb Z_{\geq 0}$ and an affine open cover $X=U_0\cup \cdots \cup U_n$, we consider the filtered ssa-$\delta$-space 
$$\Spec_{\delta}\left(\prod_{0\leq i_0<\cdots<i_{\bullet}\leq n} \mathcal O^{\infty}(U_{i_0}\cap\cdots \cap U_{i_{\bullet}})\right)$$
equipped with the natural face maps. Here and in what follows, the products of the type above are taken to be equal to the zero ring for $\bullet > n$.
However, this is {\it not} the most natural example appearing in our theory. Rather, the examples that naturally occur come from certain (filtered) {\it $\delta$-graded rings} and `$\Proj_{\delta}$' spaces attached to them. To explain this construction we need to introduce a few more concepts, partly following \cite{Bu05}, Ch. 2, Sect. 2.1, Defs. 2.18 and 2.21.

\bigskip

Let $W:=\mathbb Z[\phi]$ be the commutative ring of polynomials in an indeterminate $\phi$ with $\mathbb Z$-coefficients. If $w =\sum_{i=0}^n a_i \phi^i$ with $n\in\mathbb Z_{\geq 0}$, $a_i \in\mathbb Z$ and $a_n \neq 0$, then we set $\deg(w):=\sum_{i=0}^n a_i$ and $\ord(w):=n$. We also set $\deg(0):=0$ and $\ord(0):=0$. We endow $W$ with a partial order as follows: we have $w=\sum_{i=0}^n a_i\geq 0$ if and only if $a_i\geq 0$ for all $i\in\{0,\ldots,n\}$. We set 
$$W_+:=\{w\in W|w\geq 0\}\;\;\;\textup{and}\;\;\;W(r):=\{w\in W|\ord(w)\leq r\}\; \textup{for}\; r\in\mathbb Z_{\geq 0}.$$

\begin{df}\label{df11}
A {\it $\delta$-graded ring} is a $W$-graded $\mathbb Z_{(p)}$-algebra
$B=\bigoplus_{w\in W}B(w)$
equipped with a ring endomorphism $\phi:B\rightarrow B$ such that the following three axioms hold:
 
\medskip
{\bf (GAX1)} The ring $B$ is $p$-adically separated and flat over $\mathbb Z_{(p)}$.

\smallskip

{\bf (GAX2)} For all $w\in W$ we have $\phi(B(w))\subset B(\phi w)$.

\smallskip

{\bf (GAX3)} For all $b_1,b_2\in B(w)$ we have $b_1^pb_2^{\phi}-b_2^p b_1^{\phi}\in pB((\phi+p)w)$.

\medskip

\noindent  Here, for $b\in B$ we denote simply $b^{\phi}:=\phi(b)$. A $\delta$-graded ring as above is {\it integral} if in addition it satisfies the following two conditions:

\smallskip

{\bf (GAX4)} The ring $B/pB$ is an integral domain.

\smallskip

{\bf (GAX5)} Both $\phi$ and its reduction modulo $p$ are injective.
\end{df}

Note that, by axioms (GAX1) and (GAX4), each integral $\delta$-graded ring is an integral domain.

The $\delta$-graded rings are the objects of a category 
$\textup{\bf GrRing}_{\delta}$
whose morphisms are the morphisms
of graded rings that commute with the corresponding $\phi$s.
Our definition of integral $\delta$-graded rings is closely related to (but is not identical to) the definition of $\delta$-graded rings in \cite{Bu05}, Ch. 2, Sect. 2.1, Def. 2.18: 
the latter considers non-unitary rings and uses $W_+$ instead of $W$. With notation as in (GAX3) we write 
$$\{b_1,b_2\}_{\delta}:=\frac{b_1^pb_2^{\phi}-b_2^p b_1^{\phi}}{p}\in B((\phi+p)w),$$
where division by $p$ is well-defined by axioms (GAX1) and (GAX3). 

Note that, in the above definition $\phi$ is not assumed to be a Frobenius lift and in most applications  $\phi$ will not be a Frobenius lift.
(For a basic example when $\phi$ is not a Frobenius lift see the ring in \cite{Bu05}, Eqn. (5.10).)
Thus, in general a $\delta$-graded ring does not have a natural structure of a $\delta$-ring. However, if $B$ is an integral $\delta$-graded ring, then for each $w\in W$ and every $f\in B(w)\setminus pB(w)$, the ring of fractions
$$B_{\langle f\rangle}:=\{\frac{F}{f^v}|F\in B(vw),\ v\in W_+\}$$
has a natural structure of a $\delta$-ring, where its $p$-derivation $\delta$ is defined by the formula
$$\delta\left(\frac{F}{f^v}\right):=\frac{\{f^v,F\}_{\delta}}{f^{v(\phi+p)}},$$
as one can check easily that the rule $x\mapsto x^p+\delta(x)$ defines the ring endomorphism
of $B_{\langle f\rangle}$ induced by $\phi$. Similarly, the ring of fractions
$$B_{((p))}:=\{\frac{F}{G}|F\in B(w),\ G\in B(w)\setminus pB(w),\ w\in W\},$$
has a naturally induced structure of a $\delta$-ring; it is a discrete valuation ring with maximal ideal generated by $p$.

We say that $B$ is a $\delta$-graded ring {\it over $R$} if $R$ is a subring of $B(0)$ and $\phi$ extends the Frobenius lift on $R$.

\begin{df}\label{df12}
Let $B$ be an integral $\delta$-graded ring over $R$ and let $n\in\mathbb N$ and $w\in W$. 
By a {\it linear system} in $B$ (of length $n+1\geq 2$ and weight $w\in W$) we mean an $n+1$-tuple 
$$\f=(f_0,\ldots,f_n)\in (B(w)\setminus pB(w))^{n+1}.$$
For each such pair $(B,\f)$ one defines the connected ssa-$\delta$-space 
\begin{equation}\label{EQ004}
\Proj_{\delta}(B,\f):=\Spec_{\delta}\left(\prod_{0\leq i_0<\cdots< i_{\bullet}\leq n}
 B_{\langle f_{i_0}\cdots f_{i_{\bullet}}\rangle}\right)
\end{equation}
over $R$ with the obvious face maps.\end{df}

The ssa-$\delta$-space $X_{\delta}=\Proj_{\delta}(B,\f)$ clearly satisfies 
$$H^0(X_{\delta})=\bigcap_{i=0}^n B_{\langle f_i\rangle},$$
the intersection being taken inside $B_{((p))}$. Also $X_{\delta}$ above is easily seen to satisfy the condition in Definition
 \ref{df9} (g) so we can consider the ssa-$\delta$-space $D^{-1}X_{\delta}$ and, in particular, the $\delta$-ring of rational functions $H^0(D^{-1}X_{\delta})$ on $X_{\delta}$.

If a morphism of integral $\delta$-graded rings $u:B\rightarrow \tilde{B}$ sends some linear system $\f=(f_0,\ldots,f_n)$ in $B$ into a linear system $\tilde{\f}:=(u(f_0),\ldots,u(f_n))$ in $\tilde{B}$ (of the same length and weight), then we call $\tilde{\f}$ the image of $\f$ in $\tilde B$ and we have an induced morphism $\Proj_{\delta}(\tilde{B},\tilde{\f})\rightarrow \Proj_{\delta}(B,\f)$ in $\textup{\bf ssAff}_{\delta}$.

In algebraic geometry the linear systems are vector spaces generated by global sections of a line bundle but our linear systems are tuples and not $R$-modules, as even the $0$-th cohomology groups of the ssa-$\delta$-spaces $\Proj_{\delta}(B,\f)$ really depend on $\f$ and not only on the $R$-submodule $\sum_{i=0}^n Rf_i$ of $B(w)$. 

\smallskip

For applications we will need the `filtered version' of the above concepts:

\begin{df}\label{df13}
A {\it filtered $\delta$-graded ring} is a $\delta$-graded ring $B$ as in Definition \ref{df11} which is also equipped with a structure of a filtered graded ring, given by a filtration $(B^r)_{r\geq 0}$, for which the following two axioms hold:

\smallskip

{\bf (GAX6)} For all integers $r\geq 0$ we have $\phi(B^r)\subset B^{r+1}$.

\smallskip

{\bf (GAX7)} For all integers $r\geq 0$ and all $w\in W$ we have 
$$\{B^r(w),B^r(w)\}_{\delta}:=\{\{b_1,b_2\}_{\delta}|b_1,b_2\in B^r(w)\}\subset B^{r+1}((\phi+p)w).$$
\end{df}

Note that if (GAX6) holds and $B^r\cap pB^{r+1}=pB^r$ for all integers  $r\geq 0$, then (GAX7) follows from (GAX3).
Filtered $\delta$-graded rings are the objects of a category 
$\textup{\bf FilGrRing}_{\delta}$
whose morphisms are the morphisms of $\delta$-graded rings that are also morphisms of filtered rings.

We similarly speak about filtered $\delta$-graded rings over $R$ and about filtered integral $\delta$-graded ring.

If $B$ is a filtered integral $\delta$-graded ring over $R$, then for each $w\in W$ and every $f\in B(w)\setminus pB(w)$, the $\delta$-ring $B_{\langle f\rangle}$ has a structure of a filtered $\delta$-ring with filtration $((B_{\langle f\rangle})^r)_{r\geq 0}$ given by 
$$(B_{\langle f\rangle})^r:=\{\frac{F}{f^v}|F\in B^r(vw), v\in W_+\}.$$
Similarly, the $\delta$-ring $B_{((p))}$ has a structure of a filtered $\delta$-ring with filtration 
$((B_{((p))})^r)_{r\geq 0}$ given by
$$(B_{((p))})^r:=\{\frac{F}{G}|F\in B^r(w), G\in B^r(w)\setminus pB^r(w), w\in W\}.$$
For each linear system $\f$ in a filtered integral $\delta$-graded ring $B$ over $R$,
 the ssa-$\delta$-space $\Proj_{\delta}(B,\f)$ has an induced structure of filtered ssa-$\delta$-space.
If a  morphism $B\rightarrow\tilde B$ in $\textup{\bf FilGrRing}_{\delta}$, where $B$ and $\tilde{B}$ are filtered integral $\delta$-graded rings, maps a linear system $\f=(f_0,\ldots,f_n)$ in $B$ into a linear system $\tilde{\f}:=(u(f_0),\ldots,u(f_n))$ in $\tilde{B}$, 
and if $\lambda \in R^{\times}$,
then  we have an induced morphism  $\Proj_{\delta}(\tilde{B},\lambda\tilde{\f})\rightarrow \Proj_{\delta}(B,\f)$ in $\textup{\bf FilssAff}_{\delta}$ where $\lambda\tilde{\f}:=(\lambda u(f_0),\ldots,\lambda u(f_n))$.
 
\subsection{Concepts of $\delta$-dimension}\label{S22}

Throughout this subsection  we use the following convention: a family $(x_i)_{i\in \mathcal I}$ of elements of a ring is called {\it algebraically independent} over a subring if the members of the family are distinct
(i.e., $x_i\neq x_j$ for $i,j\in \mathcal I$, $i\neq j$) and the set $\{x_i|i\in \mathcal I \}$ is algebraically independent over the subring. We make a similar convention for linearly independent.

If $B_1\subset B_2$ is an inclusion of integral domains, let 
$$\trdeg_{B_1}(B_2)\in\mathbb N \cup \{0,\infty\}$$ 
be the supremum of all $d\in\mathbb N \cup \{0\}$ such that $B_2$ (or $\Frac(B_2)$) contains $d$ algebraically independent elements over $B_1$; here we identify the elements of 
$\mathbb N \cup \{0,\infty\}$ with cardinals.  We have
 $\trdeg_{B_1}(B_2)=\trdeg_{\Frac(B_1)}(\Frac(B_2))$.

In this subsection we take $B_0$ to be an integral domain thought of as the `ground ring'. (In Sections \ref{S2} and \ref{S3}, $B_0$ will be $R=W(k)$, while in Section \ref{S4}, $B_0$ will be an algebraically closed field.) All concepts will be relative to $B_0$ and we will often drop the reference to $B_0$ if no confusion can arise. Throughout this subsection $B$ is an integral domain that contains $B_0$. 

\begin{df}\label{df14}
Let $B$ be a filtered ring with filtration $(B^r)_{r\geq 0}$ 
such that $B_0\subset B^0$.
We say that $B$ has {\it asymptotic transcendence degree} $\alpha\in {\mathbb R}$, and we write 
$$\asytrdeg(B)=\alpha,$$
if $\trdeg_{B_0}(B^r)<\infty$ for all integers  $r\geq 0$ and 
$$\lim_{r\rightarrow \infty} \frac{\trdeg_{B_0}(B^r)}{r}=\alpha.$$
We say that $B$ has {\it asymptotic dimension} $\beta \in {\mathbb R}$, and we write 
$$\asydim(B)=\beta,$$
if $\dim(B^r \otimes_{B_0} \Frac(B_0))<\infty$ for all integers  $r\geq 0$ and 
$$\lim_{r\rightarrow \infty} \frac{\dim(B^r \otimes_{B_0} \Frac(B_0))}{r}=\beta.$$ 
\end{df}

Here $\dim$ denotes Krull dimension. If all $\Frac(B_0)$-algebras $B^r\otimes_{B_0}\Frac(B_0)$ are finitely generated, then $\asytrdeg(B)$ exists if and only
$\asydim(B)$ exists and, in case they exist, they are equal. A typical situation when $\asydim(B)$ exists but $\asytrdeg(B)$ does not arises as follows. Let $X$ be an affine smooth scheme over $\Spec(R)$ of relative dimension $d$ and we take
$$B_0=R,\ \ B^r=\mathcal O(J^r(X)),\ \ B=\mathcal O(J^{\infty}(X))=\bigcup_{r\geq 0} B^r.$$
Then by Equation (\ref{EQ002}) we have $\dim(B^r \otimes_R K)=(r+1)d$, so
\begin{equation}\label{EQ006}
\asydim(\mathcal O(J^{\infty}(X)))=d,\end{equation}
but $\asytrdeg(\mathcal O(J^{\infty}(X)))$ does not exist if $d\geq 1$ because the transcendence degrees involved in its definition are infinite.

\bigskip

We recall another concept of dimension which comes from difference algebra (see \cite{Le}) and which will play a role in our paper. 

Assume we are given an injective ring endomorphism 
$\phi:B\rightarrow B$
that sends $B_0$ onto $B_0$. So $\phi$ induces a ring endomorphism of $\Frac(B)$ that sends $\Frac(B_0)$ onto itself.
(In Sections \ref{S2} and \ref{S3}, where $B_0=R$, the restriction of $\phi$ to $B_0$ will be the unique Frobenius lift on $R$  but $\phi$ on $B$ will not necessarily be a Frobenius lift. In Subsections \ref{S42} to \ref{S44}, where $B_0$ is an algebraically closed field $\mathbb K$, the endomorphism $\phi$ will be denoted by $\sigma$ and its restriction to $\mathbb K$ will be the identity automorphism $1_{\mathbb K}$.) For 
$$(w,b)=(\sum_{i=0}^n a_i\phi^i,b)\in (W\times B^{\times})\cup (W_+\times B),$$ 
we extend the notation $b^{\phi}=\phi(b)$ by setting
$$b^w:=b^{a_0} \phi(b)^{a_1}\cdots
\phi^n(b)^{a_n}\in B.$$

Let $C$ be a $B_0$-subalgebra of $B$ such that $\phi(C)\subset C$ (i.e., a $B_0[\phi]$-subalgebra of $B$).
 A subset $\bowtie$ of $\Frac(B)$ will be called $\phi$-{\it algebraically independent} over $\Frac(C)$ if the family
$(\phi^j(b))_{(b,j)\in\bowtie\times (\mathbb N\cup\{0\})}$
of elements of $\Frac(B)$
is algebraically independent over $\Frac(C)$. A subset $\Frac(B)$ is called a $\phi$-{\it transcendence basis} of $\Frac(B)$ over $\Frac(C)$ if it is maximal in the set of $\phi$-algebraically independent subsets of $\Frac(B)$.  As such, the $\phi$-{\it transcendence degree} of $B$ 
over $C$ will be the supremum 
$\phidim_{C}(B)\in\mathbb N \cup \{0,\infty\}$
of all numbers $d\in \mathbb N \cup \{0\}$ such that $\Frac(B)$ has a subset with $d$ elements which is $\phi$-algebraically independent over $\Frac(C)$. We will use the simpler notation $\phidim(B):=\phidim_{B_0}(B)$.

 A subset $\{b_1,\ldots,b_n\}\subset \Frac(B)$ is said to $\phi$-{generate} the field $\Frac(B)$ over $\Frac(B_0)$ if the set
$\{\phi^j(b_i)|(i,j)\in\{1,\ldots,n\}\times (\mathbb N\cup\{0\})\}$
generates the field $\Frac(B)$ over $\Frac(B_0)$. We say $\Frac(B)$ is $\phi$-finitely generated over $\Frac(B_0)$ if there exist a finite subset of $B$ that $\phi$-generates $\Frac(B)$ over $\Frac(B_0)$. 
Similarly, a subset $\{b_1,\ldots,b_n\}\subset B$ is said to $\phi$-{\it generate} $B$ as a $B_0$-algebra if the set
$\{\phi^j(b_i)|(i,j)\in\{1,\ldots,n\}\times (\mathbb N\cup\{0\}\}$
generates $B$ as a $B_0$-algebra. We say $B$ is $\phi$-finitely generated over $B_0$ if there exist a finite subset of $B$ that $\phi$-generates $B$ as a $B_0$-algebra. 

We end this subsection with remarks on interrelations between the various dimension concepts introduced above.

\begin{rem}\label{R1}
The following seven properties hold:
 
 \medskip
 {\bf (a)} A $\phi$-transcendence basis of $\Frac(B)$ over $\Frac(C)$ always exists and every two such $\phi$-transcendence bases have the same cardinality (see \cite{Le}, Ch.4, Prop. 4.1.6). 
 
 
 \smallskip
 {\bf (b)} We have the additivity property (see \cite{Le}, Ch. 4, Thm. 4.1.9):
\begin{equation}\label{EQ007}
\phidim(B)=\phidim(C)+\phidim_{C}(B).\end{equation}

 \smallskip
 {\bf (c)} If $\Frac(B)$ is $\phi$-finitely generated over $\Frac(B_0)$, then $\Frac(C)$ is also $\phi$-finitely generated over $\Frac(B_0)$ (see \cite{Le}, Ch. 4, Thm. 4.4.1). 
 
 \smallskip
{\bf (d)} Each subset of $\Frac(B)$ that $\phi$-generates $\Frac(B)$ over $\Frac(C)$ contains a $\phi$-transcendence basis (see \cite{Le}, Ch. 4, Thm. 4.1.8).

 \smallskip
{\bf (e)} Assume $\{\beta_1,\ldots,\beta_n\}\subset \Frac(C)$ is $\phi$-algebraically independent over $\Frac(B_0)$ and assume $\{\gamma_1,\ldots,\gamma_m\}\subset\Frac(C)$ is a set which $\phi$-generates $\Frac(C)$ over $\Frac(B_0)$. Then $m\geq n$ and the set
$\{\gamma_1,\ldots,\gamma_m\}$ contains a subset with $n$ elements which is $\phi$-algebraically independent over $\Frac(B_0)$. To check this statement, we first remark that $\phidim(C)\geq n$. By the property (d) applied to the fields $\Frac(B_0)$ and $\Frac(C)$ we have that $\{\gamma_1,\ldots,\gamma_m\}$ contains a $\phi$-transcendence basis of $\Frac(C)$ over $\Frac(B_0)$ which, by the property (a), has $\phidim(C)$ elements and the statement follows.

 \smallskip
{\bf (f)} If $\phidim(B)=n\in \mathbb N$, then there exists a subset of $B$ with $n$ elements which is $\phi$-algebraically independent over $\Frac(B_0)$. To check this statement, let $\{\frac{b_1}{b_0},\ldots,\frac{b_n}{b_0}\}$ be a $\phi$-transcendence basis of $\Frac(B)$ over $\Frac(B_0)$, with $b_0,\ldots,b_n\in B$. Then the statement follows from the property (e) applied to $\{b_0,\ldots,b_n\}$ with $C:=B_0[\phi^j(b_i)|0\leq i\leq n, j\geq 0]$ and $m=n+1$. 

 \smallskip
{\bf (g)} Assume $B$ has a grading by an abelian semigroup in which the homogeneous components are $B_0$-submodules. Then $\phidim(B)$ is the supremum 
of all natural numbers $d$ such that $B$ has $d$ homogeneous elements that are
 $\phi$-algebraically independent over $\Frac(B_0)$. To check this statement, let $\{\frac{b_1}{b_0},\ldots,\frac{b_n}{b_0}\}\subset \Frac(B)$ be $\phi$-algebraically independent over $\Frac(B_0)$, where $b_0,\ldots,b_n\in B$ and let $\{\gamma_1,\ldots,\gamma_m\}$ be the set of all homogeneous components of $b_0,\ldots,b_n$. If $C:=B_0[\phi^j(\gamma_i)|0\leq i\leq m, j\geq 0]$, then the statement follows from the property (e).\end{rem}
 
\begin{rem}\label{R2}
We mention the following  related concept. Let $E$ be a $B_0$-module equipped with a semilinear map $\phi:E\rightarrow E$, i.e., $\phi$ on $E$ is additive and for all $\lambda\in B_0$ and all $x\in E$ we have $\phi(\lambda x)=\phi(\lambda)\phi(x)$. By the $\phi$-{\it rank} of $E$ we mean the supremum, denoted by
 $\phirank_{B_0}(E)\in\mathbb N \cup \{0,\infty\}$, of all $n\in \mathbb N\cup\{0\}$ such that there exist elements $x_1,\ldots,x_n \in E$ that are $\phi$-{\it linearly independent} over $B_0$ in the sense that the family $(\phi^j(x_i))_{(i,j)\in\{1,\ldots,n\}\times (N \cup \{0\})}$ is $B_0$-linearly independent.
 \end{rem}

 \begin{rem}\label{2.5}
 Let $B$ be a $W$-graded ring where $W=\mathbb Z[\phi]$. We assume that $B_0\subset B(0)$ and that we have $\phi(B(w))\subset B(\phi w)$ for all $w\in W$. We set
 $$W_B:=\{w\in W|B(w)\neq 0\}.$$
 If we view $W$ as a module over the semiring $W_+$, then $W_B$ is a 
 $W_+$-submodule of $W$ (recall $B$ is an integral domain and $\phi:B\to B$ is injective).
 It is easy to check that if $B\otimes_{B_0} \Frac(B_0)$ is $\phi$-finitely generated as a $\Frac(B_0)$-algebra, then $W_B$ is a finitely generated $W_+$-module.
\end{rem}

\begin{rem}\label{R3}
Assume $B$ is a filtered ring with filtration $(B^r)_{r\geq 0}$ and assume we are given   an injective ring endomorphism $\phi$ with $\phi(B_0)=B_0$
 {\it compatible} with the filtration in the sense that for all integers $r\geq 0$ we have
$$\phi(B^r)\subset B^{r+1}.$$
If $B$ has an asymptotic transcendence degree, then we have an inequality
\begin{equation}\label{EQ008}
\phidim(B)\leq\asytrdeg(B).
\end{equation}
Note that the above inequality is not an equality in general: the right-hand side effectively depends on the filtration whereas the left-hand side does not.
If $C$ is $B_0$-subalgebra of $B$ which is a 
filtered ring with filtration $(C^r)_{r\geq 0}$ such that we have $C^r\subset B^r$ and $\phi(C^r)\subset C^{r+1}$ for all integers $r \geq 0$, and if $C$ also has an asymptotic transcendence degree, then it is easy to check, using Equation (\ref{EQ007}) and Inequality (\ref{EQ008}),
 that we have an implication
\begin{equation}\label{EQ009}
\asytrdeg(C)=\asytrdeg(B)\ \ \ \Rightarrow \ \ \phidim(C)=\phidim(B).
\end{equation}
Even if either $C$ or $B$ does not have an asymptotic transcendence degree, for an integer $N\geq 0$ we have the implications $\circled{1}\Rightarrow\circled{2}\Rightarrow\circled{3}$ between the following statements:

\smallskip

\circled{1} We have $\trdeg_{C^r}(B^r)\leq N$ for all integers  $r>>0$.

\smallskip

\circled{2} We have $\trdeg_C(B)\leq N$.

\smallskip

\circled{3} We have $\phidim(C)=\phidim(B)$.
\end{rem}
 
\begin{rem}\label{R4}
The concepts of $\phi$-transcendence degree and asymptotic transcendence degree  are related as follows. If $B$ is equipped with an injective ring endomorphism $\phi$ with $\phi(B_0)= B_0$ such that $\Frac(B)$ is $\phi$-finitely generated over $\Frac(B_0)$, then it is easy to see that there exists an increasing exhausting filtration $(B^r)_{r\geq 0}$ by $B_0$-subalgebras of $B$ such that for all integers $r\geq 0$ we have $\phi(B^r)\subset B^{r+1}$ and the asymptotic transcendence degree exists and is given by
 $$\asytrdeg(B)=\phidim(B).$$
 This will not be used here as the filtrations appearing in the paper
 are given (cannot be modified) and Inequality (\ref{EQ008}) is not a priori an equality.\end{rem}
 
 \begin{rem}\label{R5}
In this remark we assume that the inclusion  $B_0\subset B$ is a morphism of $\delta$-rings. Let $\phi: B\rightarrow B$ be the Frobenius lift attached to $\delta:B\rightarrow B$. Then from Remark \ref{R1} (g) we get that $\phidim(B)$ coincides with the supremum of the set of all integers $d$ such that there exist elements $b_1,\ldots,b_d\in B$ with the property that the family
$$(\delta^j(b_i))_{(i,j)\in\{1,\ldots,d\}\times (\mathbb N \cup \{0\})}$$ 
of elements in $\Frac(B)$ is algebraically independent over $K$. A similar statement holds for the concept of finite $\phi$-generation.
 \end{rem}
 
 \begin{rem}\label{R6}
Let $B_0=R$ and let $B$ be an integral $\delta$-graded ring over $R$. We assume 
that $B(0)=R$ and 
 that there exists $w\in W\setminus \{0\}$ such that we have $B(w)\neq pB(w)$. We first state and then argue the following five properties .
 
\medskip

{\bf (a)} We have an equality
$$\phidim(B)=\phidim(B_{((p))})+1.$$

\smallskip

{\bf (b)}
The $\phi$-transcendence degree $\phidim(B)$ is the supremum of the set of all  $n\in \mathbb N\cup\{0\}$ such that $B$ contains $n$ elements that are homogeneous of the same weight and $\phi$-algebraically independent.
 
\smallskip
 
{\bf (c)} Assume  $w\geq 1$ or $w\leq -1$.  
(For instance, the inequality $w\geq 1$ holds in the example of the `projective space'
 introduced below by Equation (\ref{EQ032}) and the inequality 
 $w\leq -1$ holds in some of our main examples related to moduli schemes of abelian schemes, see Remark \ref{R23} (c) of Subsubsection \ref{S3603} or Family (\ref{EQ130}) below.)
Let $f\in B(w)\setminus pB(w)$. Then $p$ is a prime element in the ring $B_{\langle f \rangle}$ and the localization of this ring at the prime ideal $(p)$ satisfies 
\begin{equation}\label{EQ010}
(B_{\langle f \rangle})_{(p)}=B_{((p))}.\end{equation} 
In particular, we have
$$\phidim(B)=\phidim(B_{\langle f \rangle})+1.$$

 \smallskip
 
 {\bf (d)} Assume  $w\geq 1$ or $w\leq -1$ and assume in addition that $B$ is a filtered $\delta$-graded ring over $R$. Consider a  linear system $\f=(f_0,\ldots,f_n)\in (B(w)\setminus pB(w))^{n+1}$, and let
$X_{\delta}:=\Proj_{\delta}(B_{\bullet},\f)$.
Then, recalling the notation in Definition \ref{df9} (g),
defining $N_{n,\bullet}:=\left(\begin{array}{c} n+1\\ \bullet+1\end{array}\right)$ we have
 \begin{equation}\label{EQ005}
 D_{\bullet}^{-1}B_{\bullet}=(B_{((p))})^
 {N_{n,\bullet}},\end{equation}
 with the canonical coface maps, so 
 the ssa-$\delta$-space 
 $D^{-1}X_{\delta}$ is acyclic and we have
$$H^0(D^{-1}X_{\delta})=B_{((p))}.$$ 

\smallskip

{\bf (e)} Assume we are in the situation of part (d) and assume,  in addition, that $B^r(v)=0$ for all integers  $r\geq 0$ and all $v\in W$ such that $\ord(v)>r$. 
(This condition is trivially satisfied in our main examples later in the paper, see the paragraph after Definition \ref{df18} of Subsubsection \ref{S231})
Then the canonical filtrations on the two rings in Equation (\ref{EQ005}) satisfy
\begin{equation}\label{EQ011}
(D_{\bullet}^{-1}B_{\bullet})^r\subset ((B_{((p))})^r)^{N_{n,\bullet}}\subset (D_{\bullet}^{-1}B_{\bullet})^{r+\ord(w)}.\end{equation}
In particular we have the following formula for the asymptotic transcendence degree of the $\delta$-ring of rational functions on $X_{\delta}$:
$$\asytrdeg(H^0(D^{-1}X_{\delta}))=\asytrdeg(B_{((p))}).$$

In what follows we give the arguments for properties (a) to (e).

The inequality `$\geq$' in part (a) follows because each non-zero element in $B(w)$ is 
$\phi$-algebraically independent over $B_{((p))}$. To check the inequality `$\leq$' in part (a)
assume $\phidim(B)\geq n+1$ for some integer $n\geq 0$, and consider
$\phi$-algebraically independent homogeneous elements $b_0,\ldots,b_n\in B$ over $R$, $b_i\in B(w_i)\setminus pB(w_i)$, see Remark \ref{R1} (g). So the family
$(\prod_{i=0}^n b_i^{w'_i})_{(w'_0,\ldots,w'_n)\in W^{n+1}}$ is linearly independent over $R$.
As $B(0)=R$ we have $w_i\neq 0$ for all $i\in\{0,\ldots,n\}$. 
Let $v\in W$ be a non-zero common multiple of $w_0,\ldots,w_n$
and  for $i\in\{0,\ldots,n\}$ let $v_i\in W$ be such that $v=v_iw_i\in W$.
As $v_i\neq 0$ for all $i\in\{0,\ldots,n\}$, the family
$(\prod_{i=0}^n b_i^{v_iw'_i})_{(w'_0,\ldots,w'_n)\in W^{n+1}}$ of elements in $\Frac(B)$ is linearly independent over $R$.
Thus the family $(b_0^{v_i})_{i\in\{0,\ldots,n\}}$ of elements of $\Frac(B)$ is $\phi$-algebraically independent over $R$, hence the family $(b_0^{v_i}/b_0^{v_0})_{i\in\{1,\ldots,n\}}$ is algebraically independent over $R$. As the reduction of $\phi:B\rightarrow B$ modulo $p$ is injective we have
$b_i^{v_i}/b_0^{v_0}\in B_{((p))}$ for all $i\in\{1,\ldots, n\}$. Therefore $\phidim(B_{((p))})\geq n$. From this the inequality `$\leq$' follows.

To check part (b), assume again that $\phidim(B)\geq n+1$ and let $b_0,\ldots,b_n$ and the notation be as as in
 the prior paaragraph. For $i\in\{1,\ldots,n\}$, write
 $b_i^{v_i}/b_0^{v_0}=\beta_i/\gamma_i$, $\beta_i\in B(u_i)$, $\gamma_i\in B(u_i)\setminus pB(u_i)$, $u_i\in W$. We choose $u_i$ such that 
 $u:=\sum_{i=1}^n u_i\neq 0$
 (as we can multiply $\beta_1$ and $\gamma_1$ by some element in $B(w)\setminus pB(w)$). Let $\gamma:=\prod_{i=1}^n \gamma_i\in B(u)$. Then $\gamma, \beta_1/\gamma_1,\ldots,\beta_n/\gamma_n$ are $\phi$-algebraically independent over $R$
 (because, as $u\neq 0$, each element in $B(u)$ is $\phi$-algebraically independent over $B_{((p))}$). Hence $\gamma,\beta_1\gamma/\gamma_1,\ldots,\beta_n\gamma/\gamma_n\in 
 B(u)$ are $\phi$-algebraically independent over $R$ and are homogeneous of the same weight which proves part (b).
 
 To check part (c) note first that, by the axiom (GAX5), $p$ does not divide $f^{\phi^i}$  for all integers $i\geq 0$. As, by the axiom (GAX4), $p$ is prime in $B$  it follows that $p$ is  prime in $B_{\langle f\rangle}$.
 The inclusion `$\subset$' in Equation (\ref{EQ010})  is obvious. To check the inclusion `$\supset$'
 let $F/G\in B_{((p))}$ where $F\in B(v)$ and $G\in B(v)\setminus pB(v)$ for some $v\in W$ with $\ord(v)\leq r$. Write $v=v_+-v_-$ with $v_+,v_-\geq 0$ and $\ord(v_+),\ord(v_-)\leq r$. 
 Then for $w\geq 1$ we can write
 \begin{equation}\label{EQ012}
 \frac{F}{G}=\frac{F^{w}}{GF^{w-1}}=
 \frac{(F^{w}f^{v_-})/f^{v_+}}{(GF^{w-1}f^{v_-}) /f^{v_+}}\in (B_{\langle f \rangle})_{(p)}
 \end{equation}
 while for $w\leq -1$ we can write
 \begin{equation}\label{EQ012.681}
 \frac{F}{G}=\frac{F^{-w}}{GF^{-w-1}}=
 \frac{(F^{-w}f^{v_+})/f^{v_-}}{(GF^{-w-1}f^{v_+}) /f^{v_-}}\in (B_{\langle f \rangle})_{(p)}
 \end{equation}
 which implies that the inclusion `$\supset$' in Equation (\ref{EQ010}) holds.
 
 Part (d) follows directly from Equation (\ref{EQ010}).
 
To check part (e) note  that the first inclusion in Equation (\ref{EQ011}) is obvious. To check the second inclusion it is enough to note that if $F/G\in B_{((p))}$ with $F\in B^r(v)$ and $G\in B^r(v)\setminus pB^r(v)$ as in Equation (\ref{EQ012}) (resp. (\ref{EQ012.681})), then
 $F^{w}$ and $GF^{w-1}$ (resp. $F^{-w}$ and $GF^{-w-1}$) have order at most $r+\ord(w)$
 while $f^{v_+}$ (resp. $f^{v_-}$) has order at most $\ord(w)+\ord(v)\leq \ord(w)+r$.

 \end{rem}

\subsection{Concepts of $\delta$-invariant theory}\label{S23}

We review some  basic concepts introduced in \cite{Bu00, BB, Bu03, Bu05} and we introduce some new concepts that are needed in the paper. 

\subsubsection{Rings of abstract $\delta$-modular forms}\label{S231}

 We develop a formalism of $\delta$-modular functions
 and forms in an abstract (purely categorical) setting of what we will call {\it moduli data}.
 In the next subsubsection we will specialize our discussion to the case of abelian schemes and their isogenies.
 
 Fix a ring $S_0$ which we view as a `ground ring' and which, for simplicity,
  we usually drop from 
  our notation. All sets will be assumed to belong to a given universe and all categories will be small. Denote by $\textup{\bf Ring}$, 
 $\textup{\bf Cat}$, $\textup{\bf Gr}$, and $\textup{\bf Set}$ the category of rings over $S_0$, the category of categories, the category of groups, and the category of sets. 
 (So in the case of  $\textup{\bf Ring}$ the notation matches  the one in Subsection \ref{S21} if we take $S_0=\mathbb Z_{(p)}$; however  we will also be interested later in the case $S_0=R$.)
Let $\iota:\textup{\bf Gr}\rightarrow\textup{\bf Cat}$ be the natural functor and let $\mathbb G_m:\textup{\bf Ring}\rightarrow\textup{\bf Gr}$ be the multiplicative group functor. 

\begin{df}\label{df15}
A {\it moduli datum} (over $S_0$) is a quadruple $(\mathfrak M,\mathfrak D,G,D)$ that consists of:
 
 {\bf (MAX1)} a functor $\mathfrak M:\textup{\bf Ring}\rightarrow \textup{\bf Cat}$;
 
 \smallskip
 
 {\bf (MAX2)} a morphism $\mathfrak D:\mathfrak M \rightarrow \iota \circ \mathbb G_m$ of functors from $\textup{\bf Ring}$ to $\textup{\bf Cat}$;
 
 \smallskip
 
 {\bf (MAX3)} a functor $G:\textup{\bf Ring}\rightarrow\textup{\bf Gr}$ that acts on $\mathfrak M$;
 \smallskip

{\bf (MAX4)} a morphism $D:G\rightarrow \mathbb G_m$ of functors from $\textup{\bf Ring}$ to $\textup{\bf Gr}$.
\end{df}

Explicitly, (MAX1)-(MAX4) are given by a rule that functorially attaches to each ring $S$ over $S_0$ the following:
\smallskip

{\bf (1)} a category $\mathfrak M(S)$ with set of 
 objects $\pmb{\M}(S)$, set of morphisms $\pmb{\M}'(S)$, and source, target, identity, and composition maps $$\begin{array}{rcc}
 \pi_1(S),\pi_2(S):\pmb{\M}'(S) & \rightarrow & \pmb{\M}(S),\\
 \epsilon(S):\pmb{\M}(S) & \rightarrow & \pmb{\M}'(S), \\
 \mu(S):\pmb{\M}''(S) & \rightarrow & \pmb{\M}'(S)\end{array}$$
 respectively, denoted in what follows simply by $\pi_1,\pi_2,\epsilon,\mu$, where 
 $$\pmb{\M}''(S):=\{(x'_1,x'_2)\in \pmb{\M}'(S)\times 
 \pmb{\M}'(S),\ \pi_2(x'_1)=\pi_1(x'_2)\};$$ 
 
 \smallskip
 
 {\bf (2)} a multiplicative map $\mathfrak D(S):\pmb{\M'}(S)\rightarrow S^{\times}$, denoted in what follows simply by $\mathfrak D$; 
 
 \smallskip
 
 {\bf (3)} a group $G(S)$ acting on $\pmb{\M}(S)$ and $\pmb{\M}'(S)$ such that all the maps $\pi_1,\pi_2,\epsilon$ are $G(S)$-invariant, and moreover, with $GS(s)$ acting naturally on $\pmb{\M}''(S)$, $\mu$ is as well $G(S)$-invariant;

 \smallskip

 {\bf (4)} a group homomorphism $D(S):G(S)\rightarrow S^{\times}$, denoted in what follows simply by $D$. 
 
 \smallskip
 
 We write, as usual, $x'_2\circ x'_1:=\mu(x'_1,x'_2)$ for 
 $(x'_1,x'_2)\in \pmb{\M}''(S)$. Then
 multiplicativity in (2) means that 
 $\mathfrak D(x'_2\circ x'_1)=\mathfrak D(x'_2)\mathfrak D(x'_1)$. 
 
 \begin{rem}\label{R6.1}
 Assume $(\mathfrak M,\mathfrak D,G,D)$ is a moduli datum.
 Consider the equivalence relation $\sim$ on $\pmb{\M}(S)$ 
 where
 for $x,y\in \pmb{\M}(S)$ we have 
 $x \sim y$ if and only if 
 there exists $g\in G(S)$ and an isomorphism between $y$ and $g\cdot x$. 
We get a functor 
$$\mathcal A=\mathcal A_{\mathfrak M,G}:\textup{\bf Ring}\rightarrow\textup{\bf Set}$$
 that attaches to each $S\in \Ob(\textup{\bf Ring})$ 
 the quotient set
 $$\mathcal A(S):=\pmb{\M}(S)/\sim.$$
 Similarly, consider the equivalence relation $\sim'$ on $\pmb{\M}'(S)$ 
 where for $x',y'\in \pmb{\M}'(S)$ we have
 $x' \sim' y'$ if and only if 
 there exists $g\in G(S)$ and isomorphisms $z'_1,z'_2\in \pmb{\M}'(S)$ such that 
 $g\cdot x'=z'_1\circ y'\circ z'_2$.
We get a functor 
$$\mathcal A'=\mathcal A'_{\mathfrak M,G}:\textup{\bf Ring}\rightarrow \textup{\bf Set}$$ 
 that attaches to each $S\in \Ob(\textup{\bf Ring})$ the quotient set
 $$\mathcal A'(S):=\pmb{\M}'(S)/\sim'.$$
 We have a natural pair of source and target morphism of functors
 \begin{equation}
 \label{EQ012.1}
 \mathcal A'\stackrel{\pi_1,\pi_2}{\rightrightarrows} \mathcal A.
 \end{equation}
Let $\mathcal A/\mathcal A':\textup{\bf Ring}\rightarrow \textup{\bf Set}$ be the functor
 that attaches to each $S\in \Ob(\textup{\bf Ring})$ the coequalizer of the two maps $\mathcal A'(S)\rightrightarrows \mathcal A(S)$ defined by the functors (\ref{EQ012.1}). We have a natural quotient morphism of functors
 \begin{equation}
 \label{EQ012.11}
 \mathcal A\stackrel{\pi}{\rightarrow} \mathcal A/ \mathcal A'.
 \end{equation}

If the $S_0$-schemes are being identified with functors from 
 $\textup{\bf Ring}$ to $\textup{\bf Set}$, then, in all applications we have in mind
 (e.g., in the situations of Remark \ref{R6.75} and Subsubsection \ref{S232} below)
 the functors
 $\mathcal A$ and $\mathcal A'$ are `very close to being schemes': typically, `appropriate sheafifications' of them have \'{e}tale covers that are schemes (as a consequence of results on moduli spaces in `classical' algebraic geometry, see e.g., \cite{MFK}, Ch. 7, Thm. 7.9). But in these applications, the 
 functor $\mathcal A/\mathcal A'$
 is `very far from being a scheme.' To remedy this situation one is led to search for
 an `enlargement' of algebraic geometry, and our $\delta$-geometry does remedy it.
 \end{rem}
 
 Recall the partially ordered ring $W=\mathbb Z[\phi]$, the ring homomorphism $\deg:W\rightarrow \mathbb Z$, the map $\ord:W\rightarrow \mathbb N \cup \{0\}$, and the cone $W_+\subset W$ of elements $\geq 0$. 
For $S\in\Ob(\textup{\bf FilRing}_{\delta})$ and $(w,b)\in (W\times (S^m)^\times)\cup (W_+\times S^m)$, with $m\in\mathbb N\cup\{0\}$ and the writing $w=\sum_{i=0}^n a_i \phi^i$, we set
$$b^w:=b^{a_0} \phi(b)^{a_1}\cdots
\phi^n(b)^{a_n}\in S^{m+n}.$$

 Assume in what follows that the ground ring is $S_0=R$.

 \begin{df}\label{df16}
 Let $r\geq 0$ be an integer.
 A {\it $\delta$-modular function} of order $\leq r$ for a given moduli datum $(\mathfrak M,\mathfrak D,G,D)$ is a rule, $f$, that associates to each
object $S$ of $\Prol$
and to every $x \in \pmb{\M}(S^0)$
an element
$f(x,S) \in S^r$ such that $f(x,S)$ depends only on $S$ and on the isomorphism class of $x$ and such that the association $(x,S)\mapsto f(x,S)$ is functorial in $S$.
\end{df}

\noindent For each integer $r\geq 0$, denote by $M^r$ the set of $\delta$-modular functions of order $\leq r$; it has a structure of $p$-adically complete flat $R$-algebra
and there exists a natural injective $R$-algebra homomorphism $M^r\rightarrow M^{r+1}$. We set 
$$M:=\bigcup_{r\geq 0}M^r.$$
The filtered ring $M$ is an object of $\Prol$ with $p$-derivation defined by $(\delta f)(x,S):=\delta(f(x,S))$ for $f\in M$, $x\in \pmb{\M}(S^0)$, $S\in \Ob(\Prol)$. We denote by  $\phi$
the Frobenius lift on $M$ attached to $\delta$.

If in the above discussion we replace $\pmb{\M}$ by $ \pmb{\M}'$ we obtain corresponding rings $M^{\prime r}$ and $M'$. The morphisms $\pi_1,\pi_2$ in the discussion after Definition
\ref{df15}
 induce morphisms of filtered $\delta$-rings $\pi_1^*,\pi_2^*:M\rightarrow M'$, $f\mapsto \pi_i^* f$, $i\in \{1,2\}$, defined by
$$(\pi_i^*f)(x',S):=f(\pi_i(x'),S).$$
The morphisms $\pi_i^*$ and their reductions modulo $p$, $\overline{\pi_i^*}:\overline{M^r}\rightarrow \overline{M^{\prime r}}$ are injective, as one can see by taking $x'=\epsilon(x)$ in the above equality and using that $S^r$ is flat over $R$.

 \begin{df}\label{df17}
 Let $w\in W(r)$. An element $f\in M^r$ has {\it weight} $w$ if for each $g\in G(S^0)$ and every $x\in \pmb{\M}(S^0)$ the following equality holds in $S^r$:
$$
f(g\cdot x, S)= D(g)^{-w} \cdot
f(x,S).$$
\end{df}

Here  $D$ denotes the map
 $D(S^0):G(S^0)\rightarrow (S^0)^{\times}$; see point (4) after Definition \ref{df15}.
 
We denote by $M^r(w)$ the $R$-module of all elements of $M^r$ of weight $w$ and we refer to them as {\it $\delta$-modular forms} of weight $w$.
If in Definition \ref{df17} we replace $\pmb{\M}$ by $ \pmb{\M}'$ we obtain $R$-modules
$M^{\prime r}(w)$.

Finally let us we fix once and for all a rational number $\rho\in \mathbb Q$ which we refer to as the {\it slope} of the theory and consider the ideal of $W$ defined by
 $$W^{(\rho)}:=\{w\in W|\rho\cdot \deg(w)\in \mathbb Z\}.$$
 So if $\rho=a/b$ with $a,b$ relatively prime integers with $b\in\mathbb N$, then $W^{(\rho)}=W^{(\frac{1}{b})}$ consists of all $w\in W$ whose degree is divisible by $b$. In all our applications we will have $\rho\in \frac{1}{2}\mathbb Z$ so $W^{(\rho)}$ will be either $W$ or $W^{(\frac{1}{2})}$ according to $\rho$ being an integer or not.

 \begin{df}\label{df18}
 Let $w\in W(r)\cap W^{(\rho)}$. We say that an element $f\in M^r(w)$
 is {\it covariant} if for each
$S \in\Ob(\Prol)$ and every 
$x'\in \pmb{\M}'(S^0)$
the following equality holds in $S^r$:
$$
f(\pi_2(x'),S)= \mathfrak D(x')^{\rho \cdot \deg(w)} \cdot
f(\pi_1(x'),S).$$
\end{df}

Here  $\mathfrak D$ denotes the map
 $\mathfrak D(S^0):\pmb{\M}'(S^0)\rightarrow (S^0)^{\times}$; see see point (2) after Definition \ref{df15}.

For $w\in W(r)\cap W^{(\rho)}$ we denote by $\mathbb I^r(w)$ the $R$-submodule of covariant elements in $M^r(w)$. 
 If $w\in W\setminus (W(r)\cap W^{(\rho)})$ we set $\mathbb I^r(w):=0$. 
 Similarly, for $w\in W(r)\cap W^{(\rho)}$ we set $\mathbb M^r(w):=M^r(w)$ 
 and $\mathbb M^{\prime r}(w):=M^{\prime r}(w)$;
 for 
 $w\in W\setminus (W(r)\cap W^{(\rho)})$ 
 we set $\mathbb M^r(w)=\mathbb M^{\prime r}(w):=0$.

 For each integer $r\geq 0$ we consider the $W$-graded rings
 $$\mathbb I^r:=\bigoplus_{w\in W} \mathbb I^r(w), \ \ \mathbb M^r:=\bigoplus_{w\in W} \mathbb M^r(w),\ \ \mathbb M^{\prime r}:=\bigoplus_{w\in W} \mathbb M^{\prime r}(w).$$
 These rings form filtered graded systems so the 
 rings
 $$\mathbb I:=\bigcup_{r\geq 0} \mathbb I^r,\ \ 
 \mathbb M:=\bigcup_{r\geq 0} \mathbb M^r,\ \ \mathbb M':=\bigcup_{r\geq 0} 
 \mathbb M^{\prime r}$$
 have natural structures of filtered graded rings with graded pieces of weight $w$ given respectively by
 $$\mathbb I(w):=\bigcup_{r\geq 0} \mathbb I^r(w),\ \ 
 \mathbb M(w):=\bigcup_{r\geq 0} \mathbb M^r(w),\ \ \mathbb M'(w):=\bigcup_{r\geq 0} 
 \mathbb M^{\prime r}(w).$$
 These graded pieces are $0$ for $w\in W\setminus W^{(\rho)}$.
 The rings 
 $\mathbb I, \mathbb M, \mathbb M'$, equipped with the endomorphisms $\phi$ induced by the endomorphism $\phi$ on $M$, have natural structures of filtered $\delta$-graded rings, see Definitions \ref{df11} and \ref{df13}. 
 Note that for all $w\in W$ the $R$-linear map
 \begin{equation}\label{EQ012.111}
 \mathbb I(w)\rightarrow M
 \end{equation}
 has a torsion free cokernel.
 Also if we define the $R$-algebra automorphism
 \begin{equation}
 \label{EQ012.1111}
 \beta:\mathbb M'\rightarrow \mathbb M'\end{equation}
 by the formula
 $$(\beta f)(x',S):=\mathfrak D(x')^{\rho\cdot \deg(w)}\cdot f(x',S),\ \ f\in \mathbb M'(w),$$
 then an element $f\in \mathbb M$ belongs to $\mathbb I$ if and only if we have
 \begin{equation}\label{EQ019}\beta \pi_1^* f=\pi_2^*f.\end{equation}
 The natural $R$-algebra homomorphisms 
 \begin{equation}\label{EQ012.2}
 \mathbb I\rightarrow \mathbb M\stackrel{ \beta \pi_1^*, \pi_2^*}{\rightrightarrows} \mathbb M^{\prime}\end{equation}
 are trivially seen to be injective, with injective reductions modulo $p$. Also, $\mathbb I$ is the equalizer in the category of $R$-algebras (and also in the category of filtered $\delta$-graded rings over $R$) of the pair of $R$-algebra homomorphisms $(\beta \pi_1^*, \pi_2^*)$.

\smallskip
 
 \begin{rem}\label{R6.5}
 Assume, in this remark, that the filtered $\delta$-graded rings $\mathbb I, \mathbb M, \mathbb M'$ are integral.
 (This condition is satisfied 
 in our main applications, i.e., for the rings $\mathbb I,\mathbb M, \mathbb M'$ to be considered in Subsubsection \ref{S232}, see Corollary \ref{C4} (d).) 
For each linear system $\f=(f_0,\ldots,f_n)$ in $\mathbb I$
we get connected filtered ssa-$\delta$-spaces (see Equation (\ref{EQ004}))
$$\Proj_{\delta}(\mathbb I,\f):=\Spec_{\delta}\left(\prod_{0\leq i_0<\cdots< i_{\bullet}\leq n}
 \mathbb I_{\langle f_{i_0}\cdots f_{i_{\bullet}}\rangle}\right),$$
$$\Proj_{\delta}(\mathbb M,\f):=\Spec_{\delta}\left(\prod_{0\leq i_0<\cdots< i_{\bullet}\leq n}
 \mathbb M_{\langle f_{i_0}\cdots f_{i_{\bullet}}\rangle}\right),$$ 
$$\Proj_{\delta}(\mathbb M',\f):=\Spec_{\delta}\left(\prod_{0\leq i_0<\cdots< i_{\bullet}\leq n}
 \mathbb M'_{\langle f_{i_0}\cdots f_{i_{\bullet}}\rangle}\right),$$
 and we have naturally induced morphisms 
 of filtered ssa-$\delta$-spaces
 \begin{equation}\label{EQ029}
 \Proj_{\delta}(\mathbb M',\f)\stackrel{\pi_1,\pi_2}{\rightrightarrows} \Proj_{\delta}(\mathbb M,\f)\stackrel{\pi}{\rightarrow} \Proj_{\delta}(\mathbb I,\f).
 \end{equation}
 Here the elements $f_i\in \mathbb I\subset \mathbb M$ are viewed as elements of $\mathbb M'$ by identifying them 
  with the elements $\pi_1^*f_i\in \mathbb M'$; the morphism $\pi_2$ is then well-defined 
  in view of the display of Definition \ref{df18}.
 Note that the automorphism $\beta$ of $\mathbb M'$ induces the identity on
 $\Proj_{\delta}(\mathbb M',\f)$. 
 It is then easy to see that  (\ref{EQ029}) is a coequalizer diagram
 in $\textup{\bf ssAff}_{\delta}$.
 
 Assume, in addition, that
 for every integer $r\geq 0$, $w\in W(r)$, and $f\in \mathbb I^r(w)\setminus p \mathbb I^r(w)$ we have
$$(\mathbb I_{\langle f\rangle})^r=(\mathbb M_{\langle f\rangle})^r\cap \mathbb I_{\langle f\rangle}.$$
 (This condition is, again, satisfied 
 in our main applications, i.e., for the rings $\mathbb I,\mathbb M, \mathbb M'$ to be considered in Subsubsection \ref{S232}, see Corollary \ref{C5} (d).) Then
  (\ref{EQ029}) is also a coequalizer diagram 
 in $\textup{\bf FilssAff}_{\delta}$.
 The pair $(\pi_1,\pi_2)$ in Equation (\ref{EQ029}) may be viewed as an $\f$-realization, in $\delta$-geometry, of the pair of morphisms in Equation (\ref{EQ012.1}) while the morphism
 $\pi$ in Equation (\ref{EQ029}) may be viewed as an $\f$-realization of the morphism in Equation (\ref{EQ012.11}).
 \end{rem}

 \begin{rem}\label{R6.75}
 The Hecke correspondences on the stack $\mathcal A_{g,R}$ mentioned in Subsections  \ref{S11} and \ref{S12} that will be studied in this paper fit, as we shall see, into the above formalism. On the other hand our general formalism applies to Hecke correspondences on stacks defined by more general Shimura varieties (see \cite{Bu03, Bu04} and \cite{Bu05}, Ch. 8  for the case of Shimura curves) and also to correspondences arising from certain dynamical systems
 on projective spaces or abelian varieties (see \cite{Bu05}, Chs. 6 and 7 for the cases of $\mathbb P^1$ and respectively of elliptic curves). Here are two illustrations of situations coming from dynamical systems.
 
 \smallskip
 
 {\bf (a)} Let $1\in\Gamma\subset \pmb{\GL}_2(\mathbb Z_{(p)})$ be a submonoid. For each $R$-algebra $S$ one can consider the category $\mathfrak M_{\mathbb P^1,\Gamma}(S)$ with set of objects
 $$\pmb{\M}_{\mathbb P^1,\Gamma}(S):=\left\{v=\left( \begin{array}{c} v_1\\ v_2\end{array}\right),\ v_1,v_2\in S,\ v_1S+v_2S=S\right\},$$
 set of morphisms
 $$\pmb{\M}'_{\mathbb P^1,\Gamma}(S):=\Gamma \times \pmb{\M}_{\mathbb P^1,\Gamma}(S),$$
 source map given by the second projection, target map given by multiplication of matrices with column vectors, 
 composition given by
 $$\mu((g,v),(g',gv)):=(g'g,v)$$
 and unit defined by $\epsilon(v):=(1_2,v)$. Furthermore, we consider the group $G(S)=\mathbb G_m(S)=S^{\times}$ and the actions of $G(S)$ 
 on $\pmb{\M}_{\mathbb P^1,\Gamma}(S)$ and $\pmb{\M}'_{\mathbb P^1,\Gamma}(S)$
 given by usual scalar multiplication of elements of $S^2$ by units of $S$. We take the homomorphism $D$ to be the identity and define $\mathfrak D(g,v):=\det(g)$.
 We have defined a moduli datum $(\mathfrak M_{\mathbb P^1,\Gamma},\mathfrak D,\mathbb G_m,D)$. 
 Note that if $\Gamma$ has no invertible elements except the identity, then for every local ring $S$ over $R$
 we have $\mathcal A_{\mathfrak M_{\mathbb P^1,\Gamma},\mathbb G_m}(S)=\mathbb P^1(S)$.
 Take the slope to be
 $\rho=\frac{1}{2}$.
 Then the theory in \cite{Bu05}, Ch. 6 provides an analysis of the ring $\mathbb I=\mathbb I_{\mathbb P^1,\Gamma}$ in this case. In particular, if the group generated by $\Gamma$ is not solvable, then $\mathbb I$ is described in terms of certain explicit forms $f^r\in \mathbb I^r(-\phi^r-1)$, $r\geq 1$ (see \cite{Bu05}, Prop. 6.17).
 This example (which, as explained in \cite{Bu05}, Ch. 6 is directly related to the classical invariant theory of $\pmb{\SL}_2$ acting on multiple pairs of $2$-vectors) will play no role in the present paper.
 
 \smallskip
 
 {\bf (b)} Fix an elliptic curve $E$ over $R$. For each $R$-algebra $S$ we consider the category $\mathfrak M_E(S)$ with set of objects
 $$\pmb{\M}_E(S):=\{(P,\omega)|P\in E(S),\ S\omega=H^0(E_S,\Omega_{E_S/S})\},$$
 set of morphisms
 $$\pmb{\M}'_E(S):=\{(P_1,P_2,\omega,n)|(P_1,\omega),(P_2,\omega)\in \pmb{\M}(S), \ n\in \mathbb N\setminus p\mathbb N,\ 
 P_2=nP_1\},$$
 source and target given by the obvious projections, 
 composition given by
 $$\mu((P_1,P_2,\omega,n),(P_2,P_3,\omega,m)):=(P_1,P_3,\omega,nm),$$
 and unit defined by $\epsilon(P,\omega):=(P,P,\omega,1)$. Furthermore consider the group $G(S)=\mathbb G_m(S)=S^{\times}$ and the actions of $G(S)$ given by multiplication of $\omega$ by elements of $S^{\times}$. Finally consider the homomorphism $D$ given by $D(\lambda):=\lambda^{-1}$ and the last projection map $\mathfrak D(P_1,P_2,\omega,n):=n$. 
 We have defined a moduli datum $(\mathfrak M_E,\mathfrak D,\mathbb G_m,D)$. 
 Note that for every local ring $S$ over $R$ we have $\mathcal A_{\mathfrak M_E,\mathbb G_m}(S)=E(S)$.
 Take the slope to be
 $\rho=1$.
 Then the theory in \cite{Bu05}, Ch. 7 provides an analysis of the ring $\mathbb I=\mathbb I_E$ in this case. In particular, if $p\neq 2$, $E$ is not a canonical lift, and 
 $\psi_{\omega}$ is the {\it normalized $\delta$-character} of $E$ over $S$ attached to the $1$-form $\omega$ 
 (see \cite{Bu95a}, Introd. or \cite{Bu05}, Thm. 7.22, p. 197, for the case $S=R$, and similarly for a general filtered $\delta$-ring $S$), 
 then the rule $f(P,\omega):=\psi_{\omega}(P)$ defines an element $f\in \mathbb I^2(\phi^2)$;
 moreover all the elements in the ring $\mathbb I$ can be `constructed from' this element, see
 \cite{Bu05}, Thm. 7.34.
 This example plays a key role in a number of contexts \cite{Bu95a, Bu05, BM, BP} but will play no role in the present paper.
\end{rem}

\subsubsection{Rings of Siegel $\delta$-modular forms}\label{S232}

The dual of an abelian scheme $A$
over a scheme will be denoted by $\check{A}$.

An abelian scheme $A$ over $\Spec(S)$ is called {\it ordinary} if all its geometric fibers over points of $\Spec(\overline{S})$ are ordinary abelian varieties over fields. If $S=k$, let $T_p$ be the Tate module functor from the category of ordinary abelian varieties over $S$ to the category of finitely generated $\mathbb Z_p$-modules; so $T_p(A)$ is the Tate module of $A$. 

Let $\mathcal A_{g,R}$ be the stack of principally polarized abelian schemes of relative dimension $g$ over $R$-schemes (for instance, see \cite{FC}, Ch. I, Def. 4.3, p. 19). This is an algebraic stack over $\Spec(R)$ in the sense of \cite{FC}, Ch. I, Def. 4.6, p. 20 (see \cite{FC}, Ch. I, Subsect. 4.11)\footnote{In the more modern terminology, it is a Deligne--Mumford stack.}.
In particular it has an \'{e}tale cover that is a smooth scheme over $\Spec(R)$ and hence it has a well-defined relative dimension which equals to $\frac{g(g+1)}{2}$.
The stack $\mathcal A_{g,R}$ induces a functor we still denote by $\mathcal A_{g,R}$ from the category
 $\textup{\bf Ring}$ of $R$-algebras to the category $\textup{\bf Set}$ of sets: for each $R$-algebra $S$, $\mathcal A_{g,R}(S)$ is the set
of isomorphism classes of 
 principally polarized abelian schemes over $S$ of relative dimension $g$.  
 
Let $N\geq 3$ be an integer prime to $p$. The functor $\mathcal A_{g,1,N,R}$ from 
 $\textup{\bf Ring}$ to $\textup{\bf Set}$
 that attaches to each $R$-algebra $S$  the set $\mathcal A_{g,1,N,R}(S)$ of isomorphism
 classes of principally polarized abelian schemes over $S$ of relative dimension $g$ equipped with symplectic similitude level-$N$ structure 
 is represented by (and will be identified with) a smooth $R$-scheme which we denote simply by $\mathcal X:=\mathcal A_{g,1,N,R}$ (see the second footnote in Subsection \ref{S11}). 
 
 If, in the above discussion, one takes ordinary abelian schemes  instead of abelian schemes, then one obtains a functor
 $\mathcal A_{g,\ord,R}$
 from 
 $\textup{\bf Ring}$ to $\textup{\bf Set}$.
 Similarly, for $N\geq 3$ prime to $p$, we may consider the open subscheme 
 $\mathcal X_{\ord}=\mathcal A_{g,1,N,\ord,R}$ of $\mathcal X$ which is the union of $\mathcal X_K=\mathcal A_{g,1,N,K}$ and of the ordinary locus $\mathcal X_{\ord,k}$ of $\mathcal X_k=\mathcal A_{g,1,N,k}$. 

In what follows we define two moduli data 
$$(\mathfrak M_g,\mathfrak D,\pmb{\GL}_g,\det)\ \ \ \text{and}\ \ \ (\mathfrak M_{g,\ord},\mathfrak D,\pmb{\GL}_g,\det)$$
 and we apply to them the constructions of the previous subsubsection. 
The rings $\mathbb I, \mathbb M, \mathbb M'$ that correspond to these moduli data will be called 
$$\mathbb I_g, \mathbb M_g, \mathbb M'_g\ \ \ \text{and}\ \ \ \mathbb I_{g,\ord}, \mathbb M_{g,\ord}, \mathbb M'_{g,\ord},$$ 
respectively; however, if no confusion arises, we will sometimes drop the indices.

Recall that  our ground ring is $S_0=R$.
For each $R$-algebra
 $S$, let $\pmb{\M}_g(S)$ denote the set of all triples
$(A,\theta,\omega)$ where $A$ is an abelian scheme
of relative dimension $g$ over $\Spec(S)$, where $\theta:A \rightarrow
\check{A}$ is a principal polarization, and where $\omega$ is a column basis of $1$-forms on $A$ by which we mean a column vector $\omega=[\omega_1\cdots\omega_g]^{\t}$ whose entries form a basis of the $S$-module of $1$-forms $H^0(A,\Omega^1_{A/\Spec(S)})$. 

Let now
$$(A_1,\theta_1,\omega_1), (A_2,\theta_2,\omega_2)
\in \pmb{\M}_g(S)$$
 and let $u:A_1 \rightarrow A_2$ be an isogeny
(which is not assumed to be compatible
with either the principal polarizations or the column bases of $1$-forms).
Let 
$$u^{\t}:=\theta_1^{-1} \circ \check{u} \circ 
\theta_2:
A_2 \rightarrow \check{A}_2 \rightarrow \check{A}_1
 \rightarrow A_1,$$
 $$
 u^*\theta_2:=\check{u}\circ \theta_2\circ u:A_1\rightarrow A_2
 \rightarrow \check{A}_2 \rightarrow \check{A}_1.$$
 We say $u$ is a {\it polarized isogeny of degree prime to $p$
  with respect to $\theta_1$ and $\theta_2$} if there exists
 $d\in\mathbb N\setminus p\mathbb N$ such that $uu^{\t}=d:A_2\rightarrow A_2$ 
 (equivalently, $u^*\theta_2=d\theta_1:A_1\rightarrow \check{A}_1$); for such $u$ we have $\deg(u)=d^g$. If $\theta_1$ and $\theta_2$ are clear from context and  $p$ is fixed  we will not mention `of degree prime to $p$' or `with respect to $\theta_1$ and $\theta_2$' and simply say that $u$ is a {\it polarized isogeny}.
 
 We denote by $\pmb{\M}'_g(S)$ the set of all triples $(x_1,x_2,u)$ where 
 $$x_i:=(A_i,\theta_i,\omega_i)\in \pmb{\M}_g(S),\ \ i\in\{1,2\},$$
 and $u:A_1\rightarrow A_2$ is a polarized isogeny  such that $u^*\omega_2=\omega_1$.
 There exist canonical maps of sets
 $$\begin{array}{rcc}
 \pi_1,\pi_2:\pmb{\M}'_g(S) & \rightarrow & \pmb{\M}_g(S),\\
 \epsilon:\pmb{\M}_g(S) & \rightarrow & \pmb{\M}'_g(S), \\
 \mu:\pmb{\M}''_g(S) & \rightarrow & \pmb{\M}'_g(S),\end{array}$$
 where $ \pi_i(x_1,x_2,u):=x_i$, $\epsilon(x):=(x,x,\textup{id})$, 
 $\pmb{\M}''_g(S)$ is the set of all pairs $(y_1,y_2)\in \pmb{\M}'_g(S)^2$ such that 
 $\pi_2(y_1)=\pi_1(y_2)$,
 and
 $\mu$ is given by the obvious composition of isogenies. 
 We obtain a category $\mathfrak M(S)$ whose set of 
 objects is $\pmb{\M}_g(S)$, whose set of morphisms is $\pmb{\M}'_g(S)$, and where
 $\pi_1,\pi_2,\epsilon,\mu$ are the source, target, identity, and composition maps.

 We let the group $\pmb{\GL}_g(S)$ act on $\pmb{\M}_g(S)$ by the rule
$$\lambda \cdot(A,\theta,\omega):=(A,\theta,\lambda \omega)$$
for $\lambda\in \pmb{\GL}_g(S)$ and $(A,\theta,\omega)\in \pmb{\M}_g(S)$
and act on $\pmb{\M}'_g(S)$ by the rule
$$\lambda\cdot (x_1,x_2,u):= (\lambda \cdot x_1, \lambda \cdot x_2, u).$$
We take the group homomorphism $D=\det:\pmb{\GL}_g(S)\rightarrow S^{\times}$, $D(\lambda)=\det(\lambda)$.

Finally we denote by $\mathfrak D:\pmb{\M}'_g(S)\rightarrow S^{\times}$ the map defined by letting
$\mathfrak D(x_1,x_2,u)$ be the image of $d:=uu^{\t}\in \mathbb N\setminus p\mathbb N$ in $S^{\times}$. 

We have defined above a moduli datum $(\mathfrak M_g, \mathfrak D,\pmb{\GL}_g, \det)$. 

Clearly, for every local ring $S$, the set $\mathcal A_{\mathfrak M_g, \pmb{\GL}_g}(S)$ attached to this moduli datum (see Remark \ref{R6.1})
 identifies with the set $\mathcal A_{g,R}(S)$ of isomorphism classes
 of principally polarized abelian schemes $(A,\theta)$ over $S$. Similarly, the set 
 $\mathcal A'_{\mathfrak M_g, \pmb{\GL}_g}(S)$
 (see Remark \ref{R6.1})
 identifies with the set $\mathcal A'_{g,R}(S)$ of isomorphism classes of quintuples $(A_1,\theta_1,A_2,\theta_2,u)$
 where $(A_1,\theta_1),(A_2,\theta_2)\in \mathcal A_{g,R}(S)$ and $u:A_1\rightarrow A_2$ is a polarized isogeny; the set $\mathcal A'_{g,R}(S)$ may be viewed as the `union of all prime to $p$ Hecke correspondences' on $\mathcal A_{g,R}(S)$. We continue to denote by
 $\mathcal A'_{g,R}$ the obvious stack inducing the functor $\mathcal A'_{g,R}$.

In what follows, for the value of our slope, we take $\rho:=g/2$; as we shall see, this value is `forced upon us' by the compatibilities built into our specific situation (see Remark \ref{R6.95} (c)). 

By Subsubsection \ref{S231}, by adding the index $g$ we can attach to our moduli datum and to our slope $\rho=g/2$ the corresponding $R$-modules 
\begin{equation}
\label{EQ029.1}
M^r_g(w), M^{\prime r}_g(w), 
 \mathbb M^r_g(w), \mathbb M^{\prime r}_g(w),\mathbb I^r_g(w),\end{equation}
and $R$-algebras
\begin{equation}
\label{EQ029.2}
M^r_g, M_g, M^{\prime r}_g, M'_g, \mathbb M^r_g, 
\mathbb M^{\prime r}_g, 
\mathbb I^r_g, \mathbb M_g, \mathbb M'_g, \mathbb I_g.\end{equation}

 \begin{df}\label{df19}
  The elements of $M^r_g$ are called {\it Siegel $\delta$-modular functions} of genus $g$ and order $\leq r$. 
The elements of $M^r_g(w)$ are called {\it Siegel $\delta$-modular forms} of genus $g$, order $\leq r$, and weight $w$.
 A Siegel $\delta$-modular form $f \in \mathbb M^r_g(w)$ is called {\it Hecke covariant} if it is covariant, i.e., 
 if it belongs to $\mathbb I^r_g(w)$.
 \end{df} 

\begin{rem}\label{R6.85}
{\bf (a)} Corollary \ref{C4} (d) below will prove that the rings $\mathbb M_g, \mathbb M'_g, \mathbb I_g$ are filtered integral $\delta$-graded rings. So the discussion in Remark \ref{R6.75} applies to these rings.

{\bf (b)} 
 We will later on prove that we have $M^s_g(w)\cap M^r_g=0$ provided $0\leq r<\ord(w)\leq s$, see Corollary \ref{C5} (a) below. This may be seen as an a posteriori justification for our convention that
 $\mathbb M^r(w)=0$ for $w\in W\setminus W(r)$.

 {\bf (c)}
 We will later on prove that the $R$-linear maps
$\mathbb M^r_g\rightarrow M^r_g$
are injective, see Corollary \ref{C5} (b) below. On the other hand, for each $w\in W$ the cokernels of the maps
$M^r_g(w)\rightarrow M^r_g$ and $\mathbb I^r_g(w)\rightarrow M^r_g$
are trivially seen to be torsion free.

{\bf (d)} In view of the paragraph after  Definition \ref{df18}, by taking $\rho=g/2$, from very definitions we get that we have $\mathbb I^r_g(w)=\mathbb M^r_g(w)=0$ for $\frac{g\cdot \deg(w)}{2}\not\in \mathbb Z$.
\end{rem}

\smallskip

\begin{rem}\label{R6.95}
{\bf (a)} An element $f\in M^r_g$ is in $M^r_g(w)$ if and only if
 for each 
 $S\in \textup{Ob}(\Prol)$ and every
 $\lambda \in\pmb{\GL}_g(S^0)$ the following equality holds in $S^r$:
\begin{equation}\label{EQ029.3}
f(A, \theta, \lambda \omega, S)= \det(\lambda)^{-w} \cdot
f(A,\theta, \omega,S).\end{equation}
Note that if a non-zero element  $f\in M^r_g$ has a weight then the weight is unique. This 
follows from the fact that if  $w\in W$ is such that $\lambda^w=1$ for all $\lambda\in R^{\times}$ then $w=0$; the latter property follows from the injectivity of the map in Equation (\ref{EQ003}) applied to $X=\Spec(R[x,x^{-1}])$
and the element $x^w-1\in \mathcal O(J^{\ord(w)}(X))$.

 {\bf (b)}  An element $f \in \mathbb M^r_g(w)$ is in $ \mathbb I^r_g(w)$
 if and only if for each
$S \in\Ob(\Prol)$ and every 
$(x_1,x_2,u)\in \pmb{\M}'_g(S^0)$, where
$x_1=(A_1,\theta_1,\omega_1)$, $x_2=(A_2,\theta_2,\omega_2)$, $uu^{\t}=d\in\mathbb N\setminus p\mathbb N$, and $u^*\omega_2=\omega_1$, the following equality holds 
in $S^r$:
\begin{equation}\label{EQ029.4}
f(A_2,\theta_2,\omega_2,S)= d^{\frac{g\cdot \deg(w)}{2}} \cdot
f(A_1, \theta_1, \omega_1,S).\end{equation}
If $f\in\mathbb I^r_g(w)$, then by combining Equations (\ref{EQ029.3}) and (\ref{EQ029.4}) we get that 
 for every
$S \in\Ob(\Prol)$, every 
$(A_1,\theta_1,\omega_1),(A_2,\theta_2,\omega_2)\in \pmb{\M}_g(S^0)$, 
and every polarized isogeny $u:A_1\rightarrow A_2$ with
$uu^{\t}=d\in\mathbb N\setminus p\mathbb N$, the following equality holds in $S^r$:
\begin{equation}\label{EQ029.5}
f(A_2,\theta_2,\omega_2,S)= d^{\frac{g\cdot \deg(w)}{2}} \cdot \det([u])^{-w}\cdot 
f(A_1, \theta_1, \omega_1,S),\end{equation}
where $[u]\in\pmb{\GL}_g(S^0)$ is the unique matrix such that $u^*\omega_2=[u]\omega_1$.


{\bf (c)} Assume $0\neq f \in M^r_g(w)$ and $\nu\in \mathbb Z$ are such that for each 
$S \in\Ob(\Prol)$ and every 
$(x_1,x_2,u)\in \pmb{\M}'_g(S^0)$, where
$x_1=(A_1,\theta_1,\omega_1)$, $x_2=(A_2,\theta_2,\omega_2)$, $uu^{\t}=d\in\mathbb N\setminus p\mathbb N$, and $u^*\omega_2=\omega_1$, we have
\begin{equation}
\label{EQ029.6}
f(A_2,\theta_2,\omega_2,S)= d^{\nu} \cdot
f(A_1, \theta_1, \omega_1,S).\end{equation}
Then it is easy to check (using Equation (\ref{EQ029.6}) for both $u$ and $u^{\t}$ plus Equation (\ref{EQ029.3})) that
$\frac{g\cdot \deg(w)}{2}\in \mathbb Z$ and $\nu=\frac{g\cdot \deg(w)}{2}$.
\end{rem}

It turns out that, in order to study the ring $\mathbb I_g$, it is crucial to develop
 a matrix version of some of our definitions. We present this in what follows.

\begin{df}\label{df20}
Let $m\in\mathbb N$ and $r\in\mathbb N\cup\{0\}$. 
A {\it Siegel $\delta$-modular function}
 of genus $g$, size $m$, and order $\leq r$
is a rule $f$ that associates to each pair
 $(S,(A,\theta,\omega))\in \Ob(\Prol)\times\pmb{\M}_g(S^0)$
an $m \times m$ matrix
$$f(A,\theta,\omega,S) \in\pmb{\Mat}_m(S^r)$$
that depends on $S$ and the isomorphism class
of $(A,\theta,\omega)$ only and that is functorial
in $S$. \end{df}

In this paper we will use only the sizes $m=1$ and $m=g$.
The set of Siegel $\delta$-modular functions
 of genus $g$, size $1$, and order $\leq r$ coincides with $M^r_g$. 
We denote by $M^r_{gg}$ the $R$-algebra of Siegel $\delta$-modular functions
 of genus $g$, size $g$, and order $\leq r$. We have
 $M^r_{gg}=\pmb{\Mat}_g(M^r_g)$. Similarly to the case of size $1$ we set
 $M_{gg}:=\bigcup_{r\geq 0} M^r_{gg}=\pmb{\Mat}_g(M_g)$.
 
 In what follows we discuss the concepts of {\it weight} and {\it Hecke covariance} for Siegel $\delta$-modular functions of size $g$.
 Let $$\mathcal W:=\{\pm\phi^a|a\in \mathbb Z_{\geq 0}\}\subset W$$
 be the multiplicative submonoid generated by 
$-1$ and $\phi$. For an element $w\in\{\pm\phi^a\}\subset\mathcal W$ and $S\in\Ob(\Prol)$, let $\chi_w:\pmb{\GL}_m(S^0)\rightarrow \pmb{\GL}_m(S^a)$ be the homomorphism which maps $\lambda\in\pmb{\GL}_m(S^0)$ to
$$\begin{array}{rcll}
\chi_w(\lambda):=\phi^a(\lambda) & \textup{if}\ \ \ w=\phi^a,\\
\ & \ & \ & \ \\
\chi_w(\lambda):=((\phi^a(\lambda))^{\t})^{-1} &\textup{if}\ \ \ w=-\phi^a.\end{array}$$

\begin{df}\label{df21}
 Fix a pair $(w',w'')\in\mathcal W\times \mathcal W$
where $w',w''$ have order $\leq r$; so
$w',w''\in\{\pm\phi^a|a\in\{0,\ldots,r\}\}$.
A Siegel $\delta$-modular function $f\in M^r_{gg}$ is called
a {\it Siegel $\delta$-modular form} of weight $(w',w'')$
if for every $S \in\Ob(\Prol)$, for every $(A,\theta,\omega) \in \pmb{\M}_g(S^0)$
and every $\lambda \in\pmb{\GL}_g(S^0)$ the following equality holds in
$ \pmb{\Mat}_g(S^r)$:
\begin{equation}\label{EQ029.7}
f(A,\theta,\lambda \omega,S)= \chi_{w'}(\lambda) \cdot
f(A,\theta,\omega,S) \cdot \chi_{w''}(\lambda^{\t}).\end{equation}\end{df}

The $R$-module of all Siegel $\delta$-modular forms of genus $g$, size $g$,
order $\leq r$ and weight $(w',w'')$  will be denoted by $M^r_{gg}(w',w'')$.
Note that if a non-zero element of $M^r_{gg}$ has a weight $(w',w'')$ then that weight is not a priori unique (it is clearly never unique if $g=1$) but one expects that $w'+w''=0$.

\begin{df}\label{df22}
Fix a pair $(w',w'')\in\mathcal W\times \mathcal W$
where $w',w''$ have order $\leq r$; so
$w',w''\in\{\pm\phi^a|a\in\{0,\ldots,r\}\}$.
A Siegel $\delta$-modular form $f\in M^r_{gg}(w',w'')$ is {\it Hecke covariant} if 
 for each
$S \in\Ob(\Prol)$ and every 
$(x_1,x_2,u)\in \pmb{\M}'_g(S^0)$, where
$x_1=(A_1,\theta_1,\omega_1)$, $x_2=(A_2,\theta_2,\omega_2)$, $uu^{\t}=d\in\mathbb N\setminus p\mathbb N\subset \End(A_2)$,
the following equality holds,
\begin{equation}\label{EQ029.8}
f(A_2,\theta_2,\omega_2,S)= d^{-\deg(w'+w'')/2}\cdot 
f(A_1, \theta_1, \omega_1,S).
\end{equation}\end{df}

As $w',w''\in\mathcal W$ have degree $-1$ or $1$, we have $\deg(w'+w'')/2\in \{-1,0,1\}$.
We denote by $\mathbb I^r_{gg}(w',w'')$ the $R$-submodule of $M^r_{gg}(w',w'')$ of Hecke covariant Siegel $\delta$-modular forms. 

Let $f\in\mathbb I^r_{gg}(w',w'')$. Combining
Equations (\ref{EQ029.7}) and (\ref{EQ029.8}) we get that
 for every
$S \in\Ob(\Prol)$, every 
$(A_1,\theta_1,\omega_1), (A_2,\theta_2,\omega_2)\in \pmb{\M}_g(S^0)$, and every polarized isogeny $u:A_1\rightarrow A_2$ with $uu^{\t}=d\in\mathbb N\setminus p\mathbb N$
the following 
holds:
\begin{equation}\label{EQ029.9}
f(A_2,\theta_2,\omega_2,S)= d^{-\deg(w'+w'')/2}\cdot \chi_{w'}([u])\cdot 
f(A_1, \theta_1, \omega_1,S)\cdot \chi_{w''}([u]^{\t}),
\end{equation}
where $[u]\in\pmb{\GL}_g(S^0)$ is the unique matrix such that $u^*\omega_2=[u]\omega_1$.

\medskip

\begin{rem}\label{R7}
Our $R$-modules $\mathbb I^r_{gg}(\phi^a,\phi^b)$ a priori contain the $R$-modules 
$$I^r_g(\phi^a,\phi^b)\otimes_{\mathbb Z_p} R$$ defined by the $\mathbb Z_p$-modules of {\it isogeny covariant} forms $I^r_{g}(\phi^a,\phi^b)$ introduced in \cite{BB}, Subsect. 1.3. The latter spaces will play no role in our paper but note that our results here for $\mathbb I^r_{gg}(\phi^a,\phi^b)$ imply (and provide new proofs for) the corresponding results in \cite{BB} for the $\mathbb Z_p$-modules $I^r_g(\phi^a,\phi^b)$. We recall, for convenience, that the definition of the spaces 
 $I_g(\phi^a,\phi^b)$ in \cite{BB} 
 was similar to that of our spaces
$\mathbb I_{gg}(\phi^a,\phi^b)$ here, except that, for the former, one considered `{\it prolongation sequences} over $\mathbb Z_p$' rather than `filtered $\delta$-rings over $R$' and, more importantly, condition (\ref{EQ029.9}) was required to hold for all isogenies $u$ 
of degree prime to $p$, rather than just for polarized isogenies. \end{rem}

\medskip

\begin{rem}\label{R7.2}
If in  Definitions \ref{df19} to \ref{df22} we replace 
 $\pmb{\M}_g(S)$ by the set $\pmb{\M}_{g,\ord}(S)$ of all triples 
 $(A,\theta,\omega)\in \pmb{\M}_g(S)$
 with the property that $A$ is ordinary and if we replace $\pmb{\M}'_g(S)$ by the corresponding set $\pmb{\M}'_{g,\ord}(S)$, we obtain a new moduli datum $(\mathfrak M_{g,\ord},\mathfrak D,\pmb{\GL}_g,\det)$. Correspondingly, the objects $\square_g$ with the index $g$ attached to $(\mathfrak M_g,\mathfrak D,\pmb{\GL}_g,\det)$ will be replaced by objects $\square_{g,\ord}$ attached to 
 $(\mathfrak M_{g,\ord},\mathfrak D,\pmb{\GL}_g,\det)$ and will be referred to as {\it `ordinary'}. 
 For instance $M^r_g(w)$ will be replaced by $M^r_{g,\ord}(w)$ and the elements of the latter will be referred to as {\it ordinary Siegel $\delta$-modular forms} of genus $g$, order $\leq r$, and weight $w$, etc.
 All claims in Remarks \ref{R6.85} and \ref{R6.95} hold in this new context. 
 Note that here the word `ordinary' is not an adjective: each Siegel $\delta$-modular form
 defines an ordinary Siegel $\delta$-modular form, and in fact we have natural $R$-linear homomorphisms
 $$
 M^r_g(w) \rightarrow M^r_{g,\ord}(w),\ \ M^r_{gg}(w',w'') \rightarrow M^r_{gg,\ord}(w',w'')$$
 which turn out to be injective by our Serre--Tate expansion principle (see Proposition \ref{P3} (b) below).
 \end{rem}

\begin{rem}\label{R7.3}
In this remark we consider the ordinary case; an entirely similar discussion holds in the non-ordinary case.
\smallskip

{\bf (i)} The $R$-modules $\mathbb I^r_{gg,\ord}(w',w'')$ cannot be naturally assembled into a ring. However,
 for $w',w'',w'''\in\mathcal W$ we have natural maps
$$
\begin{array}{rcll}
\mathbb I^r_{gg,\ord}(w',w'') & \rightarrow & \mathbb I^r_{gg,\ord}(w'',w'), & f\mapsto f^{\t}\\
\ & \ & \ & \ \\
\mathbb I^r_{gg,\ord}(w',w'') & \rightarrow & \mathbb I^{r+1}_{gg,\ord}(\phi w',\phi w''), & f\mapsto f^{\phi}:=\phi(f),\\
\ & \ & \ & \ \\
\mathbb I^r_{gg,\ord}(w',w'')\times \mathbb I^r_{gg,\ord}(-w'',w''') & \rightarrow & \mathbb I^r_{gg,\ord}(w',w'''), &(f_1,f_2)\mapsto f_1f_2,\end{array}$$
where the first map is $R$-linear, the second one is $\mathbb Z_p$-linear, and the third one is $R$-bilinear. Also, the rule $(A,\theta,\omega, S)\mapsto 1_g$ defines an element 
of $\mathbb I^r_{gg,\ord}(-w',w')$ for each $w'\in\mathcal W$. So the set $\mathcal W$ can be seen as the set of objects
of a category whose set of morphisms from $w'\in \mathcal W$ to $w''\in \mathcal W$ identifies
with the $R$-module $\mathbb I^r_{gg,\ord}(-w',w'')$ and whose composition law is given by the above $R$-bilinear maps.

\smallskip

{\bf (ii)}
 If $f\in\mathbb I^r_{gg,\ord}(w',w'')$ has the property that we have
$$f(A,\theta,\omega,S)\in {\pmb{\GL}}_g(S^r)$$
 for all $(A,\theta,\omega,S)$, then the rule
\begin{equation}\label{EQ029.91}
(A,\theta,\omega,S)\mapsto (f(A,\theta,\omega,S))^{-1}\end{equation}
defines an element 
$$f^{-1}\in\mathbb I^r_{gg,\ord}(-w'',-w').$$
Similarly, if $f\in\mathbb I^r_{g,\ord}(w)$ has the property that we have
 $$f(A,\theta,\omega,S)\in (S^r)^{\times}$$ for all $(A,\theta,\omega,S)$, then the rule (\ref{EQ029.91})
defines an element $$f^{-1}\in\mathbb I^r_{g,\ord}(-w).$$

\smallskip

{\bf (iii)} For all $w',w''\in\mathcal W$ we have natural maps 
$$
\det: \mathbb I^r_{gg,\ord}(w',w'') \rightarrow \mathbb I^r_{g,\ord}(-w'-w''),\ \ \ f\mapsto \det(f),$$
 which for $g=1$ are canonical identifications. 
For $g=1$ we have additional canonical identifications 
$M_{11,\ord}^r(w',w'')=M_{1,\ord}^r(-w'-w'')$, etc.
\end{rem}

 \subsubsection{Rings of $\delta$-homogeneous polynomials}\label{S233}
 
 For an $n+1$-tuple of indeterminates $z:=(z_0,\ldots,z_n)$, where $n\geq 1$, 
the $\delta$-ring of $\delta$-polynomials 
 $R[z,z',z'',\ldots]$ 
  (defined by Equation (\ref{EQ001}) but with $(y,d)$ replaced by $(z,n)$) contains the elements ($i\in\{0,\ldots,n\}$):
 $$z_i^{\phi}=z_i^p+pz'_i,\ \ z_i^{\phi^2}=(z_i^p+pz_i')^p+p((z'_i)^p+pz''_i),\ \textup{etc}.$$
 We define $z':=(z_0',\ldots,z'_n)$, $z^{\phi}:=(z_0^{\phi},\ldots,z_n^{\phi})$, etc.,
 and view the ring $R[z,z',z'',\ldots]$ as a subring of the 
 $K$-algebra of polynomials $K[z,z',z'',\ldots]$ in the set of indeterminates $\{z_i,z'_i,z''_i,\ldots|i\in\{0,\ldots,n\}\}$. This $K$-algebra  is generated by the set 
$\{z_i,z^{\phi}_i,z^{\phi^2}_i,\ldots|i\in\{0,\ldots,n\}\}$, i.e.,
 \begin{equation}\label{EQ030}
 K[z,z',z'',\ldots]=K[z,z^{\phi},z^{\phi^2},\ldots]:=K[z_i,z^{\phi}_i,z^{\phi^2}_i,\ldots|i\in\{0,\ldots,n\}].\end{equation}
 On this $K$-algebra we consider the  $W_+$-grading defined by letting $z_i^{\phi^j}$ have degree $\phi^j$.
 For $w\in W_+$ we denote by $K[z,z^{\phi},z^{\phi^2},\ldots](w)$ the component of degree $w$ in the ring (\ref{EQ030}) and consider the intersection
 $$\mathbb S_n(w):=R[z,z',z'',\ldots]\cap (K[z,z^{\phi},z^{\phi^2},\ldots](w))$$
 inside the ring (\ref{EQ030}).
 
\begin{df}\label{df23}
 The ring of {\it $\delta$-homogeneous polynomials} in $z$ is the 
 $W_+$-graded ring 
 \begin{equation}\label{EQ031}
 \mathbb S_n:=\bigoplus_{w\in W_+} \mathbb S_n(w).\end{equation}
\end{df}
 
 This ring has a natural structure of $R$-algebra and is equipped with a ring endomorphism $\phi$, extending the Frobenius lift  on $R$, such that for all $w\in W_+$ we have $\phi(\mathbb S_n(w))\subset \mathbb S_n(\phi(w))$. Note however that $\phi$ is not a Frobenius lift on $\mathbb S_n$  and thus this ring does not have a natural structure of $\delta$-ring; indeed we have $z'_i=\frac{1}{p}z_i^{\phi}-\frac{1}{p}z_i^p\not\in \mathbb S_n$ because the homogeneous components $\frac{1}{p}z_i^{\phi}$ and $\frac{1}{p}z_i^p$ of $z'_i$ do not belong to $\mathbb S_n$. However, if $F_0,\ldots, F_m\in \mathbb S_n(w_1)$ and $G\in \mathbb S_m(w_2)$, then $G(F_0,\ldots,F_m)\in \mathbb S_n(w_1w_2)$. In particular,
 as for all $i,j\in\{0,\ldots,n\}$ we have
 $$\{z_i,z_j\}_{\delta}=z_i^pz'_j-z_j^pz'_i\in \mathbb S_n(\phi+p),$$
 it follows that 
 for all $F,G\in \mathbb S_n(w)$ we have $\{F,G\}_{\delta}\in \mathbb S_n((\phi+p)w)$.
 So $\mathbb S_n$ has a natural structure of a filtered integral $\delta$-graded ring over $R$ with filtration  $(\mathbb S_n\cap R[z,z',\ldots,z^{(r)}])_{r\geq 0}$. 
 
 It is sometimes useful to shift the filtration and grading on $\mathbb S_n$ in a way that it depends on a given integer $r_0\geq 0$ and on a given $w_0\in W(r_0)$. More precisely, for a given integer $r_0\geq 0$ we define a  filtration $(\mathbb S^r_n)_{r\geq 0}$ by setting
 $\mathbb S^r_n:=R$ if $0\leq r<r_0$ and 
 $$\mathbb S^r_n:=\mathbb S_n\cap R[z,z',\ldots,z^{(r-r_0)}],\ \textup{if}\ \ \ r\geq r_0,$$
and for a given $w_0\in W(r_0)$ we define
 a new grading on $\mathbb S_n$ by letting each $z_i^{\phi^j}$ have weight $\phi^jw_0$. With this new filtration and grading $\mathbb S_n$ is, again, a filtered integral $\delta$-graded ring and we refer to the pair $(r_0,w_0)$ as its {\it shift}.
 In particular, for each shift $(r_0,w_0)$ (which will not be included in our notation), we view
 the connected filtered ssa-$\delta$-space
 \begin{equation}\label{EQ032}
 \mathbb P^n_{\delta}:=\Proj_{\delta}(\mathbb S_n,z)=\Spec_{\delta}\left(
 \prod_{0\leq i_0<\cdots<i_{\bullet}\leq n} \mathbb S_{n,\langle z_{i_0}\cdots z_{i_{\bullet}}\rangle}
 \right)\end{equation}
as an analogue in the $\delta$-geometry of the $n$-th dimensional projective space.
 
 We have a $K$-algebra isomorphism (identification) of $W_+$-graded algebras
 $$\mathbb S_n \otimes_R K\simeq K[z,z^{\phi},z^{\phi^2},\ldots ].$$
 
 Note that the
 inclusion $\mathbb S_n\subset R[z,z',z'',\ldots]$ is compatible with $\phi$ acting on the two $R$-algebras and
 for all $w\in W$, the $R$-linear map 
 $$\mathbb S_n(w)\rightarrow R[z,z',z'',\ldots]$$
 has a torsion free cokernel. 
 
 The triple $(\mathbb S_n, R[z,z',z'',\ldots], z)$ with  shift $(r_0,w_0)$
 has the following universal property. 
 Consider an arbitrary triple $(B,C,\bf f)$, where
 
\medskip
 
 {\bf (i)} $C$ is a $\delta$-ring over $R$, $B$ is a $R$-subalgebra of $C$ such that $\phi(B)\subset B$ and  $B$ is equipped with an integral $\delta$-graded ring structure with respect to the restriction of $\phi$ to $B$,  with the  property that for all $w\in W$, the $R$-linear map $B(w)\rightarrow C$ has a torsion free cokernel;
 
\smallskip
 
{\bf (ii)} $\f=(f_0,\ldots,f_n)$ is a linear system in $B$ of weight $w_0\in W(r_0)$.

\medskip\noindent
 Then there exists a unique morphism of $\delta$-graded rings $\mathbb S_n\rightarrow B$ such that $z_i\mapsto f_i$ for all $i\in\{0,\ldots,n\}$, where $\mathbb S_n$ has shift $(r_0,w_0)$.
 To check this, note that there exists a unique morphism of $\delta$-rings $R[z,z',z'',\ldots]\rightarrow C$ that sends $z_i\mapsto f_i$ for all $i\in\{0,\ldots,n\}$ and this morphism sends $\mathbb S_n$ to $B$.

Let $\mathbb I$ stand for either $\mathbb I_g$ or $\mathbb I_{g,\ord}$ and let $M$ stand for either $M_g$ or $M_{g,\ord}$.
Applying the above universal property to $B=\mathbb I$ and $C=M$ (see the paragraph of Equation (\ref{EQ012.111})), we get that for every linear system $\f=(f_0,\ldots,f_n)$ in $\mathbb I$ whose all entries belong to $\mathbb I^{r_0}(w_0)$, we have an induced morphism of filtered $\delta$-graded rings 
\begin{equation}\label{EQ033}
\mathbb S_n\rightarrow \mathbb I
\end{equation}
over $R$ that maps $z_i$ to $f_i$, where $\mathbb S_n$ has shift $(r_0,w_0)$.
Then the morphism (\ref{EQ033}) induces a morphism of filtered ssa-$\delta$-spaces 
$$\Proj_{\delta}(\mathbb I,\f)\rightarrow \Proj_{\delta}(\mathbb S_n,z).$$
 
\subsubsection{Rings of invariants of $\pmb{\SL}_g$-actions}\label{S234}
 
We now define $R$-modules and $R$-algebras that will play key roles in our theory.
 Let 
 $$\begin{array}{l}
 T=T^{(0)}=(T_{ij})_{1\leq i,j\leq g}=(T^{(0)}_{ij})_{1\leq i,j\leq g},\\
 T'=T^{(1)}=(T'_{ij})_{1\leq i,j\leq g}=(T^{(1)}_{ij})_{1\leq i,j\leq g},\\
\ldots,\;
 T^{(r)}=(T^{(r)}_{ij})_{1\leq i,j\leq g},\;
\ldots\end{array}$$
 be a sequence of symmetric matrices of indeterminates over $R$. For all integers $1\leq i\leq j\leq g$ and $r\geq 0$, $T_{ij}^{(r)}=T_{ji}^{(r)}$ will have degree $1$. For every $s\in \frac{1}{2}\mathbb Z=\{\frac{m}{2}|m\in \mathbb Z\}$ such that $gs\in \mathbb Z$ we
 denote by $\mathbb H^r_g(s)$ the $R$-module of homogeneous polynomials
 $$G\in R[T,\ldots,T^{(r)}]:=R[T^{(l)}_{ij}|1\leq i\leq j\leq g, 0\leq l\leq r]$$
 of degree $gs$ that satisfy
 \begin{equation}\label{EQ034}
 G(\Lambda T\Lambda^{\t},\ldots,\Lambda T^{(r)}\Lambda ^{\t})=\det(\Lambda)^{2s}G(T,\ldots,T^{(r)})
 \end{equation}
 for all $\Lambda\in\pmb{\GL}_g(R)$. 
 Here and later, for $G$ as above and symmetric matrices $M_0=(m_{0,ij})_{1\leq i,j\leq g},\ldots, M_r=(m_{r,ij})_{1\leq i,j\leq g}$ with entries in an $R$-algebra we use the notation 
$$G(M_1,\ldots,M_r):=G(m_{0,11}, m_{0,12},\ldots,m_{0,gg},\ldots,m_{r,11}, m_{r,12},\ldots,m_{r,gg}).$$
A similar notation will be used later (see Subsection \ref{S35}) for $G$ a power series instead of a polynomial.
 
 We also set $\mathbb H^r_g(s):=0$ if $s\in \frac{1}{2}\mathbb Z$ and $gs\not\in \mathbb Z$.
 Clearly, we have $\mathbb H^r_g(s)=0$ for all $s\in \frac{1}{2}\mathbb Z$ with $s<0$. We then define the $R$-algebras
$$\mathbb H^r_g:=\bigoplus_{s\in \mathbb Z}\mathbb H^r_g(s),\ \ \ \ 
\mathbb H_{g,\tot}^r:=\bigoplus_{s\in \frac{1}{2}\mathbb Z}\mathbb H^r_g(s),$$
$$
\mathbb H^{r,\torus}_g:= \mathbb H^r_g[z_0,\ldots,z_{r-1},\frac{1}{z_0\cdots z_{r-1}}],\ \ \ \mathbb H_g^{\torus}:=
\bigcup_{r\geq 0} \mathbb H^{r,\torus}_g,$$
$$
\mathbb H_{g,\tot}^{r,\torus}:=\mathbb H_{g,\tot}^r[z_0,\ldots,z_{r-1},\frac{1}{z_0\cdots z_{r-1}}
],\ \ \ \mathbb H_{g,\tot}^{\torus}:=
\bigcup_{r\geq 0} \mathbb H_{g,\tot}^{r,\torus},$$
where $(z_n)_{n\geq 0}$ is a sequence of indeterminates. We view $\mathbb H_{g,\tot}^{\torus}$ as a filtered $W$-graded ring with filtration given by the subrings in the union that defines it and grading defined by letting each element of $\mathbb H^r_g(s)$ have weight
$-2s$ and declaring that
$$\textup{ each $z_i$ has weight $\phi^{i+1}-\phi^i$.}$$
The $R$-algebras we just introduced are analogues of algebras in the classical invariant theory of multiple quadratic forms we study in Subsections \ref{S41} to \ref{S43}.

Similarly, we consider matrices of indeterminates 
$$X_1=(X_{1,ij})_{1\leq i,j\leq g}, \ldots, X_n=(X_{n,ij})_{1\leq i,j\leq g}$$
 and we denote by $\mathbb H^n_{g,\con}$ the $R$-algebra of all 
polynomials 
$$F\in R[X_1,\ldots,X_n]:=R[X_{l,ij}|1\leq i,j\leq g,1\leq l\leq n]$$ 
such that
$$F(\Lambda X_1 \Lambda^{-1},\ldots,\Lambda X_n \Lambda^{-1})=F(X_1,\ldots,X_n)$$
 for all $\Lambda\in {\pmb{\GL}}_g(R)$. The $R$-algebras $\mathbb H^n_{g,\con}$ are analogues of the algebras that show in the classical invariant theory of multiple endomorphisms and that are studied in Subsubsections \ref{S411}, \ref{S415} and \ref{S416}.

Some of our main results will establish links between
the algebras $\mathbb I^r_{g,\ord}$ and $\mathbb I^r_g$ and the algebras $\mathbb H_{g,\tot}^{r,\torus}$ and $\mathbb H^r_{g,\con}$ (see Theorems \ref{T5} and \ref{T7} below).

\subsubsection{Rings of $\delta$-power series}\label{S235}

If $f:C\rightarrow B$ is a set theoretical map between rings and $n\in \mathbb N$, the rule $(c_{ij})_{1\leq i,j\leq n}\mapsto (f(c_{ij}))_{1\leq i,j\leq n}$ defines a function $f:\textup{Mat}_n(C)\rightarrow \textup{Mat}_n(B)$ denoted also by $f$.

Our Serre--Tate expansion theory will involve the rings
$$S^r_{\for}:=\reallywidehat{R[[T]][T',\ldots,T^{(r)}]}:=\reallywidehat{R[[T_{ij}|1\leq 
i\leq j\leq g]][T^{(l)}_{ij}|1\leq i\leq j\leq g, 1\leq l\leq r]}.$$
Their union $S_{\for}$ has a natural structure of a filtered $\delta$-ring (with $\delta T_{ij}^{(l)}=T_{ij}^{(l+1)}$ for all integers $1\leq i\leq j\leq g$ and $l\geq 0$) which is
an object of $\Prol$.
 On the other hand the Fourier expansion 
 theory developed in \cite{BB} involved the rings
 $$S^r_{\can}:=\reallywidehat{R((q))[q',\ldots,q^{(r)}]}:=\reallywidehat{R((q))[q'_{ij},\ldots,q^{(r)}_{ij}|1\leq i \leq j \leq g]},$$
where 
$$\begin{array}{l}
q=q^{(0)}=(q_{ij})_{1\leq i,j\leq g}=(q^{(0)}_{ij})_{1\leq i,j\leq g},\\
q'=q^{(1)}=(q'_{ij})_{1\leq i,j\leq g}=(q^{(1)}_{ij})_{1\leq i,j\leq g},\\
\ldots, q^{(r)}=(q^{(r)}_{ij})_{1\leq i,j\leq g}, \ldots
\end{array}
$$
 is a sequence of symmetric matrices of indeterminates  and 
 $$R((q)):=R[[q_{ij}|1\leq i\leq j\leq g]][q_{ij}^{-1}|1\leq i\leq j\leq g].$$
  Again, the union $S_{\can}$ of these rings is an object of $\Prol$, with $\delta q_{ij}^{(l)}=q_{ij}^{(l+1)}$ for all integers $1\leq i\leq j\leq g$ and $l\geq 0$.

 In order to compare later on the Fourier and Serre--Tate expansion theories 
 (see Remark \ref{R7.4} below)
 we will need to consider the $R$-subalgebras
$$S^r_{\small-can}:=\reallywidehat{R[q_{ij},q_{ij}^{-1},q'_{ij},q''_{ij},\ldots,q^{(r)}_{ij}|1\leq i\leq j\leq g]}\subset S^r_{\can}$$ 
 and the $R$-algebra monomorphisms
$S^r_{\small-can}\rightarrow S^r_{\for}$, $F\mapsto F_{\textup{for}}$ that map $q_{ij},q'_{ij},q''_{ij},\ldots$ to $1+T_{ij},1+T_{ij}',1+T_{ij}'',\ldots$ (respectively). We have an induced monomorphism $\pmb{\Mat}_g(S^r_{\small-can}) \rightarrow \pmb{\Mat}_g(S^r_{\for})$, 
$M\mapsto M_{\textup{for}}$. Hence 
$$
q_{\textup{for}}=\mathbb E_g+T,\ \ q'_{\textup{for}}=\delta(\mathbb E_g+T), \ \ q''_{\textup{for}}= \delta^2(\mathbb E_g+T), \ldots,$$
where we recall that $\mathbb E_g$ is the $g\times g$ matrix with all entries equal to $1$.

A key role in our paper will be played by the matrix of series
\begin{equation}
\label{EQ034.9}
\Psi_q=(\Psi_{q_{ij}})_{1\leq i,j\leq g}\in \pmb{\Mat}_g(S^r_{\small-can})\subset 
\pmb{\Mat}_g(S^r_{\can}),\end{equation}
where 
 $$\Psi_{q_{ij}}:=\frac{1}{p}\log\left(1+p\frac{\delta q_{ij}}{q_{ij}^p}\right)\in S^r_{\small-can},$$
 and by its image
 \begin{equation}
 \label{EQ034.99}
 \Psi:=(\Psi_q)_{\textup{for}}\in \pmb{\Mat}_g(S^r_{\for}).\end{equation}
 
 \subsubsection{Rings of $\delta$-sections of line bundles}\label{S236}
 
Let $X$ be a smooth scheme over $\Spec(R)$ and $L$ a line bundle over $X$. 
For simplicity we will assume in what follows that the reduction modulo $p$, $\overline{X}$, is connected, hence an integral scheme.
We denote by
$$\mathbb L:=\textup{\bf Spec}\left(\bigoplus_{i\in \mathbb Z} L^{\otimes i}\right)$$
the principal $\mathbb G_m$-bundle over $X$ associated to $L$. 
The $\mathbb G_m$-action on $\mathbb L$ induces an $R^{\times}$-action, denoted by $\star$, on the $p$-jet spaces $J^r(\mathbb L)$ and on their rings of global functions $\mathcal O(J^r(\mathbb L))$. For each $w\in W$ we denote by $\mathbb S^r_{X,L}(w)$ the $R$-module of all 
$f\in \mathcal O(J^r(\mathbb L))$ such that we have $\lambda \star f=\lambda^{-w}\cdot f$ for all $\lambda\in R^\times$.

\begin{df}\label{df24}
The ring of {\it $\delta$-sections} of $L$ is the ring $\mathbb S_{X,L}$ defined by
$$\mathbb S_{X,L}:=\bigcup_{r=0}^{\infty}\mathbb S^r_{X,L}, \ \ \mathbb S^r_{X,L}:=\bigoplus_{w\in W(r)} \mathbb S^r_{X,L}(w).$$\end{df}

If $X=\Spec(B)$ is affine, $L$ is trivial over $X$, and $\{z_L\}$ is a basis of $L$, then we have isomorphisms
$$\mathcal O(J^r(\mathbb L))\simeq \reallywidehat{\mathcal O(J^r(X))[z_L,z_L^{-1},z_L',\ldots,z_L^{(r)}]}=\reallywidehat{J^r(B)[z_L,z_L^{-1},z_L',\ldots,z_L^{(r)}]},$$
$$ \mathbb S^r_{X,L}(w)\simeq \mathcal O(J^r(X))\cdot z_L^w=J^r(B)\cdot z_L^w.
$$
The $R$-algebra $\mathbb S_{X,L}$ has a structure of 
filtered graded ring coming from the filtered graded system $(\mathbb S^r_{X,L})_{r\geq 0}$ and with this structure it is a
filtered integral $\delta$-graded ring over $R$. In particular, 
for each linear system $\f=(f_0,\ldots,f_n)$ in $\mathbb S_{X,L}$ we define the filtered ssa-$\delta$-space (see Equation (\ref{EQ004}))
\begin{equation}\label{EQ036}
X_{\delta,L,\f}:=\Proj_{\delta}(\mathbb S_{X,L},\f)=\Spec_{\delta}\left(
\prod_{0\leq i_0<\cdots<i_{\bullet}\leq n} \mathbb S_{X,L,\langle f_{i_0}\cdots f_{i_{\bullet}}\rangle}
\right).\end{equation}

\begin{ex}\label{EX1}
Recall the ring $\mathbb S_n$ (see Equation (\ref{EQ031})).
If we take $X=\mathbb P^n=\Proj (R[z])$, $z=(z_0,\ldots, z_n)$, and $L=\mathcal O(1)$, then
one can prove (see \cite{Bu05}, Rm. 6.9, p. 167) that
$\mathbb S_{\mathbb P^n,\mathcal O(1)}=\mathbb S_n$,
hence, with notation as in Equations (\ref{EQ036}) and (\ref{EQ032}), we have
$$\mathbb P^n_{\delta,\mathcal O(1),z}=\Proj_{\delta}(\mathbb S_{\mathbb P^n,\mathcal O(1)},z)=\Proj_{\delta}(\mathbb S_n,z)=\mathbb P^n_{\delta}.$$
\end{ex}

\medskip

We will abbreviate $X_{\delta,L,\f}$ in Equation (\ref{EQ036}) by $X_{\delta}$ if $L$ and $\f$ are clear from the context.
In view of the above example this is consistent with the notation in Equation (\ref{EQ032}) in case $X=\mathbb P^n$, $L=\mathcal O(1)$, $\f=z$.

\medskip

Going back to the case of a general $X$ and $\f$, let $P\in X(R)$ be an $R$-valued point of $X$. Then we will construct an $n$-tuple
\begin{equation}\label{EQ037}
(f_0(P),\ldots,f_n(P))\in R^{n+1},\end{equation}
well-defined up to multiplication by an element of $R^{\times}$, as follows. We can assume $X$ is affine and $L$ is trivial on $X$. Fix a basis $\{z_L\}$ of $L$, write $f_i=\varphi_i\cdot z_L^w$, with $\varphi_i\in \mathcal O(J^r(X))$, view $\varphi_i$ as a morphism $\varphi_i:J^r(X)\rightarrow \widehat{\mathbb A^1_R}$, consider the point $J^r(P)\in J^r(X)(R)$, and set $f_i(P):=\varphi_i(J^r(P))\in\widehat{\mathbb A^1_R}(R)=R$. In other words, $P\mapsto f_i(P)$ is the $\delta$-function defined by $\varphi_i$ (see Definition \ref{df7}).

\begin{df}\label{df25}
Let $P\in X(R)$. We say that $P$ is: 

\medskip
{\bf (a)} {\it $\f$-semistable} if there exists $i\in \{0,\ldots,n\}$ such that $f_i(P)\in R^\times$. 

\smallskip
{\bf (b)} {\it $\f$-unstable} if $f_0(P)=\cdots=f_n(P)=0$. \end{df}

The definition makes sense even though the vector (\ref{EQ037}) is only defined up to multiplication by a scalar in $R^{\times}$. We denote
by $X(R)_{\f}^{\textup{ss}}$ and $X(R)_{\f}^{\textup{u}}$ the sets of $\f$-semistable and $\f$-unstable (respectively) points in $X(R)$. We have
$$X(R)_{\f}^{\textup{ss}}\sqcup X(R)_{\f}^{\textup{u}}\subset X(R)$$ and in general the inclusion is strict. 

Example: for $X=\mathbb P^n$ and $\f=z=(z_0,\ldots,z_n)$, we have $\mathbb P^n(R)^{\textup{ss}}_z = \mathbb P^n(R)$. 

For general $X$ and $\f$ we have a canonical evaluation map
\begin{equation}\label{EQ038}
\mathcal E_{X,\f}:X(R)^{\textup{ss}}_{\f}\rightarrow X_{\delta}(\delta\textup{-pts})\end{equation}
defined as follows. For $P\in X(R)_{\f}^{\textup{ss}}$, we consider the non-empty subset
$$I:=\{i|f_i(P)\in R^\times\}\subset \{0,\ldots,n\},$$
and for each increasing sequence $0\leq i_0<\cdots <i_d \leq n$ 
let $R_{i_0\ldots i_d}$ be $R$ if $\{i_0,\ldots,i_d\}\subset I$ and be $0$ if $\{i_0,\ldots,i_d\}\not\subset I$.
Then consider the discrete ssa-$\delta$-space
\begin{equation}
\label{EQ038.5}
P_{\delta}:=\Spec\left( \prod_{0\leq i_0<\cdots <i_{\bullet}\leq n} R_{i_0\ldots i_{\bullet}}\right),\end{equation}
with obvious face maps; one can easily check that $P_{\delta}$ is connected and acyclic, hence a $\delta$-point. Next we consider 
 the $R$-algebra homomorphisms
$$\mathbb S_{X,L,\langle f_{i_0}\cdots 
f_{i_{\bullet}}\rangle}\rightarrow R_{i_0\ldots i_{\bullet}}$$
defined when $\{i_0,\ldots,i_{\bullet}\}\subset I$ as the evaluation at $P$. The induced morphism $P_{\delta}\rightarrow X_{\delta}$ is $\mathcal E_{X,\f}(P)$.
Note that there exists a canonical diagram
$$\begin{array}{ccccl}
X(R)_{\f}^{\textup{ss}} & \stackrel{\Phi_{\f}}{\longrightarrow} & \mathbb P^n(R)^{\textup{ss}}_z & = & \mathbb P^n(R)\\
\Big\downarrow {\mathcal E_{X,\f}} & \ & \Big\downarrow {\mathcal E_{\mathbb P^n,z}} & \ & \ \\
X_{\delta}(\delta\textup{-pts}) & \longrightarrow & \mathbb P_{\delta}^n(\delta\textup{-pts}), & \ & \ 
\end{array}$$
where
$$\Phi_{\f}(P):=(f_0(P):\cdots:f_n(P)).$$

As with $\delta$-functions (see Definition \ref{df7}) the maps $\Phi_{\f}$ are not induced by morphisms of formal schemes
from open formal subschemes of $\widehat{X}$ to $\reallywidehat{\mathbb P^n}$. Rather,
they are 
 induced by morphisms of formal schemes from open formal subschemes of 
 $J^r(X)$ to $\reallywidehat{\mathbb P^n}$ and hence they are
 given, locally in the Zariski topology of $X$, by restricted power series in local affine coordinates $x_1,\ldots,x_d$ and their $p$-derivatives $\delta x_1,\ldots,\delta x_d,\ldots,\delta^r x_1,\ldots,\delta^r x_d$; such maps are referred to in \cite{Bu95a}, Introd. and \cite{Bu05}, Ch. 2, Sect. 2.1 as {\it $\delta$-maps}. 

\begin{ex}\label{EX2}
We recall that for an integer $N\geq 3$ prime to $p$ the moduli scheme $\mathcal X=\mathcal A_{g,1,N,R}$ is smooth over $\Spec(R)$. Let $X\subset \mathcal X$ be an open subscheme with $\overline{X}$ connected and 
let $\pi:A_X\rightarrow X$ the universal abelian scheme with principal polarization $\theta_X$. Recall the 
{\it Hodge line bundle} over $X$,
$$L:=\pi_*(\wedge^g \Omega_{A_X/X}).$$
With $X$ and $L$ as above we set
\begin{equation}\label{EQ039}
\mathbb M^r_X=\bigoplus_{w\in W(r)} \mathbb M^r_X(w):=\mathbb S^r_{X,L}=\bigoplus_{w\in W(r)}\mathbb S^r_{X,L}(w),\ \ \ \mathbb M_X:=\mathbb S_{X,L}.\end{equation}
Let $\mathbb M$ stand for either $\mathbb M_g$ or $\mathbb M_{g,\ord}$; in the latter case we also assume $X\subset \mathcal X_{\ord}$.
We will define a canonical morphism of filtered $\delta$-graded rings
\begin{equation}\label{EQ040}
\mathbb M\rightarrow \mathbb M_X.
\end{equation}

First we assume that $X$ is affine and that there exists a column basis $\omega_X=(\omega_{X,1},\ldots,\omega_{X,g})^{\t}$ of $1$-forms on $A_X$.
Then for $f\in M^r_g(w)$ we consider the element
 \begin{equation}\label{EQ041}
 f_X:=f(A_X,\theta_X,\omega_X, \mathcal O^{\infty}(X))\in \mathcal O(J^r(X))\simeq \mathcal O(J^r(X))\cdot z_L^w,\end{equation}
 where we identify $z_L$ with the canonical global section of $L$ defined by the wedge product $\omega_{X,1}\wedge\ldots \wedge \omega_{X,g}$. The map $f\mapsto f_X$ defines $R$-linear maps
 \begin{equation}\label{EQ042}
 M^r_g(w)\rightarrow \mathcal O(J^r(X))\simeq \mathcal O(J^r(X))\cdot z_L^w
 \end{equation}
 \begin{equation}\label{EQ043}
 M^r_{g,\ord}(w)\rightarrow \mathcal O(J^r(X))\simeq \mathcal O(J^r(X))\cdot z_L^w
 \end{equation}
 that depend on $\omega_X$.
 For $X$ arbitrary, the above maps induce an $R$-algebra homomorphism
 (\ref{EQ040}) which is intrinsic, i.e., is independent of the choices made to get such local column bases of $1$-forms. By Corollary \ref{C4} (d) below the homomorphism (\ref{EQ040}) is injective with injective reduction modulo $p$.
\end{ex}

\subsubsection{Analogies with geometric invariant theory}\label{S237}

Before proceeding to explain our main results (see Subsection \ref{S24}) we 
give here a quick summary of the various ssa-$\delta$-spaces introduced so far and we point out some 
 analogies with geometric invariant theory as in \cite{MFK}.

Let $(\mathbb I, \mathbb M, \mathbb M')$ be one of the two triples
\begin{equation}\label{EQ043.1}
(\mathbb I_g, \mathbb M_g, \mathbb M'_g)\ \ \ \ \text{or} \ \ \ (\mathbb I_{g,\ord}, \mathbb M_{g,\ord}, \mathbb M'_{g,\ord}),
\end{equation}
see Equation (\ref{EQ029.2}) and its ordinary analogue,
 and recall that $\mathbb I, \mathbb M, \mathbb M'$ are 
 filtered integral $\delta$-graded rings. Recall also the filtered integral $\delta$-graded rings 
 $\mathbb S_n$ and $\mathbb M_X$, where $X\subset\mathcal X=\mathcal A_{g,1,N,R}$
 or $X\subset\mathcal X_{\textup{ord}}$,
 $N\geq 3$ an integer prime to $p$, $\overline{X}$ connected, see Definition \ref{df23} and Equation (\ref{EQ039}). We have homomorphisms of filtered $\delta$-graded rings
$$\mathbb S_n\rightarrow \mathbb I \rightarrow \mathbb M\rightarrow \mathbb M_X\ \ \ \text{and} \ \ \ \mathbb M\rightrightarrows \mathbb M',$$
where the homomorphism $\mathbb S_n\rightarrow \mathbb I$ is defined by a linear system $\f$ of length $n+1$ in $\mathbb I$ (see the morphism (\ref{EQ033})) and the other homomorphisms 
(see Equations (\ref{EQ012.2}) and (\ref{EQ040})) are the canonical ones which, as stated, are injective, with injective reduction modulo $p$.
The main geometric objects of our theory are then the filtered ssa-$\delta$-spaces:
 \begin{equation}\label{EQ044}
 \begin{array}{rcl}
 X_{\delta} & := & \Proj_{\delta}(\mathbb M_X,\f),\\
 \ & \ & \ \\
 Y'_{\delta} & := & \Proj_{\delta}(\mathbb M',\f),\\
 \ & \ & \ \\
 Y_{\delta} & := & \Proj_{\delta}(\mathbb M,\f),\\
 \ & \ & \ \\
 Z_{\delta} & := & \Proj_{\delta}(\mathbb I,\f),\\
 \ & \ & \ \\
 \mathbb P^n_{\delta} & := & \Proj_{\delta}(\mathbb S_n,z),\end{array}
 \end{equation}
 where we continue to denote by $\f$ the images of $\f$ in the rings $\mathbb M, \mathbb M', \mathbb M_X$.
 There exist naturally induced morphisms of filtered ssa-$\delta$-spaces
 \begin{equation}\label{EQ045}
 X_{\delta}\rightarrow Y_{\delta}\rightarrow Z_{\delta}\rightarrow \mathbb P^n_{\delta}\ \ \ \text{and}\ \ \ Y'_{\delta} \rightrightarrows Y_{\delta}.\end{equation}
 We also have an induced map of sets
\begin{equation}\label{EQ046}
\Phi_{\f}:X(R)^{\textup{ss}}_{\f} \rightarrow \mathbb P^n(R)
\end{equation}
compatible with the map $X_{\delta}(\delta\textup{-pts})\rightarrow \mathbb P^n_{\delta}(\delta\textup{-pts})$.
The morphism 
\begin{equation}\label{EQ047}
Y_{\delta}\rightarrow Z_{\delta}\end{equation}
 is a categorical quotient (coequalizer) in $\textup{\bf FilssAff}_{\delta}$ of the pair of morphisms 
\begin{equation}\label{EQ048}
Y'_{\delta} \rightrightarrows Y_{\delta};\end{equation} 
see Remark \ref{R6.5}.
By analogy with geometric invariant theory (see below), 
the pair of morphisms (\ref{EQ048}) can be viewed as an avatar (an analogue/realization in $\delta$-geometry) of the pair of morphisms
of functors from $\textup{\bf Ring}$ to $\textup{\bf Set}$
\begin{equation}
\label{EQ048.1}
\mathcal A'_{g,R}\stackrel{\pi_1,\pi_2}{\rightrightarrows} \mathcal A_{g,R}\end{equation}
 from the union $\mathcal A'_{g,R}$ of all prime to $p$ `Hecke correspondences' on $\mathcal A_{g,R}$. (A similar picture holds for the ordinary locus.)
 Note that as $N\geq 3$ is prime to $p$ the coequalizer in the category of schemes over $\Spec(R)$ of the correspondence
 \begin{equation}
\label{EQ048.2}
\mathcal A'_{g,1,N,R}\stackrel{\pi_1,\pi_2}{\rightrightarrows} \mathcal A_{g,1,N,R}\end{equation}
obtained by taking the union of all prime to $p$ `Hecke correspondences' on
$\mathcal A_{g,1,N,R}$
 is trivial, i.e., is isomorphic to $\Spec(R)$. This is so as, for an algebraically closed field $\kappa$, the generic  orbits of the Hecke correspondences are Zariski dense in $\mathcal A_{g,1,N,\kappa}$ (see \cite{Cha}, Thm. 2, for $\kappa$ of characteristic $p$).
Moreover, the morphism (\ref{EQ047}) can be viewed as an avatar of the coequalizer 

 \begin{equation}
\label{EQ049}
\mathcal A_{g,R}\rightarrow \mathcal A_{g,R}/\mathcal A'_{g,R}\end{equation}
of (\ref{EQ048.1}) in the category of functors from $\textup{\bf Ring}$ to $\textup{\bf Set}$,
or of the similar map in the ordinary case.
Similarly, we can view the morphism
\begin{equation}\label{EQ049.1}
X_{\delta}\rightarrow Z_{\delta}
\end{equation}
 as an avatar of the coequalizer of (\ref{EQ048.2}) in the category of functors from $\textup{\bf Ring}$ to $\textup{\bf Set}$.

The various properties of the ring $\mathbb I$ to be proved in Section \ref{S3}
 (such as $\phi$-transcendence degree, asymptotic transcendence degree, $\phi$-finite generation, etc.) correspond to `geometric properties' of the morphisms (\ref{EQ045}) and of the maps (\ref{EQ046}). We will not pursue this in a systematic way but some of these geometric properties will be explained in Subsection \ref{S24}. Here is an example of such a property.
 Let us say that a subset $C\subset X(R)$ is a 
{\it polarized isogeny class} 
if there exists $P\in C$ such that $C$ consists of all points $P'\in X(R)$ 
 with the property that there exists a polarized isogeny between the principally polarized abelian schemes over $\Spec(R)$ that correspond to $P$ and $P'$. 
 Directly from the definition of Hecke covariance we 
get that the maps  (\ref{EQ046})  are constant on polarized isogeny classes.
 
Here are some comments on the analogy between the situation described above and the situation in classical geometric invariant theory as in \cite{MFK}. 
Given an action $\rho:X':=G\times X\rightarrow X$
 of an algebraic group $G$ on a variety $X$ over some algebraically closed field, one says that a morphism $\pi:X\rightarrow Y$ of varieties is a categorical quotient for this action if it is a coequalizer, in the category of varieties, 
 of the two morphisms 
 \begin{equation}
 \label{EQ050.01}
 \rho,\pi_2:X' \rightrightarrows X\end{equation}
 where $\pi_2$ is the second projection. In our setting, the pair $(\ref{EQ050.01})$ is replaced by Hecke correspondences on the moduli stack $\mathcal A_{g,R}$ over $\Spec(R)$, see Equation (\ref{EQ048.1}). The fact that (\ref{EQ047}) is a coequalizer in $\textup{\bf FilssAff}_{\delta}$ of the pair of morphisms (\ref{EQ048}), suggests, following the language of the geometric invariant theory, to regard $Z_{\delta}$ as a `quotient', in $\delta$-geometry, of (an open subscheme of) the moduli scheme $\mathcal A_{g,1,N,R}$ or of its ordinary locus by the Hecke correspondences. As we shall see, $Z_{\delta}$ is a highly non-trivial object, and this is one of the main motivations here (and in \cite{Bu00, BB, Bu03, Bu05}) for enriching classical algebraic geometry with the $\delta$-geometry.

In the same vein, we view the maps $\Phi_{\f}$ in (\ref{EQ046}) as analogues, in the $\delta$-geometry of Hecke correspondences,
of maps to projective spaces defined by $G$-invariant sections of $G$-linearized bundles over a variety with $G$-action in classical geometric invariant theory, see \cite{MFK}, Ch. 1, Sect. 3. 
The role of $G$-linearization is played, in our formalism, by the consideration of the 
isomorphism $\beta$ in Equation (\ref{EQ012.1111}).
Our $\f$-semistable points are then an analogue of the $G$-semistable points in \cite{MFK}, Ch. 1, Sect. 4. 

One can ask for a $\delta$-geometric analogue of the stable points in \cite{MFK}, Ch. 1, Sect. 4. Such a concept could have an arithmetic nature as follows. If $R_{\textup{alg}}:=\{x\in R|x\;\textup{algebraic over}\; \mathbb Q\}$, then
 $\mathcal X$ comes via base change from a moduli scheme $\mathcal X_{\textup{alg}}$ over $\Spec(R_{\textup{alg}})$ and we take $X$ such that it comes via base change from an open subscheme $X_{\textup{alg}}$ of $\mathcal X_{\textup{alg}}$. Then our  tentative definition is as follows.
 A point $P\in X_{\textup{alg}}(R_{\textup{alg}})\subset X(R)$ is called {\it $\f$-stable} if it is $\f$-semistable and the intersection $\Phi_{\f}^{-1}(\Phi_{\f}(P))\cap X_{\textup{alg}}(R_{\textup{alg}})$ taken inside $X(R)$ is a finite union of polarized isogeny classes. We will not further investigate this concept.

\subsection{Main results}\label{S24}
 
\subsubsection{Serre--Tate expansions and the basic forms}
 
 Our approach is based on the following {\it Serre--Tate expansion principle} which is analogous to the Fourier expansion principle in \cite{BB}, Thm. 3.1. For the next statement fix an ordinary elliptic curve $E_0$ over $k=R/pR$ and a generator $e_0$ of (the $\mathbb Z_p$-module defined by) $T_p(E_0)$ and set $A_0:=E_0^g$. We endow $A_0$ with the standard diagonal principal polarization $\theta_0$ and with the diagonal column basis $e$ of $T_p(A_0)$ defined by $e$ by which we mean a column vector $e=[e_{0,1}\;\cdots\;e_{0,g}]^{\t}$ where each $e_{0,i}$ with $i\in\{1,\ldots,g\}$ is $e_0$ but viewed as a generator of the Tate module of the $i$-th copy $E_{0,i}$ of $E_0$ in $A_0$.
 
 \begin{thm}\label{T1}
 The following two properties hold:
 
 \medskip
 {\bf (a)}
 For every $w',w''\in\mathcal W$ there exist natural injective $R$-linear maps
 $$\mathcal E_{A_0,\theta_0,e}:\mathbb I^r_{gg}(w',w'')\rightarrow \pmb{\Mat}_g(S^r_{\textup{for}}),\ \ \ \mathcal E_{A_0,\theta_0,e}:\mathbb I^r_{gg,\ord}(w',w'')\rightarrow \pmb{\Mat}_g(S^r_{\textup{for}}),$$
 whose reductions modulo $p$ are also injective.
 
 \smallskip
 {\bf (b)} For each $w\in W$ there exist natural injective $R$-linear map
 $$\mathcal E_{A_0,\theta_0,e}:\mathbb I^r_g(w)\rightarrow S^r_{\textup{for}},\ \ \ \mathcal E_{A_0,\theta_0,e}:\mathbb I^r_{g,\ord}(w)\rightarrow S^r_{\textup{for}},$$
 whose reductions modulo $p$ are also injective. \end{thm}
 
These $R$-linear maps will be compatible
 with $R$-linear maps $\mathbb I^r_{gg}(w',w'')\rightarrow \mathbb I^r_{gg,\ord}(w',w'')$ and $\mathbb I^r_g(w)\rightarrow \mathbb I^r_{g,\ord}(w)$ 
 (which are therefore injective with injective reduction modulo $p$) and will be compatible in a natural way with multiplication.
We refer to the $R$-linear maps in Theorem \ref{T1} as {\it Serre--Tate expansion maps}. 
On our way of proving Theorem \ref{T1} we will prove that the $R$-algebra homomorphism
in Equations
(\ref{EQ040}) 
 is injective with injective reductions modulo $p$ (see Proposition \ref{P3} (b)).
 
 \medskip
 
Recall the matrix $\Psi$ in Equation (\ref{EQ034.99}). We will prove the following existence result (see Theorem \ref{T13} (a)):
 
\begin{thm}\label{T2}
There exists a generator $e_0$ of the Tate module of $E_0$ with the following property.
For every integer $r\geq 1$ there exists a unique Hecke covariant Siegel $\delta$-modular form
$$f^r:=f^r_{g,\crys}\in \mathbb I^r_{gg}(\phi^r,1)$$
whose Serre--Tate expansion equals
$$\mathcal E_{A_0,\theta_0,e} (f^r_{g,\crys})=
\Psi^{\phi^{r-1}}+p\Psi^{\phi^{r-2}}+\cdots+p^{r-1}\Psi.$$
\end{thm}

\begin{rem}\label{R7.4}
We will refer to the forms $f^r$
 as the {\it basic} Hecke covariant Siegel $\delta$-modular forms of weight $(\phi^r,1)$. These forms were introduced 
 in \cite{BB}, Subsect. 4.1 and their (crystalline) construction will be recalled in Subsection \ref{S34}.
By construction, if $p>2$ (resp. if $p=2$) these forms vanish on all quadruples $(A,\theta,\omega,R)$ in which $A$ is an abelian scheme over $\Spec(R)$ which is a canonical lift (resp. quasi-canonical lift) of its reduction $A_k$ modulo $k$. For $r=1$ the converse holds: if $p>2$ (resp. if $p=2$) and if $f^1$ vanishes at some quadruple $(A,\theta,\omega,R)$ with $A$ ordinary, then $A$ is an abelian scheme over $\Spec(R)$ which is a canonical lift (resp. quasi-canonical lift) of $A_k$.

 Recall from \cite{BB}, Subsect. 2.3 that for $w\in W$ and $w',w''\in\mathcal W$ we have natural {\it Fourier expansion maps}
 $$\mathcal E_{\infty}:\mathbb I^r_{gg}(w',w'')\rightarrow \pmb{\Mat}_g(S^r_{\can}),\ \ \ \mathcal E_{\infty}:\mathbb I^r_g(w)\rightarrow S^r_{\can}.$$
 These maps will not play any role in the present paper. Note however that
 by \cite{BB}, Eq. (4.4) 
 the Fourier expansion of $f^r_{g,\crys}$
 is given by
 \begin{equation}\label{EQ050.1}
 \mathcal E_{\infty}(f^r_{g,\crys})=\Psi_q^{\phi^{r-1}}+p\Psi_q^{\phi^{r-2}}+\cdots+p^{r-1}\Psi_q.\end{equation}
 Combining Equation (\ref{EQ050.1}) and Theorem \ref{T2} we get that the Fourier and Serre--Tate expansions of the basic forms $f^r_{g,\crys}$ coincide up to natural identification:
 $$ (\mathcal E_{\infty}(f^r_{g,\crys}))_{\textup{for}}= \mathcal E_{A_0,\theta_0,e} (f^r_{g,\crys}).$$ 
\end{rem}

We have the following rank $1$ result (see Corollaries \ref{C10} and \ref{C11} below):

\begin{thm}\label{T3}
The following two properties hold:

\medskip
{\bf (a)} The $R$-modules $\mathbb I^1_{gg}(\phi,1)$ and 
$\mathbb I^1_{gg,\ord}(\phi,1)$ coincide and are 
 free of rank $1$ generated by $f^1_{g,\crys}$.
 
 \smallskip
{\bf (b)} For all integers $n\geq 1$, the $R$-modules $\mathbb I^1_g(-n\phi-n)$ and $\mathbb I^1_{g,\ord}(-n\phi-n)$ coincide and are free of rank $1$ generated by $(\det(f^1_{g,\crys}))^n$.
\end{thm}

In addition, generalizing a construction for $g=1$ in \cite{Ba} and \cite{Bu05}, in Section 8.4 we will prove (see Theorem \ref{T14} and Corollaries \ref{C15} (e) and \ref{C17}):

\begin{thm}\label{T4}
The following three properties hold:

\medskip
{\bf (a)}
The $R$-module $\mathbb I^1_{gg,\ord}(-\phi,1)$ is free of rank $1$ and has a generator $f^{\partial}_{g,\crys}$
whose Serre--Tate expansion is the identity matrix:
$$\mathcal E_{A_0,\theta_0,e}(f^{\partial}_{g,\crys})=1_g.$$

\smallskip
{\bf (b)} 
The form $\det(f^{\partial}_{g,\crys})$ is invertible in the ring $\mathbb I^1_{g,\ord}$ and if $w\in W(r)$ is such that $\deg(w)=0$ (equivalently, if $w\in (\phi-1)W$), then the $R$-module $\mathbb I^r_{g,\ord}(w)$ is free of rank $1$ generated by
$$\det(f^{\partial}_{g,\crys})^{\frac{w}{\phi-1}}.$$

\smallskip
 {\bf (c)} We have $\mathbb I^1_g(\phi-1)=\mathbb I^1_g(1-\phi)=0$.
\end{thm}

\begin{rem}\label{R8}
We shall refer to $f^{\partial}_{g,\crys}$ as the {\it basic} ordinary Siegel $\delta$-modular form of weight $(-\phi,1)$: the free $R$-modules  $\mathbb I^r_{gg,\ord}(w',w'')$ for $w',w''\in \mathcal W$ of order $\leq r$ can be fully described using $f^{\partial}_{g,\crys}$ (see Theorem \ref{T14} below). These modules have rank $r$ if $w',w''>0$, have rank $1$ if $w'w''<0$, and vanish (have rank $0$) if $w'<0$, $w''<0$. However, the determination of the submodules $\mathbb I^r_{gg}(w',w'')\subset \mathbb I^r_{gg,\ord}(w',w'')$ for all $w',w''\in \mathcal W$ is an open problem.
\end{rem}

\subsubsection{Structure of the ring $\mathbb I_{g,\ord}$}

We next present a series of applications of our Serre--Tate expansion theory to the structure of the rings $\mathbb I_{g,\ord}$ (see Theorem \ref{T15} (b) below):

 \begin{thm}\label{T5}
 There exist natural monomorphisms of graded $R$-algebras
 $$\diamondsuit:\mathbb H_{g,\tot}^{r-1,\torus}\rightarrow \mathbb I^r_{g,\ord},\ \ \ \ 
 \diamondsuit:\mathbb H_{g,\tot}^{\torus}\rightarrow \mathbb I_{g,\ord},$$
 whose cokernels are torsion.
 \end{thm}
 
 \begin{rem}\label{R9}
 By Theorem \ref{T5} we have induced $K$-algebra isomorphisms
 $$\diamondsuit:\mathbb H_{g,\tot}^{r-1,\torus}\otimes_R
 K\rightarrow \mathbb I^r_{g,\ord}\otimes_R K,\ \ \ \ 
 \diamondsuit:\mathbb H_{g,\tot}^{\torus}\otimes_R K\rightarrow \mathbb I_{g,\ord}\otimes_R K;$$
so the study of the rings $\mathbb I^r_{g,\ord}\otimes_R K$ gets reduced to the study of the rings
 $\mathbb H^r_g\otimes_R K$, which is a topic in the geometric invariant theory of multiple quadratic forms. This invariant theory has classical roots and will be revisited, and further developed, in Subsections \ref{S41} to \ref{S43}. 
 In particular, as we shall see, the $K$-algebras $\mathbb H^r_{g,\tot}\otimes_R K$ are finitely generated, Cohen--Macauley, and unique factorization domains, thus Gorenstein, and hence so are the $K$-algebras $\mathbb I^r_{g,\ord}\otimes_R K$.
Also we will see that  (see Corollary \ref{C15} (b) and (c)):
 
 \medskip
 $\bullet$ if $g$ is odd, then we have $\mathbb H_{g,\tot}^r=\mathbb H^r_g$
 hence $\mathbb H_{g,\tot}=\mathbb H_g$ which implies $\mathbb I^r_{g,\ord}(w)=0$ for all integers  $r\geq 1$ and $w\in W(r)$ with $\deg(w)$ odd;
 
 \smallskip
 $\bullet$ if $g$ is even, then $\mathbb H_{g,\tot}\neq \mathbb H_g$ 
 and hence
 the $R$-algebras 
 $\mathbb I^r_{g,\ord}(w)$ do not vanish for all integers  $r\geq 1$ and $w\in W(r)$ with $\deg(w)$ odd. 
 
\medskip
Subsection \ref{S42} provides, in particular, a detailed study of the structure of algebra $\mathbb H_{g,\tot}^r\otimes_R K$ tensored with an algebraic closure of $K$. The results there easily translate into results
 about the algebra $\mathbb I^r_{g,\ord}\otimes_R K$;
  we leave this translation to the reader. In particular  we will introduce in 
 Subsection \ref{S4206}
 some special  `$\Theta$ and $\Upsilon$ invariants' 
and we will compare the algebras of all invariants with the algebra generated by the $\Theta$
invariants (see, for instance,  Theorems \ref{T29} and \ref{T30} and Corollary \ref{C24}). We will deduce information on the Hilbert series of the rings of invariants (see Corollaries \ref{C21} and \ref{C22}), the smallest number of homogeneous generators
(see Proposition \ref{P15} (a)), and transcendence bases (see Theorem \ref{T31}). 
For $g=2$, Subsection \ref{S43} gives finer results on Hilbert series (see Proposition \ref{P18} and Corollary \ref{C26}), the smallest number of homogeneous generators
(see Proposition \ref{P17}), and relations (see Corollary \ref{C25}).

 On the other hand, the $R$-algebra homomorphism $\diamondsuit:\mathbb H_{g,\tot}^{\torus}\rightarrow \mathbb I_{g,\ord}$ is not an isomorphism. It would be interesting to compute the cokernel of the $R$-linear map $\diamondsuit$. A related open problem is to compute the cokernel (which is, again, torsion)
 of the $R$-linear map given by the inclusion
 $$R[z,z^{\phi},z^{\phi^2},\ldots]\rightarrow \mathbb S_n.$$
 Proposition \ref{P7.5} (b) below will provide a lower bound for the torsion of the cokernel of $\diamondsuit$ in the following sense. If for an $R$-module $M$ we denote by $M[p]$ the $k$-linear space of all elements of $M$ annihilated by $p$, then for $s\in\frac{1}{2}\mathbb Z_{\geq 0}$ we will show:
 $$\dim_k\left(
 \frac{\mathbb I^{r+1}_{g,\ord}(-2(p+1)s)}{\diamondsuit(\mathbb H^r_g((p+1)s))}
 [p]\right)\geq 2 \dim_K(\mathbb H^{r-1}_g(s)_K)-3.$$
 It would be interesting to find the precise value of the left-hand side of the above inequality.
 \end{rem}
 
 \begin{rem}\label{R10}
 By combining the proofs of Theorems \ref{T2}, \ref{T4},  and \ref{T5}, the Serre--Tate expansions of all elements of $\mathbb I_{g,\ord}$ will be made explicit in Remark \ref{R23} (a).\end{rem}
 
 \begin{rem}\label{R11}
 Let $X\subset \mathcal X_{\ord}=\mathcal A_{g,1,N,R,\ord}$ be an open subscheme with $\overline{X}$ connected.
We check that for each polarized isogeny class $C\subset X(R)$ there exist `many' Hecke covariant Siegel $\delta$-modular forms that vanish on $C$. First, by Theorem \ref{T5}, there exist (many) integers $r\geq 1$ and (many) weights $w\in W(r)$
such that the (finitely generated) $R$-module $\mathbb I_{g,\ord}^r(w)$ has `large' rank. 
For a point $P\in C$ consider the principally polarized abelian scheme $(A,\theta)$ over $\Spec(R)$
that corresponds to it and the submodule $\mathbb I_{g,\ord}^{r,[P]}(w)$
of all $f\in \mathbb I_{g,\ord}^r(w)$ that vanish at $P$ in the sense that $f(A,\theta,\omega,R)=0$ for one (equivalently every) $\omega$. Note that, by Hecke covariance, each 
$f\in \mathbb I_{g,\ord}^{r,[P]}(w)$ vanishes on the whole of $C$.
We claim that the factor module $\mathbb I_{g,\ord}^r(w)/ \mathbb I_{g,\ord}^{r,[P]}(w)$ is either $0$ or (non-canonically) isomorphic to $R$. Indeed, if this factor module is non-zero, let $0\neq f_0\in \mathbb I_{g,\ord}^r(w)$ be such that the element $f_0(A,\theta,\omega,R)\in R$ has minimum $p$-adic valuation and we conclude by noting that $\mathbb I_{g,\ord}^{r,[P]}(w)$ is the kernel of the surjective $R$-linear homomorphism
$$\mathbb I_{g,\ord}^r(w)\rightarrow R,\ \ \ f\mapsto \frac{f(A,\theta,\omega,R)}{f_0(A,\theta,\omega,R)}.$$\end{rem}
 
 \begin{rem}\label{R12}
 Let $X\subset \mathcal X_{\ord}$ be, as before, an open subscheme with $\overline{X}$ connected.
By the proof of Theorem \ref{T5} it will follow that for each $w\in W$, every $f\in \mathbb I_{g,\ord}(w)$ vanishes 
 on the set of ordinary $\CM_R$-points of $X$, i.e., of $R$-valued points of $X$ that correspond to ordinary abelian schemes over $\Spec(R)$ that have {\it complex multiplication} (see Remark \ref{R23} (b)). 
 If $p>2$, then an ordinary $\CM_R$-point is the same as a $\CL_R$ point, i.e., an $R$-valued point of $\mathcal X_{\ord}$ that corresponds to an ordinary abelian scheme over $\Spec(R)$ which is the canonical lift of its reduction modulo $p$ (see \cite{V}, Thm. 1.6.3). Similarly, if $p=2$, then an ordinary $\CM_R$-point is the same as a $\QCL_R$ point, i.e., an $R$-valued point of $\mathcal X_{\ord}$ that correspond to an ordinary abelian scheme over $\Spec(R)$ which is a {\it quasi-canonical lift} of its reduction modulo $2$, which means that it defines the same filtered $F$-crystal over $k$ as the canonical lift (see \cite{V}, Thm. 1.6.3).
So for every linear system
 $\f$ in $\mathbb I_{g,\ord}$, the set of $\f$-unstable points $X(R)^{\textup{u}}_{\f}$
 in $X(R)$ contains the set of ordinary $\CM_R$-points of $X$.
 It would also be interesting to decide when, for suitable $\f$s, the locus of $\f$-stable points of $X(R)$ is non-empty and, when so, to compute it.\end{rem}

Combining Theorem \ref{T5} with results from Subsections \ref{S41} and \ref{S42} we will conclude that (see Theorem \ref{T16} (c), (e) and (g)):
 
\begin{thm}\label{T6}
  The following two properties hold:
 \medskip
 
 {\bf (a)} We have equalities:
 $$\asytrdeg(\mathbb I_{g,\ord})=\phidim(\mathbb I_{g,\ord})=\frac{g(g+1)}{2}+1.$$
 
 \smallskip
 
 {\bf (b)} The field $\Frac(\mathbb I_{g,\ord})$ is $\phi$-generated over $K$ by finitely many elements in $\mathbb I^4_{g,\ord}$.
\end{thm}

Here and from now on, $\phidim:=\phidim_R$.

\begin{rem}\label{R13}
As, by classical invariant theory, the $K$-algebras $\mathbb H^{r-1,\torus}_g\otimes_R K$
are finitely generated we also get that
$$\asydim(\mathbb I_{g,\ord})=\frac{g(g+1)}{2}+1.$$
On the other hand we do not know whether $\mathbb I_{g,\ord}\otimes_R K$ is $\phi$-finitely generated as a $K$-algebra.
\end{rem}

\begin{rem}\label{R14}
We recall the $R$-algebra monomorphism
\begin{equation}\label{EQ051}
\mathbb I_{g,\ord}\rightarrow \mathbb M_X\end{equation}
with injective reduction modulo $p$, where $X\subset \mathcal X_{\ord}$ is an arbitrary open subscheme with $\overline{X}$ connected (see Equation (\ref{EQ040})). If $X$ is affine, then by viewing $\mathbb M_X$ as filtered by the subrings $\mathbb M^r_X$, we have (see Equation (\ref{EQ006}))
$$\asydim(\mathbb M_X)=\frac{g(g+1)}{2}+1.$$
The morphism 
(\ref{EQ047}) is induced by the morphism (\ref{EQ051}) and we
recall that we view the morphism 
(\ref{EQ047}) as an avatar of the morphism
(\ref{EQ049}). Then the equality of the asymptotic dimensions of the source and target in (\ref{EQ051})
says that the `quotient map' 
(\ref{EQ049})
is a map between geometric objects `of the same dimension' in `$\delta$-geometry'. 

Note that we have $K$-algebra monomorphisms
\begin{equation}\label{EQ052}
\mathbb I^r_{g,\ord}\otimes_R K\rightarrow \mathbb M^r_X\otimes_R K\end{equation}
induced by (\ref{EQ051}). By Theorem \ref{T16} (b) below and Remark \ref{R13} the $K$-algebra 
$\mathbb I^r_{g,\ord}\otimes_R K$ is finitely generated and has dimension
 \begin{equation}\label{EQ053} 
 \dim(\mathbb I^r_{g,\ord}\otimes_R K)=
 r\frac{g(g+1)}{2}-g^2+1+r\end{equation}
 while, for $X$ affine, the (noetherian but not finitely generated)
$K$-algebra $\mathbb M^r_X\otimes_R K$ has dimension 
\begin{equation}\label{EQ054}
\dim(\mathbb M^r_X\otimes_R K)=
(r+1)\frac{g(g+1)}{2}+r.\end{equation}
The difference between the two dimensions is therefore 
$$\dim(\mathbb M^r_X\otimes_R K)-\dim(\mathbb I^r_{g,\ord}\otimes_R K)=
g^2+\frac{g(g+1)}{2}-1$$ 
which is independent of $r$; note that this difference is less than (rather than equal to) the dimension $\dim(\pmb{\GSp}_{2g,K})=2g^2+g+1$ of the group of symplectic similitudes of rank $g+1$ over $K$, as one might expect
from the analogy with the differential algebraic picture \cite{Bu95b}, Introd. and the complex analytic picture \cite{FC}, Ch. VII, Sect. 3, p. 255.
It would be interesting to study the fibers of the morphism of spectra induced by (\ref{EQ052}). One is especially interested in the structure of the rings
$\mathbb M^r_X/\mathbb I_{g,\ord}^{r,[P]}\mathbb M^r_X$ and
$\mathbb M_X/\mathbb I_{g,\ord}^{[P]}\mathbb M_X$ for $P\in X(R)$,
 where $\mathbb I_{g,\ord}^{r,[P]}\subset \mathbb I^r_{g,\ord}$ and
 $\mathbb I_{g,\ord}^{[P]}\subset \mathbb I_{g,\ord}$
 are the prime ideals (see Remark \ref{R11})
$$\mathbb I_{g,\ord}^{r,[P]}:=\bigoplus_{w\in W} \mathbb I_{g,\ord}^{r,[P]}(w),\ \ \ \mathbb I_{g,\ord}^{[P]}:=\bigcup_{r\geq 1} \mathbb I_{g,\ord}^{r,[P]}.$$
The ring $\mathbb M_X/\mathbb I_{g,\ord}^{[P]}\mathbb M_X$ can be viewed as an 
 avatar, in $\delta$-geometry, of the polarized isogeny class of $P$. 
 It has a natural structure of filtered $W$-graded ring equipped with an endomorphism induced by $\phi$.
 Note that, by Remark \ref{R11}, the ring 
 $\mathbb I_{g,\ord}/\mathbb I_{g,\ord}^{[P]}$ has a natural structure of filtered $W$-graded integral domain all of whose homogeneous components are either $0$ or (non-canonically) isomorphic to $R$ and is equipped with an endomorphism induced by $\phi$.
 One would like to know, for instance, if or when $\mathbb M_X/\mathbb I_{g,\ord}^{[P]}\mathbb M_X$ is an integral domain and one would like to compute its asymptotic dimension.
A related problem is to understand the fibers of the maps $\Phi_{\bf f}$ in (\ref{EQ046}). In the same vein one would like to understand the image of the maps 
$\Phi_{\bf f}$, e.g., the ideal of all $\delta$-homogeneous polynomials in $\mathbb S_n$ that vanish on this image. Another related open problem is to compute the transcendence degrees 
 $\textup{tr.deg}_k(\mathbb I_{g,\ord}^r\otimes_R k)$.

\end{rem}

In Section \ref{S3} we will also obtain estimates/asymptotics for the (finite) ranks of the $R$-modules 
 $\mathbb I^r_{g,\ord}(w)$ (see Theorem \ref{T15} (a)). Also, for $r\leq 2$ we will obtain $W$-graded $K$-algebra isomorphisms
$$\mathbb I^1_{g,\ord}\otimes_R K \simeq K[\Theta_0,z_0,z_0^{-1}],$$
$$\mathbb I^2_{g,\ord}\otimes_R K \simeq K[\Theta_0,\ldots,\Theta_g,z_0,z_1,\frac{1}{z_0z_1}],$$
 where $\Theta_0,\ldots,\Theta_g,z_0,z_1$ are indeterminates with weights
 $-2,\ldots,-2,\phi-1,\phi^2-\phi$ (respectively); these isomorphisms will be entirely explicit (see Corollary \ref{C16}). The structure of the rings
 $\mathbb I^r_{g,\ord}$ in case $r\geq 3$ will turn out to be rather complicated.
 However the rings $\mathbb I^r_{g,\ord}\otimes_R K$ are 
 isomorphic to the rings $\mathbb H_{g,\tot}^{r-1,\torus}\otimes_R K$ (see Theorem \ref{T5})
 and the rings $\mathbb H_{g,\tot}^{r-1}\otimes_R K$ will be studied in detail in Subsections \ref{S41} to \ref{S43}.
 
 Recalling the notation in Remark \ref{R2}, we get, for $g$ odd, an equality of $W_+$-modules 
 $$W_{\mathbb I_{g,\ord}}=\{w\in W|\deg(w)\in 2\mathbb Z,\ w\leq 0\}.$$
 The structure of $W_{\mathbb I_{g,\ord}}$ for $g$ even is more complicated. However 
 in view of the isomorphism $\mathbb I_{g,\ord}\otimes_R K\simeq \mathbb H_{g,\tot}^{\torus}\otimes_R K$ (see Theorem \ref{T5}) the problem of understanding the structure of 
 $W_{\mathbb I_{g,\ord}}$ trivially reduces to the problem of determining which homogeneous components of $\mathbb H_{g,\tot}^{\torus}$ (with respect to a natural gradation) are non-zero.
 This latter problem immediately 
 reduces to a similar problem for the ring  $\mathbb H_{g,\tot}$ which, in its turn, 
 is partially addressed (although not completely solved) in Section 4;  see Lemma \ref{L43} (a)
 and  (b) and Subsection \ref{S4206}.
 
 \subsubsection{Structure of the ring $\mathbb I_g$}

 Next we turn our attention to $\mathbb I_g$.
 We will prove that $\mathbb I_g$ has an asymptotic transcendence degree that is equal to that of $\mathbb I_{g,\textup{ord}}$ so the two rings are `close' in the sense of asymptotic transcendence degrees. Indeed, using the invariant theory of multiple endomorphisms (see Subsubsections \ref{S411}, \ref{S415} and \ref{S416}), we will prove (see Theorem \ref{T17}) the following:
 
 \begin{thm}\label{T7}
 For each integer $r\geq 4$ there exists a natural map of $W$-graded $R$-algebras
$$\square:\mathbb H^{r-2}_{g,\con}[x_1,\ldots,x_r,\xi_1,\ldots,\xi_r]\rightarrow \mathbb I^r_g$$
whose image satisfies
$$\trdeg_R(\textup{Im}(\square))\geq \trdeg_R(\mathbb I^r_{g,\ord})-\frac{g(g+1)}{2}-2.$$\end{thm}

 Here $x_i$ and $\xi_i$ are indeterminates of degree $-\phi^{i-1}-\phi^i$ and $-1-\phi^i$ (respectively)
 and the maps $\square$ are compatible as $r$ varies. 
 
 \begin{rem}\label{R15}
 Theorem \ref{T7} and Remark \ref{R3}
 imply the inequality
 \begin{equation}\label{EQ055}
 \trdeg_{\mathbb I_g} (\mathbb I_{g,\ord})\leq \frac{g(g+1)}{2}+2.\end{equation}
 It would be interesting to compute the exact value of the transcendence degree $\trdeg_{\mathbb I_g} (\mathbb I_{g,\ord})$. We emphasize  that we do not even know if the fields 
 $\Frac(\mathbb I_g)\subset \Frac(\mathbb I_{g,\ord})$ are distinct; this question is open even for $g=1$. On the other hand,
 by Theorem \ref{T7} and Equations (\ref{EQ053}) and (\ref{EQ054}) we get
 \begin{equation}\label{EQ056}
 \dim(\mathbb M_X^r\otimes_R K)-\trdeg_R(\mathbb I_g^r)\leq 2g^2+g+1=\dim(\pmb{\GSp}_{2g/K}).
 \end{equation}
 The value of $\trdeg_R(\mathbb I_g^r)$ is not known even for $g=1$.
 The analogies with the differential algebraic picture \cite{Bu95b}, Introd. and with the complex analytic picture \cite{FC}, Ch. VII, Sect. 3, p. 255, suggest that  Inequality (\ref{EQ056}) might be an equality. A related open problem is to compute the transcendence degrees $\textup{tr.deg}_k(\mathbb I_g^r\otimes_R k)$; again, this is not known even for $g=1$.

 \end{rem}
 
 Using Inequality (\ref{EQ055}) and Theorem \ref{T6} (a) we will deduce (see Corollary \ref{C18} (c)) the following result which 
 solves, in particular, a `main open problem' in \cite{BB} (see \cite{BB}, paragraph after Thm. 1.13):

 \begin{thm}\label{T8}
 We have
$$\asytrdeg(\mathbb I_g)=\phidim(\mathbb I_g)=\frac{g(g+1)}{2}+1.$$
 \end{thm}
 
 \begin{rem}\label{R16}
 As $\Frac(\mathbb I_{g,\ord})$ is $\phi$-finitely generated over $K$ so is its subfield 
 $\Frac(\mathbb I_g)$ (see Remark \ref{R1} (c)). Also, as $\mathbb I^r_g(w)\subset \mathbb I^r_{g,\ord}(w)$, it follows that for all linear systems $\f$ of weight $w$ in  
 $\mathbb I^r_{g,\ord}$
 the set of $\f$-unstable points of $\mathcal X(R)=\mathcal A_{g,1,N}(R)$ contains the set of all ordinary $\CM_R$-points. 
 Here are some questions that we leave open. Are $\mathbb I^r_g\otimes_RK$ finitely generated $K$-algebras?
 Is $\mathbb I_g\otimes_R K$ a  $\phi$-finitely generated  $K$-algebra?
What are the transcendence degrees of $\mathbb I^r_g\otimes_RK$ over $K$?
What is the structure of the rings $\mathbb M_X/\mathbb I_g^{[P]}\mathbb M_X$ where 
 $X$ is open in $\mathcal X$, with 
$\overline{X}$ connected, $P\in X(R)$, and 
$\mathbb I_g^{[P]}$ is the prime ideal of $\mathbb I_g$ defined as in Remark \ref{R14} by replacing $\mathbb I_{g,\ord}$ with $\mathbb I_g$? \end{rem}
 
 \begin{rem}\label{R17}
It would be interesting to compute the $W_+$-module $W_{\mathbb I_g}$, see Remark \ref{R2}.
As $\det(f^r)\in \mathbb I^r_g(-1-\phi^r)$ we have $-1-\phi^r\in W_{\mathbb I_g}$ for all integers  $r\geq 1$. So $W_{\mathbb I_g}$ contains the $\mathbb Z$-linear span $W_{\textup{cyc}}$ of the elements of the form $-\phi^a-\phi^b$ with $0\leq a<b$ integers. Note that $W_{\textup{cyc}}$ is a $W_+$-submodule of $W_{\mathbb I_g}$.  All examples of 
homogeneous elements in $\mathbb I_g$ 
 appearing in our paper have weights in $W_{\textup{cyc}}$.  We do not know if or when $W_{\mathbb I_g}=W_{\textup{cyc}}$. On the other hand we claim that the $W_+$-module $W_{\textup{cyc}}$ is not finitely generated. 
Indeed, if $w_1,\ldots,w_n$ generate the $W_+$-module $W_{\textup{cyc}}$, then for every integer $r\geq 1$ we can write $-1-\phi^r=\sum_{i=1}^n v_{r,i}w_i$ with $v_{r,1},\ldots,v_{r,n}\in W_+$. Taking degrees and using the fact that $\deg(w_i)\in-2\mathbb N$ we get that $-1-\phi^r=w_{i(r)}$ for some $i(r)\in \{1,\ldots,n\}$ which for $r>>0$ represents a contradiction, proving our claim. As a consequence of the claim and of Remark \ref{R2}, if for some $g\in\mathbb N$ we have
 $W_{\mathbb I_g}=W_{\textup{cyc}}$, then 
 $\mathbb I_g\otimes_R K$ is not $\phi$-finitely generated as a $K$-algebra.
 \end{rem}
 
For the next result, for a polarized abelian scheme $(A,\theta)$ over $R$, we denote by $\End(A,\theta)^{(p)}$ the subring of $\End(A)$ generated by all polarized isogenies $A\rightarrow A$. The theorem below shows that the locus in $\mathcal A_{g,1,N}(R)$ of all $(A,\theta)$s for which $\End(A,\theta)^{(p)}$ is non-commutative is contained in the vanishing locus of a non-zero Hecke covariant Siegel $\delta$-modular form (see Theorem \ref{T18} below for a more precise result):
 
 \begin{thm}\label{T9}
 There exists $w\in W$ with the property that $\mathbb I_g(w)$ contains a non-zero form
 $\mathcal D$
 such that for all
 $(A,\theta,\omega)\in \pmb{\M}_g(R)$ with $\mathcal D(A,\theta,\omega)\neq 0$,
 the ring $\End(A,\theta)^{(p)}$ is commutative.
 \end{thm}

 Let us say that the polarized isogeny class of a principally polarized abelian scheme $(A,\theta)$ over $R$ is {\it decomposable} if there exist non-zero principally polarized abelian schemes $(A_1,\theta_1)$ and $(A_2,\theta_2)$ and a polarized isogeny $A\rightarrow A_1\times A_2$, where $A_1\times A_2$ is equipped with the polarization $\theta_1\times \theta_2$. 
 We say that the polarized isogeny class of $(A,\theta)$ is {\it indecomposable} if it is not decomposable.
 The theorem below shows that the locus in $\mathcal A_{g,1,N}(R)$ of all $(A,\theta)$s whose polarized isogeny class are decomposable is contained in the vanishing locus of a non-zero Hecke covariant Siegel $\delta$-modular form (see Theorem \ref{T18.5} below for a more precise result):
 
 \begin{thm}\label{T10}
 There exists $w\in W$ with the property that $\mathbb I_g(w)$ contains a non-zero form
 $\mathcal G$
 such that for all
 $(A,\theta,\omega)\in \pmb{\M}_g(R)$ with $\mathcal G(A,\theta,\omega)\neq 0$, the polarized isogeny class of $(A,\theta)$ is indecomposable.
 \end{thm}
 
 Theorems \ref{T9} and \ref{T10} are arithmetic analogues of results in differential algebraic geometry, see \cite{Bu95b}, Prop. 6.8 and Cor. 5.11.
 
The following `$p$-adic approximation result' (see Theorem \ref{T19} (c) below) shows that $\mathbb I_g$ is `close to' $\mathbb I_{g,\ord}$ in yet another way:

\begin{thm}\label{T11}
Let $\mathbb D_g$ be the multiplicative set of all
homogeneous elements of $\mathbb I_g\backslash p \mathbb I_g$ whose weights have even degree. 
Then the natural homomorphism
$$\reallywidehat{\mathbb D_g^{-1} \mathbb I_g}
\rightarrow \reallywidehat{\mathbb D_g^{-1}\mathbb I_{g,\ord}}$$
is an isomorphism.
\end{thm}

The above theorem implies that $\mathbb I_g$ is `$p$-adically sufficiently close to' $\mathbb I_{g,\ord}$ in the sense that each element in $\mathbb I_{g,\ord}$
can be $p$-adically approximated by fractions with denominators homogeneous elements of $\mathbb I_g\setminus p\mathbb I_g$. 

\subsubsection{Flows and reciprocity}

Our next application involves the notion of $\delta$-flow introduced in \cite{Bu17}, Ch. 3, Sect. 3.6, Def. 3.122. 
Let $X$ be an affine smooth scheme over $\Spec(R)$ and let $J^1(X)$ be its first $p$-jet space.
We recall that an ideal $\mathfrak a$ in $\mathcal O(J^1(X))$ is called a {\it $\delta$-flow} on $X$ if the composite projection
$$\Spf(\mathcal O(J^1(X))/\mathfrak a) \subset J^1(X)\rightarrow J^0(X)=\widehat{X}$$
is an isomorphism. 
Our terminology here is suggested by the 
analogies between $p$-derivations and derivations (vector fields) 
that are in the background of \cite{Bu95a,Bu05}. Note that by the universal property of $J^1(X)$
every $\delta$-flow on $X$ induces a Frobenius lift on $\widehat{X}$, i.e.,  an endomorphism of $\widehat{X}$ over $\mathbb Z_p$  whose reduction modulo $p$ is  the $p$-power Frobenius.
We will prove the following result (see Theorem \ref{T20}):

\begin{thm}\label{T12}
 Let $X$ be an affine open subscheme of the moduli scheme $\mathcal X=\mathcal A_{g,1,N,R}$ over $\Spec(R)$, where $N\geq 3$ is an integer prime to $p$, such that $\overline{X}$ is contained in the ordinary locus and there exists a column basis $\omega_X$ of $1$-forms on $A_X$. Let $\mathfrak a_X$ be the ideal of $\mathcal O(J^1(X))$ locally 
 generated by the entries of $f^1_{g,\crys}(A_X,\theta_X,\omega_X)$ where $(A_X,\theta_X)$ is the universal principally polarized abelian scheme over $X$.
Then $\mathfrak a_X$ is a $\delta$-flow on $X$.
\end{thm}

\begin{rem}\label{R17.1}
It would be interesting to check whether
 the Frobenius lift on $\widehat{\mathcal A_{g,1,N,\ord}}$ induced by the $\delta$-flow in Theorem \ref{T12} coincides with the Frobenius lift 
considered in \cite{BG}, p. 3 
 that arises from the functor that attaches to every ordinary abelian scheme its quotient by the
canonical subgroup. If this were the case, then we would get the interesting consequence that the Frobenius lift in \cite{BG}, although not extendable to the whole
of $\widehat{\mathcal A_{g,1,N}}$, `comes from a global object ($f^1_{\crys}$) defined on the whole of $\mathcal A_{g,1,N}$.' \end{rem}

Using Theorem \ref{T12} and the arithmetic analogue in \cite{Bu97}, Cor. 1.7 of Manin's Theorem of the Kernel  \cite{Man}, Sect. 5, Thm. 2 we will prove the following result (see Theorem \ref{T20.5} below) which extends the `Reciprocity Theorem for $\CL_R$ points' in \cite{BP}, Thm. 3.5, to the case of arbitrary genus $g\geq 1$. We refer to the discussion before Theorem \ref{T20.5} for more details.

\begin{thm}\label{T12.12}
Assume we are given morphisms of schemes over $R$
$$\Pi:Y\rightarrow X,\ \ \ \ \Phi:Y\rightarrow A$$ 
with $X\subset \mathcal X_{\ord}$ affine, $Y$ affine, $\Pi$ \'{e}tale, and $A$ an abelian scheme  of relative dimension $d$.
Assume $Y,X,A,\Pi,\Phi$ descend to corresponding objects $Y',X',A',\Pi',\Phi'$ defined over a finite normal extension $R'$ of $\mathbb Z_p$ contained in $R$ and assume $A'$
  is non-degenerate in the sense of Definition \ref{df27} below.
Then there exists a morphism of $p$-adic formal schemes 
$$\Phi^+:\widehat{Y}\rightarrow \widehat{\mathbb G_a^d}$$
 that satisfyies the following property. For all $P_1,\ldots,P_n\in \Pi^{-1}(\CL_R)\subset Y(R)$ and all $m_1,\ldots,m_n\in \mathbb Z$ we have
$$
\sum_{i=1}^n m_i \Phi(P_i)\in A(R)_{\textup{tors}} \ \ \ \ \Longleftrightarrow \ \ \ \ 
\sum_{i=1}^n m_i \Phi^+(P_i)=0\in R^d.
$$
\end{thm}

\begin{rem}\label{R17.5}
As a consequence of Theorem \ref{T12.12} we will have (see Corollary \ref{C18.5}) 
that if $\Pi^{-1}(\CL_R)\not\subset \Phi^{-1}(A(R)_{\textup{tors}})$, then the image of the map
$$\Pi^{-1}(\CL_R)\cap \Phi^{-1}(A(R)_{\textup{tors}})\subset Y(R)\rightarrow Y(k)$$
is not Zariski dense in $\overline{Y}$. This statement is a (local) higher dimensional analogue
of the result in \cite{NS}, Thm. 1.5 according to which there exist only finitely many torsion Heegner points on each elliptic curve over $\mathbb Q$.
\end{rem}

\subsubsection{Cohomological computations}

 It would be interesting to explore the cohomology modules of the `main' ssa-$\delta$-spaces 
 (\ref{EQ044}) appearing in our theory.
 We will leave this problem open but, with notation as in Equation (\ref{EQ044}), we will present some easy computations that illustrate the flavor of the problem. 
 
 First, consider the cohomology groups $H^d(\mathbb P^1_{\delta})$.
 We will check (see Proposition \ref{P8}) that the following equalities hold:
 $$H^0(\mathbb P^1_{\delta})=R,\ \ \ \phirank_R(H^1(\mathbb P^1_{\delta}))=\infty.$$
This is in stark contrast with the situation in usual algebraic geometry where $H^1(\mathbb P^1,\mathcal O)=0$ for the projective line $\mathbb P^1$ over a ring. 
 
 Similarly, let us consider the case $g=1$ and the basic forms $f^r:=f^r_{1,\crys}\in \mathbb \mathbb I^r_{11}(\phi^r,1)=I^r_1(-1-\phi^r)$ and $f^{\partial}_{1,\crys}
 \in \mathbb I^1_{11,\ord}(-\phi,1)=\mathbb I^1_{1,\ord}(\phi-1)$ referred to above.
 Consider also the forms
 \begin{equation}\label{EQ057}
 f^{\langle r \rangle}:=(f^1)^{\phi^{r-1}} \cdot (f^{\partial})^{\frac{\phi^r+\phi^{r-1}-2}{\phi-1}}\in \mathbb I_{1,\ord}(-2);
 \end{equation}
 these forms (and their generalizations for arbitrary $g$) will play a key role in the paper.
For each integer $n\geq 2$, consider the ssa-$\delta$-space 
 $$Z_{\delta,n}:=\Proj_{\delta}(\mathbb I_{1,\ord},(f^{\langle 1 \rangle},\ldots,f^{\langle n \rangle})).$$
 We will then check (see Proposition \ref{P9}) that the following equalities hold:
 $$H^0(Z_{\delta,n})=\mathbb I_{1,\ord,\langle f^{\langle n \rangle}\rangle},\ \ \ H^d(Z_{\delta,n})=0\ \ \textup{for}\ \ d\geq 1.$$
 
 On the other hand, if we consider a matrix $\lambda=(\lambda_{ij})_{1\leq i,j\leq 2}\in \pmb{\GL}_2(R)$
 and we set $(f_1^{\lambda},f_2^{\lambda}):=(f^{\langle 1 \rangle},f^{\langle 2\rangle})\lambda^{\t}$ and
 $$Z_{\delta,2}^{\lambda}:=\Proj_{\delta}(\mathbb I_{1,\ord},(f^{\lambda}_1,f^{\lambda}_2)),$$
 then we will show (see Proposition \ref{P10}) that if each $\lambda_{ij}\neq 0$, then the following equalities hold: 
 $$H^0(Z_{\delta,2}^{\lambda})=R,\ \ \ \phirank_R(H^1(Z_{\delta,2}^{\lambda}))=\infty.$$
 
 Note that, if $\lambda=1_2$, then we have
 $Z_{\delta,2}^{\lambda}=Z_{\delta,2}$.
 So the cohomology of $Z_{\delta,2}^{\lambda}$ effectively depends on $\lambda$.
 This is, again, in stark contrast with 
 the situation in usual algebraic geometry where the {\v C}ech cohomology groups of an affine open cover of a separated scheme do not depend on the cover.
 Indeed, for each $\lambda$ the pair $(f_1^{\lambda},f_2^{\lambda})$ can be intuitively 
 viewed as defining an `affine open cover' of a `separated' ssa-$\delta$-space
 that does not depend on $\lambda$ (as the $R$-module spanned by $f_1^{\lambda},f_2^{\lambda}$ does not depend on $\lambda$). 
 
 Finally, for arbitrary $g$, we will prove (see Proposition \ref{P11}) that for a `generic choice' of a linear system $\f$ of weight $-2$ in 
 $\mathbb I_{g,\ord}$ we have $H^0(\Proj_{\delta}(\mathbb I_{g,\ord},\f))=R$.
 For a weaker version of this statement for $\mathbb I_g$ see Proposition \ref{P12} below.
 
\newpage\section{Proofs of the main results}\label{S3}

All the main results are proved in this section.
 
Subsection \ref{S31} records some basic symmetry properties of the bilinear maps over $\mathbb Z_p$ appearing in the Serre--Tate deformation theory of abelian schemes. 
Subsection \ref{S32} introduces and studies some multiplicative monoids contained in endomorphism algebras of abelian varieties over $k$;  these monoids, and some related unitary groups, will play a key role in our arguments. Subsection \ref{S33} translates 
the conditions of lifting from $k$ to $R$ suitable 
isogenies between abelian schemes into the language of matrices and relates
 them to the previously introduced unitary groups.

Subsection \ref{S34} reviews a construction of Hecke covariant Siegel $\delta$-modular forms $f^r_{g,\crys}$ introduced in \cite{BB}, Subsect. 4.1; it also introduces and studies the ordinary forms $f^{\partial}_{g,\crys}$ which generalize the ordinary form $f^{\partial}_{1,\crys}$ introduced and studied in \cite{Ba} and \cite{Bu05}, Ch. 8. Out of these forms we construct some new ordinary forms $f^{\langle r\rangle}$ of weight $-2$ that play a key role later. 

In Subsection \ref{S35} we introduce one of our main tools, the Serre--Tate expansion of a Siegel $\delta$-modular form, and we prove a {\it Serre--Tate expansion principle} (see Proposition \ref{P3} (b)). In Subsection \ref{S36} we use the Serre--Tate expansions to construct embeddings (comparison maps) $\heartsuit$ of the $R$-modules $\mathbb I^r_{g,\ord}(w)$ of Hecke covariant forms into $R$-modules $\mathbb H^r_g(s)$
of $\pmb{\SL}_g$-invariant functions on certain $R$-modules of multiple quadratic forms (see Corollary \ref{C7}). On the other hand we use the forms $f^{\langle r\rangle}$ to construct embeddings (comparison maps) $\diamondsuit$ of the $R$-modules $\mathbb H^{r-1}_g(s)$ into $R$-modules $\mathbb I^r_{g,\ord}(-2s)$. A key Lemma \ref{L17} on `elimination of logarithms' will help us use the embeddings $\heartsuit$ and $\diamondsuit$ to provide a complete description, up to torsion, of the ring $\mathbb I_{g,\ord}$ in terms of the ring of invariants $\mathbb H_g$ (see Theorem \ref{T16} (a)). 

In Subsection \ref{S37} we study the $R$-algebra $\mathbb I_g$ of (non-ordinary) Hecke covariant forms. The main idea is to use  the forms $f^r_{g,\crys}$ to find embeddings (comparison maps) $\square$ of the $R$-algebras $\mathbb H^r_{g,\textup{conj}}$ of $\pmb{\SL}_g$-invariant functions on $R$-modules of multiple endomorphisms into the ring $\mathbb I^r_g$ (see Theorem \ref{T17}). 

Subsection \ref{S38} contains the proof of our main $p$-adic approximation result Theorem \ref{T11}. In Subsection \ref{S39} we prove Theorem \ref{T12} on the $\delta$-flow structure defined by $f^1_{g,\crys}$. 
In Subsection \ref{S310.5} we prove the `reciprocity' Theorem \ref{T12.12}.
Subsection \ref{S310.66} contains a few cohomological computations. 

\subsection{Symmetry}\label{S31}

Let $A_0$ be an ordinary abelian variety over $\Spec(k)$ of dimension $g$. Let ${\mathcal M}_{A_0}=\Spf(R_g)$ be the formal deformation space of $A_0$; it is well-known that we can identify $$R_g=R[[T_{ij}|1\leq i,j\leq g]],$$ where $T_{ij}$s are algebraically independent indeterminates over $R$. 

For an abelian scheme $A$ over $\Spec(R)$ which lifts $A_0$ (thus we have $A\times_{\Spec(R)} \Spec(k)=A_0$), let 
$$q_A:T_p(A_0)\times T_p(\check{A}_0)\rightarrow 1+\mathfrak m_R=1+pR$$
be the Serre--Tate bilinear map attached to $A$, where $1+\mathfrak m_R$ is a group under multiplication. Let $e=[e_1\cdots e_g]^{\t}$ be a column basis of $T_p(A_0)$ by which we mean that its entries form a $\mathbb Z_p$-basis of $T_p(A_0)$.

We assume that there exists a principal polarization on $A$ defined by an isomorphism $\theta:A\rightarrow \check{A}$. Let $\theta_0:A_0\rightarrow \check{A}_0$ be the isomorphism obtained from $\theta$ by reduction modulo $p$. We consider the bilinear map
$$q_{A,\theta}:=q_A\circ [T_p(1_{A_0})\times T_p(\theta_0)]:T_p(A_0)\times T_p(A_0)\rightarrow 1+\mathfrak m_R$$
induced by $q_A$ (where $1_{A_0}$ is the identity of $A_0$) and the isomorphism 
$$T_p(\theta_0):T_p(A_0)\rightarrow T_p(\check{A}_0).$$ 

Let $q_{\univ,A_0}:T_p(A_0)\times T_p(\check{A}_0)\rightarrow 1+\mathfrak m_{R_g}$
be the Serre--Tate bilinear map attached to the formal deformation of $A_0$ over ${\mathcal M}_{A_0}=\Spf(R_g)$. We similarly get a universal Serre--Tate bilinear map
$$q_{\univ,\theta_0}:=q_{\univ,A_0}\circ [T_p(1_{A_0})\times T_p(\theta_0)]:T_p(A_0)\times T_p(A_0)\rightarrow 1+\mathfrak m_{R_g}.$$
Note that $q_A$ and $q_{A,\theta}$ are composites of $q_{\univ,A_0}$ and $q_{\univ,\theta_0}$ (respectively) with a uniquely determined epimorphism $1+\mathfrak m_{R_g}\rightarrow 1+\mathfrak m_R$ of multiplicative groups induced by an $R$-algebra retraction $R_g\rightarrow R$. 

Let $I_{\theta_0}$ be the ideal of $R_g$ such that $\Spf(R_g/I_{\theta_0})$ is the formal deformation space ${\mathcal M}_{A_0,\theta_0}$ of $(A_0,\theta_0)$. 

\begin{prop}\label{P1}
The bilinear map $q_{A,\theta}$ is symmetric.
\end{prop}

\noindent
{\it Proof.} 
For $(\alpha,\beta)\in T_p(A_0)\times T_p(\check{A}_0)$ we have $q_{\univ,\check{A}_0}(\beta,\alpha)=q_{\univ,A_0}(\alpha,\beta)$. This follows from Cartier duality and duality of abelian schemes, which in particular allow us to identify $A_0=\check{\check{A}}_0$ and $\Spf(R_g)={\mathcal M}_{A_0}={\mathcal M}_{\check{A}_0}$. More precisely, if $\Re$ is a local artinian $R$-algebra of residue field $k$, then a short exact sequence
$$0\rightarrow \pmb{\mu}_{\pmb{p^{\infty}},\Re}\rightarrow D_\Re\rightarrow ({\mathbb Q}_p/{\mathbb Z}_p)_\Re\rightarrow 0$$ of $p$-divisible groups over $\Spec(\Re)$ is automatically self-dual with respect to the Cartier duality and this implies the identity $q_{\univ,\check{A}_0}(\beta,\alpha)=q_{\univ,A_0}(\alpha,\beta)$.

The isomorphism $\theta_0:A_0\rightarrow\check{A}_0$ induces an isomorphism 
$$\iota_{\theta_0}:\Spf(R_g)={\mathcal M}_{A_0}\simeq \Spf(R_g)={\mathcal M}_{\check{A}_0},$$ and the ideal $I_{\theta_0}$ of $R_g$ is the smallest one so that $\iota_{\theta_0}$ modulo $I_{\theta_0}$ is the identity automorphism of $R_g/I_{\theta_0}$. In other words, the formal deformation space ${\mathcal M}_{A_0,\theta_0}$ of $(A_0,\theta_0)$ is the locus where the two maps 
$$q_{\univ,A_0},q_{\univ,\check{A}_0}\circ (T_p(\theta_0)\times T_p(\check\theta_0)^{-1}):T_p(A_0)\times T_p(\check{A}_0)\rightarrow 1+\mathfrak m_{R_g}$$ 
coincide (cf. also \cite{Ka81}, Thm. 2.1 4)). 

From the very definition of a polarization we have the symmetry property $\theta_0=\check{\theta}_0$. From this and the identity $q_{\univ,\check{A}_0}(\beta,\alpha)=q_{\univ,A_0}(\alpha,\beta)$, by composing the last two maps with the isomorphism 
$$T_p(1_{A_0})\times T_p(\theta_0):T_p(A_0)\times T_p(A_0)\rightarrow T_p(A_0)\times T_p(\check{A}_0)$$ 
we get that the resulting two maps from $T_p(A_0)\times T_p(A_0)$ to $1+\mathfrak m_{R_g}$
differ from each other via the composite with the switching involution of $T_p(A_0)\times T_p(A_0)$ defined by $(x,y)$ goes to $(y,x)$. This implies that
\begin{equation}\label{EQ058}
I_{\theta_0}=(q_{\univ,\theta_0}(e_i,e_j)-q_{\univ,\theta_0}(e_j,e_i)|1\leq i,j\leq g).
\end{equation}
Thus the symmetry of reductions of $q_{\univ,\theta_0}$ via quotients of $R_g$ holds precisely for $R_g/I_{\theta_0}$ (i.e., for ${\mathcal M}_{A_0,\theta_0}$) and hence it holds for the bilinear map $q_{A,\theta}$ which is defined via the $R$-valued point of ${\mathcal M}_{A_0,\theta_0}$ that corresponds to $(A,\theta)$.\endproof 

\medskip

The entries of $q_{A,\theta}$ satisfy certain `linear relations'. More precisely, the formal closed subscheme ${\mathcal M}_{A_0,\theta_0}=\Spf(R_g/I_{\theta_0})$ of ${\mathcal M}_{A_0}=\Spf(R_g)$ (i.e., the formal completion of the moduli stack $\mathcal A_{g,R}$ at the point corresponding to $(A_0,\theta_0)$) is a formal subtorus of ${\mathcal M}_{A_0}=\Spf(R_g)$ and one gets $g(g-1)/2$ `linear relations' for the entries of $q_{A,\theta}$ which can be read out from $\theta_0$ as follows. First, $\check{e}=[\check{e}_1\cdots \check{e}_g]^{\t}:=[T_p(\theta_0)(e_1)\cdots T_p(\theta_0)(e_g)]^{\t}$ is a column basis of $T_p(\check{A}_0)$ and we get the $T_{ij}$s via the rule $T_{ij}=q_{ij}-1=q_{\univ,\theta_0}(e_i,e_j)-1$. We have $q_{\univ,\theta_0}(e_i,e_j)-q_{\univ,\theta_0}(e_j,e_i)=T_{ij}-T_{ji}$ for all $1\leq i<j\leq g$, and, in view of Equation (\ref{EQ058}), we get:

\begin{cor}\label{C1}
We have ${\mathcal M}_{A_0,\theta_0}=\Spf(S^0_{\for})$, where
$S^0_{\for}=R[[T_{ij}|[1\leq i \leq j\leq g]]$ is identified canonically with $R_g/(T_{ij}-T_{ji}|1\leq i < j\leq g)$.
\end{cor}

If $f$ is an arbitrary column basis of $T_p(\check{A}_0)$, then the `linear relations' would involve the matrix representation of $T_p(\theta_0)$ with respect to $e$ and $f$.

\subsection{Unitary groups}\label{S32}

For the ordinary principally polarized abelian variety $(A_0,\theta_0)$ over $k$ we consider the monoid
$$\mathcal U_0:={\mathcal U}_{A_0,\theta_0}:=\{u_0\in\End(A_0)|u_0u_0^{\dagger}\in\mathbb N\setminus p\mathbb N\},$$
where $u_0\rightarrow u_0^{\dagger}:=\theta_0^{-1}\circ \check{u}_0\circ\theta_0$ is the rule that defines the Rosati involution. With $e=[e_1\cdots e_g]^{\t}$ as above, let 
$$U_0=U_{A_0,\theta_0,e}\subset\pmb{\GL}_g(\mathbb Z_p)$$
 be the image of $\mathcal U_0$ under the injective algebra homomorphism
\begin{equation}\label{EQ059}
\End(A_0)\rightarrow \End(T_p(A_0))\simeq \pmb{\Mat}_g({\mathbb Z}_p),\ \ 
v_0\mapsto M_{v_0},\end{equation}
where $\End(A_0)\rightarrow \End(T_p(A_0))$ maps $v_0\in\End(A_0)$ to $T_p(v_0)$ and where the isomorphism in (\ref{EQ059}) is defined by $e$, so we have $T_p(v_0)(e)=M_{v_0}e$.  
Then 
$U_0$ is equipped with an involution $M\mapsto M^{\dagger}$ induced by the Rosati involution, so $M_{v_0^{\dagger}}=M_{v_0}^{\dagger}$. For each $M=(m_{ij})_{1\leq i,j\leq g}\in U_0$ we write $M^{\dagger}=(m^{\dagger}_{ij})_{1\leq i,j\leq g}$. The involution on $U_0$ does not extend, in general, to an involution on $\pmb{\Mat}_g({\mathbb Z}_p)$, and is a priori unrelated to the transposition
on $\pmb{\Mat}_g({\mathbb Z}_p)$. The subgroup of $\pmb{\GL}_g(\mathbb Z_p)$ generated by $U_0$ is:
$$\mathbb Z_{(p)}^{\times}\cdot U_0:=\{\lambda M|\lambda\in\mathbb Z_{(p)}^{\times}\; \textup{and} \;M\in U_0\}\subset \pmb{\GL}_g(\mathbb Z_p).$$

\begin{ex}\label{EX3}
Consider the special case
\begin{equation}\label{EQ060}
\begin{array}{rcl}
A_0 & = &E_0^g,\\
\ & \ & \\
\theta_0 & = & \textup{standard diagonal principal polarization},\\
\ & \ & \\
 e & = & [e_1\cdots e_g]^{\t}=\textup{a diagonal basis};\end{array}
\end{equation}
so $A_0=E_0^g$ is the $g$-fold product of an ordinary elliptic curve $E_0$ over $\Spec(k)$, $\theta_0$ is obtained from (i.e., is the $g$-fold product of) the Abel--Jacobi map of $E_0$, 
we fix a basis $e_0$ of $T_p(E_0)$, 
and each $e_i$  equals $e_0$ but viewed as an element of the $i$-th direct 
summand $T_p(E_0)$ of $T_p(A_0)=T_p(E_0)^g$.
 Let ${\mathcal O}:=\End(E_0)$ and consider the canonical embedding 
 $${\mathcal O}\rightarrow\End(T_p(E_0))=\mathbb Z_p.$$ 
 Recall that $\mathcal O$ is an order in an imaginary quadratic field $\mathcal K$ and that the Rosati involution on 
 $\mathcal O$ is the complex conjugation $\lambda\mapsto \overline{\lambda}$. 
 \end{ex}

Let $M=(m_{ij})_{1\leq i,j\leq g}\in\pmb{\Mat}_g(\mathcal O)\cap\pmb{\GL}_g(\mathbb Z_p)$. Let $$u_0:A_0=E_0^g\rightarrow A_0=E_0^g$$ 
be the isogeny defined naturally by the matrix $M$. We have $u_0(P)=MP$ for all $P\in A_0(k)=E_0(k)^g$, viewed as a column vector with entries in $E_0(k)$. Note that $T_p(u_0)(e)=Me$, hence $M_{u_0}=M$ and $u_0^{\dagger}P=M^{\dagger}P$ for all $P\in A_0(k)$.
 
\begin{lemma}\label{L1}
If $(A_0,\theta_0,e)$ is as in Example \ref{EX3}, then for all $M\in U_0$ we have $M^{\dagger}=\overline{M}^{\t}$.
\end{lemma}

\noindent
{\it Proof.} The case $g=1$ is well-known. Let now $i,j\in\{1,\ldots,g\}$. Let $E_{0,i}$ be $E_0$ but viewed as the $i$-th factor of $A_0=E_0^g$. The $ij$ entry $m_{ij}$ of $M$ defines an endomorphism $m_{ij}:E_{0,j}\rightarrow E_{0,i}$ inducing $\check{m}_{ij}:\check{E}_{0,i}\rightarrow\check{E}_{0,j}$ and thus via the standard principal polarization on $E_0$ inducing an endomorphism $E_{0,i}\rightarrow E_{0,j}$ that defines the $ji$ entry $m^{\dagger}_{ji}$ of $M^{\dagger}$. Thus we have
$$m^{\dagger}_{ji}=\check{m}_{ij}=\overline{m_{ij}},$$ 
where the last equality follows from the case $g=1$ by dropping the indexes $i$ and $j$ in the domains and codomains of $m_{ij}$ and $\check{m}_{ij}$.\endproof

\medskip

An element $\lambda\in\mathcal U_0$ will be called a {\it scalar} if $M_{\lambda}$
is a scalar matrix. If $\lambda$ is a scalar, then so is $\lambda^{\dagger}$. In this case we will view $\lambda,\lambda^{\dagger}$ as elements of $\mathbb Z_p^{\times}$; so also $\lambda/\lambda^{\dagger}\in\mathbb Z_p^{\times}$. Consider the following conditions on $(A_0,\theta_0)$:

\medskip
 {\bf (UAX1)} There exists a scalar $\lambda\in\mathcal U_0$ with $\lambda/\lambda^{\dagger}$ not a root of unity.
 
 \smallskip
 
{\bf (UAX2)} The group $\mathbb Z_p^{\times}\cdot U_0$ is Zariski dense in $\pmb{\GL}_{g,\mathbb Z_p}$, i.e., each polynomial in $g^2$ indeterminates with coefficients in $\mathbb Z_p$ which vanishes on
 $\mathbb Z_p^{\times}\cdot U_0$ is the zero polynomial.
 
\medskip
 
\begin{lemma}\label{L2}
Each $(A_0,\theta_0)$ as in (\ref{EQ060}) satisfies condition \textup{(UAX1)}.
 \end{lemma}
 
\noindent
{\it Proof.} To check this we can assume $g=1$. Based on Lemma \ref{L1} it suffices to
show that there exists a non-zero $\lambda\in\mathcal K$ such that
$\lambda/\overline{\lambda}$ is not a root of unity. But this is
well-known (for instance, see \cite{Bu05}, Ch. 8, Lem. 8.21).\endproof
 
 \begin{prop}\label{P2}
Each $(A_0,\theta_0)$ as in (\ref{EQ060}) satisfies condition \textup{(UAX2)}.
\end{prop}

We first need the following general lemma:

\begin{lemma}\label{L3}
 Let $\mathcal G$ be a connected reductive group scheme over $\Spec(\mathbb Z_{(p)})$. Then $\mathcal G(\mathbb Z_{(p)})$ is Zariski dense in $\mathcal G_{\mathbb Q}$ and is dense in $\mathcal G(\mathbb R)$.
 \end{lemma}
 
\noindent
{\it Proof suggested by Offer Gabber.} Let $\kappa$ be a field equipped with a reasonable field topology: non discrete, Hausdorff, and the ring operations and inversions are continuous. This defines a topology on $X(\kappa)$ for every $\kappa$-scheme locally of finite type, compatible with morphisms, products, open and closed embeddings. Let $\mathcal B$ be a connected smooth affine group over $\Spec(\kappa)$ which is unirational (e.g., this holds if $\kappa$ is perfect, see \cite{Bore}, Ch. V, Thm. 18.2). 

We check that each non-empty open $\mathcal V$ of $\mathcal B(\kappa)$ (in the mentioned topology) is Zariski dense in $\mathcal B$. The unirationality means that there exists a non-empty open $\mathcal Y$ in an affine space $\mathbb A^n_\kappa$ and a dominant morphism $f:\mathcal Y\rightarrow\mathcal B$. By translations on the source and target we can assume that the origin $0$ of $\mathbb A^n_\kappa$ is in $\mathcal Y$ and is mapped by $f$ into the identity of $\mathcal B$. Then the Zariski density of $\mathcal V$ in $\mathcal B$ follows from the observation that non-empty opens of $\kappa^n=\mathbb A^n_\kappa(\kappa)$ are Zariski dense in $\mathbb A^n_\kappa$. 

Lemma follows by taking $\kappa=\mathbb Q$ and $\mathcal B=\mathcal G_{\mathbb Q}$ and using the fact that non-empty opens of $\mathbb Q^n=\mathbb A^n_\kappa(\kappa)$ in the $p$-adic topology are dense in $\mathbb R^n$ (in its metric topology), by the approximation theorem for independent valuations of $\mathbb Q$.\endproof
 
 \medskip
 
 \noindent
 {\it Proof of Proposition \ref{P2}.} Consider an embedding $\mathcal K=\Frac(\mathcal O)\subset \mathbb C$ and note that the group
 $\mathbb Z_{(p)}^{\times}\cdot U_0$ is contained in the subgroup $\pmb{\GL}_g(\mathcal K)$ of $\pmb{\GL}_g(\mathbb Q_p)$.
 By considering a
 $\mathcal K$-embedding of $\mathbb Q_p$ into $\mathbb C$ we get that to prove Proposition \ref{P2} it suffices to show that $\mathbb Z_{(p)}^{\times}\cdot U_0$ is Zariski dense in $\pmb{\GL}_g(\mathbb C)$.
 
Let now $\pmb{\U}_g$ be the reductive (unitary) group scheme over 
 ${\mathbb Z}_{(p)}$ whose functor of valued points attaches to each ${\mathbb Z}_{(p)}$-algebra
 $B$ the group
 $$\pmb{\U}_g(B):=\{M\in\pmb{\GL}_g(B\otimes_{\mathbb Z}{\mathcal O})|M\overline{M}^{\t}=1_g\},$$ 
 where the upper bar conjugation acts on the second factor of $B\otimes_{\mathbb Z}{\mathcal O}$. 
 From Lemma \ref{L3} we get that $\pmb{\U}_g(\mathbb Z_{(p)})$ is dense in $\pmb{\U}_g(\mathbb R)$ in the metric topology. Also from Lemma \ref{L2} we get that $\pmb{\U}_g(\mathbb Z_{(p)})\subset \mathbb Z_{(p)}^{\times}\cdot U_0$.
We conclude that Proposition \ref{P2} holds by the classical fact that
 $\pmb{\U}_g(\mathbb R)$ is dense in $\pmb{\GL}_g({\mathbb C})$ in the Zariski topology. 
 
 We include a proof of this classical fact. Let $\Gamma_g$ be the Zariski closure of $\pmb{\U}_g(\mathbb R)$ in $\pmb{\GL}_{g,\mathbb C}$. Then $\Gamma_g$ is a complex algebraic group and to check that $\Gamma_g=\pmb{\GL}_{g,\mathbb C}$ it is enough to check that its Lie algebra $\Lie(\Gamma_g)$ equals $\Lie(\pmb{\GL}_{g,\mathbb C})$ (see \cite{Bore}, Ch. I, Sect. 7.1). But $\Lie(\Gamma_g)$ contains the $\mathbb C$-span of the Lie algebra $\Lie(\pmb{\U}_g(\mathbb R))$ over $\mathbb R$ and thus it contains all the antisymmetric matrices of size $g\times g$ and all matrices of the form $\sqrt{-1}$ times a symmetric real matrix of size $g\times g$. Hence $\Lie(\Gamma_g)$ contains all the real matrices and thus all the complex matrices of size $g\times g$. We conclude that $\Gamma_g=\pmb{\GL}_{g,\mathbb C}$.
 \endproof
 
\begin{rem}\label{R20}
{\bf (a)} The set of points in the moduli stack ${\mathcal A}_{g,k}$ over $\Spec(k)$ that correspond to ordinary principally polarized abelian varieties over $\Spec(k)$ for which conditions (UAX1) and (UAX2) hold is Zariski dense. This follows from Lemma \ref{L2}, Proposition \ref{P2} and Larsen's example in \cite{Cha}, p. 443.

\smallskip
{\bf (b)} If $\lambda\in\mathcal U_0\setminus\mathbb Q$ is a scalar, then $\mathcal F:=\mathbb Q[\lambda]$ is a totally imaginary quadratic extension of $\mathbb Q$ and the Mumford--Tate group of the generic fiber over $K=\Frac(R)$ of the canonical lift of $A_0$ is naturally identified with the $2$ dimensional torus which is the Weil restriction $\textup{Res}_{\mathcal F/\mathbb Q} \mathbb G_{m,\mathcal F}$. Using this, it is easy to see that if the condition (UAX1) holds, then $A_0$ is isogeneous to $E_0^g$ for some elliptic curve $E_0$ over $\Spec(k)$ and the Rosati involution defined by the principal polarization $\theta_0$ of $A_0$ induces the complex conjugation on $\mathcal F$.

\smallskip
{\bf (c)} Similarly, if the condition (UAX2) holds, then $A_0$ is isogeneous to $E_0^g$ for some elliptic curve $E_0$ over $\Spec(k)$. One can check this by first remarking that $A_0$ is isogeneous to a power of a simple ordinary abelian variety $A_{00}$ over $\Spec(k)$ which has a model over a finite field whose Frobenius endomorphism generates a totally imaginary extension of its subfield fixed by the Rosati involution defined by $\theta_0$ and thus it is an imaginary quadratic extension of $\mathbb Q$. As $A_{00}$ is ordinary, the Honda--Serre--Tate theory implies that $A_{00}$ is an elliptic curve over $\Spec(k)$.
\end{rem}

\subsection{Lifting condition}\label{S33}
Assume $Q_1, Q_2\in\pmb{\Mat}_g(R)$ are symmetric. Let $(A_1,\theta_1)$, $(A_2,\theta_2)$ be the principally polarized abelian schemes over $\Spec(R)$ that lift $A_0$ and such that we have $(q_{A_1,\theta_1}(e_i,e_j))_{1\leq i,j\leq g}=Q_1$ and $(q_{A_2,\theta_2}(e_i,e_j))_{1\leq i,j\leq g}=Q_2$.
 By the functoriality part of Serre--Tate theory (see \cite{Ka81}, Thm. 2.1 4))
 the condition that the isogeny $u_0:A_0\rightarrow A_0$ lifts to a homomorphism $u:A_1\rightarrow A_2$ is
 that for all $1\leq i,j\leq g$ we have
 \begin{equation}\label{EQ061}
 q_{A_1,\theta_0}(e_i,T_p(u_0^{\dagger})(e_j))=q_{A_2,\theta_0}(T_p(u_0)(e_i),e_j).\end{equation}
 For all $c\in \mathbb Z_p$ and $m\in 1+\mathfrak m_R$ we let $c\star m=m\star c:=m^c\in 1+\mathfrak m_R$; hence $\star$ denotes multiplication in the (left and right) $\mathbb Z_p$-module structure of the multiplicative abelian group $1+\mathfrak m_R$. More generally, we want to consider left and right multiplication of matrices with entries in the $\mathbb Z_p$-module $1+\mathfrak m_R$ by matrices with entries in $\mathbb Z_p$. It is useful to introduce a special notation for this as follows. Let $\pmb{\Mat}_g(1+\mathfrak m_R)$ be the subset of $\pmb{\Mat}_g(R)$ formed by matrices whose all entries are in $1+\mathfrak m_R$.
 If $M\in \pmb{\Mat}_g(\mathbb Z_p)$ and 
 $X\in \pmb{\Mat}_g(1+\mathfrak m_R)$ we denote by $M\star X \in \pmb{\Mat}_g(1+\mathfrak m_R)$ the matrix with $ij$ entries given by 
 $$ (m_{i1}\star x_{1j})\cdots (m_{ig}\star x_{gj})=\prod_{l=1}^g x_{lj}^{m_{il}}$$
 and we denote by $X\star M \in \pmb{\Mat}_g(1+\mathfrak m_R)$ the matrix with $ij$ entries given by 
 $$ (x_{i1}\star m_{1j})\cdots (x_{ig}\star m_{gj})=\prod_{l=1}^g x_{il}^{m_{lj}}.$$
These operations satisfy the following four `associativity' properties: 
$$M_1\star(X\star M_2)=(M_1\star X)\star M_2,$$
$$(M_1\cdot M_2)\star X=M_1\star (M_2\star X),$$
$$X\star (M_1\cdot M_2)=(X\star M_1)\star M_2,$$
$$(X\star M)^{\t}=M^{\t}\star X^{\t}.$$
Here (and later) $M_1\cdot M_2=M_1 M_2$ denotes the usual product of the $g\times g$ matrices $M_1$ and $M_2$. 

Returning to Equation (\ref{EQ061}), set $M:=M_{u_0}$ (see Equation (\ref{EQ059})). Then Equation (\ref{EQ061}) is equivalent to the equality
 \begin{equation}\label{EQ062}Q_2=M^{-1}\star Q_1\star (M^{\dagger})^{\t}.\end{equation}

\begin{lemma}\label{L4}
For $M\in U_0$ with $MM^{\dagger}=d\cdot 1_g$, $d\in\mathbb N\setminus p\mathbb N$,
let $u_0:A_0\rightarrow A_0$ be the endomorphism corresponding to $M$; so $M=M_{u_0}$.
Assume $(A_1,\theta_1)$ is a principally polarized abelian scheme over $R$ that lifts $(A_0,\theta_0)$ and let $Q_1$ be the symmetric matrix with coefficients in $R$ that corresponds to $A_1$. 
Then the matrix $Q_2:=M^{-1}\star Q_1\star (M^{\dagger})^{\t}$ is symmetric. Moreover, if 
$(A_2,\theta_2)$ is the principally polarized abelian scheme over $\Spec(R)$ which lifts $(A_0,\theta_0)$ and corresponds to $Q_2$, then 
 there exists an isogeny $u:A_1\rightarrow A_2$ that lifts the isogeny $u_0:A_0\rightarrow A_0$ and such that
 $uu^{\t}=d$.
\end{lemma}

 \noindent
{\it Proof.} The symmetry of $Q_2$ follows from the fact that $Q_1\star (d\cdot 1_g)$ is symmetric and we have
$$Q_2=M^{-1}\star Q_1\star (M^{\dagger})^{\t}
=M^{-1}\star (Q_1\star (d\cdot 1_g)) \star (M^{-1})^{\t}.$$
As Equation (\ref{EQ062}) holds for $A_1$, $A_2$ and $M=M_{u_0}$ there exists an isogeny $u:A_1\rightarrow A_2$ that lifts the isogeny $u_0:A_0\rightarrow A_0$. As the reduction modulo $p$ map 
$\End(A_1)\rightarrow \End(A_0)$ is injective, $u_0u_0^{\dagger}=d$ implies that $uu^{\t}=d$.
\endproof

\medskip

Let $dt/t$ be the invariant form on the formal multiplicative group $\mathbb G_{m,R}^{\for}$ of the group scheme $\mathbb G_{m,R}=\Spec(R[t,t^{-1}])$ over $\Spec(R)$.

For $c\in\{1,2\}$, as $A_c$ is ordinary one has that 
the formal group $A_c^{\for}$ is isomorphic to $({\mathbb G}_{m,R}^{\for})^g$ over $R$; see \cite{Maz}, Lem. 4.27.
Let 
$\zeta_{ci}\in\Hom(A_c^{\for},{\mathbb G}_{m,R}^{\for})$ be the image via the isomorphism
\begin{equation}\label{EQ063}
T_p(\check{A}_0) \simeq\Hom(A_c^{\for},{\mathbb G}_{m,R}^{\for})\end{equation}
of the $i1$ entry of the column basis $T_p(\theta_0)(e)$ of $T_p(\check{A}_0)$. Defining
$$\omega_{ci}:=\zeta_{ci}^*(dt/t),$$
we get an $R$-basis $\omega_c=\{\omega_{c1},\ldots,\omega_{cg}\}$ of global $1$-forms on 
$A_c^{\for}$ which we identify with $1$-forms on $A_c$. We say that the $\omega_c$ is {\it attached} to the basis $e$.

\begin{lemma}\label{L5}
We assume that $M\in U_0$ and $Q_2=M^{-1}\star Q_1\star (M^{\dagger})^{\t}$ and we consider the isogeny $u:A_1\rightarrow A_2$ provided by Lemma \ref{L4}. We write $(u^{\for})^*(\omega_2)=[u]\omega_1$, where $[u]\in\pmb{\GL}_g(R)$. Then we have $[u]=M^{\dagger}$.
\end{lemma}

 \noindent
{\it Proof.} 
 By the functoriality of the isomorphism (\ref{EQ063}) we get that the following diagram is commutative:
 \begin{equation}
 \begin{array}{rcl}\label{EQ063.1}
 T_p(\check{A}_0) & \simeq & \Hom(A_1^{\for},{\mathbb G}_{m,R}^{\for})\\
 T_p(\check{u}_0) \uparrow & \ & \uparrow (u^{\for})^*\\
 T_p(\check{A}_0) & \simeq & \Hom(A_2^{\for},{\mathbb G}_{m,R}^{\for}).
 \end{array}
 \end{equation}
 Here $(u^{\for})^*(\zeta)=\zeta\circ u^{\for}$ for all $\zeta\in \Hom(A_2^{\for},{\mathbb G}_{m,R}^{\for})$. 
 Recalling that $M^{\dagger}=M_{u_0}^{\dagger}=M_{u_0^{\dagger}}$ and using Equation (\ref{EQ063.1}) we get that for all
  $i\in\{1,\ldots,g\}$ we have an equality of formal homomorphisms
 $$\zeta_{2i}\circ u^{\for}=\sum_{j=1}^g  m^{\dagger}_{ij}\zeta_{1j},$$
 hence
 $$(u^{\for})^* (\zeta_{2i}^* (dt/t))=\sum_{j=1}^g m^{\dagger}_{ij}\zeta_{1j}^*(dt/t).$$
Thus we have $[u]=M^{\dagger}$ as one can see at the level of formal Lie groups.
 \endproof
 
 \subsection{Crystalline forms}\label{S34}

 We first recall from \cite{BB}, Subsect. 4.1, the construction of some Hecke covariant Siegel $\delta$-modular forms $f^r_{g,\crys}\in\mathbb I^r_{gg}(\phi^r,1)$ indexed by $r\in\mathbb N$ which we refer to as the {\it basic} Siegel $\delta$-modular forms of weight $(\phi^r,1)$. Then,
 generalizing a construction for $g=1$ in \cite{Ba}, Constr. 3.2, p. 248 (cf. also \cite{Bu05}, Ch. 8, Sect. 8.4.3), we 
 introduce a Hecke covariant ordinary Siegel $\delta$-modular form $f^{\partial}_{g,\crys}\in\mathbb I^1_{gg,\ord}(-\phi,1)$ which we refer to as the {\it basic} Siegel $\delta$-modular form of weight $(-\phi,1)$.
 
 Let $S$ and $\pmb{\M}_g(S)$ be as in Subsection \ref{S232}, let $(A,\theta,\omega)\in \pmb{\M}_g(S)$
 and consider the de Rham $S$-modules
 $H^1_{\dR}(A/S)$ and $H^1_{\dR}(\check{A}/S)$.
 We recall from \cite{FC}, Ch. III, Sect. 9, paragraph before Prop. 9.1, that $H^1_{\dR}(A/S)$ has a canonical submodule
 $F^1H^1_{\dR}(A/S)$ canonically isomorphic to (and will be identified with) $H^0(A,\Omega_{A/S})$ and such that 
 $$H^1_{\dR}(A/S)/F^1H^1_{\dR}(A/S)\simeq H^1(A,\mathcal O)\simeq H^0(\check{A},\Omega_{\check{A}/S}^*).$$
From loc. cit. we also get that there exists a canonical perfect pairing 
\begin{equation}
\label{extra1}
\langle\ ,\ \rangle_A:H^1_{\dR}(\check{A}/S)\times H^1_{\dR}(A/S) \rightarrow S\end{equation}
which induces the canonical pairing between 
$F^1H^1_{\dR}(A/S)= H^0(A,\Omega_{A/S})$ and 
$H^1_{\dR}(\check{A}/S)/F^1H^1_{\dR}(\check{A}/S)= H^0(A,\Omega^*_{A/S})$.
In view of this compatibility, for all $(\omega,\omega')\in H^0(A,\Omega_{A/S})\times H^0(\check{A},\Omega_{\check{A}/S})$ we have $\langle \omega',\omega\rangle_A=0$.
 
We claim that the pairing $\langle \ ,\ \rangle_A$ is compatible with polarized 
isogenies in the sense that if $u:A_1\rightarrow A_2$ is a polarized isogeny of principally polarized abelian schemes, then
for all 
 $\check{\alpha}_1\in H^1_{\dR}(\check{A}_1/S)$ and $\alpha_2 \in H^1_{\dR}(A_2/S)$
 the following formula holds:
\begin{equation}\label{EQ063.5}
\langle \check{\alpha}_1,u^*\alpha_2\rangle_{A_1}=\langle \check{u}^*\check{\alpha}_1, \alpha_2\rangle_{A_2}.
\end{equation}
 To see this consider 
 first the union $\mathcal A'_{g,1,N,R}$ of all prime to $p$ `Hecke correspondences' on $\mathcal A_{g,1,N,R}$ parameterizing polarized isogenies that are compatible with level structures. Then $\mathcal A'_{g,1,N,R}$ is smooth over $R$ hence its local rings are integral domains of characteristic $0$.
 In view of compatibility with base change,  it is enough to check our claim for
 every ring $S$ with the property that $\Spec(S)$ is an open affine subscheme of $\mathcal A'_{g,1,N,R}$. So it is enough to check the claim for every integral domain $S$ of characteristic $0$; this allows us to assume $S$ is the complex field. In this case we use the analytic description of our objects as follows. If $V$ is a finite dimensional complex vector space and $\Lambda\subset V$ is a lattice, then the first Betti homology group $H_1(V/\Lambda,\mathbb Z)$ is naturally identified with $\Lambda$.
 Assume now $A$ is a complex abelian variety. Then by \cite{GH}, Ch. II, Sect. 6, pp. 327 and 328 we have canonical identifications of complex manifolds
 $$A(\mathbb C)= H^0(A,\Omega)^*/H_1(A(\mathbb C),\mathbb Z),\ \ \check{A}=H^1(A,\mathcal O)/H_1(A(\mathbb C),\mathbb Z)^*.$$
 So the natural pairing
 \begin{equation}
 \label{extra2}
 H_1(A(\mathbb C),\mathbb Z)^*\times H_1(A(\mathbb C),\mathbb Z)\rightarrow \mathbb Z
 \end{equation}
 induces a pairing on the Betti cohomology complex vector spaces:
 $$H^1(\check{A}(\mathbb C),\mathbb C)\times H^1(A(\mathbb C),\mathbb C)\rightarrow \mathbb C.$$
 The latter can be identified  with the pairing in Equation (\ref{extra1}) 
 under the natural identifications of the Betti and de Rham cohomology complex vector spaces.
 Indeed, by  \cite{FC}, Ch. III, Sect. 9, paragraph before Prop. 9.1,  the Chern class $c(\mathcal L)\in H^2_{\dR}(A\times \check{A}/\mathbb C)$ of the Poincar\'{e} bundle $\mathcal L\in \textup{Pic}(A\times \check{A})$ lies in $H^1_{\dR}(A/\mathbb C)\otimes_{\mathbb C}
  H^1_{\dR}(\check{A}/\mathbb C)\subset H^2_{\dR}(A\times \check{A}/\mathbb C)$ and it induces the perfect pairing (\ref{extra1}); on the other hand, using 
  the analytic description of the Poincar\'{e} bundle $\mathcal L$ in \cite{GH}, Ch. II, Sect. 6, p. 328 one checks that the pairing 
  between $H^1_{\dR}(A/\mathbb C)\simeq H^1(A(\mathbb C),\mathbb C)$ and $H^1_{\dR}(\check{A}/\mathbb C)\simeq H^1(\check{A}(\mathbb C),\mathbb C)$
  induced by $c(\mathcal L)$ coincides with the paring 
(\ref{extra2}).  Finally note that the pairing in Equation (\ref{extra2}) is compatible (in the obvious sense) with isogenies hence so is the pairing in Equation (\ref{extra1}).

We consider the pullback isomorphism $\theta^*: H^1_{\dR}(\check{A}/S)\rightarrow H^1_{\dR}(A/S)$
induced by $\theta:A\to\check{A}$ and define a bilinear map
$$\langle\ ,\ \rangle_{\theta}:H^1_{\dR}(A/S) \times H^1_{\dR}(A/S) \rightarrow S$$
by the rule
$$\langle\alpha,\beta\rangle_{\theta}:=\langle(\theta^*)^{-1}(\alpha),\beta\rangle_A,$$
where $\alpha,\beta\in H^1_{\dR}(A/S)$. By what was said above about the pairing $\langle\ ,\ \rangle_A$  we get that the bilinear map $\langle\ ,\ \rangle_{\theta}$ is a perfect pairing and has $H^0(A,\Omega^1)$ as an isotropic $S$-submodule. Clearly the formation of this paring commutes with  base change
$S\rightarrow S'$. Moreover, we have the following compatibility:

\begin{lemma}\label{L6}
Let $u:A_1\rightarrow A_2$ be a polarized isogeny  between principally polarized abelian schemes over $\Spec(S)$ with principal polarizations  $\theta_1$ and $\theta_2$
 (respectively). Let 
$$u^*:H^1_{\dR}(A_2/S)\rightarrow
H^1_{\dR}(A_1/S),\ \ \ u^{\t*}:H^1_{\dR}(A_1/S)\rightarrow
H^1_{\dR}(A_2/S)$$ be the induced isomorphisms. Then for all $\alpha\in H^1_{\dR}(A_1/S)$ and $\beta\in H^1_{\dR}(A_2/S)$ we have
$$\langle \alpha, u^*(\beta)\rangle_{\theta_1}=\langle u^{\t*}(\alpha),\beta\rangle_{\theta_2}.$$
\end{lemma}

\noindent
{\it Proof.} 
Using Equation (\ref{EQ063.5}) we compute
$$\langle\alpha, u^*(\beta)\rangle_{\theta_1}= \langle (\theta_1^*)^{-1}(\alpha),u^* (\beta)\rangle_{A_1}=\langle \check{u}^*\circ (\theta_1^*)^{-1}(\alpha),\beta\rangle_{A_2}$$
$$=\langle (\theta_2^*)^{-1}\circ u^{\t*}(\alpha),\beta\rangle_{A_2}=\langle u^{\t*}(\alpha),\beta\rangle_{\theta_2}.
$$
\endproof

\medskip
Let now 
 $S \in\Ob(\Prol)$ and $(A,\theta,\omega) \in \pmb{\M}_g(S^0)$.
Then for all integers $i\geq 0$, the crystalline theory provides
$\phi$-linear maps (see \cite{Bu00} or \cite{Bu05}, Sect. 8.1.9):
$$\Phi:H^1_{\dR}(A/S^0) \otimes_{S^0} S^i \rightarrow H^1_{\dR}(A/S^0) \otimes_{S^0} S^{i+1}.$$

Following \cite{BB}, Subsect. 4.1, for each integer $r\geq 1$ we have the $r$-th iterates 
$\Phi^r:H^1_{\dR}(A/S^0) \otimes_{S^0} S^i \rightarrow H^1_{\dR}(A/S^0) \otimes_{S^0} S^{i+r}$ of $\Phi$ and we define a Siegel
$\delta$-modular function $f^r_{g,\crys}$ of genus $g$, size $g$, and order $\leq r$ by the formula
$$f^r_{g,\crys}(A,\theta,\omega,S):=p^{-1}\langle\Phi^r \omega,
\omega^{\t}\rangle_{\theta}\in \pmb{\Mat}_g(S^r).$$
Writing $\omega=[\omega_1\cdots\omega_g]^{\t}$ as a column vector of $1$-forms,
 $\langle\Phi^r \omega, \omega^{\t}\rangle_{\theta}$ is the matrix whose each $ij$ entry is $\langle\Phi^r (\omega_i), \omega_j\rangle_{\theta}$.

\begin{lemma}\label{L7}
The rule $f^r_{g,\crys}$ defines a Hecke covariant Siegel $\delta$-form
of genus $g$, order $\leq r$, size $g$, and weight $(\phi^r,1)$, i.e.,
$$f^r_{g,\crys}\in\mathbb I^n_{gg}(\phi^r,1).$$
\end{lemma}

\noindent
{\it Proof.} 
This was proved in \cite{BB}, Subsect. 4.1. Alternatively one can prove this by a direct computation using Lemma \ref{L6} (see the proof of Lemma \ref{L8} below for a similar computation).
\endproof

\medskip

\begin{rem}\label{R21}
By construction, the form $f^1_{g,\crys}$ vanishes at some quadruple $(A,\theta,\omega,R)$ if and only if the filtered $F$-crystal over $\Spec(k)$ of the abelian scheme $A$ over $\Spec(R)$ is a direct sum of $2g$ filtered $F$-crystals over $\Spec(k)$ of rank $1$. Thus $(A,\theta)$ is either the canonical lift of an ordinary principally polarized abelian variety $(A_k,\theta_k)$ over $\Spec(k)$ or, if $p=2$, one of the $2^{\frac{g(g+1)}{2}}$ quasi-canonical lifts of the ordinary principally polarized abelian variety $(A_k,\theta_k)$ over $\Spec(k)$ (i.e., corresponds to a $\Spf(R)$-valued point of order $2$ of the formal torus of dimension $\frac{g(g+1)}{2}$ over $\Spf(R)$ of formal deformations of $(A_k,\theta_k)$). See \cite{Me}, App. and \cite{V}, Thm. 1.6.3 and Prop. 9.5.1.
\end{rem}

Assume now that $S \in\Ob(\Prol)$ and $(A,\theta,\omega) \in \pmb{\M}_{g,\ord}(S^0)$.
Recall that, as $A$ is ordinary,
we can consider the {\it unit root subspace} $U$ of $H^1_{\dR}(A/S^0)$, see \cite{Ka73b}, Thm. 4.1. Strictly speaking, loc. cit. is stated for the case when the reductions of $S^0$ modulo powers of $p$ are smooth over the corresponding reductions of $R$ but working in the \'etale topology, $(A,\theta)$ is the pullback of a principally polarized abelian scheme over the $p$-adic completion of a smooth $R$-algebra to which loc. cit. applies. The unit root subspace $U$ is a projective $S^0$-module which is a complement of $H^0(A,\Omega_{A/S^0})$ in $H^1_{\dR}(A/S^0)$. Assume it is free on $S^0$ and 
pick an arbitrary column vector $\eta=[\eta_1\cdots\eta_g]^{\t}$ whose entries form an $S^0$-basis of $U$. 

Consider the $g\times g$ matrix $\phi(\langle \eta,\omega^{\t}\rangle_{\theta})$ whose $ij$ entries are $\phi(\langle \eta_i,\omega_j \rangle_{\theta})\in S^1$ and the $g\times g$ matrix $\langle
\Phi \eta,\omega^{\t}\rangle_{\theta}$ whose $ij$ entries are $\langle \Phi
\eta_i,\omega_j\rangle_{\theta}\in S^1$.

Then we consider the product of invertible matrices
$$f^{\partial}_{g,\crys}(A,\theta,\omega, S):=
(\phi(\langle \eta,\omega^{\t}\rangle_{\theta}))^{-1}\cdot \langle \Phi
\eta,\omega^{\t}\rangle_{\theta}\in {\pmb{\GL}}_g(S^1).$$


\begin{lemma}\label{L8}
 The rule $f^{\partial}_{g,\crys}$ defines an ordinary Hecke covariant Siegel $\delta$-modular form of weight $(-\phi,1)$, size $g$, and order $1$, i.e.,
$$f^{\partial}_{g,\crys}\in\mathbb I^1_{gg,\ord}(-\phi,1).$$
In particular, its inverse defines a form
$$f_{\partial,g,\crys}:=(f^{\partial}_{g,\crys})^{-1}\in\mathbb I^1_{gg,\ord}(-1,\phi)$$
and the determinants of these two forms are forms
$$\det(f^{\partial}_{g,\crys})\in\mathbb I^1_{g,\ord}(\phi-1),$$
$$ \det(f_{\partial,g,\crys})=(\det(f^{\partial}_{g,\crys}))^{-1}\in\mathbb I^1_{g,\ord}(1-\phi).$$
\end{lemma}

\noindent
{\it Proof.}
The fact that $f^{\partial}_{g,\crys}$ has weight $(-\phi,1)$ follows from the following computation:
$$
\begin{array}{rcl}
f^{\partial}_{g,\crys}(A,\theta,\lambda \omega, S) & = &
(\phi(\langle \eta,\omega^{\t}\lambda^{\t}\rangle_{\theta}))^{-1}\cdot \langle \Phi \eta,\omega^{\t}\lambda^{\t}\rangle_{\theta}\\
\ & \ & \\
\ & = & ((\phi(\lambda))^{\t})^{-1}\cdot f^{\partial}_{g,\crys}(A,\theta, \omega, S) \cdot \lambda^{\t}.
\end{array}
$$
To check the Hecke covariance consider an object $S\in\Ob(\Prol)$ and  triples
$(A_1,\theta_1,\omega_1),(A_2,\theta_2,\omega_2)\in \pmb{\M}_{g,\ord}(S^0)$, let
$u:A_1\rightarrow A_2$ be a polarized isogeny with $u^*\omega_2=\omega_1$, $uu^{\t}=d\in\mathbb N\setminus p\mathbb N$, and we assume there exists a column vector $\eta_1$ whose entries form a basis of the unit root space of $A_1$.
Then, using Lemma \ref{L6} 
we have:
$$\begin{array}{rcl}
f^{\partial}_{g,\crys}(A_1,\theta_1,\omega_1,S) & = & (\phi(\langle \eta_1,\omega_1^{\t}\rangle_{\theta_1}))^{-1}\cdot \langle \Phi\eta_1,
\omega_1^{\t}\rangle_{\theta_1}\\
\ & \ & \ \\
\ & = & (\phi(\langle \eta_1,u^*\omega_2^{\t}\rangle_{\theta_1}))^{-1} \cdot \langle
\Phi \eta_1,u^*\omega_2^{\t}\rangle_{\theta_1} \\
\ & \ & \ \\
\ & = & (\phi(\langle u^{\t*}\eta_1,\omega_2^{\t}\rangle_{\theta_2}))^{-1} \cdot \langle
 u^{\t*}\Phi \eta_1,\omega_2^{\t}\rangle_{\theta_2}\\
 \ & \ & \ \\
 \ & = & (\phi(\langle u^{\t*}\eta_1,\omega_2^{\t}\rangle_{\theta_2}))^{-1} \cdot \langle
\Phi u^{\t*}\eta_1,\omega_2^{\t}\rangle_{\theta_2}\\
\ & \ & \ \\
\ & = & f^{\partial}_{g,\crys}(A_2,\theta_2,\omega_2,S).\end{array}$$
From this and  Equation (\ref{EQ029.8}), applied with $\deg(-\phi+1)=0$, we get that
$f^{\partial}_{g,\crys}\in\mathbb I^1_{gg,\ord}(-\phi,1)$.
\endproof

\medskip

Let us abbreviate $f^{\partial}:=f^{\partial}_{g,\crys}$
and for an integer $b\geq 1$ we set
$$N_b:=(f^{\partial})^{\phi^{b-1}} \cdot (f^{\partial})^{\phi^{b-2}} \cdot \cdot \cdot (f^{\partial})^{\phi}\cdot f^{\partial}\in\mathbb I^b_{gg,\ord}(-\phi^b,1).$$
We also set $N_0:=1_g$, the identity matrix.

Using the invertibility of $f^{\partial}$ we have:

\begin{cor}\label{C2}
For all integers $a,b\geq 0$ and $r\geq\max\{a,b\}$, the following $R$-linear maps are isomorphisms:
$$\begin{array}{rcll}
\mathbb I^r_{gg,\ord}(\phi^a,\phi^b) & \rightarrow & 
 \mathbb I^r_{gg,\ord}(1,1), & f\mapsto N_a^{\t}\cdot f\cdot N_b,\\
 \ & \ & \ & \ \\
 \mathbb I^r_{gg,\ord}(-\phi^a,\phi^b) & \rightarrow &
 \mathbb I^r_{gg,\ord}(-1,1), & f\mapsto N_a^{-1}\cdot f\cdot N_b,\\
 \ & \ & \ & \ \\
 \mathbb I^r_{gg,\ord}(\phi^a,-\phi^b) & \rightarrow & 
 \mathbb I^r_{gg,\ord}(1,-1), & f\mapsto N_a^{\t}\cdot f\cdot (N_b^{\t})^{-1}.\end{array}$$
\end{cor}

In particular, for each integer $a\geq 1$ we define recursively the form
\begin{equation}\label{EQ064}
f^{\langle a \rangle}=f^{\langle a \rangle}_{g,\crys}:=N_a^{\t} \cdot (f^1_{g,\crys})^{\phi^{a-1}}\cdot N_{a-1}\in\mathbb I_{gg}^a(1,1)
\end{equation}
and its transpose,
\begin{equation}\label{EQ065}
f^{\langle a \rangle t}:=(f^{\langle a \rangle}_{g,\crys})^{\t}\in\mathbb I_{gg}^a(1,1).
\end{equation}
Also we define the form
\begin{equation}\label{EQ066}
f^{[a]}=f^{[a]}_{g,\crys}:=N_a^{\t}\cdot f^a_{g,\crys}\in\mathbb I_{gg}^a(1,1)
\end{equation}
and its transpose
\begin{equation}\label{EQ067}
f^{[a]t}:=(f^{[a]}_{g,\crys})^{\t}\in\mathbb I_{gg}^a(1,1).
\end{equation}
Corollary \ref{C9} (c) below will check that $f^{\langle a \rangle t}=f^{\langle a \rangle}$ and $f^{[a]}=f^{[a]t}$. 

The following identity is a direct consequence of definitions:

\begin{equation}\label{EQ068}
f^{\langle a+1 \rangle}=(f^{\partial})^{\t}\cdot (f^{\langle a \rangle})^{\phi}\cdot f^{\partial}.
\end{equation}

As for size $g$ forms, we have the following result for size $1$ forms:

\begin{cor}\label{C3}
Let $w\in W(r)$. Let $w'\in W$ be the unique element such that
$$(\phi-1)w'=w-\deg(w).$$
 Then 
 the $R$-linear map
$$\mathbb I^r_{g,\ord}(w)\rightarrow
 \mathbb I^r_{g,\ord}(\deg(w)),\ \ \ f\mapsto \det(f^{\partial}_{g,\crys})^{-w'} \cdot f,$$
 is an isomorphism. 
\end{cor}

We end by recording the following facts proved in \cite{Bu05} based on \cite{Ba} (see \cite{Bu05}, Thm. 8.83 which relies on \cite{Bu05}, Prop. 8.22):

\begin{lemma} \label{L9}
The following three properties hold:

\medskip
{\bf (a)} We have $\mathbb I^1_1(\phi-1)=\mathbb I^1_1(1-\phi)=0$. In particular, 
$$f^{\partial}_{1,\crys}\in\mathbb I^1_{1,\ord}(\phi-1)\setminus\mathbb I^1_1(\phi-1),\ \ \ f_{\partial, 1,\crys}\in \mathbb I^1_{1,\ord}(1-\phi)\setminus \mathbb I^1_1(1-\phi).$$

\smallskip
{\bf (b)} For all distinct integers $a,b\geq 1$, the $R$-module $\mathbb I^{\max\{a,b\}}_1(-\phi^a-\phi^b)$ has rank $1$.

\smallskip
{\bf (c)} For all integers $a\geq 1$ we have $\mathbb I^a_1(-2\phi^a)=0$.
\end{lemma}
 
 \subsection{Serre--Tate expansions}\label{S35}
 
Let $T=T^{(0)},\ldots, T^{(r)}$ and their entries be as in Subsubsection \ref{S234}. Let
$$S^r_{\big-for}:=R[[T,\ldots,T^{(r)}]]=R[[T_{ij}^{(l)}|1\leq i\leq j\leq g, 0\leq l\leq r]],$$
and define
$$S^r_{\for}:=\reallywidehat{R[[T]][T',\ldots,T^{(r)}]}$$
similarly. So $S^0_{\for}=S^0_{\big-for}$ is as in Corollary \ref{C1}. 

The rings $S_{\big-for}:=\bigcup_{r\geq 0}S^r_{\big-for}$ and $S_{\for}:=\bigcup_{r\geq 0}S^r_{\for}$ are naturally in $\Ob(\Prol)$ with $\delta$ defined by the rule: $\delta(T^{(l)}_{ij})=T^{(l+1)}_{ij}$ for all integers $l\geq 0$ and $1\leq i\leq j\leq g$.

Let $(A_{S^0_{\for}},\theta_{S^0_{\for}})$ be the universal principally polarized abelian scheme over $\Spec(S^0_{\for})$ that lifts $(A_0,\theta_0)$ and that is the algebraization of the formal principally polarized abelian scheme over ${\mathcal M}_{A_0,\theta_0}$ (see Corollary \ref{C1}). Furthermore, let
$\zeta_i\in\Hom(A^{\for}_{S^0_{\for}},{\mathbb G}_m^{\for})$ be the image via the isomorphism (cf. Equation (\ref{EQ063}))
$$T_p(\check{A}_0) \simeq\Hom(A^{\for}_{S^0_{\for}},{\mathbb G}_{m}^{\for})$$
of the $i1$ entry of the column basis $T_p(\theta_0)(e)$ of $T_p(\check{A}_0)$. Defining
$\omega_i:=\zeta_i^*(dt/t)$,
we get an $R$-basis $\omega_{S^0_{\for}}=\{\omega_1,\ldots,\omega_g\}$ of global $1$-forms on  $A^{\for}_{S^0_{\for}}$ which we identify with $1$-forms on $A_{S^0_{\for}}$.

 Let $r\geq 0$ be an integer and let $w,w',w''\in W$. 
 
 For $f\in M^r_g(w)$ let 
 $$\mathcal E(f)= \mathcal E_{A_0,\theta_0,e}(f)=F(T,T',\ldots,T^{(r)}):= f(A_{S^0_{\for}},\theta_{S^0_{\for}},\omega_{S^0_{\for}},S_{\for})\in S^r_{\for}.$$
 We call $\mathcal E(f)$ the {\it Serre--Tate expansion} of $f$ at $(A_0,\theta_0,e)$.
  
 Similarly, for $f\in M^r_{gg}(w',w'')$ we let
 $$\mathcal E(f)=\mathcal E_{A_0,\theta_0,e}(f)=F(T,T',\ldots,T^{(r)})\in\pmb{\Mat}_g(S^r_{\for})\subset \pmb{\Mat}_g(S^r_{\big-for})$$
 be the matrix of formal power series obtained by evaluating $f$ at the quadruple
 $(A_{S^0_{\for}},\theta_{S^0_{\for}},\omega_{S^0_{\for}},S_{\for})$
 and we also call it the {\it Serre--Tate expansion} of $f$.

 We obtain {\it Serre--Tate expansion $R$-linear maps}
 \begin{equation}\label{EQ069}
 \mathcal E_{A_0,\theta_0,e}:M^r_g(w)\rightarrow S^r_{\for}
 \end{equation}
 and
 \begin{equation}\label{EQ070}
 \mathcal E_{A_0,\theta_0,e}:M^r_{gg}(w',w'')\rightarrow \pmb{\Mat}_g(S^r_{\for}).
 \end{equation}
As $A_{S^0_{\for}}$ is ordinary, the $R$-linear maps (\ref{EQ069}) and (\ref{EQ070}) extend to $R$-linear maps
 \begin{equation}\label{EQ071}
 \mathcal E_{A_0,\theta_0,e}:M^r_{g,\ord}(w)\rightarrow S^r_{\for}
 \end{equation}
 and
 \begin{equation}\label{EQ072}
 \mathcal E_{A_0,\theta_0,e}:M^r_{gg,\ord}(w',w'')\rightarrow \pmb{\Mat}_g(S^r_{\for}).
 \end{equation}
 
In what follows we would like to understand the restriction of these maps to Hecke covariant forms.
 
For each matrix $M\in U_0$ (so $MM^{\dagger}=d\cdot 1_g$ with $d\in\mathbb N\setminus p\mathbb N$), and for every symmetric matrices $Q_1$, $Q_2\in\pmb{\Mat}_g(R)$ that satisfy 
 $Q_2=M^{-1}\star Q_1\star (M^{\dagger})^{\t}$, let $(A_1,\theta_1)$ and $(A_2,\theta_2)$ be the corresponding principally polarized abelian schemes over $\Spec(R)$, let $\omega_1$ and $\omega_2$ be column bases of $1$-forms on $A_1$ and $A_2$ (respectively) attached to the column basis $e$ of $T_p(A_0)$, and let $u:A_1\rightarrow A_2$ be the corresponding isogeny (see Subsection \ref{S33}).
 
Let $f\in\mathbb I^r_{g,\ord}(w)$ be an ordinary Hecke covariant Siegel $\delta$-modular form of genus $g$, size $1$, order $\leq r$ and weight $w\in W$ with $\frac{g\cdot \deg(w)}{2}\in \mathbb Z$; so Equation (\ref{EQ029.5}) holds for it.
As $M$ has ${\mathbb Z}_p$-coefficients, from Lemma \ref{L5} we get that
$$\det([u])^{-w}=\det(M^{\dagger})^{-\deg(w)}.$$
Also $\det(M)\det(M^{\dagger})=d^g=\deg(u)$. 
 Hence Equation (\ref{EQ029.5}) becomes
 \begin{equation}\label{EQ073}
f(A_2,\theta_2,\omega_2,R)= \det(M)^{\deg(w)}\cdot d^{-\frac{g\cdot \deg(w)}{2}} \cdot f(A_1,\theta_1,\omega_1,R).
\end{equation}

If $f\in\mathbb I^r_{gg,\ord}(\phi^a,\phi^b)$ is a Hecke covariant Siegel $\delta$-modular form of genus $g$, size $g$, order $\leq r$ and weight $(\phi^a,\phi^b)$, then by Equation (\ref{EQ029.9}) we have:
\begin{equation}\label{EQ074}
f(A_2,\theta_2,\omega_2, R)= M^{\dagger} \cdot f(A_1,\theta_1,\omega_1, R)\cdot (M^{\t})^{-1}.
\end{equation}

If $f\in\mathbb I^r_{gg,\ord}(-\phi^a,\phi^b)$ is a Hecke covariant Siegel $\delta$-modular form of genus $g$, size $g$, order $\leq r$ and weight $(-\phi^a,\phi^b)$, then by Equation (\ref{EQ029.9}) we have:
\begin{equation}\label{EQ075}
f(A_2,\theta_2,\omega_2, R)= ((M^{\dagger})^{\t})^{-1} \cdot f(A_1,\theta_1,\omega_1, R)\cdot (M^{\dagger})^{\t}.
\end{equation}

If $f\in\mathbb I^r_{gg,\ord}(-\phi^a,-\phi^b)$ is a Hecke covariant Siegel $\delta$-modular form of genus $g$, size $g$, order $\leq r$ and weight $(-\phi^a,-\phi^b)$, then by Equation (\ref{EQ029.9}) we have:
\begin{equation}\label{EQ075.2}
f(A_2,\theta_2,\omega_2, R)= ((M^{\dagger})^{\t})^{-1} \cdot f(A_1,\theta_1,\omega_1, R)
\cdot M.
\end{equation}

Let $\mathbb Z_p[[T]]:=\mathbb Z_p[[T_{i,j}|1\leq i\leq j\leq g]]$ (recall that $T_{ij}=T_{ji}$).
For each $M\in U_0$, $MM^{\dagger}=d\cdot 1_g$, $d\in\mathbb N\setminus p\mathbb N$, 
 consider the $g\times g$ matrix with coefficients in $\mathbb Z_p[[T]]$ defined by
$$T^M:=M^{-1}\star (\mathbb E_g+T) \star (M^{\dagger})^{\t}-\mathbb E_g.$$
 Here the $\star$ operation is again defined in terms of the (left and right) $\mathbb Z_p$-module structure on the multiplicative abelian group $1+\mathbb Z_p[[T]]$.
 So if $M^{-1}=(m_{ij}^{(-1)})_{1\leq i,j\leq g}$, then the $ij$ entry of $T^M$ is
$$(T^M)_{ij}=\prod_{\tilde i=1}^g\prod_{\tilde j=1}^g (1+T_{{\tilde i}{\tilde j}})^{m^{(-1)}_{i\tilde i} m^{\dagger}_{j\tilde j}}-1\in\mathbb Z_p[[T]].$$
For every $G\in S^r_{\big-for}$ we define
$$G^M:=G(T^M,\delta(T^M),\ldots,\delta^r(T^M)).$$
Then for $F=(F_{ij})_{1\leq i,j\leq g}\in\pmb{\Mat}_g(S^r_{\big-for})$ we define 
 $$F^M:=(F_{ij}^M)_{1\leq i,j\leq g}\in\pmb{\Mat}_g(S^r_{\big-for}).$$

  Using the injectivity of the $R$-algebra homomorphism (\ref{EQ003})
 it is easy to check that if a series $F=F(T,T',\ldots,T^{(r)})\in S^r_{\for}$ has the property that
 $$F(pX,\delta (pX),\ldots,\delta^r (pX))=0$$ for all symmetric $g\times g$ matrices $X$ with coefficients in $R$, then $F=0$; indeed it is enough to note that the map
 $$X\mapsto F(pX,\delta(pX),\ldots,\delta^r(pX)$$
is a $\delta$-function (see Definition \ref{df7}) on the set of symmetric $g\times g$ matrices with entries in $R$.
 Using this remark we get that for $f\in\mathbb I^r_{g,\ord}(w)$ we have an equality
 \begin{equation}\label{EQ076}
 \mathcal E(f)^M= \det(M)^{\deg(w)}\cdot d^{-\frac{g\cdot \deg(w)}{2}} \cdot \mathcal E(f);\end{equation}
 indeed, in view of Lemma \ref{L4} and Equation (\ref{EQ073}) the difference between the right-hand side and the left-hand side of Equation 
 (\ref{EQ076}) vanishes 
 whenever $T,\ldots,T^{(r)}$ are replaced by 
 $pX,\ldots,\delta^r (pX)$ where $X$ is an arbitrary symmetric $g\times g$ matrix with entries in $R$.
 Similarly, from Equation (\ref{EQ074}) we get that for $f\in\mathbb I^r_{gg,\ord}(\phi^a,\phi^b)$ we have
 \begin{equation}\label{EQ077}
 \mathcal E(f)^M= M^{\dagger} \cdot \mathcal E(f) \cdot (M^{\t})^{-1}.\end{equation}
 Similarly, from Equation (\ref{EQ075}) we get that for 
 $f\in\mathbb I^r_{gg,\ord}(-\phi^a,\phi^b)$ we have
 \begin{equation}\label{EQ078}
 \mathcal E(f)^M= ((M^{\dagger})^{\t})^{-1}\cdot \mathcal E(f) \cdot (M^{\dagger})^{\t}.\end{equation}
 Similarly, from Equation (\ref{EQ075.2}) we get that for 
 $f\in\mathbb I^r_{gg,\ord}(-\phi^a,-\phi^b)$ we have
 \begin{equation}\label{EQ078.1}
 \mathcal E(f)^M= ((M^{\dagger})^{\t})^{-1}\cdot \mathcal E(f) \cdot M.\end{equation}
 
 Let $s\in \frac{1}{2}\mathbb Z$. If $gs\in \mathbb Z$ we consider the $R$-module
 $$\mathcal H_g^r(U_0,s)=\mathcal H_g^r(U_0,s)_R:=\{F\in S^r_{\big-for}|F^M=
 \det(M)^{-2s}d^{gs}
 \cdot F\;\;\forall\;\;M\in U_0\}.$$ 
 If $gs\not\in \mathbb Z$ we set $\mathcal H_g^r(U_0,s):=0$.

We define $R$-modules 
$$\mathcal H_{gg}^r(U_0)=\mathcal H_{gg}^r(U_0)_R, \ \ \mathcal H_{gg}^r(U_0)^*=\mathcal H_{gg}^r(U_0)^*_R,\ \ \ \mathcal H_{gg}^r(U_0)^{**}=\mathcal H_{gg}^r(U_0)^{**}_R$$
 as follows:
$$\mathcal H_{gg}^r(U_0):=\{F\in\pmb{\Mat}_g(S^r_{\big-for})|F^M=M^{\dagger}\cdot F \cdot (M^{\t})^{-1}\;\;\forall\;\;M\in U_0\},$$
$$\mathcal H_{gg}^r(U_0)^*:=\{F\in\pmb{\Mat}_g(S^r_{\big-for})|F^M=((M^{\dagger})^{\t})^{-1}\cdot F \cdot (M^{\dagger})^{\t}\;\;\forall\;\;M\in U_0\}.$$ 
$$\mathcal H_{gg}^r(U_0)^{**}:=\{F\in\pmb{\Mat}_g(S^r_{\big-for})|F^M=((M^{\dagger})^{\t})^{-1}\cdot F \cdot M\;\;\forall\;\;M\in U_0\}.$$ 

Due to Equation (\ref{EQ076}), (\ref{EQ077}), (\ref{EQ078}), (\ref{EQ078.1}) the rule $f\mapsto \mathcal E_{A_0,\theta_0,e}(f)$ defines $R$-linear maps 
\begin{equation}\label{EQ079}
\mathcal E_{A_0,\theta_0,e}:\mathbb I^r_{g,\ord}(w)\rightarrow\mathcal H^r_g(U_0,-\frac{\deg(w)}{2}),\end{equation}
\begin{equation}\label{EQ080}
\mathcal E_{A_0,\theta_0,e}:\mathbb I^r_{gg,\ord}(\phi^a,\phi^b)\rightarrow\mathcal H^r_{gg}(U_0),\end{equation}
\begin{equation}\label{EQ081}
\mathcal E_{A_0,\theta_0,e}:\mathbb I^r_{gg,\ord}(-\phi^a,\phi^b)\rightarrow\mathcal H^r_{gg}(U_0)^*,\end{equation}
\begin{equation}\label{EQ081}
\mathcal E_{A_0,\theta_0,e}:\mathbb I^r_{gg,\ord}(-\phi^a,-\phi^b)\rightarrow\mathcal H^r_{gg}(U_0)^{**}.\end{equation}
And by restriction we get $R$-linear maps
\begin{equation}\label{EQ082}
\mathcal E_{A_0,\theta_0,e}:\mathbb I^r_g(w)\rightarrow\mathcal H^r_g(U_0,-\frac{\deg(w)}{2}),\end{equation}
\begin{equation}\label{EQ083}
\mathcal E_{A_0,\theta_0,e}:\mathbb I^r_{gg}(\phi^a,\phi^b)\rightarrow\mathcal H^r_{gg}(U_0),\end{equation}
\begin{equation}\label{EQ084}
\mathcal E_{A_0,\theta_0,e}:\mathbb I^r_{gg}(-\phi^a,\phi^b)\rightarrow\mathcal H^r_{gg}(U_0)^*,\end{equation}
\begin{equation}\label{EQ084.1}
\mathcal E_{A_0,\theta_0,e}:\mathbb I^r_{gg}(-\phi^a,-\phi^b)\rightarrow\mathcal H^r_{gg}(U_0)^{**}.\end{equation}

The part (b) below is a `Serre--Tate expansion principle':

\begin{prop}\label{P3}

{\bf (a)} The $R$-linear maps in Equations (\ref{EQ040}), (\ref{EQ042}), and (\ref{EQ043}) above are injective and have torsion free cokernels. 

{\bf (b)}
The $R$-linear maps in Equations (\ref{EQ069}) to (\ref{EQ072}) and (\ref{EQ079}) to 
(\ref{EQ084.1}) above are injective and have torsion free cokernels. \end{prop}

\noindent
{\it Proof.} For part (a) we only consider the case of the $R$-linear map
$$
M^r_{g,\ord}(w)\rightarrow \mathcal O(J^r(X));
$$
 the other cases follow from this case. 
Let $f\in M^r_{g,\ord}(w)$ be such that we have $f_X=0\in \mathcal O(J^r(X))$ (see Equation (\ref{EQ041})) for some affine open subscheme $X$ of $\mathcal X_{\ord}=\mathcal A_{g,1,N,R,\ord}$ with $\overline{X}$ connected and with the property that there exists a column basis of $1$-forms on the universal abelian scheme over $X$. As the geometric fibers of the stack $\mathcal A_{g}$ over $\Spec(\mathbb Z)$ are irreducible (see \cite{FC}, Cor. 5.10), the geometric fibers of $\mathcal X_k$ have their connected components permuted transitively by the group $\pmb{\GSp}_{2g}(\mathbb Z/N\mathbb Z)$.
So if $\mathcal C$ is an arbitrary connected component of $\reallywidehat{\mathcal X}$ (the $p$-adic completion of $\mathcal X$) it follows that $f_{\Spf(C)}=0$ for some non-empty affine open formal subscheme $\Spf(C)$ of $\mathcal C$. Thus $f$ evaluated at the formal principally polarized abelian scheme over $\Spf(C)$ endowed with a column basis of $1$-forms vanishes, i.e., the element of $J^r(C)$ it defines is $0$. From the connectivity (hence irreducibility) of $\mathcal C$ we get that $f$ similarly vanishes at each formal principally polarized abelian scheme obtained from the universal one over $\mathcal C$ via a pullback and endowed with an arbitrary column basis of $1$-forms. 

Now let $B\in \Ob(\Prol)$. We want to prove that $f$ vanishes at each abelian scheme of relative dimension $g$ over $\Spec(B^0)$ which has a principal polarization and a column basis of $1$-forms, i.e., the elements of $B^r$ it defines are all $0$. There exists a finite flat \'etale homomorphism $B^0\rightarrow C^0$ such that the resulting principally polarized abelian scheme over $\Spec(C_0)$ has a symplectic similitude level-$N$ structure. By \cite{Bu05}, Lem. 3.14, the union $C:=\bigcup_{r\geq 0}(B^r\otimes_{B^0}C^0)$ has a natural structure of object in $\Prol$ and as the $B^r$-homomorphisms
\begin{equation}
\label{BntoBntomesC0}
B^r \rightarrow B^r \otimes_{B^0} C^0\end{equation}
are faithfully flat (hence injective), we can assume (by replacing $B$ with $C$) that the principally polarized abelian scheme over $\Spec(B)$ has a symplectic similitude level-$N$ structure. Thus the desired vanishing property of $f$ (i.e., of the mentioned elements of $B^r$) follows from the previous paragraph as $\mathcal C$ was an arbitrary connected component of 
$\reallywidehat{\mathcal X}$. 

The torsion freeness property for the cokernel is proved similarly by using the fact that 
the reduction modulo $p$  of the homomorphism (\ref{BntoBntomesC0}) is again faithfully flat and hence  injective.

For part (b) we only consider the case of the $R$-linear map (\ref{EQ069}): the cases of the other $R$-linear maps can be treated similarly.
Let $f\in M^r_g(w)$ be such that $\mathcal E(f)=0$.
Let $P_0\in\mathcal X(k)$ correspond to the principally polarized abelian scheme $(A_0,\theta_0)$ equipped with some symplectic similitude level-$N$ structure. Let $\Spf(C) \subset {\mathcal C}$ be an affine open formal subscheme that contains $P_0$ and sufficiently small so that column bases of $1$-forms on the corresponding formal abelian scheme over $\Spf(C)$ exist. As $\Spf(C)$ is connected and formally smooth over $\Spf(R)$, the $R$-algebra homomorphism $C\rightarrow S^0_{\for}$ that corresponds to the point $P_0$ is injective with torsion free cokernel and thus also the $R$-algebra homomorphism
$$J^r(C)\rightarrow S^r_{\big-for}$$
is injective with torsion free cokernel (see \cite{Bu05}, Ch. 4, Subsect. 4.5.1, Prop. 4.43). 
Thus $f$ evaluated at the formal principally polarized abelian scheme over $\Spf(C)$ endowed with a column basis of $1$-forms, i.e., the element of $J^r(C)$ it defines is $0$ as $\mathcal E(f)$ is so. From part (a) we get that $f=0$. The assertion about the torsion free cokernel is proved similarly.
\endproof

\medskip

Let $\mathbb I\in \{\mathbb I_g,\mathbb I_{g,\ord}\}$ and similarly consider the rings
 $\mathbb M,\mathbb M'$ (see Equation (\ref{EQ043.1})).
From Proposition \ref{P3} (b),  we get directly:

\begin{cor}\label{C4}
The following four properties hold:

\medskip
{\bf (a)} 
If $f_1\in\mathbb I^r(w_1)$ and $f_2\in\mathbb I^r(w_2)$ with $w_1,w_1\in W$ are such that
$f_1f_2$ is $0$ (respectively divisible by $p$) in $\mathbb I^r(w_1+w_2)$, then
either $f_1$ or $f_2$ is $0$ (respectively is divisible by $p$ in $\mathbb I^r(w_1)$
or $\mathbb I^r(w_2)$). 

\smallskip
{\bf (b)} The ring $\mathbb I$ is a $p$-adically separated integral domain and $p$ is a prime element of it, hence the local ring $\mathbb I_{(p)}$ is a discrete valuation ring. Moreover, we have $\mathbb I^r\cap p\mathbb I^{r+1}=p\mathbb I^r$ for all integers $r\geq 0$.

\smallskip
{\bf (c)} The endomorphism $\phi:\mathbb I\rightarrow \mathbb I$ is injective with injective reduction modulo $p$.

\smallskip
{\bf (d)} The analog statements (a) to (c) hold for the rings $\mathbb M,\mathbb M',\mathbb M_X$ in place of $\mathbb I$. Thus the rings $\mathbb I,\mathbb M,\mathbb M',\mathbb M_X$ have natural structures of filtered integral $\delta$-graded rings. Moreover, the morphism in (\ref{EQ040}) is injective with injective reduction modulo $p$.\end{cor}

Property (a) of the next corollary says  that no non-zero Siegel $\delta$-modular form of size $1$
can have an order smaller than the order of its weight. 

\begin{cor}\label{C5}
For all $w\in W$ and integers $r\geq 0$ the following five properties hold:

\medskip
{\bf (a)}
For all integers $s,r$ with $0\leq s<\ord(w)\leq r$ we have $M^r_g(w)\cap M^s_g=0$.

\smallskip
{\bf (b)} The $R$-algebra homomorphism $\mathbb M^r_g\rightarrow M^r_g$ is injective.

\smallskip
{\bf (c)} For all $v\in W$, if $0\neq F\in \mathbb I_g(w)$ and $G\in M^r_g$ are such that $FG\in \mathbb I_g(w+v)$, then $G\in \mathbb I^r_g(v)$.

\smallskip
{\bf (d)} If $f\in \mathbb I^r_g(w)\setminus p \mathbb I^r_g(w)$, then
$$(\mathbb I_{g,\langle f\rangle})^r=(\mathbb M_{g,\langle f\rangle})^r\cap \mathbb I_{g,\langle f\rangle}.$$

\smallskip
{\bf (e)} Properties (a) to (d) hold in the ordinary case as well.
\end{cor}

\noindent
{\it Proof.} We consider only the non-ordinary case as the ordinary care is similar. 
To prove part (a) we check that the assumption that there exists a non-zero $f\in M^r_g(w)\cap M^s_g$ leads to a contradiction. For each $\lambda\in (S^0_{\for})^{\times}$ we have
\begin{equation}\label{EQ085}
\mathcal E^{\lambda}(f):=f(A_{S^0_{\for}},\theta_{S^0_{\for}},\lambda \omega_{S^0_{\for}},S_{\for})\in S^s_{\for},\end{equation}
\begin{equation}\label{EQ086}
\mathcal E(f):=f(A_{S^0_{\for}},\theta_{S^0_{\for}}, \omega_{S^0_{\for}},S_{\for})\in S^s_{\for}.\end{equation}
As $f$ has weight $w$, 
we have $\mathcal E^{\lambda}(f)=\lambda^{-gw} \mathcal E(f)$, see Equation (\ref{EQ029.3}). As by the Serre--Tate expansion principle (Proposition
\ref{P3} (b)) we have $\mathcal E(f)\neq 0$,
we conclude that 
$\lambda^{gw}\in \Frac(S^s_{\for})$.
But, as $s<\ord(w)$, this is clearly false, say, for $\lambda:=1_g+T_{11}$.

To prove part (b) we check that if $(f_w)_{w\in W(r)}\in \mathbb M^r_g=\oplus_{w\in W(r)}M^r_g(w)$ and $\sum_{w\in W(r)}f_w=0$ in $M^r_g$, then $f_w=0$ for all $w\in W(r)$. Using the notation in (\ref{EQ085}) and (\ref{EQ086}) we have $\sum_{w\in W(r)}\mathcal E^{\lambda}(f_w)=0$ for all $\lambda\in R^{\times}$ hence $\sum_{w\in W(r)}\lambda^{-w}\mathcal E(f_w)=0$ in $S^r_{\for}$. For each monomial $\mu$ in the indeterminates $T_{ij},T'_{ij},\ldots,T^{(r)}_{ij}$ with $i,j\in\{1,\ldots,g\}$ let $a_{w,\mu}\in R$ be the coefficient of $\mu$ in $\mathcal E(f_w)$. We get that $\sum_{w\in W(r)}a_{w,\mu}\lambda^w=0$ for all $\mu$ and
all $\lambda\in R^{\times}$.
The latter easily implies that $a_{w,\mu}=0$ for all $\mu$ and $w\in W(r)$, hence $f_w=0$ for all $w\in W(r)$.

For  part (c) we use the fact that when checking the conditions defining weight and Hecke covariance it suffices to do so in the `universal situation' in which case the filtered $\delta$-rings  under consideration are integral domains. 
So it is enough to check the conditions for every
 quadruple $(A,\theta,\omega,S)$ with $S$ an integral domain 
 and $F(A,\theta,\omega,S)\neq 0$.
 For every such quadruple and each $\lambda\in (S^0)^{\times}$ we have
$$\begin{array}{rcl}
F(A,\theta,\lambda \omega,S)G(A,\theta,\lambda \omega,S) & = & (FG)(A,\theta,\lambda \omega,S) \\
\ & \ & \ \\
\ & = & \det(\lambda)^{-w-v}F(A,\theta,\omega,S)G(A,\theta,\omega,S)\\
\ & \ & \ \\
\ & = & 
\det(\lambda)^{-v}F(A,\theta,\lambda\omega,S)G(A,\theta,\omega,S).\end{array}$$
As $S$ is an integral domain we can cancel the non-zero factor $F(A,\theta,\lambda\omega,S)$
and get that $G$ has weight $v$. A similar argument shows that $G$ is Hecke covariant.

To prove part (d) let $F/f^v=G/f^u\in (\mathbb M_{g,\langle f\rangle})^r\cap \mathbb I_{g,\langle f\rangle}$ with $v,u\in W$, $F\in \mathbb I_g(wv)$, and $G, f^u\in \mathbb M^r_g(wu)$. Then $f^v G\in \mathbb I_g(wu+wv)$. By part (c) we have $G\in \mathbb I_g(wu)$, hence
$G\in \mathbb I^r_g(wu)$ and therefore $G/f^u\in (\mathbb I_{g,\langle f\rangle})^r$.

The proofs of parts (a) to (d) in the ordinary case are similar.\endproof

\subsection{Comparison maps: the ordinary case}\label{S36}

In this subsection we first define some new modules and algebras and basic `comparison maps' $\spadesuit$, $\clubsuit$, $\diamondsuit$, and $\heartsuit$ between them. Then we will use these `comparison maps', in conjunction with our Serre--Tate expansion maps $\mathcal E$ and with the invariant theory 
of multiple quadratic forms (see Subsections \ref{S41} to \ref{S43}), to prove our basic computations for different dimensions of $K$-algebras and dimensions of
$K$-vector spaces of Hecke covariant ordinary Siegel $\delta$-modular
forms.

\subsubsection{Basic modules of forms}\label{S3601}

We start by considering some variants of prior modules of formal power series (and new modules, consisting of polynomials) over rings other than $R$.

First, if in the definition of 
$$\mathcal H^r_g(U_0,s),\ \ \mathcal H^r_{gg}(U_0),\ \ \mathcal H_{gg}^r(U_0)^*,\ \ \mathcal H_{gg}^r(U_0)^{**}$$
we replace $R$ by an arbitrary integral domain $B$ that contains $\mathbb Z_p$, we get
 corresponding $B$-modules
$$\mathcal H^r_g(U_0,s)_B,\ \ \mathcal H^r_{gg}(U_0)_B,\ \ \mathcal H_{gg}^r(U_0)^*_B,\ \ 
\mathcal H_{gg}^r(U_0)^{**}_B.$$
Next, for each polynomial 
$$G\in B[T,\ldots,T^{(r)}]:=B[T_{ij},T'_{ij},\ldots,T^{(r)}_{ij}|1\leq i,j\leq g]$$ 
(recall that we have $T_{ij}=T_{ji}$) and each $M\in U_0$ we set
 $$G_M:=G(M^{-1}\cdot T\cdot (M^{\dagger})^{\t},\ldots,M^{-1}\cdot T^{(r)}\cdot (M^{\dagger})^{\t}).$$
Recall, here and below, $\cdot$ is matrix multiplication. 

Let $s\in \frac{1}{2}\mathbb Z$.  If $gs\not\in \mathbb Z$, then we set $H^r_g(U_0,s)_B:=0$. If $gs\in \mathbb Z$, then we define $B$-modules
 $$H^r_g(U_0,s)_B:=\{G\in B[T,\ldots,T^{(r)}]|G_M=\det(M)^{-2s}d^{gs}\cdot G\;\forall\;
 M\in U_0\},$$ 
 $$H^r_{gg}(U_0)_B:=\{G\in\pmb{\Mat}_g(B[T,\ldots,T^{(r)}])|G_M=M^{\dagger}\cdot G\cdot (M^{\t})^{-1}\;\forall\;M\in U_0\},$$ 
 $$H^r_{gg}(U_0)^*_B:=\{G\in\pmb{\Mat}_g(B[T,\ldots,T^{(r)}])|G_M=((M^{\dagger})^{\t})^{-1} \cdot G \cdot (M^{\dagger})^{\t}\;\forall\; M\in U_0\},$$
 $$H^r_{gg}(U_0)^{**}_B:=\{G\in\pmb{\Mat}_g(B[T,\ldots,T^{(r)}])|G_M=((M^{\dagger})^{\t})^{-1} \cdot G \cdot M\;\forall\; M\in U_0\}.$$

If $D\subset B$ is a multiplicative set, then we have
 $$H^r_g(U_0,s)_{D^{-1}B}=D^{-1}H^r_g(U_0,s)_B,\ \ 
 H^r_{gg}(U_0)_{D^{-1}B}=D^{-1}H^r_{gg}(U_0)_B,$$
 $$ H^r_{gg}(U_0)^*_{D^{-1}B}=D^{-1}H^r_{gg}(U_0)^*_B,\ \ \ H^r_{gg}(U_0)^{**}_{D^{-1}B}=D^{-1}H^r_{gg}(U_0)^{**}_B.$$
 Also, if $B\subset B'$ is a field extension, then we have
 $$H^r_g(U_0,s)_{B'}=H^r_g(U_0,s)_B\otimes_B B',\ \ H^r_{gg}(U_0)_{B'}=H^r_{gg}(U_0)_B\otimes_B B',$$
 $$H^r_{gg}(U_0)^*_{B'}=H^r_{gg}(U_0)^*_B\otimes_B B',\ \ \ H^r_{gg}(U_0)^{**}_{B'}=H^r_{gg}(U_0)^{**}_B\otimes_B B'.$$
Last two properties do not hold if $H$ is replaced by $\mathcal H$ as we shall see below.

For an infinite integral domain $B$  and for $s\in \frac{1}{2}\mathbb Z$ with $gs\in \mathbb Z$, let $\mathbb H^r_g(s)_B$ be the $B$-module of homogeneous polynomials $G$ in the ring $B[T,\ldots,T^{(r)}]$ of degree $gs$ that satisfy
 \begin{equation}\label{EQ087}
 G(X\cdot T\cdot X^{\t},\ldots,X\cdot T^{(r)}\cdot X^{\t})=\det(X)^{2s}G\end{equation}
 for all $X\in\pmb{\Mat}_g(B)$. 
 Here the entries  of the matrices $T,\ldots,T^{(r)}$ are given degree $1$.
 Note that $\mathbb H^r_g(s)_R$ coincides with the modules $\mathbb H^r_g(s)$ introduced in Equation (\ref{EQ034}).
 If $B$ is a field, then $\mathbb H^r_g(s)_B$ is the $B$-vector space of $\pmb{\SL}_g(B)$-invariant homogeneous polynomials of degree $gs$ in $B[T,\ldots,T^{(r)}]$ under the left action 
 defined by the rule: 
 $$(X,F)\mapsto F(X^{-1}\cdot T\cdot (X^{-1})^{\t},\ldots,X^{-1}\cdot T^{(r)}\cdot (X^{-1})^{\t}).$$
 Similarly, for an infinite integral domain $B$ 
 let $\mathbb H^r_{gg}(1)_B$ be the $B$-module of homogeneous matrices 
 $$G\in \pmb{\Mat}_g(B[T,\ldots,T^{(r)}])$$ of degree $1$
 such that we have
 \begin{equation}\label{EQ088}
 G(X\cdot T\cdot X^{\t},X\cdot T'\cdot X^{\t},\ldots,X\cdot T^{(r)}\cdot X^{\t})=X \cdot G(T,T',\ldots,T^{(r)})\cdot X^{\t}
 \end{equation}
 for all $X\in\pmb{\Mat}_g(B)$.  Again, for each multiplicative set $D\subset B$ we have
 $$\mathbb H^r_g(s)_{D^{-1}B}=D^{-1}\mathbb H^r_g(s)_B\ \ \ \textup{and}\ \ \ \mathbb H^r_{gg}(1)_{D^{-1}B}=D^{-1}\mathbb H^r_{gg}(1)_B.$$
 Also, if $B\subset B'$ is an extension of infinite fields, then
 $$\mathbb H^r_g(s)_{B'}=\mathbb H^r_g(s)_B\otimes_B B'\ \ \ \textup{and}\ \ \ \mathbb H^r_{gg}(1)_{B'}=\mathbb H^r_{gg}(1)_B\otimes_B B'.$$
 
 \medskip
 
 In what follows we will repeatedly use conditions (UAX1) and (UAX2) introduced in Subsection \ref{S32}.
 
\begin{lemma}\label{L10}
Assume $(A_0,\theta_0)$ satisfies condition \textup{(UAX1)} and let $B$ be an integral domain that contains $\mathbb Z_p$. 

\medskip
{\bf (a)} Each polynomial in $H^r_g(U_0,s)_B$ is homogeneous of degree $gs$. In particular, we have $H^r_g(U_0,s)_B=0$ for all $s<0$ and $H^r_g(U_0,0)_B=B$.

\smallskip
{\bf (b)} Each matrix in $H^r_{gg}(U_0)_B$ is homogeneous of degree $1$.

\smallskip
{\bf (c)} We have $H^r_{gg}(U_0)^*_B\subset \pmb{\Mat}_g(B)$, i.e., each matrix in $H^r_{gg}(U_0)^*_B$ is homogeneous of degree $0$.

\smallskip
{\bf (d)} We have $H^r_{gg}(U_0)^{**}_B=0$.
 \end{lemma}
 
 \noindent
 {\it Proof.}
As $(A_0,\theta_0)$ satisfies condition (UAX1), there exists $\lambda\in\mathcal U_0$ such that 
$\lambda/\lambda^{\dagger}$ is not a root of unity. For $F\in H^r_g(U_0,s)$, we write $F=\sum_{n\geq 0} F_n$, with each $F_n$ homogeneous of degree $n$, and we compute
 $$F_{\lambda^2 \cdot 1_g}=\sum_{n\geq 0} (\lambda/\lambda^{\dagger})^{-2n}\cdot F_n.$$
 On the other hand, as $F\in H^r_g(U_0,s)$, we have
 $$F_{\lambda^2 \cdot 1_g}=\lambda^{-4gs}d^{2gs}\cdot F=(\lambda/\lambda^{\dagger})^{-2gs}\cdot F=\sum_{n\geq 0} (\lambda/\lambda^{\dagger})^{-2gs}F_n.$$
From the two expressions of $F_{\lambda^2 \cdot 1_g}$ we get that $F_n=0$ if $n\neq gs$. Thus part (a) holds. Parts (b), (c), (d) are checked similarly. \endproof

 \begin{prop}\label{P4}
 Assume $(A_0,\theta_0)$ satisfies conditions \textup{(UAX1)} and \textup{(UAX2)}.
 Then 
 $$H^r_g(U_0,s)_{\mathbb Z_p}=\mathbb H^r_g(s)_{\mathbb Z_p},\ \ 
 H^r_{gg}(U_0)_{\mathbb Z_p}=\mathbb H^r_{gg}(1)_{\mathbb Z_p},\ \ 
 H^r_{gg}(U_0)^*_{\mathbb Z_p}=\mathbb Z_p\cdot 1_g
 $$ and hence for every integral domain $B$ that contains $\mathbb Z_p$ we have
 $$H^r_g(U_0,s)_B=\mathbb H^r_g(s)_B,\ \ H^r_{gg}(U_0)_B=\mathbb H^r_{gg}(1)_B,
 \ \ 
 H^r_{gg}(U_0)^*_B=B\cdot 1_g.$$
 \end{prop}
 
 \noindent
 {\it Proof.} Let $G\in H^r_g(U_0,s)_{\mathbb Z_p}$
 and we prove that $G\in\mathbb H^r_g(s)_{\mathbb Z_p}$. (The converse is 
 proved similarly.) We need to show that Equation (\ref{EQ087}) holds
 for all $X\in\pmb{\Mat}_g(\mathbb Z_p)$. By Lemma \ref{L10} (a), $G$ is homogeneous of degree $gs$. As $\mathbb Z_{(p)}^{\times} \cdot U_0$ is Zariski dense in $\pmb{\GL}_{g,\mathbb Z_p}$, it is enough to show that Equation (\ref{EQ087}) holds for all $X\in\mathbb Z_{(p)}^{\times} \cdot U_0$. By homogeneity of (\ref{EQ087}) it is enough to show that Equation  (\ref{EQ087})
 holds for all $X\in U_0$. We know that we have
 \begin{equation}\label{EQ089}
 G(Y^{-1}T(Y^{\dagger})^{\t},\ldots,Y^{-1}T^{(r)}(Y^{\dagger})^{\t})=
 \det(Y)^{-2s}d^{gs}\cdot G(T,\ldots,T^{(r)})\end{equation}
 for all $Y\in U_0$. Write $YY^{\dagger}=d\cdot 1_g$. We get
 $Y^{-1}=d^{-1}Y^{\dagger}$. From this and the homogeneity part of Lemma \ref{L10} (a), as $\det(Y)\det(Y^{\dagger})=d^g$, we get that the left-hand side of Equation  (\ref{EQ089}) equals
 $$
 \det(Y^{\dagger})^{2s}G(Y^{\dagger}T(Y^{\dagger})^{\t},\ldots,
 Y^{\dagger}T^{(r)}(Y^{\dagger})^{\t}).$$
 Setting $X=Y^{\dagger}$ we get that Equation  (\ref{EQ087}) holds for all $X\in U_0$.
 
 Similarly, let $G\in H^r_{gg}(U_0)_{\mathbb Z_p}$ and let us prove that $G\in\mathbb H^r_{gg}(1)_{\mathbb Z_p}$; the converse is trivial (and does not need condition (UAX2)). 
 Indeed, by Lemma \ref{L10} (b), $G$ is homogeneous of degree $1$. We know that 
 $$G(M^{-1}T(M^{\dagger})^{\t},\ldots,M^{-1}T^{(r)}(M^{\dagger})^{\t})=M^{\dagger}\cdot G \cdot (M^{\t})^{-1}$$
 holds
 for all $M\in U_0$. Using the fact that $MM^{\dagger}=d\cdot 1_g$ plus the homogeneity of $G$ and writing $X=M^{-1}$
 we get that
 Equation (\ref{EQ088}) holds
 for all $X\in\mathbb Z_{(p)}^{\times}\cdot U_0$. Using condition (UAX2) we get that Equation
(\ref{EQ088}) holds for all $X\in\pmb{\GL}_g(\mathbb Z_p)$, 
so for all $X\in\pmb{\Mat}_g(\mathbb Z_p)$,
hence $G\in\mathbb H^r_{gg}(1)_{\mathbb Z_p}$.

If $G\in H^r_{gg}(U_0)^*_{{\mathbb Z}_p}\subset\pmb{\Mat}_g(\mathbb Z_p)$ (see Lemma \ref{L10} (c)), we get that $G$ commutes with all matrices in the set $\{(M^{\dagger})^{\t}|M\in U_0\}$.
By condition (UAX2) we get that $G$ commutes with all the elements of $\pmb{\GL}_g(\mathbb Z_p)$, hence $G\in\mathbb Z_p\cdot 1_g$.
\endproof

\medskip
 
 In what follows we consider integers $s, r\geq 0$. Let
 $$\Delta(g,r):=\{(m_0,\ldots,m_r)\in\mathbb Z_{\geq 0}^{r+1}|m_0+\cdots+m_r=g\}.$$
Let $B$ be an infinite integral domain.
Note that $\mathbb H^r_g(s)_B$ contains all homogeneous polynomials of degree $gs$ which are products of $s$ polynomials of the form
\begin{equation}\label{EQ090}
\Theta_{m_0,\ldots,m_r}\in\mathbb H^r_g(1)_B\ \ \ \textup{\for} \ (m_0,\ldots, m_r)\in\Delta(g,r)
\end{equation}
that are defined by the identity
$$\det(y_0 T+y_1T'+\cdots+ y_rT^{(r)})=\sum_{(m_0,\ldots, m_r)\in\Delta(g,r)} \Theta_{m_0,\ldots,m_r} y_0^{m_0}y_1^{m_1}\ldots y_r^{m_r}$$
in which $y_0,\ldots,y_r$ are  indeterminates; so $\Theta_{m_0,\ldots,m_r}$ are partial polarizations of the invariant polynomial $\det$. For $r=1$ the classical notation is $\Theta_i=\Theta_{g-i,i}$ for $0\leq i\leq g$; so $\Theta_0=\det(T)$.

For an infinite field $B$ of characteristic $\mathfrak p\geq 0$ and $s\in \frac{1}{2}\mathbb Z$, the dimension
\begin{equation}
\label{Dpgrs}
D_{\mathfrak p}(g,r,s):=\dim_B(\mathbb H^r_g(s)_B),\end{equation}
depends (as we shall see) only on $g,r,s,\mathfrak p$, and equals $0$ for $s<0$. For $B$ an infinite integral domain consider the integral extension of $B$-algebras
$$\mathbb H^r_{g,B}:=\bigoplus_{s\in \mathbb Z}\mathbb H^r_g(s)_B\subset 
\mathbb H^r_{g,\tot,B}:=\bigoplus_{s\in \frac{1}{2}\mathbb Z}\mathbb H^r_g(s)_B.$$
For an infinite field $B$ set
$$D(g,r):=\dim(\mathbb H^r_{g,B})=\dim(\mathbb H^r_{g,\tot,B}),$$
 which, as we shall see shortly, does indeed depend only on $g,r$.

We have the following result (see Theorem \ref{T29} (a) in Subsubsection \ref{S4205}) which, at least in characteristic $0$, is classical:
 
\begin{lemma}\label{L11}
Let $B$ be an infinite field of characteristic $\mathfrak p\geq 0$, $\mathfrak p\neq 2$.
The $B$-algebra $\mathbb H^1_{g,\tot, B}$ is generated by the algebraically independent elements $\Theta_0,\ldots,\Theta_g$.
 Hence, for $s\in \mathbb N\cup \{0\}$, we have
$$D_{\mathfrak p}(g,1,s)=\frac{(g+s)!}{g!s!}\;\;\;\textup{and}\;\;\;D_{\mathfrak p}(g,0,s)=1.$$
In addition $D_{\mathfrak p}(g,1,s)=0$ if $g$ is even and $s\notin \mathbb Z$.
\end{lemma}

Also we have:

\begin{lemma}\label{L12}
The following three properties hold:

\smallskip

{\bf (a)}
Let $B$ be either an infinite field or $R$. Then
for each integer $r\geq 1$, the $B$-algebra $\mathbb H^r_{g,\tot, B}$ is finitely generated and a unique factorization domain.

\smallskip
{\bf (b)} For $B$ an infinite field we have
$$D(g,r)=(r+1)\cdot \frac{g(g+1)}{2}-g^2+1.$$

\smallskip
{\bf (c)} For a field $B$ of characteristic $0$ there exists a polynomial $D_{g,r}(x)\in\mathbb Q[x]$ of degree 
$$D(g,r)-1=\frac{g(rg+r-g+1)}{2}$$
such that for all $s\in \mathbb Z_{\geq 0}$ big enough we have
$$D_0(g,r,s)=D_{g,r}(s).$$
Moreover, if $g$ is even, then there exists a polynomial $D_{g,\odd,r}(x)\in\mathbb Q[x]$ of degree at most $D(g,r)-1$
such that for all $s\in \mathbb N$ big enough we have
$$D_0(g,r,\frac{1}{2}+s)=D_{g,\odd,r}(s).$$
\end{lemma}

For detailed information on the numbers $D_0(g,r,s)$ (not needed for the present section) we refer to results in Section 4 (e.g.,
Theorem \ref{T29}
in Subsubsection \ref{S4205}, Proposition \ref{P15} in Subsubsection \ref{S4207},
and Corollaries \ref{C21} and \ref{C22} in Subsubsection \ref{S4210}).

\medskip

\noindent
{\it Proof.}
If $B$ is an infinite field or $R$, then, by Lemma \ref{L43} (a) and (b) of Subsubsection \ref{S4202}, 
 $\mathbb H^r_{g,\tot, B}$ is the $B$-algebra of 
 invariants $B[T_{ij}|1\leq i\leq j\leq g]^{\pmb{\SL}_g}$. 
To prove part (b) or part (a) in case $B$ is a field, we can assume $B$ is algebraically closed in which case we conclude by Corollary \ref{C20} of Subsubsection \ref{S4203}.
It remains to prove part (a) in case $B=R$. Recalling that the upper bar signifies the reduction modulo $p$ (or in the present context, tensoring with $k$ over $R$), note that the homomorphism $\overline{\mathbb H^{r}_{g,\tot,R}}\rightarrow \mathbb H^{r}_{g,\tot,k}$ 
is injective because if $p$ times a polynomial is invariant then the polynomial is invariant.
Hence $p$ is a prime element in $\mathbb H^r_{g,\tot,R}$. As $\mathbb H^r_{g,\tot,R} \otimes_R K=\mathbb H^r_{g,\tot,K}$ 
is a unique factorization domain it follows that $\mathbb H^{r}_{g,\tot,R}$ is a unique factorization domain (see \cite{Mat}, Thm. 20.2), and from a general result on finite generation of rings of invariants over universally Japanese rings such as $R$ (see \cite{Se}, Thm. 2; cf. also Subsubsection \ref{S4203}) we get that $\mathbb H^r_{g,\tot,R}$ is a finitely generated $R$-algebra. However, for convenience we include a direct argument. 

Let $\mathbb H^{r,\Theta}_{g,\tot,R}$ be the $R$-subalgebra of 
$\mathbb H^{r}_{g,\tot,R}$
generated by the elements in Equation (\ref{EQ090}): it is
of finite type over $R$. As $R$ is universally Japanese, the
normalization $\mathbb B$ of $\mathbb H^{r,\Theta}_{g,\tot,R}$ in $\Frac(\mathbb H^{r}_{g,\tot,R})$ is a
finitely generated $R$-algebra.
The rings $\mathbb H^{r,\Theta}_{g,\tot,R}$ and $\mathbb H^r_{g,\tot,R}$ are $\frac{1}{2}\mathbb Z_{\geq 0}$-graded. Then
 $\mathbb B$ is a $\frac{1}{2}\mathbb Z_{\geq 0}$-graded $R$-subalgebra of
$\mathbb H^{r}_{g,\tot,R}$ as the latter is normal.
From Theorem \ref{T30} (a) of Subsubsection \ref{S4208} we get that $\mathbb B\otimes_R
K=\mathbb H^r_{g,\tot,K}=\mathbb H^r_{g,\tot,R}\otimes_R K$ and the image $\tilde{\mathbb B}$ of $\overline{\mathbb B}$ in
$\overline{\mathbb H^r_{g,\tot,R}}\subset \mathbb H^r_{g,\tot,k}$ is such that the $k$-homomorphism $\tilde{\mathbb B}\to \mathbb H^r_{g,\tot,k}$
is finite. Thus $\overline{\mathbb H^r_{g,\tot,R}}$ is a finitely generated $k$-algebra (being a
submodule of the finite module $\mathbb H^r_{g,\tot,k}$ over the finitely generated $k$-algebra $\tilde{\mathbb B}$). We pick a set of homogeneous generators for this $k$-algebra, lift them to a set of homogeneous elements
of $\mathbb H^r_{g,\tot,R}$, and enlarge this set to generate a $\frac{1}{2}\mathbb Z_{\geq 0}$-graded
$R$-subalgebra $\mathbb A$ of $\mathbb H^r_{g,\tot,R}$ with the properties that we have $\mathbb A\otimes_R
K=\mathbb H^r_{g,\tot,K}=\mathbb H^r_{g,\tot,R}\otimes_R K$ and the $k$-homomorphism $\overline{\mathbb A}\to \mathbb \overline{\mathbb H^r_{g,\tot,R}}$ is
surjective.
We claim that $\mathbb A= \mathbb H^r_{g,\tot,R}$ and this will complete the argument that $\mathbb H^r_{g,\tot,R}$ is a finitely generated $R$-algebra. 
To check the claim, let $s\in \frac{1}{2}\mathbb Z_{\geq 0}$ and let $\mathbb A(s)$ and 
$\mathbb H^r_{g,\tot,R}(s)$ be the corresponding graded $R$-module pieces. Then $\mathbb A(s)\subset \mathbb H^r_{g,\tot,R}(s)$ 
becomes an isomorphism after tensoring with $K$ and it becomes surjective after tensoring with $k$. As $\mathbb A(s)$ is a finitely generated $R$-module, we get that $\mathbb A(s)=\mathbb H^r_{g,\tot,R}(s)$. As this is true for all $s\in \frac{1}{2}\mathbb Z_{\geq 0}$, we conclude that our claim holds.

Part (c) follows from Corollary \ref{C21} of Subsubsection \ref{S4210}.\endproof

\medskip

For explicit formulas for $D_{2,2}(t)$, $D_{2,3}(t)$ and $D_{2,4}(t)$ see Corollary \ref{C26} (a), (b) and (c) (respectively) of Subsubsection \ref{S433}. 

We continue by providing a construction of elements in some of the modules introduced above. Let
$$\psi(t,t'):=\frac{1}{p}\log \frac{1+t^p+pt'}{(1+t)^p}\in R[[t,t']],$$
where $t,t'$ as above are two indeterminates. Let 
$$\Psi\in \pmb{\Mat}_g(R[[t,t']])$$ 
be such that its $ij$ entry is 
$$\Psi_{ij}:=\psi(T_{ij},T'_{ij})=\frac{\log(\frac{1+T_{ij}^{\phi}}{(1+T_{ij})^p})}{p}
$$
 and set 
$$\mathcal F:=\det(\Psi)\in R[[t,t']].$$
Then $\mathcal F$ satisfies 
$$\mathcal F^M=\det(M)^{-2}d^g \cdot \mathcal F$$
and thus $\mathcal F\in\mathcal H^1_g(U_0,1)=\mathcal H^1_g(U_0,1)_R$.
On the other hand, recall that $K=\Frac(R)$ and consider the logarithm of the 
multiplicative formal group
$$\ell(t):=\log(1+t)=\sum_{n=1}^{\infty}(-1)^{n+1}\frac{t^n}{n}\in K[[t]].$$
Let $\ell(T):=(\ell(T_{ij}))_{1\leq i,j\leq g}\in\pmb{\Mat}_g(K[[t]])$. Then $\tilde{\mathcal F}:=\det(\ell(T))$ satisfies 
$$\tilde{\mathcal F}^M=\det(M)^{-2}d^g\cdot \tilde{\mathcal F}$$
for all $M\in U_0$ and thus $\tilde{\mathcal F} \in\mathcal H^0_g(U_0,1)_K\subset \mathcal H^r_g(U_0,1)_K$ for all integers $r\geq 0$. But for all integers $r\geq 0$ we have $\tilde{\mathcal F}\not\in\mathcal H^r_g(U_0,1)_R$.

More generally, consider the formal power series
$$\Psi_{m_0,m_1\ldots,m_r}\in K[[T,\ldots,T^{(r)}]]:=K[[T^{(l)}_{ij}|1\leq i\leq j\leq g, 0\leq l\leq r]]$$
defined by the equality
$$\det(z_0 \ell(T)+z_1 \Psi+z_2 \Psi^{\phi}+\cdots+z_r\Psi^{\phi^{r-1}})$$
$$=\sum_{(m_0,\ldots, m_r)\in\Delta(g,r)} z_0^{m_0}z_1^{m_1}\ldots z_r^{m_r}\Psi_{m_0,m_1\ldots,m_r}$$
in which $z_0,\ldots,z_r$ are $r+1$ indeterminates. In other words
$$\Psi_{m_0,m_1\ldots,m_r}:=\Theta_{m_0,m_1,\ldots,m_r}(\ell(T),\Psi,\ldots, \Psi^{\phi^{r-1}}),$$
where $\Theta_{m_0,m_1,\ldots,m_r}\in\mathbb H^r_g(1)$ was defined in Equation (\ref{EQ090}).
Then clearly
$$\Psi_{m_0,m_1,\ldots,m_r}\in\mathcal H^r_g(U_0,1)_K.$$
Note that, in the notation above, we also have
\begin{equation}\label{EQ091}
\Psi,\Psi^{\phi},\ldots,\Psi^{\phi^{r-1}}\in\mathcal H^r_{gg}(U_0)_R,\end{equation}
\begin{equation}\label{EQ092}
\ell(T)\in\mathcal H^r_{gg}(U_0)_K\setminus \mathcal H^r_{gg}(U_0)_R.\end{equation}

\subsubsection{The  maps $\spadesuit$, $\clubsuit$, and $\heartsuit$}\label{S3602}

\begin{prop}\label{P5}
There exist well-defined $K$-linear maps
$$\spadesuit:\mathbb H^r_g(s)_K\rightarrow \mathcal H^r_g(U_0,s)_K$$
given by the formula
$$F\mapsto F^{\spadesuit}:=F(\ell(T),\Psi,\ldots, \Psi^{\phi^{r-1}}).$$
\end{prop}

\noindent
{\it Proof.}
This follows using Equations (\ref{EQ091}) and (\ref{EQ092}) via a computation similar to the one in the proof of Proposition \ref{P4}.
\endproof

\medskip
 
In order to define our next `comparison map' $\clubsuit$ we need some preparation.

\begin{lemma}\label{L13}
Let $n\in\mathbb N$. Let $a_1,\ldots,a_n\in\mathbb Z_p$ and $x=\{x_1,\ldots,x_n\}$ be a set of indeterminates. Then for every integer $r\geq 0$ the following identity holds 
$$\delta^r (\sum_{m=1}^n a_m x_m)-\sum_{m=1}^n a_m x_m^{(r)}\in (x,x',\ldots,x^{(r)})^2$$
in the ring $\mathbb Z_p[x,x',\ldots,x^{(r)}]:=\mathbb Z_p[x_m,x'_m,\ldots,x_m^{(r)}|1\leq m\leq n]$.
\end{lemma}

The proof of the lemma is an easy induction on $n\in\mathbb N$. 

With $T$ the universal symmetric $g\times g$ matrix as above and $r\geq 0$ an integer we immediately get:

\begin{lemma}\label{L14}
For all $M\in U_0$ we have an equality of $g\times g$ matrices
$$\delta^r(T^M)=M^{-1} \cdot T^{(r)} \cdot (M^{\dagger})^{\t}+\C^M_2,$$
where the entries of the `correction factor' $\C^M_2$ belong to the squared ideal $(T^{(l)}_{ij}|1\leq i\leq j\leq g, 0\leq l\leq r)^2$ of the $\mathbb Z_p$-subalgebra $\mathbb Z_p[[T,\ldots,T^{(r)}]]$ (i.e., $\mathbb Z_p[[T^{(l)}_{ij}|1\leq i\leq j\leq g, 0\leq l\leq r]]$) of $S^r_{\big-for}$.
\end{lemma}

For each commutative ring $B$, integer $a\geq 0$, and formal power series $F\in B[[T,\ldots,T^{(r)}]]:=B[[T^{(l)}_{ij}|1\leq i\leq j\leq g, 0\leq l\leq r]]$, let the degree $a$ component $F_a$ of $F$ be the sum of all terms of $F$ of total degree $a$. If $F\neq 0$, the initial component $F^{\clubsuit}$ of $F$ will be the non-zero component of $F$ which is of smallest total degree, and we set $0^{\clubsuit}:=0$. More generally, for $F$ a matrix whose entries are formal power series in $B[[T,\ldots,T^{(r)}]]$, let $F_a$ be the matrix whose entries are the $a$-components of the entries of the matrix $F$, and we define $F^{\clubsuit}$ similarly. 

 Directly from Lemma \ref{L13} we get:
 
\begin{lemma}\label{L15}
For an integral domain $B$ that contains $\mathbb Z_p$ and all $s\in \frac{1}{2}\mathbb Z$ with $gs\in \mathbb Z$ the following four properties hold:
 
\medskip
 {\bf (a)}
 If $F\in\mathcal H^r_g(U_0,s)_B$, then $F^{\clubsuit}\in H^r_g(U_0,s)_B$.
 
 \smallskip
 {\bf (b)}
 If $F\in\mathcal H^r_{gg}(U_0)_B$, then $F^{\clubsuit}\in H^r_{gg}(U_0)_B$.
 
 \smallskip
 {\bf (c)}
 If $F\in\mathcal H^r_{gg}(U_0)^*_B$, then $F^{\clubsuit}\in H^r_{gg}(U_0)^*_B$.
 
 \smallskip
 {\bf (d)}
 If $F\in\mathcal H^r_{gg}(U_0)^{**}_B$, then $F^{\clubsuit}\in H^r_{gg}(U_0)^{**}_B$.

 \end{lemma}
 
 \medskip
 
 From Lemmas \ref{L10} and \ref{L15} we get:
 
\begin{cor}\label{C6}
Assume $(A_0,\theta_0)$ satisfies condition \textup{(UAX1)} and let $B$ be an integral domain that contains $\mathbb Z_p$. Then for all $s\in \frac{1}{2}\mathbb Z$ with $gs\in \mathbb Z$
the following four properties hold:

\medskip
{\bf (a)}
For each non-zero $F\in\mathcal H^r_g(U_0,s)_B$, the polynomial $F^{\clubsuit}\in H^r_g(U_0,s)_B$ is the component of degree $gs$ of $F$ (i.e., $\deg(F^{\clubsuit})=gs$). Thus 
 the $B$-linear map 
 $$\clubsuit: \mathcal H^r_g(U_0,s)_B\rightarrow H^r_g(U_0,s)_B,\ \ F\mapsto F^{\clubsuit}$$
 is injective
 and hence $\mathcal H^r_g(U_0,s)_B$ has finite rank bounded by the rank of $H^r_g(U_0,s)_B$; in particular, we have $\mathcal H^r_g(U_0,s)_B=0$ for $s<0$ and $\mathcal H^r_g(U_0,0)$ has rank $\leq 1$. Moreover, for every $F_1\in\mathcal H^r_g(U_0,s_1)_B$ and $F_2\in\mathcal H^r_g(U_0,s_2)_B$ we have 
 $F_1F_2\in\mathcal H^r_g(U_0,s_1+s_2)_B$ and
 $(F_1F_2)^{\clubsuit}=F_1^{\clubsuit}F_2^{\clubsuit}$.
 
 \smallskip
 {\bf (b)} For each non-zero $F\in\mathcal H^r_{gg}(U_0)_B$, the matrix $F^{\clubsuit}\in H^r_{gg}(U_0)_B$ is the component of degree $1$ of $F$ (i.e., we have $\deg(F^{\clubsuit})=1$). Thus 
 the $B$-linear map 
 $$\clubsuit:\mathcal H^r_{gg}(U_0)_B\rightarrow H^r_{gg}(U_0)_B,\ \ \ F\mapsto F^{\clubsuit}$$ is injective and hence $\mathcal H^r_{gg}(U_0)_B$ has finite rank bounded by the rank of $H^r_{gg}(U_0)_B$. 
 
 \smallskip
 {\bf (c)} For each $F\in \mathcal H^r_{gg}(U_0)^*_B$, the matrix $F^{\clubsuit}\in H^r_{gg}(U_0)^*_B$ is the component of degree $0$ of $F$. Thus 
 the $B$-linear map 
 $$\clubsuit:\mathcal H^r_{gg}(U_0)^*_B\rightarrow H^r_{gg}(U_0)^*_B,\ \ F\mapsto F^{\clubsuit}$$ is injective.
 
 \smallskip
 {\bf (d)} We have $\mathcal H^r_{gg}(U_0)^{**}_B=0.$
 \end{cor}
 
 As we shall see later, the $B$-linear map $\mathcal H^r_g(U_0,s)_B\rightarrow H^r_g(U_0,s)_B$ is not surjective in general (even if $g=1$). 
 
From Propositions \ref{P3} (b) and \ref{P4} and Corollary \ref{C6} (a) we get directly:

\begin{cor}\label{C7}
Let $w\in W$. If $\deg(w)>0$, then $\mathbb I^r_{g,\ord}(w)=0$. If $\deg(w)=0$, then $\mathbb I^r_{g,\ord}(w)$ has rank at most $1$. For $w\in W(r)$, if $s:=-\deg(w)/2\in \frac{1}{2}\mathbb N$ and $gs\in \mathbb Z$, then we have a natural injective $R$-linear map
 $$\heartsuit:\mathbb I^r_{g,\ord}(w)\stackrel{\mathcal E}{\rightarrow}
 \mathcal H^r_g(U_0,s)\stackrel{\clubsuit}{\rightarrow} \mathbb H^r_g(s)$$
and thus the rank of the free $R$-module $\mathbb I^r_{g,\ord}(w)$ is less or equal to $D_0(g,r,s)$. 
\end{cor}

We recall that, by definition,  $\mathbb I^r_{g,\ord}(w)=0$ if $\frac{g\cdot \deg(w)}{2}\not\in \mathbb Z$, see
Remarks \ref{R6.85} (d)  and \ref{R7.2}.

\begin{lemma}\label{L16}
For $s\in \frac{1}{2}\mathbb N$ with $gs\in \mathbb Z$ and $r\geq 1$ an integer, the composite $K$-linear map
$$\mathbb H^r_g(s)_K\stackrel{\spadesuit}{\rightarrow} \mathcal H^r_g(U_0,s)_K\stackrel{\clubsuit}{\rightarrow} \mathbb H^r_g(s)_K$$
is given by the formula
$$F\mapsto (F^{\spadesuit})^{\clubsuit}=F(T,T'-T,p(T''-T'),\ldots, p^{r-1}(T^{(r)}-T^{(r-1}))).$$
In particular, this composite is injective, hence a $K$-linear isomorphism.
Thus the $K$-linear maps $\clubsuit$ and $\spadesuit$ are isomorphisms.
\end{lemma}

\noindent
{\it Proof.}
This follows from the fact that the linear part of $\ell(T)$ is $T$ and the linear part 
of 
$\Psi^{\phi^i}$ equals $p^i(T^{(i+1)}-T^{(i)})$.
\endproof

\medskip

The following lemma on `elimination of logarithms' will play a key role later on:

\begin{lemma}\label{L17}
Let $r\geq 1$ be an integer. Let $F\in R[T,\ldots,T^{(r)}]$ be a homogeneous polynomial such that we have
\begin{equation}\label{EQ093}
F(\ell(T),\Psi,\Psi^{\phi},\ldots,\Psi^{\phi^{r-1}})\in S^r_{\big-for}\otimes_R K.\end{equation}
Then $F\in R[T',\ldots,T^{(r)}]:=R[T^{(l)}_{ij}|1\leq i\leq j\leq g, 1\leq l\leq r]$.
\end{lemma}

\noindent
{\it Proof.} Let $\mathcal S_r$ be the set 
of all homogeneous polynomials
$$F\in R[T,\ldots,T^{(r)}]\backslash R[T',\ldots,T^{(r)}]$$
 such that Equation (\ref{EQ093}) holds. We want to show $\mathcal S_r=\emptyset$.
 Assume $\mathcal S_r\neq \emptyset$ for some $r\geq 1$ and seek a contradiction.
 Let $r$ be the minimum of all positive integers for which $\mathcal S_r\neq \emptyset$ and let $F\in \mathcal S_r$ have degree
 $$\deg(F)=c:=\min\{\deg(G)|G\in \mathcal S_r\}.$$
 Write $F=F_0+F_1$ where $F_0\in R[T',\ldots,T^{(r)}]$ and $F_1$ is in the ideal
 $(T)$ of $R[T,\ldots,T^{(r)}]$ generated by the entries of $T$. As
 $F\in \mathcal S_r$ we get that $F_1\in \mathcal S_r$. Replacing $F$ by $F_1$ we can assume that, in addition,
 $F\in (T)$. For $1\leq i,j\leq g$ we compute
 $$\begin{array}{rcl}
 \frac{\partial}{\partial T_{ij}^{(r)}} (\Psi_{ij}^{\phi^{r-1}}) & = & \frac{1}{p}\cdot 
 \frac{\partial}{\partial T_{ij}^{(r)}} \left(\log\left( \frac{1+T_{ij}^{\phi^r}}{(1+T_{ij}^{\phi^{r-1}})^p}
 \right)\right)\\
 \ & \ & \ \\
 \ & = & \frac{1}{p} \cdot \left( \frac{1+T_{ij}^{\phi^r}}{(1+T_{ij}^{\phi^{r-1}})^p}
 \right)^{-1} \cdot \frac{\partial}{\partial T_{ij}^{(r)}} \left( \frac{1+T_{ij}^{\phi^r}}{(1+T_{ij}^{\phi^{r-1}})^p}
 \right)\\
 \ & \ & \ \\
 \ & = & \frac{1}{p} \cdot (1+T_{ij}^{\phi^r})^{-1} \cdot \frac{\partial}{\partial T_{ij}^{(r)}}
 (T_{ij}^{\phi^r}) \ \\
 \ & \ & \ \\
 \ & = & p^{r-1} \cdot (1+T_{ij}^{\phi^r})^{-1}. \end{array}$$
 Hence
 $$
 \frac{\partial}{\partial T_{ij}^{(r)}}(F(\ell(T),\Psi,\ldots,\Psi^{\phi^{r-1}}))=
 \frac{\partial F}{\partial T_{ij}^{(r)}}(\ell(T),\Psi,\ldots,\Psi^{\phi^{r-1}})\cdot p^{r-1} \cdot (1+T_{ij}^{\phi^r})^{-1}.
 $$
As the right hand side of the above equality belongs to $S^r_{\big-for}\otimes_R K$ it follows that
 $$
 \frac{\partial F}{\partial T_{ij}^{(r)}}(\ell(T),\Psi,\ldots,\Psi^{\phi^{r-1}})\in S^r_{\big-for}\otimes_R K.
 $$
As $\frac{\partial F}{\partial T_{ij}^{(r)}}$ is homogeneous of degree $c-1$, by the minimality of $F$, we must have
 $$\frac{\partial F}{\partial T_{ij}^{(r)}}\in R[T',\ldots,T^{(r)}].$$
 On the other hand, as $r\geq 1$ and $F\in (T)$ we also have 
 $$\frac{\partial F}{\partial T_{ij}^{(r)}}\in (T)\subset R[T,\ldots,T^{(r)}].$$
 We conclude that for all integers $1\leq i,j\leq g$ we have 
 $$\frac{\partial F}{\partial T_{ij}^{(r)}}=0.$$
 If $r\geq 2$ we get $F\in \mathcal S_{r-1}$, contradicting the minimality of the positive integer $r$.
 If $r=1$ we get $F\in R[T]$. As $F(\ell(T))\in R[[T]]\otimes_R K$ it follows that
 $$
 \frac{\partial}{\partial T_{ij}}(F(\ell(T)))=
 \frac{\partial F}{\partial T_{ij}}(\ell(T))\cdot (1+T_{ij})^{-1}\in R[[T]]\otimes_R K.
 $$
 Hence
 $$\frac{\partial F}{\partial T_{ij}}(\ell(T))\in R[[T]]\otimes_R K.
 $$
 By the minimality of $c$ we get
 $$\frac{\partial F}{\partial T_{ij}}\in R[T']\cap R[[T]]=R,$$
 hence $c=1$. So $F(\ell(T))=\sum \lambda_{ij}\ell(T_{ij})\in R[[T]]\otimes_RK$, with $\lambda_{ij}\in R$, hence $F=0$, a contradiction.
 \endproof

\subsubsection{Applications to size $g$ forms}\label{S3603}

Next we briefly turn our attention to size $g$ forms. 
For $B$ an infinite integral domain we recall that $\mathbb H^0_{gg}(1)_B$ is the $B$-module of all matrices 
$F\in\pmb{\Mat}_g(B[T])$ such that for all $X\in\pmb{\Mat}_g(B)$, inside $\pmb{\Mat}_g(B[T])$ we have an identity 
\begin{equation}\label{EQ094}
F(XTX^{\t})=X\cdot F(T)\cdot X^{\t}.\end{equation}
 
\begin{lemma}\label{L18}
Let $B$ be a field of characteristic $0$. Then the $B$-vector space $\mathbb H^0_{gg}(1)_B$
has dimension $1$ with basis the matrix $T$.
\end{lemma}

\noindent
{\it Proof.}
We can assume that $B$ is algebraically closed. Setting $T=1_g$ and taking $X\in\pmb{\O}_g(B)=\{X\in\pmb{\Mat}_g(B)|XX^{\t}=1_g\}$ in Equation (\ref{EQ094}), we get
$$F(1_g)=X \cdot F(1_g)\cdot X^{-1}.$$
So $F(1_g)$ is a matrix that commutes with all matrices in $\pmb{\O}_g(B)$.
We claim that $F(1_g)$ must be a scalar matrix. Indeed the matrix $F(1_g)$ induces an endomorphism of the $\pmb{\O}_g(B)$-module $B^g$. This module is known to be simple;
this is trivial to check for $g\in\{1,2\}$ while for $g\geq 3$ this follows from the fact that the $\pmb{so}_{g,B}$-module $B^g$ is simple (for instance, see \cite{FH}, Thms. 19.2 and 19.4). Then our claim follows from Schur's Lemma. 

In view of our claim we can write $F(1_g)=\lambda\cdot 1_g$ with $\lambda\in B$. Setting
 $T=1_g$ in Equation (\ref{EQ094}), we get that for all $X\in\pmb{\GL}_g(B)$ we have
$$F(XX^{\t})=\lambda X X^{\t}.$$
By Sylvester's theorem, each symmetric matrix in $\pmb{\GL}_g(B)$ can be written in the form $XX^{\t}$ for some 
$X\in\pmb{\GL}_g(B)$. It follows that $F(Y)=\lambda Y$ for every symmetric matrix $Y\in\pmb{\GL}_g(B)$
hence $F=\lambda \cdot T$.
\endproof

\begin{prop}\label{P6}
Assume $B$ is a field of characteristic $0$. Then $\mathbb H^r_{gg}(1)_B$
has a basis consisting of the matrices $T,T',\ldots,T^{(r)}$.
\end{prop}

 \noindent
{\it Proof.} Let $G\in\mathbb H^r_{gg}(1)_B$; so, for all $M\in\pmb{\GL}_g(B)$,
 Equation (\ref{EQ088}) holds.
As $G$ is homogeneous of degree $1$ we have
 \begin{equation}\label{EQ095}
G(T,T',\ldots,T^{(r)})=G(T,0,\ldots,0)+G(0,T',0,\ldots,0)+\cdots+G(0,0,\ldots,T^{(r)}).\end{equation}
Setting all but one of the matrices $T,T',\ldots,T^{(r)}$ in Equation (\ref{EQ088}) equal to $0$, from Lemma \ref{L18} we get that the terms in the right-hand side of Equation (\ref{EQ095}) have the form
$b_0T,b_1T',\ldots,b_rT^{(r)}$
respectively, with $b_0,\ldots, b_r\in B$, and we are done by Equation (\ref{EQ095}).
\endproof

\begin{cor}\label{C8}
Assume conditions \textup{(UAX1)} and \textup{(UAX2)} hold.
Then for every integer $r\geq 0$ the $R$-module $\mathcal H_{gg}^r(U_0)_R$ has rank $r$ and hence for each element $F$ in it 
there exists a non-zero element $\lambda\in R$ such that $\lambda F$ is an $R$-linear combination of 
$\Psi,\Psi^{\phi},\ldots,\Psi^{\phi^{r-1}}$.
\end{cor}

\noindent
{\it Proof.}
The elements $\ell(T),\Psi,\ldots,\Psi^{\psi^{r-1}}$ of the $K$-vector space $\mathcal H^r_{gg}(U_0)_K$ are linearly independent.
By Corollary \ref{C6} (b) we have
$$\dim_K(\mathcal H^r_{gg}(U_0)_K)\leq \dim_K(H^r_{gg}(U_0)_K).$$
By Proposition \ref{P4} we have
$$\dim_K(H^r_{gg}(U_0)_K)=\dim_K(\mathbb H^r_{gg}(1)_K).$$
By Proposition \ref{P6} we have
$$\dim_K(\mathbb H^r_{gg}(1)_K)=r+1.$$
So $\ell(T),\Psi,\ldots,\Psi^{\psi^{r-1}}$ is a basis of $\mathcal H^r_{gg}(U_0)_K$.
Let $F\in\mathcal H_{gg}^r(U_0)_R$. Then there exists a non-zero $\lambda\in R$ and there exist
$\lambda_0,\ldots,\lambda_r\in R$ such that
$$\lambda F=\lambda_0\ell(T)+\lambda_1\Psi+\cdots+\lambda_r\Psi^{\phi^{r-1}}.$$
So $\lambda_0\ell(T)\in S^0_{\for}$. This clearly implies 
 $\lambda_0=0$ which ends our proof.
\endproof

\medskip

{\it From now on throughout the paper, we will assume, for simplicity, that
$(A_0,\theta_0,e)$ is as in Example \ref{EX3} and $e$ is chosen such that
$\mathcal E_{E_0,\theta_0,e}(f^1_{1,\crys})=\psi$. This is possible by \cite{Bu05}, Ch. 8, Sect. 8.1, pp. 239--240. In particular conditions \textup{(UAX1)} and \textup{(UAX2)} hold.}

\medskip

For the next result we note that, if $t$ is an indeterminate, then the base change of the universal formal abelian scheme of relative dimension $g$ over $\Spf(R_g/(T_{ij}-T_{ji}|1\leq i<j\leq g))$ which lifts $(A_0,\theta_0)$ via the morphism $$\Spf(R[[t]])\rightarrow \Spf(R_g/(T_{ij}-T_{ji}|1\leq i<j\leq g))$$ 
induced by $T_{ij}\mapsto t\cdot \delta_{ij}$ is the $g$-fold product of the universal formal elliptic curve
over $\Spf(R[[t]])$ which lifts $E_0$ with its product polarization. 
Then substituting
\begin{equation}\label{EQ096}
(T,T',\ldots,T^{(r)})\;\;\textup{by}\;\; (t\cdot 1_g,t'\cdot 1_g,\ldots,t^{(r)}\cdot 1_g)\end{equation}
 in the Serre--Tate expansion matrix $\mathcal E_{A_0,\theta_0,e}(f^r_{g,\crys})$, where $t',t'',\ldots$ are indeterminates that play the role of $T', T'',\ldots$ for $g=1$, we get
$$\mathcal E_{E_0,\theta_0,e}(f^r_{1,\crys})\cdot 1_g\in\pmb{\Mat}_g(R[[t,t',\ldots,t^{(r)}]]).$$

\begin{thm}\label{T13}
For all integers $r\geq 1$ the following two properties hold:

\medskip
{\bf (a)} The form $f^r_{g,crys}\in\mathbb I^r_{gg}(\phi^r,1)$ has the Serre--Tate expansion 

$$\mathcal E_{A_0,\theta_0,e}(f^r_{g,crys})=\Psi^{\phi^{r-1}}+p\Psi^{\phi^{r-2}}+\cdots+p^{r-1}\Psi.$$

\smallskip
{\bf (b)} The free $R$-modules $\mathbb I^1_{gg}(\phi,1)$ and $\mathbb I^1_{gg,\ord}(\phi,1)$ are equal and have rank $1$, with basis $f^1_{g,\crys}$.

\smallskip
\end{thm}

\noindent
{\it Proof.}
Part (a) was proved for $g=1$ in \cite{Bu05}, Ch. 8, Subsect. 8.4.3, Prop. 8.61. For arbitrary $g\geq 1$, from Corollary \ref{C8} we get that there exists $\lambda_0,\ldots,\lambda_r\in R$, such that we have an equality
$$\mathcal E_{A_0,\theta_0,e}(f^r_{g,crys})=\lambda_0\Psi+\cdots+\lambda_r\Psi^{\phi^{r-1}}.$$
Making the substitution (\ref{EQ096}) in this equality
 we get that
$$\mathcal E_{A_0,\theta_0,e}(f^r_{1,crys})=\lambda_0\psi+\cdots+\lambda_r\psi^{\phi^{r-1}}.$$
By the genus $1$ case of the theorem we get that 
 $$(\lambda_0,\ldots,\lambda_r)= (p^{r-1},\ldots,p,1),$$
 from which part (a) follows.
 
To prove part (b) we note that, by Corollary
\ref{C8}, $\mathcal H^1_{gg}(U_0)$ has rank $1$ and we conclude by Proposition \ref{P3} (b), and part (a). \endproof

\medskip

Similarly we have the following basic properties:

\begin{thm}\label{T14}
The following three properties hold:

\medskip
{\bf (a)} The Serre--Tate expansion of $f^{\partial}_{g,\crys}\in\mathbb I^1_{gg,\ord}(-\phi,1)$ equals $1_g$.

\smallskip
{\bf (b)} The free $R$-module $\mathbb I^1_{gg,\ord}(-\phi,1)$ has basis $f^{\partial}_{g,\crys}$.

\smallskip

{\bf (c)} Moreover for all $a,b\in\{0,1,\ldots,r\}$ we have the
following rank results:

\smallskip
{\bf (c.i)} the $R$-module $\mathbb I^r_{gg,\ord}(\phi^a,\phi^b)$
  has rank $r$;

  \smallskip
{\bf (c.ii)} the $R$-modules $\mathbb I^r_{gg,\ord}(-\phi^a,\phi^b)$
and $\mathbb I^r_{gg,\ord}(\phi^a,-\phi^b)$ have rank $1$;

  \smallskip

  {\bf (c.iii)} the $R$-module $\mathbb I^r_{gg,\ord}(-\phi^a,-\phi^b)$
  is $0$.

\end{thm}

\smallskip

\noindent
{\it Proof.} 
From Propositions \ref{P3} (b) and \ref{P4} and Corollary \ref{C8} we get that the free $R$-module 
$\mathbb I^1_{gg,\ord}(-\phi,1)$
 has rank at most $1$.
As it contains $f^{\partial}_{g,\crys}$, by the torsion freeness part in Proposition \ref{P3} (b) we get that it is equal to $\mathbb I^1_{gg}(-\phi,1)$ and this latter $R$-module has rank $1$. Also, as $\mathcal H^1_{gg}(U_0)^*$ contains 
$1_g$, by Propositions \ref{P3} (b) and \ref{P4} plus Corollary \ref{C6} (c), we have $\mathcal H^1_{gg}(U_0)^*=R1_g$. So the Serre--Tate expansion of $f^{\partial}_{g,\crys}$ equals $\eta \cdot 1_g$ for some $\eta\in R$.
We claim that $\eta=1$. It suffices to check this for $g=1$, and this was checked in \cite{Bu05}, Prop. 8.59. So parts (a) and (b) are proved.

 Part (c.i) can be proved by a similar argument: one notes that, by Corollary
\ref{C8}, $\mathcal H^r_{gg}(U_0)$ has rank $r$ and one concludes by Proposition \ref{P3} (b), Corollary \ref{C2}, and Theorem \ref{T13} (a). 

Part (c.ii) follows similarly by combining part (a), 
 Propositions \ref{P3} (b) and \ref{P4}, Corollary \ref{C6} (c), and Lemma \ref{L8}.

Part (c.iii) follows similarly 
by combining Proposition \ref{P3} (b) and 
Corollary \ref{C6} (d).\endproof

\begin{cor}\label{C9}
 For all integers $a\geq 1$, with the forms $f^{\langle a \rangle}$ and $f^{[a]}$ as in Equations (\ref{EQ064}) and (\ref{EQ066}), the following four properties hold:

\medskip
{\bf (a)} The Serre--Tate expansion of $f^{\langle a \rangle}$ equals $\Psi^{\phi^{a-1}}$.

\smallskip
{\bf (b)} The Serre--Tate expansion of $f^{[a]}$ equals $\Psi^{\phi^{a-1}}+p\Psi^{\phi^{a-2}}+\cdots+p^{a-1}\Psi$.

\smallskip
{\bf (c)} The identities $f^{\langle a \rangle t}=f^{\langle a \rangle}$ and $f^{[a]t}=f^{[a]}$ hold.

\smallskip
{\bf (d)} If $a\geq 2$, then we have an equality
$$f^{[a]}=f^{\langle a \rangle}+pf^{\langle a-1 \rangle}+\cdots+p^{a-1}f^{\langle 1 \rangle}.$$
\end{cor}

\noindent
{\it Proof.} Parts (a) and (b) follow from Theorems \ref{T13} (a) and \ref{T14} (a).
To check part (c), note that the forms $f^{\langle a \rangle t}$ and $f^{\langle a \rangle}$ both belong to $\mathbb I^a_{gg}(1,1)$ and have the same Serre--Tate expansion due to part (a) and the symmetry of the matrix $\Psi=(\Psi_{ij})_{1\leq i,j\leq g}$. We conclude that part (c) holds by the Serre--Tate expansion principle (Proposition \ref{P3} (b)). The same argument applies to $f^{[a]}$. Part (d) follows from parts (a) and (b) and the Serre--Tate expansion principle.
\endproof

\subsubsection{The maps $\diamondsuit$ and $\Omega$}\label{S3604}

For all integers $r\geq 0$ we recall the $R$-algebra $\mathbb H^r_{g,\tot}=\bigoplus_{s\in \frac{1}{2}\mathbb Z} \mathbb H^r_g(s)$, and we consider the $R$-algebra
$$\mathbb I^{r,\mathbb Z}_{g,\ord}:=\bigoplus_{s\in \frac{1}{2}\mathbb Z} \mathbb I^r_{g,\ord}(-2s)=\bigoplus_{m\in \mathbb N\cup \{0\}} \mathbb I^r_{g,\ord}(-m).$$
 Denoting $R[T,T',T'',\ldots]:=\cup_{r\geq 0} R[T,\ldots,T^{(r)}]$, we consider the $R$-algebra monomorphism
 $$\sigma:R[T,T',T'',\ldots]\rightarrow R[T,T',T'',\ldots],\ \ \ F\mapsto F^{\sigma},$$
 where
 $$T^{\sigma}:=T',\ (T')^{\sigma}:=T'',\ \textup{etc}.,\ \ \ a^{\sigma}:=a,\ \ \ a\in R,$$
 and the $\mathbb Z_p$-algebra automorphism
 $$\tau:R[T,T',T'',\ldots]\rightarrow R[T,T',T'',\ldots],\ \ \ F\mapsto F^{\tau},$$
 where
 $$T^{\tau}:=T, \ (T')^{\tau}:=T',\ \textup{etc}.,\ \ \ a^{\tau}:=\phi(a),\ \ a\in R;$$
 so for all integers $l\geq 0$ and $1\leq i\leq j\leq g$ we have $\sigma(T_{ij}^{(l)})=T_{ij}^{(l+1)}$ and $\tau(T_{ij}^{(l)})=T_{ij}^{(l)}$.
We have $\sigma\tau=\tau\sigma$. Moreover, for $r\geq 1$ we have induced $\mathbb Z_p$-linear maps
\begin{equation}\label{EQ097}
\mathbb H^{r-1}_g(s)\rightarrow \mathbb H^r_g(s),\ \ \ F\mapsto F^{\sigma},\end{equation}
\begin{equation}\label{EQ097.2}
\mathbb H^{r-1}_g(s)\rightarrow \mathbb H^{r-1}_g(s),\ \ \ F\mapsto F^{\tau}.\end{equation}
and a $\mathbb Z_p$-bilinear map
\begin{equation}
\label{EQ097.5}
\mathbb H^{r-1}_g(s)\times \mathbb H^{r-1}_g(s)\rightarrow\mathbb H^r_g((p+1)s),\end{equation}
defined by
$$ (F,G)\mapsto \{F,G\}_{\sigma,\tau}:=G^{\sigma\tau}F^p-F^{\sigma\tau}G^p.$$

\begin{prop}\label{P7}
For all integers $r\geq 1$ and all $s\in \frac{1}{2}\mathbb Z$ the following three properties hold:

\medskip
{\bf (a)} We have a well-defined $R$-linear map
$$\diamondsuit:\mathbb H^{r-1}_g(s)\rightarrow \mathbb I^r_{g,\ord}(-2s),\ \ \ F\mapsto F^{\diamondsuit}:=F(f^{\langle 1 \rangle},\ldots,f^{\langle r \rangle}).$$

\smallskip
{\bf (b)} For $F\in\mathbb H^{r-1}_g(s)$ we have the following formula:
\begin{equation}\label{EQ098}
(F^{\diamondsuit})^{\phi}=\det(f^{\partial})^{-2s}\cdot (F^{\sigma \tau })^{\diamondsuit}.\end{equation}

\smallskip
{\bf (c)} For $F,G\in\mathbb H^{r-1}_g(s)$ we have the following formula:
\begin{equation}\label{EQ099}
(\{F,G\}_{\sigma,\tau})^{\diamondsuit}=
p\cdot \det(f^{\partial})^{2s}\cdot \{F^{\diamondsuit},G^{\diamondsuit}\}_{\delta}
\end{equation}
(see the part after Definition \ref{df11} of Subsection 2.1 for the notation $\{F^{\diamondsuit},G^{\diamondsuit}\}_{\delta}$). 
\end{prop}

\noindent
{\it Proof.}
The expression $F(f^{\langle 1 \rangle},\ldots,f^{\langle r \rangle})$ is well-defined as $f^{\langle 1 \rangle},\ldots,f^{\langle r \rangle}$ are symmetric,
see Corollary \ref{C9} (c).
This expression is trivially seen to define an element in $\mathbb I^r_{g,\ord}(-2s)$. Thus part (a) holds. 

Equation (\ref{EQ098}) reads 
\begin{equation}\label{EQ100}
F(f^{\langle 1 \rangle},\ldots,f^{\langle r \rangle})^{\phi}=\det(f^{\partial})^{-2s}\cdot F^{\tau}(f^{\langle 2 \rangle},\ldots,f^{\langle r+1\rangle});\end{equation}
as Equation (\ref{EQ100}) 
follows immediately from Equation (\ref{EQ068}), part (b) holds. Part (c) holds as Equation (\ref{EQ099}) follows directly from Equation (\ref{EQ098}).
\endproof

\begin{thm}\label{T15}
For all $r\in \mathbb N$ and all $s\in \frac{1}{2} \mathbb Z$ the $R$-linear map
$$\diamondsuit:\mathbb H^{r-1}_g(s)\rightarrow \mathbb I^r_{g,\ord}(-2s)$$ 
 is injective and has a torsion cokernel. In particular:
 
 \medskip
 {\bf (a)} For $w\in W$, if $s:=-\deg(w)/2$, then we have
 $$\rank_R(\mathbb I^r_{g,\ord}(w))= D_0(g,r-1,s).$$
 
 \smallskip
 {\bf (b)} The $R$-algebra homomorphism
 $$\diamondsuit: \mathbb H^{r-1}_{g,\tot}\rightarrow \mathbb I^{r,\mathbb Z}_{g,\ord}$$
 is injective and it induces a $K$-algebra isomorphism
 $$\diamondsuit: \mathbb H^{r-1}_{g,\tot}\otimes_R K\rightarrow \mathbb I^{r,\mathbb Z}_{g,\ord}\otimes_R K.$$\end{thm}

\noindent
{\it Proof.} We have the following diagram of $R$-linear maps:
$$
\begin{array}{ccccccc}
\mathbb H^{r-1}_g(s)_K & \stackrel{\sigma}{\rightarrow} & \mathbb H^r_g(s)_K &
\stackrel{\spadesuit}{\rightarrow} & \mathcal H^r_g(U_0,s)_K & \stackrel{\clubsuit}{\rightarrow}
& \mathbb H^r_g(s)_K\\
\alpha \uparrow & \ &\ &\ & \beta\uparrow & \ & \gamma\uparrow\\
\mathbb H^{r-1}_g(s) & \stackrel{\diamondsuit}{\rightarrow} & \mathbb I^r_{g,\ord}(-2s) & 
\stackrel{\mathcal E}{\rightarrow} & \mathcal H^r_g(U_0,s) & \stackrel{\clubsuit}{\rightarrow}
& \mathbb H^r_g(s)
\end{array}
$$
where $\alpha$, $\gamma$ and $\beta$ are inclusions.
As the $R$-linear maps on the upper row are injective it follows that $\diamondsuit$ on the lower row is injective. Now for each $f\in\mathbb I^r_{g,\ord}(-2s)$ we have that 
$$\beta(\mathcal E(f))\in S^r_{\big-for}.$$
As $\spadesuit$ is an $R$-linear isomorphism (see Lemma \ref{L16}), there exists $F\in\mathbb H^r_g(s)_K$ such that $F^{\spadesuit}=\beta(\mathcal E(f))$.
By Lemma \ref{L17} we must have $F\in R[T',\ldots,T^{(r)}]$. But then, clearly $F=F_0^{\sigma}$ for some $F_0\in \mathbb H^{r-1}_g(s)_K$. We conclude that the image of the injective $R$-linear map $\beta\circ \mathcal E$ is contained in the image of the $R$-linear map $\spadesuit\circ \sigma$. Hence we have an inequality
$$\dim(\mathbb I^r_{g,\ord}(-2s)\otimes_R K)\leq \dim(\mathbb H^{r-1}_g(s)_K)$$
which, as $\diamondsuit$ is injective, is an equality. Hence the cokernel of $\diamondsuit$ is torsion. This proves part (b). Part (a) follows from part (b), Equation (\ref{Dpgrs}), and Corollary \ref{C3}.
\endproof

\begin{cor}\label{C10}
For all integers $n\geq 1$, the $R$-modules $\mathbb I^1_g(-n\phi-n)$ and $\mathbb I^1_{g,\ord}(-n\phi-n)$ coincide and are free of rank $1$ with basis $(\det(f^1_{g,\crys}))^n$.
\end{cor}

\noindent
{\it Proof.}
Theorem \ref{T15} (a) and Lemma \ref{L11} imply immediately that the statement is true
after tensorization with $K$. So to prove the corollary it is enough to check that  $(\det(f^1_{g,\crys}))^n$ is not divisible by $p$ in $\mathbb I^1_{g,\ord}(-n\phi-n)$. The latter condition follows
from the torsion freeness part of the Serre--Tate expansion principle (see Proposition \ref{P3} (b)) and the fact  that the Serre--Tate expansion of $(\det(f^1_{g,\crys}))^n$ equals $1$.\endproof

\medskip

For ordinary forms of size $g$ we have:

\begin{cor}\label{C11}
For all integers $a,b\geq 0$ and $r\geq \max\{a,b\}$, the $R$-module $\mathbb I^r_{gg,\ord}(\phi^a,\phi^b)$ has rank $r$.
\end{cor}

\noindent
{\it Proof.} Recall that by Corollary \ref{C8} the $R$-module $\mathcal H_{gg}^r(U_0)_R$ has rank $r$. 
Hence by Proposition \ref{P3} (b) the $R$-module $\mathbb I^r_{gg,\ord}(1,1)$ has rank $\leq r$. As $f^{\langle 1 \rangle},\ldots,f^{\langle r \rangle}$ belong to the latter module and are $R$-linearly independent (as their Serre--Tate expansions are, see Corollary \ref{C9} (a)) it follows that
$\mathbb I^r_{gg,\ord}(1,1)$ has rank $r$. We conclude by Corollary \ref{C2}.
\endproof

\medskip

By Theorem \ref{T15} (a) the $R$-linear maps
$\diamondsuit:\mathbb H^{r-1}_g(s)\rightarrow \mathbb I^r_{g,\ord}(-2s)$
 become surjective after inverting $p$;  but they are not  surjective in general (see Proposition \ref{P7.5} (b) below for a more precise statement). 
 
 Indeed let us consider  the torsion $R$-modules
 \begin{equation}\label{EQ100.5}
 \mathbb T^r_g(s):=\frac{\mathbb I^r_{g,\ord}(-2s)}{\diamondsuit(\mathbb H^{r-1}_g(s))},\ \ \ 
 \mathbb T^r_g:=\frac{\mathbb I^{r,\mathbb Z}_{g,\ord}}{\diamondsuit(\mathbb H^{r-1}_{g,\tot})}=\bigoplus_{s\in \frac{1}{2}\mathbb Z_{\geq 0}} \mathbb T^r_g(s).\end{equation}
 Denote by 
 $\mathbb T^r_g(s)[p]$ the set of all elements of $\mathbb T^r_g(s)$ annihilated by $p$. So
 $\mathbb T^r_g(s)[p]$ has a natural structure of $k$-vector space. 
 We will denote by $\mathbb T^r_g(s)[p]^{\tau}$ the $k$-vector space obtained 
 from $\mathbb T^r_g(s)[p]$ by restricting the scalars via the Frobenius automorphism of $k$.
 
 Consider the $\mathbb Z_p$-bilinear map
 $$\Omega:\mathbb H^{r-1}_g(s)\times \mathbb H^{r-1}_g(s)\rightarrow \mathbb T^{r+1}_g((p+1)s)[p]$$
 defined (see Proposition \ref{P7} (c)) by:
 $$\begin{array}{rcl}
 \Omega(F,G) &:= & (\det(f^{\partial}))^{2s}\{F^{\diamondsuit},G^{\diamondsuit}\}_{\delta}+\diamondsuit(\mathbb H^r_g((p+1)s))\\
 \ & \ & \ \\
 \ & = & \frac{1}{p}(\{F,G\}_{\sigma,\tau})^{\diamondsuit}+\diamondsuit(\mathbb H^r_g((p+1)s)).
 \end{array}$$
Clearly, the map $\Omega$ induces an
 antisymmetric
 $k$-bilinear map
 $$
 \overline{\Omega}:\overline{\mathbb H^{r-1}_g(s)}\times \overline{\mathbb H^{r-1}_g(s)}\rightarrow \mathbb T^{r+1}_g((p+1)s)[p]^{\tau}.$$
 For instance, additivity in the first argument follows from the formula (see Proposition \ref{P7} (b)):
 $$
 \{F_1^{\diamondsuit}+F_2^{\diamondsuit},G^{\diamondsuit}\}_{\delta}=
 \{F_1^{\diamondsuit},G^{\diamondsuit}\}_{\delta}+\{F_2^{\diamondsuit},G^{\diamondsuit}\}_{\delta}-\det(f^{\partial})^{-2s}
 (G^{\sigma\tau}\cdot C_p(F_1,F_2))^{\diamondsuit}.
 $$
To state our next result and for future use it is convenient to make the following definition. 

\begin{df}\label{df26}
Let $U,V$ be two vector spaces over an arbitrary field. A bilinear antisymmetric map
$\Omega:U\times U\rightarrow V$ is called {\it totally non-degenerate} if for all $x,y\in U$, we have $\Omega(x,y)=0$ if and only if $x$ and $y$ are linearly dependent.
\end{df}
 
 \begin{prop}\label{P7.5}
For an integer $r\geq 2$ and $s\in \frac{1}{2}\mathbb Z_{\geq 0}$ the following three properties hold:
 
 \smallskip
 
 {\bf (a)} The above bilinear map $\overline{\Omega}$  is totally non-degenerate.
 
 \smallskip
 
 {\bf (b)} We have an inequality
 $$ \dim_k(\mathbb T^{r+1}_g((p+1)s)[p])\geq 2D_0(g,r-1,s)-3.$$
 In particular, if $D_0(g,r-1,s)\geq 2$, then the map 
 $$\diamondsuit:\mathbb H^r_g((p+1)s)\rightarrow \mathbb I^{r+1}_{g,\ord}(-2(p+1)s)$$
 is not surjective.
 
 \smallskip
 
 {\bf (c)} The
annihilator of the $R$-module $\mathbb T^r_g$ is $0$. \end{prop}
 
 \noindent
 {\it Proof.} To prove part (a) assume 
 $F,G\in \mathbb H^{r-1}_g(s)$ 
 and $\overline{\Omega}(\overline{F},\overline{G})=0$; we need to show that $\overline{F}$ and $\overline{G}$ are $k$-linearly dependent. We have $\Omega(F,G)=0$ so there exists $H\in \mathbb H^r_g((p+1)s)$ such that 
 $(\{F,G\}_{\sigma,\tau})^{\diamondsuit}=p H^{\diamondsuit}$.
 By the injectivity of $\diamondsuit$ we have
 $\{F,G\}_{\sigma,\tau}=pH$.
 This implies that either $\overline{G}=0$ (in which case we are done) or $\overline{G}\neq 0$ and 
 $(\overline{F}/\overline{G})^{\sigma \tau}=(\overline{F}/\overline{G})^p$ in the field $k(T,\ldots,T^{(r)})$. The latter clearly implies $\overline{F}/\overline{G}\in k$ and we are done again.
 
 To check part (b) 
 note that by part (a) plus
 a general fact of linear algebra (see Lemma \ref{L19} below)
 the $k$-linear span of the image of $\overline{\Omega}$ has dimension at least $2\dim_k(\overline{\mathbb H^{r-1}_g(s)})-3$. We then conclude by
 the fact that
$$\dim_k(\overline{\mathbb H^{r-1}_g(s)})=\dim_K(\mathbb H^{r-1}_g(s)_K)=D_0(g,r-1,s).$$

To check part (c) note that, by part (b), $\mathbb T^r_g[p]\neq 0$ hence there exists a prime element $F$ of the unique factorization domain $\mathbb H^{r-1}_{g,\tot}$ with $\frac{F}{p}\in \mathbb I^{r,\mathbb Z}_{g,\ord}\setminus\diamondsuit(\mathbb H^{r-1}_{g,\tot})$. For each $n\in\mathbb N$ the image of $\frac{F^n}{p^n}$ in $\mathbb T_g$ has annihilator $p^nR$ in $R$. Thus the annihilator of $\mathbb T_g$ in $R$ is contained in $\cap_{n\geq 1} p^n R=0$, hence is $0$.
\endproof

\medskip

In the above proof we applied over $k$ the following general fact:

\begin{lemma}\label{L19}
Let $U$ and $V$ be finite dimensional vector spaces over an algebraically closed field $B$ and let $b:U\times U \rightarrow V$ be a $B$-bilinear antisymmetric map. Assume $b$ is totally non-degenerate. Then the following inequality holds:
$$\dim_B(\textup{Span}(\Im(b)))\geq 2 \dim_B(U)-3.$$
\end{lemma}

\noindent
{\it Proof.} Let $u_1,\ldots,u_n$ be a basis of $U$ and set $v_{ij}:=b(u_i,u_j)\in V$ for $1\leq i<j\leq n$.
Let $L$ be the vector space over $B$ of all $(n^2-n)/2$-tuples $c:=(c_{ij})_{(i,j)\in \{1,\ldots, n\}^2\setminus\{(l,l)|l\in\{1,\ldots,n\}\}}\in B^{(n^2-n)/2}$ and let 
$$L':=\{c\in L|\sum_{1\leq i<j\leq n} c_{ij}v_{ij}=0\}.$$ 
For $c\in L\setminus\{0\}$ let $[c]\in \mathbb P(L):=(L\setminus \{0\})/B^{\times}$ be the corresponding point in the projective space; similarly,
for $\lambda:=(\lambda_i)_{1\leq i \leq n}\in B^n\setminus\{0\}$ let $[\lambda]\in \mathbb P(B^n)$ be the corresponding point.
Let 
$$\varphi:(\mathbb P(B^n)\times \mathbb P(B^n))\setminus \Delta \rightarrow \mathbb P(L)$$ 
be the morphism
$$\varphi([\lambda],[\mu]):=[(\lambda_i\mu_j-\mu_i\lambda_j)_{1\leq i<j\leq n}],$$
where $\Delta$ is the diagonal of $\mathbb P(B^n)\times \mathbb P(B^n)$.
The image $\Im(\varphi)$ is the Grassmannian of $2$-planes in an $n$-dimensional vector space, so $\dim(\Im(\varphi))=2n-4$. 
If there exists $[(c_{ij})]\in \Im(\varphi)\cap \mathbb P(L')$, then $c_{ij}=\lambda_i\mu_j-\lambda_j\mu_i$ for some $(\lambda_i)_{1\leq i\leq n},(\mu_i)_{1\leq i\leq n}\in B^n\setminus \{0\}$ and for $x:=\sum_{i=1}^n \lambda_i u_i$ and $y:=\sum_{i=1}^n \mu_i u_i$ we get that $x$ and $y$ are $B$-linearly independent and $b(x,y)=0$, a contradiction. Thus we have $\Im(\varphi)\cap \mathbb P(L')=\emptyset$
and hence
$$\dim_B(L')\leq \dim_B(L)-1-\dim(\Im(\varphi))\leq \frac{n(n-1)}{2}-2n+3,$$
which implies that the lemma holds.
\endproof

\begin{rem}\label{R22}
Lemma \ref{L19} is optimal in the sense that for every vector space $U$ over $B$ of dimension at least $2$ there exists a vector space $V$ and a totally non-degenerate antisymmetric bilinear map
$b:U\times U\rightarrow V$ such that the inequality of Lemma \ref{L19} is an equality. To construct such a $V$ and $b$ fix a basis $u_1,u_2,\ldots,u_n$ of $U$ with $n\geq 2$,
let $L$ and $\varphi$ be as in the above proof, let $L'\subset L$ be a subspace of dimension
$ \frac{n(n-1)}{2}-2n+3$ such that $\mathbb P(L')\cap \Im(\varphi)=\emptyset$ (such an $L'$ exists
because $\dim(\Im(\varphi))+\frac{n(n-1)}{2}-2n+3\leq \dim_B(L)-1$)
 and take $V:=(\wedge^2 U)/U'$ where $U'$ is the subspace of $\wedge^2 U$ consisting of all expressions of the form $\sum_{i<j} c_{ij} u_i\wedge u_j$ with $(c_{ij})_{1\leq i<j\leq n}\in L'$. Finally take $b$ to be the natural quotient composition $U\times U\rightarrow \wedge^2 U\rightarrow V$. \end{rem}

In what follows we consider the rings:
$$\mathbb H_{g,\tot}:=\bigcup_{r=0}^{\infty} \mathbb H^r_{g,\tot},\ \ \ \mathbb I^{\mathbb Z}_{g,\ord}:=\bigcup_{r=0}^{\infty} \mathbb I^{r,\mathbb Z}_{g,\ord}.$$

In view of Lemma \ref{L12} (b), we get:

\begin{cor}\label{C12}
For all integers $r\geq 2$ we have
$$\trdeg_R(\mathbb I^{r,\mathbb Z}_{g,\ord})=r\frac{g(g+1)}{2}-g^2+1.$$
\end{cor}

\begin{cor}\label{C13}
We have
$$\asytrdeg(\mathbb I^{\mathbb Z}_{g,\ord})=\frac{g(g+1)}{2}.$$\end{cor}

On the other hand recall the rings
$\mathbb I^r_{g,\ord}$ 
and their union $\mathbb I_{g,\ord}$. Also recall the filtered $\delta$-ring $\mathbb I_{g,\ord,((p))}$.
 For $r\in\mathbb N$ let $z_{<r}:=\{z_0,\ldots,z_{r-1}\}$ be a set of indeterminates and we consider the $R$-algebra monomorphism between rings of Laurent polynomials in the set of indeterminates $z_{<r}$,
 $$\mathbb H_{g,\tot}^{r-1,\torus}:=\mathbb H_{g,\tot}^{r-1}[z_{<r},\frac{1}{z_0\cdots z_{r-1}}]\rightarrow \mathbb I^{r,\mathbb Z}_{g,\ord}[z_{<r},\frac{1}{z_0\cdots z_{r-1}}],$$
 that is induced by $\diamondsuit$ and sends $z_i$ to $z_i$ for all $i\in\{0,\ldots,r-1\}$.
The ring of Laurent polynomials 
 $$\mathbb H_{g,\tot}^{\torus}:=\cup_{r\geq 1} \mathbb H_{g,\tot}^{r-1,\torus}$$
in the set of indeterminates $z:=\cup_{r\in\mathbb N}z_{<r}$ is a filtered $W$-graded ring with the grading defined by letting the elements
 of $\mathbb H_g^r(s)$ have weight $-2s$ and letting $z_i$ have weight $\phi^{i+1}-\phi^i$ for each integer $i\geq 0$, see Subsubsection \ref{S234}.
 
\begin{thm}\label{T16}
The following eight properties hold:

\medskip
{\bf (a)} For all integers $r\geq 1$ we have an $\mathbb I^{r,\mathbb Z}_{g,\ord}$-algebra isomorphism 
$$\mathbb I^{r,\mathbb Z}_{g,\ord}[z_{<r},\frac{1}{z_0\cdots z_{r-1}}]\simeq \mathbb I^r_{g,\ord},$$
and hence monomorphisms of $R$-algebras, 
$$\diamondsuit:\mathbb H_{g,\tot}^{r-1,\torus}\rightarrow \mathbb I^r_{g,\ord},\ \ \ 
\diamondsuit:\mathbb H_{g,\tot}^{\torus}\rightarrow \mathbb I_{g,\ord},$$
whose cokernels are torsion.

\smallskip
{\bf (b)} For all integers $r\geq 2$ we have
$$ \trdeg_R(\mathbb I^r_{g,\ord})=r\frac{g(g+1)}{2}-g^2+1+r.$$

\smallskip
{\bf (c)} We have
$$\asytrdeg(\mathbb I_{g,\ord})=\frac{g(g+1)}{2}+1.$$

\smallskip
{\bf (d)} We have an equality $\mathbb I_{g,\ord,((p))}=\mathbb I^{\mathbb Z}_{g,\ord,((p))}$.
In particular, we have
$$\asytrdeg(\mathbb I_{g,\ord,((p))})=\frac{g(g+1)}{2}.$$

\smallskip
{\bf (e)} We have
$$\phidim(\mathbb I_{g,\ord})=\frac{g(g+1)}{2}+1.$$

\smallskip
{\bf (f)} We have
$$\phidim(\mathbb I_{g,\ord,((p))})=\frac{g(g+1)}{2}.$$

\smallskip
{\bf (g)} The field $\Frac(\mathbb I_{g,\ord})$ is $\phi$-generated over $K$ by a finite set of elements in $\mathbb I^4_{g,\ord}$.

\smallskip
{\bf (h)} For each integer $r\geq 1$, the form $\det(f^r_{g,crys})$ belongs to the image of the $R$-algebra monomorphism $\diamondsuit:\mathbb H_g^{r-1,\torus}\rightarrow \mathbb I^r_{g,\ord}$.
\end{thm}

\noindent
{\it Proof.} To check part (a) we define a $\mathbb I^{r,\mathbb Z}_{g,\ord}$-algebra homomorphism
$$\mathbb I^{r,\mathbb Z}_{g,\ord}[z_{<r},\frac{1}{z_0\cdots z_{r-1}}]\rightarrow \mathbb I^r_{g,\ord}$$
 by the rule
$$z_i\mapsto \det(f^{\partial}_{g,\crys})^{\phi^i}\in\mathbb I^r_{g,\ord}(\phi^i(\phi-1)),\ \ \ i\in\{0,\ldots,r-1\}.$$
Hence, for $s\in \frac{1}{2}\mathbb Z$, $f\in\mathbb I^r_{g,\ord}(-2s)$, and $w'=\sum_{i=0}^{r-1}a_i\phi^i\in W$ we have 
$$f\cdot z_0^{a_0}\cdots z_{r-1}^{a_{r-1}}\mapsto f\cdot (f^{\partial}_{g,\crys})^{w'}\in\mathbb I^r_{g,\ord}(w),$$
where 
$$w:=-2s+(\phi-1)w'\in W,\ \ \ \deg(w)=-2s.$$
By Corollary \ref{C3}, we have an induced $R$-linear isomorphism 
$$\mathbb I^r_{g,\ord}(-2s)\cdot z_0^{a_0}\cdots z_{r-1}^{a_{r-1}}\simeq \mathbb I^r_{g,\ord}(w).$$
Based on this, as $w$ determines the pair $(s,w')$, we conclude that the part (a) and the first sentence of part (d) hold. 

Part (b) follows from part (a) and Corollary \ref{C10}. Part (c) follows from part (b). The second sentence of part (d) follows from part (c).

For the proof of parts (e) and (f) recall our  convention (see the beginning of Subsection \ref{S22})
on algebraically independent families. Furthermore, by the cardinality of a family $(x_i)_{i\in \mathcal I}$ we will mean the cardinality of the set $\{x_i|i\in \mathcal I\}$.
Also, by the union of two families $(x_i)_{i\in \mathcal I}$ and $(y_j)_{j\in \mathcal J}$  we will understand the family $(z_k)_{k\in \mathcal I\coprod \mathcal J}$, indexed by the coproduct $\mathcal I\coprod \mathcal J$ in the category of sets where, if $\kappa_1:\mathcal I\rightarrow \mathcal I\coprod \mathcal J$ and $\kappa_2:\mathcal J\rightarrow \mathcal I\coprod \mathcal J$ are the natural injections, then 
$z_{\kappa_1(i)}:=x_i$ for all $i\in \mathcal I$ and $z_{\kappa_2(j)}:=y_j$ for all $j\in \mathcal J$. 
 
In order to check parts (e) and (f), based on part (d), Inequality (\ref{EQ008}) and Remark \ref{R6} (a) it suffices to check that $\mathbb I_{g,\ord}$ contains $\frac{g(g+1)}{2}+1$ elements that are $\phi$-algebraically independent. To check this we use the fact that, by Corollary \ref{C23} (a) of Subsection \ref{S42}, the ring $\mathbb H_g$ contains homogeneous polynomials 
$F_1,\ldots,F_{\frac{g(g+1)}{2}}$ of degree $gs$, for some $s\in \mathbb N$,
with coefficients in $\mathbb Z$, such that the family
$$(\sigma^jF_i)_{(i,j)\in\{1,\ldots,\frac{g(g+1)}{2}\}\times (\mathbb N\cup\{0\})}$$
is algebraically independent over $\mathbb Q$, therefore over $K$.
 In particular, $\sigma$ and $\tau$ act as the identity on the coefficients of these polynomials. Hence, for each integer $N\geq 0$, the family
$$((F_i^{\sigma^j})^{\diamondsuit})_{(i,j)\in\{1,\ldots,\frac{g(g+1)}{2}\}\times\{0,\ldots,N\}}$$
of elements in $\mathbb I_{g,\ord}$ 
is algebraically independent over $K$. 
 By Proposition \ref{P7} (b), for each integer $j\geq 1$, we have 
\begin{equation}\label{EQ101}
(F_i^{\sigma^j})^{\diamondsuit}=(\det(f^{\partial}))^{2s(1+\phi+\cdots+\phi^{j-1})}
\cdot ((F_i)^{\diamondsuit})^{\phi^j}.\end{equation}
Thus the family of elements of weight $-2s$
\begin{equation}\label{EQ102}
\left((\det(f^{\partial}))^{2s(1+\phi+\cdots+\phi^{j-1})}
\cdot ((F_i)^{\diamondsuit})^{\phi^j}\right)_{(i,j)\in\{1,\ldots,\frac{g(g+1)}{2}\}\times\{0,\ldots,N\}}\end{equation}
is algebraically independent over $K$.
By weight considerations one easily sees that the union of the family (\ref{EQ102}) with the family
\begin{equation}\label{EQ103}
(\det(f^{\partial})^{\phi^j})_{j\in\{0,\ldots,N-1\}}\end{equation}
is algebraically independent over $K$. We claim that the union of the family (\ref{EQ103})
with the family
\begin{equation}\label{EQ104}
(((F_i)^{\diamondsuit})^{\phi^j})_{(i,j)\in\{1,\ldots,\frac{g(g+1)}{2}\}\times\{0,\ldots,N\}}\end{equation}
is algebraically independent over $K$. Indeed this follows from the following two facts:

\medskip
{\bf (i)} The $K$-algebra generated by 
the members of the union of the families (\ref{EQ102}) and (\ref{EQ103}) is contained in the $K$-algebra generated by the members of the union of the families (\ref{EQ104}) and (\ref{EQ103}).

\smallskip
{\bf (ii)} The cardinality of the union of the families (\ref{EQ102}) and (\ref{EQ103}) is the same as the cardinality of the union of the families (\ref{EQ104}) and (\ref{EQ103}).

\medskip\noindent Hence the family
\begin{equation}\label{EQ105}
(b_i)_{i\in \{0,\ldots,{\frac{g(g+1)}{2}}\}}\end{equation}
where
\begin{equation} 
b_0:=\det(f^{\partial}),\ b_1:=(F_1)^{\diamondsuit},\ldots,b_{\frac{g(g+1)}{2}}:=(F_{\frac{g(g+1)}{2}})^{\diamondsuit},\end{equation}
is $\phi$-algebraically independent and parts (e) and (f) follow. 

To check part (g), we note that by Corollary \ref{C23} (b) and (d) of Subsection \ref{S42}, the field
$\Frac(\mathbb H_{g,\tot, \mathbb K}^3)$ contains elements 
$G_1,\ldots,G_{N_3}$ such that the family
\begin{equation}\label{EQ106}
(G_i^{\sigma^j})_{(i,j)\in\{1,\ldots,N_3\}\times (\mathbb N\cup\{0\})}\end{equation}
generates the field $\Frac(\mathbb H_{g,\tot,\mathbb K})$ over $\mathbb K$.
By enlarging $N_3\in\mathbb N$ we can assume that $G_1,\ldots,G_{N_3}$ are homogeneous elements of $\mathbb H_{g,\tot,\mathbb Q}^3$. Then the family (\ref{EQ106})
generates the field $\Frac(\mathbb H_{g,\tot,\mathbb Q})$ over $\mathbb Q$
and hence the field $\Frac(\mathbb H_{g,\tot})$ over $K$. As $\sigma$ acts trivially on the coefficients of $G_i$ we get that Equation (\ref{EQ101}) holds with $F_i$ replaced by $G_i$.
We conclude that the Hecke covariant ordinary Siegel $\delta$-modular forms
$$ \det(f^{\partial}),(G_1)^{\diamondsuit},\ldots,(G_{N_3})^{\diamondsuit}\in \mathbb I^4_{g,\ord}$$
and their images via $\phi^j$, $j\geq 0$, generate the field $\Frac(\mathbb I_{g,\ord})$ over $K$.

Part (h) follows from the fact that $\det(f^r_{g,crys})$ equals a $W$-power of $\det(f^{\partial}_{g,\crys})$ times $\det(f^{[r]})$ and thus (see Corollary \ref{C9} (d) and Equation (\ref{EQ090})) times a $\mathbb Z$-linear combination of elements of the form $\Theta_{m_0,\ldots,m_r}(f^{\langle 1 \rangle},\ldots,f^{\langle r \rangle})$. 
\endproof

\begin{cor}\label{C14}
We have an inequality
$$ \trdeg_k(\mathbb I^r_{g,\ord}\otimes_R k)\leq r\frac{g(g+1)}{2}-g^2+1+r.
$$
\end{cor}

\noindent
{\it Proof.} As $p$ is a prime element in the ring $\mathbb I^r_{g,\ord}$ (see Corollary 4 (a)), it follows by Cohen's `dimension inequality' in [Ma], Thm. 15.5 applied to the extension $R\subset
\mathbb I^r_{g,\ord}$ and to the prime ideal $p\mathbb I^r_{g,\ord}$
that $\trdeg_k(\mathbb I^r_{g,\ord}\otimes_R k)\leq \trdeg_R(\mathbb I^r_{g,\ord})$.
We conclude by Theorem \ref{T16} (b).
\endproof



\begin{cor}\label{C15}
The following seven properties hold:

\medskip
{\bf (a)} For each integer $r\geq 1$, the $K$-algebra $\mathbb I_{g,\ord}^r\otimes_R K$ is finitely generated, Cohen--Macauley, unique factorization domain, hence Gorenstein, and moreover $\mathbb I_{g,\ord}^r$ is a unique factorization domain.
  
\smallskip
{\bf (b)} If $g$ is odd, then for all $w\in W$ with $\deg(w)$ odd we have $\mathbb I_{g,\ord}^r(w)=0$.

\smallskip
{\bf (c)} If $g$ is even, then there exists $r\in\{0,1,2,3,4\}$ and $w\in W(r)$ with $\deg(w)$ odd such that $\mathbb I_{g,\ord}^r(w)\neq 0$. 

\smallskip
{\bf (d)} We have $\mathbb I^r_{g,\ord}(w)=0$ for all $w\in W$ with $\deg(w)=-1$.

\smallskip
{\bf (e)} Let $w\in W(r)$ with $\deg(w)=0$, so $w\in (\phi-1)W$. Then the $R$-module $\mathbb I^r_{g,\ord}(w)$ is free of rank $1$ with basis 
$\det(f^{\partial}_{g,\crys})^{\frac{w}{\phi-1}}$.

\smallskip
{\bf (f)} All non-zero elements in $\mathbb I_{g,\ord}(-2)\otimes_R K$ generate prime ideals in the ring
$\mathbb I_{g,\ord}\otimes_R K$. 

\smallskip
{\bf (g)} Every two $K$-linearly independent elements in $\mathbb I_{g,\ord}(-2)\otimes_R K$ generate two distinct ideals in $\mathbb I_{g,\ord}\otimes_R K$. \end{cor}

\noindent
{\it Proof.} As we have a $K$-algebra isomorphism $\mathbb H_{g,\tot}^{r-1,\torus}\otimes_R K\simeq \mathbb I^r_{g,\ord}\otimes_R K$ (see Theorem \ref{T16} (a)) to prove the $K$-algebra statement of part (a) it suffices to show that $\mathbb H_{g,\tot}^{r-1}\otimes_R K$ has all the desired properties and this follows from Subsection \ref{S4203}. 

To prove that $\mathbb I_{g,\ord}^r$ is a unique factorization domain, we first recall that $p$ is a prime element of $\mathbb I_{g,\ord}^r$ (see Corollary \ref{C4} (b)) and we will use this to check that we have an identity $\mathbb I_{g,\ord}^r=\mathbb I_{g,\ord}^r[\frac{1}{p}]\cap (\mathbb I_{g,\ord}^r)_{(p)}$ between subrings of $\Frac(\mathbb I_{g,\ord}^r)$. The inclusion "$\subset$" is clear and for the reversed inclusion we consider an element $x\in \mathbb I_{g,\ord}^r[\frac{1}{p}]\cap (\mathbb I_{g,\ord}^r)_{(p)}$. Let $i\in\mathbb N\cup\{0\}$ be the smallest integer such that $y:=p^ix\in \mathbb I_{g,\ord}^r$. It suffices to show that the assumption that $i>0$ leads to a contradiction. As $x\in (\mathbb I_{g,\ord}^r)_{(p)}$, there exists $a\in\mathbb I_{g,\ord}^r\setminus p\mathbb I_{g,\ord}^r$ such that $z:=sx\in\mathbb I_{g,\ord}^r$. From the identity $p^iz=sy$ between elements of $\mathbb I_{g,\ord}^r$, as $i>0$ we get that $p$ divides $sy$. So, as $p$ is a prime element of $\mathbb I_{g,\ord}^r$ which does not divide $s$, we get that $p$ divides $y$ and therefore $p^{i-1}x=\frac{y}{p}\in\mathbb I_{g,\ord}^r$ which contradicts the minimality property of $i>0$. 

As $\mathbb I_{g,\ord}^r[\frac{1}{p}]\simeq \mathbb I_{g,\ord}^r\otimes_R K$ is a unique factorization domain, it is a Krull ring. Moreover, $(\mathbb I_{g,\ord}^r)_{(p)}$ is a discrete valuation ring (see Corollary \ref{C4} (b)). From the last two sentences and the prior paragraph we get that $\mathbb I_{g,\ord}^r$ is a Krull ring (see \cite{Sa}, Cor. (a) after Prop. 4.1). Based on this and \cite{Sa}, Cor. after Thm. 6.3 we get that $\mathbb I_{g,\ord}^r$ is a unique factorization domain. Thus part (a) holds.

Part (b) follows from Lemma \ref{L43} (a) of Subsection \ref{S4202}, part (c) follows from Corollary \ref{C23} (c) of Subsection \ref{S4211} and
part (d) follows from Proposition \ref{P14} of Subsection \ref{S4204}, these references being applied over an algebraic closure of $K$. 

Part (e) follows directly from the proof of Theorem \ref{T16} (a).

To check part (f), 
let $f\in \mathbb I^r_{g,\ord}(-2)\otimes_R K$ be a non-zero element.
By the unique factorization property in part (a), it is enough to check that 
$f$ is irreducible in  $\mathbb I_{g,\ord}^r\otimes_R K$. But if $f$ is reducible in $\mathbb I_{g,\ord}^r\otimes_R  K$, it is easy to see that we have a product decomposition $f=f'f''$ with $f'\in \mathbb I_{g,\ord}^r(w')\otimes_R  K$ and $f''\in \mathbb I_{g,\ord}^r(w'')\otimes_R K$ for some $w',w''\in W$ and $f',f''$ non-units of $\mathbb I^r_{g,\ord}\otimes_R  K$.
 Hence $\deg(w')+\deg(w'')=-2$. As $\deg(w')$ and $\deg(w'')$ are $\leq 0$ (see Corollary \ref{C7})
 it follows by part (d) that either
 $\deg(w')$ or $\deg(w'')$
 is $0$. But for $\deg(w)=0$ the $K$-vector space $\mathbb I_{g,\ord}(w)\otimes_R  K$ is one dimensional with basis an invertible element of $\mathbb I_{g,\ord}\otimes_R  K$, see part (e). So either $f'$ or $f''$ is invertible in $\mathbb I^r_{g,\ord}\otimes_R  K$, a contradiction.
 
 Finally part (g) follows form the fact that $\mathbb I_{g,\ord}(0)\otimes_R K=K$.
\endproof

\begin{rem}\label{R23}
{\bf (a)} Now we can explicitly compute the Serre--Tate expansions of all elements of the $R$-module
$\mathbb I^r_{g,\ord}(w)$, $w\in W$. Indeed, by the proof of Theorem
\ref{T16} (a),
 each element in $\mathbb I^r_{g,\ord}(w)$ is a $K$-multiple of an element of the form 
\begin{equation}\label{EQ107}
f=F^{\diamondsuit}\cdot \det(f^{\partial})^{w'}\end{equation}
for some $F\in \mathbb H^{r-1}_g(s)$, where $s:=-\deg(w)/2$ and $w':=\frac{w+2s}{\phi-1}\in W$. 
On the other hand it follows by Theorems \ref{T13} (a) and \ref{T14} (a) that the Serre--Tate expansion of $f$ above is given by the formula:
$$\mathcal E_{A_0,\theta_0,e}(f)=F(\Psi,\Psi^{\phi},\ldots,\Psi^{\phi^{r-1}}).$$

\smallskip
{\bf (b)} For $w\in W$, each element of $\mathbb I_{g,\ord}(w)$ vanishes at every
ordinary $\CM_R$-point of $\mathcal X_{\ord}$: it has the form (\ref{EQ107}) and $F^{\diamondsuit}$ vanishes at every ordinary $\CM_R$-point by Proposition \ref{P7} (a) and Remark \ref{R21}. In particular, for all linear systems $\f$ of weight $w$ in
$\mathbb I_{g,\ord}$
the set of $\f$-unstable points $\mathcal X_{\ord}(R)^{\textup{u}}_{\f}$ in $\mathcal X_{\ord}(R)$ contains the set of ordinary $\CM_R$-points. More generally, we can describe the set of $\f$-unstable points $\mathcal X_{\ord}(R)^{\textup{u}}_{\f}$
 in terms of the unstable locus of $\pmb{\SL}_g$ acting on spaces of multiple quadratic forms as follows. The $K$-algebra $\mathbb H^{r-1}_{g,\tot}\otimes_R K$ is known (by classical invariant theory) to be finitely generated. Let $s\in \frac{1}{2}\mathbb N$ be such that if
 $F_1,\ldots,F_{D_0(g,r-1,s)}$ is a $K$-basis of the homogeneous component $\mathbb H^{r-1}_g(s)_K$, then the radical ideal generated by $F_1,\ldots,F_{D_0(g,r-1,s)}$ in $\mathbb H^{r-1}_g\otimes_R K$ contains a set of generators of the $K$-algebra $\mathbb H^{r-1}_g\otimes_R K$. We can assume $F_1,\ldots,F_{D_0(g,r-1,s)}\in \mathbb H^r_g(s)$, i.e., the coefficients of the $F_i$s are in $R$.
 With the notation of Subsubsection \ref{S234}, let $\pmb{\V}_{g,R}^r:=\Spec(R[T,T',\ldots,T^{(r-1)}])$ and we consider the closed subscheme $\Sigma^r_{\pmb{\SL}_g}$ of
 $\pmb{\V}_{g,R}^r$ defined by $F_1,\ldots,F_{D_0(g,r-1,s)}$; for an algebraic closure $\mathbb K$ of $K$, the $\mathbb K$-valued points of $\Sigma^r_{\pmb{\SL}_g}$, are precisely the unstable points
 for the action of $\pmb{\SL}_g(\mathbb K)$ on $\pmb{\V}_{g,R}^r(\mathbb K)$ in the sense of \cite{MFK}, App. to Ch. I, Sect. B, Def. 
 By the proof of Theorem \ref{T16} (a), a basis of the $K$-linear space $\mathbb I^r_{g,\ord}(-2s)$
 is given by the forms
 $$F_1^{\diamondsuit},\ldots,F_{D_0(g,r-1,s)}^{\diamondsuit}.$$
 Let $X\subset \mathcal X_{\ord}$ be an affine open subscheme with $\overline{X}$ connected and such that the universal polarized abelian scheme $(A_X,\theta_X)$ has a column basis $\omega_X$ of $1$-forms on $A_X$ and consider the map
 $$f^{\langle 1,\ldots, r\rangle}:X(R)\rightarrow \pmb{\V}_{g,R}^r(R)$$
 which maps $P\in X(R)$ to 
 $$(f^{\langle 1 \rangle}(P),\ldots, f^{\langle r \rangle}(P)):=(f^{\langle 1 \rangle}(A_P,\theta_P,\omega_P,R),\ldots, f^{\langle r \rangle}(A_P,\theta_P,\omega_P,R))$$
 where $A_P,\theta_P,\omega_P$ over $\Spec(R)$ are as usual induced by $A_X,\theta_X,\omega_X$
 via base change by $P:\Spec(R)\rightarrow X$. Then, one can check, by directly using the definitions, that the set 
 $\mathcal X_{\ord}(R)^{\textup{u}}_{\f}$ of $\f$-unstable points for $\f$ defined by an arbitrary $R$-basis of 
 $\mathbb I_{g,\ord}^r(-2s)$, such as the one above, is given by the formula
 $$\mathcal X_{\ord}(R)^{\textup{u}}_{\f}=(f^{\langle 1, \ldots, r \rangle})^{-1}(\Sigma^r_{\pmb{\SL}_g}(R)).$$
 
 \smallskip
 
 {\bf (c)} In the proof of Theorem \ref{T16} (e) and (f) we considered elements of $\mathbb I_{g,\ord}$, $b_0=\det(f^{\partial})$ of weight $\phi-1$ and $b_i= (F_i)^{\diamondsuit}$ for $i\in\{1,\ldots,\frac{g(g+1)}{2}\}$ of weight $-2s$ for some $s\in\mathbb N$, and we checked that the family $(b_i)_{i\in\{0,\ldots,\frac{g(g+1)}{2}\}}$ of Equation (\ref{EQ105}) is $\phi$-algebraically independent over $R$. Based on this, the argument in the proof of Remark \ref{R6} (a) implies that
 the family $$(\tilde{b}_i)_{i\in \{0,\ldots,\frac{g(g+1)}{2}\}}$$
 of elements of the same weight $-sg(g+1)$,
 $$\tilde{b}_0:=b_1b_2 \cdots b_{\frac{g(g+1)}{2}},\ \ \tilde{b}_1:=b_0^{2s}b_1^{\phi}b_2\cdots b_{\frac{g(g+1)}{2}}, \ldots,
 \tilde{b}_{\frac{g(g+1)}{2}}:=b_0^{2s}b_1b_2\cdots b_{\frac{g(g+1)}{2}}^{\phi},$$
 is $\phi$-algebraically independent.
 \end{rem}

The following is a direct consequence of Lemma \ref{L11} and Theorem \ref{T16} (a):

\begin{cor}\label{C16}
We have a $K[\Theta_0]$-algebra isomorphism
$$
K[\Theta_0,z_0,z_0^{-1}] \simeq \mathbb I^1_{g,\ord}\otimes_R K .
$$
that maps $z_0$ to $\det(f^{\partial}_{g,\crys})$ and that extends to a $ K[\Theta_0,\ldots,\Theta_g]$-algebra isomorphism
 $$
 K[\Theta_0,\ldots,\Theta_g,z_0,z_1,\frac{1}{z_0z_1}]\simeq \mathbb I^2_{g,\ord}\otimes_R K,
$$
which maps $z_1$ to $\det(f^{\partial}_{g,\crys})^{\phi}$.
\end{cor}

The structure of $\mathbb I^r_{g,\ord}\otimes_R K$ is more complicated as $r\geq 3$ increases. A `good picture' of what happens in this case can be obtained by combining Theorem \ref{T16}
(a) with the results of Subsections \ref{S42} and \ref{S43}; we leave this to the reader.

In view of Theorem \ref{T16} (a) (and the work in Subsection \ref{S42}) the structure of the ordinary rings $\mathbb I^r_{g,\ord}$ is more or less well understood. By contrast, 
the structure of the non-ordinary rings $\mathbb I^r_g$ remains somewhat mysterious. We will address these rings in the next two subsections. Before that, we record the following generalization of Lemma \ref{L9} (a):

\begin{cor}\label{C17}
We have 
$\mathbb I^1_g(\phi-1)=\mathbb I^1_g(1-\phi)=0$ and thus $\mathbb I_g\neq \mathbb I_{g,\ord}$.
\end{cor}

\noindent
{\it Proof.}
We will only show that the assumption $\mathbb I^1_g(\phi-1)\neq 0$ leads to contradiction; the case of $\mathbb I^1_g(1-\phi)$ is similar (and actually easier) and is left to the reader. 
As the $R$-module $\mathbb I^1_{g,\ord}(\phi-1)/\mathbb I^1_g(\phi-1)$ is torsion free (see Proposition \ref{P3} (b)), from Corollary \ref{C15} (e) and our assumption we get that $\det(f^{\partial}_{g,\crys})\in \mathbb I^1_g(\phi-1)=\mathbb I^1_{g,\ord}(\phi-1)$.
Let $Y\subset Y_1(N)$ be an affine open subset of  the affine modular curve of level $\Gamma_1(N)$ over $\Spec(R)$ with $N\geq 4$ not divisible by $p$, such that the reduction modulo $p$ of $Y$ is not contained in the ordinary locus of the reduction modulo $p$ of $Y_1(N)$. Let $E_Y$ be the universal elliptic curve over $Y$, let $\theta_{E_Y}$ be the Abel--Jacobi map (polarization) of $E_Y$ and let $\omega_{E_Y}$ be a column basis of $1$-forms on $E_Y$. Let $(E_{1,R},\theta_1),\ldots,(E_{g-1,R},\theta_{g-1})$ be ordinary elliptic curves over $\Spec(R)$ endowed with their Abel--Jacobi maps (polarizations) such that we have $f^{\partial}_{1,\crys}(E_{i,R},\theta_i,\omega_i, R)\in R^{\times}$ for $i\in\{1,\ldots,g-1\}$, where each $\omega_i$ is a column basis of $1$-forms on $E_{i,R}$. Let $(E_{i,Y},\theta_{i,Y},\omega_{i,Y}):=(E_{i,R},\theta_i,\omega_i)\times_{\Spec(R)}Y$. Let $(A_Y,\theta_{A_Y},\omega_{A_Y})\in\pmb{\M}_g(\mathcal O(Y))$ be defined by
$$(A_Y,\theta_{A_Y}):=(E_Y,\theta_{E_Y})\times_Y (E_{1,Y},\theta_{1,Y})\times_Y\cdots \times_Y (E_{g-1,Y},\theta_{g-1,Y})$$
and by the column basis of $1$-forms on $A_Y$ whose $11$ entry is the pullback of $\omega_Y$ and whose $i1$ entry for $i\in\{2,\ldots,g\}$ is the pullback of $\omega_{i-1,Y}$.
We have that
$$\det(f^{\partial}_{g,\crys})(A_Y,\theta_{A_Y},\omega_{A_Y},\mathcal O^{\infty}(Y))\in \mathcal O(J^1(Y)).$$
On the other hand the above element equals
$$f^{\partial}_{1,\crys}(E_Y,\theta_{E_Y},\omega_{E_Y},\mathcal O^{\infty}(Y))\cdot
\prod_{i=1}^{g-1} f^{\partial}_{1,\crys}(E_{i,R},\theta_i,\omega_i,R).$$
So we have
$$f_Y:=f^{\partial}_{1,\crys}(E_Y,\theta_{E_Y},\omega_{E_Y},\mathcal O^{\infty}(Y))\in 
\mathcal O(J^1(Y)).$$
For a prime $\ell$ that does not divide $pN$, let  $Y_1(N,\ell)$ be the curve  over $\Spec(R)$ that parametrizes isogenies of degree $\ell$ between elliptic curves with level $\Gamma_1(N)$ structure that respect the level $\Gamma_1(N)$ structures.
 We consider the two canonical projections $\pi_1,\pi_2:Y_1(N,\ell)\rightarrow Y_1(N)$ and the intersection $Y_{\ell}:=\pi_1^{-1}(Y)\cap \pi_2^{-1}(Y)$. Let $Y_{\ord}\subset Y$ and $Y_{\ell,\ord}\subset Y_{\ell}$ be the open subschemes which over $K$ are isomorphisms and whose fibers over $k$ are ordinary loci.  
 For $i\in \{1,2\}$ let $\pi_i^*E_Y/Y_{\ell}$ be the pullback of $E_Y/Y$ via $\pi^*_i:Y_{\ell}\rightarrow Y$ and consider the $1$-forms $\pi_i^*\omega_{E_Y}$ on $\pi_i^*E_Y$. 
 Let
 $$\lambda:=\frac{u^*(\pi_1^*\omega_{E_Y})}{\pi_2^*\omega_{E_Y}}\in \mathcal O(Y_{\ell})^{\times}$$
 with $u:\pi_2^*E_Y\rightarrow \pi_1^*E_Y$  the `universal isogeny' over $Y_{\ell}$; therefore we have $\lambda^{\phi-1}\in \mathcal O(J^1(Y_{\ell}))^{\times}$. By the Hecke covariance of $f^{\partial}_{1,\crys}$ 
  we have the following equality in 
 $\mathcal O(J^1(Y_{\ell,\ord}))$:
\begin{equation}
\label{relationbetweenpullbacks}
\pi_2^*f_Y=\lambda^{\phi-1}\cdot \pi_1^*f_Y.\end{equation}
 By the injectivity of the $R$-algebra homomorphism $\mathcal O(J^1(Y_{\ell}))\rightarrow \mathcal O(J^1(Y_{\ell,\ord}))$ the equality (\ref{relationbetweenpullbacks}) holds in $\mathcal O(J^1(Y_{\ell}))$. The existence of such a function $f_Y$ in $\mathcal O(J^1(Y))$ contradicts \cite{Bu05}, Thm. 8.83 3).\endproof

\begin{rem}\label{R24}
It may be reasonable to expect that $\mathbb I^r_g(w)=0$ for all $w\neq 0$ with $\deg(w)=0$.
By Corollary \ref{C15} (e) this is equivalent to the statement that 
\begin{equation}\label{EQ108}
\det(f^{\partial}_{g,\crys})^{\frac{w}{\phi-1}}\not\in \mathbb I^r_g(w)\end{equation}
for all such $w$.
Corollary \ref{C17} gives a partial result in this direction. 
One can show that Statement (\ref{EQ108}) holds for $g=1$ provided $w=\sum_{i=0}^r a_i\phi^i\in W$ with 
$a_0,\ldots,a_r$ integers such that $\sum_{i=0}^r a_i=0$ and 
$\sum_{i=0}^ra_ip^i<0$. This follows from the fact that the reduction modulo $p$ of $f^{\partial}_{1,\crys}$ is the Hasse invariant, see \cite{Bu05}, Prop. 8.57.
On the other hand note that we do not know, for instance, if Statement (\ref{EQ108}) holds even for $g=1$ and $w=n(\phi-1)$, with $n\geq 2$ an integer. This is so as, for $n\geq 2$, from the fact that $f^{\partial}_{1,\crys}\not\in \mathbb I^1_1(\phi-1)$
it does not follow a priori that $(f^{\partial}_{1,\crys})^n\not\in \mathbb I^1_1(n(\phi-1))$. Concretely, the following phenomenon occurs. 
Let $S$ be a smooth $R$-algebra, let $h\in S\backslash pS$, and let $f\in \widehat{S_h}$ be such that $f^n\in \widehat{S}$ for some $n\geq 2$; then in general it does not follow that $f \in \widehat{S}$. We include a simple example. Let $S:=R[x]$ with $x$ an indeterminate, 
 let $n\geq 2$ be an integer that does not divide $p$, and let $h:=x$. The equation $y^n=x^{n-1}(x-p)$ in $y$ has a solution 
 $$z:=x-p/n+\sum_{m=2}^{\infty} p^mz_m\in\widehat{R[x,x^{-1}]},$$ with each $z_m\in R[x,x^{-1}]$.
 We have $z^n=x^{n-1}(x-p)\in \widehat{R[x]}$ but $z\not\in \widehat{R[x]}$. \end{rem}

\subsection{Comparison maps: the non-ordinary case}\label{S37}

In this subsection we define some further new modules and algebras and we introduce some further `comparison maps' $\triangle$, $\square$, between them. We will use these maps, in conjunction with our Serre--Tate expansion maps $\mathcal E$ and with the invariant theory of multiple endomorphisms (see Subsubsections \ref{S411}, \ref{S415} and \ref{S416}), to derive our basic estimates from below for the dimensions of $K$-vector spaces of Hecke covariant (non-ordinary) Siegel $\delta$-modular forms.

\subsubsection{The maps $\square$}\label{S3701}

Let $n\in \mathbb N$ and let $X_1=(X_{1,jl})_{1\leq j,l\leq g},\ldots,X_n=(X_{n,jl})_{1\leq j,l\leq g}$ be universal $g\times g$ matrices whose entries are indeterminates.
Let $B$ be an integral domain which is a $\mathbb Z_p$-algebra. Consider the left action of ${\pmb{\GL}}_n(B)$ on the polynomial $B$-algebra $B[X_1,\ldots,X_n]:=B[X_{i,jl}|1\leq i\leq n, 1\leq j,l\leq g]$ defined by the rule:
$$(M,F)\mapsto F(M^{-1}X_1M,\ldots,M^{-1}X_nM).$$
We denote by $\mathbb H^n_{g,\con,B}$ the $B$-algebra of all $\pmb{\GL}_g(B)$-invariant polynomials in $B[X_1,\ldots,X_n]$. 
Clearly 
$$\mathbb H^n_{g,\con,D^{-1}B}=\mathbb H^n_{g,\con}\otimes_B D^{-1}B$$
for each multiplicative set $D\subset B$. Also 
$$\mathbb H^n_{g,\con,B'}=\mathbb H^n_{g,\con}\otimes_B B'$$ for each extension $B\subset B'$ of infinite fields. Moreover, we have an identity $\mathbb H^n_{g,\con,R}=\mathbb H^n_{g,\con}$.
For $B$ an algebraically closed field and $n\geq 2$ we have the formula (see Theorem \ref{T26} of Subsubsection \ref{S411}):
$$\trdeg_R(\mathbb H^n_{g,\con,B})=D_{\con}(g,n):=(n-1)g^2+1.$$
Also recall that the `first main theorem' in this setting 
 (see Theorem \ref{T28} of Subsection \ref{S416})
implies that if $B$ an infinite field, then the $B$-algebra $\mathbb H^n_{g,\con,B}$ is generated by finitely many polynomials of the form
$$\tau_{j;i_1,i_2,\cdots, i_N}:=c_j(X_{i_1}X_{i_2}\cdots X_{i_N})$$
where $j\in \{1,\ldots,g\}$, $N\in\mathbb N$, and $i_1,\ldots,i_N\in\{1,\ldots,n\}$, and where
for a matrix $M\in\pmb{\Mat}_g(C)$ with $C$ a ring, we let $c_0(M),\ldots,c_g(M)\in C[t]$ be such that
 $\det(t\cdot 1_g-M)=\sum_{j=0}^g (-1)^jc_j(M)t^{g-j}\in C[t]$
 is its characteristic polynomial; henceforth $c_0(M)=1$, $c_1(M)$ is the trace $\Trace(M)$ of $M$, and $c_g(M)=\det(M)$.
There exists a `second main theorem' describing relations between these generators (see \cite{P}, Introd.) but we will not need these relations in what follows. 

Note that the natural $k$-algebra homomorphisms
$$\mathbb H^n_{g,\con,R}\otimes_R k\rightarrow \mathbb H^n_{g,\con,k}$$
are isomorphisms.
Indeed they are surjective by the first main theorem and they are injective because
for $F\in R[X_1,\ldots,X_n]$, if $pF$ is $\pmb{\GL}_g(R)$-invariant, then $F$ itself is $\pmb{\GL}_g(R)$-invariant. By Nakayama's Lemma applied to the homogeneous components of $\mathbb H^n_{g,\con,R}$ it follows that the latter $R$-algebra is generated
by finitely many elements of the form $\tau_{j;i_1,i_2,\cdots, i_N}$; we will refer to this statement as the `first main theorem over $R$.'

More generally, let $x_1,\ldots,x_{n+2},\xi_1,\ldots,\xi_{n+2}$ be indeterminates and view \begin{equation}\label{EQ109}
R[X_1,\ldots,X_n,x_1,\xi_1,\ldots,x_{n+2},\xi_{n+2}]:=R[X_1,\ldots,X_n][x_1,\xi_1,\ldots,x_{n+2},\xi_{n+2}]
\end{equation}
 as a $W$-graded $R$-algebra
by letting 
$$\deg(x_i):=-\phi^{i-1}-\phi^i\;\;\textup{and}\;\;\deg(\xi_i):=-1-\phi^i\;\;\forall i\in\{1,\ldots,n+2\}$$
and
$$\deg(X_{i,jl}):=-1-\phi^i-\phi^{i+1}-\phi^{i+2}\;\;\forall i\in\{1,\ldots,n\}\;\;\textup{and}\;\;\forall j,l\in\{1,\ldots,g\}.$$ 
Then by definition
$$\mathbb H^n_{g,\con}[x_1,\ldots,x_{n+2},\xi_1,\ldots,\xi_{n+2}]$$
 is a $W$-graded $R$-subalgebra of (\ref{EQ109}). 
We will prove the following: 

\begin{thm}\label{T17}
For all integers $n\geq 2$ there exists a natural homomorphism
$$\square=\square_n:\mathbb H^n_{g,\con}[x_1,\ldots,x_{n+2},\xi_1,\ldots,\xi_{n+2}]\rightarrow \mathbb I^{n+2}_g,$$
of $W$-graded $R$-algebras whose image satisfies
$$\trdeg_R(\textup{Im}(\square))\geq (n-1)\frac{g(g+1)}{2}+n+g+1.$$
\end{thm}

The homomorphisms $\square_n$ in the above theorem will be compatible as $n\geq 2$ varies and will be entirely explicit, as we will see shortly. 

\begin{cor}\label{C18}
The following four properties hold:
 
 \medskip
 {\bf (a)} For all integers $r\geq 4$ we have inequalities
$$
(r-3)\frac{g(g+1)}{2}+r+g-1
\leq \trdeg_R (\mathbb I^r_g)\leq 
r\cdot \frac{g(g+1)}{2}-g^2+1+r.$$

\smallskip
{\bf (b)} We have an inequality
$$\trdeg_{\mathbb I_g} (\mathbb I_{g,\ord})\leq \frac{g(g+1)}{2}+2.$$

\smallskip
{\bf (c)} We have equalities
$$
\asytrdeg(\mathbb I_g)= \phidim(\mathbb I_g)=\frac{g(g+1)}{2}+1.
$$

\smallskip
{\bf (d)} We have an equality
$$
\phidim(\mathbb I_{g,((p))})=\frac{g(g+1)}{2}.
$$
 \end{cor}
 
 \noindent
 {\it Proof.}
Part (a) follows from Theorems \ref{T17} and \ref{T16} (b).
 To check part (b) note that the difference between the upper and the lower bounds 
 for $\trdeg_R (\mathbb I^r_g)$ in part (a)
 equals $\frac{g(g+1)}{2}+2$. From this and Theorem \ref{T16} (b) we get that
$$\trdeg_{\mathbb I^r_g} (\mathbb I^r_{g,\ord})=\trdeg_{R} (\mathbb I^r_{g,\ord})-\trdeg_{R} (\mathbb I^r_{g})\leq \frac{g(g+1)}{2}+2$$
for all $r\geq 4$. Thus part (b) follows from the part of Remark \ref{R3} after
 Equation (\ref{EQ009}). The equality 
 $\asytrdeg(\mathbb I_g)=\frac{g(g+1)}{2}+1$ in part 
 (c) follows from part (a).
 The equality 
 $\phidim(\mathbb I_g)=\frac{g(g+1)}{2}+1$ in part (c) follows from the implication in Equation (\ref{EQ009}) plus Theorem \ref{T16} (c); alternatively one can use part (b) plus Theorem \ref{T16} (c) and the additivity 
 of the $\phi$-transcendence degree (see Equation (\ref{EQ007})). 
 Part (d) follows from part (c) and Remark \ref{R6} (a).
 \endproof
 
\begin{rem}\label{R25}
{\bf (a)} By an argument similar to the one in the proof of Corollary \ref{C14} we have that
$$\trdeg_k(\mathbb I^r_g\otimes_R k)\leq \trdeg_R(\mathbb I^r_g).$$


{\bf (b)}
By Remark \ref{R23} it follows that for all linear systems $f$ weight $w\in W$ in
$\mathbb I_{g,\ord}$
the set of $f$-unstable points
$\mathcal X(R)^{\textup{u}}_f$ in $\mathcal X(R)$ contains the set of ordinary $\CM_R$-points.
\end{rem}

\subsubsection{Cyclic products and the maps $\triangle$}\label{S3702}

In order to proceed with our proof we need to recall (and slightly extend) the {\it cyclic product construction} from \cite{BB}, Subsect. 4.2.
 Recall we already introduced our {\it basic} forms 
$f^r=f^r_{g,\crys}\in\mathbb I^r_{gg}(\phi^r,1)$, $r\in \mathbb N$. 
For integers $a>b\geq 0$ we set 
\begin{equation}
\label{deffab}
f^{ab}:=\phi^b f^{a-b}\in\mathbb I^a_{gg}(\phi^a,\phi^b)\;\;\textup{and}\;\;f^{ba}:=(f^{ab})^{\t}\in \mathbb I^a_{gg}(\phi^b,\phi^a).\end{equation}
Let $\cyc^r$ be the set of all $2s$-tuples $(a_1,a_2,\ldots,a_{2s-1},a_{2s})$ of arbitrary even length $2s\in 2\mathbb N$, such that the following three properties hold:

\medskip
{\bf (i)} For all $i\in\{1,\ldots,2s\}$, we have $a_i\in\{0,1,\ldots,r\}$.

\smallskip
{\bf (ii)} For all $i\in\{1,\ldots,2s-1\}$, we have $a_i\neq a_{i+1}$.

\smallskip
{\bf (iii)} We have $a_{2s}\neq a_1$.

\medskip

\noindent 
For each $(a_1,a_2,\ldots,a_{2s-1},a_{2s})\in\cyc^r$ we define the 
{\it matrix cyclic product}
$$
F^{a_1,a_2,\ldots, a_{2s-1},a_{2s}}\in M^r_{gg}=\pmb{\Mat}_g(M^r_g)$$
by the formula
\begin{equation}
\label{defFa1}F^{a_1,a_2,\ldots, a_{2s-1},a_{2s}}:=f^{a_1a_2}(f^{a_3a_2})^*f^{a_3a_4}(f^{a_5a_4})^*\cdots f^{a_{2s-1}a_{2s}}(f^{a_1a_{2s}})^*.\end{equation}
Note that for each
$S \in\Ob(\Prol)$ and every 
$(x_1,x_2,u)\in \pmb{\M}'_g(S^0)$, where
$x_1=(A_1,\theta_1,\omega_1)$, $x_2=(A_2,\theta_2,\omega_2)$, $uu^{\t}=d\in\mathbb N\setminus p\mathbb N$,
the following equality holds (see Equation (\ref{EQ029.9})):
  \begin{equation}\label{EQ110}
 F^{a_1,\ldots,a_{2s}}(x_2,S)=d^{-gs}(\det([u]))^{-w} [u]^{\phi^{a_1}} F^{a_1,\ldots,a_{2s}}(x_1,S) ([u]^{\phi^{a_1}})^{-1},
 \end{equation}
 where $w:=-\sum_{i=1}^{2s}\phi^{a_i}\in W$ and $[u]\in\pmb{\GL}_g(R)$ is defined by the identity $u^*\omega_2=[u]\omega_1$.
 Then, for $j\in \{1,\ldots,g\}$ we define the (scalar)
{\it cyclic products}
$$f_{j;a_1,a_2,\ldots, a_{2s-1},a_{2s}}
:=c_j(F^{a_1,a_2, \ldots ,a_{2s-1},a_{2s}})\in M^r_g.$$
For $j=1$ we abbreviate
\begin{equation}
\label{deffa1}
f_{a_1,a_2,\ldots, a_{2s-1},a_{2s}}:=f_{1;a_1,a_2,\ldots, a_{2s-1},a_{2s}}
:=\Trace(F^{a_1,a_2, \ldots ,a_{2s-1},a_{2s}})\in M^r_g.\end{equation}
Note that, for $a\in\mathbb N$ we have 
$$f_{0,a}=\Trace(F^{0,a})=\Trace(f^{0a}(f^{0a})^*)=\Trace(\det(f^{0a})\cdot 1_g)=g\cdot \det(f^a).$$
Clearly the $K$-algebra generated by the members of the family 
$$(f_{j;a_1,a_2,\ldots, a_{2s-1},a_{2s}})_{(j,(a_1,a_2,\ldots, a_{2s-1},a_{2s}))\in \{1,\ldots g\}\times \textup{Cyc}^r}$$ is the same as the $K$-algebra generated by the members of the  family 
$$(f_{a_1,a_2,\ldots, a_{2s-1},a_{2s}})_{(a_1,a_2,\ldots, a_{2s-1},a_{2s})\in \textup{Cyc}^r};$$
however, the $R$-algebras generated by the members of these two families are distinct if $p\le g$.
By Equation (\ref{EQ110}) we have
\begin{equation}\label{finIr}
f_{j;a_1,a_2,\ldots, a_{2s-1},a_{2s}}\in\mathbb 
I^r_g(-j(\sum_{i=1}^{2s}\phi^{a_i})).\end{equation}

On the other hand, for each integer $r\geq 1$ we consider the $R$-subalgebra
$$R[Q^{(1)},\ldots,Q^{(r)}]\subset S^r_{\for}$$
generated by the entries of the matrices 
\begin{equation}
\label{defQ(s)}
Q^{(s)}:=\Psi^{\phi^{s-1}}\in \pmb{\Mat}_g(S^r_{\for}),\ \ s\in \mathbb N
\end{equation} (see Equations (\ref{EQ034.9}) and (\ref{EQ034.99})). The 
entries of these symmetric matrices that are above and on the diagonal are (distinct and) algebraically independent over $R$ (because, by definition, the same is true for the matrices $\Psi,\delta \Psi,\ldots, \delta^{r-1}\Psi$ and, on the other hand, the $K$-algebra generated by the entries of $Q^{(1)},\ldots,Q^{(r)}$ coincides with the $K$-algebra generated by the entries of $\Psi,\delta \Psi,\ldots, \delta^{r-1}\Psi$);
 we will view these entries of 
 $Q^{(1)},\ldots,Q^{(r)}$
  as indeterminates over $R$ and also 
 as indeterminates in the ring
 $R[Q^{(1)},\ldots,Q^{(r)}]/(p)$. 
We consider the reduction modulo $p$ 
homomorphism
$$\triangle:R[Q^{(1)},\ldots,Q^{(r)}]\rightarrow k[Q^{(1)},\ldots,Q^{(r)}]:=R[Q^{(1)},\ldots,Q^{(r)}]/(p),\ \ G\mapsto G^{\triangle}.$$
We note that the natural $k$-algebra homomorphism
$$k[Q^{(1)},\ldots,Q^{(r)}]\rightarrow S^r_{\for}\otimes_R k=k[[T]][T',\ldots,T^{(r)}]$$
is not injective. 

For each $(a_1,a_2,\ldots,a_{2s-1},a_{2s})\in\cyc^r$
we define a matrix with entries in the ring $R[Q^{(1)},\ldots,Q^{(r)}]$ by the formula
$$Y_{a_1,a_2,\ldots, a_{2s-1},a_{2s}}:=Q^{(m_1)}(Q^{(m_2)})^*Q^{(m_3)}(Q^{(m_4)})^*\ldots Q^{(m_{2s-1})}(Q^{(m_{2s})})^*,$$
where
$$m_1:=\max\{a_1,a_2\},\ \ m_2:=\max\{a_2,a_3\},\ldots, m_{2s}:=\max\{a_{2s},a_1\}.$$
For each such matrix and $j\in \{1,\ldots,g\}$ we define a polynomial
$$y_{j;a_1,a_2,\ldots, a_{2s-1},a_{2s}}:=c_j(Y_{a_1,a_2,\ldots, a_{2s-1},a_{2s}})\in R[Q^{(1)},\ldots,Q^{(r)}]$$
and we set
$$y_{a_1,a_2,\ldots, a_{2s-1},a_{2s}}:=y_{1;a_1,a_2,\ldots, a_{2s-1},a_{2s}}.$$

\begin{lemma}\label{L20}
 For each $(a_1,a_2,\ldots,a_{2s-1},a_{2s})\in\cyc^r$ we have
 $$(\mathcal E_{A_0,\theta_0,e}(f_{j;a_1,a_2,\ldots, a_{2s-1},a_{2s}}))^{\triangle}= y^{\triangle}_{j;a_1,a_2,\ldots, a_{2s-1},a_{2s}}\neq 0.$$
\end{lemma}

\noindent
 {\it Proof.} This follows directly from Theorem \ref{T13} (a).\endproof

\medskip

With $\Lambda_0,\ldots,\Lambda_r$ as indeterminates, let 
the polynomial algebras 
$$R[Q^{(1)},\ldots,Q^{(r)},\Lambda_0,\ldots,\Lambda_r]:= R[Q^{(1)},\ldots,Q^{(r)}][\Lambda_0,\ldots,\Lambda_r],$$ 
$$k[Q^{(1)},\ldots,Q^{(r)},\Lambda_0,\ldots,\Lambda_r]:= k[Q^{(1)},\ldots,Q^{(r)}][\Lambda_0,\ldots,\Lambda_r].$$ 
We extend the homomorphism $\triangle$ above to a reduction modulo $p$ homomorphism between these polynomial
algebras
$$\triangle:R[Q^{(1)},\ldots,Q^{(r)},\Lambda_0,\ldots,\Lambda_r]\rightarrow k[Q^{(1)},\ldots,Q^{(r)},\Lambda_0,\ldots,\Lambda_r],$$
 by setting $\Lambda_j^{\triangle}:=\Lambda_j$. 

\begin{rem}\label{R26}
For $F_1,\ldots,F_m\in R[Q^{(1)},\ldots,Q^{(r)},\Lambda_0,\ldots,\Lambda_r]$, if
$F_1^{\triangle},\ldots, F_m^{\triangle}$ are algebraically independent in $k[Q^{(1)},\ldots,Q^{(r)},\Lambda_0,\ldots,\Lambda_r]$, then $F_1,\ldots,F_m$ are algebraically independent in $R[Q^{(1)},\ldots,Q^{(r)},\Lambda_0,\ldots,\Lambda_r]$.\end{rem}

We define
$$ z_{j;a_1,a_2,\ldots, a_{2s-1},a_{2s}}:=y_{j;a_1,a_2,\ldots , a_{2s-1},a_{2s}} \Lambda_{a_1}^j\Lambda_{a_2}^j\cdots\Lambda_{a_{2s-1}}^j\Lambda_{a_{2s}}^j,$$ 
$$
z_{a_1,a_2,\ldots, a_{2s-1},a_{2s}}:=z_{1;a_1,a_2,\ldots, a_{2s-1},a_{2s}}.
$$
So for $\mathcal E=\mathcal E_{A_0,\theta_0,e}$ we have
\begin{equation}\label{EQ113}
(\mathcal E(f_{j;a_1,a_2,\ldots, a_{2s-1},a_{2s}})
\Lambda_{a_1}^j\cdots\Lambda_{a_{2s}}^j)^{\triangle}= z^{\triangle}_{j;a_1,a_2,\ldots, a_{2s-1},a_{2s}}.\end{equation}

\begin{lemma}\label{L21}
Consider the unique $R$-algebra homomorphism
$$\square:R[X_1,\ldots,X_n,x_1,\ldots,x_{n+2},\xi_1,\ldots,\xi_{n+2}]\rightarrow M^{n+2}_g, \ \ F\mapsto F^{\square},$$
defined by the rules
$$\begin{array}{rcll}
X_{i}^{\square}& := & F^{0,i,i+1,i+2}, & i\in \{1,\ldots,n\},\\
\ & \ & \ & \ \\
x_i^{\square} & := & f_{i-1,i}, & i\in \{1,\ldots,n+2\},\\
\ & \ & \ & \ \\
\xi_i^{\square}& :=&f_{0,i}, & i\in \{1,\ldots,n+2\}.\end{array}$$
Then $\square$ induces a $W$-graded $R$-algebra homomorphism 
$$\square:\mathbb H^n_{g,\con}[x_1,\ldots,x_{n+2},\xi_1,\ldots,\xi_{n+2}]\rightarrow \mathbb I^{n+2}_g.$$
\end{lemma}

\noindent
{\it Proof.}
By the `first main theorem over $R$' it is enough to check that for all $i_1,\ldots,i_N\in\{1,\ldots,g\}$ the element $\tau_{j;i_1,\ldots , i_N}^{\square}$ belongs to $\mathbb I^{n+2}_g(w)$, where $w\in\mathbb Z[\phi]$ is the sum of the degrees of $X_{i_1},\ldots,X_{i_N}$; but this follows from the fact that
$$\tau_{j;i_1,\ldots, i_N}^{\square}=f_{j;0,i_1,i_1+1,i_1+2,0,i_2,i_2+1,i_2+2,\ldots,0,i_N,i_N+1,i_N+2}.$$
\endproof

\noindent
{\it Proof of Theorem \ref{T17}.} 
Consider the polynomial $k$-algebra $$k[X_1,\ldots,X_n,x_1,\ldots,x_{n+2}]$$ in the indeterminates 
$x_1,\ldots,x_{n+2}$ and 
the entries of the matrices $X_1,\ldots,X_n$.
Consider the $k$-algebra homomorphism
\begin{equation}\label{EQ114}
k[X_1,\ldots,X_n,x_1,\ldots,x_{n+2}]\rightarrow k[Q^{(1)},\ldots,Q^{(n+2)}], \ \ F\mapsto \widetilde{F},\end{equation}
defined by
$$
\widetilde{X_i}:=Q^{(i)}(Q^{(i+1)})^*\in \pmb{\Mat}_g(S^{i+1}_{\for}),\ \ i\in\{1,\ldots,n\},$$
$$
\widetilde{x_i}:=\det(Q^{(i)}), \ \ i\in\{1,\ldots, n+2\}.$$
We let $\reallywidetilde{\mathbb H^n_{g,\con,k}}$ be the image of $\mathbb H^n_{g,\con,k}$ under the $k$-algebra homomorphism (\ref{EQ114}). For integers $n\geq 2$, by Theorem \ref{T27} of Subsubsection \ref{S411} applied with $(Q^{(1)},\ldots,Q^{(n+1)})\longleftrightarrow (Q_0,\ldots,Q_n$) we have
\begin{equation}\label{EQ115}
\dim(\reallywidetilde{\mathbb H^n_{g,\con,k}})= (n-1)\frac{g(g+1)}{2}+g.
\end{equation}
By Equation (\ref{EQ115})
there exists a set 
$$I\subset \bigcup_{N=1}^{\infty}(\{1,\ldots,g\}\times\{1,\ldots,n\}^N)$$ whose cardinality $|I|$ is given by
\begin{equation}\label{EQ116}
 |I|=(n-1)\frac{g(g+1)}{2}+g\end{equation}
such that   
 the family
\begin{equation}\label{EQ117}
(c_j(Q^{(i_1)}(Q^{(i_1+1)})^*\cdots Q^{(i_N)}(Q^{(i_N+1)})^*))_{(j,i_1,\ldots,i_N)\in I},\end{equation}
is algebraically independent over $k$ in the ring $k[Q^{(1)},\ldots,Q^{(n+1)}]$.
 But the members of the family (\ref{EQ117}) can be written as
\begin{equation}\label{EQ118}
\frac{c_j(Q^{(i_1)}(Q^{(i_1+1)})^*Q^{(i_1+2)}(Q^{(i_1+2)})^*\cdots Q^{(i_N)}(Q^{(i_N+1)})^*Q^{(i_N+2)}(Q^{(i_N+2)})^*)}{\det(Q^{(i_1+2)})^j\cdots \det(Q^{(i_N+1)})^j\det(Q^{(i_N+2)})^j}\end{equation}
which in their turn can be written as
$$\left(\frac{y_{j;0,i_1, i_1+1, i_1+2,0,i_1,i_1+1,i_1+2,\ldots, 0,i_N, i_N+1, i_N+2}}{(y_{i_1+1,i_1+2} \cdots y_{i_N+1,i_N+2})^j}\right)^{\triangle},\ \ (j, i_1,\ldots,i_N)\in I.$$
So if we consider the two families
\begin{equation}\label{EQ119}
(y^{\triangle}_{i-1,i})_{i\in\{1,\ldots,n+2\}},
\end{equation}
\begin{equation}\label{EQ120}
(y^{\triangle}_{j;0,i_1, i_1+1, i_1+2,\ldots,0,i_N, i_N+1, i_N+2})_{(j,i_1,\ldots,i_N)\in I},
\end{equation}
then the members of the family (\ref{EQ119})$\cup$(\ref{EQ120}) generate a field extension of $k$ of transcendence degree $\geq |I|$.
If we consider the family
\begin{equation}\label{EQ121}
(\Lambda_{i})_{i\in \{0,\ldots,n+2\}},
\end{equation} we get that 
the field extension of $k$ generated by 
 the members of the family
(\ref{EQ119})$\cup$(\ref{EQ120})$\cup$(\ref{EQ121})
has transcendence degree $\geq |I|+n+3$. 
Consider the four families
\begin{equation}\label{EQ122}
(z^{\triangle}_{j;0,i_1, i_1+1, i_1+2,\ldots,0,i_N, i_N+1, i_N+2})_{(j,i_1,\ldots,i_N)\in I},
\end{equation}
\begin{equation}\label{EQ123}
(z^{\triangle}_{i-1,i})_{i\in\{1,\ldots,n+2\}},
\end{equation}
\begin{equation}\label{EQ124}
(z^{\triangle}_{0,i})_{i\in\{2,\ldots,n+2\}},
\end{equation}
\begin{equation}\label{EQ125}
(\Lambda_i)_{i\in \{0,n+2\}}. 
\end{equation}
We claim that 
the field extension of $k$ generated by
 the members of the family (\ref{EQ122})$\cup$(\ref{EQ123})$\cup$(\ref{EQ124})$\cup$(\ref{EQ125})
 has transcendence degree $\geq |I|+n+3$.
Indeed 
the field extension of $k$ generated by the members of the family (\ref{EQ119})$\cup$(\ref{EQ120})$\cup$(\ref{EQ121}) is contained in the field extension of $k$ generated  by the members of the family 
(\ref{EQ122})$\cup$(\ref{EQ123})$\cup$(\ref{EQ124})$\cup$(\ref{EQ125}); for if $i\in\{1,\ldots,n+2\}$ then
 $\max\{0,i\}=\max \{i-1,i\}=i$, so 
$Y_{0,i}=Q^{(i)}(Q^{(i)})^*=Y_{i-1,i}$, hence
$y_{0,i}=y_{i-1,i}$, and therefore
$$\Lambda_{i-1}=\Lambda_0\frac{z_{i-1,i}}{z_{0,i}}.$$
Hence the 
field extension of $k$ generated by the 
members of the family (\ref{EQ122})$\cup$(\ref{EQ123})$\cup$(\ref{EQ124})
has transcendence degree $\geq |I|+n+1$.
Consider the three families
\begin{equation}\label{EQ126}
(\mathcal E(f_{j;0,i_1, i_1+1, i_1+2,\ldots,0,i_N, i_N+1, i_N+2})\Lambda_0^j\cdots \Lambda_{i_N+2}^j)_{(j,i_1,\ldots,i_N)\in I},
\end{equation}
\begin{equation}\label{EQ127}
(\mathcal E(f_{i-1,i})\Lambda_{i-1}\Lambda_i)_{i\in\{1,\ldots,n+2\}},
\end{equation}
\begin{equation}\label{EQ128}
(\mathcal E(f_{0,i})\Lambda_0\Lambda_i)_{i\in\{2,\ldots,n+2\}}.
\end{equation}
By Equation (\ref{EQ113}) and Remark \ref{R26} it follows that the 
field extension of $K$ generated by the
members of the family (\ref{EQ126})$\cup$(\ref{EQ127})$\cup$(\ref{EQ128})
has transcendence degree $\geq |I|+n+1$.
By the Serre--Tate expansion principle (Proposition \ref{P3} (b)) the $R$-algebra homomorphism
\begin{equation}\label{EQ129}
\bigoplus_{w\leq 0}\mathbb I^{n+2}_g(w)\rightarrow R[[T,\ldots,T^{(n+2)}]][\Lambda_0,\ldots,\Lambda_{n+2}]\end{equation}
defined by 
$$f\mapsto \mathcal E(f)\cdot \Lambda_0^{c_0}\cdots \Lambda_{n+2}^{c_{n+2}},\ \ \ f\in \mathbb I^{n+2}_g(w),\ \ w=-\sum_{j=0}^{n+2} c_j\phi^j\leq 0,$$
is injective. Consider the three families
 \begin{equation}\label{EQ130}
(f_{j;0,i_1, i_1+1, i_1+2,\ldots,0,i_N, i_N+1, i_N+2})_{(j,i_1,\ldots,i_N)\in I},
\end{equation}
\begin{equation}\label{EQ131}
(f_{i-1,i})_{i\in\{1,\ldots,n+2\}},
\end{equation}
\begin{equation}\label{EQ132}
(f_{0,i})_{i\in\{2,\ldots,n+2\}}.
\end{equation}
We get that the 
field extension of $K$ generated by the 
members of the family (\ref{EQ130})$\cup$(\ref{EQ131})$\cup$(\ref{EQ132})
has transcendence degree $\geq |I|+n+1$  and we conclude by Equation (\ref{EQ116}).
 \endproof
 
 \subsubsection{Applications to some loci in $\mathcal A_g$}\label{S3703}
 
 We end by showing that certain remarkable Zariski dense loci in $\mathcal A_g$ are contained in zero loci of Hecke covariant Siegel $\delta$-modular forms.
 
We say that two matrices $M_1,M_2\in\pmb{\Mat}_g(R)$ {\it pseudo-commute} if there exists $c\in R^{\times}$ such that $M_1M_2=cM_2M_1$. For $(A,\theta,\omega)\in \pmb{\M}_g(R)$, we recall that $\textup{End}(A,\theta)^{(p)}$ is the subring of $\textup{End}(A)$ generated by all polarized isogenies $u:A\rightarrow A$. Consider the map
 \begin{equation}
 \label{onemoreextraeqn}
 [ \ ]_{\omega}:\textup{End}(A,\theta)^{(p)}\rightarrow \pmb{\Mat}_g(R)\end{equation} 
defined by the rule $u\mapsto [u]_{\omega}$, where the matrix $[u]_{\omega}$ is uniquely determined by the identity $u^*\omega=[u]_{\omega} \omega$. We claim that the map (\ref{onemoreextraeqn}) is a ring monomorphism. Indeed, embedding $R$ into the complex field $\mathbb C$ we are reduced to checking the claim for $R$ replaced by $\mathbb C$; however, in this case, our claim follows from the fact that for each abelian variety $A$ over $\mathbb C$, the endomorphism ring $\textup{End}(A)$
is the ring of $\mathbb C$-linear endomorphisms of $H^0(A(\mathbb C),\Omega)^*$ 
 that send the lattice $H_1(A(\mathbb C),\mathbb Z)$ into itself.
 
\begin{lemma}\label{L22}
 Let $F^{a_1,\ldots,a_{2s}}$ be an arbitrary matrix cyclic product and recall that  $f_{a_1,\ldots,a_{2s}}\in\mathbb  I^r_g(-(\sum_{i=1}^{2s}\phi^{a_i}))$ is its trace (see Equation (\ref{finIr})).
 Then for every $(A,\theta,\omega)\in \pmb{\M}_g(R)$ the following two properties hold:
 
 \medskip
 {\bf (a)} The matrix $F^{a_1,\ldots,a_{2s}}(A,\theta,\omega)\in \pmb{\Mat}_g(R)$ pseudo-commutes with each matrix of the form $M^{\phi^{a_1}}$ with $M\in\Im([ \ ]_{\omega})$.
 
 \smallskip
 {\bf (b)} If $f_{a_1,\ldots,a_{2s}}(A,\theta,\omega)\neq 0$, then part (a) holds with `pseudo-commutes' replaced by `commutes'.
 \end{lemma}
 
 \noindent
 {\it Proof.} Part (a) follows from Equation (\ref{EQ110}). Part (b) also follows by taking traces in Equation (\ref{EQ110}) and using part (a).\endproof
 
 \begin{thm}\label{T18}
 For each integer $s\geq 2$ and every non-negative distinct integers $a_1,\ldots,a_{2s}$ there exists a non-zero form
 $$\mathcal D^{a_1,\ldots,a_{2s}}\in \mathbb I_g(-(g^2-g+1)\sum_{i=1}^{2s}\phi^{a_i})$$
 such that for all
 $(A,\theta,\omega)\in \pmb{\M}_g(R)$, if
 $\mathcal D^{a_1,\ldots,a_{2s}}(A,\theta,\omega)\neq 0$,
 then the ring $\End(A,\theta)^{(p)}$ is commutative.
 \end{thm}
 
 \noindent
 {\it Proof.}
 For a monic polynomial $C(t)=t^g+c_1t^{g-1}+\cdots+c_g\in B[t]$ with coefficients in a ring $B$ let
 $\textup{Disc}_g^0(C(t)):=\textup{Disc}_g(1,c_1,\ldots,c_g)\in B$, where $\textup{Disc}_g$ is the discriminant of the generic binary (i.e., in $2$ indeterminates) form of degree $g$. Also for $M\in\pmb{\Mat}_g(B)$ we let $\chi_M(t):=\det(t\cdot 1_g-M)\in B[t]$ be its characteristic polynomial. We can assume that $a_1<\cdots <a_{2s}$.
 Define 
 $$D^{a_1,\ldots,a_{2s}}:=\textup{Disc}_g^0(\chi_{F^{a_1,\ldots,a_{2s}}}(t))\in M_g$$
 where $F^{a_1,\ldots,a_{2s}}$ is the corresponding matrix cyclic product. The coefficients 
 $c_1,\ldots,c_g$ of $\chi_{F^{a_1,\ldots,a_{2s}}}(t)$ belong to $M_g$ and, by Newton's formulas, for each $j\in \{1,\ldots,g\}$, $c_j$ multiplied by a suitable power of $p$ is a weighted homogeneous polynomial of degree $j$, with integer coefficients,
 in the forms $\Trace((F^{a_1,\ldots,a_{2s}})^a)$, $a\in\{1,\ldots,j\}$, with respect to the weights
 $1,\ldots,j$. So for each $j$ we have $c_j\in \mathbb I_g(-j\sum_{i=1}^{2s}\phi^{a_i})$. 
 Also $\text{Disc}^0_g$ is a weighted homogeneous polynomial 
 of degree $g(g-1)$, with respect to the weights $1,\ldots,g$. Hence it follows that
 $D^{a_1,\ldots,a_{2s}}\in \mathbb I_g(-g(g-1)\sum_{i=1}^{2s}\phi^{a_i})$. 
 Set 
 $$\mathcal D^{a_1,\ldots,a_{2s}}:=f_{a_1,\ldots,a_{2s}}D^{a_1,\ldots,a_{2s}}.$$
 Note that, by Lemma \ref{L22} (b), if $\mathcal D^{a_1,\ldots,a_{2s}}$ does not vanish at some triple $(A,\theta,\omega)\in \pmb{\M}_g(R)$, then every element of $\phi^{a_1}(\textup{Im}([\ ]_{\omega}))$ commutes with the matrix $F^{a_1,\ldots,a_{2s}}(A,\theta,\omega)$ that has distinct eigenvalues. Hence 
 $\phi^{a_1}(\textup{Im}([\ ]_{\omega}))$
 is contained in an abelian $R$-subalgebra of $\pmb{\Mat}_g(R)$.
In order to conclude the proof we need to check that $D^{a_1,\ldots,a_{2s}}\neq 0$.  It suffices to show that the image $\mathcal E(D^{a_1,\ldots,a_{2s}})^{\triangle}$ of $\mathcal E(D^{a_1,\ldots,a_{2s}})$ via  the reduction modulo $p$ homomorphism between  rings of polynomials
$$\triangle:R[Q^{(1)},\ldots,Q^{(2s)}]\rightarrow k[Q^{(1)},\ldots,Q^{(2s)}]$$ is non-zero. 
 But, by Theorem \ref{T13} (a), 
 we have 
 $$\mathcal E(f^r_{g,\crys})\equiv \Psi^{\phi^{r-1}}\ \ \textup{mod}\ \ p$$
 hence, using Equations (\ref{deffab}), (\ref{defFa1}), (\ref{deffa1}), and (\ref{defQ(s)}) 
 and the equalities
 $$\max\{a_1,a_2\}=a_2,\ \ldots, \max\{a_{2s-1},a_{2s}\}=\max\{a_{2s},a_1\}=a_{2s}$$
 we get that 
 the mentioned image equals the  polynomial
\begin{equation}
\label{mathcalEis}
\textup{Disc}^0_g(\chi_{Q^{(a_2)}(Q^{(a_3)})^*\cdots Q^{(a_{2s-2})}(Q^{(a_{2s-1})})^*Q^{(a_{2s})}(Q^{(a_{2s})})^*}(t)).\end{equation}
 As $s\geq 2$, evaluating the matrices of indeterminates $Q^{(a_3)},\ldots,Q^{(a_{2s})}$ at $1_g$ we get the non-zero polynomial $\textup{Disc}^0_g(\chi_{Q^{(a_2)}}(t))$, so $D^{a_1,\ldots,a_{2s}}\neq 0$.
\endproof
 
\begin{thm}\label{T18.5}
Let $s\geq 2$ be an integer and let
 $$a_1<\cdots<a_{2s}\ \ \ \text{and}\ \ \ b_1<\cdots<b_{2s}$$
 be non-negative integers such that $a_1=b_1$ and $\{a_2,a_3\}\cap \{b_2,b_3\}=\emptyset$.
 Define
 $$m_g:=\sum_{j=1}^{g-1} j\binom{g}{j},\ \ \ w_a:=\sum_{i=1}^{2s}\phi^{a_i},\ \ \ w_b:=\sum_{i=1}^{2s}\phi^{b_i}.$$
 Then there exists a non-zero form
 $$\mathcal G\in \mathbb I_g(-(m_g+1)(w_a+w_b))$$
 such that for all
 $(A,\theta,\omega)\in \pmb{\M}_g(R)$ with $\mathcal G(A,\theta,\omega)\neq 0$, the polarized isogeny class of $(A,\theta)$ is indecomposable.
\end{thm}
 
 \noindent 
 {\it Proof.}
 For a ring $C$, matrix $M\in\pmb{\Mat}_g(C)$ and integer $j\in\{1,\ldots,g-1\}$
let $M^{\wedge j}\in\pmb{\Mat}_{\binom{g}{j}}(C)$ be the matrix representation of the $j$-th wedge power of the $C$-linear endomorphism of $C^g$ defined by $M$ with respect to the standard basis 
$$\{e_{i_1}\wedge\cdots\wedge e_{i_j}|1\leq i_1<\cdots<i_j\leq g\}$$ of $\bigwedge^{\raisebox{-0.0ex}{\scriptsize $j$}} C^g$ ordered lexicographically, where $\{e_1,\ldots,e_g\}$ is the standard basis of $C^g$. For $M_1,M_2\in\pmb{\Mat}_g(C)$ we have $(M_1M_2)^{\wedge j}=M_1^{\wedge j}M_2^{\wedge j}$ and we set
$$\Phi_j(M_1,M_2):=\det(M_1^{\wedge j}M_2^{\wedge j}-M_2^{\wedge j}M_1^{\wedge j})\in C.$$
 Consider the matrix cyclic products
 $$F_a:=F^{a_1,\ldots,a_{2s}},\ \ \ F_b:=F^{b_1,\ldots,b_{2s}}.$$
As $a_1=b_1$ the products $F_aF_b$ and $F_bF_a$ are also matrix cyclic products:
 $$F_aF_b=F^{a_1,\ldots, a_{2s},b_1,\ldots,b_{2s}},\ \ \ F_bF_a=F^{b_1,\ldots, b_{2s},a_1,\ldots,a_{2s}}.$$
 Define
 $G_j:=\Phi_j(F_a,F_b)\in \mathbb I_g$, $j\in \{1,\ldots,g-1\}$. Using Equation (\ref{EQ029.9})
 one easily checks that
 $$G_j\in \mathbb I_g(-j\binom{g}{j}(w_a+w_b)).$$
 Thus $\prod_{j=1}^{g-1} G_j\in \mathbb I_g(-m_g(w_a+w_b))$. From this  and Equation (\ref{finIr}) we get:
 $$\mathcal G:=f_{a_1,\ldots,a_{2s}}f_{b_1,\ldots,b_{2s}}\prod_{j=1}^{g-1} G_j\in \mathbb I_g(-(m_g+1)(w_a+w_b)).$$
 We claim that $\mathcal G$ has the desired properties. We first prove that $\mathcal G$ is non-zero. It is enough to check that each $G_j$ is non-zero. As in the proof of Theorem 
 \ref{T18} it is enough to show that the image $\mathcal E(G_j)^{\triangle}$ of $\mathcal E(G_j)$ in the ring $k[Q^{(1)},\ldots,Q^{(\max\{a_{2s},b_{2s}\})}]$ is non-zero. 
  Evaluating the matrices of indeterminates 
  $Q^{(a_4)},\ldots,Q^{(a_{2s})}$ and $Q^{(b_4)},\ldots,Q^{(b_{2s})}$
   at $1_g$
  one gets the expression
 $$\widetilde{G_j}:=\Phi_j(Q^{(a_2)}(Q^{(a_3)})^*,Q^{(b_2)}(Q^{(b_3)})^*).$$
As $\{a_2,a_3\}\cap \{b_2,b_3\}=\emptyset$ the entries of the symmetric matrices 
$Q^{(a_2)}$, $(Q^{(a_3)})^*$, $Q^{(b_2)}$, $(Q^{(b_3)})^*$ that are above and on the diagonal are algebraically independent.
 So it is enough to find symmetric matrices $A_2,A_3,B_2,B_3\in \pmb{\Mat}_g(k)$ such that
 $\Phi_j(A_2A_3^*,B_2B_3^*)\neq 0$. Each diagonalizable matrix in $\pmb{\GL}_g(k)$
 is a product of two symmetric matrices (see Remark \ref{R37} in Subsubsection \ref{S415}).
 So it is enough to find two diagonalizable matrices $D_1,D_2\in \pmb{\GL}_g(k)$ such that
 $\Phi_j(D_1,D_2)\neq 0$. This follows from Lemma \ref{L36} in Subsubsection \ref{S414} and from the fact that the set of diagonalizable matrices is Zariski dense in the set of all matrices.
 
 To end our proof we need to show that if the polarized isogeny class $(A,\theta)$ is decomposable then 
 $\mathcal G(A,\theta,\omega)=0$ for all  $\omega$. By Hecke covariance
 we may assume $(A,\theta)=(A_1,\theta_1)\times_{\textup{Spec}(R)} (A_2,\theta_2)$  and we may take $\omega=(\textup{pr}_1^*(\omega_1),\textup{pr}_2^*(\omega_2))^{\t}$ where $\textup{pr}_1,\textup{pr}_2$ are the projections of $A$ onto $A_1,A_2$ (respectively) and $\omega_1,\omega_2$ are column bases of $1$-forms on $A_1,A_2$ (respectively). For $m\in\mathbb Z$ let $[m]_i:A_i\rightarrow A_i$, $i\in \{1,2\}$ be the multiplication by $m$ on $A_i$. Then we may consider the polarized isogenies,
 $$[1]_1\times [1]_2, [1]_1\times [-1]_2:A\rightarrow A.$$
 Their sum $u:=[2]_1\times [0]_2$ belongs to $\End(A,\theta)^{(p)}$. 
 Let $g_1$ and $g_2$ be the relative dimensions of $A_1$ and $A_2$ so the relative dimension of $A$ is $g=g_1+g_2$.
 Note that the matrix
 $[u]_{\omega}$ of $u$ with respect to $\omega$ has the block form 
 $\left(\begin{array}{cc} 2 \cdot 1_{g_1}& 0_{g_1\times g_2}\\ 0_{g_2\times g_1} & 0_{g_2\times g_2}\end{array}\right)$ where, for $i,j\in \{1,2\}$, we denoted by $0_{g_i\times g_j}$  the $g_i\times g_j$ matrix with zero entries.
 By Lemma \ref{L22} (b) this matrix commutes with $F_a(A,\theta,\omega)$ and $F_b(A,\theta,\omega)$ hence the latter matrices have entries in $R$ and a block form 
 $\left(\begin{array}{cc} \star & 0\\ 0 & \star\end{array}\right)$ so the endomorphisms over $K$ defined by these matrices have an invariant subspace of dimension $g_1\in \{1,\ldots,g-1\}$ in common. By Lemma \ref{L37} in Subsubsection \ref{S414} we get $\Phi_j(F_a(A,\theta,\omega),F_b(A,\theta,\omega))=0$, hence $\mathcal G(A,\theta,\omega)=0$.\endproof

 \subsection{$p$-adic approximation}\label{S38}

In this subsection we show how to $p$-adically approximate elements of $\mathbb I^r_{g,\ord}$ by quotients of elements in $\mathbb I^r_g$. We begin by introducing some new concepts.

Let 
$$W_{\even}(r):=\{w\in W(r)|\deg(w)\in 2\mathbb Z\}$$ 
and let
$\mathbb D\subset \mathbb I_g$ be a multiplicative set such that the following two conditions hold:

\medskip
\noindent
{\bf (1)} We have $\mathbb D\subset \cup_{r=1}^{\infty}\cup_{w\in W_{\even}(r)} \mathbb I_g^r(w)\backslash p \mathbb I^r_g(w)$.

\smallskip
\noindent
{\bf (2)} We have $\mathbb D^{\phi}\subset\mathbb D$.

\medskip

\noindent An element $f\in\mathbb D^{-1}\mathbb I_g$ is said to have weight $w\in W$
if there exist $r\in\mathbb N$, $v\in W_{\even}(r)$, $F\in\mathbb I^r_g(v+w)$, and $G\in\mathbb D\cap  \mathbb I^r_g(v)$,
such that $f=G^{-1}F$. Denote by $(\mathbb D^{-1}\mathbb I_g)(w)$ the $R$-submodule of $\mathbb D^{-1}\mathbb I_g$ of elements of weight $w$. Similarly, an element $f\in\mathbb D^{-1}\mathbb I_{g,\ord}$ is said to have weight $w\in W$
if there exist $r\in\mathbb N$, $v\in W_{\even}(r)$, $F\in\mathbb I^r_{g,\ord}(v+w)$, and $G\in\mathbb D\cap \mathbb I^r_{g,\ord}(v)$,
such that $f=G^{-1}F$. Denote by $(\mathbb D^{-1}\mathbb I_{g,\ord})(w)$ the $R$-submodule of $\mathbb D^{-1}\mathbb I_{g,\ord}$ of elements of weight $w$. Note that if a non-zero  element
of $\mathbb D^{-1}\mathbb I_{g}$ or $\mathbb D^{-1}\mathbb I_{g,ord}$ has a weight $w$, then $w$ is uniquely determined. 
 
Fix a pair $(w',w'')\in\mathcal W\times \mathcal W$.
We say that an element $f\in\mathbb D^{-1} M_{gg}$ 
has weight $(w',w'')$ if there exists $r\in\mathbb N$ and a triple
 $$(w,F,G)\in W_{\even}(r)\times M^r_{gg}\times (\mathbb I_g^r(w)\backslash p \mathbb I^r_g(w))$$ 
 such that  the following two conditions hold:
 
 \medskip
 
 {\bf (a)} We have $f=G^{-1}F$.
 
 \smallskip
 
{\bf (b)} For each $S \in\Ob(\Prol)$, for all $(A,\theta,\omega) \in \pmb{\M}_g(S^0)$
and all $\lambda \in\pmb{\GL}_g(S^0)$ we have
\begin{equation}\label{EQ133}
F(A,\theta,\lambda \omega,S)= \det(\lambda)^{-w}\cdot \chi_{w'}(\lambda) \cdot
F(A,\theta,\omega,S) \cdot \chi_{w''}(\lambda)^{\t}.\end{equation}

\medskip

Note that if the above is the case and  $(\tilde{w},\tilde{F},\tilde{G})$ is another triple in the set
$W_{\even}(r)\times M^r_{gg}\times (\mathbb I_g^r(w)\backslash p \mathbb I^r_g(w))$ such that condition (a) holds for $(\tilde{F},\tilde{G})$, then condition (b) holds for $(\tilde{w},\tilde{F})$.
The $R$-submodule of $\mathbb D^{-1}M_{gg}$ of all elements of 
 weight $(w',w'')$ will be denoted by 
 $$(\mathbb D^{-1}M_{gg})(w',w'').$$
Similar definitions can be made for $\mathbb D^{-1}M_{gg,\ord}$ in place of $\mathbb D^{-1}M_{gg}$. As in the case of Definition \ref{df21}, a non-zero element of $\mathbb D^{-1} M_{gg}$ or $\mathbb D^{-1}M_{gg,\ord}$ could (at least for $g=1$) have two different weights $(w',w'')$ and one expects that $w'+w''=0$.

We remark that if $f\in (\mathbb D^{-1}M_{gg})(w',w'')\cap (\mathbb D^{-1}M_{gg})^{\times}$, then we have $f^{-1}\in (\mathbb D^{-1}M_{gg})(-w'',-w')$.

Let $n\in\mathbb N$. Let $\{f_1,\ldots,f_n\}$ be a set of elements in $(\mathbb D^{-1}M_{gg,\ord})^{\times}$ that have weights in $\mathcal W\times \mathcal W$ and let $(w',w'')\in\mathcal W\times\mathcal W$. By a {\it coherent product} of elements in $\{f_1,\ldots,f_n\}$ of weight $(w_1,w_2)$ we mean a product of the form
$$f^{\bullet}_1 f_2^{\bullet}\cdots f^{\bullet}_N\in (\mathbb D^{-1}M_{gg,\ord})(w',w'')\cap
(\mathbb D^{-1}M_{gg,\ord})^{\times}$$ for which the following two conditions hold:

\medskip

\noindent {\bf (i)}  For each $i\in\{1,\ldots,N\}$ there exist $j\in\{1,\ldots,n\}$ and an integer $s\geq 0$ with the property that $f^{\bullet}_i$ equals either $f_j^{\phi^s}$ or $(f_j^{\phi^s})^{\t}$.

\smallskip

\noindent {\bf (ii)} For each $i\in\{1,\ldots,N\}$ there exists $(w'_i,w''_i)\in\mathcal W\times\mathcal W$ such that the form $f^{\bullet}_i$ has weight 
 $(w'_i,w''_i)$ and moreover $w_1'=w'$, $w''_N=w''$, and for all $i\in\{1,\ldots,N-1\}$ we have $w''_i+w'_{i+1}=0$.

\medskip

\noindent Note that if $f_1,\ldots,f_n\in (\mathbb D^{-1}M_{gg})^{\times}$, then, in the notation above, we have: 
$$f^{\bullet}_1 f_2^{\bullet}\cdots f^{\bullet}_N\in (\mathbb D^{-1}M_{gg})(w',w'')\cap
 (\mathbb D^{-1}M_{gg})^{\times}.$$
 
 For $(w',w'')\in\mathcal
W\times\mathcal W$, an element $\pi$ in the $p$-adic completion
$\reallywidehat{\mathbb D^{-1} M_{gg,\ord}}$ will be called a {\it
coherent combination} of elements in $\{f_1,\ldots,f_n\}$ of weight
$(w',w'')$ if we have a series representation $\pi=\sum_{i=0}^{\infty} a_i \pi_i$ where the
sequence $(a_i)_{i\geq 0}$ in $R$ converges to $0$ in the $p$-adic
topology and each $\pi_i\in\mathbb D^{-1}
M_{gg,\ord}(w',w'')$ is a coherent product of elements in
$\{f_1,\ldots,f_n\}$ of weight $(w',w'')$.  We will say $\pi$ is {\it regular} if in the series representation of $\pi$ we can choose $a_0\not\in\mathfrak m_R$ and $a_i\in\mathfrak m_R$ for all $i\geq 1$. We will say $\pi$ is {\it finite} if in the series representation of $\pi$ we can choose  $a_i=0$ for all but finitely many $i$.
If $f_1,\ldots,f_n$ belong to $(\mathbb D^{-1}M_{gg})^{\times}$, then $\pi$ comes from an element in $\reallywidehat{\mathbb D^{-1} M_{gg}}$; if, in addition, $\pi$ is finite, then $\pi$ comes from an element in $\mathbb D^{-1} M_{gg}$.

In the notation above we have:

\begin{lemma}\label{L23}
We assume $n\geq 2$. Let $\pi_1$ be a coherent product of elements in the set
$$\Sigma:=\{f_2,\ldots,f_n, f_2^{-1},\ldots,f_n^{-1}\},$$
 of weight $(w',w'')$, and let $\pi_2$ be a coherent combination of elements in 
 the set $\Sigma\cup\{f_1\}$,
 of the same weight $(w',w'')$.
Then 
$$(\pi_1+p\cdot \pi_2)^{-1}=\pi_1^{-1}+p\cdot \pi_3$$
where $\pi_3$ is a coherent combination of elements in the set $\Sigma\cup\{f_1\}$, of the same weight
$(-w'',-w')$ as $\pi_1^{-1}$. In particular, $(\pi_1+p\cdot\pi_2)^{-1}$ is regular.
\end{lemma}

\noindent
{\it Proof.}
Based on the identity
$$(\pi_1+p\cdot \pi_2)^{-1}=\sum_{i=0}^{\infty} (-p)^i(\pi_1^{-1}\pi_2)^i\pi_1^{-1},$$
we can take $a_i:=(-p)^i$ for $i\in\mathbb N\cup\{0\}$ and the lemma follows.\endproof

\medskip

In what follows set $f^r:=f^r_{g,\crys}$ for $r\in\mathbb N$ and assume that, in addition to the conditions (1) and (2) above, the multiplicative set $\mathbb D$ also satisfies the third condition:

\medskip

\noindent
{\bf (3)} $\det(f^1),\det(f^2)\in\mathbb D$.

\medskip

\noindent For instance, as $\mathbb D$ we could take the minimal set
\begin{equation}
\label{Dmin}
\mathbb D_{\min}:=\{\det(f^1)^{w'}\det(f^2)^{w''}|w',w''\in W_+\}\end{equation}
(contained in the image of
the $R$-algebra monomorphism $\diamondsuit:\mathbb H_g^{\torus}\rightarrow \mathbb I_{g,\ord}$, see Theorem \ref{T16} (h)) or the maximal set
$$\mathbb D_g:= \cup_{r=1}^{\infty} \cup_{w\in W_{\even}(r)} \mathbb (\mathbb I^r_g(w)\backslash p \mathbb I^r_g(w)).$$
The key fact in our discussion below is that we have the equality:
\begin{equation}
\label{keyeq}
f^2=(f^1)^{\phi}\cdot f^{\partial}+p\cdot (((f^{\partial})^{\phi})^{-1})^{\t} \cdot f^1.\end{equation}
Indeed, by Theorem \ref{T13} (a) and Theorem \ref{T14} (a), the Serre--Tate expansions of the left and right hand sides of Equation (\ref{keyeq}) coincide. On the other hand the two sides have the same weight.  So the equality in Equation (\ref{keyeq}) holds by the Serre--Tate expansion principle (Proposition \ref{P3} (b)).
In particular, by setting $f_{\partial}:=(f^{\partial})^{-1}$, in $\mathbb D^{-1} M_{gg,\ord}$ we get an equality
\begin{equation}\label{EQ134}
f^{\partial}=((f^1)^{\phi})^{-1}\cdot f^2-
p\cdot ((f^1)^{\phi})^{-1}\cdot ((f_{\partial})^{\phi}))^{\t} \cdot f^1
\end{equation}
where $((f^1)^{\phi})^{-1}\cdot f^2$ is a coherent product of elements in the set
\begin{equation}\label{EQ135}
\Pi:=\{f^1, (f^1)^{-1},f^2,(f^2)^{-1}\}\end{equation}
of weight $(-\phi,1)$ and $((f^1)^{\phi})^{-1}\cdot ((f_{\partial})^{\phi}))^{\t} \cdot f^1$ is a coherent product of elements in the set
\begin{equation}\label{EQ136}
\Pi':=\{f_{\partial},f^1, (f^1)^{-1},f^2,(f^2)^{-1}\}=\Pi\cup\{f_{\partial}\}
\end{equation}
of the same weight.
In particular, by Lemma \ref{L23}, in $\mathbb D^{-1}M_{gg,\ord}$ we have 
\begin{equation}\label{EQ137}
f_{\partial}=(f^2)^{-1}(f^1)^{\phi}+p\pi_1,
\end{equation}
where $(f^2)^{-1}(f^1)^{\phi}$ is a coherent product of elements in the set $\Pi$ and $\pi_1$ is a coherent combination of elements in the set $\Pi'$, both of them having weight $(-1,\phi)$.

\begin{lemma}\label{L24}
For every integer $n\geq 2$ there exist a finite regular coherent combination $\eta_n\in\mathbb D^{-1}M_{gg}$ of elements in Equation (\ref{EQ135}) and a coherent combination $\pi_n\in\reallywidehat{\mathbb D^{-1}M_{gg,\ord}}$
of elements in Equation (\ref{EQ136}), both of weight $(-1,\phi)$, such that in $\reallywidehat{\mathbb D^{-1}M_{gg,\ord}}$ we have an identity
\begin{equation}\label{EQ138}
f_{\partial}=\eta_n+p^n\pi_n.\end{equation}
\end{lemma}

\noindent
{\it Proof.}
We proceed by induction on $n$. The base of the induction for $n=2$ holds due to the Equation (\ref{EQ137}).
Assuming that the Equation (\ref{EQ138}) holds for some $n\geq 2$, the passage from $n$ to $n+1$ is achieved by substituting
 all occurrences of $f_{\partial}$ in $\pi_n$ in that equation by the right-hand side of Equation (\ref{EQ137}) and regrouping the terms to get the desired form. 
\endproof

\begin{cor}\label{C18.1}
For every integer $n\geq 2$ there exist a finite regular coherent combination $\eta^*_n\in\mathbb D^{-1}M_{gg}$ of elements in Equation (\ref{EQ135}) and a coherent combination $\pi^*_n\in\reallywidehat{\mathbb D^{-1}M_{gg,\ord}}$
of elements in Equation (\ref{EQ136}), both of weight $(-\phi,1)$, such that in $\reallywidehat{\mathbb D^{-1}M_{gg,\ord}}$ we have an identity
$$f^{\partial}=\eta^*_n+p^n\pi^*_n.$$
\end{cor}

\noindent
{\it Proof.}
Replace $f_{\partial}$ in the right-hand side of Equation (\ref{EQ134}) by the right-hand side of Equation (\ref{EQ138}).
\endproof

\medskip

Recall the forms $f^{\langle a \rangle}$ from Equation (\ref{EQ064}) indexed by $a\in\mathbb N$; then Corollary \ref{C18.1} immediately implies the following:

\begin{cor}\label{C18.2}
For all integers $a\geq 1$ and $n\geq 2$ there exist a finite regular coherent combination $\eta^{(a)}_n\in\mathbb D^{-1}M_{gg}$ of elements in Equation (\ref{EQ135}) and a coherent combination $\pi^{(a)}_n\in\reallywidehat{\mathbb D^{-1}M_{gg,\ord}}$
of elements in Equation (\ref{EQ136}), both of weight $(1,1)$, such that in $\reallywidehat{\mathbb D^{-1}M_{gg,\ord}}$ we have an identity
$$f^{\langle a \rangle}=\eta^{(a)}_n+p^n\pi^{(a)}_n.$$
\end{cor}

Recall that we have $R$-algebra monomorphisms
$$\mathbb H_{g,\tot}^{\torus}\stackrel{\diamondsuit}\rightarrow
\mathbb I_{g,\ord}
\rightarrow \mathbb D^{-1}\mathbb I_{g,\ord}\rightarrow \mathbb D^{-1}M_{g,\ord} \rightarrow\reallywidehat{\mathbb D^{-1}M_{g,\ord}}$$
which, for $s\in \frac{1}{2}\mathbb Z$, induce $R$-linear maps 
$$\mathbb H_g(s) \rightarrow \mathbb I_{g,\ord}(-2s)\rightarrow (\mathbb D^{-1} \mathbb I_{g,\ord})(-2s)\rightarrow \mathbb D^{-1} \mathbb I_{g,\ord},$$
where 
$$\mathbb H_g(s):=\bigcup_{r\geq 1} \mathbb H^r_g(s),\ \ \ \mathbb I_{g,\ord}(-2s):=\bigcup_{r\geq 1} \mathbb I^r_{g,\ord}(-2s).$$
Also, we have $R$-algebra homomorphisms
$$\mathbb D^{-1}\mathbb I_g \rightarrow 
\mathbb D^{-1}\mathbb I_{g,\ord} \rightarrow 
\mathbb D^{-1}M_{g,\ord} \rightarrow\reallywidehat{\mathbb D^{-1}M_{g,\ord}},$$
that induce for each $w\in W$ the $R$-linear maps
$$(\mathbb D^{-1}\mathbb I_g)(w)\rightarrow 
 (\mathbb D^{-1}\mathbb I_{g,\ord})(w)\rightarrow\reallywidehat{\mathbb D^{-1}M_{g,\ord}}.$$

\begin{lemma}\label{L25}
For every $w_0\in W$ the reductions of the homomorphism
$$(\mathbb D^{-1}\mathbb I_{g,\ord})(w_0)\rightarrow\reallywidehat{\mathbb D^{-1}M_{g,\ord}}$$
modulo powers of $p$ are injective.
\end{lemma}

\noindent
{\it Proof.}
Let $n\in\mathbb N$ and let $G^{-1}F\in (\mathbb D^{-1}\mathbb I_{g,\ord})(w_0)$ where $G\in \mathbb D \cap\mathbb I_g(w)$ and $F\in\mathbb I_{g,\ord}(w+w_0)$ with $w\in W$ of even degree.
We assume that the image of $G^{-1}F$ in 
$\reallywidehat{\mathbb D^{-1}M_{g,\ord}}$ belongs to the submodule $p^n \reallywidehat{\mathbb D^{-1}M_{g,\ord}}$. Then the image of $G^{-1}F$ in $\mathbb D^{-1}M_{g,\ord}$ belongs to $p^n \mathbb D^{-1}M_{g,\ord}$. So there exist $G_1\in\mathbb D\cap \mathbb I_g(w_1)$, $H\in\mathbb D\cap \mathbb I_g(v)$ for some $v,w_1\in W$ of even degree, and $F_1\in M_{g,\ord}$ such that the image of $HG_1F$ in $M_{g,\ord}$ equals the image of 
$p^nHF_1G$. By the flatness over $R$ condition in the definition of the objects of $\Prol$ (see Definition \ref{df5}),
as $HG_1F$ is Hecke covariant of weight $v+w+w_1+w_0$, we get that
$HF_1G$ is also Hecke covariant of the same weight. By Corollary \ref{C5} (c) we conclude that $F_1$ is Hecke covariant of weight $w_1+w_0$. So we have $G_1^{-1}F_1\in (\mathbb D^{-1}\mathbb I_{g,\ord})(w_0)$. As $HG_1F$ and $p^n HGF_1$ coincide in $M_{g,\ord}$, they coincide in $\mathbb I_{g,\ord}(v+w+w_1+w_0)$. So $G^{-1}F=p^n G_1^{-1}F_1\in p^n(\mathbb D^{-1}\mathbb I_{g,\ord})(w_0)$.\endproof

\medskip
The proof of the following elementary lemma is left as an exercise.

\begin{lemma}\label{L26}
Let $C\subset B$ be an extension of flat $\mathbb Z_{(p)}$-algebras whose reduction modulo $p$, $C/pC\rightarrow B/pB$, is injective. If $b\in B$ is such that there exists $N\in\mathbb N$ with the property that $p^Nb$ belongs to the $p$-adic closure of $C$ in $B$, then $b$ itself belongs to the $p$-adic closure of $C$ in $B$.
\end{lemma}

Here by the `$p$-adic closure of $C$ in $B$' we mean the closure of $C$ in $B$ where the latter is equipped with its   $p$-adic topology.

\begin{lemma}\label{L27}
The following two properties hold:

\medskip
{\bf (a)} The $R$-algebra homomorphism
$$\mathbb D^{-1}\mathbb I_g\rightarrow \mathbb D^{-1} \mathbb I_{g,\ord}$$
 is injective and its reduction modulo $p$ is injective. 

\smallskip
{\bf (b)} The $R$-algebras 
 $\mathbb D^{-1}\mathbb I_g$ and $\mathbb D^{-1}\mathbb I_{g,\ord}$ are $p$-adically separated integral domains having $p$ as a prime element.
\end{lemma}

\noindent
{\it Proof.} Part (a) follows the Serre--Tate expansion principle (Proposition \ref{P3} (b)) via localization which is an exact functor. Part (b) follows directly from Corollary \ref{C4} (b). 
\endproof

\medskip

In view of Lemma \ref{L27} (a) we will speak about inclusions $\mathbb D^{-1}\mathbb I_g\subset \mathbb D^{-1} \mathbb I_{g,\ord}$ and
 $\reallywidehat{\mathbb D^{-1}\mathbb I_g}
\subset \reallywidehat{\mathbb D^{-1} \mathbb I_{g,\ord}}$. Here is our main approximation result:

\begin{thm}\label{T19}
The following three properties hold:

\medskip
{\bf (a)} For every integer $s\in \frac{1}{2}\mathbb Z$ with $gs\in \mathbb Z$, the image of the $R$-linear map
$\mathbb H_g(s)\rightarrow (\mathbb D^{-1}\mathbb I_{g,\ord})(-2s)$ is contained in the $p$-adic closure of the image of the $R$-linear map
$$(\mathbb D^{-1} \mathbb I_g)(-2s) \rightarrow 
(\mathbb D^{-1}\mathbb I_{g,\ord})(-2s).$$

\smallskip
{\bf (b)} For every $w\in W$, the image of $\det(f^{\partial})^w
\in\mathbb I_{g,\ord}((\phi-1)w)$ in $(\mathbb D^{-1}\mathbb I_{g,\ord})((\phi-1)w)$ is contained in the $p$-adic closure of the image of
$$(\mathbb D^{-1} \mathbb I_g)((\phi-1)w) \rightarrow 
(\mathbb D^{-1}\mathbb I_{g,\ord})((\phi-1)w).$$ 

\smallskip
{\bf (c)} We have
$$\reallywidehat{\mathbb D^{-1} \mathbb I_g}
= \reallywidehat{\mathbb D^{-1}\mathbb I_{g,\ord}}.$$
\end{thm}

\noindent
{\it Proof.} We start by proving part (a). Part (b) can be proved by an entirely similar argument. 

In view of Lemma \ref{L25} and Theorem \ref{T15} (b) it is enough to show that the image in $\reallywidehat{\mathbb D^{-1}M_{g,\ord}}$ of each element of $\mathbb I^r_{g,\ord}(-2s)$ 
of the form $F(f^{\langle 1 \rangle},\ldots,f^{\langle r \rangle})$,
where $F\in\mathbb H^{r-1}_g(s)$, is a $p$-adic limit of images of elements
in $(\mathbb D^{-1}\mathbb I_g)(-2s)$. 

Fix such an $F$ and an integer $n\geq 1$. We seek an element 
$F[n]\in (\mathbb D^{-1} \mathbb I_g)(-2s)$ whose image in the ring $\reallywidehat{\mathbb D^{-1}M_{g,\ord}}$ differs from the image of $F(f^{\langle 1 \rangle},\ldots,f^{\langle r \rangle})$ by an element 
in $p^n\reallywidehat{\mathbb D^{-1}M_{g,\ord}}$.

By Corollary \ref{C18.2} there exist
finite regular coherent combinations 
$$\eta^{(1)}_n,\ldots,\eta^{(r-1)}_n\in\mathbb D^{-1}M_{gg}=\pmb{\Mat}_g(\mathbb D^{-1}M_g)$$
 of elements in the set (\ref{EQ135}) of weight $(1,1)$ such that
\begin{equation}\label{EQ139}
F(f^{\langle 1 \rangle},\ldots,f^{\langle r-1\rangle})\equiv F(\eta^{(1)}_n,\ldots,\eta^{(r-1)}_n)\ \ \ \mod\ \ \ p^n\end{equation}
in the ring $\reallywidehat{\mathbb D^{-1}M_{g,\ord}}$. The elements $\eta^{(a)}_n$ can be written as
\begin{equation}
\label{etaequalsGinversePhi}\eta_n^{(a)}=G^{-1}\Phi_{n,a}\end{equation}
where $G\in\mathbb D_{\textup{min}}$ (see Equation (\ref{Dmin}))  is of some weight $w\in W$ of even degree and $\Phi_{n,a}\in M_{gg}$ satisfy the following two conditions:

\medskip

\noindent {\bf (*)} For
each triple $(S, (A,\theta,\omega),\lambda)\in\Ob(\Prol)\times\pmb{\M}_g(S^0)\times\pmb{\GL}_g(S^0)$ we have:
\begin{equation}
\label{Phina0}
\Phi_{n,a}(A,\theta,\lambda \omega,S)= \det(\lambda)^{-w}\cdot \lambda \cdot
\Phi_{n,a}(A,\theta,\omega,S) \cdot \lambda^{\t}.\end{equation}

\smallskip

\noindent {\bf (**)} For each
$S \in\Ob(\Prol)$ and every triple $(x_1,x_2,u)\in \pmb{\M}'_g(S^0)$, where
$x_1=(A_1,\theta_1,\omega_1)$, $x_2=(A_2,\theta_2,\omega_2)$, $uu^{\t}=d\in\mathbb N\setminus p\mathbb N$, and $u^*\omega_2=\omega_1$, the following equality holds:
\begin{equation}
\label{Phina}
\Phi_{n,a}(A_2,\theta_2,\omega_2,S)= d^{-1+\frac{g\cdot \deg(w)}{2}}\cdot 
\Phi_{n,a}(A_1, \theta_1, \omega_1,S).\end{equation}

\medskip

\noindent Condition (*) 
follows from the fact that $f^1$ and $f^2$ have weights $(\phi,1)$ and $(\phi^2,1)$, respectively.
Condition 
 (**) follows from the Hecke covariance of $f^1$ and $f^2$ (see Equation (\ref{EQ029.4})).
As $F$ is homogeneous of degree $gs$, by Equation (\ref{etaequalsGinversePhi}) we get that
\begin{equation}\label{EQ140}
F(\eta^{(1)}_n,\ldots,\eta^{(r-1)}_n)=(G^{gs})^{-1}\Phi,\ \ \ \Phi:=F(\Phi_{n,1},\ldots,\Phi_{n,r-1})\in M_g.\end{equation}
Using Equation (\ref{Phina0}) 
we get  that for every triple $(S,(A,\theta,\omega),\lambda)$ as in (*) above we have
$$
\Phi(A,\theta,\lambda \omega,S) =\det(\lambda)^{2s-gsw}\cdot \Phi(A,\theta,\omega,S).$$
Thus $\Phi$ has weight $gsw-2s$. Also, by Equation (\ref{Phina}), for every $S$ and $(x_1,x_2,u)$ as in (**) above we have
$$\Phi(A_2,\theta_2,\omega_2,S)=d^{gs(-1+\frac{g\cdot \deg(w)}{2})} \Phi(A_1,\theta_1,\omega_1,S).$$
By Equation (\ref{EQ029.4}) we get that $\Phi\in\mathbb I_g^{r'}(gsw-2s)$ for some $r'\in \mathbb N$. 
(Note that a priori one cannot take $r'=r$.)
Thus we have $(G^{gs})^{-1}\Phi\in (\mathbb D^{-1}\mathbb I_g)(-2s)$.
Then the element $F[n]:=F(\eta^{(1)}_n,\ldots,\eta^{(r-1)}_n)$ has the desired properties in view of Equations (\ref{EQ139}) and (\ref{EQ140}).

By parts (a) and (b), the image $\Im(\mathbb H_{g,\tot}^{\torus}\rightarrow \mathbb D^{-1}\mathbb I_{g,\ord}$) (see Theorem \ref{T16} (a)) is contained in the $p$-adic closure of $\mathbb D^{-1} \mathbb I_g$ in 
$\mathbb D^{-1}\mathbb I_{g,\ord}$. By Lemmas \ref{L26} and \ref{L27} (a),
we get that every element of $\mathbb D^{-1}\mathbb I_{g,\ord}$ is contained in the $p$-adic closure of $\mathbb D^{-1} \mathbb I_g$ in 
$\mathbb D^{-1}\mathbb I_{g,\ord}$. Now, by Lemma \ref{L27} (a) again, each Cauchy sequence in $\mathbb D^{-1}\mathbb I_{g,\ord}$ which is contained in $\mathbb D^{-1} \mathbb I_g$ is as well a Cauchy sequence in $\mathbb D^{-1} \mathbb I_g$. We conclude that every element
of $\reallywidehat{\mathbb D^{-1}\mathbb I_{g,\ord}}$ that belongs to 
 $\mathbb D^{-1}\mathbb I_{g,\ord}$ is contained in $\reallywidehat{\mathbb D^{-1} \mathbb I_g}$.
This implies that part (c) holds.
\endproof

\subsection{Flows}\label{S39}

In this subsection we show that the form $f^1_{g,\crys}$ induces $\delta$-flows in the sense recalled below on affine open subschemes of $\mathcal X_{\ord}$ whose reductions modulo $p$ are non-empty.

Let $X$ be an affine smooth scheme over $\Spec(R)$ and let $J^1(X)$ be its first $p$-jet space.
Following \cite{Bu17}, Ch. 3, Sect. 3.6, Def. 3.122, an ideal $\mathfrak a$ of $\mathcal O(J^1(X))$ will be called a {\it $\delta$-flow} on $X$ if the composite projection
$$Z^1(\mathfrak a):=\Spf(\mathcal O(J^1(X))/\mathfrak a) \subset J^1(X)\rightarrow \widehat{X}$$
is an isomorphism.

By the universal property 
of $p$-jet spaces, giving a $\delta$-flow on $X$
is equivalent to giving a $p$-derivation on $\mathcal O(\widehat{X})$.
(The notion of {\it $\delta$-flow} was suggested by the analogy between $p$-derivations and vector fields which was in the background of \cite{Bu95a,Bu05}.)

\medskip

\begin{rem}\label{R27}
If $\mathfrak a$ is a $\delta$-flow on $X$, then for all integers $n\geq 1$ the projection
$$Z^n(\mathfrak a):=\Spf(\mathcal O(J^n(X))/(\mathfrak a,\delta \mathfrak a,\ldots,\delta^{n-1}\mathfrak a)) \subset J^n(X)\rightarrow \widehat{X}$$
is an isomorphism. To check this we show that the $R$-algebra homomorphism 
$\mathcal O(\widehat{X})\rightarrow \mathcal O(Z^n(\mathfrak a))$ is an isomorphism. 
To check injectivity note that, by the universal property of $p$-jet spaces, the $p$-derivation on $\mathcal O(\widehat{X})$ defines a right inverse $\mathcal O(J^n(X))\rightarrow \mathcal O(\widehat{X})$ for
the $R$-algebra homomorphism $\mathcal O(\widehat{X}) \rightarrow \mathcal O(J^n(X))$; but clearly this right inverse factors through an $R$-algebra homomorphism
 $\mathcal O(Z^n(\mathfrak a))\rightarrow \mathcal O(\widehat{X})$ which is a right inverse for $\mathcal O(\widehat{X})\rightarrow \mathcal O(Z^n(\mathfrak a))$. 
 Hence $\mathcal O(\widehat{X})\rightarrow \mathcal O(Z^n(\mathfrak a))$ is injective.
We check that the $R$-algebra homomorphism $\mathcal O(\widehat{X})\rightarrow \mathcal O(Z^n(\mathfrak a))$ is surjective by induction on $n\in\mathbb N$. 
 The case $n=1$ is our $\delta$-flow hypothesis. Now assume the statement for some $n\geq 1$ and take an element $u$ of $\mathcal O(Z^{n+1}(\mathfrak a))$. So $u$ is a $p$-adic limit of elements of the form $v$ and $\delta w$ with $v,w\in\mathcal O(Z^n(\mathfrak a))$. By the induction hypothesis $v$ and $w$ come from elements of $\mathcal O(\widehat{X})$. Hence $\delta w$ comes from an element of $\mathcal O(Z^1(\mathfrak a))$. By case $n=1$, $\delta w$ comes then from an element of $\mathcal O(\widehat{X})$. Thus $u$ comes from an element of $\mathcal O(\widehat{X})$. This ends the induction.
\end{rem}

\medskip

Our next purpose is to show that Hecke covariant forms of weight $(\phi,1)$ naturally define $\delta$-flows on $\mathcal X_{\ord}\subset\mathcal X$. We first introduce a few extra concepts on $\delta$-flows.

Let $X$ be a smooth affine scheme over $\Spec(R)$ of relative dimension $d$.
 Assume $\bar X$ is irreducible.
Let $f=(f_1,\ldots,f_d)\in\mathcal O(J^1(X))^d$. 
Assume furthermore there exist \'etale coordinates $\tau=
(\tau_1,\ldots,\tau_d)\in\mathcal O(X)^d$, i.e., their corresponding morphism $\tau:X\rightarrow \mathbb A^d_R$ is \'etale. For a $d$-tuple of indeterminates $\tau':=(\tau'_1,\ldots,\tau'_d)$, let $\mathcal O(X)[\tau']:=\mathcal O(X)[\tau'_1,\ldots,\tau'_n]$ be the polynomial $\mathcal O(X)$-algebra in the entries of $\tau'$; we also view $\tau'$ as an element of $\mathcal O(X)[\tau']^d$.
 By \cite{Bu95a}, Prop. 1.4 we have a $\mathcal O(X)$-algebra isomorphism
 $$\reallywidehat{\mathcal O(X)[\tau']}\simeq  \mathcal O(J^1(X))$$
 that maps $\tau_i'$ to $\delta \tau_i$ for $i\in \{1,\ldots,d\}$.
 
We say $f$ is {\it quasi-linear with respect to $\tau$} if there exist
$$(L_{ij})_{1\leq i,j\leq d}\in\pmb{\GL}_d(\mathcal O(\widehat{X})), \ \ (E_1,\ldots,E_d)\in \mathcal O(\widehat{X})^d, \ \ (J_1,\ldots,J_d) \in\mathcal O(J^1(X))^d$$
 such that for all $i\in\{1,\ldots,d\}$, inside $\mathcal O(J^1(X))$ we have
 \begin{equation}\label{EQ141}
f_i=\sum_{j=1}^d L_{ij}\tau_j'+E_i+pJ_i.
\end{equation}

Next consider a point $P\in X(R)$, $P:\Spec(R)\rightarrow X$, and an isomorphism between the completion of $X$ along the image of $P$ (i.e., the completion of the local ring of $X$ at its $k$-valued closed point defined by $P$) and the formal power series ring
$R[[\underline{z}]]:=R[[z_1,\ldots,z_d]]$ 
 where $\underline{z}:=(z_1,\ldots,z_d)$ is a $d$-tuple of indeterminates to be called {\it formal coordinates} at $P$. For simplicity we view this isomorphism as an identification.
 We also assume that the image of $P$ under the  map $\tau(R):X(R)\rightarrow \mathbb A^d(R)$ is $0\in \mathbb A^d(R)$.
 Let $R[[\underline{z}]][\underline{z}']:=R[[\underline{z}]][z_1',\ldots,z_d']$ 
 where $\underline{z}':=(z'_1,\ldots,z'_d)$ is a $d$-tuple of new indeterminates.
 We have an $R$-algebra isomorphism (which we view also as an identification) $R[[\tau]]:=R[[\tau_1,\ldots,\tau_d]]\simeq  R[[\underline{z}]]$ induced by the composite $R$-algebra homomorphism $R[\tau_1,\ldots,\tau_d]\to\mathcal O(X)\to R[[\underline{z}]]$; 
 let $S_i\in R[[\tau]]$ be the image of $z_i$ for $i\in \{1,\ldots,d\}$.
 Let $R[[\tau]][\tau']:=R[[\tau]][\tau_1',\ldots,\tau_d']$ and let $\delta:R[[\tau]]\rightarrow \reallywidehat{R[[\tau]][\tau']}$ be the unique $p$-derivation relative to the inclusion such that $\delta \tau_i=\tau_i'$ for $i\in \{1,\ldots,d\}$.
 So we have a natural $R$-algebra composite homomorphism
 \begin{equation}
 \label{alghommm}
 \mathcal O(J^1(X))\rightarrow \reallywidehat{\mathcal O(X)[\tau']}\rightarrow 
\reallywidehat{R[[\tau]][\tau']}\rightarrow \reallywidehat{R[[\underline{z}]][\underline{z}']}\end{equation}
where the first and the third homomorphisms are isomorphisms, 
the inverse of the last isomorphism mapping $z'_i$ to $\delta S_i$ for $i\in \{1,\ldots,d\}$.
The homomorphism (\ref{alghommm}) is injective as $\overline{J^1(X)}$ is irreducible and smooth over $k$.

We say $f$ is {\it formally quasi-linear at $P$ with respect to $\underline{z}$} if there exist 
$$(u_{ij})_{a\leq i,j\leq d}\in\pmb{\GL}_d(R[[\underline{z}]]),\ \ (V_1,\ldots,V_d)\in R[[\underline{z}]]^d,\ \ (W_1,\ldots,W_d)\in\reallywidehat{R[[\underline{z}]][\underline{z}']}^d$$
 such that for all $i\in\{1,\ldots,d\}$ we have the following equality in the ring $\reallywidehat{R[[\underline{z}]][\underline{z}']}$:
$$f_i=\sum_{j=1}^d u_{ij}\tau_j'+V_i+pW_i.$$
Our terminology comes from the analogy with classical differential equations where
{\it quasi-linearity} means `linearity in the derivatives'.

Clearly, if $f$ is quasi-linear with respect to $\tau$, it is also formally quasi-linear at $P$ with respect
to $\underline{z}$.
We have the following  partial converse of this which is an extension of  \cite{BP}, Lem. 4.66:

\begin{lemma}\label{L28}
Let $X,P,f,\tau,\underline{z}$ be as above.
Assume $f$ is formally quasi-linear at $P$ with respect to $\underline{z}$. Then there exists an  affine open subscheme $Y$ of $X$ that contains $P$ and such that the image of $f$ in $\mathcal O(J^1(Y))^d$ is quasi-linear with respect to the image of $\tau$ in $\mathcal O(Y)^d$.
\end{lemma}

\noindent
{\it Proof.} Recall the injective homomorphism (\ref{alghommm}) which we view  as an inclusion. 
For $\alpha=(\alpha_1,\ldots,\alpha_d)\in \mathbb Z_{\geq 0}^d$ 
and $b=(b_1,\ldots,b_d)\in B^d$ with $B$ a ring,
we write 
$b^{\alpha}:=b_1^{\alpha_1}\cdots b_d^{\alpha_d}\in B$.
For $i\in\{1,\ldots,d\}$ we write
$$S_i:=S_i(\tau)=z_i=\sum_{\alpha\in\mathbb Z_{\geq 0}^d} c_{\alpha i}\tau^{\alpha}\in R[[\tau]],$$
where each $c_{\alpha i}\in R$. Let $\tau^p:=(\tau_1^p,\ldots,\tau_d^p)\in\mathcal O(X)^d$
and consider the $d$-tuple $\tau^p+p\tau'\in \mathcal O(X)[\tau']^d$.
Then
$$ z'_i=\frac{1}{p}\left[ \sum_{\alpha\in\mathbb Z_{\geq 0}^d} c_{\alpha i}^{\phi}(\tau^p+p\tau')^{\alpha}-\left( \sum_{\alpha\in\mathbb Z_{\geq 0}^d} c_{\alpha i} \tau^{\alpha} \right)^p \right]=\sum_{j=1}^d(\partial S_i/\partial \tau_j)^p \tau'_j+a_i+pb_i
$$
for some $a_i \in R[[\tau]]$ and $b_i \in\reallywidehat{R[[\tau]][\tau']}$. Denoting by $\bar{f}_i$ the reduction modulo $p$ of $f_i$, 
 we get the identity
 \begin{equation}
 \label{theidentityy}
 \overline{f}_i=\sum_{j=1}^d \overline{l}_{ij}\tau'_j+ \overline{a}_i\end{equation}
 in $k[[\tau]][\tau']$, where $\overline{a}_i\in k[[\tau]]$ is the reduction modulo $p$ of $a_i$ and, for $(u_{ij})_{1\leq i,j\leq d}$ and its reduction 
 $$(\overline{u}_{ij})_{1\leq i,j\leq d}\in\pmb{\GL}_d(k[[\underline{z}]])=\pmb{\GL}_d(k[[\tau]])$$modulo $p$,
 $$\overline{l}_{ij}:=\sum_{m=1}^d \overline{u}_{im}\cdot (\overline{\partial S_m/\partial \tau_j})^p\in k[[\tau]].$$
 Here $\overline{\partial S_m/\partial \tau_j}\in k[[\tau]]$ is the reduction modulo $p$ of $\partial S_m/\partial \tau_j\in R[[\tau]]$. 
 On the other hand we have
 $$\overline{f}_i\in\mathcal O(\overline{X})[\tau']\subset k[[\tau]][\tau'].$$
As $\overline{l}_{ij},\ \overline{a}_i\in k[[\tau]]$, from Equation (\ref{theidentityy})
 we get that 
 $$\overline{l}_{ij},\ \overline{a}_i\in\mathcal O(\overline{X}).$$
 Choose $L_{ij},E_i\in\mathcal O(X)$ that lift $\overline{l}_{ij},\overline{a_i}$ respectively.
 Then Equation (\ref{EQ141}) holds for some $J_i\in\mathcal O(J^1(X))$ as $\overline{\mathcal O(J^1(X))}\rightarrow k[[\tau]][\tau']$ is a $k$-algebra monomorphism. As the matrices 
 $$(u_{ij}(0))_{1\leq i,j\leq d},\ ((\partial S_i/\partial \tau_j)(0))_{1\leq i,j\leq d}\in\pmb{\Mat}_d(R)$$
 are invertible so is
 $(L_{ij}(P))_{1\leq i,j\leq d}\in\pmb{\Mat}_d(R)$. We conclude by taking $Y$ to be the non-zero locus of the determinant of the matrix $(L_{ij})_{1\leq i,j\leq d}\in\pmb{\Mat}_d(\mathcal O(X))$.\endproof

\begin{rem}\label{R28}
With $X,f,\tau$ as in the discussion preceding Lemma \ref{L28}, we consider additional \'etale coordinates $\tilde{\tau}=
(\tilde{\tau}_1,\ldots,\tilde{\tau}_d)\in\mathcal O(X)^d$. If $f$ is quasi-linear with respect to $\tilde{\tau}$, then it is also quasi-linear with respect to $\tau$. Indeed, let $P\in X(R)$.
Compositing $\tau$ and $\tilde{\tau}$ with translations of $\mathbb A^d$ we may assume 
that the maps $X(R)\rightarrow \mathbb A^d(R)=R^d$ induced by $\tau$ and $\tilde{\tau}$ map $P$ into $0$. Then we take $[\underline{z}]$ such that we have $R[[\underline{z}]]=R[[\tilde{\tau}]]:=R[[\tilde{\tau}_1,\ldots,\tilde{\tau}_d]]$. As $f$ is quasi-linear with respect to $\tilde{\tau}$ it is formally quasi-linear at $P$ with respect to $\underline{z}$ at $P$. Based on the proof of Lemma \ref{L28} we have $\bar Y=\bar X$, as the matrix $(\bar l_{ij})_{1\leq i,j\leq n}\in\pmb{\Mat}_d(\mathcal O(\bar X))$ is invertible.
\end{rem}

\medskip

The link between quasi-linearity and $\delta$-flows is given by the following lemma:

\begin{lemma}\label{L29}
Let $X,f,\tau$ be as in the discussion preceding Lemma \ref{L28}.
Assume $f$ is quasi-linear with respect to $\tau$.
Then the ideal
$\mathfrak a$ of $\mathcal O(J^1(X))$ generated by the entries of $f$ is
a $\delta$-flow on $X$. \end{lemma}

\noindent
{\it Proof.} Recall that $Z^1(\mathfrak a)=\Spf(\mathcal O(J^1(X)/\mathfrak a))$ and that $\overline{*}$ means the reduction modulo $p$ of $*$. It is enough to check that the projection
$$\overline{Z^1(\mathfrak a)}\subset \overline{J^1(X)}\rightarrow \overline{X}$$
is an isomorphism; but this is clear from the definition of quasi-linearity.
\endproof

\medskip

Let $(A_{\mathcal X},\theta_{\mathcal X})$ be the universal principally polarized abelian scheme over $\mathcal X=\mathcal A_{g,1,N,R}$. For a scheme $X$ over $\mathcal X$, let $(A_X,\theta_X):=X\times_{\mathcal X} (A_{\mathcal X},\theta_{\mathcal X})$. Until Subsection 3.10 we let $d:=\frac{g(g+1)}{2}$ be the relative dimension of $\mathcal X$ over $R$. 

For each ring $B$ and matrix $M=(m_{ij})_{1\leq i,j\leq g}\in\pmb{\Mat}_g(B)$, we define the $d$-tuple
$$\overrightarrow{M^{\textup{upp}}}:=(m_{11},m_{12},\ldots, m_{1g},m_{22},\ldots,m_{2g},m_{33},\ldots,m_{gg})\in B^{d}$$ 
whose entries $m_{ij}$ with $1\leq i\leq j\leq g$ are listed in the lexicographic order.

\begin{lemma}\label{L30}
There exists an affine open subscheme 
 $\mathcal Y \subset \mathcal X$ which admits \'etale coordinates $\tau$, for which $\overline{\mathcal Y}$ is non-empty, for which the abelian scheme $A_{\mathcal Y}$ has a column basis $\omega_{\mathcal Y}$ of $1$-forms on $A_{\mathcal Y}/{\mathcal Y}$ and such that for the matrix
 $$f_{\mathcal Y}:=f^1_{g,crys}(A_{\mathcal Y},\theta_{\mathcal Y},\omega_{\mathcal Y},\mathcal O^{\infty}(\mathcal Y))\in\pmb{\Mat}_g(\mathcal O(J^1({\mathcal Y})))$$
 the following three properties hold:
 
\medskip 
{\bf (a)} The $d$-tuple 
$\overrightarrow{f_{\mathcal Y}^{\textup{upp}}}\in\mathcal O(J^1(\mathcal Y))^d$
is quasi-linear with respect to $\tau$.

\smallskip
{\bf (b)} The ideal $\mathfrak a$ generated by the entries of the matrix $f_{\mathcal Y}$ is equal to the ideal $\mathfrak a^{\textup{upp}}$ generated by the entries of the $d$-tuple $\overrightarrow{f_{\mathcal Y}^{\textup{upp}}}$.\footnote{It is also equal to the ideal generated by the entries of the
symmetric matrix $f^{\langle 1\rangle}_{g,crys}(A_{\mathcal Y},\theta_{\mathcal
Y},\omega_{\mathcal Y})$, provided $\mathcal Y$ modulo $p$ is contained in the ordinary locus.}

\smallskip
{\bf (c)} The ideal $\mathfrak a$ is a $\delta$-flow on $\mathcal Y$. 
\end{lemma}

\noindent
{\it Proof.}
We recall that $(A_0,\theta_0,e)$ is as in Example \ref{EX3}. Equipping $(A_0,\theta_0)$ with a symplectic similitude level-$N$ structure, we get a $k$-valued point $P_0$ of $\mathcal X$. Let $P\in\mathcal X(R)$ be the $R$-valued point which is the canonical lift of $P_0$, i.e., is defined by the canonical lift $A_P$ of $A_0$.

We take $X\subset \mathcal X$ to be an affine open subscheme that admits \'etale coordinates $\tau$ and contains $P$. We denote  $T_{11},T_{12},\ldots,T_{1g},T_{22},\ldots, T_{2g},T_{33},\ldots,T_{gg}$ by 
$z_1,\ldots,z_d$
and for $1\leq j <i\leq g$ we will use as before $T_{ij}:=T_{ji}$; so $T^{\t}=T$ and $R[[\underline{z}]]=R[[T]]=S^0_{\for}$ is as before. 
Let $\omega_{S^0_{\for}}$ be the column basis of $1$-forms on the corresponding abelian scheme $A_{S^0_{\for}}/\Spec(S^0_{\for})$ that is attached to $e$. The restriction $\omega_P$ of $\omega_{S^0_{\for}}$ to $A_P$
is obtained from $\omega_{S^0_{\for}}$ by setting $T=0$. Replacing $X$ by a smaller open subscheme we can assume there exists a column basis $\omega_X$ of $1$-forms on $A_X/X$ whose restriction to $A_P$
equals $\omega_P$. So we have
$$\omega_X=\lambda\cdot \omega_{S^0_{\for}},\ \ \lambda=\lambda(T)\in\pmb{\GL}_g(S^0_{\for}),\ \ \lambda(0)=1_g.$$
Write $\lambda=(\lambda_{ij})_{1\leq i,j\leq g}$, $\lambda_{ij}=\lambda_{ij}(T)$.
Replacing $f$ by a $R^{\times}$-multiple of it and using the homogeneity property of $f=f^1_{g,crys}$ and Theorem \ref{T13} (a), we get that
$$f_X:=f(A_X,\theta_X,\omega_X)= \lambda\cdot f(A_{S^0_{\for}},\theta_{S^0_{\for}},\omega_{S^0_{\for}})
\cdot (\lambda^{\t})^{\phi}=\lambda \cdot \Psi \cdot (\lambda^{\t})^{\phi}.$$
We have 
$$\overrightarrow{f_X^{\textup{upp}}}=(f_{X,11},f_{X,12},\ldots, f_{X,1g},f_{X,22},\ldots,f_{X,2g},f_{X,33},\ldots,f_{X,gg}),$$
 where for $1\leq i,j\leq g$ we have (recall that $\Psi$ is symmetric)
$$\begin{array}{rcl}
f_{X,ij} & = & \sum_{m,s=1}^g \lambda_{im}\Psi_{ms}\lambda_{sj}^{\phi}\\
\ & \ & \ \\
\ & = & \sum_{1\leq m<s\leq g} (\lambda_{im}\lambda_{sj}^{\phi}+\lambda_{is}\lambda_{mj}^{\phi})\Psi_{ms}+\sum_{s=1}^g \lambda_{is}\lambda_{sj}^{\phi}\Psi_{ss}.
\end{array}
$$
The $d$-tuple $\overrightarrow{\Psi^{\textup{upp}}}$ is easily seen to be formally quasi-linear at $P$ with respect to $\underline{z}$. 
Recall the ring
 $S^1_{\for}:=\reallywidehat{S^0_{\for}[T']}$ in Subsubsection \ref{S235} and write
$$\overrightarrow{f^{\textup{upp}}_X}=Q\cdot \overrightarrow{\Psi^{\textup{upp}}},\ \ Q=(q_{ij,ms})_{1\leq i\leq j\leq g,1\leq m\leq s\leq g}\in\pmb{\Mat}_d(S^1_{\for}),$$
where a $d$-tuple is identified with a matrix of size $d\times 1$ and
$$\begin{array}{rcl}
q_{ij,ms} & := & \lambda_{im}\lambda_{sj}^{\phi}+\lambda_{is}\lambda_{mj}^{\phi}\ \ \textup{for}\ \ m<s,\\
\ & \ & \ \\
q_{ij,ss} & := & 
\lambda_{is}\lambda_{sj}^{\phi}.\end{array}
$$

\noindent
\begin{claim}\label{Claim1}
We have $Q=Q_0+pQ_1$, $Q_0\in\pmb{\GL}_d(S^0_{\for})$ and $Q_1\in\pmb{\Mat}_d(S^1_{\for})$. 
\end{claim}

Indeed the reduction $\overline{q}_{ij,ms}\in\overline{S^1_{\for}}=k[[T]][T']$ of $q_{ij,ms}$ modulo $p$ is actually in $\overline{S^0_{\for}}=k[[T]]$, so we can write $\overline{q}_{ij,ms}=\overline{q}_{ij,ms}(T)$ and, using the Kronecker delta,  we have
$$
\overline{q}_{ij,ms}(0)= \delta_{im}\delta_{js}+\delta_{is}\delta_{jm}\ \ \ \for\ \ m<s,
$$
and
$$\overline{q}_{ij,ss}(0)=\delta_{is}\delta_{js}.$$
So 
$\overline{q}_{ij,ms}(0)=0$ for $(i,j)\neq (m,s)$ and $\overline{q}_{ms,ms}=1$; in other words, the reduction modulo $p$ of the matrix $Q$ has entries in $\overline{S^0_{\for}}$ and its evaluation at $0$ is the identity matrix. This proves Claim \ref{Claim1}.

By Claim \ref{Claim1} and the fact that $\overrightarrow{\Psi^{\textup{upp}}}$ is formally quasi-linear at $P$ with respect to $\underline{z}$ it follows that
$\overrightarrow{f^{\textup{upp}}_X}$ is formally quasi-linear at $P$ with respect to $\underline{z}$. By Lemma \ref{L28} there exists an affine open subscheme $\mathcal Y\subset X$ as in the part (a) of the lemma, and we let $f_{\mathcal Y}$ and $\overrightarrow{f_{\mathcal Y}^{\textup{upp}}}$ be the analogues of $f_X$ and $\overrightarrow{f_X^{\textup{upp}}}$ (respectively).

We prove part (b) of the lemma. Clearly $\mathfrak a^{\textup{upp}}\subset \mathfrak a$. Thus we have closed embeddings
$$Z^1(\mathfrak a)=\Spf(\mathcal O(J^1(\mathcal Y))/\mathfrak a)\subset Z^1(\mathfrak a^{\textup{upp}})=\Spf(\mathcal O(J^1(\mathcal Y))/\mathfrak a^{\textup{upp}})\subset J^1(\mathcal Y)$$
between formal schemes.
By part (a), the 
projection
$$\overline{Z^1(\mathfrak a^{\textup{upp}})}\subset \overline{J^1(\mathcal Y)}\rightarrow \overline{\mathcal Y}$$
is an isomorphism. To prove that $\mathfrak a=\mathfrak a^{\textup{upp}}$ it is enough to show that
the closed embedding $Z^1(\mathfrak a)\subset Z^1(\mathfrak a^{\textup{upp}})$ is an isomorphism, hence that the closed embedding
$\overline{Z^1(\mathfrak a)}\subset\overline{Z^1(\mathfrak a^{\textup{upp}})}$ is an isomorphism. As $\overline{\mathcal Y}$ is reduced it is enough to show that the image of the projection
\begin{equation}\label{EQ142}
\overline{Z^1(\mathfrak a)}\subset \overline{Z^1(\mathfrak a^{\textup{upp}})}\rightarrow \overline{\mathcal Y}\end{equation}
contains a Zariski open sense subset. Now the ordinary locus $\overline{\mathcal Y}_{\ord}$
of $\overline{\mathcal Y}$ is open and dense. Let $P_0'\in\overline{\mathcal Y}_{\ord}(k)$ be a $k$-valued point
and let $P'\in \mathcal Y(R)$ be the $R$-valued point corresponding to the canonical lift of the abelian variety over $\Spec(k)$ that corresponds to $P_0'$. Let $J^1(P')\in J^1(\mathcal Y)(R)$ be the first jet of $P'$. Then all the components $f_{\mathcal Y,ij}$ of $\overrightarrow{f_{\mathcal Y}^{\textup{upp}}}$ vanish at $J^1(P')$ hence $\overline{J^1(P')}\in\overline{Z^1(\mathfrak a)}(k)$,
and thus $\overline{J^1(P')}$ maps to $P_0'$. This shows that the image of the projection (\ref{EQ142})
contains a dense open subscheme which ends our proof of part (b).

Part (c) follows from parts (a) and (b) plus Lemma \ref{L29}.
\endproof

\begin{lemma}\label{L31}
Let $\pi:X'\rightarrow X$ be an \'etale morphism between affine smooth schemes over $\Spec(R)$ and let $\mathfrak a \subset \mathcal O(J^1(X))$. Let $\mathfrak a':=\mathfrak a \cdot \mathcal O(J^1(X'))$. 
If $\mathfrak a$ is a $\delta$-flow on $X$ then 
$\mathfrak a'$ is a $\delta$-flow on $X'$. The converse also holds provided $\pi$ is surjective.
\end{lemma}

\noindent
{\it Proof.}
As $\pi$ is \'etale, based on \cite{Bu95a}, Prop. 1.4,
we have an isomorphism $J^1(X')\simeq J^1(X)\times_{\widehat{X}} \widehat{X'}$. Hence $Z^1(\mathfrak a')\simeq Z^1(\mathfrak a)\times_{\widehat{X}}\widehat{X'}$. Clearly, if $Z^1(\mathfrak a)\rightarrow \widehat{X}$ is an isomorphism, then
$Z^1(\mathfrak a')\rightarrow \widehat{X'}$ is an isomorphism; the converse also holds if $\pi$ is surjective because then $\pi$ is faithfully flat. \endproof

\begin{cor}\label{C18.3}
Let $X_1,X_2,Y$ be affine smooth schemes over $\Spec(R)$ and let
$$X_1\stackrel{\pi_1}{\longleftarrow} Y \stackrel{\pi_2}{\longrightarrow} X_2$$ be \'etale morphism with $\pi_1$ surjective.
Let $m\geq 1$ be an integer and let $g^{(1)}\in\mathcal O(J^1(X_1))^m$ and 
$g^{(2)}\in\mathcal O(J^1(X_2))^m$ be $m$-tuples of elements (viewed as $m\times 1$ matrices) such that we have 
$$\pi_1^*g^{(1)}=M_{12}\cdot \pi^*_2g^{(2)}$$
for some matrix $M_{12}\in\pmb{\GL}_m(\mathcal O(J^1(Y)))$.
Let $\mathfrak a_1$ and $\mathfrak a_2$ be the ideals in $\mathcal O(J^1(X_1))$ and $\mathcal O(J^1(X_2))$ generated by the entries of $g^{(1)}$ and $g^{(2)}$ (respectively).
If $\mathfrak a_2$ is a $\delta$-flow on $X_2$, then $\mathfrak a_1$ is a $\delta$-flow on $X_1$.\end{cor}

\noindent
{\it Proof.} This follows directly from Lemma \ref{L31} applied to $\pi_2$ and $\pi_1$.\endproof

\medskip
 
If $X$ is an affine open subscheme of $\mathcal X$  such that there exists a column basis $\omega_X$ of $1$-forms on $A_X/X$, then the ideal $\mathfrak a_X$ of $\mathcal O(J^1(X))$ generated by the entries of 
$$f_X:=f^1_{g,\crys}(A_X,\theta_X,\omega_X,\mathcal O^{\infty}(X))\in\pmb{\Mat}_g(\mathcal O(J^1(X)))$$
does not depend on the choice of the basis $\omega_X$. 
 
We now take $X$ to be an arbitrary affine open subscheme of $\mathcal X$. For each affine
 cover $X=\cup X_{\alpha}$ with every $X_{\alpha}$ such that there exists a column basis $\omega_{X_{\alpha}}$ of $1$-forms
 on $A_{X_{\alpha}}/X_{\alpha}$, the ideals $\mathfrak a_{X_{\alpha}}$ glue together to give a well-defined ideal $\mathfrak a_X$ in $\mathcal O(J^1(X))$ that does not depend on the affine cover.
 
 \begin{thm}\label{T20}
 Let $X$ be an affine open subscheme of $\mathcal X$ whose reduction modulo $p$ is non-empty and contained in the ordinary locus $\overline{\mathcal X}_{\ord}$ of $\overline{\mathcal X}$. Then the ideal $\mathfrak a_X$ is a $\delta$-flow on $X$.
 \end{thm}

\noindent
{\it Proof.}
Let $\mathcal Y\subset X$ and $\omega_{\mathcal Y}$ be as in Lemma \ref{L30} and its proof. It suffices to show that $X$ can be covered with affine open subschemes $X_{\alpha}$ such that $\mathfrak a_{X_{\alpha}}$ is a $\delta$-flow on $X_{\alpha}$. So we can assume that there exists a column basis $\omega_X$ of $1$-forms on $A_X/X$ and it suffices to show that each point $P_0\in X(k)$ is contained in an affine open subscheme $X_0$ of $X$ such that $\mathfrak a_{X_0}$ is a $\delta$-flow on $X_0$. 
 
Let $\ell$ be a prime distinct from $p$. By the density of the $\ell$-power polarized isogeny classes of ordinary points (see \cite{Cha}, Thm. 2) there exists an isogeny
$$u_0:A_{P_0}\rightarrow A_{P_0'},$$
 where $(A_{P_0},\theta_{P_0})$ is the principally polarized abelian variety over $\Spec(k)$ that corresponds to $P_0$ and where $(A_{P_0'},\theta_{P_0'})$ is the principally polarized abelian variety  over $\Spec(k)$ that corresponds to some point $P_0'\in \mathcal Y(k)$, such that 
 $u_0$ preserves the principal polarizations up to multiplication by $\ell^{s_0}$ for some 
 $s_0\in\mathbb N\cup\{0\}$ (i.e., the pullback of $\theta_{P_0'}$ by $u_0$ equals $\theta_{P_0}$ times $\ell^{s_0}$). Thus $\deg(u_0)=\ell^{2gs_0}$ is prime to $p$.

Let 
$$\mathcal X \stackrel{\pi_1}{\leftarrow} \mathcal X' \stackrel{\pi_2}{\rightarrow} \mathcal X$$
be the Hecke correspondence that parameterizes pairs of principally polarized abelian schemes of relative dimension $g$ 
with symplectic similitude level-$N$ structures together with an isogeny 
of degree $\ell^{2gs_0}$ which preserves the principal polarizations up to multiplication by $\ell^{s_0}$ and let 
$$u_{\mathcal X'}:\pi_1^*A_{\mathcal X}\rightarrow \pi_2^*A_{\mathcal X}$$
be the corresponding universal isogeny of abelian schemes over $\mathcal X'$ (compatible with the principal polarizations). 

It is well-known that $\mathcal X'$ is smooth and that the two projections $\pi_1,\pi_2$ are finite \'etale. We recall the argument for this. We fix a perfect alternating form $\psi_{2g}:\mathbb Z_{\ell}^{2g}\times \mathbb Z_{\ell}^{2g}\rightarrow \mathbb Z_{\ell}$. Working in the finite \'etale topology of $\mathcal X$, we can assume that $\ell^{s_0}$ divides $N$. But in this case $\mathcal X'$ is the disjoint union of copies of $\mathcal X$ indexed by $\mathbb Z_{\ell}$-lattices $\mathcal L$ of $\mathbb Q_{\ell}^{2g}$ which satisfy the following property:

\medskip\noindent
{\bf (i)} {\it We have $\ell^{s_0}\mathbb Z_{\ell}^{2g}\subset\mathcal L\subset \ell^{-s_0}\mathbb Z_{\ell}^{2g}$, $s_0$ is the smallest whole number with this property and we get a perfect alternating form $\psi_{2g}:\mathcal L\times\mathcal L\rightarrow \mathbb Z_{\ell}.$}

\medskip

\noindent
Moreover, the two projections $\pi_1$ and $\pi_2$ induce isomorphisms when restricted to each $\mathcal X$ copy of $\mathcal X'$. 

We view $u_0:A_{P_0}\rightarrow A_{P_0'}$ as a $k$-valued point of $\mathcal X'$. Setting
$$X':=\pi_1^{-1}(X)\cap \pi_2^{-1}(\mathcal Y) \subset \mathcal X'$$
we have that $u_0$ belongs to $X'(k)$. 
As $\pi_1$ is flat, $\pi_1(X')$ is open in $X$. Hence there exists an affine open subscheme $X_0$ of $\pi_1(X')$ that contains $P_0$.
 Setting $X_0':=\pi_1^{-1}(X_0)$ we have a 
diagram
$$X_0\stackrel{\Pi_1}{\leftarrow} X_0'\stackrel{\Pi_2}{\rightarrow} \mathcal Y$$
with $\Pi_1,\Pi_2$ induced by $\pi_1,\pi_2$, respectively, 
and an induced
isogeny
$$u_{X_0'}:\Pi_1^*A_{X_0}\rightarrow \Pi_2^*A_{\mathcal Y}$$
of abelian schemes over $X_0'$.
Consider the matrix
$$f_{X_0}:=f^1_{g,crys}(A_{X_0},\theta_{X_0},\omega_{X_0},\mathcal O^{\infty}(X_0))\in\pmb{\Mat}_g(\mathcal O(J^1(X_0)))$$
where $\theta_{X_0}$ is induced by $\theta_X$.
By Hecke covariance we have
$$\Pi_1^*f_{X_0}=\beta_1 \cdot \Pi^*_2f_{\mathcal Y} \cdot \beta_2$$
for some $\beta_1,\beta_2\in\pmb{\GL}_g(\mathcal O(J^1(X_0')))$.
We apply Corollary \ref{C18.3} by taking the quadruple  $(X_1,X_2,Y,m)$ to be equal to $(X_0,\mathcal Y,X_0',g^2)$ and with $$Q_{12}\in\pmb{\GL}_{g^2}(\mathcal O(J^1(X_0')))$$ as the matrix representation (with respect to some standard bases) of the $\mathcal O(J^1(X_0'))$-linear 
automorphism defined by the rule $\gamma\mapsto \beta_1\cdot \gamma \cdot \beta_2$, $\gamma\in\pmb{\Mat}_g(\mathcal O(J^1(X_0')))$; as the ideal $\mathfrak a_{\mathcal Y}$ generated by the entries of $f_{\mathcal Y}$ is a $\delta$-flow on $\mathcal Y$ (see Lemma \ref{L30} (c)), it follows that the ideal $\mathfrak a_{X_0}$ generated by the entries of $f_{X_0}$ is as well a $\delta$-flow on $X_0$.\endproof

\subsection{Reciprocity for canonical lifts}\label{S310.5}

The purpose of this subsection is to derive an application of our flow theorem (see Theorem \ref{T20})
which extends to arbitrary genus $g$ the `Reciprocity Theorem for CL points' in \cite{BP}, Thm. 3.5. 

\begin{df}\label{df27}
Consider a finite normal extension $R'$ of $\mathbb Z_p$ contained in $R$ and $A'$ an abelian scheme over $R'$. 
Let $p^{\nu}$ be the cardinality of the field $R'/pR'$ and 
let $A$ be the abelian scheme over $R$ obtained from $A'$ via base change.
We say that $A'$ is  {\it non-degenerate} if the following two conditions hold:

\smallskip

{\bf (NAX1)} The characteristic polynomial of the $p^{\nu}$-power Frobenius endomorphism
of $\overline{A'}$ has simple roots only.

\smallskip

{\bf (NAX2)} The Serre--Tate matrix $q(A):=(q_A(e_i,\check{e}_j))_{1\leq i,j\leq d}$ of $A$ with respect to some (equivalently every) bases $e_1,\ldots, e_d$ and $\check{e}_1,\ldots, \check{e}_d$ of $T_p(A_0)$ and $T_p(\check{A}_0)$ respectively satisfies the condition
$$\det\left(p^{-1}(q(A)-\mathbb E_d)\right) \not\equiv 0\ \ \ \text{mod}\ \ \ p$$
(recall that $\mathbb E_d$ is the $d\times d$ matrix all of whose entries are $1$).
\end{df}

For $A$ and $A'$ as above we denote by $A(\mathbb Z_p^{\textup{ur}})$ the set $A'(\mathbb Z_p^{\textup{ur}})$ of $\mathbb Z_p^{\textup{ur}}$-points of $A'$ where $\mathbb Z_p^{\textup{ur}}$ is the maximal unramified extension extension of $\mathbb Z_p$ in 
$R$; this makes sense as $R'\subset \mathbb Z_p^{\textup{ur}}$. We identify $R=\widehat{\mathbb Z_p^{\textup{ur}}}$. 
Note that the set $A'(\mathbb Z_p^{\textup{ur}})$ is independent of the choice of the model $A'/R'$ of $A/R$. This is so as the pullback functor from the category of abelian schemes over $\Spec(\mathbb Z_p^{\textup{ur}})$ to the category of abelian schemes over $\Spec(R)$ is full faithful.

Recall that $\mathcal X_{\ord}$ denotes the ordinary locus of $\mathcal X:=\mathcal A_{g,1,N,R}$ with $N\geq 3$. Let $\CL_R\subset \mathcal X_{\ord}(R)$ be the set of canonical lift points (see Remark \ref{R12}). Also, for an abelian group $\Gamma$ we denote by $\Gamma_{\textup{tors}}$ its torsion subgroup.

The main result of this subsection is the following:

\begin{thm}\label{T20.5}
Assume we are given morphisms of schemes over $R$
$$\Pi:Y\rightarrow X,\ \ \ \ \Phi:Y\rightarrow A$$ 
with $X\subset \mathcal X_{\ord}$ affine, $Y$ affine, $\Pi$ \'{e}tale, and $A$ an abelian scheme  of relative dimension $d$. Assume that $X,Y,A,\Pi,\Phi$ descend to corresponding objects 
$X',Y',A',\Pi',\Phi'$ defined over a  finite normal extension $R'$ of $\mathbb Z_p$ contained in $R$ and assume $A'$ is non-degenerate. 
Then there exists a morphism of $p$-adic formal schemes 
$$\Phi^+:\widehat{Y}\rightarrow \widehat{\mathbb G_a^d}$$
such that for all integers $n\geq 1$, all points $P_1,\ldots,P_n\in \Pi^{-1}(\CL_R)\subset Y(R)$ and all integers $m_1,\ldots,m_n\in \mathbb Z$, the following equivalence holds
$$
\sum_{i=1}^n m_i \Phi(P_i)\in A(R)_{\textup{tors}} \ \ \ \ \Longleftrightarrow \ \ \ \ 
\sum_{i=1}^n m_i \Phi^+(P_i)=0\in R^d.
$$
\end{thm}

Above the maps $\Pi(R):Y(R)\rightarrow X(R)$,
$\Phi(R):Y(R)\rightarrow A(R)$, and $\Phi^+(R):Y(R)=\widehat{Y}(R)\rightarrow \widehat{\mathbb G_a^d}(R)=R^d$
were denoted simply by $\Pi, \Phi$, and $\Phi^+$ (respectively); here $\widehat{Y}(R)$ denotes, as usual, the set of $\textup{Spf}(R)$-points of $\widehat{Y}$. 
A map such as $\Phi^+$ above was called in \cite{BP} a {\it reciprocity map}; we refer to \cite{BP}, Rm. 3.7
for comments on this terminology. In particular we get:

\begin{cor}\label{C18.5}
In the situation  of Theorem \ref{T20.5} assume in addition that there exists $Q\in \Pi^{-1}(\CL_R)$ such that
$\Phi(Q)\not\in A(R)_{\textup{tors}}$. Then the image of 
the map
\begin{equation}
\label{EQ142.2}
\Pi^{-1}(\CL_R)\cap \Phi^{-1}(A(R)_{\textup{tors}})\rightarrow Y(R)\rightarrow Y(k)
\end{equation}
is not Zariski dense in $\overline{Y}$. (Here the first arrow in (\ref{EQ142.2}) is the inclusion and 
$Y(k)=\overline{Y}(k)$.)
\end{cor}

\noindent
{\it Proof.}  Taking $n=1$ and $m_1=1$ in Theorem \ref{T20.5} we get that for all $P\in  \Pi^{-1}(\CL_R)\subset Y(R)$ we have 
\begin{equation}
\label{extraaa}
\Phi(P)\in A(R)_{\textup{tors}}\ \ \Longleftrightarrow  \ \Phi^+(P)=0.\end{equation}
By the existence of $Q$ there exists $j\in \{1,\ldots,d\}$ such that the $j$th component $\Phi^+_j\in \mathcal O(\widehat{Y})$ of $\Phi^+$ is non-zero. Write $\Phi^+_j=p^N \Phi^+_{j,0}$ with $N\geq 0$ and $\Phi^+_{j,0}\in \mathcal O(\widehat{Y})\backslash p\mathcal O(\widehat{Y})$. Hence
the  image $\overline{\Phi^+_{j,0}}$ of $\Phi^+_{j,0}$ 
in $\mathcal O(\overline{Y})$ is non-zero. By the equivalence in Equation (\ref{extraaa}) the image of the map (\ref{EQ142.2}) is contained in the zero locus of $\overline{\Phi^+_{j,0}}\in \mathcal O(\overline{Y})$.
\endproof

\medskip

In order to prove  Theorem \ref{T20.5} recall from \cite{Bu95a}, Introd. that a {\it $\delta$-character} of order $\leq r$ of an abelian scheme $A$ over $R$ is a $\delta$-function $\psi:A(R)\rightarrow R=\mathbb G_a(R)$ of order $\leq r$ in the sense of Definition \ref{df7} which is also a group homomorphism. We need the following result from \cite{Bu97}, Cor. 1.7:

\begin{thm}\label{T20.7}
Consider a  finite normal extension $R'$ of $\mathbb Z_p$ contained in $R$, let 
 $A'$ be a non-degenerate abelian scheme over $R'$ of relative dimension $d$, and let $A$ be the abelian scheme over $R$ obtained via base change from $A'$. Then there exist $\delta$-characters $\psi_1,\ldots,\psi_d:A(R)\rightarrow R$ of order $\leq 2$ such that
the map
$$\psi:=(\psi_1,\ldots,\psi_d):A(R)\rightarrow R^d=\mathbb G_a^d(R)$$
is surjective and satisfies
$$A(\mathbb Z_p^{\textup{ur}}) \cap \textup{Ker}(\psi)=A(\mathbb Z_p^{\textup{ur}})_{\textup{tors}}.$$
\end{thm}

The theorem above can be viewed as an arithmetic analogue of Manin's Theorem of the Kernel \cite{Man}, Sect. 5, Thm. 2.

\medskip

\noindent {\it Proof of Theorem \ref{T20.5}}.
Let $\psi_1,\ldots,\psi_d$ be as in Theorem \ref{T20.7}. For $j\in \{1,\ldots,d\}$ we 
consider the $\delta$-function $\psi_j\circ \Phi$ on $Y(R)$, identified with an element of 
$\mathcal O(J^2(Y))$. On the other hand,
by Theorem \ref{T20} the ideal $\mathfrak a_X\subset \mathcal O(J^1(X))$ is a $\delta$-flow on $X$. By Lemma \ref{L31} the ideal $\mathfrak a_X \cdot \mathcal O(J^1(Y))$ 
(generated by the image of $\mathfrak a_X$ via the homomorphism $\Pi^*:\mathcal O(J^1(X))
\rightarrow \mathcal O(J^1(Y))$ induced by $\Pi$)
is a $\delta$-flow on $Y$. By Remark \ref{R27} the algebra map
$$\mathcal O(\widehat{Y})\rightarrow \mathcal O(J^2(Y))/(\mathfrak a_X,\delta \mathfrak a_X)$$
is an isomorphism which, for simplicity we view as an equality. So for each $j\in \{1,\ldots,d\}$ we have an equality of functions $Y(R)\rightarrow R$,
\begin{equation}
\label{mixedemo}
\psi_j\circ \Phi=\Phi^+_j+\sum_{m=1}^M a_{j,m}\cdot (\varphi_{j,m}\circ \Pi)+\sum_{m=1}^M b_{j,m} \cdot (\delta \circ \varphi_{j,m}\circ \Pi),\end{equation}
with $\Phi^+_j\in \mathcal O(\widehat{Y})$, $M\in \mathbb N$, where for $m\in \{1,\ldots,M\}$ $a_{j,m}, b_{j,m} \in \mathcal O(J^2(Y))$ are identified with functions $Y(R)\rightarrow R$ and $\varphi_{j,m}\in \mathfrak a_X$ is identified with a function $X(R)\rightarrow R$. These identifications are obtained by viewing the injective $R$-algebra homomorphism in Equation (\ref{EQ003}) as an inclusion. Also the dots in Equation (\ref{mixedemo}) denote usual multiplication of $R$-valued functions.
Consider the  morphism of $p$-adic formal schemes
$$\Phi^+:=(\Phi_1^+,\ldots,\Phi_d^+):\widehat{Y}\rightarrow \widehat{\mathbb G_a^d}.$$
The functions $\varphi_{j,m}\circ \Pi$ vanish on every point $P\in \Pi^{-1}(\CL_R)$. So for each such $P$ we have 
\begin{equation}
\label{EQ142.5}
\psi_j(\Phi(P))=\Phi_j^+(P).\end{equation}
Assume now that $P_1,\ldots,P_n\in \Pi^{-1}(\CL_R)$ and $m_1,\ldots,m_n\in \mathbb Z$.

\begin{claim}\label{Claim2} 
We have $P_1,\ldots,P_n\in Y'(\mathbb Z_p^{\textup{ur}})$. \end{claim}

Indeed as $\Pi$ is \'{e}tale it is sufficient to show that for each $j\in\{1,\ldots,n\}$, the image $\Pi(P_j)\in X(R)$ in fact belongs to $X'(\mathbb Z_p^{\textup{ur}})$. Recall now the Shimura--Taniyama theorem that abelian varieties with $\CM$ over an algebraically closed  field of characteristic $0$ are defined over number fields (for instance, see \cite{CCO}, Thm. 1.7.2.1). Note that all polarizations and symplectic similitude level-$N$ structures on such
abelian varieties are then also defined over  (possibly larger) number fields. We deduce that $\Pi(P_j)\in X'(F)$
 for some finite extension $F$ of the fraction field of $R'$ contained in an algebraic closure of $K$. As $F\cap R\subset \mathbb Z_p^{\textup{ur}}$ our claim follows.
 
By Claim \ref{Claim2}, we have
$\Phi(P_1),\ldots,\Phi(P_n)\in A(\mathbb Z_p^{\textup{ur}})$. By Theorem \ref{T20.7}, and by the fact that $\psi_j$ are group homomorphisms,
the condition 
$$\sum_{i=1}^n m_i \Phi(P_i)\in A(R)_{\textup{tors}}$$
is equivalent to the condition 
\begin{equation}
\label{EQ142.7}
\sum_{i=1}^n m_i \psi_j(\Phi(P_i))=0,\ \ \ \text{for}\ \ \ j\in \{1,\ldots,d\}.\end{equation}
By Equation (\ref{EQ142.5}) we get that the condition (\ref{EQ142.7}) is equivalent to the condition
$$\sum_{i=1}^n m_i \Phi^+(P_i)=0.$$
\endproof

\subsection{Some computations of cohomology groups}\label{S310.66}

In this subsection we present a few basic computations of cohomology 
groups of some simple ssa-$\delta$-spaces that appeared naturally in our paper.
We will not pursue here more general
cases.

\begin{prop}\label{P8}
For $\mathbb P^1_{\delta}:=\Proj_{\delta}(\mathbb S_1, (z_0,z_1))$ we have
$$H^0(\mathbb P^1_{\delta})=R\ \ \ \textup{and}\ \ \ \phirank_R(H^1(\mathbb P^1_{\delta}))=\infty.$$
\end{prop}

In what follows we recall the forms $f^{\langle r \rangle}$ in Equation (\ref{EQ057}).

\begin{prop}\label{P9}
For $Z_{\delta,n}:=\Proj_{\delta}(\mathbb I_{1,\ord},(f^{\langle 1 \rangle},\ldots,f^{\langle n \rangle}))$, $n\geq 2$ we have
 $$H^0(Z_{\delta,n})=\mathbb I_{1,\ord,\langle f^{\langle n \rangle}\rangle}\ \ \ \textup{and} \ \ 
 \ H^d(Z_{\delta,n})=0\ \ \ \textup{for}\ \ d\geq 1.$$
\end{prop}

\begin{prop}\label{P10}
For $\lambda=(\lambda_{ij})_{1\leq i,j\leq 2}\in \pmb{\GL}_2(R)$
 set $\begin{bmatrix}f_1^{\lambda} \\
 f_2^{\lambda} \\
 \end{bmatrix}:=\lambda\cdot \begin{bmatrix} f^{\langle 1 \rangle}\\ 
 f^{\langle 2 \rangle}\\
 \end{bmatrix}$ and
 $$Z_{\delta,2}^{\lambda}:=\Proj_{\delta}(\mathbb I_{1,\ord},(f^{\lambda}_1,f^{\lambda}_2)).$$
If all entries $\lambda_{ij}$ of $\lambda$ are non-zero, then we have
 $$H^0(Z_{\delta,2}^{\lambda})=R\ \ \ \textup{and}\ \ \ \phirank_R(H^1(Z_{\delta,2}^{\lambda}))=\infty.$$
\end{prop}

To state a partial generalization of Proposition \ref{P10} we recall the maps $\sigma,\tau:\mathbb H_g(s)\rightarrow \mathbb H_g(s)$ in Equation (\ref{EQ097}) and $\diamondsuit:\mathbb H_g(s)\rightarrow \mathbb I_{g,\ord}(-2s)$ in Proposition \ref{P7} (a). 

\begin{df}\label{df28}
The elements of the sets
$$\mathbb H_g(1)\setminus\sigma(\mathbb H_g(1))\ \ \textup{and}\ \ 
\diamondsuit(\mathbb H_g(1)\setminus\sigma(\mathbb H_g(1)))$$ will be
called primitive.\end{df}

\begin{prop}\label{P11}
Assume $\f=(f_0,\ldots,f_n)$ is a linear system of weight $-2$ in $\mathbb I_{g,\ord}$ such that $f_0$ and $f_1$ are primitive and $R$-linearly independent. Then we have
$$H^0(\Proj_{\delta}(\mathbb I_{g,\ord},\f))=R.$$
\end{prop}

For the  non-ordinary case we will prove:

\begin{prop}\label{P12}
There exists a linear system $(\varphi_1,\varphi_2)$ of weight $-1-\phi-\phi^2-\phi^3$ in $\mathbb I_g$ such that
$$H^0(\Proj_{\delta}(\mathbb I_g,(\varphi_1,\varphi_2)))=R.$$
\end{prop}
 
Proposition \ref{P9} follows from the following lemma applied with $B$ equal to $\mathbb I_{1,\ord}$ (see Theorem \ref{T4} (b) and Equation (\ref{EQ057}) for the units part):
 
\begin{lemma}\label{L32}
Let $B$ be an integral $\delta$-graded ring. Let 
$f_0,\ldots,f_n\in B(w)\setminus pB(w)$, $w\in W$.
 If there exist units $g_1,\ldots,g_n\in B^\times$ such that we have $f_i=g_i\dot (f_0)^{\phi^{i}}$ for all $i\in\{1,\ldots, n\}$, then for 
 $Z_{\delta}:=\Proj_{\delta}(B,(f_0,\ldots,f_n))$ we have
 $$H^0(Z_{\delta})=B_{\langle f_n\rangle}\ \ \textup{and}\ \  H^d(Z_{\delta})=0\ \ \textup{for}\ \ d\geq 1.$$
 \end{lemma}
 
 \noindent
 {\it Proof.}
 The hypotheses 
 imply that for all integers $d\geq 1$ and $0\leq i_0<\cdots<i_d\leq n$ we have
 $$B_{\langle f_{i_0}\cdots f_{i_d}\rangle}=B_{\langle f_{i_0}\rangle}\supset \cdots \supset B_{\langle f_{i_d}\rangle}.$$
 The $H^0$ statement follows.
 Now let $(a_{i_0\cdots i_d})_{0\leq i_0<\cdots <i_d\leq n}$ be a cocycle in 
 $$B_d:=\prod_{0\leq i_0<\cdots< i_{\bullet}\leq d}
 B_{\langle f_{i_0}\cdots f_{i_{\bullet}}\rangle}.$$
  For integers $0\leq i_0<\cdots<i_{d-1}\leq n$ define 
 $$b_{i_0\cdots i_{d-1}}:=a_{i_0\cdots i_{d-1} n}\ \ \ \textup{if}\ \ \ i_{d-1}<n,$$
 $$b_{i_0\cdots i_{d-1}}:=0\ \ \ \textup{if}\ \ i_{d-1}=n.$$
 One easily checks, using the cocycle condition for $(a_{i_0\cdots i_d})_{0\leq i_0<\cdots <i_d\leq n}$, that the latter is the image by $\partial_{d-1}:B_{d-1}\rightarrow B_d$  of $(b_{i_0\cdots i_{d-1}})_{0\leq i_0<\cdots <i_{d-1}\leq n}\in B_{d-1}$ hence $H^d(Z_{\delta})=0$.
 \endproof
 
 \medskip
 
The following lemma will be used in the proofs of Propositions \ref{P8} and \ref{P10}:
 
\begin{lemma}\label{L33}
 Let $B$ be an integral $\delta$-graded ring with $B(0)=R$ and with $p$ as a prime element. Let $w\in W$ and let $(f_i)_{i\in\{0,\ldots,n\}}$ be a family of elements of $B(w)\setminus pB(w)$. Assume that the principal ideals generated in $B\otimes_R K$ by each of the 
 members of the family $(f_i^{\phi^j})_{(i,j)\in \{0,1\}\times (\mathbb N\cup\{0\})}$
are prime and distinct. Then, for $Z_{\delta}:=\Proj_{\delta}(B,(f_0,\ldots,f_n))$ we have
$H^0(Z_{\delta})=R$.
 \end{lemma}
 
\noindent
 {\it Proof.}
 Every element in $H^0(Z_{\delta})$ can be written as 
 $$\frac{F_0}{f_0^{w_0}}=\frac{F_1}{f_1^{w_1}}$$ with $w_i\in W_+$ and $F_i\in B(ww_i)$ for $i\in \{0,1\}$.
 Clearing out the denominators we get that $f_0^{w_0}$ divides $F_0$ in $B\otimes_R K$ hence in $B$ (because $f_0$ is not divisible in $B$ by the prime element $p$). So
 $$\frac{F_0}{f_0^{w_0}}\in B(0)=R.$$
 \endproof
 
 For the proof of Proposition \ref{P8} we will need an extra lemma:

\begin{lemma}\label{L34}
Let $(x_1),(y_1),(x_2),(y_2)$ be distinct principal ideals of an integral domain $B$ such that $y_1x_2\not\in (x_1,y_2)$ and let $f_1,f_2\in B$. Then for all non-negative integers $n_1,m_1,n_2,m_2$ we have
\begin{equation}\label{EQ143}
\frac{y_1 x_2}{x_1y_2}\neq \frac{f_1}{x_1^{n_1}y_1^{m_1}} - \frac{f_2}{x_2^{n_2}y_2^{m_2}}.
\end{equation}
\end{lemma}

\noindent
{\it Proof.}
Assume there exist $n_1,m_1,n_2,m_2\in\mathbb N\cup\{0\}$ such that we have equality in (\ref{EQ143}). We can assume the sum $n_1+m_1+n_2+m_2$ is minimal.
Clearing out the denominators and using the fact that the principal prime ideals 
$(x_1),(y_1),(x_2),(y_2)$ are distinct, 
we immediately get that $n_1=1$, $m_1=0$, $n_2=0$, $m_2=1$, hence
$y_1x_2=y_2f_1-x_1f_2\in (x_1,y_2)$, a contradiction.
\endproof

\medskip

We view Lemma \ref{L34} as an abstraction of the fact that if $B$ is a field and $P,Q\in \mathbb P_B^1(B)$ are two distinct points, then we have $H^1(X,\mathcal O_X)\neq 0$ for
$$X:=(\mathbb P_B^1 \times \mathbb P_B^1) \setminus \{(P,Q),(Q,P)\}.$$

\medskip

\noindent
{\it Proof of Proposition \ref{P8}}.
Its identity $H^0(\mathbb P^1_{\delta})=R$ follows from Lemma \ref{L33} applied with $B=\mathbb S_1$.
For $H^1(\mathbb P^1_{\delta})$ it is enough to check that there exists no non-trivial $R$-linear combination 
 \begin{equation}\label{EQ144}
 \sum_{b,a\in \mathbb Z, b> a \geq 0} \lambda_{ab}
 \frac{z_1^{\phi^a}z_0^{\phi^b}}{z_0^{\phi^a}z_1^{\phi^b}}\in \mathbb S_{1,\langle z_0z_1\rangle},\ \ \lambda_{ab}\in R,
 \end{equation}
which can be written as a difference
 \begin{equation}\label{EQ145}
 \frac{F_0}{z_0^{w_0}}-\frac{F_1}{z_1^{w_1}},\ \ \ F_0\in \mathbb S_1(w_0), F_1\in \mathbb S_1(w_1),\ w_0,w_1\in W_+.\end{equation}
 To check this, we assume there exists an equality 
 $$ \sum_{b,a\in \mathbb Z, b> a \geq 0} \lambda_{ab}
 \frac{z_1^{\phi^a}z_0^{\phi^b}}{z_0^{\phi^a}z_1^{\phi^b}}= \frac{F_0}{z_0^{w_0}}-\frac{F_1}{z_1^{w_1}}
 $$
 between the elements appearing in Equations (\ref{EQ144}) and (\ref{EQ145})
 such that we have $\lambda_{a_0b_0}\neq 0$ for some integers $b_0>a_0\geq 0$.
 We view this equality as an identity between rational functions in the indeterminates $z_i^{\phi^c}$ and 
  substitute $z_i^{\phi^c}$ by $1$ for $c\in (\mathbb N\cup \{0\})\setminus \{a_0,b_0\}$ and $i\in\{0,1\}$ in this identity. Under this substitution, if $(a,b)$ is a pair of integers with $b>a\geq 0$ and 
  $(a,b)\neq (a_0,b_0)$, then  the term $\lambda_{ab} \frac{z_1^{\phi^a}z_0^{\phi^b}}{z_0^{\phi^a}z_1^{\phi^b}}$ becomes $\lambda_{a_0b} \frac{z_1^{\phi^{a_0}}}{z_0^{\phi^{a_0}}}$ if $a=a_0$, becomes $\lambda_{ab_0} \frac{z_0^{\phi^{b_0}}}{z_1^{\phi^{b_0}}}$ if $b=b_0$, and becomes $\lambda_{ab}$ if $a\neq a_0$ and $b\neq b_0$. 
  Moreover, the difference $\frac{F_0}{z_0^{w_0}}-\frac{F_1}{z_1^{w_1}}$ becomes a difference $\frac{F'_0}{z_0^{w'_0}}-\frac{F'_1}{z_1^{w'_1}}$ in which
   $$F'_0,F'_1\in R[z_0^{\phi^{a_0}},z_0^{\phi^{b_0}},z_1^{\phi^{a_0}},z_1^{\phi^{b_0}}]$$
 and $w'_0,w'_1$ are linear combinations of $\phi^{a_0}$ and $\phi^{b_0}$ with non-negative integers coefficients .
We get an equality 
  $$ 
  \frac{z_1^{\phi^{a_0}}z_0^{\phi^{b_0}}}{z_0^{\phi^{a_0}}z_1^{\phi^{b_0}}}= \frac{G_0}{z_0^{v_0}}-\frac{G_1}{z_1^{v_1}}$$  where 
   $v_0:=w'_0+\phi^{a_0}$, $v_1:=w_1'+\phi^{b_0}$,    
   $$
    G_0:=\lambda_{a_0 b_0}^{-1}z_0^{v_0}\left(\frac{F'_0}{z_0^{w'_0}}-\sum_{b\neq b_0}
    \lambda_{a_0b} \frac{z_1^{\phi^{a_0}}}{z_0^{\phi^{a_0}}}-\sum_{a\neq a_0,b\neq b_0} \lambda_{ab}
    \right)\in K[z_0^{\phi^{a_0}},z_0^{\phi^{b_0}},z_1^{\phi^{a_0}},z_1^{\phi^{b_0}}],
   $$
   $$
    G_1:=\lambda_{a_0 b_0}^{-1}z_1^{v_1}\left(\frac{F'_1}{z_1^{w'_1}}+\sum_{a\neq a_0}
    \lambda_{ab_0} \frac{z_0^{\phi^{b_0}}}{z_1^{\phi^{b_0}}}
    \right)\in K[z_0^{\phi^{a_0}},z_0^{\phi^{b_0}},z_1^{\phi^{a_0}},z_1^{\phi^{b_0}}].
   $$
   We get a contradiction based on Lemma \ref{L34} applied with $$(B,x_1,y_1,x_2,y_2):=(\mathbb S_1,z_0^{\phi^{a_0}},z_0^{\phi^{b_0}},z_1^{\phi^{a_0}},z_1^{\phi^{b_0}}).$$\endproof
 
 \medskip
 
 \noindent
 {\it Proof of Proposition \ref{P10}.}
 By Theorem \ref{T16} (a), we have 
 $$\mathbb I_{1,\ord}\otimes_R K=K[f^{\langle 1 \rangle},f^{\langle 2 \rangle},f^{\langle 3 \rangle},\ldots, f^{\partial}, (f^{\partial})^{-1}, (f^{\partial})^{\phi}, (f^{\partial})^{-\phi},(f^{\partial})^{\phi^2},(f^{\partial})^{-\phi^2}\ldots],$$
 where the elements $f^{\langle 1 \rangle},f^{\langle 2 \rangle},f^{\langle 3 \rangle},\ldots, f^{\partial}, (f^{\partial})^{\phi},(f^{\partial})^{\phi^2},\ldots$
 are $K$-algebraically independent. For all integers $a\geq 0$ and $i\geq 1$ we have 
 $$(f^{\langle i \rangle})^{\phi^a}=(f^{\partial})^{-2\frac{\phi^a-1}{\phi-1}}\cdot f^{\langle i+a\rangle}.$$
Setting $f_i:=f_i^{\lambda}$ for $i\in\{1,2\}$ 
 and 
 $$u_{ia}:=\lambda_{i1}^{\phi^a} f^{\langle a+1\rangle}+\lambda_{i2}^{\phi^a} f^{\langle a+2\rangle}$$
 we have
 $$f_i^{\phi^a}=(f^{\partial})^{-2\frac{\phi^a-1}{\phi-1}}\cdot u_{ia}.$$
 
 \begin{claim}\label{Claim3}
 The classes in $H^1(Z_{\delta}^{\lambda})$ of the elements
 $$\eta_{ab}:=\frac{f_2^{\phi^a} f_1^{\phi^b}}{f_1^{\phi^a}f_2^{\phi^b}}=\frac{u_{2a}u_{1b}}{u_{1a}u_{2b}}\in \mathbb I_{1,\ord,\langle f_1 f_2\rangle},\ \ \ 0\leq a< b,$$
 are $R$-linearly independent; in other words, no non-trivial $R$-linear combination of the elements $\eta_{ab}$ can be written in the form
 \begin{equation}\label{EQ146}
 \frac{F_1}{f_1^{w_1}}-\frac{F_2}{f_2^{w_2}}\end{equation}
 with $w_i\in W$ and $F_i\in \mathbb I_{1,\ord}(-2w_i)$ for $i\in\{1,2\}$.
 \end{claim}
 
Claim \ref{Claim3} implies the equality $\phirank_R(H^1(Z_{\delta,2}^{\lambda}))=\infty$ of Proposition \ref{P10}. 
 
To prove Claim \ref{Claim3} assume there exists a non-trivial $R$-linear combination of the elements $\eta_{ab}$ which can be written in the form (\ref{EQ146}) and seek a contradiction. Fix a pair $(a_0,b_0)$ of integers with $0\leq a_0< b_0$
 such that the coefficient of $\eta_{a_0b_0}$ in this linear combination is non-zero. We take $$B:=K[f^{\langle n \rangle}|n\in\mathbb N]$$ and let $\Sigma\subset B$ be the multiplicative monoid generated by all $u_{ic}$ with $i\in\{1,2\}$ and integers $c\geq 0$, 
 $c\not\in \{a_0,b_0\}$.
 Setting  $(f^{\partial})^{\phi^j}=1$ in $F_1, F_2$ for all $j\in \mathbb N\cup\{0\}$
 we get (exactly as in the proof of Proposition \ref{P8}) an equality in the ring of fractions $\Sigma^{-1}B$ of the form
 $$\frac{u_{2a_0}u_{1b_0}}{u_{1a_0}u_{2b_0}}=\frac{H_1}{u_{1a_0}^{n_1} u_{1b_0}^{m_1}}-\frac{H_2}{u_{2a_0}^{n_2} u_{2b_0}^{m_2}}$$
 where $n_1,m_1, n_2, m_2\in\mathbb N\cup\{0\}$ and $H_1,H_2\in \Sigma^{-1}B$.
 By Lemma \ref{L34} we get a contradiction if we prove the following:
 
 \begin{claim}\label{Claim4}
 For $0\leq a<b$ the principal ideals $(u_{1a}), (u_{1b}), (u_{2a}), (u_{2b})$ of $\Sigma^{-1}B$ are prime and distinct and we have $u_{1b}u_{2a}\not\in (u_{1a},u_{2b})\cdot \Sigma^{-1}B$.
\end{claim}
 
It is an easy exercise to show that Claim \ref{Claim4} follows from the following:
 
 \begin{claim}\label{Claim5}
For $\mathbb O:=\{1,2\}\times\mathbb Z_{\geq 0}$ the following three properties hold:
 
 \medskip
 {\bf (a)} The principal ideals $(u_{ij})$ of $B$ indexed by $(i,j)\in\mathbb O$ are prime and pairwise distinct.
 
 \smallskip
 {\bf (b)} For every two distinct pairs $(i_1,j_1),(i_2,j_2)\in\mathbb O$ the ideal $(u_{i_1 j_1}, u_{i_2 j_2})B$ is prime in $B$.
 
 \smallskip
 {\bf (c)} For every three distinct pairs $(i_1,j_1), (i_2,j_2), (i_3,j_3)\in\mathbb O$ we have
 $u_{i_3j_3}\not\in (u_{i_1j_1},u_{i_2j_2})B$.
\end{claim}
 
To check Claim \ref{Claim5} we view $B=K[f^{\langle n \rangle}|n\in\mathbb N]$ as a graded $K$-algebra by giving  $f^{\langle n\rangle}$  degree $1$ for all $n\geq 1 $. Parts (a) and (b) of Claim \ref{Claim5} are clear because $u_{ij}$ are pairwise $K$-linearly independent linear forms in $B$; here we are using that all entries of the matrix $\lambda$ are non-zero. Now using the grading on $B$ we see that part (c) of Claim \ref{Claim5} is equivalent to showing that for all three distinct pairs $(i_1,j_1), (i_2,j_2), (i_3,j_3)\in\mathbb O$ the linear forms
 $u_{i_1j_1}, u_{i_2j_2}, u_{i_3j_3}$ are $K$-linearly independent. Let $\Gamma$ be the $3\times \infty$ matrix whose $(k,l)$-entries are the coefficients of $f^{\langle l \rangle}$ in $u_{i_kj_k}$.
 We have to show that $\Gamma$ has a non-zero $3\times 3$ minor.
 We may assume $j_1\leq j_2\leq j_3$. As the case $j_1=j_2=j_3$ does not occur (as $i_1,i_2,i_3\in\{1,2\})$ are not distinct), the $3\times 3$ minor corresponding to the columns $l_1,l_2,l_3$ is non-zero, where, depending on the three disjoint possible cases, the triple $(l_1,l_2,l_3)$ is as follows:
 
 \smallskip
 {\it Case 1}: If $j_1<j_2<j_3$, then $(l_1,l_2,l_3):=(j_1+1, j_2+1, j_3+1)$.
 
 \smallskip
 {\it Case 2}: If $j_1=j_2<j_3$, then $(l_1,l_2,l_3):=(j_1+1, j_2+2, j_3+2)$. 
 
 \smallskip
 {\it Case 3}: If $j_1<j_2=j_3$, then $(l_1,l_2,l_3):=(j_1+1, j_2+1, j_3+2)$.
 
 \smallskip
In all three cases we used, again, the fact that all entries of the matrix $\lambda$ are non-zero. At this point Claim \ref{Claim5} is proved. Thus the statement about $H^1$ in Proposition \ref{P10} is proved. The statement about $H^0$ follows from Lemma \ref{L33} and Claim \ref{Claim5} (a).\endproof
 
 \medskip
 \noindent
 {\it Proof of Proposition \ref{P11}}.
 Let $f_i:=(F_i)^{\diamondsuit}$, $i\in \{0,1\}$, $F_i\in \mathbb H_g(1)$ primitive. Let $B:=\mathbb I_{g,\ord}$. By Lemma 
 \ref{L33} it is enough
 to show that the principal ideals generated in $B\otimes_R K$ by each of the 
 members of the family $(f_i^{\phi^j})_{(i,j)\in \{0,1\}\times (\mathbb N\cup\{0\})}$
are prime and distinct.
 By Equation (\ref{EQ098}) it is enough to show that the principal ideals 
 $$\mathfrak a_{ij}=((F_i^{\tau^j\sigma^j})^{\diamondsuit})$$
  generated in $B\otimes_R K$ by the members of the family 
 \begin{equation}\label{EQ147}
 ((F_i^{\tau^j\sigma^j})^{\diamondsuit})_{(i,j)\in \{0,1\}\times (\mathbb N\cup\{0\})}\end{equation}
 are prime and distinct. By
 Corollary \ref{C15} (e) and (f) the ideals $\mathfrak a_{ij}$ are prime and it is enough to show that the members of the family (\ref{EQ147}) are pairwise $R$-linear independent.
 Assume this is not the case. As the $\mathbb Z_p$-algebra monomorphism $\sigma\circ \tau$ fixes no polynomial of degree $\geq 1$ and as $\diamondsuit$ is an $\mathbb Z_p$-algebra monomorphism there exist two indices $j_0,j_1\in \mathbb N\cup\{0\}$ with $j_1\geq j_0$ and $\lambda\in K$ such that
 $$F_{i_0}^{\tau^{j_0}\sigma^{j_0}}=\lambda\cdot F_{i_1}^{\tau^{j_1}\sigma^{j_1}}=
 (\phi^{-j_1}(\lambda)\cdot F_{i_1})^{\tau^{j_1}\sigma^{j_1}}$$
 where $\{i_0,i_1\}=\{0,1\}$.
So we have
 $$F_{i_0}=(\phi^{-j_1}(\lambda)\cdot F_{i_1})^{\tau^j\sigma^j}$$
 where $j:=j_1-j_0\geq 0$. As $F_0$ and $F_1$ are primitive we get
 $j=0$. But then $F_0$ and $F_1$ are $K$-linearly dependent, a contradiction.
 \endproof
 
 \medskip
 
 In order to prove Proposition \ref{P12} we need the following lemma in which we use the forms $f^{[a]}$ introduced by Equation (\ref{EQ066}):
 
 \begin{lemma}\label{L35}
 The forms $\det(f^{[1]}), \det(f^{[2]})\in \mathbb I_{g,\ord}(-2)$ are primitive and $R$-linearly independent.
 \end{lemma}
 
\noindent
 {\it Proof.}
 The linear independence follows from the Serre--Tate expansion principle, see Corollary \ref{C9} (b). To prove that $\det(f^{[1]})$ and $\det(f^{[2]})$ are primitive,
 recall (see Equation (\ref{EQ090})) the definition of the partial polarizations $\Theta_0,\ldots,\Theta_g\in\mathbb H_g(1)$ of the invariant polynomial $\det$ via the identity
 $$\det(y_0T+y_1T')=\sum_{i=0}^g \Theta_i (T,T') y_0^{g-i}y_1^i$$
 and (see Corollary \ref{C9} (d)) the formula
 $f^{[2]}=pf^{\langle 1\rangle}+f^{\langle 2\rangle}$.
 Consider the polynomial $\Theta^{[2]}:=\sum_{i=0}^g p^{g-i}\Theta_i (T,T')$.
As $\det(f^{[1]})=\Theta_0(f^{\langle 1 \rangle},0)=(\Theta_0(T,0))^{\diamondsuit}$, $\det(f^{[1]})$ is primitive.
On the other hand we have:
$$\det(f^{[2]})=\det(pf^{\langle 1\rangle}+f^{\langle 2\rangle})=\sum_{i=0}^g\Theta(f^{\langle 1 \rangle}, f^{\langle 2 \rangle})p^{g-i}=(\Theta^{[2]})^{\diamondsuit}.$$
So we are left to check that $\Theta^{[2]}$ is primitive. But the coefficient of the monomial $T_{11}\cdots T_{g-1,g-1}T'_{gg}$ in $\Theta_i(T,T')$ is $1$ if $i=1$ and $0$ if $i\in \{0,2,\ldots,g\}$. So $\Theta^{[2]}\not\in R[T']$ which ends our proof.
 \endproof
 
 \medskip
 \noindent
 {\it Proof of Proposition \ref{P12}}.
 For $\varphi_1:=\det(f^1)^{1+\phi^2}$, $\varphi_2:=\det(f^2)^{1+\phi}$
 we have
 $$R\subset H^0(Z_{\delta},(\varphi_1,\varphi_2)) =\mathbb I_{g,\langle \varphi_1 \rangle}\cap \mathbb I_{g,\langle \varphi_2\rangle}\subset \mathbb I_{g,\ord,\langle \varphi_1 \rangle}\cap \mathbb I_{g,\ord,\langle \varphi_2\rangle}$$
 $$\subset \mathbb I_{g,\ord,\langle \det(f^1) \rangle}\cap \mathbb I_{g,\ord,\langle \det(f^2) \rangle}\subset \mathbb I_{g,\ord,\langle \det(f^{[1]}) \rangle}\cap \mathbb I_{g,\ord,\langle \det(f^{[2]}) \rangle}=R,
 $$
 where the last equality follows from Proposition \ref{P11} and Lemma \ref{L35}.
 \endproof
 
\newpage\section{Background and required results in GIT}\label{S4}

We recall that this section is completely self-contained, can be read out at any time, and contains results in geometric invariant theory (GIT) that are often used and cited in the prior sections. Its aim is to study the structure of algebras of invariants of the group $\pmb{\SL}_g$ acting
on tuples of quadratic forms and on tuples of endomorphisms in arbitrary (mixed) characteristics. Though this theory has classical roots, unless explicitly stated, all results in this section are new. 

We will start, in Subsection \ref{S41},
by giving formulas (see Subsubsection \ref{S411}) for the (Krull) dimensions of algebras of invariants of $\pmb{\SL}_g$ acting on tuples of quadratic forms (see Theorem \ref{T25}) and on tuples of endomorphisms (see Theorem \ref{T26}); these formulas may be classical but we could not find a reference for them. We continue by giving an estimate (see Theorem \ref{T27}) for the dimensions of the images of certain basic homomorphisms between the two types of algebras of invariants; this estimate plays a key role in our applications to Sections \ref{S2} and \ref{S3}. The basic tool in the proofs of the above three theorems carried on respectively in Subsubsections \ref{S413} to \ref{S415} is a study of Hilbert--Mumford stability in the corresponding contexts (see Subsubsection \ref{S412}).

Subsection \ref{S42} offers a detailed study of the invariants of $\pmb{\SL}_g$ acting
on tuples of quadratic forms for arbitrary $g$. 
In particular  we will introduce some special  `$\Theta$ and $\Upsilon$ invariants' 
(see Subsubsection \ref{S4206})
and we will compare the algebras of all invariants with the algebras generated by the $\Theta$
invariants (see, for instance,  Theorems \ref{T29} and \ref{T30} and Corollary \ref{C24}). We will deduce information on the Hilbert series of the algebras of invariants (see Corollaries \ref{C21} and \ref{C22}), the smallest number of homogeneous generators
(see Proposition \ref{P15} (a)), and transcendence bases (see Theorem \ref{T31} and Corollary \ref{C23} (b)). 

Subsection \ref{S43} gives finer results in the case $g=2$
on Hilbert series (see Proposition \ref{P18} and Corollary \ref{C26}), the smallest number of homogeneous generators
(see Proposition \ref{P17}), and relations (see Corollary \ref{C25}).

In Subsection \ref{S44} we prove a general result about the $\sigma$-finite generation of algebras of invariants for arbitrary reductive groups.

Throughout this section we adopt some general notation as follows.

Let $B^\times$ be the group of units of a commutative ring $B$. 

Let $\mathbb K$ be an arbitrary algebraically closed field. Let $p:=\ch(\mathbb K)$; so either $p=0$ or $p>0$ is a prime number. Product of schemes over $\Spec(\mathbb K)$ will be denoted by $\times$. 

For a finite dimensional vector space $V$ over $\mathbb K$, let
$$\pmb{\V}:=\Spec(\Sym(V^*)),$$ 
where $\mathcal O(\pmb{\V}):=\Sym(V^*)$ is the symmetric algebra of the dual $V^*$ of $V$. Thus, for a $\mathbb K$-algebra $\mathbb B$ we have $\pmb{\V}(\mathbb B)=\mathbb B\otimes_{\mathbb K} V$, and in particular $\pmb{\V}(\mathbb K)=V$. Let $r\in\mathbb N\cup\{0\}$. Let $V^{r+1}$ be the direct sum of $r+1$-copies of $V$. 

Let $g\in\mathbb N$. Let $\pmb{\GL}_g:=\pmb{\GL}_{g,\mathbb K}$ be the reductive group over $\Spec(\mathbb K)$ of $g\times g$ invertible matrices. Let $\pmb{\SL}_g:=\pmb{\SL}_{g,\mathbb K}$ be the semisimple group over $\Spec(\mathbb K)$ which is the kernel of the determinant homomorphism $\det:\pmb{\GL}_g\rightarrow\mathbb G_{m,\mathbb K}$. Let $\pmb{\SO}_g$ be the closed subgroup scheme of $\pmb{\SL}_g$ of orthogonal matrices of determinant $1$. If $p=2$, then $\pmb{\SO}_g$ is smooth if and only if $g=1$. If $p\neq 2$, then $\pmb{\SO}_1$ is trivial, $\pmb{\SO}_2\simeq \mathbb G_{m,\mathbb K}$, and $\pmb{\SO}_g$ is semisimple for $g\geq 3$. Let $\pmb{\GT}_g$ be the standard maximal torus of $\pmb{\GL}_g$ of invertible diagonal matrices. Thus $\pmb{\T}_g:=\pmb{\GT}_g\cap \pmb{\SL}_g$ is the standard maximal torus of $\pmb{\SL}_g$. Let $1_g$ be the identity element of $\pmb{\GL}_g(\mathbb K)$ and of its subgroups.

Let $G$ be a reductive group over $\Spec(\mathbb K)$. Let $G^0$ be the connected component of its identity element and let $Z(G)$ be its center. All semisimple groups over $\mathbb K$ will be connected. Each $G$-module $V$ will be a finite dimensional left $G$-module given by a representation $\rho_V:G\to\pmb{\GL}_V$. Thus $V^*$ is also a $G$-module given by the (composite homomorphism) representation $\rho_{V^*}:G\xrightarrow{\textup{inv}} G\xrightarrow{\rho_V}\pmb{\GL}_V\simeq \pmb{\GL}_{V^*}$, where $G\xrightarrow{\textup{inv}} G$ is the involutory inverse antiisomorphism and $\pmb{\GL}_V\simeq \pmb{\GL}_{V^*}$ is the canonical antiisomorphism. For $(h,v,w)\in G(\mathbb K)\times V\times V^*$, we have $(h\cdot w)(v)=w(h^{-1}\cdot v)\in\mathbb K$. 

\subsection{Stable points and dimensions of quotients}\label{S41}

Let $\Mat_g$ be the $\mathbb K$-vector space of $g\times g$ matrices with entries in $\mathbb K$, and let 
$$\V_g:=\{M\in\Mat_g|M^{\t}=M\}\subset\Mat_g$$ 
be its subspace of symmetric matrices, where the upper right index $\t$ denotes transpose. Let $\Z_g$ be the subspace of $\V_g$ of all diagonal matrices. For $\lambda_1,\ldots,\lambda_g\in\mathbb K$ let 
$$\Diag(\lambda_1,\ldots,\lambda_g)\in\Z_g$$ 
be such that for $i\in\{1,\ldots,g\}$ its $ii$ entry is $\lambda_i$. Let 
$$\L_g:=\mathbb K1_g\times\Z_g\subset\V_g^2;$$ 
we have $\dim_{\mathbb K}(\L_g)=g+1$.
We recall that $M\in\Mat_g$ is called {\it alternate} if $M^{\t}=-M$ and, in case $p=2$, all its diagonal entries are $0$. As such, if $p=2$, $Q\in\V_g$ is called {\it non-alternate} if at least one of its diagonal entries is non-zero. 

Let $\Trace:\Mat_g\to\mathbb K$ be the trace $\mathbb K$-linear map. 

By our convention above we have 
$$\pmb{\Z}_g=\Spec(\Sym(\Z_g^*))\subset\pmb{\V}_g=\Spec(\Sym(\V_g^*))\subset\pmb{\Mat}_g=\Spec(\Sym(\Mat_g^*)).\footnote{In Sections \ref{S1} to \ref{S3}, $\pmb{\Mat}_g$ is a ring scheme over $\Spec(\mathbb Z)$, but in this subsection, $\pmb{\Mat}_g$ will be a shorter notation for $\pmb{\Mat}_{g,\mathbb K}$ viewed just as a vector group variety over $\Spec(\mathbb K)$.}$$
In what follows by $\pmb{\SL}_g$-module we will mean a left $\pmb{\SL}_g$-module.

We view $\Mat_g$ as an $\pmb{\SL}_g$-module via the classical congruence action 
$$\pmb{\SL}_g\times\pmb{\Mat}_g\xrightarrow{\cong}\pmb{\Mat}_g$$ 
given on $\mathbb K$-valued points by the rule:
$(\Lambda,M)\mapsto \Lambda M\Lambda^{\t}$. Then $\V_g$ is an $\pmb{\SL}_g$-submodule of $\Mat_g$ and thus $\V_g^{r+1}$ and its dual $(\V_g^{r+1})^*=(\V_g^*)^{r+1}$ naturally become $\pmb{\SL}_g$-modules. 

\subsubsection{Main results on dimensions of quotients}\label{S411}
Our first goal is to prove the following formula which is likely to be classical but for which we could not find a reference:

\begin{thm}\label{T25}
For all integers $r\geq 1$ we have
$$\dim(\Sym((\V_g^{r+1})^*)^{\pmb{\SL}_g})=(r+1)\cdot \frac{g(g+1)}{2}-g^2+1.$$
\end{thm}

We will also consider the classical conjugation action of $\pmb{\SL}_g$ on $\Mat_g$ 
(when not clear from context, it will be specified which one of the two structures of $\Mat_g$ as a $\pmb{\SL}_g$-module is used). We view $\Mat_g$ as a left $\pmb{\SL}_g$-module via the action 
$$\pmb{\SL}_g\times\pmb{\Mat}_g\xrightarrow{\con}\pmb{\Mat}_g$$ 
given on $\mathbb K$-valued points by the rule:
$(\Lambda,M)\mapsto \Lambda M\Lambda^{-1}$. Then $\Mat_g^{r+1}$ naturally becomes an $\pmb{\SL}_g$-module isomorphic to its dual $(\Mat_g^*)^{r+1}$. 

\medskip
We will prove the following formula which, also, is likely to be classical:

\begin{thm}\label{T26}
For all integers $r\geq 1$, for the conjugation action we have
$$\dim(\Sym((\Mat_g^{r+1})^*)^{\pmb{\SL}_g})=r\cdot g^2+1.$$
\end{thm}

Let $n\in\mathbb N$. We consider now the morphism 
\begin{equation}\label{EQ148}
\pmb{\pi}_n:\pmb{\V}_g^{n+1}\rightarrow \pmb{\Mat}_g^n\end{equation}
 given on $\mathbb K$-valued points by the map
\begin{equation}\label{EQ149}
\pi_n:=\pmb{\pi}_n(\mathbb K):\V_g^{n+1}\rightarrow \Mat_g^n,\ \ \ \pi_n(Q_0,\ldots,Q_n)= (Q_0Q_1^*,Q_1Q_2^*,\ldots, Q_{n-1}Q_n^*),\end{equation}
where we recall that $M^*\in\Mat_g$ denotes the adjugate of $M\in\Mat_g$; we have $M^*M=MM^*=\det(M)1_g$. The morphism (\ref{EQ149})
is $\pmb{\SL}_g$-equivariant 
(with $\pmb{\SL}_g$ acting on $\V_g$ by congruence and on $\Mat_g$ by conjugation)
and hence
 induces a $\mathbb K$-algebra homomorphism
$$\Sym((\Mat_g^n)^*)^{\pmb{\SL}_g}\rightarrow \Sym((\V_g^{n+1})^*)^{\pmb{\SL}_g}$$
whose image we denote by $\reallywidetilde{\Sym((\Mat_g^n)^*)^{\pmb{\SL}_g}}$. 
 In Sections \ref{S2} and \ref{S3}, the analogues of $\mathbb K$-algebras $\Sym((\Mat_g^n)^*)^{\pmb{\SL}_g}$ and $\reallywidetilde{\Sym((\Mat_g^n)^*)^{\pmb{\SL}_g}}$ but over the field of fractions $K$ of the ring of $p$-typical Witt vectors with coefficients in an algebraic closure of $\mathbb F_p$ were denoted by $\mathbb H^n_{g,\con,K}$ and $\reallywidetilde{\mathbb H^n_{g,\con,K}}$ (respectively). 
We will prove the following formula used in Section \ref{S3}:

\begin{thm}\label{T27}
For all integers $n\geq 2$ we have
$$\dim(\reallywidetilde{\Sym((\Mat_g^n)^*)^{\pmb{\SL}_g}})= (n-1)\frac{g(g+1)}{2}+g.$$
\end{thm}

For the proofs of Theorems \ref{T25}, \ref{T26} and \ref{T27}
we need to first recall from \cite{MFK} some basic facts of geometric invariant theory.

\subsubsection{Stable points}\label{S412}

Let $V$ be a $G$-module. We recall that a point $v\in V=\pmb{\V}(\mathbb K)$ is called {\it stable} if its $G$-orbit $G\cdot v$ is closed and its stabilizer $\Stab_G(v)$ is finite (see \cite{MFK}, App. to Ch. I, Sect. B, Def.). 
Recall that the $\mathbb K$-algebra of invariants $\mathcal O(\pmb{\V})^G$ is finitely generated (see \cite{MFK}, App. to Ch. I, Sect. A, Thm. A.1.0) and normal; if in addition $G$ is semisimple,
$\mathcal O(\pmb{\V})^G$ is a unique factorization domain.
Write $\pmb{\V}/G:=\Spec
(\mathcal O(\pmb{\V})^G)$ and denote by 
$$\varpi:\pmb{\V}\rightarrow \pmb{\V}/G$$ 
the induced projection which is submersive and a universal categorical quotient if $p=0$ (see \cite{MFK}, Ch. I, Sect. 2, Thm. 1.1) and is submersive and a uniform categorical quotient if $p>0$ (see \cite{MFK}, App. to Ch. I, Sect. C, Thm. A.1.1).\footnote{In the literature one also uses the notation $\pmb{\V}//G$ in place of $\pmb{\V}/G$; for simplicity we will stick here with the notation $\pmb{\V}/G$.} 

Let $V^{\s}\subset V$ be the set of stable points. Then $V^{\s}$ is the set of $\mathbb K$-valued points of a $G$-invariant and open (possibly empty) subscheme $\pmb{\V}^{\s}$ of $\pmb{\V}$ which is also $\mathbb K^\times$-invariant (see \cite{MFK}, App. to Ch. I, Sect. B, Def.). If $\pmb{\X}$ is a locally closed subscheme of $\pmb{\V}$, let $\pmb{\X}^{\s}:=\pmb{\X}\cap \pmb{\V}^{\s}$.

\begin{Fact}\label{F1}
 For each $v\in V^{\s}=\pmb{\V}^{s}(\mathbb K)$, the reduced fiber $(\varpi^{-1}(\varpi(v)))_{\red}$ is $G\cdot v$ and is a closed subscheme of $\pmb{\V}$ of dimension $\dim(G)-\dim(\Stab_G(v))=\dim(G)$. Moreover, $\varpi(\pmb{\V}^{\s})$ is an open subvariety of $\pmb{\V}/G$.
 \end{Fact}
 
\noindent
{\it Proof.} As $\varpi$ is a categorical quotient, we have $G\cdot v\subset \varpi^{-1}(\varpi(v))$. Thus to prove that we have $(\varpi^{-1}(\varpi(v)))_{\red}=G\cdot v$ it suffices to show that the assumption that $(\varpi^{-1}(\varpi(v)))_{\red}\setminus G\cdot v$ has a $\mathbb K$-valued point $v'$ leads to a contradiction. The Zariski closure $\overline{G\cdot v'}$ of $G\cdot v'$ in $\pmb{\V}$ is contained in $(\varpi^{-1}(\varpi(v)))_{\red}$ as well as in $V\setminus G\cdot v$. [Argument: if $\overline{G\cdot v'}\cap G\cdot v\neq\emptyset$, then $\dim(G\cdot v')>\dim(G\cdot v)=\dim(G)$, contradiction.] From \cite{MFK}, Ch. I, Sect. 2, Cor. 1.2 and App. to Ch. I, Sect. C, Cor. A.1.3 we get that there exists a global function on $\pmb{\V}/G$ which is $0$ on $\overline{G\cdot v}=G\cdot v$ and is $1$ on $\overline{G\cdot v'}$. This contradicts that $\overline{G\cdot v'}\cup G\cdot v$ maps to $\varpi(v)\in (\pmb{\V}/G)(\mathbb K)$. 

From the first part we get that $\varpi^{-1}(\varpi(\pmb{\V}^{\s}))=\pmb{\V}^{\s}$. So $\varpi^{-1}(\varpi(\pmb{\V}^{\s}))$ is open in $\pmb{\V}$. As $\varpi$ is submersive, we get that $\varpi(\pmb{\V}^{\s})$ is open in $\pmb{\V}/G$.\endproof

\medskip
For uniform or universal geometric quotients we refer to \cite{MFK}, Ch. 0, Defs. 0.6 and 07.
 
\begin{Fact}\label{F2}
If $\pmb{\W}$ is a non-empty open subvariety of $\varpi(\pmb{\V}^{\s})$, then the morphism $\varpi_{\pmb{\W}}:\varpi^{-1}(\pmb{\W})\rightarrow \pmb{\W}$ is a uniform geometric quotient and we have
$\dim(\pmb{\V}/G)=\dim_{\mathbb K}(V)-\dim(G)$. Moreover, if $p=0$, then $\varpi_{\pmb{\W}}:\varpi^{-1}(\pmb{\W})\rightarrow \pmb{\W}$ is a universal geometric quotient.
\end{Fact}
 
\noindent
{\it Proof.} To check the first and the third (last) part we can assume that $\pmb{\W}$ is affine and this case is a particular case of \cite{MFK}, Ch. I, Sect. 2, Amplific. 1.3 as the geometric fibers of $\varpi_{\pmb{\W}}$ are closed (see Fact \ref{F1}), see also \cite{MFK}, Ch. I, Prop. 1.9 or Thm. 1.10 for the `uniform' part. 

The second part follows from the dimension formula for the generic fiber of the dominant morphism $\varpi_{\pmb{\W}}$ (see Fact \ref{F1}).\endproof

\begin{df}\label{df29}
A $G$-invariant homogeneous polynomial $F\in\mathcal O(\pmb{\V})$ will be called a $G$-{\it separating polynomial on $V$} (or on $\pmb{\V}$) if $F\neq 0$ and for all $v\in V=\pmb{\V}(\mathbb K)$ with $F(v)\neq 0$ the stabilizer $\Stab_G(v)$ is finite. 
\end{df}

We have the following interpretation of stable points in terms of $G$-separating polynomials on $V$ which is a slight improvement of a particular case of \cite{MFK}, Ch. 1, Amplific. 1.11:

\begin{Fact}\label{F3}
The stable locus $V^{\s}$ is non-empty if and only if there exists a $G$-separating polynomial $F$ on $V$. More precisely, for a point $w\in V$ the following two statements are equivalent:

\medskip\noindent
\circled{\textup{1}} We have $w\in V^{\s}$.

\smallskip\noindent
\circled{\textup{2}} There exist a $G$-separating polynomial $F$ on $V$ such that $F(w)\neq 0$. 

\medskip\noindent
Thus, if $F$ is a $G$-separating polynomial on $V$ and $D(F)=(\pmb{\V}/G)\setminus Z(F)$ is its principal open of $\pmb{\V}/G$ (so $Z(F)$ is the zero locus of $F$), then we have $\varpi^{-1}(D(F))\subset \pmb{\V}^{\s}$, the morphism $\varpi^{-1}(D(F))\rightarrow D(F)$ is a uniform geometric quotient and we have $\dim(\pmb{\V}/G)=\dim_{\mathbb K}(V)-\dim(G)$. Moreover, if $p=0$, then the morphism $\varpi^{-1}(D(F))\rightarrow D(F)$ is a universal geometric quotient.
\end{Fact}

\noindent
{\it Proof.} We first prove that $\circled{\textup{2}}\Rightarrow\circled{\textup{1}}$. Let $v\in\varpi^{-1}(D(F))(\mathbb K)$ and we consider $v'\in V=\pmb{\V}(\mathbb K)$ in the Zariski closure of the orbit $G\cdot v$. Then $F(v')=F(v)\neq 0$ and thus $\Stab_{G}(v')$ is finite; therefore the orbits $G\cdot v$ and $G\cdot v'$ have the same dimension equal to $\dim(G)$. This forces $v'\in G\cdot v$, i.e., the orbit $G\cdot v$ is closed, and we get that $v\in V^{\s}$. Thus $\varpi^{-1}(D(F))\subset \pmb{\V}^{\s}$, hence $w\in V^{\s}$, i.e., $\circled{\textup{2}}\Rightarrow\circled{\textup{1}}$, and the second part of the fact follows from Fact \ref{F2}.

To end the proof we are left to prove that $\circled{\textup{1}}\Rightarrow\circled{\textup{2}}$. If $w\in V^{\s}$, then the closed $G$-invariant subsets $G\cdot w$ and $\pmb{\V}\setminus\pmb{\V}^{\s}$ of $\pmb{\V}$ are disjoint. From this and \cite{MFK}, Ch. I, Sect. 2, Cor. 1.2 and App. to Ch. I, Sect. C, Cor. A.1.3 we get that there exists a global function $F$ of $\pmb{\V}/G$ which is $1$ on $G\cdot w$ and is $0$ on $\pmb{\V}\setminus\pmb{\V}^{\s}$. Clearly $F\neq 0$. As $\pmb{\V}\setminus\pmb{\V}^{\s}$ is $\mathbb K^\times$-invariant, we can assume that $F$ is homogeneous. If $v\in V$ is such that $F(v)\neq 0$, then $v\in V^{\s}$ and thus $F$ is a $G$-separating polynomial on $V$ with $F(w)\neq 0$, i.e., $\circled{\textup{1}}\Rightarrow\circled{\textup{2}}$.\endproof
 
\medskip
We have the following trivial functorial property which is also a particular case of \cite{MFK}, Ch. 1, Prop. 1.18 and App. to Ch. I, end of Sect. C:

\begin{Fact}\label{F4}
Let $f:V\rightarrow V'$ be a surjective morphism of $G$-modules. If $F\in\mathcal O(\pmb{\V}')$ is a $G$-separating polynomial on $V'$, then $F\circ f\in\mathcal O(\pmb{\V})$ is a $G$-separating polynomial on $V$ (equivalently, we have $f^{-1}((V')^{\s})\subset V^{\s}$). 
\end{Fact}

\begin{Fact}\label{F5}
We assume that $V$ is a faithful $G$-module such that under the closed embedding homomorphism $\rho_V:G\rightarrow\pmb{\GL}_{V}$ we have $\rho_V(G^0)\subset \pmb{\SL}_{V}$. Then there exists $m\in\{1,\ldots,\dim_{\mathbb K}(V)\}$ such that $(V^m)^{\s}$ is non-empty.
\end{Fact}

\noindent
{\it Proof.} As the sets $(V^m)^{\s}$ do not change if $G$ is replaced by $G^0$, we can assume that $\rho_V$ is an inclusion and $G=G^0\subset\pmb{\SL}_{V}$. Let $d:=\dim_{\mathbb K}(V)$. We identify $V=\mathbb K^d$, $\pmb{\SL}_{V}=\pmb{\SL}_d$ and $V^d=\Mat_d$. Let $\Lambda\in V^d=\Mat_d$ be such that it has determinant $1$. The orbit $\pmb{\SL}_d\cdot \Lambda$ is the closed subset of $\pmb{\Mat}_d$ formed by all matrices of determinant $1$ and the stabilizers $\Stab_{\pmb{\SL}_d}(\Lambda)$ and $\Stab_G(\Lambda)$ are trivial. Thus the orbit $G\cdot \Lambda$ inside $\pmb{\SL}_d\cdot \Lambda$ can be identified with the closed subgroup $G$ of $\pmb{\SL}_d$ and therefore it is closed. We conclude that $\Lambda\in (V^d)^{\s}$.\endproof

\subsubsection{Proof of Theorem \ref{T25}}\label{S413}

To prove Theorem \ref{T25} we can assume that $g\geq 2$ and we first review a few classical facts for $g\geq 2$ (which for $g=1$ are trivial). 

\begin{Fact}\label{F6}
Let $Q_0\in\V_g$ be invertible. If $p=2$, we assume $Q_0$ is non-alternate. Then for each element $\nu\in\mathbb K$ such that $\nu^g=\det(Q_0)$, $Q_0$ is $\pmb{\SL}_g$-equivalent to $\nu 1_g$  (i.e., $Q_0\in\pmb{\SL}_g\cdot \nu 1_g$).
\end{Fact}

\noindent
{\it Proof.} 
As $Q_0$ is non-alternate (no matter what $p$ is), it is congruent to a diagonal matrix (see \cite{A}, Thm. 6). As $Q_0$ is invertible and $\mathbb K$ is algebraically closed, we get that $Q_0$ is congruent to $1_g$. Thus $Q_0$ is $\pmb{\SL}_g$-equivalent to $\lambda 1_g$ where $\lambda\in\mathbb K$ is such that $\lambda^g=\det(Q_0)$. We can assume $\nu\neq\lambda$. As $\lambda^g=\nu^g$, $\lambda 1_g$ and $\nu 1_g$ belong to the same $\pmb{\T}_g$-orbit but not to the same $\pmb{\SO}_g$-orbit and not necessarily to the same $Z(\pmb{\SL}_g)$-orbit: concretely, if $a:=\sqrt[2g]{\frac{\nu}{\lambda}}\in\mathbb K^\times$, then $a^g\in\{-1,1\}$ and for $\Lambda\in\{a1_g,\Diag(-a,a,\ldots,a)\}\cap\pmb{\T}_g(\mathbb K)$ we have $\Lambda \lambda 1_g\Lambda^{\t}=\nu 1_g$.\endproof

\begin{Fact}\label{F7}
Let the pair $(Q_0,Q_1)\in\V_g^2$ be such that $Q_0$ is invertible and $Q_0^{-1}Q_1$ has $g$ distinct eigenvalues. If $p=2$ we also assume that $\Trace(Q_0^{-1}Q_1)\neq 0$. Then $(Q_0,Q_1)$ is $\pmb{\SL}_g$-equivalent to a pair $(\lambda 1_g,D_1)\in \L_g$ with $\lambda\in\mathbb K^\times$.
\end{Fact}

\noindent
{\it Proof.} We will first show that the assumption that $p=2$ and $Q_0$ is alternate leads to a contradiction. As the rank of an alternate matrix is even and as $Q_0$ is invertible, this assumption implies that $g$ is even and from \cite{A}, Thm. 4 we get that $Q_0$ is congruent to a block diagonal matrix $\mathcal J_g:=\Diag(\mathcal J,\ldots,\mathcal J)$ whose $g/2$ diagonal blocks are the elementary matrix $\mathcal J:=\begin{bmatrix} 0 & 1\\ 1& 0 \\ \end{bmatrix}$ of type I. Thus $Q_0$ is $\pmb{\SL}_g$-equivalent to $\nu\mathcal J_g$ with $\nu\in\mathbb K^\times$ such that $\nu^g=\det(Q_0)$ and as such we can assume that $Q_0=\nu\mathcal J_g$. For $a,b,c\in\mathbb K$ we have $\Trace\Big(\mathcal J\begin{bmatrix} a & b\\ b& c\\ \end{bmatrix}\Big)=2b=0$. From this, as $Q_1$ is symmetric and $g$ is even, decomposing $Q_1=(q_{ij})_{1\leq i,j\leq g/2}$ into blocks $q_{ij}=q_{ji}^{\t}\in\Mat_2$ we compute $\Trace(Q_0^{-1}Q_1)=\nu^{-1}\sum_{i=1}^{g/2}\Trace(\mathcal Jq_{ii})=0$, contradiction. 

Thus, if $p=2$, then $Q_0$ is non-alternate, and from Fact \ref{F6} we get that we can assume that $Q_0=\lambda 1_g$ for some $\lambda\in\mathbb K^\times$. Thus, as $Q_0^{-1}Q_1$ has $g$ distinct eigenvalues, the same holds for $Q_1$. We write $Q_1=\Lambda_2D_2\Lambda_2^{-1}$, where $\Lambda_2\in\pmb{\GL}_g(\mathbb K)$ and $D_2\in\Z_g$ has $g$ distinct eigenvalues. As $Q_1\in\V_g$, from the identity $Q_1^{\t}=Q_1$ we easily get that $\Lambda_2^{\t}\Lambda_2$ commutes with $D_2$ and hence, as $D_2$ has $g$ distinct eigenvalues, $\Lambda_2^{\t}\Lambda_2\in\Z_g$ is diagonal, and we write $\Lambda_2^{\t}\Lambda_2=\Lambda_1^{\t}\Lambda_1$ with $\Lambda_1\in\pmb{\GL}_g(\mathbb K)\cap\Z_g$. For $\Lambda:=\Lambda_1\Lambda_2^{-1}\in\pmb{\SO}_g(\mathbb K)$ we compute 
$$D_1:=\Lambda Q_1\Lambda^{\t}=\Lambda Q_1\Lambda^{-1}=\Lambda_1\Lambda_2^{-1}Q_1\Lambda_2\Lambda^{-1}=\Lambda_1D_2\Lambda_1^{-1}\in \Z_g.$$
Thus the pair $(Q_0,Q_1)$ is $\pmb{\SL}_g$-equivalent to $\Lambda\cdot (Q_0,Q_1)=(\lambda 1_g,D_1)$.\endproof

\medskip
We recall that if 
$$\Theta_0,\ldots,\Theta_g\in\mathcal O(\pmb{\V}_g^2)^{\pmb{\SL}_g}=\mathcal O(\pmb{\V}_g^2/\pmb{\SL}_g)$$
 are the partial polarizations of the invariant polynomial $\det$ defined by the identity
\begin{equation}\label{EQ149.25}
\det(y_0Q_0+y_1Q_1)=\sum_{i=0}^g y_0^{g-i}y_1^i\Theta_i(Q_0,Q_1),
\end{equation}
where $y_0,y_1$ are two indeterminates, 
then the 
Cayley's {\it tact-invariant} (see \cite{Dol}, Subsect. 2.3.2, Thm. 2.3.12) is
$$J:=\textup{Disc}_g(\Theta_0,\ldots,\Theta_g)\in\mathcal O(\pmb{\V}_g^2)^{\pmb{\SL}_g}=\mathcal O(\pmb{\V}_g^2/\pmb{\SL}_g),$$
where $\textup{Disc}_g$ is the discriminant of the generic binary form of degree $g$; so $D_g$ is a homogeneous polynomial in $g+1$ indeterminates with coefficients in $\mathbb K$ and of degree $2g-2$. 

\medskip\noindent
 {\it Proof of Theorem \ref{T27}.} 
If $p\neq 2$, let $F_1:=1$. If $p=2$, let $F_1:=\Theta_1$. For an integer $r\geq 1$ we consider
 the projection $\V_g^{r+1}\rightarrow\V_g^2$
onto the direct sum of the first two $\V_g$ copies of $\V_g^{r+1}$. Based on Facts \ref{F3} and \ref{F4}, to end the proof it suffices to prove that 
$$F_0:=\Theta_0\cdot F_1\cdot J\in\mathcal O(\pmb{\V}_g^2)^{\pmb{\SL}_g}=\mathcal O(\pmb{\V}_g^2/\pmb{\SL}_g)$$ 
is an $\pmb{\SL}_g$-separating polynomial on $\V_g^2$.

To check this, let
 $(Q_0,Q_1)\in\V_g^2$ be such that $F_0(Q_0,Q_1)\neq 0$, i.e., $\det(Q_0)\neq 0$ and $J(Q_0,Q_1)\neq 0$ and, in case $p=2$, $\Theta_1(Q_0,Q_1)\neq 0$. We have to show that $\Stab_{\pmb{\SL}_g}((Q_0,Q_1))$ is finite. If $p=2$, we also have $\Trace(Q_0^{-1}Q_1)=\det(Q_0)^{-1}\Theta_1(Q_0,Q_1)\neq 0$. The pair $(Q_0,Q_1)$ is $\pmb{\SL}_g$-equivalent to a pair $(\lambda 1_g,D_1)\in\L_g\subset \V_g^2$ with $\lambda\in\mathbb K^\times$ (see Fact \ref{F7}).

As $J(\lambda 1_g,D_1)=J(Q_0,Q_1)\neq 0$ and
\begin{equation}\label{EQ149.5}
\det(y_01_g+y_1\lambda^{-1}D_1)=\lambda^{-g}\Big(\sum_{i=0}^g y_0^{g-i}y_1^i\Theta_i(\lambda 1_g,D_1)\Big),
\end{equation}
we get that $\lambda^{-1}D_1$ has distinct eigenvalues. Thus $D_1$ has distinct diagonal entries and hence its centralizer in $\pmb{\SL}_g$ is $\pmb{\T}_g$. Therefore we have $\Stab_{\pmb{\SO}_g}(D_1)=\pmb{\SO}_g\cap \pmb{\T}_g$. 
 As $\Stab_{\pmb{\SL}_g}((Q_0,Q_1))$ is $\pmb{\SL}_g(\mathbb K)$-conjugate to 
 $$\Stab_{\pmb{\SL}_g}((\lambda 1_g,D_1))=\Stab_{\Stab_{\pmb{\SL}_g}(\lambda 1_g)}(D_1)=\Stab_{\pmb{\SO}_g}(D_1)=\pmb{\SO}_g\cap \pmb{\T}_g\simeq\pmb{\mu}_{2,\mathbb K}^{g-1},$$
 we get that $\Stab_{\pmb{\SL}_g}((Q_0,Q_1))$ is a finite group scheme over $\Spec(\mathbb K)$.\endproof
 
 
\subsubsection{Proof of Theorem \ref{T26}}\label{S414}

To prove Theorem \ref{T26} we need a few preparations. Let $e_1:=[1\; 0\; 0\;\cdots 0]^{\t}$, $\ldots,e_g:=[0\; 0\; 0\;\cdots 1]^{\t}$ be the standard ordered basis of $\mathbb K^g$. For $M\in\Mat_g$ and $q\in\{1,\ldots,g-1\}$
let $M^{\wedge q}\in\Mat_{\binom{g}{q}}$ denote the matrix of the $q$-th wedge power of the endomorphism of $\mathbb K^g$ defined by $M$ with respect to the standard basis $\{e_{i_1}\wedge\cdots\wedge e_{i_q}|1\leq i_1<\cdots<i_q\leq g\}$ ordered lexicographically. For $D:=\Diag(\lambda_1,\ldots,\lambda_g)$, $D^{\wedge q}\in \Z_{\binom{g}{q}}$ has diagonal entries $\lambda_{i_1}\cdots\lambda_{i_q}$s with $1\leq i_1<\cdots <i_q\leq g$. For $M_0,M_1\in\Mat_g$ we have $(M_0M_1)^{\wedge q}=M_0^{\wedge q}M_1^{\wedge q}$ and we set
$$\Phi_q(M_0,M_1):=\det(M_0^{\wedge q}M_1^{\wedge q}-M_1^{\wedge q}M_0^{\wedge q}).$$ 

We recall that $\pmb{\SL}_g$ acts on $\Mat^2_g$ by conjugation. 

\begin{lemma}\label{L36}
For every $q\in\{1,\ldots,g-1\}$, $\Phi_q(M_0,M_1)$ is a non-zero ${\pmb{\SL}}_g$-invariant homogeneous polynomial in the entries of $M_0,M_1$.
\end{lemma}

\noindent
{\it Proof.}
The only non-obvious thing is that $\Phi_q\neq 0$. To check that $\Phi_q\neq 0$, it suffices to find $\Lambda_0,\Lambda_1\in\pmb{\GL}_g(\mathbb K)$ such that 
$$\det(\Lambda_0^{\wedge q}-(\Lambda_1\Lambda_0\Lambda_1^{-1})^{\wedge q})\neq 0.$$
We will take $\Lambda_0\in \Z_g$ such that $\Lambda_0^{\wedge q}$ has distinct non-zero eigenvalues and $\Lambda_1$ to be the permutation matrix attached to the $g$-cycle $(1\cdots g)$. With these choices, 
$(\Lambda_1\Lambda_0\Lambda_1^{-1})^{\wedge q}$ becomes a diagonal matrix obtained from $\Lambda_0^{\wedge q}$ by permuting its diagonal entries via a permutation without fixed points, hence $\Lambda_0^{\wedge q}-(\Lambda_1\Lambda_0\Lambda_1^{-1})^{\wedge q}\in \Z_{\binom{g}{q}}$ has non-zero entries on the diagonal and thus is invertible. To give a concrete example of such a $\Lambda_0$, let $N\geq 2^g$ be an integer not divisible by $p$, let $\zeta$ be a primitive $N$-th root of unity in $\mathbb K$, and choose
$$\;\;\;\;\;\;\;\;\;\;\;\;\;\;\;\;\;\;\;\;\;\;\;\;\;\;\;\Lambda_0:=\Diag(\zeta^{2^0}, \zeta^{2^1}, \zeta^{2^2},\ldots,\zeta^{2^{g-1}})\in\Z_g.\;\;\;\;\;\;\;\;\;\;\;\;\;\;\;\;\;\;\;\;\;\;\;\;\;\;\;\,\hfill{\square}$$

\begin{lemma}\label{L37}
Let $q\in\{1,\ldots,g-1\}$ and let $M_0,M_1\in\Mat_g$ be such that they have a common invariant subspace $W$ of dimension $q$. Then $\Phi_q(M_0,M_1)=0$.
\end{lemma}

\noindent
{\it Proof.}
As each non-zero vector $w_q$ of the line $\bigwedge^{\raisebox{-0.0ex}{\scriptsize $q$}}(W)$ of $\bigwedge^{\raisebox{-0.0ex}{\scriptsize $q$}}(\mathbb K^g)=\mathbb K^{\binom{g}{q}}$ 
is an eigenvector of both $M_0^{\wedge q}$ and $M_1^{\wedge q}$, we have
$$(M_0^{\wedge q}M_1^{\wedge q}-M_1^{\wedge q}M_0^{\wedge q})w_q=0.$$
This implies that $\Phi_q(M_0,M_1)=0$.\endproof

\medskip\noindent
 {\it Proof of Theorem \ref{T26}.} The $\pmb{\SL}_g$-invariant homogeneous polynomial
 $$\Phi(M_0,M_1):=\prod_{q=1}^{g-1} \Phi_q(M_0,M_1)\in\mathcal O(\pmb{\Mat}_g^2)$$
is non-zero (see Lemma \ref{L36}). Similar to the proof of Theorem \ref{T25}, by Facts \ref{F3} and \ref{F4}, to prove Theorem \ref{T26} it suffices to prove that $\Phi$ is an ${\pmb{\SL}}_g$-separating polynomial on $\Mat_g^2$. 
Let $(M_0,M_1)\in\Mat_g^2$ with $\Phi(M_0,M_1)\neq 0$. We know that $M_0$ and $M_1$ do not have a common non-trivial invariant subspace, see Lemma \ref{L37}. If $\Lambda\in\Stab_{\pmb{\SL}_g}((M_0,M_1))(\mathbb K)$, then each eigenspace of $\Lambda$, being $M_0$-invariant and $M_1$-invariant, must be $\mathbb K^g$. This implies that $\Stab_{\pmb{\SL}_g}((M_0,M_1))(\mathbb K)=Z(\pmb{\SL}_g)(\mathbb K)\simeq\pmb{\mu}_g(\mathbb K)$ is finite. Thus $\Phi$ is indeed an ${\pmb{\SL}}_g$-separating polynomial on $\Mat_g^2$.\endproof

\subsubsection{Proof of Theorem \ref{T27}}\label{S415}

We now know that $(\pmb{\Mat}_g^2)^{\s}$ is a non-empty open subscheme of $\pmb{\Mat}_g^2$ but to prove Theorem \ref{T27} we need that in fact the closed subscheme $\pmb{\V}_g^2$ and the two open subschemes $\pmb{\GL}_g^2$ and $(\pmb{\Mat}_g^2)^{\s}$ of $\pmb{\Mat}_g^2$ have a non-empty intersection.

\begin{lemma}\label{L38}
We have $(\pmb{\GL}_g^2)^{\s}\cap(\pmb{\Z}_g\times\pmb{\V}_g)\neq \emptyset$.
\end{lemma}

\noindent
{\it Proof.} If $p=0$ let $\mathbb K_0$ be the algebraic closure of $\mathbb Q$ in $ \mathbb K$ and let $\mathbb K_1:=\mathbb C$. If $p>0$ let $\mathbb K_0$ be the algebraic closure of $\mathbb F_p$ in $\mathbb K$ and let $\mathbb K_1$ be the completion of an algebraic closure of $\mathbb K_0((x))$ with respect to its standard topology (here $x$ is an indeterminate). The lemma holds for $\mathbb K$ if and only if it holds for $\mathbb K_0$ and if and only if it holds for $\mathbb K_1$. Thus we can assume that $\mathbb K=\mathbb K_1$ is complete with respect to an absolute value. The set of $\mathbb K$-valued points of a variety over $\mathbb K$ has a natural topology induced from that on $\mathbb K$ to be referred as the metric topology.

For $\nu_1,\ldots,\nu_g\in\mathbb K^{\times}$ distinct let $D_0:=\Diag(\nu_1,\ldots,\nu_g)\in\pmb{\GL}_g(\mathbb K)$. Let $\mathbb E_g\in\Mat_g$ be the matrix whose $g^2$ entries are all $1$. For $\lambda\in\mathbb K$ the condition that $\lambda\mathbb E_g+D_0$ is an invertible matrix with $g$ distinct eigenvalues is expressed by $\lambda$ not being a zero of a polynomial $\chi(x)\in\mathbb K[x]$. As $\chi(0)\neq 0$, there exists $a\in\mathbb K^\times$ such that $\chi(a)\neq 0$. This implies that the matrix $Q_1:=\mathbb E_g+a^{-1}D_0$ is invertible and has $g$ distinct eigenvalues. 

Let $\overline{\bf o}$ be the closure in $\pmb{\Mat}_g^2$ of the $\pmb{\SL}_g$-orbit ${\bf o}$ of $(D_0,Q_1)$. We have $(D_0,Q_1)\in \pmb{\GL}_g^2(\mathbb K)\cap (\Z_g\times\V_g)$. As $\Stab_{\pmb{\SL}_g}(D_0)=\pmb{\T}_g$ we compute
$$\Stab_{\pmb{\SL}_g}((D_0,Q_1))=\Stab_{\pmb{\T}_g}(Q_1)=\Stab_{\pmb{\T}_g}(\mathbb E_g)=Z(\pmb{\SL}_g)\simeq \pmb{\mu}_{g,\mathbb K}.$$
Thus $(D_0,Q_1)\in (\Mat_g^2)^{\s}$ if and only if $\overline{\bf o}={\bf o}$. Hence to prove the lemma it suffices to check that $\overline{\bf o}={\bf o}$. Let $(Q_0,Q'_1)\in\overline{\bf o}(\mathbb K)$. We check that there exists $\Lambda\in\pmb{\SL}_g(\mathbb K)$ such that $(\Lambda D_0\Lambda^{-1},\Lambda Q_1\Lambda^{-1})=(Q_0,Q'_1)$. The characteristic polynomials of $Q_0$ and $D_0$ are the same and hence $Q_0$ is similar to $D_0$. Thus to prove the existence of $\Lambda$ we can assume that $Q_0=D_0$.

Let $(\Lambda_n)_{n\geq 1}$ be a sequence in $\pmb{\SL}_g({\mathbb K})$ such that
$(\Lambda_nD_0\Lambda_n^{-1},\Lambda_nQ_1\Lambda_n^{-1})$ converges in the metric topology to $(D_0,Q_1')$. Its existence can be checked by considering an irreducible closed affine curve $C\subset\overline{\bf o}$ which passes
through $(D_0,Q'_1)$ and is such that $C\cap (\overline{\bf o}\setminus {\bf o})$ is finite and using the fact that the complement of finitely many points
in $\tilde{C}(\mathbb K)$, where $\tilde{C}$ is the normalization of $C$, is dense in $\tilde{C}(\mathbb K)$ in the metric topology.

If $\pmb{\U}_g^+$ (resp. $\pmb{\U}_g^-$) is the unipotent subgroup of $\pmb{\SL}_g$ of upper (resp. lower) triangular matrices with all diagonal entries $1$, then the product morphism $\pmb{\U}_g^-\times\pmb{\U}_g^+\times \pmb{\T}_g\rightarrow \pmb{\SL}_g$ is an open embedding (by the standard Bruhat decomposition, see \cite{Bore}, Ch. IV, Sect. 14, Subsect. 14.12, Thm. (a)) and thus we can identify $\pmb{\U}_g^-\times\pmb{\U}_g^+$ with an open subscheme of the quotient scheme $\pmb{\Q}_q:=\pmb{\SL}_g/\pmb{\T}_g$. As $\Lambda_nD_0\Lambda_n^{-1}$ converges to $D_0$, $\Lambda_n\pmb{\T}_g\in \pmb{\Q}_g(\mathbb K)$ converges to $\pmb{\T}_g\in \pmb{\Q}_g(\mathbb K)$. Thus there exists $N\in\mathbb N$ such that for $n\geq N$ there exists $\Sigma_n\in (\pmb{\U}_g^-\times\pmb{\U}_g^+)(\mathbb K)\subset\pmb{\SL}_g(\mathbb K)$ with the properties that we have $\Lambda_n\pmb{\T}_g=\Sigma_n\pmb{\T}_g$ and the sequence $(\Sigma_n)_{n\geq N}$ converges to $1_g$. For $n\geq N$, let $D_n:=\Sigma_n^{-1}\Lambda_n$; we have $D_n\in\pmb{\T}_g(\mathbb K)$. As $(\Sigma_nD_nD_0D_n^{-1}\Sigma_n^{-1},\Sigma_nD_nQ_1D_n^{-1}\Sigma_n^{-1})_{n\geq N}$ converges to $(D_0,Q'_1)$ and both $(\Sigma_n)_{n\geq N}$ and $(\Sigma_n^{-1})_{n\geq N}$ converge to $1_g$, the sequence $(D_nD_0D_n^{-1},D_nQ_1D_n^{-1})_{n\geq N}$ is equal to the sequence $(D_0,D_nQ_1D_n^{-1})_{n\geq N}$ and converges to $(D_0,Q'_1)$. 

As $Z(\pmb{\GL}_g)$ centralizes $Q_1$ we can change the $D_n$s so that they are in $\pmb{\GT}_g(\mathbb K)$. Thus for $n\geq N$ there exists $D_n=\Diag(1,\lambda_{2,n},\ldots,\lambda_{g,n})\in\pmb{\GT}_g(\mathbb K)$ such that $(D_nQ_1D_n^{-1})_{n\geq N}$ converges to $Q'_1$.
For an integer $1<i\leq g$, the $i1$ entry of $D_nQ_1D_n^{-1}$ is $\lambda_{i,n}$ and thus
$(\lambda_{i,n})_{n\geq N}$ converges to the $i1$ entry $\lambda_i$ of $Q'_1$. Similarly,
$(\lambda_{i,n}^{-1})_{n\geq N}$ converges to the $1i$ entry of $Q'_1$. The last two sentences imply that $\lambda_i\neq 0$.
We conclude that $(D_n)_{n\geq N}$ converges to $D:=(1,\lambda_2,\ldots,\lambda_g)\in\pmb{\GT}_g(\mathbb K)$ and therefore $DQ_1D^{-1}=Q'_1$.
For $\Lambda:=\sqrt[g]{\det(D^{-1})}D\in \pmb{\T}_g(\mathbb K)$, we have
$(\Lambda D_0\Lambda^{-1},\Lambda Q_1\Lambda^{-1})=(D_0,Q_0')$. This implies that $\overline{\bf o}={\bf o}$. \endproof

\medskip

For the next two lemmas we recall the morphism $\pmb{\pi}_n$ in Equation (\ref{EQ148}).

\begin{lemma}\label{L39}
We assume that $n\geq 2$. Then there exist invertible symmetric matrices $Q_0,\ldots,Q_n\in\V_g$ such that we have $\pi_n(Q_0,\ldots,Q_n)\in(\pmb{\GL}_g^n)^{\s}(\mathbb K)$.
\end{lemma}

\noindent
{\it Proof.} By Lemma \ref{L38} we can find invertible symmetric matrices $Q_0,Q_1'\in\V_g$ such that $(Q_0,Q_1')\in (\Mat^2_g)^{\s}$. Let $Q_1:=1_g$. Write $Q_1'=Q_2^*$, so $Q_2\in\V_g$. By Fact \ref{F3} there exists an $\pmb{\SL}_g$-separating homogeneous polynomial $F$ on $\Mat^2_g$ such that $F(Q_0,Q_1')=F(Q_0Q_1^*,Q_1Q_2^*)\neq 0$. The function $\Mat_g^n\rightarrow \mathbb K$ defined by
$(M_0,\ldots,M_{n-1})\mapsto F(M_0,M_1)$ is a separating polynomial (see Fact \ref{F4}) and does not vanish at a point of the form $\pi_n(Q_0,Q_1,Q_2,\ldots,Q_n)$ with $Q_3\ldots,Q_n\in\V_g$. Hence the latter points are stable by Fact \ref{F3} and thus belong to $(\pmb{\GL}_g^n)^{\s}(\mathbb K)$ provided $Q_3,\ldots,Q_n$ are invertible.
\endproof

\begin{lemma}\label{L40}
We assume that $n\geq 2$. There exists a Zariski dense open subscheme $\pmb{\W}_{g,n}$ of $\pmb{\GL}_g^{n+1}\cap\pmb{\V}_g^{n+1}$ such that for each point $(Q_0,\ldots,Q_n)$ in $\pmb{\W}_{g,n}(\mathbb K)$ the fiber
$\pmb{\pi}_n^{-1}(\pmb{\pi}_n(Q_1,\ldots,Q_{n+1}))$ has dimension $1$.
\end{lemma}

\noindent
{\it Proof.}
Let $(Q_0,\ldots,Q_n)\in\V_g^{n+1}$ and $\lambda\in \mathbb K^{\times}$. Let
$\lambda_n,\lambda_{n-1},\ldots,\lambda_0\in \mathbb K^{\times}$ be defined recursively by $\lambda_n:=\lambda$
and $\lambda_i=\lambda_{i+1}^{1-g}$ for $i\in\{n-1,n-2,\ldots,0\}$.
Then we have 
$$\pi_n(Q_0,\ldots,Q_n)=\pi_n(\lambda_0 Q_0,\ldots,\lambda_nQ_n);$$
hence $\pmb{\pi}_n^{-1}(\pmb{\pi}_n(Q_0,\ldots,Q_n))$ has dimension $\geq 1$ if $(Q_0,\ldots,Q_n)$ is not the zero vector of $\V_g^{n+1}$ and it is easy to see that the same holds even if $(Q_0,\ldots,Q_n)$ is the zero vector of $\V_g^{n+1}$.

We claim that there exists $(Q_0,\ldots,Q_n)\in\pmb{\GL}_g^{n+1}(\mathbb K)\cap\V_g^{n+1}$ such that
$\pmb{\pi}_n^{-1}(\pmb{\pi}_n(Q_0,\ldots,Q_n))$ has dimension $\leq 1$ and thus has dimension $1$; as the dimension of fibers of $\pmb{\pi}_n$ is upper semicontinuous on the source (see \cite{E}, Thm. 14.8), this will hold for 
$\mathbb K$-valued points $(Q_0,\ldots,Q_n)$ of a Zariski dense open subscheme $\pmb{\W}_{g,n}$ of $\pmb{\GL}_g^{n+1}\cap\pmb{\V}_g^{n+1}$.

To check our claim, recall from Lemma \ref{L38} that $(\pmb{\GL}_g^2)^{\s}\cap\pmb{\V}_g^2\neq \emptyset$. Let 
 $Q_0:=1_g$, let $(Q_1,Q_2)\in (\pmb{\GL}_g^2)^{\s}(\mathbb K)\cap\V_g^2$, and, if $n\geq 3$, let $(Q_3,\ldots,Q_n)\in (\pmb{\GL}_g^{n-2})(\mathbb K)\cap\V_g^{n-2}$. Let $\widetilde{Q}_0,\ldots,\widetilde{Q}_n\in\V_g$ be such that
$$\pi_n(\widetilde{Q}_0,\ldots,\widetilde{Q}_n)=\pi_n(Q_0,\ldots,Q_n).$$ 
We have
\begin{equation}\label{EQ150}
\widetilde{Q}_i\widetilde{Q}_{i+1}^*=Q_iQ_{i+1}^*\;\;\; \forall i\in\{0,\ldots,n-1\}.
\end{equation}
It follows that $\widetilde{Q}_0,\ldots,\widetilde{Q}_{n-1}$ are invertible and thus there exist matrices $\Lambda_0,\ldots,\Lambda_n\in \pmb{\GL}_g(\mathbb K)$ such that we have $\widetilde{Q}_0=Q_0\Lambda_0=\Lambda_0$, $\widetilde{Q}_1=Q_1\Lambda_1,\ldots$, $\widetilde{Q}_n=Q_n\Lambda_n$. From Equation (\ref{EQ150}) we get that
\begin{equation}\label{EQ151}
\Lambda_i\Lambda_{i+1}^*=1_g\;\;\;\forall i\in\{0,\ldots,n-1\}.
\end{equation}
As $\widetilde{Q}_0=\Lambda_0$ is symmetric, by induction on $i\in\{0,\ldots,n\}$ we get from Equation (\ref{EQ151}) that $\Lambda_i\in\V_g$.
As $\widetilde{Q}_1,\widetilde{Q}_2\in\V_g$, we have $Q_1\Lambda_1=(Q_1\Lambda_1)^{\t}=\Lambda_1Q_1$ and $Q_2\Lambda_2=(Q_2\Lambda_2)^{\t}=\Lambda_2Q_2$. Thus $\Lambda_1$ commutes with $Q_1$ and $\Lambda_2$ commutes with $Q_2$.
 Hence $\Lambda_2^{-1}$ commutes with $Q_2$ and hence so does $\Lambda_1$. Therefore $\Lambda_1$ is a $\mathbb K$-valued point of the one dimensional subvariety $Z(\pmb{\GL}_g)\Stab_{\pmb{\SL}_g}((Q_1,Q_2))$ of $\pmb{\GL}_g$ and thus, as $\Lambda_0,\Lambda_2,\ldots,\Lambda_n$ are determined by $\Lambda_1$ up to a finite number of possibilities (see Equation (\ref{EQ151})), our claim is proved. 
\endproof

\begin{rem}\label{R37}
The morphism $\pmb{\pi}_1:\pmb{\V}_g^2\rightarrow \pmb{\Mat}_g$ is dominant and thus for a generic point $(Q_0,Q_1)\in\V_g^2$, $\pmb{\pi}_1^{-1}(\pmb{\pi}_1(Q_0,Q_1))$ has dimension $g$. Even more is true, namely: the morphism $\pmb{\V}_g^2\rightarrow \pmb{\Mat}_g\times_{\Spec(\mathbb K)} \mathbb A^1_{\mathbb K}$ defined on $\mathbb K$-valued points by the rule $(Q_0,Q_1)\mapsto (Q_0Q_1^*,\det(Q_0))$ is dominant.
To check this it suffices to show that for each pair $(\Lambda,\lambda)\in \pmb{\GL}_g(\mathbb K)\times \mathbb K^\times$ with $\Lambda$ diagonalizable there exist $Q_0,Q_1\in\V_g$ such that $Q_0Q_1^*=\Lambda$ and $\det(Q_0)=\lambda$. We write 
 $\Lambda=\Lambda_1D\Lambda_1^{-1}$, where $D\in \Z_g$ is invertible and $\det(\Lambda_1)=1$. Let $\nu:=\sqrt[\leftroot{-2}\uproot{2}g]{\lambda}$. Then we can write $\Lambda=Q_0Q_2$, where $Q_0:=\nu\cdot \Lambda_1\Lambda_1^{\t}$, $Q_2:=\nu^{-1}\cdot (\Lambda_1^{-1})^{\t}D\Lambda_1^{-1}$. 
Clearly $Q_0,Q_2\in\V_g$ are invertible and $\det(Q_0)=\nu^g=\lambda$. We can write $Q_2=Q_1^*$, where $Q_1\in\V_g$ (as it is a scalar multiple of $Q_2^{-1}$), so $Q_0Q_1^*=\Lambda$.\footnote{We have $\Mat_g=(\V_g\cap\pmb{\GL}_g(\mathbb K))\cdot\V_g$ (the proof of \cite{Bos}, Thm. 1 works for all square matrices with entries in a ring which are similar to a Jordan matrix).}\end{rem}

\medskip\noindent
 {\it Proof of Theorem \ref{T27}.} 
Consider the morphisms
$$\pmb{\V}_g^{n+1}
\stackrel{\pmb{\pi}_n}{\longrightarrow} \pmb{\Mat}_g^n \stackrel{\varpi_n}{\longrightarrow} \pmb{\Mat}_g^n/\pmb{\SL}_g,$$
the Zariski closure $\pmb{\X}_n$ of $\pmb{\pi}_n(\pmb{\V}_g^{n+1})$ in $\pmb{\Mat}_g^n$, and the Zariski closure $\pmb{\Y}_n$ of $\varpi_n(\pmb{\X}_n)$ (or of $(\varpi_n\circ\pmb{\pi}_n)(\pmb{\V}_g^{n+1})$) in $\pmb{\Mat}_g^n/\pmb{\SL}_g$.
 By Lemma \ref{L40} we have
 $$\dim(\pmb{\X}_n)= (n+1)\frac{g(g+1)}{2}-1.$$
 By Lemma \ref{L39}, $\pmb{\X}_n\cap (\pmb{\Mat}_g^n)^{\s}\neq \emptyset$. 
As $\pmb{\X}_n$ is $\pmb{\SL}_g$-invariant and all fibers
 of $(\pmb{\Mat}_g^n)^{\s} \stackrel{\varpi_n}{\longrightarrow} (\pmb{\Mat}_g^n)^{\s}/\pmb{\SL}_g$
 are $\pmb{\SL}_g$-orbits of
 dimension $g^2-1$, we get that 
 $$\dim(\pmb{\Y}_n)= \dim(\pmb{\X}_n)-(g^2-1)= (n+1)\frac{g(g+1)}{2}-g^2
 =(n-1)\frac{g(g+1)}{2}+g.$$
 This equality implies that Theorem \ref{T27} holds.
\endproof

\subsubsection{On `Main Theorems'}\label{S416}

We recall the `first main theorem' for the conjugation action of ${\pmb{\SL}}_g$ on $\pmb{\Mat}^{r+1}_g$ (conjectured by Artin, proved for $p=0$ in
 \cite{G}, Thm. 17.4 and \cite{P}, Thm. 1.3, and for $p>0$ in \cite{Don92, Don93}; cf. also \cite{DCP}, Thm. 1.10). 
For a $\mathbb K$-algebra $\mathbb B$ and $M\in\pmb{\Mat}_g(\mathbb B)$, the characteristic polynomial
 $$\chi_M(x):=\det(x\cdot 1_g-M)=\sum_{j=0}^g (-1)^jc_j(M)x^{g-j}\in\mathbb B[x]$$
satisfies $c_0(M)=1$, $c_1(M)=\Trace(M)$, $c_g(M)=\det(M)$.
We identify 
$$\Sym((\Mat_g^{r+1})^*)=\mathbb K[X_0,\ldots,X_r],$$
where $\mathbb K[X_0,\ldots,X_r]:=\mathbb K[X_{l,ij}|0\leq l\leq r,\ 1\leq i,j\leq g]$
with each $X_l=(X_{l,ij})_{1\leq i,j\leq g}$ a $g\times g$ matrix with indeterminates as entries. 
Here we view every $X_{l,ij}$ as the element of $(\Mat_g^{r+1})^*$ that sends each $(M_0,\ldots,M_r)\in\Mat_g^{r+1}$ to the $ij$ entry of $M_l$. 

As the identifications $X_{l,ij}=X_{l,ij}(X_0,\ldots,X_r)$ make sense for all $1\leq i, j\leq g$ and all $l\in\{0,\ldots,r\}$, for $F\in \Sym((\Mat_g^{r+1})^*)$ we can write 
$$F=F(X_0,\ldots,X_r)=F(X_{0,11},X_{0,12},\ldots,X_{r,gg})$$
and we can similarly evaluate $F$ at elements of $\Mat^{r+1}_g\otimes_{\mathbb K} \Sym((\Mat_g^{r+1})^*)$. 

The action of $\Lambda\in\pmb{\SL}_g(\mathbb K)$ on $F$ is given by the rule 
 $$(\Lambda\cdot F)(X_0,\ldots,X_r):=F(\Lambda^{-1} X_0\Lambda,\ldots,\Lambda^{-1} X_r\Lambda).$$
 To check this we can assume $r=0$ and $F=X_{0,ij}$. We need to check that $\Lambda \cdot X_{0,ij}$ is the $ij$ entry of $\Lambda^{-1}X_0\Lambda$. But
 this follows from the fact that for all $M\in \Mat_g$, $(\Lambda \cdot X_{0,ij})(M)$ is the $ij$ entry of the matrix $\Lambda^{-1}M\Lambda$, confirm the very definition of the duals of $\pmb{\SL}_g$-modules. 

 The `first main theorem' is the following:

\begin{thm}\label{T28}
The $\mathbb K$-algebra $\Sym((\Mat_g^{r+1})^*)^{\pmb{\SL}_g}$ is generated by the following set of homogeneous polynomials:
$$\{c_j(X_{l_1}\cdots X_{l_N})|j\in \{1,\ldots,g\},N\in\mathbb N,l_1,\ldots,l_N\in\{0,\ldots, r\}\}.$$
 If $p\nmid g!$, then the $\mathbb K$-algebra $\Sym((\Mat_g^{r+1})^*)^{\pmb{\SL}_g}$ is generated by the subset defined by $j=1$.\footnote{If $p=0$, one can take $N\in\{1,\ldots,2^g-1\}$ (see \cite{P}, Thm. 3.4 (a)). Also, there exists a `second main theorem' (see \cite{P}, Thm. 4.5 (a) for $p=0$ and \cite{DCP}, Introd. for $p>0$) which describes the relations among these generators but it will not be used in the paper.}
\end{thm}

\begin{rem}\label{R38}
We recall the following general result (see \cite{MFK}, App. to Ch. I, Sect. F, Thm.)
If $G$ is semisimple and acts on $\pmb{\X}=\Spec(\mathbb B)$ with $\mathbb B$ a finitely generated $\mathbb K$-algebra which is a unique factorization domain, then the following two statements are equivalent:

\medskip\noindent
\circled{\textup{s1}} There exists a non-empty $G$-invariant open subscheme $\pmb{\Y}$ of $\pmb{\X}$ such that for each $x\in\pmb{\Y}(\mathbb K)$, the orbit of $x$ is closed in $\pmb{\X}$.

\smallskip\noindent
\circled{\textup{s2}} There exists a non-empty $G$-invariant open subscheme $\pmb{\S}$ of $\pmb{\X}$ such that for each $x\in\pmb{\S}(\mathbb K)$, the reduced stabilizer of $x$ in $G$ is reductive.

\medskip\noindent 
This result can be used to get less explicit proofs of Theorems \ref{T25} and \ref{T26}. Our explicit proofs of Theorems \ref{T25} and \ref{T26} were used, however, to prove Theorems \ref{T10} and \ref{T27} and will be repeatedly used in what follows. 
\end{rem}

\subsection{General results on quadratics for arbitrary $g\geq 1$}\label{S42}

We return to the invariant theory of quadratic forms. 

\subsubsection{The $\pmb{\SL}_g$-module $\V_g$}\label{S4201}

Let $\varepsilon_g\in\{0,1\}$ be the remainder of the division of $g$ by $2$. As we have $\Ker(\pmb{\SL}_g\to\pmb{\GL}_{V_g})\simeq\pmb{\mu}_{2,\mathbb K}^{1-\varepsilon_g}$, we will use $\mathbb W_{\varepsilon_g}$, where 
$$\mathbb W_0:=\frac{1}{2}\mathbb Z_{\geq 0}\;\;\;\textup{and}\;\;\;\mathbb W_1:=\mathbb Z_{\geq 0},$$
and $\pmb{\SL}'_g:=\pmb{\SL}_g/\Ker(\pmb{\SL}_g\to\pmb{\GL}_{V_g})$. If $g$ is odd we have $\pmb{\SL}'_g=\pmb{\SL}_g$ and if $g$ is even we have a short exact sequence $1\to\pmb{\mu}_{2,\mathbb K}\to\pmb{\SL}_g\to\pmb{\SL}'_g\to 1$.

Let $\tau_g:\pmb{\SL}_g\simeq\pmb{\SL}_g$ be the involutory automorphism which on $\mathbb K$-valued points is given by the rule: $\Lambda\mapsto (\Lambda^{\t})^{-1}$. If $g\geq 3$, then $\tau_g$ is the `standard' outer automorphism of $\pmb{\SL}_g$. For a finite dimensional $\pmb{\SL}_g$-module $V$ given by a representation $\rho_V:\pmb{\SL}_g\to\pmb{\GL}_{V}$, let $\tau_g^*(V)$ be the $\pmb{\SL}_g$-module given by the representation $\rho_V\circ\tau_g$; thus we have $\tau_g^*(V)=V$ as $\mathbb K$-vector spaces. If $g=2$, then $\tau_2$ is an inner automorphism of $\pmb{SL}_2$ and therefore $\tau_2^*(V)$ and $V$ are isomorphic $\pmb{\SL}_2$-modules.

We consider the standard $g$-dimensional $\pmb{\SL}_g$-module $\W_g:=\mathbb K^g$, where the action $\pmb{\SL}_g\times\pmb{\W}_g\rightarrow\pmb{\W}_g$ is given on $\mathbb K$-valued points by the rule: $(\Lambda,v)\mapsto \Lambda\cdot v$. If $p>0$, we recall that $\W_g^{(p)}$ is the $\pmb{\SL}_g$-module given on $\mathbb K$-valued points by the rule: $(\Lambda,v)\mapsto \Lambda^{[p]}\cdot v$ where $\Lambda^{[p]}\in\pmb{\SL}_g(\mathbb K)$ is the matrix obtained from $\Lambda\in\pmb{\SL}_g(\mathbb K)$ by replacing each entry by its $p$-th power. 

Let $\{e_1^*,\ldots,e_g^*\}$ be the standard basis of $\W_g^*$: it is the dual of the standard basis $\{e_1,\ldots,e_g\}$ of $\W_g$. Using the Kronecker delta, for $i,j\in\{1,\ldots,g\}$ we have $e_j^*(e_i)=\delta_{ij}$. The standard basis of the $\mathbb K$-vector space 
$$\Mat_g=\Hom_{\mathbb K}(\W_g,\W_g)=\W_g\otimes_{\mathbb K} \W_g^*$$
is $\{e_i\otimes e_j^*|1\leq i,j\leq g\}$. For $m\in\{1,\ldots,g\}$ we have $(e_i\otimes e_j^*)(e_m)=\delta_{jm}e_i$. The standard basis of $\V_g$ is 
$$\{e_i\otimes e_i^*|1\leq i\leq g\}\cup\{e_i\otimes e_j^*+e_j\otimes e_i^*|1\leq i<j\leq g\}.$$

For $g\geq 2$ and $i,j\in\{1,\ldots,g\}$ distinct let 
\begin{equation}\label{EQ151.3}
\Lambda_{i,j}\in\pmb{\SO}_g(\mathbb K)\end{equation} 
be such that it maps $(e_i,e_j)$ to $(e_j,-e_i)$ and fixes each $e_{i'}$ with $i'\in\{1,\ldots,g\}\setminus\{i,j\}$. If $p=2$, we have $\Lambda_{i,j}=\Lambda_{j,i}$.

The classification of finite dimensional simple modules over a simply connected semisimple group $G$ over $\mathbb K$ in terms of weights associated to a maximal torus $T_G$ of $G$ (see \cite{Ja}, part II, CG. 2 and 5) implies: 

\medskip\noindent
 {\it ($\sharp$)} If $V$ and $V'$ are two finite dimensional $G$-modules which are isomorphic as $T_G$-modules, then they are isomorphic as $G$-modules provided either both are semisimple $G$-modules or one is a simple $G$-module. 

\medskip
If $p=2$, let $\V_g^{\alt}$ be the $\mathbb K$-vector subspace of $\V_g$ of alternating symmetric matrices and let $\V_g^{\Diag}:=\V_g/\V_g^{\alt}$; we have 
$$\dim_{\mathbb K}(\V_g^{\alt})=\frac{g(g-1)}{2}\;\;\;\textup{and}\;\;\;\dim_{\mathbb K}(\V_g^{\Diag})=g.$$
The $\mathbb K$-vector subspace $\bigtriangleup_g\subset \Mat_g$ of upper triangular matrices is a direct supplement of $\V_g^{\alt}$ in $\Mat_g$ and we have a canonical $\mathbb K$-linear isomorphism
\begin{equation}\label{EQ151.4}\bigtriangleup_g\to\Mat_g/V_g^{\alt}.
\end{equation}

The $\mathbb K$-linear map $\Mat_g\to\V_g$ defined by the rule $M\mapsto M+M^{\t}$ is surjective if and only if $p\neq 2$: if $p=2$, the image is $\V_g^\alt$. 

\begin{lemma}\label{L41}
If $p=2$, then the following ten properties hold:

\medskip
{\bf (a)} We have a short exact sequence of $\pmb{\SL}_g$-modules
\begin{equation}\label{EQ151.5}
0\rightarrow\V_g^{\alt}\rightarrow\V_g\rightarrow\V_g^{\Diag}\rightarrow 0.
\end{equation}

{\bf (b)} The $\pmb{\SL}_g$-module $\V_g^{\Diag}$ is simple and isomorphic to $\W_g^{(2)}$.

\smallskip
{\bf (c)} The $\pmb{\SL}_g$-module $\V_g^{\alt}$ is simple (if $g=2$, then $\V_2^{\alt}=\mathbb K\mathcal J$ is a trivial simple $\pmb{\SL}_2$-module).

\smallskip
{\bf (d)} The short exact sequence (\ref{EQ151.5}) does not split.

\smallskip
{\bf (e)} The $\pmb{\SL}_g$-modules $\V_g^{\alt}$, $\bigwedge^{\raisebox{-0.0ex}{\scriptsize $2$}}(\W_g)$ and $\Mat_g/\V_g$ are isomorphic.

\smallskip
{\bf (f)} The $\pmb{\SL}_g$-module $\V_g$ is not isomorphic to either $\V_g^*$ or $\tau_g^*(\V_g^*)$.

\smallskip
{\bf (g)} The $\pmb{\SL}_g$-modules $\V_g$ and $\Sym^2(\W_g)$ are not isomorphic. 

\smallskip
{\bf (h)} The only $\pmb{\SL}_g$-submodules of $\Mat_g$ are $0$, $\V_g^{\alt}$, $\V_g$ and $\Mat_g$. 

\smallskip
{\bf (i)} The $\pmb{\SL}_g$-modules $\Mat_g^*$ and $\tau_g^*(\Mat_g)$ are isomorphic. 

\smallskip
{\bf (j)} The $\pmb{\SL}_g$-modules $\Mat_g/\V_g^{\alt}$ and $\tau_g^*(\V_g^*)$ are isomorphic.
\end{lemma}

\noindent
{\it Proof.} Part (a) holds as from \cite{A}, Thm. 1 we get that $\V_g^{\alt}$ is an $\pmb{\SL}_g$-submodule of $\V_g$. Part (b) follows from ($\sharp$) applied with $G=\pmb{\SL}_g$ as we have a natural identification $\V_g^{\Diag}=\W_g^{(2)}$ as $\pmb{\T}_g$-modules. 

To prove parts (c) and (d), we first remark that the $\pmb{\SL}_g$-submodule $\V_g^M$ of $\Mat_g$ generated by a non-zero $M\in\Mat_g$ is the $\mathbb K$-span of the congruence class of $M$. As the congruence relation on $\V_g^{\alt}\setminus\{0\}$ has precisely $\lfloor \frac{g}{2}\rfloor$ congruence classes given by the (even) rank (see \cite{A}, Thm. 4), it is easy to see that for $Q\in V_g^{\alt}$ we have $\V_g^Q=\V_g^{\alt}$. From this part (c) follows. If $Q\in\V_g$ is non-alternating, then $\V_g^Q\cap\Z_g\neq 0$ (see \cite{A}, Thm. 6) and using this we easily get that $\Z_g\subset\V_g^Q$ and thus (again based on loc. cit.) $\V_g\setminus\V_g^{\alt}\subset\V_g^Q$. This implies that $\V_g^Q=\V_g$ and that (d) holds.

Part (e) follows from ($\sharp)$ as the $\pmb{\T}_g$-modules $\V_g^{\alt}$, $\bigwedge^{\raisebox{-0.0ex}{\scriptsize $2$}}(\W_g)$ and $\Mat_g/\V_g$ are isomorphic. Part (f) follows from parts (a) to (d) and their duals. 

As we have an $\pmb{\SL}_g$-invariant epimorphism $\Sym^2(\W_g)\rightarrow \bigwedge^{\raisebox{-0.0ex}{\scriptsize $2$}}(\W_g)$, part (g) follows from parts (d) and (e). 

Based on parts (a) to (d) (or on the proof of parts (c) and (d)), to prove part (h) it suffices to show that if $M\in\Mat_g\setminus\V_g$, then $\V_g^M=\Mat_g$. Based on parts (c) and (e), we have an $\pmb{\SL}_g$-invariant epimorphism $\V_g^M\to\Mat_g/\V_g$. Thus we can assume that $M\in\V_g+\mathbb Ke_1\otimes e_2^*$ and using the decomposing of $\V_g+\mathbb Ke_1\otimes e_2^*$ into $\pmb{\T}_g$-modules on which $\pmb{\T}_g$ acts via one character, we can assume that in fact $M=\nu e_2\otimes e_1^*+(e_1\otimes e_2^*)$ with $\nu\in\mathbb K\setminus\{1\}$. As 
$$(1_g+e_2\otimes e_1^*)M(1_g+e_2\otimes e_1^*)^{\t}=(M+e_2\otimes e_2^*)(1_g+e_1\otimes e_2^*)=M+(1+\nu)(e_2\otimes e_2^*)$$ is in $\V_g^M$, we get that $(1+\nu)(e_2\otimes e_2^*)\in V_g^M\cap (\Z_g\setminus\{0\})$ and from the proof of parts (c) and (d) we get that $\V_g^M$ contains $\V_g$. From this, as $\V_g^M$ surjects onto $\Mat_g/\V_g$, we conclude that $\V_g^M=\Mat_g$. Thus part (h) holds.

Part (i) follows from the fact that the non-degenerate trace bilinear map 
$$\Mat_g\times \tau_g^*(\Mat_g)\to\mathbb K$$ 
is $\pmb{\SL}_g$-invariant, i.e., for $\Lambda\in\pmb{\SL}_g(\mathbb K)$ and $M_1,M_2\in\Mat_g$, we have 
$$\Trace(\Lambda M_1\Lambda^{\t}(\Lambda^{\t})^{-1}M_2\Lambda^{-1})=\Trace(\Lambda M_1M_2\Lambda^{-1})=\Trace(M_1M_2).$$

The dual $\Mat_g^*\to V_g^*$ of the $\pmb{\SL}_g$-invariant monomorphism $\V_g\to\Mat_g$ induces an $\pmb{\SL}_g$-invariant epimorphism $\tau_g^*(\Mat_g^*)\to \tau_g^*(V_g^*)$ which, based on part (i), can be viewed as an $\pmb{\SL}_g$-invariant epimorphism $\Mat_g\to \tau_g^*(V_g^*)$. Thus, based on part (h) and reasons of dimensions, we get that $\tau_g^*(V_g^*)$ is isomorphic to the quotient $\Mat_g/\V_g^{\alt}$ of $\Mat_g$. Thus part (j) holds. \endproof 

\begin{lemma}\label{L42}
If $p\neq 2$, then the following two properties hold:

\medskip
{\bf (a)} The $\pmb{\SL}_g$-modules $\V_g$ and $\Sym^2(\W_g)$ (resp. $\V_g^*$ and $\Sym^2(\W_g^*)$) are simple and isomorphic. 

\smallskip
{\bf (b)} The $\pmb{\SL}_g$-modules $\V_g$ and $\tau_g^*(\V_g^*)$ (or $\tau_g^*(\V_g)$ and $\V_g^*$) are isomorphic. 
\end{lemma}

\noindent
{\it Proof.} We only prove part (a) for $\V_g$ and $\Sym^2(\W_g)$ as the proof for $\V_g^*$ and $\Sym^2(\W_g^*)$ is the same. The $\pmb{\T}_g$-modules $\V_g$ and $\Sym^2(\W_g)$ are isomorphic under the natural $\mathbb K$-linear map that maps $e_i\otimes e_j^*+e_j\otimes e_i^*$ to the image of $e_i\otimes e_j\in\W_g\otimes_{\mathbb K} \W_g$ into $\Sym^2(\W_g)$ for all $1\leq i\leq j\leq g$. The proof of Lemma \ref{L41} (c) and (d) can be easily adapted to give that $\V_g$ is a simple $\pmb{\SL}_g$-module. From the last two sentences and ($\sharp$) we get that the $\pmb{\SL}_g$-modules $\V_g$ and $\Sym^2(\W_g)$ are isomorphic. Thus part (a) holds. Part (b) follows from part (a) as the $\pmb{\SL}_g$-modules $\W_g$ and $\tau_g^*(\W_g^*)$ are isomorphic.\endproof

\subsubsection{Rings and vector spaces of invariants: the $\mathbb K$ case}\label{S4202} 
With arbitrary characteristic $p\geq 0$, we begin by giving an explicit description of the $\pmb{\SL}_g(\mathbb K)$-action on the $\mathbb K$-algebra $\Sym((\V_g^{r+1})^*)$.
Let $T_{ij}\in\V_g^*$ be defined by letting $T_{ij}(Q)$ be the $ij$ entry 
of $Q\in\V_g$. Hence the set $\{T_{ij}|1\leq i\leq j\leq g\}$ is a basis of $\V_g^*$. 
We write 
$$\Sym((\V_g^{r+1})^*)=\mathbb K[T,\ldots,T^{(r)}]:=\mathbb K[T_{ij},\ldots,T^{(r)}_{ij}|1\leq i\leq j\leq g],$$ where $T=T^{(0)}=(T_{ij})_{1\leq i,j\leq g}=(T^{(0)}_{ij})_{1\leq i,j\leq g},\ldots, T^{(r)}=(T^{(r)}_{ij})_{1\leq i,j\leq g}$ are symmetric matrices of size $g\times g$ of indeterminates and for $l\in\{0,\ldots,r\}$, $\{T_{ij}^{(l)}|1\leq i\leq j\leq g\}$ is $\{T_{ij}|1\leq i\leq j\leq g\}$ but viewed as a $\mathbb K$-basis of the $l$-th factor $\V_g^*$ of $(\V_g^{r+1})^*$. As the identifications $T_{ij}^{(l)}=T_{ij}^{(l)}(T,\ldots,T^{(r)})$ make sense for all $1\leq i\leq j\leq g$ and all $l\in\{0,\ldots,r\}$, for every polynomial $F\in \mathbb K[T_{ij},\ldots,T^{(r)}_{ij}|1\leq i\leq j\leq g]$ we can write
$$F=F(T,\ldots,T^{(r)})=F(T_{11},T_{12},\ldots,T^{(r)}_{gg})$$
as well as we can similarly evaluate $F$ at elements of $\V_g\otimes_{\mathbb K}\Sym((\V_g^{r+1})^*)$.

Recall the $\pmb{\SL}_g$-action on $\V_g$ given by the rule $(\Lambda,Q)\mapsto \Lambda Q \Lambda^{\t}$, with $\Lambda\in\pmb{\SL}_g(\mathbb K)$ and $Q\in\V_g$. We claim that we have an identity
 $$(\Lambda\cdot F)(T,\ldots,T^{(r)})=F(\Lambda^{-1} T(\Lambda^{-1})^{\t},\ldots,(\Lambda^{-1} T^{(r)}(\Lambda^{-1})^{\t}).$$
 To check this we can assume $r=0$ and $F=T_{ij}$ and we need to prove that $\Lambda \cdot T_{ij}$ is the $ij$ entry of the symmetric matrix $\Lambda^{-1} T(\Lambda^{-1})^{\t}$. But
 this follows from the fact that for all $Q\in \V_g$, $(\Lambda \cdot T_{ij})(Q)$ is the $ij$ entry of the matrix $\Lambda^{-1}Q (\Lambda^{-1})^{\t}$, by the very definition of the duals of $\pmb{\SL}_g$-modules. 
 
In matrix form, for $l\in\{0,\ldots,r\}$ we have
\begin{equation}\label{EQ151.6}
(\Lambda\cdot T_{ij}^{(l)})_{1\leq i\leq j\leq g}=\Lambda^{-1}T^{(l)}(\Lambda^{-1})^{\t}.
 \end{equation}
Recalling the notation in Equation (\ref{EQ151.3}), for $g\geq 2$ and $i,j\in\{1,\ldots,g\}$ distinct we have
\begin{equation}\label{EQ151.7}
\Lambda_{ij}\cdot T_{ij}^{(l)}=T_{ij}^{(l)},\;\;\;\Lambda_{ij}\cdot T_{ii}^{(l)}=T_{jj}^{(l)}\;\;\;\textup{and}\;\;\;\Lambda_{ij}\cdot T_{jj}^{(l)}=T_{ii}^{(l)}.
\end{equation}
For $g\geq 3$, $i,j\in\{1,\ldots,g\}$ distinct, and $i',j'\in\{1,\ldots,g\}\setminus\{i,j\}$ we have
\begin{equation}\label{EQ151.8}
\Lambda_{ij}\cdot T_{i'j'}^{(l)}=T_{i'j'}^{(l)}.
\end{equation}
For $g\geq 3$ and $i,j,i'\in\{1,\ldots,g\}$ distinct we have
\begin{equation}\label{EQ151.9}
\Lambda_{ij}\cdot T_{ji'}^{(l)}=-T_{ii'}^{(l)}\;\;\;\textup{and}\;\;\;\Lambda_{ij}\cdot T_{ii'}^{(l)}=T_{ji'}^{(l)}.
\end{equation}
 
We view 
$$\Sym((\V_g^{r+1})^*)=\oplus_{s\in\mathbb Z_{\ge
0}} \Sym((\V_g^{r+1})^*)_s$$ 
as a $\mathbb Z_{\geq 0}$-graded $\mathbb K$-algebra by declaring each $T_{ij}^{(l)}$ with $0\leq l\leq r$ and $1\leq i\leq j\leq g$ to be homogeneous of degree $1$. For a non-zero polynomial $h\in \Sym((\V_g^{r+1})^*)_s$ we will speak about its $T$-degree $s$. If for each $l\in\{0,\ldots,r\}$, $h$ has partial homogeneous degree $s_l$ in the upper index ${}^{(l)}$, then we also say that $h$ has partial $T$-degree $(s_0,\ldots,s_l)$; we have $s=\sum_{l=0}^r s_l$. 

We consider the $\mathbb K$-monomorphisms (inclusions) 
$$\Sym((\V_g^1)^*)\subset \Sym((\V_g^2)^*)\subset\cdots\subset\Sym((\V_g^r)^*)\subset\Sym((\V_g^{r+1})^*)\subset\cdots;$$
they correspond to natural projections $\cdots\rightarrow\V_g^{r+1}\rightarrow\V_g^r\rightarrow\cdots\rightarrow\V_g^2\rightarrow\V_g^1$ on first coordinates and induce $\mathbb Z_{\geq 0}$-graded $\mathbb K$-monomorphisms (inclusions)
\begin{equation}\label{EQ152}
\Sym((\V_g^1)^*)^{\pmb{\SL}_g}\subset \Sym((\V_g^2)^*)^{\pmb{\SL}_g}\subset\cdots
\subset\Sym((\V_g^r)^*)^{\pmb{\SL}_g}\subset\cdots.\end{equation}
As $\pmb{\SL}_g$ is semisimple,
the rings in Equation (\ref{EQ152}) are unique factorization domains with $\mathbb K^\times$ as their groups of units.

For $s\in\mathbb W_{\varepsilon_g}$ let 
$$\mathbb H^r_g(s)_{\mathbb K}:=
 \Sym((\V_g^{r+1})^*)_{gs}^{\pmb{\SL}_g}.$$
 Thus we get a 
 $\mathbb W_{\varepsilon_g}$-graded $\mathbb K$-algebra 
 $$\mathbb H_{g,\tot,\mathbb K}^r:=\bigoplus_{s\in\mathbb W_{\varepsilon_g}} \mathbb H^r_g(s)_{\mathbb K}$$
contained in $\Sym((\V_g^{r+1})^*)$ and a $\mathbb Z_{\geq 0}$-graded $\mathbb K$-algebra 
 $$\mathbb H_{g,\mathbb K}^r:=\bigoplus_{s\in\mathbb Z_{\ge
0}} \mathbb H^r_g(s)_{\mathbb K}$$
contained in either $\mathbb H_{g,\tot,\mathbb K}^r$ or $\Sym((\V_g^{r+1})^*)$. If $g$ is even, we consider the $\mathbb H_{g,K}^{r}$-module
$$\mathbb H^{r,\perp}_{g,\mathbb K}:=\bigoplus_{s\in \mathbb W_0\setminus\mathbb Z_{\geq 0}} \mathbb H^r_g(s)_{\mathbb K}=\bigoplus_{s\in \frac{1}{2}+\mathbb Z_{\geq 0}} \mathbb H^r_g(s)_{\mathbb K}.$$

\begin{lemma}\label{L43}
The following three properties hold:

\medskip
{\bf (a)} If $g$ is odd, then we have $\mathbb H^r_{g,\mathbb K}=\mathbb H_{g,\tot,\mathbb K}^r=\Sym((\V_g^{r+1})^*)^{\pmb{\SL}_g}$. 

\smallskip
{\bf (b)} We assume that $g$ is even. Then we have $\mathbb H_{g,\tot,\mathbb K}^r=\Sym((\V_g^{r+1})^*)^{\pmb{\SL}_g}$ and the square of each homogeneous polynomial of $\mathbb H_{g,\tot,\mathbb K}^r$ belongs to $\mathbb H^r_{g,\mathbb K}$ and hence the $\mathbb K$-monomorphism (inclusion) $\mathbb H^r_{g,\mathbb K}\rightarrow \mathbb H_{g,\tot,\mathbb K}^r$ is finite. Moreover, $\mathbb H^r_{g,\mathbb K}$ is normal and exactly one of the following two possibilities hold:

\smallskip\noindent
{\bf (b.i)} either $\mathbb H^{r,\perp}_{g,\mathbb K}=0$, i.e., we have an identity $\mathbb H^r_{g,\mathbb K}=\mathbb H_{g,\tot,\mathbb K}^r$;

\smallskip\noindent
{\bf (b.ii)} or $\mathbb H^{r,\perp}_{g,\mathbb K}\neq 0$ and then $[\Frac(\mathbb H_{g,\tot,\mathbb K}^r):\Frac(\mathbb H^r_{g,\mathbb K})]=2$ and there exists a homogeneous prime element 
$$\mathfrak h_{g,r,\odd}\in\mathbb H^{r,\perp}_{g,\mathbb K}$$ 
of the smallest $T$-degree such that we have three identities
$$\begin{array}{rcl}
\Frac(\mathbb H_{g,\tot,\mathbb K}^r) & = & \Frac(\mathbb H^r_{g,\mathbb K})\oplus \Frac(\mathbb H^r_{g,\mathbb K})\mathfrak h_{g,r,\odd}\\
\ & \ & \\
\ & = & \Frac(\mathbb H^r_{g,\mathbb K})[x]/[x^2-(\mathfrak h_{g,r,\odd})^2],\\
\ & \ & \\
\Frac(\mathbb H_{g,\tot,\mathbb K}^r)\otimes_{\mathbb H_{g,\mathbb K}^r} \mathbb H_{g,\mathbb K}^{r,\perp} & = & \Frac(\mathbb H_{g,\tot,\mathbb K}^r)\mathfrak h_{g,r,\odd},\end{array}$$
hence the finite field extension $\Frac(\mathbb H^r_{g,\mathbb K})\to \Frac(\mathbb H^r_{g,\tot,\mathbb K})$ is separable if and only if $p\neq 2$.

\smallskip
{\bf (c)} If $\mathbb H_{g,\tot,\mathbb K}^r$ is Cohen--Macaulay, then $\mathbb H^r_{g,\mathbb K}$ is also Cohen--Macaulay and $\mathbb H_{g,\tot,\mathbb K}^r$ is in fact Gorenstein.

\smallskip
{\bf (d)} If $p=0$, then $\mathbb H_{g,\tot,\mathbb K}^r$ is Cohen--Macaulay.\end{lemma}

\noindent
{\it Proof.}
If there exists a non-zero polynomial $h\in\Sym((\V_g^{r+1})^*)_s^{\pmb{\SL}_g}$, then for a $\mathbb K$-algebra $\mathbb B$, an element $\lambda\in\pmb{\mu}_g(\mathbb B)=Z(\pmb{\SL}_g)(\mathbb B)$ acts on $h$ as the multiplication by $\lambda^{2s}=1$ and thus $g$ divides $2s$. From this the part (a) and the first two parts of (b) follow. As $\pmb{\SL}_g$ is semisimple, from \cite{HR}, Cor. 1.9 (see \cite{Mur}, Thm. for the Gorenstein conclusion) we get that part (d) holds and part (c) holds for $\mathbb H_{g,\tot,\mathbb K}^r$. Thus parts (b) and (c) hold if $\mathbb H^r_{g,\mathbb K}=\mathbb H_{g,\tot,\mathbb K}^r$.

To end the proof it suffices to show that if $\mathbb H^r_{g,\mathbb K}\neq\mathbb H_{g,\tot,\mathbb K}^r$ (so $g$ is even), then $[\Frac(\mathbb H_{g,\tot,\mathbb K}^r):\Frac(\mathbb H^r_{g,\mathbb K})]=2$, there exists $\mathfrak h_{g,r,\odd}\in\mathbb H_{g,\mathbb K}^{r,\perp}$ such that (b.ii) holds, $\mathbb H^r_{g,\mathbb K}$ is normal, and, if $\mathbb H_{g,\tot,\mathbb K}^r$ is Cohen--Macaulay, so is $\mathbb H^r_{g,\mathbb K}$. 

Let $\mathfrak h_{g,r,\odd}\in\mathbb H_{g,\mathbb K}^{r,\perp}$ be a non-zero homogeneous element of the smallest (positive) $T$-degree. As all prime divisors of $\mathfrak h_{g,r,\odd}$ in the unique factorization domain $\mathbb H_{g,\tot,\mathbb K}^r$ are homogeneous of positive $T$-degree, from the smallest $T$-degree property of $\mathfrak h_{g,r,\odd}$, we get that $\mathfrak h_{g,r,\odd}$ is a prime element of $\mathbb H_{g,\tot,\mathbb K}^r$.

As $\mathbb H_{g,\tot,\mathbb K}^r$ is an integral domain, $1$ and $\mathfrak h_{g,r,\odd}$ generate an $\mathbb H^r_{g,\mathbb K}$-submodule of $\mathbb H_{g,\tot,\mathbb K}^r$ which is free of rank $2$. If the $\mathbb K$-algebra $\mathbb H^r_{g,\tot,\mathbb K}$ is generated by homogeneous elements $h_1,\ldots,h_m\in\mathbb H^r_{g,\mathbb K}$ and $\mathfrak h_{g,r,\odd},h_1^{\perp},\ldots,h_{m^\perp}^{\perp}\in\mathbb H_{g,\mathbb K}^{r,\perp}$, then the $\mathbb K$-algebra $\mathbb H^r_{g,\mathbb K}$ is generated by $h_1,\ldots,h_m$ and by products of two elements in the set $\{\mathfrak h_{g,r,\odd},h_1^{\perp},\ldots,h_{m^\perp}^{\perp}\}$, and moreover $\mathfrak h_{g,r,\odd},h_1^{\perp},\ldots,h_{m^\perp}^{\perp}\in\mathbb H_{g,\mathbb K}^{r,\perp}$ are homogeneous generators of the $\mathbb H^r_{g,\mathbb K}$-module $\mathbb H_{g,\mathbb K}^{r,\perp}$. Thus, as 
$$(\mathfrak h_{g,r,\odd})^2,\mathfrak h_{g,r,\odd}h_1^{\perp},\ldots,\mathfrak h_{g,r,\odd}h_{m^\perp}^{\perp}\in\mathbb H^r_{g,\mathbb K},$$ 
the three identities of part (b) follow from the identities
\begin{equation}\label{EQ153}
\begin{array}{rcl}
\Frac(\mathbb H_{g,\tot,\mathbb K}^r) & = & \Frac(\mathbb H^r_{g,\mathbb K})(\mathfrak h_{g,r,\odd},h_1^{\perp},\ldots,h_m^{\perp})\\
\ & \ & \ \\
\ & 
= & \Frac(\mathbb H^r_{g,\mathbb K})(\mathfrak h_{g,r,\odd})\\
\ & \ & \ \\
\ & = & \Frac(\mathbb H^r_{g,\mathbb K})[x]/(x^2-(\mathfrak h_{g,r,\odd})^2)\\
\ & \ & \ \\
\ & 
= & \Frac(\mathbb H^r_{g,\mathbb K})\oplus \Frac(\mathbb H^r_{g,\mathbb K})\mathfrak h_{g,r,\odd}.
\end{array}
\end{equation}

As $\mathbb H^r_{g,\tot,\mathbb K}=\mathbb H^r_{g,\mathbb K}\oplus \mathbb H^{r,\perp}_{g,\mathbb K}$ is normal (being a unique factorization domain) and as each non-zero element of $\mathbb H^{r,\perp}_{g,\mathbb K}$ belongs to $\Frac(\mathbb H^r_{g,\mathbb K})\mathfrak h_{g,r,\odd}$, $\mathbb H_{g,\mathbb K}^r$ equals to its normalization in $\mathbb H^r_{g,\tot,\mathbb K}$ and hence is normal.

We will use standard equivalent definitions of the notion Cohen--Macaulay ring for finitely generated $\mathbb Z_{\geq 0}$-graded algebras over fields (see \cite{Cho}, Thm. 1). As $\mathbb H_{g,\tot,\mathbb K}^r$ is Cohen--Macaulay of dimension $(r+1)\frac{g(g+1)}{2}-g^2+1$ (see (a) and Theorem \ref{T25}), there exists a homogeneous system of parameters $h_1,\ldots,h_{(r+1)\frac{g(g+1)}{2}-g^2+1}$ which is a regular sequence of $\mathbb H_{g,\tot,\mathbb K}^r$. Part (b) implies that $h_1^2,\ldots,h^2_{(r+1)\frac{g(g+1)}{2}-g^2+1}$ is a homogeneous system of parameters which is a regular sequence of $\mathbb H^r_{g,\mathbb K}$, hence $\mathbb H^r_{g,\mathbb K}$ is Cohen--Macaulay.\endproof

\subsubsection{Rings and vector spaces of invariants: the general case}\label{S4203}

To list more applications needed in Sections \ref{S2} and \ref{S3}, let 
$$\mathbb H^{r,\mathbb Z}_{g,\tot}:=\mathbb H^r_{g,\tot,\mathbb C}\cap (\mathbb Z[T,\ldots,T^{(r)}])\;\;\;\textup{and}\;\;\;\mathbb H^{r,\mathbb Z}_{g}:=\mathbb H^r_{g,\mathbb C}\cap (\mathbb Z[T,\ldots,T^{(r)}]).$$
For $s\in\mathbb W_{\varepsilon_g}$ let $\mathbb H^{r,\mathbb Z}_{g,\tot}(s):=\mathbb H^r_{g,\tot}(s)_{\mathbb C}\cap \mathbb H^{r,\mathbb Z}_{g,\tot}$. Also, for $s\in\mathbb Z_{\geq 0}$ we define $\mathbb H^{r,\mathbb Z}_{g}(s):=\mathbb H^r_{g}(s)_{\mathbb C}\cap \mathbb H^{r,\mathbb Z}_{g}=\mathbb H^{r,\mathbb Z}_{g,\tot}(s)$. We have direct sum decompositions
$$\mathbb H^{r,\mathbb Z}_{g}=\oplus_{s\in\mathbb Z_{\geq 0}}\mathbb H^{r,\mathbb Z}_{g}(s)\subset\mathbb H^{r,\mathbb Z}_{g,\tot}=\oplus_{s\in\mathbb W_{\varepsilon_g}}\mathbb H^{r,\mathbb Z}_{g,\tot}(s).$$
For an arbitrary $\mathbb Z$-algebra $B$ we set 
$$\mathbb H^{r,\mathbb Z}_{g,\tot,B}:=\mathbb H^{r,\mathbb Z}_{g,\tot}\otimes_{\mathbb Z} B,\;\;\;\mathbb H^{r,\mathbb Z}_{g,\tot}(s)_{B}:=\mathbb H^{r,\mathbb Z}_{g,\tot}(s)\otimes_{\mathbb Z} B,\;\;\;\mathbb H^{r,\mathbb Z}_{g,B}:=\mathbb H^{r,\mathbb Z}_{g}\otimes_{\mathbb Z} B\;\;\;\textup{and}$$
$$\mathbb H^{r,\mathbb Z}_{g}(s)_{B}:=\mathbb H^{r,\mathbb Z}_{g}(s)\otimes_{\mathbb Z} B.$$
We get direct sum decompositions
$$\mathbb H^{r,\mathbb Z}_{g,B}=\oplus_{s\in\mathbb Z_{\geq 0}}\mathbb H^{r,\mathbb Z}_{g}(s)_{B}\subset \mathbb H^{r,\mathbb Z}_{g,\tot,B}=\oplus_{s\in\mathbb W_{\varepsilon_g}}\mathbb H^{r,\mathbb Z}_{g,\tot}(s)_{B}.$$
From \cite{Se}, Thm. 2 we get that the $\mathbb Z$-algebra $\mathbb H^{r,\mathbb Z}_{g,\tot}$ is finitely generated and arguments similar to the ones of the proof of Lemma \ref{L43} (a) and (b) show that the same holds for $\mathbb H^{r,\mathbb Z}_g$. Thus the $B$-algebras $\mathbb H^{r,\mathbb Z}_{g,\tot,B}$ and $\mathbb H^{r,\mathbb Z}_{g,B}$ are finitely presented. Moreover, if $B$ is a field of characteristic $0$, then the $B$-algebras $\mathbb H^{r,\mathbb Z}_{g,\tot,B}$ and $\mathbb H^{r,\mathbb Z}_{g,B}$ are Cohen--Macaulay unique factorization domains and hence Gorenstein (see \cite{HR}, Cor. 1.9) and their groups of units are $B^\times$.

If $p=0$, then we have canonical identifications
\begin{equation}\label{EQ153.2}
\mathbb H^r_{g,\tot,\mathbb K}=\mathbb H^{r,\mathbb Z}_{g,\tot,\mathbb K}\;\;\;\textup{and}\;\;\;\mathbb H^r_{g,\mathbb K}=\mathbb H^{r,\mathbb Z}_{g,\mathbb K}
\end{equation} 
compatible with the gradings but for $p>0$ we only have 
canonical inclusions ($\mathbb K$-algebra monomorphisms compatible with the gradings)
\begin{equation}\label{EQ153.3}
\mathbb H^{r,\mathbb Z}_{g,\tot,\mathbb K}\subset \mathbb H^r_{g,\tot,\mathbb K}\;\;\;\textup{and}\;\;\;\mathbb H^{r,\mathbb Z}_{g,\mathbb K}\subset\mathbb H^r_{g,\mathbb K}.
\end{equation}

Based on the last two paragraphs, we can set
$$D(g,r):=(r+1)\frac{g(g+1)}{2}-g^2+1,$$
$$\begin{array}{rclll}
D_{p}(g,r,s) & := & \dim_{\mathbb K}(\mathbb H^r_g(s)_{\mathbb K}) & \textup{for} & s\in\mathbb W_{\varepsilon_g},\\
\ & \ & \ & \ & \ \\
D(g,r,s) & := & D_0(g,r,s) & \textup{for} & s\in\mathbb W_{\varepsilon_g}.\end{array}$$
For instance, $D(1,r)=r+1$, $\mathbb H_{1,\mathbb K}^r$ is a polynomial $\mathbb K$-algebra in $r+1$ indeterminates of degree $1$ and therefore $D_{p}(1,r,s)=\binom{r+s}{s}$ for all $s\in\mathbb Z_{\geq 0}$. Similarly, for $g\geq 2$ we have $\mathbb H_{g,\mathbb K}^0=\mathbb K$. From inclusions (\ref{EQ153.3}) we get that 
$$D_{p}(g,r,s)\geq D(g,r,s).$$ 
Once an identity $D_{p}(g,r,s)=D(g,r,s)$ is proved, afterwards we will use the simpler notation $D(g,r,s)$ instead of $D_{p}(g,r,s)$.

From Theorem \ref{T25} and Lemma \ref{L43} (a) and (b) we get directly:

\begin{cor}\label{C20}
For all integers $r\geq 1$ we have
$$\dim(\mathbb H^r_{g,\tot,\mathbb K})=\dim(\mathbb H^r_{g,\mathbb K})=D(g,r).$$
\end{cor}

\subsubsection{The case $s=\frac{1}{2}$}\label{S4204}

To compute the numbers $D_{p}(g,r,s)$ in the case $s=\frac{1}{2}$ we recall the computation of all $\pmb{\T}_g$-invariants of $\Sym((\V_g^{r+1})^*)=\mathbb K[T,\ldots,T^{(r)}]$. 

\begin{rem}\label{R41}
The torus $\pmb{\T}_g$ leaves invariant the $\mathbb K$-span of each monomial 
$$\mu:=\prod_{1\leq i\leq j\leq g}\prod_{l=0}^r (T_{ij}^{(l)})^{n_{ij,l}}$$ 
in the $T_{ij}^{(l)}$s, where each $n_{ij,l}\in\mathbb N\cup\{0\}$. For $m\in\{1,\ldots, g\}$, let
$$\Index_m(\mu):=\sum_{1\leq i\leq j\leq g}\sum_{l=0}^r n_{ij,l}(\delta_{mi}+\delta_{mj}).$$ 
We have 
\begin{equation}\label{EQ153.35}
\sum_{m=1}^g \Index_m(\mu)=2\deg(\mu).
\end{equation}

The monomial $\mu$ is $\pmb{\T}_g$-invariant if and only if $\Index_m(\mu)$ is independent on $m\in\{1,\ldots, g\}$, in which case it will be denoted simply as $\Index(\mu)$ and from Equation (\ref{EQ153.35}) we get that $$2\deg(\mu)=g\Index(\mu);$$ 
thus, if $g$ is even, we have $\deg(\mu)=\frac{g}{2}$ if and only if $\Index(\mu)=1$. 
\end{rem}

\begin{prop}\label{P14}
If $g$ is even, then we have $D_{p}(g,r,\frac{1}{2})=0$. 
\end{prop}

\noindent
{\it Proof.} If $g=2$, we have $\mathbb H_2^r(\frac{1}{2})_{\mathbb K}\subset \oplus_{l=0}^r \mathbb KT_{12}^{(l)}$ (see Remark \ref{R41}) and $T_{12}^{(l)}$ is not an $\pmb{\SL}_2$-invariant (for $p=2$, see Lemma \ref{L41} (a) to (d) in the dual context or see Lemma \ref{L41} (j)). Thus $\mathbb H_2^r(\frac{1}{2})_{\mathbb K}=0$ and $D_{p}(2,r,\frac{1}{2})=0$.


For $g\geq 4$ we show that the assumption that there exists $h\in\mathbb H_g^r(\frac{1}{2})_{\mathbb K}\setminus\{0\}$ leads to a contradiction. We know $h$ is in the span of monomials $\prod_{m=1}^{\frac{g}{2}} T_{i_m,j_m}^{(l_m)}$ with $i_m< j_m$, $\{i_1,j_1,\ldots,i_{\frac{g}{2}},j_{\frac{g}{2}}\}=\{1,\ldots,g\}$ and $l_1,\ldots,l_{\frac{g}{2}}\in\{0,\ldots,r\}$, see Remark \ref{R41}. We fix such a monomial $\mu:=\prod_{m=1}^{\frac{g}{2}} T_{i_m,j_m}^{(l_m)}$ which appears with a non-zero coefficient $c$ in $h$. 

For $p\neq 2$, the coefficient of $\mu$ in $\Lambda_{i_1,j_2}\cdot h$ is $-c$ (see Equations (\ref{EQ151.7}) to (\ref{EQ151.9})) which contradicts the fact that $h$ is an $\pmb{\SL}_g$-invariant. 

If $p=2$, let $\Lambda\in\pmb{\SL}_g(\mathbb K)$ be such that it normalizes $\Span(\{e_{i_1},e_{j_1}\})$, fixes each $e_{i'}$ with $i'\in\{1,\ldots,g\}\setminus\{i_1,j_1\}$ and (by the case $g=2$) does not fix $T^{(l_1)}_{i_1,j_1}$. Based on Equation (\ref{EQ151.6}), similar to Equations (\ref{EQ151.7}) to (\ref{EQ151.9}) we argue that for each $l\in\{0,\ldots,r\}$, $\Lambda$ leaves invariant $\Span(\{T_{i_1,i_1}^{(l)},T_{i_1,j_1}^{(l)},T_{j_1,j_1}^{(l)}\})$ and $\Span(\{T_{ij}^{(l)}|i\notin\{i_1,j_1\}\;\textup{or}\; j\notin\{i,j\}\})$ and fixes each $T_{ij}^{(l)}$ with $i,j\in\{1,\ldots,g\}\setminus\{i_1,j_1\}$. From the last two sentences we get that $\Lambda\cdot h\neq h$ which contradicts the fact that $h$ is an $\pmb{\SL}_g$-invariant.\endproof

\subsubsection{The case $r=1$}\label{S4205}

To compute the numbers $D_{p}(g,r,s)$ in the case $r=1$ and $p\neq 2$ we first recall the following well-known result.

\begin{lemma}\label{L44}
Let $y_0,\ldots,y_g$ be indeterminates. Let $\sigma_{g,0}(y_1,\ldots,y_g):=1$ and
$$\sigma_{g,1}(y_1,\ldots,y_g):=\sum_{m=1}^g y_m,\,\cdots,\,\sigma_{g,g}(y_1,\ldots,y_g):=\prod_{m=1}^g y_m$$ 
be the elementary symmetric polynomials in the indeterminates $y_1,\ldots,y_g$. Then the $\mathbb K$-algebra endomorphism 
$$\theta_g:\mathbb K[y_0,\ldots,y_r]\to \mathbb K[y_0,\ldots,y_r]$$ 
given for $m\in\{1,\ldots,g\}$ by the rule 
$$\theta_g(y_m):=y_0^{g-m}\sigma_{g,m}(y_1,\ldots,y_g)$$
becomes integral after inverting $y_0$ (i.e., the induced endomorphism $(\theta_g)_{y_0}$ 
of $\mathbb K[y_0,\ldots,y_r]_{y_0}$ is integral).
\end{lemma}

\noindent
{\it Proof.} Let $z_1:=y_0^{-1}y_1,\ldots,z_g:=y_0^{-1}y_g\in \mathbb K[y_0,\ldots,y_r][y_0^{-1}]$. The $\mathbb K$-algebra monomorphism $\mathbb K[\sigma_1(z_1,\ldots,z_g),\ldots,\sigma_g(z_1,\ldots,z_g)]\rightarrow \mathbb K[z_1,\ldots,z_g]$
is finite as $z_1,\ldots,z_g$ are roots of the monic polynomial
$$\sum_{m=0}^g (-1)^mz^{g-m}\sigma_m(z_1,\ldots,z_m)$$and by tensoring it with $\mathbb K[y_0,y_0^{-1}]$ over $\mathbb K$ we get that $(\theta_g)_{y_0}$ is integral.\endproof

\medskip

The $p=0$ case of the next theorem is classical (for instance, it is proved in \cite{T}, Ch. 20, Sect. 6, p. 304).

\begin{thm}\label{T29}
The polynomials $\Theta_0,\ldots,\Theta_g$ 
are algebraically independent over $\mathbb K$ and hence the $\mathbb K$-subalgebra $\mathbb H_{g,\mathbb K}^{1,\Theta}:=\mathbb K[\Theta_0,\ldots,\Theta_g]$ of $\mathbb H_{g,\tot,\mathbb K}^1$ generated by them is a polynomial $\mathbb K$-algebra in $g+1$ indeterminates. Moreover, we have:

\medskip
{\bf (a)} If $p\neq 2$, then we have
$$\mathbb H_{g,\mathbb K}^{1,\Theta}=\mathbb H_{g,\mathbb K}^1=\mathbb H_{g,\tot,\mathbb K}^1=\Sym((\V_g^2)^*)^{\pmb{\SL}_g}=\mathcal O(\pmb{\V}_g^2)^{\pmb{\SL}_g}$$
(thus, if $g$ is even, we have $\mathbb H_{g,\mathbb K}^{1,\perp}=0$ and $D_{p}(g,1,s)=0$ if $s\in\mathbb W_0\setminus\mathbb Z_{\geq 0}$) and $D_{p}(g,1,s)=\binom{g+s}{s}$ for all $s\in\mathbb Z_{\geq 0}$.

\smallskip
{\bf (b)} If $p=2$, then the $\mathbb K$-monomorphisms $\mathbb H_{g,\mathbb K}^{1,\Theta}\to \mathbb H_{g,\mathbb K}^1\to\mathbb H_{g,\tot,\mathbb K}^1$ are integral and purely inseparable (but $\mathbb H_{g,\mathbb K}^{1,\Theta}\to \mathbb H_{g,\mathbb K}^1$ is not an isomorphism for $g>1$).\end{thm}

\noindent
{\it Proof.} With $F_0\in\mathbb H_{g,\mathbb K}^{1,\Theta}$ as in the proof of Theorem \ref{T25}, we claim that by inverting $F_0$ the $\mathbb K$-algebra monomorphism $\mathbb H_{g,\mathbb K,F_0}^{1,\Theta}\to\mathbb H_{g,\tot,\mathbb K,F_0}^1$ is integral and purely inseparable for $p=2$ and is an isomorphism if $p\neq 2$. 

The fibers of the morphism $\varpi_{F_0}:\Spec(\mathbb K[T,T^{(1)}]_{F_0})\to\Spec(\mathbb H_{g,\mathbb K,F_0}^{1,\Theta})$ were in essence computed in the proof of Theorem \ref{T25}. To explain this we consider a point $P\in\Spec(\mathbb H_{g,\mathbb K,F_0}^{1,\Theta})(\mathbb K)$ given by a surjective $\mathbb K$-homomorphism $\mathbb H_{g,\mathbb K,F_0}^{1,\Theta}\to\mathbb K$ which for $i\in\{0,\ldots,g\}$ maps $\Theta_i$ to some $\nu_i\in\mathbb K$. Let $\lambda_0\in\mathbb K$ be such that $\lambda_0^g=\nu_0$ and let $\lambda_1,\ldots,\lambda_g\in\mathbb K$ be such that (see Equation (\ref{EQ149.5}))
$$\nu_0\prod_{i=1}^n (x-\lambda_i)=\sum_{i=0}^g x^{g-i}\nu_i.$$ 
As $F_0(\nu_0,\ldots,\nu_g)\neq 0$, $\lambda_0\neq 0$, $\lambda_1,\ldots,\lambda_g$ are distinct, and, if $p=2$, we have $\sum_{i=1}^g \lambda_i\neq 0$. Let $D_0:=\lambda_0 1_g$, $D_1:=\Diag(\lambda_0 \lambda_1,\ldots,\lambda_0 \lambda_g)\in\Z_g$. The pair $(D_0,D_1)$ defines a $\mathbb K$-valued point of $\varpi_{F_0}^{-1}(P)$ and let ${\bf o}$ be its orbit in $\pmb{\V}_g^2$; it is a smooth, closed subscheme of $\pmb{\V}_g^2$ of dimension $g^2-1$. We have an $\pmb{\SL}_g$-invariant closed embedding ${\bf o}\subset \varpi_{F_0}^{-1}(P)$.

In this paragraph we check that ${\bf o}=(\varpi_{F_0}^{-1}(P))_{\red}$, i.e., that we have ${\bf o}(\mathbb K)=(\varpi_{F_0}^{-1}(P))(\mathbb K)$. It suffices to show that if $\lambda'_0,\lambda'_1,\ldots,\lambda'_g\in\mathbb K$ are such that $(\lambda'_0)^g=\nu_0=\lambda_0^g$ and $\{\lambda_0 \lambda_1,\ldots,\lambda_0 \lambda_g\}=\{\lambda'_0 \lambda'_1,\ldots,\lambda'_0 \lambda'_g\}$, then $(D_0,D_1)$ is $\pmb{\SL}_g$-equivalent to $(\lambda'_01_g,\Diag(\lambda'_0 \lambda'_1,\ldots,\lambda'_0 \lambda'_g))$. Based on the proof of Theorem \ref{T25}, we can assume that $\lambda'_0=\lambda_0$ and it suffices to show that $D_1$ and $\Diag(\lambda_0 \lambda'_1,\ldots,\lambda_0 \lambda'_g)$ are $\pmb{\SO}_g$-equivalent. But this follows by repeatedly applying the fact that for $g\geq 2$ and integers $1\leq i<j\leq g$, similarly to Equation (\ref{EQ151.7}) we get that 
$\Lambda_{ij}\cdot D_1$ is obtained from $D_1$ by interchanging its $ii$ and $jj$ entries. 

To prove that $\varpi_{F_0}^{-1}(P)={\bf o}$ is reduced if $p\neq 2$, it suffices to show that if $\mathbb B:=\mathbb K[x]/(x^2)=\mathbb K\oplus \mathbb K\epsilon$, with $\epsilon:=x+(x^2)$, then the injection 
$${\bf o}(\mathbb B)\rightarrow \varpi_{F_0}^{-1}(P)(\mathbb B)$$ 
is in fact a bijection. To check this, it suffices to show that if $(D'_0,D'_1)\in \mathbb B\otimes_{\mathbb K}\V_g$ is such that modulo $\epsilon$ is $(D_0,D_1)$, $\det(D_0')=\det(D_0)$ and $(D'_0)^{-1}D'_1$ has eigenvalues $\lambda_1,\ldots,\lambda_g$, then there exists $\Lambda'\in\Ker(\pmb{\SL}_g(\mathbb B)\rightarrow \pmb{\SL}_g(\mathbb K))$ such that we have $(\Lambda'D_0(\Lambda')^{\t},\Lambda'D_1(\Lambda')^{\t})=(D'_0,D'_1)$. To check this, as $p\neq 2$ each $Q\in\V_g$ of trace zero is of the form $\lambda_0(M+M^{\t})$ with $M\in\Mat_g$ of trace zero (for instance, we can take $M=\frac{1}{2\lambda_0}Q$) and we can assume that $D'_0=D_0$\footnote{If $g\geq 2$, for $p=2$ we cannot assume that $D'_0=D_0$.}, and we show that there exists $\Lambda'\in\Ker(\pmb{\SO}_g(\mathbb B)\rightarrow \pmb{\SO}_g(\mathbb K))$ such that we have $\Lambda'D_1(\Lambda')^{\t}=\Lambda'D_1(\Lambda')^{-1}=D_1'$. As $D_1$ and $D'_1$ have the same distinct eigenvalues, there exists a matrix $M\in \Mat_g$ such that $D'_1=(1_g+\epsilon M)D_1(1_g-\epsilon M)=D_1+\epsilon(MD_1-D_1M)$. As $MD_1-D_1M$ and $D_1$ are symmetric, $D_1$ commutes with $M+M^{\t}$. As $D_1\in \Z_g$ has distinct eigenvalues, we get that $M+M^{\t}\in \Z_g$. This implies that there exists $M_1\in\Mat_g$ with $M_1^{\t}+M_1=0$ (so $M_1\in\Lie(\pmb{\SO}_g)$) such that $M-M_1\in \Z_g$. We have $MD_1-D_1M=M_1D_1-D_1M_1$ and hence we can take $\Lambda'=1_g+\epsilon M_1$. 

The morphism $\Spec(\mathbb H_{g,\tot,\mathbb K,F_0}^1)\rightarrow \Spec(\mathbb H_{g,\mathbb K}^{1,\Theta,F_0})$ induces a bijection at the level of $\mathbb K$-valued points and has generically reduced fibers if $p\neq 2$, see the previous paragraph. This implies that this morphism is integral, purely inseparable if $p=2$ and birational if $p\neq 2$, and from Zariski's Main Theorem we get that it is an isomorphism if $p\neq 2$, i.e., the claim holds. 

From previous paragraph and Corollary \ref{C20} we get an equality $\dim(\mathbb H_{g,\mathbb K}^{1,\Theta})=g+1$. As $\mathbb H_{g,\mathbb K}^{1,\Theta}$ is a quotient of a polynomial $\mathbb K$-algebra in $\dim(\mathbb H_{g,\mathbb K}^{1,\Theta})=g+1$ indeterminates, it is a polynomial $\mathbb K$-algebra in $g+1$ indeterminates and $\Theta_0,\ldots,\Theta_g$ are algebraically independent over $\mathbb K$. Thus $\mathbb H_{g,\mathbb K}^{1,\Theta}\subset\mathbb H_{g,\tot,\mathbb K}^1$ is an inclusion between two noetherian unique factorization (hence normal) domains such that for $p\neq 2$ we have $\Frac(\mathbb H_{g,\mathbb K}^{1,\Theta})=\Frac(\mathbb H_{g,\tot,\mathbb K}^1)$.

Recall that $\L_g=\mathbb K1_g\times \Z_g\subset\V_g^2$ and we consider the composite morphism
$$\theta_g^*:\pmb{\L}_g\to\pmb{\V}_g^2\to (\pmb{\V}_g^2)^{\pmb{\SL}_g}=\Spec(\mathbb H^r_{g,\tot,\mathbb K})\rightarrow \Spec(\mathbb K[\Theta_0,\ldots,\Theta_g]).$$
As for $(D_0,D_1):=(\lambda_0 1_g,\Diag(\lambda_1,\ldots,\lambda_g))\in\L_g$ and $m\in\{1,\ldots,g\}$ we have $\Theta_m(D_0,D_1)=\lambda_0^{g-m}\sigma_m(\lambda_1,\ldots,\lambda_g)$ (see Equation (\ref{EQ149.25})), the notation matches with the one of Lemma \ref{L44}, i.e., $\theta_g^*$ can be identified with the morphism between spectra defined by $\theta_g$. From this and Lemma \ref{L44}, as $\mathbb H_{g,\mathbb K,
\Theta_0}^{1,\Theta}\subset\mathbb H_{g,\tot,\mathbb K,\Theta_0}^1\subset\Frac(\mathbb H_{g,\mathbb K,\Theta_0}^{1,\Theta})$, we get an integral $\mathbb K$-algebra monomorphism $\mathbb H_{g,\mathbb K,\Theta_0}^{1,\Theta}\rightarrow\mathbb H_{g,\tot,\mathbb K,\Theta_0}^1$ which for $p=2$ is purely inseparable and for $p\neq 2$ is an isomorphism as we know it is birational (see previous paragraph) and as $\mathbb H_{g,\mathbb K,\Theta_0}^{1,\Theta}$ is normal. Switching the roles of $T^{(0)}$ and $T^{(1)}$, we similarly get an integral $\mathbb K$-algebra monomorphism $\mathbb H_{g,\mathbb K,\Theta_g}^{1,\Theta}\to\mathbb H_{g,\tot,\mathbb K,\Theta_g}^1$ which for $p=2$ is purely inseparable and for $p\neq 2$ is an isomorphism. Thus $\mathbb H_{g,\tot,\mathbb K}^1\subset \mathbb H_{g,\mathbb K,\Theta_0}^{1,\Theta}\cap \mathbb H_{g,\mathbb K,\Theta_g}^{1,\Theta}=\mathbb H_{g,\mathbb K}^{1,\Theta}$ if $p\neq 2$. Hence $\mathbb H_{g,\mathbb K}^{1,\Theta}=\mathbb H_{g,\tot,\mathbb K}^1$ if $p\neq 2$ from which part (a) follows. If $p=2$, then for the normalization $\mathbb H_{g,\mathbb K}^{1,\Theta,\n}$ of $\mathbb H_{g,\mathbb K}^{1,\Theta}$ in $\mathbb H_{g,\tot,\mathbb K}^1$ we have $\mathbb H_{g,\mathbb K}^{1,\Theta,\n}=\mathbb H_{g,\mathbb K,\Theta_0}^{1,\Theta,\n}\cap\mathbb H_{g,\mathbb K,\Theta_g}^{1,\Theta,\n}$ and similarly we argue that $\mathbb H_{g,\mathbb K}^{1,\Theta,\n}=\mathbb H_{g,\tot,\mathbb K}^1$ from which part (b) follows. 
\endproof

\subsubsection{Basic $\pmb{\SL}_g$-invariants: $\Theta$s and $\Upsilon$s}\label{S4206}

In order to present a generalization of the first part of Theorem \ref{T29} we first introduce the index set 
$$\Delta(g,r):=\{(m_0,\ldots,m_r)\in\mathbb Z_{\geq 0}^{r+1}|\sum_{i=0}^r m_i=g\}$$
that will index the polynomials generalizing $\Theta_0,\ldots,\Theta_g$. If 
$$\Delta^{+1}(g,r):=\{(n_0,\ldots,n_r)\in\mathbb N^{r+1}|\sum_{i=0}^r n_i=g+r+1\},$$ 
then the rule $(m_0,\ldots,m_r)\rightarrow (m_0+1,\ldots,m_r+1)$ defines a bijection
$$\Delta(g,r)\rightarrow\Delta^{+1}(g,r).$$ 
It is well-known that $\Delta^{+1}(g,r)$ has as many elements as the number of ways of inserting $r$ commas among the $g+r$ spaces (places) between two consecutive $1$s of the sequence $1\;1\;\cdots 1\; 1$ formed by $g+r+1$ numbers $1$, and thus it is equal to $\binom{g+r}{r}$. For instance, the $4$-tuple $(1,2,1,2)\in\Delta^{+1}(2,3)$ corresponds to the following insertion of commas $1,\; 1\; 1,\; 1,\; 1\; 1$. We conclude that $\Delta(g,r)$ has $\binom{g+r}{r}$ elements. For $r\in\{1,2\}$ or $g=1$ we have $\binom{g+r}{r}=D(g,r)$ and for $r\geq 3$ and $g\geq 2$ we have $\binom{g+r}{r}>D(g,r)$. Moreover, for $r\geq 3$ we have
$$\lim_{g\rightarrow\infty} \frac{\binom{g+r}{r}}{D(g,r)}=\infty.$$

We define polynomials
$$\Theta_{m_0,\ldots,m_r}\in\mathbb H^{r,\mathbb Z}_g(1)\subset\mathbb H_{g}^{r,\mathbb Z}\subset \mathbb H_{g,\tot}^{r,\mathbb Z}$$
indexed by $(m_0,\ldots,m_r)\in\Delta(g,r)$ via
 the identity
$$\det(\sum_{i=0}^r y_iQ_i)=\sum_{(m_0,\ldots,m_r)\in\Delta(g,r)} (\prod_{i=0}^r y_i^{m_i})\Theta_{m_0,\ldots,m_r}(Q_0,Q_2,\ldots,Q_r),$$
where $y_0,\ldots,y_r$ are $r+1$ indeterminates; so the $\Theta_{m_0,\ldots,m_r}$s are partial polarizations of the invariant polynomial $\det$.

We also denote by $\Theta_{m_0,\ldots,m_r}$ its image in $\mathcal O(\pmb{\V}_g^{r+1})^{\pmb{\SL}_g}=\mathcal O(\pmb{\V}_g^{r+1}/\pmb{\SL}_g)$. Note that for all $i\in\{0,\ldots,g\}$ we have $\Theta_i=\Theta_{g-i,i}$ and for $r\geq 2$ inside $\mathcal O(\pmb{\V}_g^{r+1})^{\pmb{\SL}_g}$ we identify it as well with $\Theta_{g-i,i,0,0,\ldots,0}$ (the number of $0$s being $r-1$), via the inclusions (\ref{EQ152}).

Let $\mathbb H_g^{r,\Theta}$ be the $\mathbb Z_{\geq 0}$-graded $\mathbb Z$-subalgebra 
of $\mathbb H_g^{r,\mathbb Z}$ (or $\mathbb H_{g,\tot}^{r,\mathbb Z}$) generated by all $\Theta_{m_0,\ldots,m_r}$s with $(m_0,\ldots,m_r)\in\Delta(g,r)$ (the case $r=1$ was first introduced and used in Theorem \ref{T29}). For a $\mathbb Z$-algebra $B$, we define its $\Theta$-rings of invariants $\mathbb H_{g,B}^{r,\Theta}:=\mathbb H_g^{r,\Theta}\otimes_{\mathbb Z} B$. We have $B$-algebra monomorphisms compatible with the gradings
\begin{equation}\label{EQ153.4}
\mathbb H_{g,B}^{r,\Theta}\subset \mathbb H_{g,B}^{r,\mathbb Z}\subset \mathbb H_{g,\tot,B}^{r,\mathbb Z}.
\end{equation}

For an integer $l\geq 0$, let 
$$c(T^{(l)}):=[T^{(l)}_{11}\; T^{(l)}_{12}\; \cdots\; T^{(l)}_{1g}\;T^{(l)}_{22}\; T^{(l)}_{23}\;\cdots\; T^{(l)}_{gg}]^{\t}$$
be the column vector whose entries are $T^{(l)}_{ij}$ with $1\leq i\leq j\leq g$ listed increasingly in the lexicographic order of the indexes $ij$. 
If the inequality $r\geq \frac{g(g+1)}{2}-1$ holds, then for each increasing sequence $q_1<q_2<\cdots<q_{\frac{g(g+1)}{2}}$ of integers in the set $\{0,1,\ldots,r\}$, let 
$$\Upsilon_{q_1,\ldots,q_{\frac{g(g+1)}{2}}}:=\det([c(T^{(q_1)})|\cdots |c(T^{(q_{\frac{g(g+1)}{2}})})]).$$
If $g$ is odd, then we have
$$\Upsilon_{q_1,\ldots,q_{\frac{g(g+1)}{2}}}\in\mathbb H^{q_{\frac{g(g+1)}{2}}}_g(\frac{g(g+1)}{2})_{\mathbb K}\subset\mathbb H_{g,\mathbb K}^{q_{\frac{g(g+1)}{2}}}\subset \mathbb H_{g,\mathbb K}^r.$$
If $g$ is even, then we have
$$\Upsilon_{q_1,\ldots,q_{\frac{g(g+1)}{2}}}\in\Sym((V_g^{1+q_{\frac{g(g+1)}{2}}})^*)_{\frac{g(g+1)}{2}}^{\pmb{\SL}_g}\subset\mathbb H_{g,\mathbb K}^{q_{\frac{g(g+1)}{2}},\perp}\subset\mathbb H_{g,\mathbb K}^{r,\perp}.$$
For instance, if $g=2$ (so $r\geq 2$), then we can take $\mathfrak h_{2,r,\odd}=\Upsilon_{0,1,2}$. 
The number of such $\Upsilon_{q_1,\ldots,q_{\frac{g(g+1)}{2}}}$ with $q_{\frac{g(g+1)}{2}}\leq r$ is $\binom{r+1}{\frac{g(g+1)}{2}}$ (which makes sense even if $r< \frac{g(g+1)}{2}-1$, being $0$ in such a case).

\begin{rem}\label{R42}
As for integers $r>l\geq 1$ the $\pmb{\SL}_g$-invariant projections $\V_g^{r+1}\rightarrow\V_g^{l+1}$ admit $\pmb{\SL}_g$-invariant sections, for the $\mathbb Z_{\geq 0}$-graded $\mathbb K$-monomorphisms of (\ref{EQ152}) we have retractions $\mathbb H_{g,\tot,\mathbb K}^r\rightarrow\mathbb H_{g,\tot,\mathbb K}^l$ which map $\mathbb H_{g,\mathbb K}^{r}$ onto $\mathbb H_{g,\mathbb K}^l$ and $\mathbb H_{g,\mathbb K}^{r,\Theta}$ onto $\mathbb H_{g,\mathbb K}^{l,\Theta}$. Thus, if $\mathbb H_{g,\mathbb K}^r=\mathbb H_{g,\tot,\mathbb K}^r$, then we also have $\mathbb H_{g,\mathbb K}^l=\mathbb H_{g,\tot,\mathbb K}^l$. Similarly, if $\mathbb H_{g,\mathbb K}^{r}=\mathbb H_{g,\mathbb K}^{r,\Theta}$, then we also have $\mathbb H_{g,\mathbb K}^l=\mathbb H_{g,\mathbb K}^{l,\Theta}$, and if $\mathbb H_{g,\mathbb K}^{r,\Theta}=\mathbb H_{g,\tot,\mathbb K}^r$, then we also have $\mathbb H_{g,\mathbb K}^{l,\Theta}=\mathbb H_{g,\tot,\mathbb K}^l$.
\end{rem}




\subsubsection{The case $s=1$}\label{S4207}

Let $\Gamma_{p}(g,r)$ (resp. $\Gamma_{p}(g,\tot,r)$) be the smallest number of homogeneous generators of the $\mathbb Z_{\geq 0}$-graded $\mathbb K$-algebra $\mathbb H^r_{g,\mathbb K}$ (resp. $\mathbb W_{\varepsilon_g}$-graded $\mathbb K$-algebra $\mathbb H^r_{g,\tot,\mathbb K}$). 

\begin{prop}\label{P15}
For all $g,r\in\mathbb N$ the following two properties hold:

\medskip
{\bf (a)} We have an inclusion
$$\mathbb H_g^r(1)_{\mathbb K}\supset \Span(\{\Theta_{m_0,\ldots,m_r}|(m_0,\ldots,m_r)\in\Delta(g,r)\})$$
as well as inequalities
\begin{equation}\label{EQ153.5}
\min\{\Gamma_{p}(g,r),\Gamma_{p}(g,\tot,r)\}\geq D_{p}(g,r,1)\geq\binom{g+r}{r}\geq D(g,r).
\end{equation} 
Thus, for $r\geq 3$ we have
$$\lim_{g\rightarrow\infty} \frac{\min\{\Gamma_{p}(g,r),\Gamma_{p}(g,\tot,r)\}}{D(g,r)}=\infty.$$

\smallskip
{\bf (b)} If $p\nmid g!$, then we have $$\mathbb H_g^r(1)_{\mathbb K}= \Span(\{\Theta_{m_0,\ldots,m_r}|(m_0,\ldots,m_r)\in\Delta(g,r)\})$$
as well as 
\begin{equation}\label{EQ153.75}
D_{p}(g,r,1)=\binom{g+r}{r}
\end{equation}
\end{prop}

\noindent
{\it Proof.} As $\Theta_{m_0,\ldots,m_r}$ has partial $T$-degree $(m_0,\ldots,m_r)$, the set of polynomials $\{\Theta_{m_0,\ldots,m_r}|(m_0,\ldots,m_r)\in\Delta(g,r)\}$ has $\binom{g+r}{r}$ elements and is linearly independent over $\mathbb K$. Thus we have $D_{p}(g,r,1)\geq\binom{g+r}{r}$. As the first inequality of (\ref{EQ153.5}) is obvious and the last inequality of (\ref{EQ153.5}) was checked above, to end the proof it suffices to show that $D_{p}(g,r,1)\le\binom{g+r}{r}$ if $p\nmid g!$. 

To check this, it suffices to show that for all $(m_0,\ldots,m_r)\in\Delta(g,r)$, the $\mathbb K$-vector space of polynomials in $\mathbb H_{g,\tot,\mathbb K}^r$ of partial $T$-degree $(m_0,\ldots,m_r)$ is one dimensional (generated by $\Theta_{m_0,m_1,\ldots,m_r}$). Equivalently, we have to show that for all $(m_0,\ldots,m_r)\in\Delta(g,r)$, $(\otimes_{l=0}^r \Sym^{m_l}(V_g^*))^{\pmb{\SL}_g}$ has dimension at most $1$. As we have natural $\pmb{\SL}_g$-epimorphisms 
$$(V_g^*)^{\otimes m_l}\rightarrow \Sym^{m_l}(V_g^*)$$ 
which admit sections (as $m_l$, if non-zero, is invertible in $\mathbb K$) it suffices to show that $\dim_{\mathbb K}(((V_g^*)^{\otimes g})^{\pmb{\SL}_g})\leq 1$ (recall that $\sum_{l=0}^r m_l=g$). 

The case $g=1$ is trivial and hence we can assume $g\geq 2$. This implies that $p\neq 2$. Therefore we can identify canonically $\V_g=\Sym^2(\W_g)$ (see Lemma \ref{L42} (a)) with the direct summand of $\W_g\otimes\W_g$ formed by vectors fixed by the natural `symmetrization' action of the group $(\mathbb Z/2\mathbb Z)$ on $\W_g\otimes\W_g$ (which interchanges the two factors). Taking duals we can also identify $\V_g^*=\Sym^2(\W_g^*)$ with the direct summand of $\W_g^*\otimes\W_g^*$ fixed by the natural `symmetrization' action of the group $(\mathbb Z/2\mathbb Z)$ on $\W_g^*\otimes\W_g^*$. Implicitly, we identify $(\V_g^*)^{\otimes g}$ with the direct summand of $(\W_g^*\otimes_{\mathbb K}\W_g^*)^{\otimes g}$ formed by elements fixed by the natural diagonal action of the group $(\mathbb Z/2\mathbb Z)^g$ on $(\W_g^*\otimes_{\mathbb K}\W_g^*)^{\otimes g}$. 

It is well-known that we can identify (see \cite{HR}, Sect 1., Applic. (a), Case (ii) and Ex. after Cor. 1.9)
$$\Sym((\W_g^*)^{2g})=\mathbb K[\underline{x}_1,\ldots,\underline{x}_{2g}]:=\mathbb K[x_{j,i}|1\leq i\leq g,1\leq j\leq 2g]$$ with $\underline{x}_j:=[x_{j,1}\; x_{j,2}\;\cdots x_{j,g}]^{\t}$, in such a way that $\Sym((\W_g^*)^{2g})^{\pmb{\SL}_g}$ is the $\mathbb K$-subalgebra generated by the determinants $d_{\underline{j}}$ of the $g\times g$ matrices $\underline{x}_{\underline{j}}$ whose columns are $\underline{x}_{j_1},\ldots,\underline{x}_{j_g}$, where $\underline{j}=(j_1,\ldots,j_g)$ is an arbitrary $g$-tuple of integers satisfying $1\leq j_1<j_2<\cdots<j_g\leq 2g$. Let $\nabla(g)$ be the set with $\binom{2g}{g}$ elements formed by such $g$-tuples $\underline{j}=(j_1,j_2,\ldots,j_g)$. Let $\nabla_1(g)$ be the subset of $\nabla(g)$ formed by all $g$-tuples $\underline{j}=(j_1,j_2,\ldots,j_g)$ with $j_1=1$; it has $\binom{2g-1}{g}$ elements. For $\underline{j}\in\nabla_1(g)$, let $\underline{j}'=(j_1',j_2',\ldots,j_g')\in\nabla(g)$ be such that $\sqcup_{i=1}^g \{j_i,j_i'\}=\{1,\ldots,2g\}$. Let $\underline{j}_0:=(1,3,\ldots,2g-1)\in\nabla_1(g)$.

We know that 
$$((\W_g^*\otimes_{\mathbb K} \W_g^*)^{\otimes g})^{\pmb{\SL}_g}=\{\sum_{\underline{j}\in\nabla_1(g)} c_{\underline{j}}d_{\underline{j}}d_{\underline{j'}}|c_{\underline{j}}\in\mathbb K\;\forall \underline{j}\in\nabla_1(g)\}.$$
An element $\underline{\varepsilon}=(\varepsilon_1,\ldots,\varepsilon_g)\in(\mathbb Z/2\mathbb Z)^g$ acts naturally on $((\W_g^*\otimes_{\mathbb K} \W_g^*)^{\otimes g})^{\pmb{\SL}_g}$ via the $\mathbb K$-automorphism $\iota_{\underline{\varepsilon}}$ of $\mathbb K[\underline{x}_1,\ldots,\underline{x}_{2g}]$ which interchanges $\underline{x}_{2l-1}\longleftrightarrow\underline{x}_{2l}$ for all $l$ with $\varepsilon_l=1+2\mathbb Z$ and which fixes both $\underline{x}_{2l-1}$ and $\underline{x}_{2l}$ if $\varepsilon_l=2\mathbb Z$. To end the proof it suffices to show that the $\mathbb K$-vector space
$$(((\W_g^*\otimes_{\mathbb K} \W_g^*)^{\otimes g})^{\pmb{\SL}_g})^{(\mathbb Z/2\mathbb Z)^g}$$
has dimension at most $1$. 

If $\theta=\sum_{\underline{j}\in\nabla_1(g)} c_{\underline{j}}d_{\underline{j}}d_{\underline{j'}}\in (((\W_g^*\otimes_{\mathbb K} \W_g^*)^{\otimes g})^{\pmb{\SL}_g})^{(\mathbb Z/2\mathbb Z)^g}$, then clearly we have $c_{\underline{j}}=0$ if there exists $l\in\{1,2,\ldots,g\}$ such that both $2l-1$ and $2l$ are entries of either $\underline{j}$ or $\underline{j}'$; we call such a $\underline{j}\in\nabla_1(g)$ {\it unmixed}. As the set of products $d_{\underline{j}}d_{\underline{j'}}$ indexed by all mixed (i.e., not unmixed) $\underline{j}\in\nabla_1(g)$ is permuted transitively by $(\mathbb Z/2\mathbb Z)^g$, we get that $\theta$ is uniquely determined by $c_{\underline{j}_0}$, hence $(((\W_g^*\otimes_{\mathbb K} \W_g^*)^{\otimes g})^{\pmb{\SL}_g})^{(\mathbb Z/2\mathbb Z)^g}$ has dimension at most $1$.\endproof

\subsubsection{On singular loci}\label{S4208}

We call a subspace $W$ of $\V_g$ {\it singular} if for all $C\in W$ we have $\det(C)=0$.

\begin{lemma}\label{L45}
Let $Q_1,Q_2\in V_g$ be such that $\Span(\{Q_1,Q_2\})$ is singular. Then the intersection $\Ker(Q_1)\cap\Ker(Q_2)$ is non-trivial (i.e., it is not $0$) and there exists $Q\in\Span(\{Q_1,Q_2\})$ such that $\Ker(Q)=\Ker(Q_1)\cap\Ker(Q_2)$.
\end{lemma}

\noindent
{\it Proof.} Let $d_0:=\dim_{\mathbb K}(\Ker(Q_1)\cap\Ker(Q_2))$, $d_1:=\dim_{\mathbb K}(\Ker(Q_1))$, and $d_2:=\dim_{\mathbb K}(\Ker(Q_2))$. We have $d_1,d_2\in\{1,\ldots,g\}$ and $d_0\in\{0,1,\ldots,\min\{d_1,d_2\}\}$. We can assume that $d_1<g$ as the case $d_1=g$, i.e., $Q_1=0$, is trivial.

Up to $\pmb{\SL}_g$-equivalence, we can assume that $Q_1=\Diag(0,\ldots,0,1,\ldots,1)$ (the number of $0$s being $d_1$), that 
$$\Ker(Q_2)=\Span(\{e_1,\ldots,e_{d_0},e_{d_1+1},\ldots,e_{d_1+d_2-d_0}\}),$$ 
and that $Q_2$ has four (possible empty) diagonal blocks which are $0_{d_0}$, $1_{d_1-d_0}$, $0_{d_2-d_0}$, $1_{g-d_1-d_2+d_0}$, where $0_{m}$ denotes the zero matrix of size $m\times m$. 

If $d_0=0$, then $\det(y_1Q_1+Q_2)$ is a monic polynomial in $y_1$ of degree $g-d_1>0$ and thus there exists $y_1\in\mathbb K$ such that $y_1Q_1+Q_2$ is non-singular. This contradicts the fact that $\Span(\{Q_1,Q_2\})$ is singular. Thus $d_0>0$, i.e., the intersection $\Ker(Q_1)\cap\Ker(Q_2)$ is non-trivial.

If $M^{\wr}$ is the matrix obtained from a square $g\times g$ matrix $M$ by removing its first $d_0$ rows and columns, then a similar argument shows that there exist $y_1\in\mathbb Q$ such that $y_1Q_1^{\wr}+Q_2^{\wr}$ is a non-singular matrix of size $(g-d_0)\times (g-d_0)$
and therefore the rank of $Q:=y_1Q_1+Q_2$ is $g-d_0$, i.e., $\Ker(Q)$ has dimension $d_0$ and is equal to $\Ker(Q_1)\cap\Ker(Q_2)=\Span(\{e_1,\ldots,e_{d_0}\})$.\endproof

\begin{prop}\label{P16}
Let $(Q_0,\ldots,Q_r)\in\V_g^{r+1}=\pmb{\V}_g^{r+1}(\mathbb K)$ be such that the $\mathbb K$-vector subspace $\Span(\{Q_0,\ldots,Q_r\})$ of $\V_g$ is singular. Then the following two properties hold:

\medskip
{\bf (a)} The intersection $\cap_{i=0}^r \Ker(Q_i)$ is non-trivial.

\smallskip
{\bf (b)} The closure in $\pmb{\V}_g^{r+1}$ of the $\pmb{\SL}_g$-orbit of $(Q_0,\ldots,Q_r)$ contains the zero vector of $\V_g^{r+1}=\pmb{\V}_g^{r+1}(\mathbb K)$.
\end{prop}

\noindent
{\it Proof.} We first prove (a). As the case $r=0$ is trivial we can assume that $r>0$. We prove by induction on $j\in\{1,\ldots,r\}$ that there exists $\widetilde{Q}_j\in\Span(\{Q_0,\ldots,Q_j\})$ such that we have $0\neq \Ker(\widetilde{Q}_j)=\cap_{i=0}^j \Ker(Q_i)$. The base of the induction for $j=1$ follows from Lemma \ref{L45} applied to $Q_0,Q_1\in\V_g$. For $j\in\{1,\ldots,r\}$, assuming that $\widetilde{Q}_j$ exists, from Lemma \ref{L45} applied to $\widetilde{Q}_j,Q_{j+1}\in\V_g$, we get that there exists a symmetric matrix $\widetilde{Q}_{j+1}\in\Span(\{\widetilde{Q}_j,Q_{j+1}\})\subset\Span(\{Q_0,\ldots,Q_{j+1}\})$ such that we have an identity $0\neq \Ker(\widetilde{Q}_{j+1})=\Ker(\widetilde{Q}_j)\cap\Ker(Q_{j+1})=\cap_{i=0}^{j+1} \Ker(Q_i)$. This ends the induction. Part (a) follows from the existence of $\widetilde{Q}_r$.

To prove part (b), as $\cap_{i=0}^r \Ker(Q_i)$ is non-trivial, up to an $\pmb{\SL}_g$-equivalence, we can assume that $e_g\in\cap_{i=0}^r \Ker(Q_i)$. This implies that the last rows (and columns) of each $Q_i$ are $0$. For $\lambda\in\mathbb K^\times$, for the diagonal matrix 
$$\Lambda_{\lambda}:=\Diag(\lambda,\lambda,\ldots,\lambda,\lambda^{-g+1})\in\pmb{\SL}_g(\mathbb K)$$ 
we compute
$$\Lambda_{\lambda}\cdot (Q_0,\ldots,Q_g)=(\Lambda_{\lambda}Q_0\Lambda_{\lambda}^{\t},\ldots,\Lambda_{\lambda}Q_g\Lambda_{\lambda}^{\t})=(\lambda^2Q_0,\ldots,\lambda^2Q_r),$$
and by taking $\lambda\rightarrow 0$ we conclude that part (b) holds.\endproof

\subsubsection{On $\Theta$-rings of $\pmb{\SL}_g$-invariants}\label{S4209}

The following theorem complements and extends Theorem \ref{T29}.

\begin{thm}\label{T30}
For all integers $r\geq 1$ the following five properties hold:

\medskip
{\bf (a)} The $\mathbb W_{\varepsilon_g}$-graded $\mathbb K$-algebra 
$$\mathbb J^r_{g,\mathbb K}:=\mathbb H_{g,\tot,\mathbb K}^r/(\Theta_{m_0,\ldots,m_r}|(m_0,\ldots,m_r)\in
\Delta(g,r))$$
is local artinian and therefore we have integral $\mathbb K$-monomorphisms
$$\mathbb H_{g,\mathbb K}^{r,\Theta}\rightarrow\mathbb H_{g,\tot,\mathbb K}^r\;\;\;\;\textup{and}\;\;\;\;\mathbb H_{g,\mathbb K}^{r,\Theta}\rightarrow \mathbb H_{g,\mathbb K}^r.$$ 

\smallskip
{\bf (b)} The $\mathbb K$-algebra monomorphisms (\ref{EQ153.3}) are integral.

\smallskip
{\bf (c)} There exists a subset with $D(g,r)$ elements 
$$\mathcal T_0\subset \mathbb H_g^r(1)_{\mathbb K}=\Span(\{\Theta_{m_0,\ldots,m_r}|(m_0,\ldots,m_r)\in\Delta(g,r)\})$$
which is algebraically independent, i.e., the $\mathbb K$-subalgebra 
$$\mathbb K[\mathcal T_0]:=\mathbb K[\Theta|\Theta\in\mathcal T_0]$$ 
of $\mathbb H_{g,\mathbb K}^r$ is a polynomial $\mathbb K$-algebra in $D(g,r)$ indeterminates, and for which the $\mathbb K$-monomorphisms 
$\mathbb K[\mathcal T_0]\rightarrow\mathbb H_{g,\tot,\mathbb K}^r$
and $\mathbb K[\mathcal T_0]\rightarrow\mathbb H_{g,\mathbb K}^r$ are finite.

\smallskip
{\bf (d)} Referring to {\bf (b)}, if $\mathbb H_{g,\tot,\mathbb K}^r$ is Cohen--Macaulay, then the $\mathbb K[\mathcal T_0]$-modules $\mathbb H_{g,\tot,\mathbb K}^r$ and $\mathbb H_{g,\mathbb K}^r$ are free of finite rank. 

\smallskip
{\bf (e)} 
If $r\in\{1,2\}$, then we have $\mathbb H_{g,\mathbb K}^{r,\Theta}=\mathbb K[\mathcal T_0]$ and thus $\mathbb H_{g,\mathbb K}^{r,\Theta}$ is a polynomial $\mathbb K$-algebra (equal for $r=1$ and $p\neq 2$ with $\mathbb H_{g,\mathbb K}^1=\mathbb H_{g,\tot,\mathbb K}^1$, see Theorem \ref{T29} (a)).
\end{thm}

\noindent
{\it Proof.} 
We recall that $\pmb{\V}_g^{r+1}=\Spec(\mathbb K[T,T',\ldots,T^{(r)}])$ and we define
$$\mathbb B_{g,\mathbb K}^r:=\mathbb K[T,T',\ldots,T^{(r)}]/(\Theta_{m_0,\ldots,m_r}|(m_0,\ldots,m_r)\in\Delta(g,r)),$$ 
$$\mathbb N_{g,\mathbb K}^r:=\Ker(\mathbb J_{g,\mathbb K}^r\to \mathbb B_{g,\mathbb K}^r).$$
We consider the $\pmb{\SL}_g$-invariant closed subscheme 
$$\pmb{\V}_g^{r+1,\Theta}:=\Spec(\mathbb B_{g,\mathbb K}^r)\subset \pmb{\V}_g^{r+1}.$$

If $p=0$, then, as $\pmb{\V}_g^{r+1}\rightarrow
\pmb{\V}_g^{r+1}/\pmb{\SL}_g$ is a universal categorical quotient, we have $(\pmb{\V}_g^{r+1,\Theta})^{\pmb{\SL}_g}=\Spec(\mathbb J^r_{g,\mathbb K})$ and thus $\mathbb N_{g,\mathbb K}^r=0$. 

If $p>0$, then, as $\pmb{\V}_g^{r+1}\rightarrow
\pmb{\V}_g^{r+1}/\pmb{\SL}_g$ is surjective (being submersive), the natural morphism $\pmb{\V}_g^{r+1,\Theta}\to \Spec(\mathbb J^r_{g,\mathbb K})$ is also surjective and therefore $\mathbb N_{g,\mathbb K}^r$ is a nilpotent ideal of $\mathbb J_{g,\mathbb K}^r$.

Note that if $(Q_0,\ldots,Q_r)\in\pmb{\V}_g^{r+1,\Theta}(\mathbb K)\subset\V_g^{r+1}$, then, with $y_0,\ldots,y_r$ as indeterminates, we have $\det(\sum_{i=0}^r y_iQ_i)=0$, and hence the $\mathbb K$-vector subspace $\Span(\{Q_0,\ldots,Q_r\})$ of $\V_g$ is singular. Therefore, if $h\in\mathbb J^r_{g,\mathbb K}$ and if $\nu\in\mathbb K$ is its value at the zero vector of $\pmb{\V}_g^{r+1,\Theta}(\mathbb K)\subset\V_g^{r+1}$, from Proposition \ref{P16} (b) we get that $h-\nu$ has value $0$ at each $(Q_0,\ldots,Q_r)\in\pmb{\V}_g^{r+1,\Theta}(\mathbb K)$ and thus its image in $\mathbb J^r_{g,\mathbb K}/\mathbb N^r_{g,\mathbb K}\subset \mathbb B^r_{g,\mathbb K}$ is a nilpotent element, hence $h-\nu$ itself is a nilpotent element. This implies that part (a) holds. 

Part (b) follows from part (a) and inclusions (\ref{EQ153.4}) with $B=\mathbb B$.

As the field $\mathbb K$ is infinite and all $\Theta_{m_0,\ldots,m_r}$s are homogeneous of the same $T$-degree $g$, from part (a) and the homogeneous form of Noether's normalization theorem (see \cite{E}, Thm. 13.3), we get that part (c) holds. 

If $\mathbb H_{g,\tot,\mathbb K}^r$ is Cohen--Macaulay, then so is $\mathbb H_{g,\mathbb K}^r$ (see Lemma \ref{L43} (c)) and this implies that the $\mathbb K[\mathcal T_0]$-modules $\mathbb H_{g,\tot,\mathbb K}^r$ and $\mathbb H_{g,\mathbb K}^r$ are free (see either \cite{Ke}, Cor. 4.2, or \cite{Cho}, Thm.). Thus part (d) holds. 

For $r\in\{1,2\}$ we have $D(g,r)=D(g,r,1)=\binom{g+r}{r}$ and thus we have $\mathbb H_{g,\mathbb K}^{r,\Theta}=\mathbb K[\mathcal T_0]$; therefore part (e) follows from part (c).\endproof

\subsubsection{Applications of the Cohen--Macaulay property}\label{S4210}

We recall that the partially-ordered semiring of Hilbert (Poincar\'e) series is the semiring $\mathbb Z_{\geq 0}[[x]]$ endowed with the following partial order $\preceq$: for two series $\wp_1(x),\wp_2(x)\in\mathbb Z_{\geq 0}[[x]]$ we say that $\wp_1(x)\preceq \wp_2(x)$ if there exists $\wp_0(x)\in\mathbb Z_{\geq 0}[[x]]$ such that $\wp_2(x)=\wp_1(x)+\wp_0(x)$. 

\begin{cor}\label{C21}
Let $g,r\in\mathbb N$. We assume that $\mathbb H_{g,\tot,\mathbb K}^r$ is Cohen--Macaulay. Then the following two properties hold:

\medskip
{\bf (a)} There exists a uniquely determined non-decreasing sequence 
$$(c_p(g,r,l))_{l=1}^{\varrho^-_{p}(g,r)}$$ 
of integers $\geq 2$ with $\varrho^-_{p}(g,r)\in\mathbb N\cup\{0\}$ such that the Hilbert series
$$\wp_{p}(g,r)(x):=\sum_{s=0}^{\infty} D_{p}(g,r,s)x^s=\sum_{s=0}^{\infty} \dim_{\mathbb K}(\mathbb H^r_g(s)_{\mathbb K})x^s$$
of $\mathbb H^r_{g,\mathbb K}$ is the rational polynomial function
$$\frac{1+[D_{p}(g,r,1)-D(g,r)]x+\sum_{l=1}^{\varrho^-_{p}(g,r)} x^{c_p(g,r,l)}}{(1-x)^{D(g,r)}}.$$
Thus 
$$\wp_{p}(g,r)(x)\succeq 
\frac{1+[D_{p}(g,r,1)-D(g,r)]x}{(1-x)^{D(g,r)}}:=\wp_{p}(g,r;1)(x)$$
$$\succeq\frac{1+\left[\binom{g+r}{r}-D(g,r)\right]x}{(1-x)^{D(g,r)}}=\wp_0(g,r;1)(x),$$
and the equality $\wp_{p}(g,r;1)(x)=\wp(g,r;1)(x)$ holds if $p\nmid g!$.
Moreover, for all integers $s\geq \max\{0,c_p(g,r,\varrho^-_{p}(g,r))-D(g,r)+1\}$ we have an inequality 
$$D_{p}(g,r,s)-\binom{D(g,r)+s-1}{D(g,r)-1}-\left[D_{p}(g,r,1)-D(g,r)\right]\binom{D(g,r)+s-2}{D(g,r)-1}\geq 0$$
whose left-hand side is equal to the sum
$$\sum_{l=1}^{\varrho^-_{p}(g,r)}\binom{D(g,r)+s-1-c_p(g,r,l)}{D(g,r)-1}.$$

\medskip
{\bf (b)} Let 
$$\varrho_{p}(g,r):=1+\varrho^-_{p}(g,r)+D_{p}(g,r,1)-D(g,r)\in\mathbb N.$$
If $\mathbb H_{g,\mathbb K}^{r,\perp}\neq 0$, then there exists a uniquely determined non-decreasing sequence 
$$(c_p(g,\odd,r,l))_{l=1}^{\varrho_{p}(g,r)}$$ 
of integers $\geq 1$ such that the Hilbert series (in $\mathbb Z_{\geq 0}[[x^{\frac{1}{2}}]]$)
$$\wp_{p}(g,\tot,r)(x^{\frac{1}{2}}):=\sum_{s\in\mathbb W_0} D_{p}(g,r,s)x^s=\sum_{s\in\mathbb W_0} \dim_{\mathbb K}(\mathbb H^r_g(s)_{\mathbb K})x^s$$
of $\mathbb H^r_{g,\tot,\mathbb K}$ is the rational polynomial function in $x^{\frac{1}{2}}$
$$\wp_{p}(g,r)(x)+\frac{x^{\frac{1}{2}}\sum_{l=1}^{\varrho_{p}(g,r)} x^{c_p(g,\odd,r,l)}}{(1-x)^{D(g,r)}}.$$
\end{cor}

\noindent
{\it Proof.} Let the subset $\mathcal T_0=\{h_1,\ldots,h_{D(g,r)}\}\subset \mathbb H_g^r(1)_{\mathbb K}$ be such that the $\mathbb K$-homomorphism $\mathbb K[\mathcal T_0]=\mathbb K[h_1,\ldots,h_{D(g,r)}]\rightarrow \mathbb H_{g,\tot,\mathbb K}^r$ is finite, see Theorem \ref{T30} (c). Let $h_{D(g,r)+1},\ldots,h_{D_{p}(g,r,1)}\in \mathbb H_g^r(1)_{\mathbb K}$ be such that $\{h_l|1\leq l\leq D_{p}(g,r,1)\}$ is a $\mathbb K$-basis of $\mathbb H_g^r(1)_{\mathbb K}$. We take the non-decreasing sequence $(c_p(g,r,l))_{l=1}^{\varrho^-_{p}(g,r)}$ in $\mathbb N\setminus\{1\}$ such that there exist non-zero polynomials 
$$h_{D_{p}(g,r,1)+l}\in \mathbb H_g^r(c_p(g,r,l))_{\mathbb K}$$ 
indexed by $l\in\{1,\ldots,\varrho^-_{p}(g,r)\}$, with the property that we have a direct sum decomposition of $\mathbb Z_{\geq 0}$-graded $\mathbb K[\mathcal T_0]$-modules 
\begin{equation}\label{EQ155}
\mathbb H^r_{g,\mathbb K}=\mathbb K[\mathcal T_0]\oplus\left( \bigoplus_{l=D(g,r)+1}^{D_p(g,r,1)+\varrho^-_{p}(g,r)} \mathbb K[\mathcal T_0]h_{l}\right),
\end{equation}
by the proof of \cite{Ke}, Cor. 4.2 or \cite{Cho}, Thm. 1 and the fact that $\mathbb H^r_{g,\mathbb K}$ is Cohen--Macaulay (see Lemma \ref{L43} (c)). 

If $\mathbb H_{g,\mathbb K}^{r,\perp}\neq 0$, i.e., $g$ is even and $[\Frac(\mathbb H_{g,\tot,\mathbb K}^r):\Frac(\mathbb H^r_{g,\mathbb K})]=2$, then the $\mathbb K[\mathcal T_0]$-module $\mathbb H^{r,\perp}_{g,\mathbb K}$ is also free (by the Quillen--Suslin's theorem, see \cite{La}, Ch. XXI, Sect. 3, Thm. 3.5) and it has the same rank $\varrho_{p}(g,r)$ as $\mathbb H^r_{g,\mathbb K}$ (see Equations (\ref{EQ153})) and thus there exist non-zero homogeneous polynomials 
$$h_{D_{p}(g,r,1)+\varrho^-_{p}(g,r)+1}\in\mathbb H^r_g(\frac{1}{2}+c_p(g,\odd,r,1))_{\mathbb K}\subset \mathbb H_{g,\mathbb K}^{r,\perp},\ldots,$$
$$h_{D_{p}(g,r,1)+\varrho^-_{p}(g,r)+\varrho_{p}(g,r)}\in\mathbb H^r_g(\frac{1}{2}+c_p(g,\odd,r,\varrho_{p}(g,r)))_{\mathbb K}\subset \mathbb H_{g,\mathbb K}^{r,\perp}$$ 
such that we have a direct sum decomposition

\begin{equation}\label{EQ156}
\mathbb H_{g,\tot,\mathbb K}^r=\mathbb K[\mathcal T_0]
\oplus \left(
\bigoplus_{l=D(g,r)+1}^{D_{p}(g,r,1)+\varrho^-_{p}(g,r)+\varrho_{p}(g,r)} \mathbb K[\mathcal T_0]h_{l}\right)
\end{equation}
of $\mathbb W_0$-graded $\mathbb K[\mathcal T_0]$-modules.

The corollary follows from \cite{Ke}, Cor. 4.2 and the last two direct sum decompositions. \endproof

\medskip

From the direct sum decompositions (\ref{EQ155}) and (\ref{EQ156}) and from Theorem \ref{T30} (d) and (e) we get directly (for $r=1$, cf. also Theorem \ref{T29} (a), and for $r=2$, cf. also Proposition \ref{P15} (b)):

\begin{cor}\label{C22}
For all $g,r\in\mathbb N$, if $\mathbb H_{g,\tot,\mathbb K}^r$ is Cohen--Macaulay, then the following two properties hold:

\medskip
{\bf (a)} The integers $\varrho^-_{p}(g,r)\geq 0$ and $\varrho_{p}(g,r)\geq 1$ of Corollary \ref{C21} have the property that for each $\mathbb K$-algebra $\mathbb K[\mathcal T_0]$
as in Theorem \ref{T30} (c) we have
$$\varrho_{p}(g,r)=\rank_{\mathbb K[\mathcal T_0]}(\mathbb H_{g,\mathbb K}^r)=1+D_{p}(g,r,1)-D(g,r)+\varrho^-_{p}(g,r).$$
If $p\neq 2$ or $g=1$, then for $r=1$ we have
$$\varrho_{p}(g,1)=\rank_{\mathbb H_{g,\mathbb K}^{1,\Theta}}(\mathbb H_{g,\mathbb K}^1)=1.$$
If $p\nmid g!$, then for $r=2$ we have
$$\varrho_{p}(g,2)=\rank_{\mathbb H_{g,\mathbb K}^{2,\Theta}}(\mathbb H_{g,\mathbb K}^2)=1+\varrho^{-}_{p}(g,2).$$

\smallskip
{\bf (b)} If $g$ is even and $\mathbb H_{g,\tot,\mathbb K}^r\neq\mathbb H^r_{g,\mathbb K}$, then we have
$$\rank_{\mathbb K[\mathcal T_0]}(\mathbb H_{g,\tot,\mathbb K}^r)=2\rank_{\mathbb K[\mathcal T_0]}(\mathbb H_{g,\mathbb K}^r)=2\rank_{\mathbb K[\mathcal T_0]}(\mathbb H_{g,\mathbb K}^{r,\perp})=2\varrho_p(g,r).$$
If $p\nmid g!$, then for $r=2$ we have
$$\rank_{\mathbb H_{g,\mathbb K}^{2,\Theta}}(\mathbb H_{g,\tot,\mathbb K}^2)=2\rank_{\mathbb H_{g,\mathbb K}^{2,\Theta}}(\mathbb H_{g,\mathbb K}^2)=2\varrho_{p}(g,2)=2+2\varrho^-_{p}(g,2).$$
\end{cor}

\subsubsection{Transcendence bases}\label{S4211}

It is convenient to work this subsubsection with $\pmb{\SL}'_g$ instead of with 
$\pmb{\SL}_g$. 

Let $(\pmb{\V}_g^{r+1}/\pmb{\SL}'_g)_{\smooth}$ be the largest open subscheme of $\pmb{\V}_g^{r+1}/\pmb{\SL}'_g$ with the property that the morphism $\pmb{\V}_g^{r+1}\rightarrow \pmb{\V}_g^{r+1}/\pmb{\SL}'_g$ is smooth above it, i.e., if $\pmb{\V}_{g,\smooth}^{r+1}$ is the inverse image of $(\pmb{\V}_g^{r+1}/\pmb{\SL}'_g)_{\smooth}$ in $\pmb{\V}_g^{r+1}$, then the induced morphism $\pmb{\V}_{g,\smooth}^{r+1}\to (\pmb{\V}_g^{r+1}/\pmb{\SL}'_g)_{\smooth}$ is smooth. As $\Frac(\mathbb H_{g,\tot,\mathbb K}^r)$ is equal to its purely inseparable closure in $\Frac(\mathbb K[T,\ldots,T^{(r)}])$ (see \cite{DK}, Thm. 2.3.12), the field extension $\Frac(\mathbb H_{g,\tot,\mathbb K}^r)\to\Frac(\mathbb K[T,\ldots,T^{(r)}])$ is separably generated (see \cite{Mat}, Thm. 26.2 and 26.4). Thus $(\pmb{\V}_g^{r+1}/\pmb{\SL}'_g)_{\smooth}$ is non-empty.

We recall that $(\pmb{\V}^{r+1}_{g,\smooth})^{\s}$ is the stable locus of $\pmb{\V}^{r+1}_{g,\smooth}$ and we have $\pmb{\V}_{g,\smooth}^{r+1}/\pmb{\SL}'_g=(\pmb{\V}_g^{r+1}/\pmb{\SL}'_g)_{\smooth}$. Let $F_0=\Theta_0\cdot F_1\cdot J\in\mathcal O(\pmb{\V}_g^2)^{\pmb{\SL}'_g}=\mathcal O(\pmb{\V}_g^2/\pmb{\SL}'_g)$ be as in the proof of Theorem \ref{T25}. 

\begin{lemma}\label{L46}
For integers $r\geq l\geq 2$ we consider the natural $\pmb{\SL}'_g$-morphism 
$$\varpi_g^{r+1,l}:\pmb{\V}_g^{r+1}\rightarrow \pmb{\V}_g^{l}\times_{(\pmb{\V}_g^{l}/\pmb{\SL}'_g)} (\pmb{\V}_g^{r+1}/\pmb{\SL}'_g).$$
The following three properties hold:

\medskip
{\bf (a)} The $\pmb{\SL}'_g$-morphism $\varpi_g^{3,2}$ is generically finite of degree $2^{g-2+\varepsilon_g}$ and it is generically \'etale if and only if $p\neq 2$.
 
 \smallskip
 {\bf (b)} For $l\geq 3$, the $\pmb{\SL}'_g$-morphism $\varpi_g^{r+1,l}$ is generically finite of degree $1$. More precisely, if $\pmb{\X}_g^{l}$ is the intersection of $(\pmb{\V}_{g,\smooth}^{l})^{\s}$ with the non-empty $\pmb{\SL}'_g$-invariant open subscheme of $(\pmb{\V}_g^{l})^{\s}$ whose $\mathbb K$-valued points are those $l$-tuples $(Q_0,\ldots,Q_{l-1})\in\V_g^{l}=\pmb{\V}_g^{l}(\mathbb K)$ for which $F_0(Q_0,Q_1)\neq 0$ and the group $\Stab_{\pmb{\SL}'_g}((Q_0,\ldots,Q_{l-1}))$ is trivial, and if $\pmb{\X}_g^{l}\times\pmb{\V}_g^{r+1-l}$ is the open subscheme of $\pmb{\V}_g^{r+1}$ which is the pullback of $\pmb{\X}_g^{l}$ via the projection morphism $\pmb{\V}_g^{r+1}\rightarrow\pmb{\V}_g^{l}$, then the $\pmb{\SL}'_g$-morphism
$$\varpi_g^{r+1,l,-}:\pmb{\X}_g^{l}\times\pmb{\V}_g^{r+1-l}\rightarrow \pmb{\X}_g^{l}\times_{(\pmb{\V}_g^{l}/\pmb{\SL}'_g)} (\pmb{\V}_{g}^{r+1})/\pmb{\SL}'_g)$$
induced naturally by $\varpi_g^{r+1,l}$ is an isomorphism.\footnote{The $\pmb{\SL}'_g$-invariant morphism $\pmb{\X}_g^{l}\to \pmb{\X}_g^{l}/\pmb{\SL}'_g$ is a geometric quotient (see Fact \ref{F2}) and the $\pmb{\SL}'_g$-action on $\pmb{\X}_g^{l}$ is free. Thus \cite{MFK}, Ch. 0, Prop. 0.9 implies that $\pmb{\X}_g^{l}\to \pmb{\X}_g^{l}/\pmb{\SL}'_g$ is a fiber bundle in the sense of \cite{MFK}, Ch. 0, Def. 0.10.}

\smallskip
{\bf (c)} If $l\geq 3$, the field extension $\Frac(\pmb{\V}_g^l/\pmb{\SL}'_g)\rightarrow \Frac(\pmb{\V}_g^{r+1}/\pmb{\SL}'_g)$ (induced naturally by the projection morphism $\pmb{\V}_g^{r+1}\rightarrow\pmb{\V}_g^{l}$) is purely transcendental of transcendence degree equal to $(r+1-l)\frac{g(g+1)}{2}$. Moreover, a subset of $\Frac(\pmb{\V}_g^{r+1}/\pmb{\SL}'_g)$ with $(r+1-l)\frac{g(g+1)}{2}$ elements is a generating transcendence basis (resp. a transcendence basis) of $\Frac(\pmb{\V}_g^{r+1}/\pmb{\SL}'_g)$ over $\Frac(\pmb{\V}_g^l/\pmb{\SL}'_g)$ if and only if it is a generating transcendence basis (resp. a transcendence basis) of $\mathbb K(T,T',\ldots,T^{(r)}):=\Frac(\mathbb K[T,T',\ldots,T^{(r)}])$ over $\mathbb K(T,T',\ldots,T^{(l)}):=\Frac(\mathbb K[T,T',\ldots,T^{(l)}])$.
\end{lemma}

\noindent
{\it Proof.} 
If $P_2:=(Q_0,Q_1,Q_2)\in (\pmb{\V}^3_{g,\smooth})^{\s}(\mathbb K)$ maps to a $P_1:=(Q_0,Q_1)$ in $(\pmb{\V}^2_{g,\smooth})^{\s}(\mathbb K)$ with $F_0(Q_0,Q_1)\neq 0$, then the fiber $(\varpi_g^{3,2})^{-1}(\varpi_g^{3,2}(P_2))$ is naturally identified with the quotient group scheme 
$$\Stab_{\pmb{\SL}'_g}(P_1)/\Stab_{\pmb{\SL}'_g}(P_2)\simeq\Stab_{\pmb{\SL}_g}(P_1)/\Stab_{\pmb{\SL}_g}(P_2).$$
Based on the proof of Theorem \ref{T25} we can assume that $P_1$ is in $\mathbb K1_g\times \Z_g$ with $\Stab_{\pmb{\SL}_g}(P_1)\simeq \pmb{\mu}_{2,\mathbb K}^{g-1}$. But choosing $Q_2$ to be generic (such as all its $ij$ entries with $1\leq i\leq j\leq g$ are non-zero), we have $\Stab_{\pmb{\SL}'_g}(P_2)=\Spec(\mathbb K)$, equivalently we have $\Stab_{\pmb{\SL}_g}(P_2)=\Ker(\pmb{\SL}_g\to\pmb{\GL}_{\V_g})\simeq\pmb{\mu}_{2,\mathbb K}^{1-\varepsilon_g}$. We conclude that for such generic $Q_2$, the fiber $(\varpi_g^{3,2})^{-1}(\varpi_g^{3,2}(P_2))$ is naturally identified with $\pmb{\mu}_{2,\mathbb K}^{g-2+\varepsilon_g}$ from which part (a) follows. 

Let $P_r=(Q_0,\ldots,Q_r)\in (\pmb{\X}_g^{l}\times\pmb{\V}_g^{r+1-l})(\mathbb K)$; so $P_{l-1}:=(Q_0,Q_1,\ldots,Q_{l-1})$ is in $(\pmb{\X}^l_g)(\mathbb K)$. If $P_r\in (\pmb{\V}^{r+1}_{g,\smooth})^{s}(\mathbb K)$, then the fiber $(\varpi_g^{r+1,l})^{-1}(\varpi_g^{r+1,l}(P_r))$ is naturally identified with the quotient group scheme $\Stab_{\pmb{\SL}'_g}(P_l)/\Stab_{\pmb{\SL}'_g}(P_r)$ and thus it is $\Spec(\mathbb K)$. 

Similarly, if $P_r\notin (\pmb{\V}^{r+1}_{g,\smooth})^{s}(\mathbb K)$, then $[(\varpi_g^{r+1,l})^{-1}(\varpi_g^{r+1,l}(P_r))]_{\red}$ is $\Spec(\mathbb K)$. 

From the last two paragraphs we get that the morphism $\varpi_g^{r+1,l,-}$ is birational and induces a bijection at the level of $\mathbb K$-valued points. As the target of $\varpi_g^{r+1,l,-}$ is normal, being smooth over the normal scheme $(\pmb{\V}_{g}^{r+1})/\pmb{\SL}'_g$, from Zariski's Main Theorem we get that $\varpi_g^{r+1,l,-}$ is an isomorphism. Thus part (b) holds. 

Part (c) is a standard `descent and vector bundle' consequence of part (b). We will only recall the `descent and vector bundle' argument for the first part of part (c). Let $H_l\in\mathbb H_{g,\tot,\mathbb K}^{l-1}$ be such that the pullback to $\pmb{\V}_g^l$ of the principal open $D(H_l)$ of $\Spec(\mathbb H_{g,\tot,\mathbb K}^{l-1})=\pmb{\V}_g^l/\pmb{\SL}'_g$ is contained in $\pmb{\X}_g^l$ and thus it is $(\pmb{\X}_g^l)_{H_l}$. Inverting $H_l$, the $\mathbb H_{g,\tot,\mathbb K,H_l}^{l-1}$-algebra $\mathbb H_{g,\tot,\mathbb K,H_l}^r$ is $\mathbb Z_{\geq 0}$-graded by the total degree in just the upper indexes ${}^{(l)},\ldots,{}^{(r)}$. The fact that this $\mathbb Z_{\geq 0}$-graded structure makes the morphism 
$$\Spec(\mathbb H_{g,\tot,\mathbb K,H_l}^r)\rightarrow \Spec(\mathbb H_{g,\tot,\mathbb K,H_l}^{l-1})$$ 
get the structure of a locally free vector bundle over $\Spec(\mathbb H_{g,\tot,\mathbb K,H_l}^{l-1})$ follows after pullback via the (smooth) faithfully flat morphism $(\pmb{\X}_g^l)_{H_l}\to \Spec(\mathbb H_{g,\tot,\mathbb K,H_l}^{l-1})$ from part (b). This implies that the first part of (c) holds.\endproof

\medskip

For an integer $r\geq 2$ let $\Delta_3(g,r)$ be the subset of $\Delta(g,r)$ formed by all $r+1$-tuples of the form $(a,b,0,\ldots,0,\ldots,0,c,0,\ldots,0)$ and let $\Delta_3^*(g,r)$ be the subset of $\Delta(g,r)$ formed by all $r+1$-tuples of the form $(a,b,0,\ldots,0,\ldots,0,c)$ with $c>0$. We have $\Delta_3(g,2)=\Delta(g,2)$. 

\begin{thm}\label{T31}
Let $r\geq 2$ be an integer. Then the set 
$$\mathcal T_{g,r}:=\{\Theta_{a,b,0,\ldots,0,\ldots,0,c,0,\ldots,0}|(a,b,0,\ldots,0,\ldots,0,c,0,\ldots,0)\in\Delta_3(g,r)\}$$ is a transcendence basis of $\Frac(\mathbb H_{g,\tot,\mathbb K}^r)$ over $\mathbb K$.\end{thm}

\noindent
{\it Proof.} As $D(g,2)=D(g,2,1)=\binom{g+2}{2}$ is the number of elements of the set $\Delta_3(g,2)=\Delta(g,2)$, for $r=2$ we have $\Span(\mathcal T_{g,r})=\Span(\mathcal T_0)$ and the theorem follows from Theorem \ref{T30} (c). 

We assume now that $r>2$. The finite field extension 
\begin{equation}\label{EQ157}
\eta_1:=\Frac(\mathbb H^1_{g,\tot,\mathbb K})\supset\mathbb K(\Theta_0,\Theta_1,\ldots,\Theta_g)=\mathbb K(\Theta_{g,0},\Theta_{g-1,1},\ldots,\Theta_{0,g})
\end{equation} 
is an isomorphism if $p\neq 2$ (see Theorem \ref{T29} (a)).
Let $\eta_2:=\Frac(\mathbb H^2_{g,\tot,\mathbb K})$. As the theorem holds for $r=2$, the field extension $$\eta_1(\Theta_{a,b,c}|(a,b,c)\in\Delta_3^*(g,2))\rightarrow\eta_2$$ 
is finite. From this and Lemma \ref{L46} (a), by denoting for $2\leq l\leq r$
$$\kappa_l:=\mathbb K(T_{ij}, T^{(1)}_{i,j}|1\leq i\leq j\leq g)(\Theta_{a,b,0,\ldots,0,c}|(a,b,0,\ldots,0,c)\in\Delta_3^*(g,l)),$$ 
we get a finite field extension 
$$\kappa_2\rightarrow\mathbb K(T_{ij}, T^{(1)}_{ij},T^{(2)}_{ij}|1\leq i\leq j\leq g)=\Frac(\mathcal O(\pmb{\V}_g^3)).$$ 
Replacing the index $^{(2)}$ by the index $^{(l)}$ with $2<l\leq r$, we get a finite field extension
$$\kappa_l\rightarrow\mathbb K(T_{ij}, T^{(1)}_{ij},T^{(l)}_{ij}|1\leq i\leq j\leq g).$$
From the last two sentences, if $\kappa$ denotes the composite field of $\kappa_2,\ldots,\kappa_r$, we get a finite field extension
$$\kappa\rightarrow\mathbb K(T_{ij}, T^{(1)}_{ij},\ldots,T^{(r)}_{ij}|1\leq i\leq j\leq g).$$ 
As $\kappa$ is $\mathbb K(T_{ij}, T^{(1)}_{i,j}|1\leq i\leq j\leq g)$-generated by the family of $(r-1)\frac{g(g+1)}{2}$ polynomials
$$(\Theta_{\Box})_{\Box\in\sqcup_{l=2}^r \Delta_3^*(g,l)},$$ by reasons of transcendence degrees, we get that this family is algebraically independent over $\mathbb K(T_{ij}, T^{(1)}_{i,j}|1\leq i\leq j\leq g)$ and thus also over $\eta_1$. This implies that the theorem holds for $r>2$.\endproof

\medskip

We define
$$\mathbb K_{g,r}:=\Frac(\mathbb H_{g,\tot,\mathbb K}^r)\supset \Frac(\mathbb H_{g,\mathbb K}^r)=:\mathbb K_{g,\even,r},$$
$$\mathbb K_{g,\infty}:=\bigcup_{r\geq 1} \mathbb K_{g,r}\supset \bigcup_{r\geq 1} \mathbb K_{g,\even,r}=:\mathbb K_{g,\even,\infty},$$
$$\mathbb H_{g,\tot,\mathbb K}:=\bigcup_{r\geq 1} \mathbb H_{g,\tot,\mathbb K}^r\;\;\;\textup{and}\;\;\;\mathbb H_{g,\mathbb K}:=\bigcup_{r\geq 1} \mathbb H_{g,\mathbb K}^r,$$
the unions being taken inside the field of fractions of
$$\mathbb T_g:=\bigcup_{r\geq 0}\mathbb K[T^{(0)},T^{(1)},\ldots,T^{(r)}].$$

We consider the shift $\mathbb K$-monomorphism
$$\sigma:\mathbb T_g\rightarrow\mathbb T_g$$
that sends $\sigma(T^{(r)})=T^{(r+1)}$ for all integers $r\geq 0$ (i.e., $\sigma(T_{ij}^{(r)})=T_{ij}^{(r+1)}$ for all integers $1\leq i\leq j\leq g$ and $r\geq 0$). 

\begin{df}\label{df30}
We say that a subset $\bowtie$ of a subfield $\kappa$ of $\mathbb K_{g,\infty}$ is a {\it $\sigma$-generating set} of $\kappa$ if the subset $\kappa\cap \{\sigma^i(h)|i\geq 0,\; h\in\bowtie\}$ generates $\kappa$ over $\mathbb K$.\footnote{This is a more general notion than the `$\phi$-generate' notion of Subsection \ref{S22} as we do not requite that $\sigma(\kappa)\subset\kappa$.} In such a case, we will also say that $\kappa$ is $\sigma$-generated by $\bowtie$.
\end{df}

\begin{cor}\label{C23}
The following four properties hold:

\medskip
{\bf (a)} The set of distinct polynomials
$$\{\sigma^i(\Theta_{a,b,c})|(a,b,c)\in\Delta_3^*(g,2),i\geq 0\}\subset \mathbb T_g$$
is algebraically independent over $\mathbb K$, i.e., the set $\{\Theta_{a,b,c})_{(a,b,c)\in\Delta_3^*(g,2)}\}$ is $\sigma$-algebraically independent in the sense of Subsection \ref{S22}. 

\smallskip
{\bf (b)} The `big field' $\mathbb K_{g,\infty}$ is $\sigma$-generated by a finite subset of $\mathbb K_{g,3}$, i.e., there exists $N\in\mathbb N$ and $h_1,\ldots,h_{N}\in\mathbb K_{g,3}$ such that we have 
$$\mathbb K_{g,\infty}=\mathbb K(\sigma^i(h_j)|1\leq j\leq N,i\geq 0).$$

\smallskip
{\bf (c)} If $g$ is even, then we have $\mathbb H_{g,\mathbb K}^{3,\perp}\neq 0$. 

\smallskip
{\bf (d)} The `big even field' $\mathbb K_{g,\even,\infty}$ is $\sigma$-generated by a finite subset of $\mathbb K_{g,\even,4}$,
 i.e., there exists $\tilde N\in\mathbb N$ and $\tilde h_1,\ldots,\tilde h_{\tilde N}$ in $\mathbb K_{g,\even,4}$ such that we have
$$\mathbb K_{g,\even,\infty}=\mathbb K(\sigma^i(\tilde h_j)|1\leq j\leq\tilde N,i\geq 0).$$\end{cor}

\noindent
{\it Proof.} For part (a) it suffices to show that for all integers $i\geq 0$
$$\{\sigma^i(\Theta_{a,b,c})|(a,b,c)\in\Delta_3^*(g,2)\}=\{\Theta_{0,0,\ldots,0,a,b,c}|(a,b,c)\in\Delta_3^*(g,2)\}$$ 
(where the number of consecutive $0$s is $i$) is a transcendence basis of $\mathbb K_{g,i+2}$ over $\mathbb K_{g,i+1}$. For $i=0$ this follows from Theorem \ref{T31} applied with $r=2$ and Equation (\ref{EQ157}). For $i\geq 2$, by considering a $\mathbb K$-automorphism $\iota_{0,1;i+1,i+2}$ of 
$$\mathbb K[T^{(0)},T^{(1)},\ldots,T^{(i+2)}]$$ 
which interchanges $T_{ij}^{(0)}\longleftrightarrow T^{(i)}_{ij}$ and $T^{(1)}_{ij}\longleftrightarrow T^{(i+1)}_{ij}$ and which fixes each $T^{(l)}_{ij}$ with $l\in\{2,3,\ldots,i-1,i+2\}$, it suffices to show that, when inserting $i$ consecutive $0$s, the family of polynomials 
$$(\Theta_{a,b,0,\ldots,0,c})_{(a,b,c)\in\Delta_3^*(g,2)}=(\iota_{0,1;i+1,i+2}(\Theta_{0,\ldots,0,a,b,c}))_{(a,b,c)\in\Delta_3^*(g,2)}$$
is a transcendence basis of $\mathbb K_{g,i+2}$ over $\mathbb K_{g,i+1}$, but this follows from Theorem \ref{T31} applied with $r=i+2$. For $i=1$, the same argument applies but by considering the $\mathbb K$-automorphism $\iota_{0,1,2}$ of $\mathbb K[T^{(0)},T^{(1)},T^{(2)},T^{(3)}]$ that circularly maps $T_{ij}^{(0)}\mapsto T^{(1)}_{ij}\mapsto T^{(2)}_{ij}\mapsto T^{(0)}_{ij}$ and fixes each $T_{ij}^{(3)}$.

To prove part (b), let $\{h_1,\ldots,h_{N'}\}\subset\mathbb K_{g,2}$ be a $\sigma$-generating set of $\mathbb K_{g,2}$. We take $N\in\{N',\ldots,N'+\frac{g(g+1)}{2}\}$ to be the smallest number such that there exists a generating transcendence basis of $\mathbb K_{g,3}$ over $\mathbb K_{g,2}$ of the form 
$$\mathcal T_3:=\{h_{N'+1},\ldots,h_N\}\cup\left(
\bigcup_{q=N+1}^{N'+\frac{g(g+1)}{2}}\{\sigma^{i_q}(h_{j_q})\}\right)$$ 
with every $i_q\in\{0,1,2,3\}$ and $j_q\in\{1,\ldots,N'\}$, see Lemma \ref{L46} (c). The fact that $\{h_1,\ldots,h_N\}\subset \mathbb K_{g,3}$ is a $\sigma$-generating set of $\mathbb K_{g,\infty}$ follows from the fact (see Lemma \ref{L46} (c)) that for each integer $i\geq 0$, $\sigma^i(\mathcal T_3)$ is a generating transcendence basis of $\mathbb K_{g,3+i}$ over $\mathbb K_{g,2+i}$.

Let now $g$ be even. Part (c) follows from the inclusion $\{h_1,\ldots,h_N\}\subset \mathbb K_{g,3}$ and the fact 
that $\mathbb H_{g,\mathbb K}^{-1+\frac{g(g+1)}{2},\perp}\neq 0$ (by the existence of $\Upsilon_{1,2,\ldots,\frac{g(g+1)}{2}}$). 

To prove part (d), based on part (c), let $r\in\{3,4\}$ be the smallest integer such that
$\mathbb H_{g,\mathbb K}^{r-1,\perp}\neq 0$. Let $\{\tilde h_1,\ldots,\tilde h_{\tilde N'}\}\subset\mathbb H_{g,\mathbb K}^{r-1}$ be a $\sigma$-generating set of $\Frac(\mathbb H_{g,\mathbb K}^{r-1})$. A `descent and vector bundle' argument similar to the one of the proof of Lemma \ref{L46} (c) shows that for each integer $i\geq 0$, as $\mathbb K_{g,r+i-1}\rightarrow\mathbb K_{g,r+i}$ is a purely transcendental field extension, the field extension $\Frac(\mathbb H_{g,\mathbb K}^{r-1+i})\rightarrow\Frac(\mathbb H_{g,\mathbb K}^{r+i})$ is also purely transcendental, and a subset of $\Frac(\mathbb H_{g,\mathbb K}^{r+i})$ with $\frac{g(g+1)}{2}$ elements is a generating transcendence basis (or a transcendence basis) over $\Frac(\mathbb H_{g,\mathbb K}^{r-1+i})$ if and only if it is so over $\mathbb K_{g,r+i-1}$. The rest of the argument is as in the proof of part (b) and gives that we can choose $\tilde N\in\{\tilde N',\ldots,\tilde N'+\frac{g(g+1)}{2}\}$.\endproof

\begin{rem}\label{R43}
The fact that the `big field' and the `big even field' of Corollary \ref{C23} (b) and (d) are $\sigma$-finitely generated follows from the following general fact of difference algebra applied to the subfields $\mathbb K_{g,\infty}$ and $\mathbb K_{g,\even,\infty}$ of $\Frac(\mathbb T_g)$ (see \cite{Le}, Ch. 4, Thm. 4.4.1): each $\sigma$-subfield of a $\sigma$-finitely generated field is $\sigma$-finitely generated. The novel part of Corollary \ref{C23} (b) and (d) are the bounds $3$ and $4$ for the `order' of the $\sigma$-generators of $\mathbb K_{g,\infty}$ and $\mathbb K_{g,\even,\infty}$ (respectively).
\end{rem}

\subsubsection{Additional invariants}\label{S4212}

Similar to Lemma \ref{L46}, for $p\neq 2$ we consider the natural $\pmb{\SL}_g$-morphism 
$$\theta_g^{3,2}:\pmb{\V}_g^3\rightarrow \pmb{\V}_g^2\times_{(\pmb{\V}_g^2/\pmb{\SL}_g)} \Spec(\mathbb H_{g,\mathbb K}^{2,\Theta}).$$
This makes sense as (see Theorem \ref{T29} (a)) $\pmb{\V}_g^2/\pmb{\SL}_g=\Spec(\mathbb H_{g,\mathbb K}^1)=\Spec(\mathbb H_{g,\mathbb K}^{1,\Theta})$. From Theorem \ref{T30} (c) and Lemma \ref{L46} (a) we get that $\theta_g^{3,2}$ is generically finite of some degree 
$$b_p(g,2)\in\mathbb N,$$
which we can compute as follows. 

\begin{lemma}\label{L47}
We assume that $p\nmid g!$ and $p\neq 2$. If $\mathbb H_{g,\tot,\mathbb K}^2$ is Cohen--Macaulay, then the following two properties hold:

\medskip
{\bf (a)} If $g$ is odd or if $g$ is even but $\mathbb H^{2,\perp}_{g,\mathbb K}\neq 0$, then we have $b_p(g,2)=2^{g-1}\varrho_p(g,2)$.

\smallskip
{\bf (b)} If $g$ is even and $\mathbb H^{2,\perp}_{g,\mathbb K}=0$ (so $g\geq 4$ as $\Upsilon_{1,2,3}\in \mathbb H^{2,\perp}_{2,\mathbb K}$ is non-zero), then we have $b_p(g,2)=2^{g-2}\varrho_p(g,2)$.
\end{lemma}

\noindent
{\it Proof.} We have $b_p(g,2)=2^{g-2+\varepsilon_g}\rank_{\mathbb H_{g,\mathbb K}^{2,\Theta}}(\mathbb H_{g,\tot,\mathbb K}^2)$, see Lemma \ref{L46} (a). From this and Corollary \ref{C22} (a) and (b) applied with $r=2$ we get that lemma holds.\endproof

\medskip
For $p\neq 2$, if the finite field extension $\Frac(\mathbb H_{g,\mathbb K}^{2,\Theta})\to \Frac(\mathbb H_{g,\mathbb K}^2)$ is separable, then we can interpret $b_p(g,2)$ as follows.\footnote{One can introduce the separable version $b^{\sep}_p(g,r)$ of $b_p(g,r)$ which is a divisor of $b_p(g,r)$ and then with no assumption on the finite field extension $\Frac(\mathbb H_{g,\mathbb K}^{2,\Theta})\to \Frac(\mathbb H_{g,\mathbb K}^2)$ we get the same interpretation for $b_p^{\sep}(g,r)$. Moreover, by considering the smooth locus of the morphism $\Spec(\mathbb Z[T,\ldots,T^{(r)}])\rightarrow\Spec(\mathbb H_{g,\tot}^{r,\mathbb Z})$, it is a standard piece of algebraic geometry to get that for $p>>0$ we have $b_p(g,r)=b_p^{\sep}(g,r)$.} For each triple $(D_0,D_1,D_2)\in \mathbb K 1_g\times \Z_g^2$
such that $F_0(D_0,D_1)\neq 0$ (see proof of Theorem \ref{T25} for $F_0$), $D_1=\Diag(\nu_1,\ldots,\nu_g)$ generic in a sense which can be made precise (mainly just that $\prod_{i=1}^g\nu_i\prod_{1\leq i<j\leq g} (\nu_i+\nu_j)\neq 0$), and $D_2$ is also generic in a sense which can be made precise, then $b_p(g,2)$ is the number of $\mathbb K$-valued points of 
$$(\theta_g^{3,2})^{-1}(\theta_g^{3,2}(D_0,D_1,D_2))$$ 
and thus it is the number of solutions $Q_2\in\V_g$ of the system of equations
$$\Theta_{a,b,c}(D_0,D_1,Q_2)=\Theta_{a,b,c}(D_0,D_1,D_2)\;\;\textup{for all} \;\; (a,b,c)\in\Delta(g,2).$$
This system of equations is equivalent to the fact that the identity
\begin{equation}\label{EQ158}
\det(y_0D_0+y_1D_1+y_2Q_2)=\det(y_0D_0+y_1D_1+y_2D_2)
\end{equation}
holds for all $y_0,y_1,y_2\in\mathbb K$. This equation clearly holds if $y_2=0$ and thus we can assume that $y_2\neq 0$ and we denote $\lambda_0:=y_0y_2^{-1}$ and $\lambda_1:=y_0y_2^{-1}$. 

Equation (\ref{EQ158}) is equivalent to the characteristic polynomials equation 
\begin{equation}\label{EQ158.5}
\chi_{\lambda_1D_1+Q_2}(x)=\chi_{\lambda_1D_1+D_2}(x)\;\;\textup{for all}\;\; \lambda_1\in\mathbb K,
\end{equation}
which implies the equation
\begin{equation}\label{EQ159}
\Trace((\lambda_1D_1+Q_2)^i)=\Trace((\lambda_1D_1+D_2)^i)\;\;\textup{for all}\;\; i\in\{1,\ldots,g\}.
\end{equation}
If moreover $p\nmid g!$, then in fact Equations (\ref{EQ158.5}) and (\ref{EQ159}) are equivalent. 

Writing 
\begin{equation}\label{EQ160}
\Trace((\lambda_1D_1+Q_2)^v)=\sum_{m=0}^v \lambda_1^mP_{g,v,m}(\nu_1,\ldots,\nu_g,Q_2),
\end{equation}
each $P_{g,v,m}(w_1,\ldots,w_g,T^{(2)})=P_{g,v,m}(w_1,\ldots,w_g,T^{(2)}_{11},\ldots, T^{(2)}_{gg})$ is a homogeneous polynomial in $g+\frac{g(g+1)}{2}$ indeterminates of total degree $v$ and of partial degree $v-m$ in the index $^{(2)}$. For $m=v$ we have 
$$P_{g,v,v}(w_1,\ldots,w_g,T^{(2)}_{11},\ldots, T^{(2)}_{gg})=\sum_{i=1}^g w_i^v.$$ 
The identities (\ref{EQ159}) are equivalent to the system 
\begin{equation}\label{EQ161}
P_{g,v,m}(\nu_1,\ldots,\nu_g,T^{(2)})=P_{g,v,m}(\nu_1,\ldots,\nu_g,D_2)\;\;\textup{for all}\;\; 0\leq m<v\leq g
\end{equation}
of $\frac{g(g+1)}{2}$ equations in the $\frac{g(g+1)}{2}$ indeterminates $T^{(2)}_{ij}$ with $1\leq i\leq j\leq g$. 

\begin{thm}\label{T32}
If $p=0$, then the following three properties hold:

\medskip
{\bf (a)} For all integers $g\geq 1$ we have $b_0(g,2)=\prod_{i=1}^g i^{1+g-i}.$

\smallskip
{\bf (b)} If $g\geq 3$ is odd or if $g$ is even and $\mathbb H^{2,\perp}_{g,\mathbb K}\neq 0$, then we have
$$\varrho_0(g,2)=\prod_{i=3}^g i^{1+g-i}.$$

\smallskip
{\bf (c)} If $g$ is even and $\mathbb H^{2,\perp}_{g,\mathbb K}=0$ (so $g\geq 4$), then we have
 $$\varrho_0(g,2)=2\prod_{i=3}^g i^{1+g-i}.$$
\end{thm}

\noindent
{\it Proof.} To check this, we can assume that $\mathbb K=\mathbb C$ and that the distinct $\nu_1,\ldots,\nu_g$ are in fact non-zero real numbers with $\prod_{1\leq i<j\leq g} (\nu_i+\nu_j)\neq 0$. 

We claim that if $D_2$ is the zero matrix $0_g$ of size $g\times g$, then the system of equations (\ref{EQ161}) has the unique solution $Q_2=0_g$. 

We first check that the only solution $Q_2=(q_{ij})_{1\leq i,j\leq g}\in\V_g$ of the system of equations (\ref{EQ161}) which has all entries $q_{ij}=q_{ji}$ in $\mathbb R$ is $Q_2=0_g$. For $\lambda_1\in\mathbb R^\times$, as 
$$\chi_{\lambda_1D_1+Q_2}(x)=\chi_{\lambda_1D_1}(x)=\prod_{j=1}^g (x-\lambda_1\nu_i)$$ has $g$ distinct real zeros, we get that there exists an ordered orthonormal basis $(v_{1,\lambda_1},\ldots,v_{g,\lambda_1})$ of $\mathbb R^g$ formed by eigenvectors of $\lambda_1D_1+Q_2$ which correspond to the $g$-tuple of eigenvalues $(\lambda_1\nu_1,\ldots,\lambda_1\nu_g)$. By the compactness of the $g-1$ dimensional sphere of $\mathbb R^g$ there exists a sequence $(\lambda_{1,n})_{n\geq 1}$ of non-zero real numbers that converges to $0$ and such that $(v_{1,\lambda_{1,n}},\ldots,v_{g,\lambda_{1,n}})$ converge to an ordered orthonormal basis $(v_{1,0},\ldots,v_{g,0})$ of $\mathbb R^g$. Passing to the limit $n\rightarrow\infty$ in the identity 
$$(\lambda_{1,n}D_1+Q_2)(v_{j,\lambda_{1,n}})=\lambda_{1,n}\nu_jv_{j,\lambda_{1,n}},$$ 
we get that we have $Q_2(v_{j,0})=0$ for all $j\in\{1,\ldots,g\}$. Thus $Q_2=0_g$. 

If $Q_2=Y+iZ\in\V_g$, with $Y$ and $Z$ real symmetric matrices and with $i\in\mathbb K$ as the usual square root of $-1$ (and not as an index) is a solution of the system of equations (\ref{EQ161}), then the block real matrix
$$
Q_2^+:=\begin{bmatrix}
Y & -Z \\
-Z & -Y
\end{bmatrix}
$$
is a solution of the system of equations (\ref{EQ161}) in which $g$ and $D_1$ are replaced by $2g$ and the block diagonal real matrix $D_1^+:=\left[\begin{array}{rr}
D_1& 0_g\\ 0_g&-D_1
\end{array}
\right]$ (i.e., for $\lambda_1\in\mathbb R^\times$, as $\lambda_1D_1+Q_2$ has $g$ distinct eigenvalues $\lambda_1\nu_1,\ldots,\lambda_1\nu_g$, the matrix $\lambda_1D_1^++Q_2^+$ has $2g$ distinct eigenvalues $\lambda_1\nu_1,\ldots,\lambda_1\nu_g,-\lambda_1\nu_1,\ldots,-\lambda_1\nu_g$ and so the polynomial Equation (\ref{EQ158.5}) holds) and from the prior paragraph we get that $Q_2^+=0_{2g}$, i.e., $Q_2=0_g$. This proves the claim. 

The above claim implies that the system of equations (\ref{EQ161}) has the same number of $\mathbb K$-valued points as 
the projective scheme in the $1+\frac{g(g+1)}{2}$ indeterminates $T^{(2)}_{00}$ and $T^{(2)}_{i,j}$ with $1\leq i\leq j\leq g$ which is given by the homogeneous equations for all $0\leq m<v\leq g$
\begin{equation}\label{EQ162}
P_{g,v,m}(\nu_1,\ldots,\nu_g,T^{(2)})=(T^{(2)}_{00})^{v-m}P_{g,v,m}(\nu_1,\ldots,\nu_g,D_2).
\end{equation}
Choosing $D_2$ to be generic so that the scheme 
$(\theta_g^{3,2})^{-1}(\theta_g^{3,2}(D_0,D_1,D_2))$ is \'etale isomorphic to $\Spec(\mathbb K^{b_0(g,2)})$, from B\'ezout's Theorem (for instance, see \cite{F}, Ch. 12, Thm. 12.3 or see \cite{H}, Ch. 1, Thm. 7.7 applied repeatedly), we get that $b_0(g,2)$ is the product of the degrees of the homogeneous forms of the system of Equation (\ref{EQ162}) and thus part (a) holds. 

Parts (b) and (c) follow from part (a) and Lemma \ref{L47}.\endproof

\begin{cor}\label{C24}
We assume $g\geq 3$. Then for all integers $r\geq 2$, the finite morphism $\Spec(\mathbb H_{g,\mathbb K}^r)\rightarrow\Spec(\mathbb H_{g,\mathbb K}^{r,\Theta})$ is not an isomorphism and thus we have an inequality $\Gamma_p(g,r)>D(g,r,1)$.
\end{cor}

\noindent
{\it Proof.} We first assume that $p=0$ As $g\geq 3$, we have $\varrho_0(g,2)\geq 2$ (see Theorem \ref{T32} (a) and (b)). If $r=2$, then $\mathbb H_{g,\mathbb K}^{r,\Theta}=\mathbb K[\mathcal T_0]$ (see Theorem \ref{T30} (e)) and thus the case $r=2$ follows from the direct sum decomposition (\ref{EQ155}). The passage from $r=2$ to $r>2$ follows from Remark \ref{R42}. 

The case $p>0$ follows from the previous paragraph and the existence of the $\mathbb K$-algebra monomorphisms (\ref{EQ153.3}) and (\ref{EQ153.3}).\endproof

\begin{ex}\label{EX4}
We assume that $g=3$. Our six polynomials in the nine indeterminates $w_1$, $w_2$, $w_3$, $T^{(2)}_{11}$, $T^{(2)}_{12}$, $T^{(2)}_{13}$, $T^{(2)}_{22}$, $T^{(2)}_{23}$, and $T^{(2)}_{33}$ are:
$$
\begin{array}{rcl}
P_{3,1,0} & = & T^{(2)}_{11}+T^{(2)}_{22}+T^{(2)}_{33}=\Trace(T^{(2)}),\\
\ & \ & \ \\
P_{3,2,0} & = & \Trace((T^{(2)})^2),\\
\ & \ & \ \\
P_{3,2,1} & = & 2(w_1T^{(2)}_{11}+w_2T^{(2)}_{22}+w_3T^{(2)}_{33}),\\
\ & \ & \ \\
P_{3,3,0} & = & \Trace((T^{(2)})^3),\\
\ & \ & \ \\
P_{3,3,1} & = & 3w_1[(T^{(2)}_{11})^2+(T^{(2)}_{12})^2+(T^{(2)}_{13})^2]+3w_2[(T^{(2)}_{12})^2+(T^{(2)}_{22})^2+(T^{(2)}_{23})^2]\\
\ & \ & \ \\
\ & \ & 
+3w_3[(T^{(2)}_{13})^2+(T^{(2)}_{23})^2+(T^{(2)}_{33})^2],\\
\ & \ & \ \\
P_{3,3,2} & = & 3(w_1^2T^{(2)}_{11}+
w_2^2T^{(2)}_{22}+w_3^2T^{(2)}_{33}).\end{array}$$
We check that if $D_2=(d_{ij})_{1\leq i,j\leq 3}$ is generic, then for $p\notin\{2,3\}$ the system of equations (\ref{EQ161}) has $12$ solutions. As $\nu_1$, $\nu_2$, and $\nu_3$ are distinct, the Vandermonde $3\times 3$ matrix $V_{\nu_1,\nu_2,\nu_3}$ whose $ij$ entry for $1\leq i,j\leq 3$ is $\nu_i^{j-1}$ is invertible. Thus, as the $3$ equations in $Q_2=(q_{ij})_{1\leq i,j\leq 3}\in\V_g$
$$P_{3,1+i,i}(\nu_1,\nu_2,\nu_3,Q_2)=P_{3,1+i,i}(\nu_1,\nu_2,\nu_3,D_2)$$ 
for $i\in\{0,1,2\}$ can be substituted (as $p\notin\{2,3\}$) by one matrix equation 
$$[q_{11}\;q_{22}\;q_{33}]V_{\nu_1,\nu_2,\nu_3}=[d_{11}\;d_{22}\;d_{33}]V_{\nu_1,\nu_2,\nu_3},$$
we conclude that $q_{11}=d_{11}$, $q_{22}=d_{22}$ and $q_{33}=d_{33}$. We are left to show that by denoting $\alpha:=d_{12}=d_{21}$, $\beta:=d_{13}=d_{31}$, $\gamma:=d_{23}=d_{32}$, $x:=q_{12}=q_{21}$, $y:=q_{13}=q_{31}$, and $z:=q_{23}=q_{32}$, the system of $3$ equations
\begin{equation}\label{EQ163}
x^2+y^2+z^2=\alpha^2+\beta^2+\gamma^2,
\end{equation}
\begin{equation}\label{EQ164}
\nu_1(y^2+z^2)+\nu_2(x^2+z^2)+\nu_3(x^2+y^2)=\nu_1(\beta^2+\gamma^2)+\nu_2(\alpha^2+\gamma^2)+\nu_3(\alpha^2+\beta_2),
\end{equation}
\begin{equation}\label{EQ164.5}
2xyz-d_{33}x^2-d_{22}y^2-d_{11}z^2=2\alpha\beta\gamma-d_{33}\alpha^2-d_{22}\beta^2-d_{11}\gamma^2
\end{equation}
in the $3$ indeterminates $x$, $y$ and $z$ has $12$ solutions. Note that here we have reinterpreted the equations $P_{3,i,0}(\nu_1,\nu_2,\nu_3,Q_2)=P_{3,i,0}(\nu_1,\nu_2,\nu_3,D_2)$ for $i\in\{0,1,2\}$ as a single polynomial identity $\chi_{Q_2}(x)=\chi_{D_2}(x)$ and thus have replaced the equation $P_{3,3,0}(\nu_1,\nu_2,\nu_3,Q_2)=P_{3,3,0}(\nu_1,\nu_2,\nu_3,D_2)$ by the equation $\det(Q_2)-q_{11}q_{22}q_{33}=\det(D_2)-d_{11}d_{22}d_{33}$ (in view of the identities $q_{jj}=d_{jj}$ for $j\in\{1,2,3\}$).
Based on Equation (\ref{EQ163}), Equation (\ref{EQ164}) can be replaced by the equation
\begin{equation}\label{EQ165}
y^2+\nu z^2=\beta^2+\nu\gamma^2,
\end{equation}
where $\nu:=(\nu_3-\nu_1)(\nu_2-\nu_1)^{-1}\in\mathbb K\setminus\{0,1\}$. Based on Equations (\ref{EQ163}) and (\ref{EQ165}), by taking $D_2$ generic so that the $d_{11}$, $d_{22}$, $d_{33}$ are distinct, the Equation (\ref{EQ164.5}) can be replaced, after a scalar multiplication of $(x,y,z,\alpha,\beta,\gamma)$ by an element of $\mathbb K^\times$, by the equation
\begin{equation}\label{EQ166}
xyz-z^2=\alpha\beta\gamma-\gamma^2.
\end{equation}
By denoting 
$$\mathfrak P:=xyz,$$ 
the system of $3$ Equations (\ref{EQ163}), (\ref{EQ165}), (\ref{EQ166}) gives identities
$$z^2=\mathfrak P-\alpha\beta\gamma+\gamma^2,$$
\begin{equation}\label{EQ167}
y^2=-\nu\mathfrak P+\nu\alpha\beta\gamma-\nu\gamma^2+\beta^2,
\end{equation}
\begin{equation}\label{EQ168}
x^2=(\nu-1)\mathfrak P+\alpha^2+\alpha\beta\gamma(1-\nu)+\nu\gamma^2
\end{equation}
and thus $\mathfrak P$ is a root of the polynomial 
\begin{equation}\label{EQ169}
w^2-(w-\alpha\beta\gamma+\gamma^2)[-\nu w+\nu\alpha\beta\gamma-\nu\gamma^2+\beta^2][(\nu-1)w+\alpha^2+\alpha\beta\gamma(1-\nu)+\nu\gamma^2],
\end{equation}
which is a polynomial of degree $3$ in $w$ as $\nu\in\mathbb K^\times\setminus\{0,1\}$. 
It is easy to see that for generic $\alpha$, $\beta$, $\gamma$, $\nu\in\mathbb K$, polynomial (\ref{EQ169}) has $3$ distinct and non-zero solutions and that Equations (\ref{EQ168}) and (\ref{EQ167}) and $z=(xy)^{-1}\mathfrak P$ produce $4$ solutions for each value of $\mathfrak P$. For $p=0$, this implies that $b_0(3,2)=12$. From this and Lemma \ref{L46} (a) we get that for $p=0$ the finite $\mathbb K$-monomorphism 
$$\mathbb H_{3,\mathbb K}^{2,\Theta}\rightarrow\mathbb H_{3,\mathbb K}^2=\mathbb H_{3,\tot,\mathbb K}^2=\Sym((\V_3^3)^*)^{\pmb{\SL}_3}$$
induces at the level of field of fractions a field extension of degree $3=12/4$.
\end{ex}

\subsection{Finer results on quadratics for the case $g=2$}\label{S43}

In this subsection we assume that $g=2$. After we exemplify all prior notions we will include generalizations of the polynomials $\Theta_{2,0},\Theta_{1,1,},\Theta_{0,2},\Upsilon_{0,1,2}$ in order to present a more detailed description of the $\mathbb K$-algebra monomorphism $\mathbb H_{2,\mathbb K}^r\rightarrow\mathbb H_{2,\tot,\mathbb K}^r$ for each integer $r\geq 1$. 

\begin{ex}\label{EX5}
In this paragraph we assume $p\neq 2$. The three polynomials in the five indeterminates $\nu_1$, $\nu_2$, $T^{(2)}_{11}$, $T^{(2)}_{12}$, $T^{(2)}_{22}$ of the proof of Theorem \ref{T32} (a) are:
$$P_{2,1,0}=T^{(2)}_{11}+T^{(2)}_{22}=\Trace(T^{(2)}),$$
$$P_{2,2,0}=(T^{(2)}_{11})^2+2(T^{(2)}_{12})^2+(T^{(2)}_{22})^2=\Trace((T^{(2)})^2)=\Trace(T^{(2)})^2-2\det(T^{(2)}),$$
$$P_{2,2,1}=2(\nu_1T^{(2)}_{11}+\nu_2T^{(2)}_{22}).$$
Thus, if $D_2=(d_{ij})_{1\leq i,j\leq 2}$ is generic in the sense that $d_{12}=d_{21}\neq 0$, then the system of equations (\ref{EQ161}) has $2$ solutions $Q_2=(q_{ij})_{1\leq i,j\leq 2}\in\V_2$ given by $q_{11}=d_{11}$, $q_{22}=d_{22}$, and $q_{12}=q_{21}\in\{\sqrt{d_{12}},-\sqrt{d_{12}}\}$. We get that $b_0(2,2)=2$.\footnote{More generally, we get that $b_p^{\sep}(2,2)=2$ for $p>2$.} As $0\neq \Upsilon_{0,1,2}\in\mathbb H_{2,\mathbb K}^{2,\perp}$, the $\mathbb K$-monomorphism $\mathbb H_{2,\mathbb K}^2\rightarrow \mathbb H_{2,\tot,\mathbb K}^2$ is not birational (see Lemma \ref{L43} (b)). From this and the equality $b_0(2,2)=2$ we get that for $p=0$ the finite $\mathbb K$-monomorphism $\mathbb H_{2,\mathbb K}^{2,\Theta}\rightarrow \mathbb H_{2,\mathbb K}^2$ is birational and thus (as $\mathbb H_{2,\mathbb K}^{2,\Theta}$ is normal being a polynomial $\mathbb K$-algebra, see Theorem \ref{T30} (e)) it is an isomorphism, i.e., for $p=0$ we have
\begin{equation}\label{EQ170}
\mathbb H_{2,\mathbb K}^2=\mathbb H_{2,\mathbb K}^{2,\Theta}=\mathbb K[\Theta_{2,0,0},\Theta_{1,1,0},\Theta_{1,0,1},\Theta_{0,2,0},\Theta_{0,1,1},\Theta_{0,0,2}].
\end{equation}
As $\mathbb H_{2,\mathbb K}^{2,\perp}\neq 0$, its rank as a free $\mathbb H_{2,\mathbb K}^{2,\Theta}$-module for $p=0$ is also $1$ (see Corollary \ref{C22} (b)) and thus is generated by each non-zero homogeneous polynomial of the smallest degree, i.e., for $p=0$ we have
$$\mathbb H_{2,\mathbb K}^{2,\perp}=\mathbb H_{2,\mathbb K}^2\Upsilon_{0,1,2}=\mathbb H_{2,\mathbb K}^{2,\Theta}\Upsilon_{0,1,2}.$$
This implies that for $p=0$ we have
$$\mathbb H_{2,\tot,\mathbb K}^2=\mathbb H_{2,\mathbb K}^2\oplus\mathbb H_{2,\mathbb K}^2\Upsilon_{0,1,2}=\mathbb H_{2,\mathbb K}^{2,\Theta}\oplus \mathbb H_{2,\mathbb K}^{2,\Theta}\Upsilon_{0,1,2}.$$
As $\sigma(\Theta_{2,0,0})=\Theta_{0,2,0}$, $\sigma^2(\Theta_{2,0,0})=\Theta_{0,0,2}$ and $\sigma(\Theta_{1,0,1})=\Theta_{0,1,1}$, the set $\{\Theta_{2,0,0},\Theta_{1,1,0},\Theta_{1,0,1}\}$ is a $\sigma$-generating set of $\mathbb K_{2,\even,2}$ for $p=0$. The proof of Corollary \ref{C23} (d) for $g=2$ and $p=0$ gives that we can take $\tilde N=4$ as it is easy to check that $\{\Theta_{1,0,0,1},\sigma(\Theta_{1,0,1}),\sigma^2(\Theta_{1,1})\}$ is a generating transcendence basis of $\Frac(\mathbb K[T,T',T'',T'''])$ over $\Frac(\mathbb K[T,T',T''])$.
\end{ex}

\subsubsection{Singular loci and homogeneous rings of Grassmannians}\label{S431}

Let $S_2^{r+1}:=\prod_{l=0} ^r S_2^{(l)}$ be the singular locus of $\pmb{\V}_2^{r+1}$. We have 
$$S_2^{(l)}=\Spec(\mathbb K[T_{11}^{(l)},T_{22}^{(l)},T_{12}^{(l)}]/(T_{11}^{(l)}T_{22}^{(l)}-(T_{12}^{(l)})^2))=\Spec(\mathbb K[u_l^2,v_l^2,u_lv_l]),$$
where $u_0,v_0, u_1,v_1,\ldots$ are indeterminates and where the identification 
\begin{equation}\label{EQ171}
\mathbb K[T_{11}^{(l)},T_{22}^{(l)},T_{12}^{(l)}]/(T_{11}^{(l)}T_{22}^{(l)}-(T_{12}^{(l)})^2)
=\mathbb K[u_l^2,v_l^2,u_lv_l]
\end{equation} 
is given by First Isomorphism Theorem applied to the $\mathbb K$-homomorphism 
$$\jmath_l:\mathbb K[T_{11}^{(l)},T_{22}^{(l)},T_{12}^{(l)}]\rightarrow\mathbb K[u_l,v_l]$$ 
that maps $T_{11}^{(l)},T_{22}^{(l)},T_{12}^{(l)}$ respectively to $u_l^2,v_l^2,u_lv_l$. We recall that we have $T_{11}^{(l)}T_{22}^{(l)}-(T_{12}^{(l)})^2=\Theta_{(0,\dots,0,2,0,\ldots,0)}$, where the $2$ in $(0,\dots,0,2,0,\ldots,0)$ is on the $l$-th position. The group $\pmb{\SL}_2$ acts on $\mathbb K[u_l,v_l]$ via its action on $\Span(u_l,v_l)$ given by the rule: $\Lambda\in\pmb{\SL}_2(\mathbb K)$ maps $[u_l\; v_l]^{\t}$ to $\Lambda[u_l\; v_l]^{\t}$. 

We consider the $\mathbb Z_{\geq 0}$-graded $\mathbb K$-algebra 
$$\mathbb Y_r:=\mathbb K[u_0,v_0,\ldots,u_r,v_r],$$ 
where the degree of each $u_l$ and $v_l$ is $1$. The total degree will be called the $uv$-degree and we will speak about the partial $uv$-degree $(m_0,m_1,\ldots,m_r)$ of non-zero polynomials in the span of monomials $\prod_{l=0}^r u_l^{a_l}v_l^{m_l-a_l}$ with each $a_l\in\{0,\ldots,m_l\}$.

The identifications (\ref{EQ171}) are $\pmb{\SL}_2$-invariant and combine together to induce a finite $\pmb{\SL}_2$-morphism of degree $2^{r+1}$
$$\Spec(\mathbb Y_r)\rightarrow S_2^{r+1}=\Spec(\mathbb Y_r^{\alleven})$$
that identifies $S_2^{r+1}$ with the spectrum of the $\mathbb Z_{\geq 0}$-graded $\mathbb K$-subalgebra $\mathbb Y_r^{\alleven}$ of $\mathbb K[u_0,v_0,\ldots,u_r,v_r]$ generated by homogeneous polynomials whose partial $uv$-degrees are $r+1$-tuples of $(2\mathbb Z_{\geq 0})^{r+1}$ and that is associated to the image of the $\mathbb K$-homomorphism
$$\jmath:\mathbb K[T,\ldots,T^{(r)}]\rightarrow \mathbb Y_r$$
which extends each $\jmath_l$ with $l\in\{0,\ldots,r\}$ and which maps homogeneous polynomials of partial $T$-degree $(m_0,\ldots,m_r)$ to homogeneous polynomials of partial $uv$-degree $(2m_0,\ldots,2m_r)$. We have $\Ker(\jmath)=(\det(T),\ldots,\det(T^{(r)}))$. 

For distinct integers $i,j\in\{0,\ldots,r\}$, let 
$$y_{i,j}:=u_iv_j-v_iu_j\in\mathbb Y_r.$$ 
Thus $y_{i,j}=-y_{j,i}$. If $0\leq i<j\leq r$, the $\mathbb K$-homomorphism $\jmath$ maps 
$$\Theta_{(0,\dots,1,0,\ldots,0,1,0,\ldots,0)}=T_{11}^{(i)}T_{22}^{(j)}+T_{22}^{(i)}T_{11}^{(j)}-2T_{12}^{(i)}T_{12}^{(j)},$$ 
where the $1$s in $(0,\dots,1,0,\ldots,0,1,0,\ldots,0)$ are on the $i+1$ and $j+1$, to $$u_i^2v_j^2+v_i^2u_j^2-2u_iv_iu_jv_j=y_{ij}^2.$$

We define
$$\mathcal I_r:=\{(i,j)\in\mathbb Z^2|0\leq i<j\leq r\},$$
$$\mathcal J_r:=\{(i,j,n,s)\in\mathbb Z^4|0\leq i<j<n<s\leq r\},$$
$$\Det_{2,r}:=\{\det(T^{(i)})|i\in\{0,\ldots,r\}\}\;\;\textup{and}$$
$$I_{\Delta(2,r)}:=\{\Theta_{m_0,\ldots,m_r}|(m_0,\ldots,m_r)\in\Delta(2,r)\}.$$ 

We get an identity of $\mathbb Z_{\geq 0}$-graded $\mathbb K$-algebras induced by $\jmath$:
\begin{equation}\label{EQ172}
\mathbb K[T,\ldots, T^{(r)}]/(I_{\Delta(2,r)})=\mathbb Y_r^{\alleven}/(y_{i,j}^2|(i,j)\in\mathcal I_r).
\end{equation}

It is well-known that the $\mathbb K$-subalgebra $\mathbb Y_r^{\pmb{\SL}_2}$ of $\mathbb Y_r$ is generated by the $y_{i,j}$s and is the homogeneous ring of the Grassmannian $\Gr(2,r+1)$ of $\mathbb P^1_{\mathbb K}$ lines in $\mathbb P^r_{\mathbb K}$ (see \cite{HR}, Sect 1., Applic. (a), Case (ii) and Ex. after Cor. 1.9). From this and the Pl\"ucker relations we get that 
$$\mathbb Y_r^{\pmb{\SL}_2}=\mathbb K[y_{i,j}|(i,j)\in\mathcal I_r]
=\mathbb K[Y_{i,j}|(i,j)\in\mathcal I_r]/(Y_{i,j,n,s}|(i,j,n,s)\in\mathcal J_r),$$
where $Y_{i,j}$s are indeterminates, where for all $(i,j,n,s)\in\mathcal J_r$ we define 
$$Y_{i,j,n,s}:=Y_{i,j}Y_{n,s}-Y_{i,n}Y_{j,s}+Y_{i,s}Y_{j,n},$$ and where $y_{i,j}=Y_{i,j}+(Y_{i',j',n,s}|(i',j',n,s)\in\mathcal J_r)$ for all $(i,j)\in\mathcal I_r$. It is convenient to add extra indeterminates $Y_{i,i}$ for all $i\in\{0,\ldots,r\}$. We will speak about (total) $Y$-degrees: each $Y_{i,j}$ has $Y$-degree $1$. Hence each $Y_{i,j,n,s}$ has $Y$-degree $2$. 

For $g=2$ to generalize the polynomials $\Theta_{2,0},\Theta_{1,1,},\Theta_{0,2},\Upsilon_{0,1,2}$ we first introduce the appropriate index sets. For an integer $m$ such that $2\leq m\leq r+1$ let $\mathcal C_{r,m}$ be the set of all $m$-cycles of the permutation group $\Perm(\{0,\ldots,r\})$. Let $\mathcal C_r:=\cup_{m=2}^{r+1} \mathcal C_{r,m}$. We call an $m$-cycle $\omega=(q_1,\ldots,q_m)\in\mathcal C_{r,m}$ increasing or unimodular if it has exactly $m-1$ inversions, i.e, if up to circular rearrangement we have $q_1<q_2<\cdots<q_m$. Let $\mathcal C_{r,m}^{\uni}$ be the subset of $\mathcal C_{r,m}$ formed by all increasing $m$-cycles. Let $\mathcal C_r^{\uni}:=\cup_{m=2}^{r+1} \mathcal C_{r,m}^{\uni}$. Let $\mathcal C_r^{\uni,\even}$ be the subset of $\mathcal C_r^{\uni}$ formed by all cycles of even orders (lengths) and let $\mathcal C_r^{\uni,\odd}$ be the subset of $\mathcal C_r^{\uni}$ formed by all cycles of odd orders (lengths). The inclusion $\mathcal C_{r,m}^{\uni}\subset \mathcal C_{r,m}$ is an equality precisely when $m=2$. 

Let $\mathcal P_r$ be the power set of $\{0,\ldots,r\}$ and let $\mathcal P_{r,m}$ be its subset formed by all subsets with $m$ elements. We have a natural bijection $\mathcal P_{r,m}\rightarrow\mathcal C_{r,m}^{\uni}$ which maps a subset $\{q_1,q_2,\ldots,q_m\}\in\mathcal P_{r,m}$ with $q_1<q_2<\cdots< q_m$ to the increasing $m$-cycle $(q_1,\ldots,q_m)$. This implies that $\mathcal C_{r,m}^{\uni}$ has $\binom{r+1}{m}$ elements and that $\mathcal C_r^{\uni}$ has $\sum_{m=2}^{r+1} \binom{r+1}{m}=2^{r+1}-r-2$ elements. 

For $\omega=(q_1,\ldots,q_m)\in\mathcal C_{r,m}$ and $q_{m+1}:=q_1$, we check that there exists a unique homogeneous polynomial
$$\Xi_{\omega}\in\mathbb H_{2,\tot,\mathbb K}^r\subset\mathbb K[T,\ldots, T^{(r)}]$$
which is of of partial $T$-degree $(m_0,\ldots,m_r)$ with each $m_l\in\{0,1\}$, $m_l$ being $1$ if and only if $l\in\{q_1,\ldots,q_m\}$ (hence it is of $T$-degree $m$), such that the element $\jmath(\Xi_{\omega})\in\mathbb Y_r$ is the circular product
$$\xi_{\omega}:=-\prod_{i=1}^m (u_{q_i}v_{q_{i+1}}-v_{q_i}u_{q_{i+1}})=-\prod_{i=1}^m y_{q_i,q_{i+1}}\in (\mathbb Y_r^{\alleven})^{\pmb{\SL}_2}.$$
The uniqueness follows from the fact that the ideal $(\Det_{2,r})=\Ker(\jmath)$ contains no element of partial $T$-degree $(m_0,\ldots,m_r)$. 

For $p=0$, the existence of $\Xi_{\omega}$ follows from the fact that $\pmb{\V}_2^{r+1}\to \pmb{\V}_2^{r+1}/\pmb{\SL}_2$ is a universal categorical quotient (see Subsubsection \ref{S412}) which implies that
\begin{equation}\label{EQ172.2}
(\mathbb Y_r^{\alleven})^{\pmb{\SL}_2}=\{\mathbb K[T,\ldots,T^{(r)}]/(\Det_{2,r})\}^{\pmb{\SL}_2}=\mathbb H_{2,\tot,\mathbb K}^r/(\Det_{2,r}).
\end{equation} 

Taking $\mathbb K$ to be an algebraic closure of $\mathbb Q$, due to its uniqueness, $\Xi_{\omega}$ is fixed by the Galois group of $\mathbb Q$ and thus, with the notation of Subsection \ref{S4203}, we have $\Xi_{\omega}\in\mathbb H_{2,\tot,\mathbb Q}^{r,\mathbb Z}$. Considering the $\mathbb Z$ version 
$$\jmath^{\mathbb Z}:\mathbb Z[T,\ldots,T^{(r)}]\to\mathbb Y_r^{\mathbb Z}:=\mathbb Z[u_0,v_0,\ldots,u_r,v_r]$$ of $\jmath:\mathbb K[T,\ldots,T^{(r)}]\to\mathbb Y_r$, the only homogeneous polynomials of $\mathbb Q[T,\ldots,T^{(r)}]$ of partial $T$-degree $(m_0,\ldots,m_r)$ which map to $\mathbb Y_r^{\mathbb Z}$ are the ones of $\mathbb Z[T,\ldots,T^{(r)}]$. This implies that we have $\Xi_{\omega}\in\mathbb H_{2,\tot}^{r,\mathbb Z}$. Hence, based on Equation (\ref{EQ153.3}), using scalar extensions we get that $\Xi_{\omega}$ exists even if $p>0$. 

For instance, we have
$$\Xi_{(0,1)}=\Theta_{1,1}\;\; \textup{and}\;\; \Xi_{(0,1,2)}=-\Xi_{(0,2,1)}=\Upsilon_{0,1,2}$$
as well as
$$\xi_{(0,1)}=y_{0,1}^2\;\; \textup{and}\;\; \xi_{(0,1,2)}=-\xi_{(0,2,1)}=y_{0,1}y_{1,2}y_{0,2}.$$

For $\omega=(q_1,\ldots,q_m)\in\mathcal C_{r,m}$, let 
$$\xi_{\omega}^Y:=-\prod_{i=1}^m Y_{q_i,q_{i+1}}\in\mathbb K[Y_{i,j}|(i,j)\in\mathcal I_r].$$ 

It is easy to see that a monomial in the $y_{ij}$s has a partial $uv$-degree in $(2\mathbb Z_{\geq 0})^{r+1}$ if and only if, up to sign, it is a monomial in the $\xi_{\omega}$s with $\omega\in\mathcal C_r$. From this we get that 
\begin{equation}\label{EQ173}
\begin{split}
(\mathbb Y_r^{\alleven})^{\pmb{\SL}_2}=\mathbb K[\xi_{\omega}|\omega\in\mathcal C_r]=\mathbb K[y_{i,j}^2,\xi_{\omega}|(i,j)\in\mathcal I_r,\omega\in\mathcal C_r\setminus\mathcal C_{r,2}].\\
=\mathbb K[\xi_{\omega}^Y|\omega\in\mathcal C_r]/\{(Y_{i,j,n,s}|(i,j,n,s)\in\mathcal J_r)\cap\mathbb K[\xi_{\omega}^Y|\omega\in\mathcal C_r]\}.
\end{split}
\end{equation}

From Equations (\ref{EQ172}) 
 and (\ref{EQ173}) we get natural $\mathbb K$-algebra homomorphisms 
\begin{equation}\label{EQ174}
\begin{split}
\mathbb H_{2,\tot,\mathbb K}^r/(I_{\Delta(2,r)})\xrightarrow{\zeta}\mathbb K[\xi_{\omega}|\omega\in\mathcal C_r]/(y_{i,j}^2|(i,j)\in\mathcal I_r)\;\;\;\;\;\;\;\;\;\;\;\;\;\;\;\;\;\;\;\;\\
=\mathbb K[\xi_{\omega}^Y|\omega\in\mathcal C_r]/\{[(Y_{i,j}^2|(i,j)\in\mathcal I_r)+(Y_{i,j,n,s}|(i,j,n,s)\in\mathcal J_r)]\cap\mathbb K[\xi_{\omega}^Y|\omega\in\mathcal C_r]\}.
\end{split}
\end{equation}
As for $\omega=(q_1,\ldots,q_m)\in\mathcal C_{r,m}$ we have $\Xi_{\omega}\in\mathbb H_{2,\tot,\mathbb K}^r$, $\zeta$ is an epimorphism. Moreover, if $p=0$, then from universal categorical quotient property of the morphism $\pmb{\V}_2^{r+1}\rightarrow\pmb{\V}_2^{r+1}/\pmb{\SL}_2$, we get that $\zeta$ is an isomorphism.

Let $\bar{\xi}_{\omega}\in \mathbb H_{2,\tot,\mathbb K}^r/(I_{\Delta(2,r)})$ be the image of $\Xi_{\omega}$; so $\zeta(\bar{\xi}_{\omega})$ is the image of $\xi_{\omega}$.

\begin{lemma}\label{L48}
Let $m\geq 4$ be an integer. The following three properties hold:

\medskip
{\bf (a)} For an integer $r\geq 3$, $\{q_1,\ldots,q_m\}\in\mathcal P_{r,m}$, and $s\in\{3,\ldots,m-1\}$ we have an identity
$$\xi_{(q_1,q_2,\ldots,q_m)}+\xi_{(q_1,q_{s+1},q_{s+2},\ldots,q_m)}\xi_{(q_2,q_3,\ldots,q_s)}=(-1)^{s}\xi_{(q_1,q_s,q_{s-1},q_{s-2},\ldots,q_2,q_{s+1},q_{s+2},\ldots,q_m)}.$$

\smallskip
{\bf (b)} If $r$ and $\{q_1,\ldots,q_m\}$ are as in part (a), then for $p\neq 2$ we have $\zeta(\bar{\xi}_{(q_1,q_2,q_3,q_4,\ldots,q_m)})=0$.

\smallskip
{\bf (c)} If $r\geq 2$ is an integer and $\{q_1,q_2,q_3\},\{q_1',q_2',q_3'\}\in\mathcal P_{r,3}$, then for $p\neq 2$ we have $\zeta(\bar{\xi}_{(q_1,q_2,q_3)}\bar{\xi}_{(q_1',q_2',q_3')})=0$.
\end{lemma}

\noindent
{\it Proof.} Part (a) just translates in the language of $\xi$s the fact that the product 
$$(y_{q_1,q_2}y_{q_{s+1},q_s}-y_{{q_1},q_{s+1}}y_{q_2,q_s}+y_{q_1,q_s}y_{q_2,q_{s+1}})y_{q_2,q_3}\cdots y_{q_{s-1},q_s}y_{q_{s+1},q_{s+2}}\cdots y_{q_m,q_1}$$ 
is $0$. From part (a) by taking $s=3$ and working modulo $I_{\Delta(2,r)}$ and $(y_{i,j}^2|(i,j)\in\mathcal I_r)$ and from the identity $\bar{\xi}_{(q_2,q_3)}=0$ we get that 
\begin{equation}\label{EQ175}
\zeta(\bar{\xi}_{(q_1,q_2,q_3,q_4,\ldots,q_m)})=-\zeta(\bar{\xi}_{(q_1,q_3,q_2,q_4,\ldots,q_m)}).
\end{equation} 
Applying Equation (\ref{EQ175}) repeatedly by shifting $q_2$ we get that:
\begin{equation}\label{EQ176}
\zeta(\bar{\xi}_{(q_1,q_2,q_3,q_4,\ldots,q_m)})=-\zeta(\bar{\xi}_{(q_1,q_3,q_2,q_4,\ldots,q_m)})=\cdots=(-1)^{m-1}\zeta(\bar{\xi}_{(q_1,q_2,q_3,q_4,\ldots,q_m)}).
\end{equation}
If $m$ is even, from Equation (\ref{EQ176}) we get that $\zeta(\bar{\xi}_{(q_1,q_2,q_3,q_4,\ldots,q_m)})=0$. If $m\geq 5$ is odd, then from part (a) by taking $s=4$, by working modulo $I_{\Delta(2,r)}$ and $(y_{i,j}^2|(i,j)\in\mathcal I_r)$, and by using that part (b) holds in the even case (so $\zeta(\bar{\xi}_{(q_1,q_5,q_6,\ldots,q_m)})=0)$, we get that 
$\zeta(\bar{\xi}_{(q_1,q_2,q_3,q_4,\ldots,q_m)})=\zeta(\bar{\xi}_{(q_1,q_4,q_3,q_2,q_5,q_6,\ldots,q_m)}).$ From Equation (\ref{EQ175}) applied thrice we get that 
$$\zeta(\bar{\xi}_{(q_1,q_4,q_3,q_2,q_5,q_6,\ldots,q_m)})=-\zeta(\bar{\xi}_{(q_1,q_2,q_3,q_4,\ldots,q_m)}).$$ 
Thus, as $p\neq 2$, we get that $\zeta(\bar{\xi}_{(q_1,q_2,q_3,q_4,\ldots,q_m)})=0$. Therefore part (b) holds.

To prove (c), we can assume that $r\geq 5$ and the set $\{q_1,q_2,q_3,q_1',q_2',q_3'\}$ has $6$ elements. Not to carry the upper indexes $'$, to prove (c) we will take $m=6$ and it will suffice to check that $\zeta(\bar{\xi}_{(q_1,q_5,q_6)}\bar{\xi}_{(q_2,q_3,q_4)})=0$. But this follows from part (a) applied with $(s,m)=(4,6)$ and from part (b) applied with $m=6$. Thus part (c) holds.\endproof

\subsubsection{Generators}\label{S432}

\begin{prop}\label{P17}
If $\mathbb H_{2,\tot,\mathbb K}^r/(I_{\Delta(2,r)})\xrightarrow{\zeta}\mathbb K[\xi_{\omega}|\omega\in\mathcal C_r]/(y_{i,j}^2|(i,j)\in\mathcal I_r)$ is an isomorphism and $p\neq 2$ (for instance, if $p=0$), then the following two properties hold:

\medskip
{\bf (a)} The smallest number of homogeneous generators of the $\mathbb Z_{\geq 0}$-graded $\mathbb K$-algebra $\mathbb H_{2,\tot,\mathbb K}^r$ is 
$$\Gamma_p(2,\tot,r)=\frac{(r+1)(r^2+2r+6)}{6}.$$
More precisely, the set 
$$\{\Theta_{m_0,\ldots,m_r}|(m_0,\ldots,m_r)\in\Delta(2,r)\}\cup \{\Xi_{\omega}|\omega\in\mathcal C_{r,3}^{\uni}\}$$ 
is a minimal homogeneous set of generators of $\mathbb H_{2,\tot,\mathbb K}^r$.

\smallskip
{\bf (b)} The set $\{\Theta_{m_0,\ldots,m_r}|(m_0,\ldots,m_r)\in\Delta(2,r)\}$ is a minimal homogeneous set of generators of $\mathbb H_{2,\mathbb K}^r$, i.e., we have $\mathbb H_{2,\mathbb K}^{r,\Theta}=\mathbb H_{2,\mathbb K}^r$ and (with the notation of Proposition \ref{P15}) $\Gamma_{p}(2,r)=D(2,r,1)=\binom{r+2}{2}$.
\end{prop}

\noindent
{\it Proof.} The case $r=1$ follows from Theorem \ref{T29} and the fact that $\mathcal C_{r,3}^{\uni}$ is the empty set. Thus we can assume that $r\geq 2$. 

The set $\Delta(2,r)$ has $D(2,r,1)=\binom{r+2}{2}$ elements and the set $C_{r,3}^{\uni}$ has $\binom{r+1}{3}$ elements. Thus indeed $\Gamma_p(2,\tot,r)$ has $\frac{(r+1)(r^2+2r+6)}{6}$ elements. 

We get a minimal homogeneous set $\Gamma_{\min}$ of generators of $\mathbb H_{2,\tot,\mathbb K}^r$ as follows. We write $\Gamma_{\min}=\sqcup_{s=0}^{\infty} \Gamma_{\min,s}$, where $\Gamma_{\min,s}$ is the subset of $\Gamma_{\min}$ formed by all homogeneous elements of $T$-degree $s$ and is obtained (constructed) inductively as follows. First, $\Gamma_{\min,1}$ is a basis of 
$\Sym((\V_2^{r+1})^*)_1^{\pmb{\SL}_{2}}$ and thus is the empty set (see Proposition \ref{P14}). Second, once $\Gamma_{\min,1},\ldots,\Gamma_{\min,s}$ are obtained (constructed), $\Gamma_{\min,s+1}$ is a basis of a subspace of $\Sym((\V_2^{r+1})^*)_{s+1}^{\pmb{\SL}_{2}}$ which is a direct supplement in $\Sym((\V_2^{r+1})^*)_{s+1}^{\pmb{\SL}_{2}}$ of the intersection $\Sym((\V_2^{r+1})^*)_{s+1}^{\pmb{\SL}_{2}}\cap \mathbb K[\sqcup_{i=1}^s\Gamma_{\min,s}]$. Thus $\Gamma_{\min,2}$ is a basis of $$\Sym((\V_2^{r+1})^*)_2^{\pmb{\SL}_{2}}=\mathbb H_2^r(1)_{\mathbb K}=\Span(\{\Theta_{m_0,\ldots,m_r}|(m_0,\ldots,m_r)\in\Delta(2,r)\})$$ (see Proposition \ref{P15} (b) and the fact that $p\neq 2$ does not divide $2!=2$) and therefore it has $\binom{r+2}{2}$ elements. From Equation (\ref{EQ174}) we get that we can assume that $\Gamma_{\min,s}\subset\{\Xi_{\omega}|\omega\in\mathcal C_{r,s}\}$. From this and Lemma \ref{L48} (b) and (c) we get that for $s>3$ we have $\Gamma_{\min,s}=\emptyset$ (this is automatically so if $r=2$). As $\Xi_{(q_1,q_2,q_3)}=-\Xi_{(q_1,q_3,q_2)}$, we get that $\Gamma_{\min,3}\subset\{\Xi_{\omega}|\omega\in\mathcal C_{r,3}^{\uni}\}$. Clearly 
the natural $\mathbb K$-linear map from
$\Span(\{\xi_{\omega}^Y|\omega\in\mathcal C_{r,3}^{\uni}\})$ to $$\mathbb K[\xi_{\omega}^Y|\omega\in\mathcal C_r]/\{[(Y_{i,j}^2|(i,j)\in\mathcal I_r)+(Y_{i,j,n,s}|(i,j,n,s)\in\mathcal J_r)]\cap\mathbb K[\xi_{\omega}^Y|\omega\in\mathcal C_r]\}$$
is injective. From this and Equation (\ref{EQ174}) we get that $\Gamma_{\min,3}=\{\Xi_{\omega}|\omega\in\mathcal C_{r,3}^{\uni}\}$. So part (a) holds. Part (b) follows from part (a) and Lemma \ref{L48} (c).\endproof

\medskip
For each $(i,j)\in\mathcal I_r$ let $\underline{(i,j)}\in\Delta(2,r)$ be the unique $r+1$-tuple which has $1$ as the $i+1$ and $j+1$ entries and let $\Theta_{\underline{(i,j)}}^Y$ be an indeterminate of degree $1$. 
For each $i\in\{0,\ldots,r\}$ let $\underline{i}\in\Delta(2,r)$ be the unique $r+1$-tuple which has $2$ as the $i+1$ entry and let $\Theta_{\underline{i}}^Y$ be an indeterminate. 

If $\mathbb H_{2,\tot,\mathbb K}^r/(I_{\Delta(2,r)})\xrightarrow{\zeta}\mathbb K[\xi_{\omega}|\omega\in\mathcal C_r]/(y_{i,j}^2|(i,j)\in\mathcal I_r)$ is an isomorphism and $p\neq 2$ (for instance, if $p=0$), then from Proposition \ref{P17} (b) we get that we have a surjective homomorphism
$$\mathfrak q_r:\mathbb K[\Theta_{\underline{\star}}^Y|\star\in\{0,\ldots,r\}\cup\mathcal I_r]\rightarrow \mathbb H_{2,\mathbb K}^r$$
of $\mathbb Z_{\geq 0}$-graded $\mathbb K$-algebras that maps each $\Theta_{\underline{\star}}^Y$ to $\Theta_{\underline{\star}}$. We will speak about the $\Theta^Y$-degrees of homogeneous elements of $\mathbb K[\Theta_{\underline{\star}}^Y|\star\in\{0,\ldots,r\}\cup\mathcal I_r]$. 

\subsubsection{On Hilbert series}\label{S433}

We refer to Corollary \ref{C21} for the uniquely determined non-decreasing sequence $(c_p(2,r,l))_{l=1}^{\varrho_p^-(2,r)}$ of integers $\geq 2$ with $\varrho_p^-(2,r)\in\mathbb N\cup\{0\}$ and for $\varrho_p(2,r)\in\mathbb N$. We recall that $\varrho_0(2,2)=1$, see Theorem \ref{T32} (b).

\medskip

\begin{prop}\label{P18}
With the notation of Corollary \ref{C21}, for all integers $r\geq 3$ the following three properties hold for $p=0$:

\medskip
{\bf (a)} We have $$\varrho_0(2,r)=2^{r-1}(r-2)^{-1}\binom{2r-3}{r}-1-\binom{r-1}{2}$$ (thus $\varrho_0(2,3)=2$, $\varrho_0(2,4)=16$, $\varrho_0(2,5)=105$).

\smallskip
{\bf (b)} We have
\begin{equation}\label{EQ177}
\begin{split}
\frac{(1+x^3)[1+\frac{(r-1)(r-2)x^2}{2}+\sum_{l=1}^{\varrho_0^-(2,r)} x^{2c_0(2,r,l)}]\left[\sum_{a=0}^{\lfloor \frac{r+1}{2}\rfloor} \binom{r+1}{2a}x^a\right]}{(1-x^2)^{2r-1}}\\ 
\preceq \frac{\frac{(1+x)^{2r-1}}{r-1}\left[\sum_{j=1}^{r-1}
\binom{r-1}{j}\binom{r-1}{j-1}x^{j-1}\right]}{(1-x^2)^{2r-1}}.\;\;\;\;\;\;\;\;\;\;\;\;\;\;\;\;\;\;\;\;\;\;\;\;\;
\end{split}
\end{equation}

\smallskip
{\bf (c)} We have inequalities
$$c_0(2,r,\varrho(2,r)-1)\leq \lfloor \frac{3r-4}{2}\rfloor\;\;\;\textup{and}\;\;\;c_0(2,r,\varrho(2,r))\leq \lfloor \frac{3(r-1)}{2}\rfloor$$ (thus $c_0(2,3,1)=2$, $c_0(2,3,2)\in\{2,3\}$, $c_0(2,4,16)\leq 4$, $c_0(2,5,104)\leq 5$, and $c_0(2,5,105)\leq 6$).
\end{prop}

\noindent
{\it Proof.} Let $\mathcal T_0$ be as in Theorem \ref{T30} (c): it is a subset with $3r$ elements of $\mathbb H_2^r(1)_{\mathbb K}=\Span(\{\Theta_{m_0,\ldots,m_r}|(m_0,\ldots,m_r)\in\Delta(2,r)\})$
which is algebraically independent and for which the $\mathbb K$-monomorphisms 
$\mathbb K[\mathcal T_0]\rightarrow\mathbb H_{2,\tot,\mathbb K}^r$
and $\mathbb K[\mathcal T_0]\rightarrow\mathbb H_{2,\mathbb K}^r$ are finite and flat, i.e., the $\mathbb Z_{\geq 0}$-graded quotient $\mathbb K$-algebra 
$\mathbb H_{2,\tot,\mathbb K}^r/(\mathcal T_0)$ has dimension $0$ (equivalently, it is a local artinian ring). Based on Corollary \ref{C22} (b) and the identities 
$$D(2,r,1)-D(2,r)=\binom{r+2}{2}-3r=\binom{r-1}{2},$$ to prove part (a) it suffices to show that 
\begin{equation}\label{EQ178}
\rank_{\mathbb K[\mathcal T_0]}(\mathbb H_{2,\tot,\mathbb K}^r)=2^r(r-2)^{-1}\binom{2r-3}{r}. 
\end{equation}

To check this, we first check that we can assume that we have
$$\mathcal T_0\supset\mathcal D_0:=\Det_{2,r}=\{\Theta_{\underline{i}}|i\in \{0,\ldots,r\}\}=\{\det(T^{(i)})|i\in \{0,\ldots,r\}\}.$$
As $p=0$, $\pmb{\V}_2^{r+1}\rightarrow\pmb{\V}_2^{r+1}/\pmb{\SL}_2$ is a universal categorical quotient, hence:
$$\Spec(\mathbb H_{2,\tot,\mathbb K}^r/(\mathcal D_0))=\Spec((\Sym((\V_2^{r+1})^*)/(\mathcal D_0))^{\pmb{\SL}_2})=(S_2^{r+1})^{\pmb{\SL}_2}.$$
Recall that $\dim((S_2^{r+1})^{\pmb{\SL}_2})$ is equal to $$\dim(\mathbb K[u_l^2,u_lv_l,v_l^2|l\in\{0,\ldots,r\}]^{\pmb{\SL}_2})=\dim(\mathbb Y_r^{\pmb{\SL}_2})=2(r+1)-3=2r-1$$
(one can easily check this based on Fact \ref{F2}).
From this and Theorem \ref{T30} (a) we get that there exists a subset $\mathcal E_0$ with $2r-1$ elements of $\mathbb H_2^r(1)_{\mathbb K}$ such that by taking $\mathcal T_0=\mathcal D_0\cup\mathcal E_0$, the $\mathbb Z_{\geq 0}$-graded $\mathbb K$-algebra 
$\mathbb H_{2,\tot,\mathbb K}^r/(\mathcal T_0)$ has dimension $0$; so the proof of Theorem \ref{T30} (c) applies to give us that the chosen $\mathcal T_0$ has the desired properties.

We identify $\mathbb K[\mathcal T_0]/(\mathcal D_0)=\mathbb K[\mathcal E_0]$. Based on the previous paragraph, to prove that Equation (\ref{EQ178}) holds it suffices to show that 
\begin{equation}\label{EQ179}
\rank_{\mathbb K[\mathcal E_0]}(\mathbb H_{2,\tot,\mathbb K}^r/(\mathcal D_0))=2^r(r-2)^{-1}\binom{2r-3}{r}. 
\end{equation}

We consider the reductive group 
$$\pmb{\G}_{2,r}:=[\pmb{\SL}_2\times (\mathbb Z/2\mathbb Z)_{\mathbb K}^{r+1}]/(\mathbb Z/2\mathbb Z)_{\mathbb K}$$
over $\Spec(\mathbb K)$, where $(\mathbb Z/2\mathbb Z)_{\mathbb K}$ is embedded diagonally in the center of the product $\pmb{\SL}_2\times (\mathbb Z/2\mathbb Z)_{\mathbb K}^{r+1}$. It acts naturally on $\mathbb Y_r$, with the non-trivial element of the $l+1$-th factor of $(\mathbb Z/2\mathbb Z)_{\mathbb K}^{r+1}$ mapping $(u_l,v_l)$ to $(-u_l,-v_l)$ and fixing each $u_i$ and $v_i$ with $i\in\{0,\ldots,l-1,l+1,\ldots,r\}$. From the identity $\mathbb Y_r^{(\mathbb Z/2\mathbb Z)_{\mathbb K}^{r+1}}=\mathbb K[u_l^2,u_lv_l,v_l^2|l\in\{0,\ldots,r\}]$ we deduce new identifications
$$\mathbb Y_r^{\pmb{\G}_{2,r}}=(\mathbb Y_r^{\pmb{\SL}_2})^{(\mathbb Z/2\mathbb Z)_{\mathbb K}^r}=(\mathbb Y_r^{(\mathbb Z/2\mathbb Z)_{\mathbb K}^{r+1}})^{\pmb{\PGL}_2}$$
$$=(\mathbb Y_r^{\alleven})^{\pmb{\SL}_2}=\mathbb H_{2,\tot,\mathbb K}^r/(\mathcal D_0),$$
where we identify $(\mathbb Z/2\mathbb Z)_{\mathbb K}^r$ with the quotient of $(\mathbb Z/2\mathbb Z)_{\mathbb K}^{r+1}$ by its diagonal $(\mathbb Z/2\mathbb Z)_{\mathbb K}$ subgroup. We have a natural direct sum decomposition
\begin{equation}\label{EQ180}
\mathbb Y_r^{\pmb{\SL}_2}=\oplus_{\chi\in\Hom((\mathbb Z/2\mathbb Z)^r,\{-1,1\})} \mathbb Y^{\pmb{\SL}_2}_{r,\chi},
\end{equation}
where the abstract group $(\mathbb Z/2\mathbb Z)^r=(\mathbb Z/2\mathbb Z)^r_{\mathbb K}(\mathbb K)$ acts on $\mathbb Y^{\pmb{\SL}_2}_{r,\chi}$ via the character $\chi$. We translate this decomposition in terms of the parity of the entries of partial $uv$-degrees as follows. For $a\in\{0,\ldots,\lfloor \frac{r+1}{2}\rfloor\}$, let $\Par_{r,a}$ be the set of $r+1$-tuples $(o_0,\ldots,o_r)\in\{\even,\odd\}^{r+1}$ with the property that the number of times $\odd$ shows up in the sequence $o_0,\ldots,o_r$ is exactly $2a$; it has $\binom{r+1}{2a}$ elements. The set 
$$\Par_r:=\bigcup_{a=0}^{\lfloor \frac{r+1}{2}\rfloor} \Par_{r,a}$$
has exactly $\sum_{a=0}^{\lfloor \frac{r+1}{2}\rfloor} \binom{r+1}{2a}=2^r$ elements. For $\underline{o}=(o_0,\ldots,o_r)\in\Par_r$, let 
$$\mathbb Y^{\pmb{\SL}_2}_{r,\underline{o}}$$ 
be the direct summand of $\mathbb Y^{\pmb{\SL}_2}$ formed by homogeneous polynomials of partial $uv$-degree $(m_0,\ldots,m_r)$ such that for each index $l\in\{0,\ldots,r\}$, the parity of $m_l$ is exactly $o_l$. Then for each $\chi\in\Hom((\mathbb Z/2\mathbb Z)^r,\{-1,1\})$, there exists a unique element $o_{\chi}\in\Par_r$ such that we have
$$\mathbb Y^{\pmb{\SL}_2}_{r,\chi}=\mathbb Y^{\pmb{\SL}_2}_{r,o_{\chi}}.$$
For instance, if $\chi_0\in\Hom((\mathbb Z/2\mathbb Z)^r,\{-1,1\})$ is the trivial character, then $o_{\chi_0}=(\even,\even,\even,\ldots,\even)$ and we have
$$\mathbb H_{2,\tot,\mathbb K}^r/(\mathcal D_0)=\mathbb Y^{\pmb{\SL}_2}_{r,\chi_0}=\mathbb Y^{\pmb{\SL}_2}_{r,(\even,\even,\ldots,\even)}.$$
Each $\mathbb Y^{\pmb{\SL}_2}_{r,\chi}$ is a non-zero torsion free $\mathbb H_{2,\tot,\mathbb K}^r/(\mathcal D_0)$-module. For instance, 
$$0\neq \prod_{i=0}^{a-1} y_{2i,2i+1}\in\mathbb Y^{\pmb{\SL}_2}_{r,(\odd,\odd,\ldots,\odd,\even,\ldots,\even)},$$
where $\odd$ shows up exactly $2a$ times, and there exist no homogeneous elements of $\mathbb Y^{\pmb{\SL}_2}_{r,(\odd,\odd,\ldots,\odd,\even,\ldots,\even)}$ of positive $uv$-degree less than $2a$ (equivalently, $Y$-degree less than $a$). This implies that 
\begin{equation}\label{EQ181}
\dim_{\Frac(\mathbb H_{2,\tot,\mathbb K}^r/(\mathcal D_0))}(\Frac(\mathbb H_{2,\tot,\mathbb K}^r/(\mathcal D_0))\otimes_{\mathbb H_{2,\tot,\mathbb K}^r/(\mathcal D_0)} \mathbb Y^{\pmb{\SL}_2}_{r,\chi})=1.
\end{equation}
Quillen--Suslin's theorem implies that each $\mathbb Y^{\pmb{\SL}_2}_{r,\chi}$ is a free $\mathbb K[\mathcal E_0]$-module. From this and Equation (\ref{EQ181}) we get that for $\chi\in\Hom((\mathbb Z/2\mathbb Z)^r,\{-1,1\})$ we have
\begin{equation}\label{EQ182}
\rank_{\mathbb K[\mathcal E_0]}(\mathbb H_{2,\tot,\mathbb K}^r/(\mathcal D_0))=\rank_{\mathbb K[\mathcal E_0]}(\mathbb Y^{\pmb{\SL}_2}_{r,\chi}).
\end{equation}
Let $q(r):=\rank_{\mathbb K[\mathcal E_0]}(\mathbb Y_r^{\pmb{\SL}_2})$. From Equations (\ref{EQ180}) and (\ref{EQ182}) we get that 
$$q(r)=2^r\rank_{\mathbb K[\mathcal E_0]}(\mathbb H_{2,\tot,\mathbb K}^r/(\mathcal D_0)).$$
Thus to prove that Equation (\ref{EQ179}) holds it suffices to show that we have
$$q(r)=2^{2r}(r-2)^{-1}\binom{2r-3}{r}.$$

To check this we recall (see \cite{GW}, Sect. 6 and \cite{Br}, Thm. 1) that the Hilbert series of the $\mathbb Z_{\geq 0}$-graded $\mathbb K$-algebra $\mathbb Y_r^{\pmb{\SL}_2}$ which is the homogeneous coordinate ring of the Grassmannian $\Gr(2,r+1)$ and in which all $y_{i,j}$s with $0\leq i<j\leq r$ are of degree $1$ (to be called as well $Y$-degree), is
\begin{equation}\label{EQ183}
\wp_1(2,r)(x):=\frac{\sum_{j=1}^{r-1}
(r-1)^{-1}\binom{r-1}{j}\binom{r-1}{j-1}x^{j-1}}{(1-x)^{2r-1}}.
\end{equation}
The composite $\mathbb K$-monomorphism $\mathbb K[\mathcal E_0]\rightarrow\mathbb H_{2,\tot,\mathbb K}^r/(\mathcal D_0)\rightarrow\mathbb Y_r^{\pmb{\SL}_2}$ is integral and thus the $\mathbb K[\mathcal E_0]$-module $\mathbb Y_r^{\pmb{\SL}_2}$ is free (see \cite{Ke}, Cor. 4.2). If we have an ordered $\mathbb K[\mathcal E_0]$-basis $w_1,\ldots,w_{q(r)}$ of $\mathbb Y_r^{\pmb{\SL}_2}$ formed by homogeneous elements such that the sequence $(\deg(w_l))_{1\leq l\leq q(r)}$ of $Y$-degrees is non-decreasing, then we have 
\begin{equation}\label{EQ184}
\wp_1(2,r)(x)=\frac{\sum_{l=1}^{q(r)} x^{\deg(w_l)}}{(1-x^2)^{2r-1}}.
\end{equation}
[We recall that for $(i,j)\in\mathcal I_r$ we have $\jmath(\Theta_{\underline{(ij)}})=y_{ij}^2$, and thus each element of $\mathcal E_0$ has $Y$-degree $2$; also, it has $T$-degree $2$.]

From Equations (\ref{EQ183}) and (\ref{EQ184}) we get that 
\begin{equation}\label{EQ185}
\sum_{l=1}^{q(r)} x^{\deg(w_l)}=\left[\sum_{j=1}^{r-1}(r-1)^{-1}\binom{r-1}{j}\binom{r-1}{j-1}x^{j-1}\right](1+x)^{2r-1}.
\end{equation}
By taking $x=1$ in this equation we get that
$$q(r)=2^{2r-1}(r-1)^{-1}(\sum_{j=1}^{r-1}\binom{r-1}{j}\binom{r-1}{j-1})=2^{2r}[2(r-1)]^{-1}\binom{2(r-1)}{r-2},$$
where we used that $\sum_{j=1}^{r-1}\binom{r-1}{j}\binom{r-1}{j-1}=\sum_{j=1}^{r-1}\binom{r-1}{j}\binom{r-1}{r-j}$ is the coefficient of $x^r$ in the product $(1+x)^{r-1}(1+x)^{r-1}=(1+x)^{2(r-1)}$ and hence it is $\binom{2(r-1)}{r}$. Thus $q(r)=2^{2r}(r-2)^{-1}\binom{2r-3}{r}$, i.e., Equation (\ref{EQ178}) holds. So part (a) holds.

To check that parts (b) and (c) hold, we first remark that the leading term of the right-hand side of the Equation (\ref{EQ185}) is $x^{2r-1}x^{r-2}=x^{3r-3}$. Thus $\deg(w_{q(r)})=3r-3$ and $\deg(w_{q(r)-1})\leq 3r-4$. 

As $\Perm(\{0,\ldots,r\})$ acts naturally on $\mathbb Y_r$ in such a way that for all $a\in\{0,\ldots,\lfloor \frac{r+1}{2}\rfloor\}$ it permutes transitively $\Par_{r,a}$ as well as the set $\{\mathbb Y^{\pmb{\SL}_2}_{r,\chi}|o_{\chi}\in\Par_{r,a}\}$, if $o_{\chi}\in\Par_{r,a}$, then the Hilbert series $\wp_{1,a}(2,r)(x)$ of $\mathbb Y^{\pmb{\SL}_2}_{r,\chi}$ depends, as the notation suggests, only on $a$. We conclude that:
$$\wp_1(2,r)(x)=\sum_{a=0}^{\lfloor \frac{r+1}{2}\rfloor} \binom{r+1}{2a}\wp_{1,a}(2,r)(x).$$
If $o_{\chi}\in\Par_{r,a}$, then 
each element of a homogeneous $\mathbb K[\mathcal E_0]$-basis of 
$\mathbb Y^{\pmb{\SL}_2}_{r,\chi}$ has $Y$-degree at least $a$ and therefore, we have
\begin{equation}\label{EQ186}
x^a\wp_{1,0}(2,r)(x)\preceq\wp_{1,a}(2,r)(x).
\end{equation}
The Hilbert series of $\mathbb H^r_{2,\mathbb K}/(\mathcal D_0)$ but viewed as a $2\mathbb Z_{\geq 0}$-graded $\mathbb K$-algebra (as the $y_{i,j}^2$s have $Y$-degree $2$ or equivalently as the $\Theta$s have $T$-degree $2$) is $(1-x^2)^{r+1}\wp_0(2,r)(x^2)$. The Hilbert series $\wp_{1,0}(2,r)(x)$ is the Hilbert series of $\mathbb H^r_{2,\tot,\mathbb K}/(\mathcal D_0)$ (viewed with the $T$-degree grading which is compatible with the one just mentioned for $\mathbb H^r_{2,\mathbb K}/(\mathcal D_0)$) and thus is the sum of $(1-x^2)^{r+1}\wp_0(2,r)(x^2)$ and of the Hilbert series $\wp_{1,0}^{\perp}(2,r)(x)$ of $\mathbb H^{r,\perp}_{2,\mathbb K}/\mathcal D_0H^{r,\perp}_{2,\mathbb K}$. Similarly to the inequality (\ref{EQ186}) we argue that
\begin{equation}\label{EQ187}
x^3(1-x^2)^{r+1}\wp_0(2,r)(x^2)\preceq\wp_{1,0}^{\perp}(2,r)(x).
\end{equation}

Combining Equation (\ref{EQ185}) and inequalities (\ref{EQ186}) and (\ref{EQ187}) with Corollary \ref{C21} we get that part (b) holds. From part (b) and the two relations $\deg(w_{q(r)})=3r-3$ and $\deg(w_{q(r)-1})\leq 3r-4$ we get the inequalities $2c_0(2,r,\varrho_0^-(2,r))\leq 3r-3$ and $2c_0(2,r,\varrho_0^-(2,r)-1)\leq 3r-4$; so part (c) holds.\footnote{Due to the presence of the denominator $(1-x^2)^ {2r-1}$ in part (b), one cannot conclude that the inequality $2c_0(2,r,\varrho^-_0(2,r))+3+\lfloor \frac{r+1}{2}\rfloor\leq 3r-3$ holds, and in fact this inequality is false if $r\in\{3,4\}$.}\endproof

\begin{ex}\label{EX6}
We take $r=3$ and continue to assume that $p=0$. As we have $\dim(\mathbb H_{2,\mathbb K}^3)=9$ and as $\mathbb H_{2,\mathbb K}^3$ is generated (see Proposition \ref{P17} (b)) by the $10$ polynomials $\Theta_{(m_0,m_1,m_2,m_3)}\in \mathbb H_2^3(1)_{\mathbb K}$ indexed by the elements $(m_0,m_1,m_2,m_3)\in\Delta(2,3)$, the $\mathbb K$-epimorphism $\mathfrak q_r$ induces an identification
$$\mathbb H_{2,\mathbb K}^3\simeq\mathbb K[\Theta^Y_{(m_0,m_1,m_2,m_3)}|(m_0,m_1,m_2,m_3)\in\Delta(2,3)]/(\Pi_4(0,1,2,3))$$
where $\Pi_4(0,1,2,3)$ is a homogeneous polynomial in the $10$ indeterminates of $\Theta^Y$-degree equal to $\rank_{\mathbb H_{2,\mathbb K}^{3,\Theta}}(\mathbb H_{2,\mathbb K}^3)=2+\varrho_0(2,3)=4$, uniquely determined up to multiplication by elements of $\mathbb K^\times$. Thus, for $r=3$ we get that 
$$\Ker(\mathfrak q_3)=(\Pi_4(0,1,2,3)).$$ 
The Hilbert series of a hypersurface of degree $4$ in $\mathbb P^9_{\mathbb K}$ is 
$$\lambda(2,3)(x)=\frac{1-x^4}{(1-x)^{10}}=\frac{1+x+x^2+x^3}{(1-x)^9}$$
and thus we have $c_0(2,3,1)=2$ and $c_0(2,3,2)=3$ (as suggested by Proposition \ref{P18} (b)). 
\end{ex}

\begin{cor}\label{C25}
For all integers $r\geq 3$ and $m\in\{4,\ldots,r+1\}$ the following two types of relations hold for $p=0$:

\medskip
{\bf (a) ($\Xi$ cyclic relations)} For all $(\{q_1,\ldots,q_m\},s)\in\mathcal P_{r,m}\times\{3,\ldots,m-1\}$ we have the following identity
$$\Xi_{(q_1,q_2,\ldots,q_m)}+\Xi_{(q_1,q_{s+1},q_{s+2},\ldots,q_m)}\Xi_{(q_2,q_3,\ldots,q_s)}=(-1)^{s}\Xi_{(q_1,q_s,q_{s-1},\ldots,q_2,q_{s+1},q_{s+2},\ldots,q_m)}$$
between elements of $H^r_{2,\tot,\mathbb K}$.

\smallskip
{\bf (b) ($\Xi$ Pl\"ucker relations)} For $(i,j,n,s)\in\mathcal J_r$, let $\Pi^{(2)}_4(i,j,n,s)$ be
$$(\Theta_{\underline{(i,j)}}^Y)^2(\Theta_{\underline{(n,s)}}^Y)^2+(\Theta_{\underline{(i,n)}}^Y)^2(\Theta_{\underline{(j,s)}}^Y)^2+(\Theta_{\underline{(i,n)}}^Y)^2(\Theta_{\underline{(j,s)}}^Y)^2$$
and let $\Pi^{(4)}_4(i,j,n,s)$ be 
$$2\Theta_{\underline{(i,j)}}^Y\Theta_{\underline{(i,n)}}^Y\Theta_{\underline{(n,s)}}^Y\Theta_{\underline{(j,s)}}^Y+2\Theta_{\underline{(i,j)}}^Y\Theta_{\underline{(i,n)}}^Y\Theta_{\underline{(n,s)}}^Y\Theta_{\underline{(j,n)}}^Y+2\Theta_{\underline{(i,n)}}^Y\Theta_{\underline{(i,s)}}^Y\Theta_{\underline{(j,s)}}^Y\Theta_{\underline{(j,n)}}^Y.$$
Then there exists a homogeneous polynomial 
$$\Pi^{\det}_4(i,j,n,s)\in(\Theta_{\underline{\star}}^Y|\star\in\{i,j,n,s\})$$ 
of the same partial degrees as $\Pi_4^{(2)}(i,j,n,s)$ or $\Pi_4^{(4)}(i,j,n,s)$ and such that 
$$\Pi_4(i,j,n,s):=\Pi_4^{(2)}(i,j,n,s)+\Pi_4^{(4)}(i,j,n,s)+\Pi^{\det}_4(i,j,n,s)\in\Ker(\mathfrak q_r).$$
\end{cor}

\noindent
{\it Proof.} Part (a) follows from Lemma \ref{L48} (a). Part (b) is a direct consequence of the Pl\"ucker relations. More precisely, starting from the Pl\"ucker relation $y_{i,j}y_{n,s}-y_{i,n}y_{j,s}+y_{i,s}y_{j,n}=0$ we first get that 
$$y_{i,n}^2y_{j,s}^2=y_{i,j}^2y_{n,s}^2+2y_{i,j}y_{n,s}y_{i,s}y_{j,n}+y_{i,s}^2y_{j,n}^2,$$
and second that 
$$4y_{i,j}^2y_{n,s}^2y_{i,s}^2y_{j,n}^2=(y_{i,n}^2y_{j,s}^2-y_{i,j}^2y_{n,s}^2-y_{i,s}^2y_{j,n}^2)^2,$$
and thus finally that 
$$y_{i,j}^4y_{n,s}^4+y_{i,n}^4y_{j,s}+y_{i,n}^4y_{j,s}-2y_{i,j}^2y_{i,n}^2y_{n,s}^2y_{j,s}^2-2y_{i,j}^2y_{i,n}^2y_{n,s}^2y_{j,n}-2y_{i,n}^2y_{i,s}^2y_{j,s}^2y_{q_1,n}=0.$$
To prove part (b) we can assume that $r=3$; so $(i,j,n,s)=(0,1,2,3)$. For $r=3$, part (b) follows once we recall that $y_{i',j'}^2=\xi_{(i',j')}$ which implies that $\Pi^{(2)}_4(i,j,n,s)+\Pi^{(4)}_4(i,j,n,s)$ modulo the ideal $(\Theta_{\underline{\star}}^Y|\star\in\{0,1,2,3\})+\Ker(\mathfrak q_r)$ is $0$ (cf. also Equations (\ref{EQ172.2}) and (\ref{EQ173})).
\endproof

\begin{ex}\label{EX7}
We take $r\geq 3$. For $(i,j,n,s)\in\mathcal J_r$, let $\Pi_4(i,j,n,s)$ be as in Corollary \ref{C25} (b): their $T$-degrees are $8$, their partial $T$-degrees $(m_0,\ldots,m_r)$ are such that $m_l\in\{0,2\}$, and $m_l$ is $0$ if and only if $l\in\{0,\ldots,r\}\setminus\{i,j,n,s\}$, and their $\Theta^Y$-degrees are $4$. As $\Ker(\mathfrak q_3)=(\Pi_4(0,1,2,3))$, it follows that $\Ker(\mathfrak q_r)$ contains no non-zero homogeneous polynomials of $T$-degree $4$ (so of $\Theta^Y$-degree $2$) and all homogeneous polynomials of $T$-degree $6$ or $8$ (so of $\Theta^Y$-degrees $3$ or $4$) must either involve $5$ indeterminates among the $T,\ldots,T^{(r)}$s or must belong to the $\mathbb K$-vector space $\Span(\{\Pi_4(i,j,n,s)|(i,j,n,s)\in\mathcal J_r\})$. 
This implies that for all $r\geq 3$ we have
$$D(2,r,2)=\binom{1+\frac{(r+1)(r+2)}{2}}{2},$$
$$D(2,r,3)\leq \binom{2+\frac{(r+1)(r+2)}{2}}{3},$$
$$D(2,r,4)\le\binom{3+\frac{(r+1)(r+2)}{2}}{4}-\binom{r+1}{4}.$$
\end{ex} 
In particular, as $\varrho_0(2,4)=16$ and $c_0(2,4,16)\leq 4$ (see Proposition \ref{P18} (c)) and as $D(2,4,1)=15$, $D(2,4,2)=120$, $D(2,4,3)=680-\alpha$, and $D(2,4,4)=3,055-\beta$, with $\alpha,\beta\in\mathbb N\cup\{0\}$, we easily get that 
$$\wp_0(2,4)(x)=\frac{1+3x+6x^2+(10-\alpha)x^3+\alpha x^4}{(1-x)^{12}},$$
and thus $c_0(2,4,1)=\cdots=c_0(2,4,6)=2,\;\; c(2,4,7)=\cdots =c_0(2,4,16-\alpha)=3$, $c_0(2,4,17-\alpha)=\cdots =c_0(2,4,16)=4$, $\alpha\in\{0,1,\ldots,10\}$ and $\beta=10+11\alpha$. 

Representation theory gives an explicit way to compute Hilbert series. For instance, \cite{DK}, Cor. 4.6.9 and Rm. 4.6.10 allow us to compute $\wp_0(2,r)(x^2)$ as the even part of the coefficient of $z^0$ of the fraction
$$\frac{1-z^2}{(1-xz^{-2})^{r+1}(1-x)^{r+1}(1-xz^{2})^{r+1}}$$
and thus of the fraction
$$\frac{\sum_{i=0}^{\infty} [x^{2i}\binom{r+i}{i}^2- x^{2i+1}\binom{r+i}{i}\binom{r+i+1}{i+1}]}{(1-x)^{r+1}}.$$
A simple calculation shows that for $r=4$ this fraction is a Hilbert series which modulo $x^9\mathbb Z_{\geq 0}[[x]]$ is equal to
$$1+15x^2+10x^3+12-x^4+126x^5+680x^6+855x^7+3,045x^8.$$
In other words, we have $\alpha=0$, $\beta=10$, and
$$\wp_0(2,4)(x)=\frac{1+3x+6x^2+10x^3}{(1-x)^{12}}.$$

As $\beta=10$, for each $(q_0,q_1;q_2,q_3,q_4)$ with $\{q_0,\ldots,q_4\}=\{0,1,2,3,4\}$, $q_0<q_1$, $q_2<q_3<q_4$, there exists a non-zero homogeneous polynomial 
$$\Pi_4(q_0,q_1;q_2,q_3,q_4)\in\Ker(\mathfrak q_5)$$ 
of partial $T$-degrees $(m_0,m_1,m_2,m_3,m_4)$ such that $m_{q_0}=m_{q_1}=1$ and $m_{q_2}=m_{q_3}=m_{q_4}=2$. 

\medskip

By combining Corollary \ref{C21} with Proposition \ref{P18} (a) and (c), the equality $\varrho_0^-(2,2)=0$ (see Equation (\ref{EQ170})), and Examples \ref{EX6} and \ref{EX7}, we get directly:

\begin{cor}\label{C26}
The following four properties hold:

\medskip
{\bf (a)} For all integers $s\geq 0$ we have $D(2,1,s)=\binom{s+2}{2}$ and
$$D(2,2,s)=\binom{s+5}{5}=\frac{(s+1)(s+2)(s+3)(s+4)(s+5)}{120}.$$

\smallskip
{\bf (b)} For all integers $s\geq 0$ we have $D(2,3,s)=\sum_{i=5}^8 \binom{s+i}{8}$ and thus
$$D(2,3,s)=\frac{(s+1)(s+2)(s+3)(s+4)(s+5)(s^3+9s^2+46s+84)}{2(7!)}.$$

\smallskip
{\bf (c)} For all integers $s\geq 0$ we have 
$$D(2,4,s)=\binom{s+11}{11}+3\binom{s+10}{11}+6\binom{s+9}{11}+10\binom{s+8}{11}$$
$$=\frac{[\prod_{i=1}^8 (s+i)](20s^3+103s^2+535s+990)}{11!}.$$

\smallskip
{\bf (d)} Let $r\geq 5$ be an integer. Then for all integers $s\geq 0$, $D(2,r,s)$ is 
$$\binom{3r+s-1}{3r-1}+\binom{r-1}{2}\binom{3r+s-2}{3r-1}+\sum_{l=1}^{\varrho^-(2,r)}\binom{3r+s-1-c_l(2,r)}{3r-1}$$
and the rational number 
$$\frac{D(2,r,s)-\binom{3r+s-1}{3r-1}-\binom{r-1}{2}\binom{3r+s-2}{3r-1}}{\frac{2^{r-1}}{(r-2)}\binom{2r-3}{r}-1-\binom{r-1}{2}}$$
belongs to the interval of the real line
$$\Big[\binom{\lfloor \frac{3r}{2}\rfloor+s+1}{3r-1},\binom{3r+s-3}{3r-1}\Big].$$
Thus the unique polynomial $D_{2,r}(x)\in\mathbb Q[x]$ such that for all integers $s\geq 0$ we have $D(2,r,s)=D_{2,r}(s)$, is divisible by 
$\prod_{i=1}^{3r-1-\lfloor \frac{3(r-1)}{2}\rfloor} (x+i)$,
has leading coefficient $\frac{2^{r-1}}{(r-2)(3r-1)!}\binom{2r-3}{r}$ and satisfies $D_{2,r}(0)=1$ and $D_{2,r}(1)=\binom{r+1}{2}$.
\end{cor}

\subsection{A general result on $\sigma$-finite generation}\label{S44}

For the sake of future references we include here the following generalization of Theorem \ref{T31} and Corollary \ref{C23}. Let $d\in\mathbb N$ and let 
$$\underline{x}:=\underline{x}^{(0)}:=[x_1\cdots x_d]^{\t},\ \underline{x}':=\underline{x}^{(1)}:=[x_1'\cdots x_d']^{\t},\ldots,\underline{x}^{(r)}:=[x_1^{(r)}\cdots x_d^{(r)}]^{\t}, \ldots$$ 
be a sequence of column vectors with entries 
 indeterminates of degree $1$. We consider the polynomial $\mathbb Z_{\geq 0}$-graded $\mathbb K$-algebras
$$\mathbb X:=\cup_{r\geq 0} \mathbb X_r,\;\;\;\;\mathbb X_r:=\mathbb K[x_i^{(l)}|1\leq i\leq d,0\leq l\leq r].$$ 
We consider the shift $\mathbb K$-algebra monomorphism 
$$\sigma:\mathbb X\rightarrow \mathbb X$$
that sends $x_i^{(l)}$ to $x_{i}^{(l+1)}$ for all $(i,l)\in \{1,\ldots,d\}\times\mathbb Z_{\geq 0}$.

If $G$ is a reductive group over $\Spec(\mathbb K)$ equipped with a faithful representation $\rho:G\rightarrow\pmb{\GL}_{d,\mathbb K}$, we let $G$ act 
on $\mathbb X$ via the matrix rules $(h,\underline{x}^{(l)})\mapsto \rho(h)\underline{x}^{(l)}$ for all $(h,l)\in G(\mathbb K)\times \mathbb Z_{\geq 0}$; it restricts to an action of $G$ on $\mathbb X_r$ for each $r\in\mathbb Z_{\geq 0}$.
Clearly, if $\vartheta\in\mathbb X^G$, then $\sigma(\vartheta)\in\mathbb X^G$. 

For $j,t\in\mathbb Z_{\geq 0}$, we consider the $G$-invariant $\mathbb K$-automorphism 
$$\iota_{j;j+t}:\mathbb X\rightarrow\mathbb X$$ 
that for $i\in\{1,\ldots,d\}$ interchanges $x_i^{(j)}\longleftrightarrow x_i^{(j+t)}$ and fixes each $x_i^{(l)}$ with $l\in\mathbb Z_{\geq 0}\setminus\{j,j+t\}$; each $\iota_{j;j}$ is an identity automorphism.

\begin{thm}\label{T33}
We assume that there exists a positive integer $m$ such that we have a $G$-separating polynomial on $V^m$ (i.e., the stable locus $(V^m)^{\s}$ is non-empty).\footnote{If the identity component of $G$ is semisimple, then $m\in\{1,\ldots,d\}$ exists (see Fact \ref{F5}).} Let $\mathcal T_0$ be an arbitrary transcendence basis of $\Frac(\mathbb X_{m-1}^G)$ over $\mathbb K$ and let $\mathcal T=\{\vartheta_1,\ldots,\vartheta_d\}\subset\mathbb X_m^G$ be an arbitrary transcendence basis of $\Frac(\mathbb X_m^G)$ over $\Frac(\mathbb X_{m-1}^G)$. Then the following four properties hold:

\medskip
{\bf (a)} For each $l\in\mathbb Z_{\geq 0}$ the union $\mathcal T_0\cup\left(\bigcup_{j=0}^l\iota_{m;m+j}(\mathcal T)\right)$ is a transcendence basis of $\Frac(\mathbb X_{m+l}^G)$ over $\mathbb K$.

\smallskip
{\bf (b)} The set $\mathcal T$ is $\sigma$-algebraically independent over $\mathbb K$, i.e., the set 
$$\{\sigma^j(\vartheta_i)|(i,j)\in\{1,\ldots,d\}\times\mathbb Z_{\geq 0}\}=\cup_{j\geq 0} \sigma^j(\mathcal T)$$
of polynomials of $\mathbb X^G$ is algebraically independent over $\mathbb K$.

\smallskip
{\bf (c)} The field $\Frac(\mathbb X^G)$ is $\sigma$-generated by a finite subset of $\Frac(\mathbb X_{m+1}^G)$. 

\smallskip
{\bf (d)} We assume that for each integer $l\geq 0$ we have a direct sum decomposition
$$\mathbb X_{m+l-1}^G=\mathbb X_{m+l-1,\even}^G\oplus \mathbb X_{m+l-1,\odd}^G$$
that corresponds to homogeneous degrees in $2a\mathbb Z_{\geq 0}$ and $a(2\mathbb N-1)$ (respectively), where $a\in\mathbb N$ is fixed. We also assume that for $l>> 0$, $\mathbb X_{m+l-1,\odd}^G\neq 0$. Then $\mathbb X_{m+1,\odd}^G\neq 0$ and the field 
$$\Frac(\mathbb X^G_{\even}):=\cup_{l\geq 0} \Frac(\mathbb X_{m+l-1,\even}^G)$$ is $\sigma$-generated by a finite subset of $\Frac(\mathbb X_{m+2,\even}^G)$. 
\end{thm}

\noindent
{\it Proof.} As the representation $\rho:G\rightarrow \pmb{\GL}_{d,\mathbb K}$ is faithful, for a generic $m+1$-tuple $(v_0,\ldots,v_m)\in V^{m+1}$, the stabilizer $\Stab_G((v_0,\ldots,v_m))$ is trivial. Based on this, the proof of the theorem is the same as the proofs of Theorem \ref{T31} and Corollary \ref{C23}. We include here only details of the arguments of the last two results which used concrete notation that pertained exclusively to quadratic forms. 

As in the proof of Lemma \ref{L46} (b) and (c) one checks that for all integers $j>t\geq 0$ we have the following property:

\medskip\noindent
{\bf (i)} the natural $G$-morphism

$$\rho_{m+j,m+t}:\pmb{\V}^{m+j}\rightarrow \pmb{\V}^{m+t}\times_{(\pmb{\V}^{m+t}/G)} (\pmb{\V}^{m+j}/G)$$ 
is generically finite. For $t>0$, $\rho_{m+j,m+t}$ is generically finite of degree $1$, the field extension $\Frac(\mathbb X_{m+t-1}^G)\rightarrow\Frac(\mathbb X_{m+j-1}^G)$ is purely transcendental of transcendence degree $d(j-t)$, and a subset with $d(j-t)$ elements of $\Frac(\mathbb X_{m+j-1}^G)$ is a generating transcendence basis (resp. a transcendence basis) of $\Frac(\mathbb X_{m+j-1}^G)$ over $\Frac(\mathbb X_{m+t-1}^G)$ if and only if it is a generating transcendence basis (resp. transcendence basis) of $\Frac(\mathbb X_{m+j-1})$ over $\Frac(\mathbb X_{m+t-1})$.

\medskip
As in the proof of Theorem \ref{T31} we check that the following two properties hold:

\medskip\noindent
{\bf (ii)} for all integers $l\geq 0$, the set $\iota_{m;m+l}(\mathcal T)$ is a transcendence basis of 
$$\Frac((\mathbb K_{m-1}[x_1^{(m+l)},\ldots,x_d^{(m+l)}])^G)$$ over $\Frac(\mathbb X_{m-1}^G)$;

\smallskip\noindent
{\bf (iii)} the union $\mathcal T\cup(\bigcup_{t\geq 1} \iota_{m;m+t}(\mathcal T))$ is an algebraically independent subset of $\Frac(\mathbb X)$ over $\Frac(\mathbb X_{m-1})$ and thus it is a transcendence basis of $\Frac(\mathbb X^G)$ over $\Frac(\mathbb X_{m-1}^G)$.

\medskip
For $(t,r)\in\mathbb N\times\mathbb Z_{\geq 0}$ let $\iota_{0,1,\ldots,t-1;t+r,t+1,\ldots,2t-1+r}$ be the $G$-invariant $\mathbb K$-automorphism of $\mathbb X$ which for each $i\in\{1,\ldots,d\}$ interchanges $x_i^{(l)}\longleftrightarrow x_i^{(l+t+r)}$ for all $l\in\{0,\ldots,t-1\}$ and fixes each $x_i^{(l)}$ with $l$ in the difference set $\mathbb Z_{\geq 0}\setminus\{0,1,\ldots,t-1,t+r,t+1,\ldots,2t-1+r\}$. 

From the property (iii) we first get that part (a) holds and second that (as in the proof of Corollary \ref{C23}):

\medskip\noindent
{\bf (iv)} for each integer $t\geq m$, $\sigma^t(\mathcal T)=\iota_{0,1,\ldots,m-1;t,t+1,\ldots,t+m-1}(\iota_{m;m+t}(\mathcal T))$ is a transcendence basis of 
$\Frac(\mathbb X_{m+t}^G)$ over $\Frac(\mathbb X_{m+t-1}^G)$;

\smallskip\noindent
{\bf (v)} if $t\in\{1,\ldots,m-1\}$, then $\sigma^t(\mathcal T)=\iota_{\pi_{m+t}}(\iota_{m;m+t}(\mathcal T))$ is a transcendence basis of 
$\Frac(\mathbb X_{m+t}^G)$ over $\Frac(\mathbb X_{m+t-1}^G)$, where $\pi_{m+t}$ is the shift by $t$ permutation of $\{0,\ldots,m+t-1\}$ with the property that $\pi_{m+t}(l)$ is congruent to $l+t$ modulo $m+t$ and where $\iota_{\pi_{m+t}}$ is the $G$-invariant $\mathbb K$-automorphism of $\mathbb X$ which for each $i\in\{1,\ldots,d\}$ maps $x_i^{(l)}$ to $x_i^{(\pi_{m+t}(l))}$ for all $l\in\{0,\ldots,m+t-1\}$ and fixes $x_i^{(l)}$ for all integers $l\geq m+t$. 

\medskip
Part (b) follows from the property (iv). From the second part of the property (i) we get that the proofs of Corollary \ref{C23} (b) to (d) apply entirely to give us that parts (c) and (d) hold as well.\endproof

\begin{rem}\label{R44}
The fact that the fields in parts (c) and (d) of the theorem above are $\sigma$-finitely generated follows, again, from the general fact of difference algebra in \cite{Le}, Ch. 4, Thm. 4.4.1. The novel part of parts (c) and (d) is, again, the bounds $m+1$ and $m+2$ (respectively) for the `order' of the $\sigma$-generators.
\end{rem}

\medskip
\noindent
{\bf Acknowledgment.}
The first author acknowledges support from the Simons Foundation through grants 615356 and 311773. The second author would like to thank Binghamton University for good working conditions and Ofer Gabber for the example of Footnote 3 and for suggesting the proof of Lemma \ref{L3}.

\hbox{}
\hbox{Alexandru Buium,\;\;\;E-mail: buium@math.unm.edu}
\hbox{Address: Department of Mathematics and Statistics,} 
\hbox{University of New Mexico, Albuquerque, NM 8713, U.S.A.}

\bigskip\smallskip

\hbox{Adrian Vasiu,\;\;\;E-mail: adrian@math.binghamton.edu}
\hbox{Address: Department of Mathematical Sciences, Binghamton University,}
\hbox{P. O. Box 6000, Binghamton, New York 13902-6000, U.S.A.}

\end{document}